\documentclass[12pt]{tlpbook}
\usepackage{times,mathrsfs,amscd,amsmath,amsfonts,amssymb,graphicx,tikz,enumerate,xcolor,fancyhdr}
\usepackage[T1]{fontenc}
\usepackage{helvet}
\usepackage[most, listings, theorems]{tcolorbox}
\usepackage{babel}
\usepackage{thmtools}
\definecolor{grey}{rgb}{0.5,0.5,0.5}
\definecolor{orange}{rgb}{1,0.5,0}
\definecolor{purple}{rgb}{0.7,0,1}
\definecolor{darkgreen}{rgb}{0,0.5,0.2}
\definecolor{darkblue}{rgb}{0,0,1}
\definecolor{linkc}{rgb}{0.4,0.2,0.4}
\usepackage{hyperref}
\hypersetup{
    colorlinks=true,
    citecolor=darkblue,
    linkcolor=linkc,
    urlcolor=red,
    }

\newcounter{myquestion}
\newcommand{\atc}{\addtocounter{myquestion}{1}}

\newtcolorbox[auto counter, number within = chapter]
{theorem}[2][]{%
fonttitle= \bfseries\upshape, fontupper= \upshape,  arc=0mm,
colback=red!5,colframe=red!80!black,
title={Theorem \thetcbcounter\; #2},#1}

\newtcolorbox[use counter from = theorem, number within = chapter]
{proposition}[2][]{%
fonttitle= \bfseries\upshape, fontupper= \upshape, arc=0mm,
colback=red!5,colframe=red!80!black,
title={Proposition \thetcbcounter\; #2},#1}

\newtcolorbox[use counter from = theorem, number within = chapter]
{lemma}[2][]{%
fonttitle= \bfseries\upshape, fontupper= \upshape,
arc=0mm, colback=green!10,colframe=green!50!black,
title={Lemma \thetcbcounter\; #2},#1}

\newtcolorbox[use counter from = theorem, number within = chapter]
{corollary}[2][]{%
fonttitle= \bfseries\upshape, fontupper= \upshape, arc=0mm,
colback=green!10,colframe=green!50!black,
title={Corollary \thetcbcounter\; #2},#1}

 \newtcolorbox{myproof}[1]{%
fonttitle= \bfseries\upshape, fontupper= \upshape, arc=0mm,
colback=grey!10,colframe=grey!50!black, colbacktitle=grey!75!black,
enhanced, attach boxed title to top center={yshift=-2mm}, title={#1}}

\newtcolorbox[auto counter, number within = chapter]
{definition}[2][]{%
fonttitle= \bfseries\upshape, fontupper= \upshape, arc=0mm,
colback=blue!5,colframe=blue!75!black, 
title={Definition  \thetcbcounter\; #2},#1}

 \newtcolorbox[auto counter, number within = chapter]{remark}[2][]{%
fonttitle= \bfseries\upshape, fontupper= \upshape, arc=0mm,
colback=yellow!10,colframe=yellow!20!red,
title={Remark \thetcbcounter\; #2},#1}

\newtcolorbox[auto counter, number within = chapter]{example}[2][]{%
fonttitle= \bfseries\upshape, fontupper= \upshape, arc=0mm,
colback=orange!15,colframe=orange!55!black,
title={Example \thetcbcounter\; #2},#1}

 \newtcolorbox{example2}[1]{%
fonttitle= \bfseries\upshape, fontupper= \upshape, arc=0mm,
colback=orange!15,colframe=orange!55!black,  title=#1}

 \newtcolorbox{solution}[1]{%
fonttitle= \bfseries\upshape, fontupper= \upshape, arc=0mm,
colback=grey!10,colframe=grey!50!black, colbacktitle=grey!75!black,
enhanced, attach boxed title to top center={yshift=-2mm}, title={#1}}

 \newtcolorbox{highlight}[1]{%
fonttitle= \bfseries\upshape, fontupper= \upshape, arc=0mm,
colback=yellow!10,colframe=yellow!20!red, colbacktitle=yellow!20!red,
enhanced, attach boxed title to top center={yshift=-2mm}, title={#1}}

\newtcolorbox[auto counter, number within = section]{question}[2][]{%
fonttitle= \bfseries\upshape, fontupper= \upshape, arc=0mm,
colback=purple!5,colframe=purple!75!black, title={\linkt Question #2\linko},#1}

\definecolor{covercolor}{rgb}{0.7,0,0.2}
\definecolor{covertext}{rgb}{1,1,0.3}
\definecolor{darkblue}{rgb}{0,0,0.5}

\newcommand{\di}{\displaystyle}

\newcommand{\bookauthor}{Teo Lee Peng}
\newcommand{\booktitle}{Mathematical Analysis\\Volume I}
\newcommand{\fa}{$\hspace{0cm}$}
\newcommand{\bp}{\end{myproof}\begin{myproof}{}}
\newcommand{\bs}{\end{solution}\begin{solution}{}}
\newcommand{\be}{\end{example}\begin{example2}{}}
\newcommand{\vp}{\vfill\pagebreak}
\newcommand{\sbr}{\vspace{0.8cm}\hrule\vspace{0.5cm}}
\newcommand{\linkt}{\hypersetup{linkcolor=white}}
\newcommand{\linko}{\hypersetup{linkcolor=linkc}}

\begin{document}
  \begin{coverpage}
~\vspace{2cm}
\begin{center}
{\fontfamily{phv}\fontseries{mc}\fontsize{24}{28}\selectfont \textcolor{covertext}{  {\bfseries {\booktitle}} }

\vspace{6cm}

\textcolor{covertext}{{ \bfseries{\bookauthor}}}

}
\end{center}
 \end{coverpage}

\pagecolor{white}
\title{\uppercase{\booktitle}}

\author{\bookauthor}

\dedication{}

\date{\today}

\maketitle

 \frontmatter
\setcounter{page}{1}
\hypersetup{linkcolor=darkblue}

\tableofcontents

\chapter*{Preface}
 
Mathematical analysis is a standard course which introduces students to rigorous reasonings in mathematics, as well as the theories needed for advanced analysis courses. It is a compulsory course for all mathematics majors. It is also strongly recommended for students that major in computer science, physics, data science, financial analysis,  and other areas that require a lot of analytical skills. Some standard textbooks in mathematical analysis include the classical one by  Apostol \cite{Apostol} and Rudin \cite{Rudin}, and the modern one by Bartle \cite{Bartle}, Fitzpatrick \cite{Fitzpatrick}, Abbott \cite{Abbott}, Tao \cite{Tao_1, Tao_2} and Zorich \cite{Zorich_1, Zorich_2}.
 
 This book is the first volume of the textbooks intended for a one-year course in mathematical analysis.   We introduce the fundamental concepts in a pedagogical way. Lots of examples are given to illustrate the theories. 
We assume that students are familiar with the material of calculus such as those in the book \cite{Stewart}. Thus, we do not emphasize on the computation techniques. Emphasis is put on building up  analytical skills through rigorous reasonings. 

Besides calculus, it is also assumed that students have taken introductory courses in discrete mathematics and linear algebra, which covers topics such as logic, sets, functions, vector spaces, inner products, and quadratic forms. Whenever needed, these concepts would be briefly revised. 

In this book, we have defined all the mathematical terms we use carefully. While most of the terms have standard definitions, some of the terms may have definitions  defer from authors to authors. The readers are advised to check the definitions of the terms used in this book when they encounter them. This can be easily done by using the search function provided by any PDF viewer. The readers are also encouraged to fully utilize the hyper-referencing provided.

 \vspace{0.9cm}
~\hfill\bookauthor

%-----------------------------------------------------------------------------
% End of preface.tex
%-----------------------------------------------------------------------------

\mainmatter

\chapter{The Real Numbers}\label{ch1}

 \section{Logic, Sets and Functions}
 
In this section, we give a brief review of propositional logic, sets and functions.  It is assumed that students have taken an introductory course which covers these topics, such as a course in discrete mathematics \cite{Rosen}. 

\begin{definition}{Proposition} 
A {\bf proposition}, usually denoted by $p$, is a declarative sentence that is either true or false, but not both.
\end{definition}

\begin{definition}{Negation of a Proposition}
 If $p$ is a proposition, $\neg p$ is the {\bf negation} of $p$. The proposition $p$ is true if and only if the negation $\neg p$ is false.
\end{definition}

From two propositions $p$ and $q$, we can apply logical operators and   obtain a compound proposition. 

\begin{definition}{Conjunction of   Propositions}
 If $p$ and $q$ are propositions, $p\wedge q$ is the {\bf conjunction} of $p$ and $q$, read as "$p$ and $q$". The proposition $p\wedge q$ is true if and only if both $p$ and $q$ are true.
\end{definition}

\begin{definition}{Disjunction of   Propositions}
 If $p$ and $q$ are propositions, $p\vee q$ is the {\bf disjunction} of $p$ and $q$, read as "$p$ or $q$". The proposition $p\vee q$ is true if and only if either $p$ is true or $q$ is true.
\end{definition}

\begin{definition}{Implication of   Propositions}
 If $p$ and $q$ are propositions, the proposition $p\to q$ is read as "$p$ {\bf implies} $q$". It is false if and only if $p$ is true but $q$ is false.
\end{definition}

$p\to q$ can also be read as "if $p$ then $q$ or "$p$ only if $q$". In mathematics, we usually write $p\implies q$ instead of $p\to q$. 

\begin{definition}{Double Implication}
 If $p$ and $q$ are propositions, the proposition $p\longleftrightarrow q$ is read as "$p$ {\bf if and only if} $q$". It is the conjunction of $p\to q$ and $q\to p$. Hence, it is true if and only if both $p$ and $q$ are true, or both $p$ and $q$ are false.
\end{definition}

The stament ``$p$ if and only if $q$'' is often expressed as $p\iff q$.

Two compound propositions $p$ and $q$ are said to be logically equivalent, denoted by $p\equiv q$, provided that $p$ is true if and only if $q$ is true.

Logical equivalences are important for working with mathematical proofs. Some equivalences such as commutative law, associative law, distributive law are obvious. Other important equivalences are listed in the theorem below.
\begin{theorem}{Logical Equivalences}

Let $p$, $q$, $r$ be propositions. 
\begin{enumerate}[1.]
\item $p\to q \;\equiv \;\neg p\vee q$
\item De Morgan's Law
\begin{enumerate}[(i)]
\item $\neg(p\vee q)\;\equiv\;\neg p\,\wedge\,\neg q$
\item $\neg(p\wedge q)\;\equiv\;\neg p\,\vee\,\neg q$
\end{enumerate}

\end{enumerate}
\end{theorem}

A very important equivalence is the equivalence of an implication with its contrapositive.
\begin{theorem}{Contraposition}
If $p$ and $q$ are propositions, $p\to q$ is equivalent to $\neg q\to\neg p$.
\end{theorem}

In mathematics, we are often dealing with statements that depend on variables. Quantifiers are used to specify the extent to which such a statement is true. Two commonly used quantifiers are "for all" ($\forall$) and "there exists" ($\exists$). 

For negation of statements with quantifiers, we have the following generalized De Morgan's law.

\begin{theorem}{Generalized De Morgan's Law}
\begin{enumerate}[1.]
\item $\neg\left(\forall x\;P(x)\right)\;\equiv\;\exists x\;\neg P(x)$
\item $\neg\left(\exists x\;P(x)\right)\;\equiv\;\forall x\;\neg P(x)$
\end{enumerate}
\end{theorem}

For nested quantifiers, the ordering is important if different types of quantifiers are involved. For example, the statement
\[\forall x\;\exists y\; \;x+y=0\] is not equivalent to the statement
\[\exists y\;\forall x\;\;x+y=0.\]When the domains for $x$ and $y$ are both the set of real numbers, the first statement is true, while the second statement is false.

For a set $A$, we use the notation $x\in A$ to denote $x$ is an element of the set $A$; and the notation $x\notin A$ to denote $x$ is not an element of $A$.

\begin{definition}{Equal Sets}
Two sets $A$ and $B$ are equal if they have the same elements. In logical expression, $A=B$ if and only if
\[x\in A\iff x\in B.\]
\end{definition}

\begin{definition}
{Subset}
If $A$ and $B$ are sets, we say that $A$ is a {\bf subset} of $B$,  denoted by $A\subset B$, if every element of $A$ is an element of $B$. In logical expression, $A\subset B$ means that
\[x\in A\implies x\in B.\]
\end{definition}

When $A$ is a subset of $B$, we will also say that  $A$ is contained in $B$, or $B$ contains $A$.

We say that $A$ is a {\bf proper subset} of $B$ if $A$ is a subset of $B$ and $A\neq B$. 
In some textbooks,  the symbol "$\subseteq$" is used to denote subset, and the symbol "$\subset$" is reserved for proper subset. In this book, we will not make such a distinction. Whenever we write $A\subset B$, it means $A$ is a subset of $B$, not necessary a proper subset.

There are   operations that can be defined on sets, such as union, intersection, difference and complement. 

\begin{definition}
{Union of Sets}
If $A$ and $B$ are sets, the {\bf union} of $A$ and $B$ is the set $A\cup B$ which contains all elements that are either in $A$ or in $B$. In logical expression,
\[x\in A\cup B\iff (x\in A)\;\vee\;(x\in B).\]
\end{definition}

\begin{definition}
{Intersection of Sets}
If $A$ and $B$ are sets, the {\bf intersection} of $A$ and $B$ is the set $A\cap B$ which contains all elements that are  in both $A$ and  $B$. In logical expression,
\[x\in A\cap B\iff (x\in A)\;\wedge\;(x\in B).\]
\end{definition}

\begin{definition}
{Difference of Sets}
If $A$ and $B$ are sets, the difference of $A$ and $B$ is the set $A\setminus B$  which contains all elements that are  in  $A$ and not in  $B$. In logical expression,
\[x\in A\setminus B\iff (x\in A)\;\wedge\;(x\notin B).\]
\end{definition}

\begin{definition}
{Complement of a Set}
If $A$ is a set that is contained in a universal set $U$, the {\bf complement} of $A$ in $U$ is the set $A^C$ which contains all elements that are in $U$ but not in $A$. In logical expression,
\[x\in A^C\iff (x\in U)\;\wedge\;(x\notin A).\]
\end{definition}

Since a universal set can vary from context to context, we will usually avoid using the notation $A^C$ and use $U\setminus A$ instead for the complement of $A$ in $U$. The advantage of using the notation $A^C$ is that De Morgan's law takes a more succint form.

\begin{proposition}{De Morgan's Law for Sets}

If $A$ and $B$ are sets in a universal set $U$, and $A^C$ and $B^C$ are their complements in $U$, then
\begin{enumerate}[1.]
\item
$(A\cup B)^C=A^C\cap B^C$
\item $(A\cap B)^C=A^C\cup B^C$
\end{enumerate}
\end{proposition}

\begin{definition}{Functions}
When $A$ and $B$ are sets, a {\bf function} $f$ from $A$ to $B$, denoted by $f:A\rightarrow B$, is a correspondence that assigns every element of $A$ a unique element in $B$. If $a$ is in $A$, the {\bf image} of $a$ under the function $f$ is denoted by $f(a)$, and it is an element of $B$.

$A$ is called the {\bf domain} of $f$, and $B$ is called the {\bf codomain} of $f$.
\end{definition}
\begin{definition}
{Image of a Set}
If $f:A\to B$ is a function and $C$ is a subset of $A$, the image of $C$ under $f$ is the set 
\[f(C)=\left\{f(c)\,|\,c\in C\right\}.\]
$f(A)$ is called the range of $f$.

\end{definition}

\begin{definition}
{Preimage of a Set}
If $f:A\to B$ is a function and $D$ is a subset of $B$, the preimage of $D$ under $f$ is the set 
\[f^{-1}(D)=\left\{a\in A\,|\,f(a)\in D\right\}.\]
 
\end{definition} 
Notice that $f^{-1}(D)$ is a notation, it does not mean that the function $f$ has an inverse.

Next, we turn to discuss injectivity and surjectivity of functions.
\begin{definition}{Injection}
We say that a function $f:A\to B$ is an {\bf injection}, or the function $f:A\rightarrow B$ is {\bf injective}, or the function $f:A\to B$ is {\bf one-to-one}, if no pair of distinct elements of $A$ are mapped to the same element of $B$.
Namely,
\[a_1\neq a_2 \implies f(a_1)\neq f(a_2).\]

\end{definition}

Using contrapositive, a function is injective provided that
\[f(a_1)=f(a_2)\implies a_1=a_2.\]

\begin{definition}{Surjection}
We say that a function $f:A\to B$ is a {\bf surjection}, or the function $f:A\rightarrow B$ is {\bf surjective}, or the function $f:A\to B$ is {\bf onto}, if every element of $B$ is the image of some element in $A$.
Namely,
\[\forall b\in B, \exists a\in A, f(a)=b.\]
Equivalently, $f:A\to B$ is surjective if the range of $f$ is $B$. Namely, $f(A)=B$.

\end{definition}

\begin{definition}{Bijection}
We say that a function $f:A\to B$ is a {\bf bijection}, or the function $f:A\rightarrow B$ is bijective, if it is both injective and surjective.

A bijection is also called a  {\bf one-to-one correspondence}.

\end{definition}

Finally, we would like to make a remark about some notations. If $f:A\rightarrow B$ is a function with domain $A$, and $C$ is a subset of $A$, the restriction of $f$ to $C$ is the function $f|_C:C\rightarrow B$ defined by $f|_C(c)=f(c)$ for all $c\in C$. When no confusion arises, we will often denote this function simply as $f:C\rightarrow B$.
\vp
\section{The Set of Real Numbers and Its Subsets}\label{sec1.2}
In this section, we introduce the set of real numbers using an intuitive approach. 

\begin{definition}{Natural Numbers}
The set of {\bf natural numbers} $\mathbb{N}$ is the set that contains the counting numbers,    1, 2, 3 $\ldots$, which are also called positive integers. 
\end{definition}

$\mathbb{N}$ is an inductive set. The number 1 is the smallest element of this set. If $n$ is a natural number, then $n+1$ is also a natural number.  

The number 0 corresponds to nothing. 

For every positive integer $n$, $-n$ is a number which   produces 0 when adds to $n$. This number $-n$ is called the negative of $n$, or the additive inverse of $n$.

$-1$, $-2$, $-3$, $\ldots$, are called negative integers.

\begin{definition}{Integers}
The set of {\bf integers} $\mathbb{Z}$ is the set that contains all positive integers, negative integers and 0.
\end{definition}

We will also use the notation $\mathbb{Z}^+$ to denote the set of positive integers. 

\begin{definition}{Rational Numbers}
 The set of {\bf rational numbers} $\mathbb{Q}$ is the set defined as
\[\mathbb{Q}=\left\{\left.\frac{m}{n}\,\right|\,m,n\in\mathbb{Z}, n\neq 0\right\}.\]
\end{definition}

Each rational number is a quotient of two integers, where the denominator is nonzero. The set of integers $\mathbb{Z}$ is a subset of the set of rational numbers $\mathbb{Q}$.

Every rational number $m/n$ has a decimal expansion. For example, 
\[-\frac{23}{4}=-5.75,\]
\[\frac{27}{7}=3.857142857142\ldots=3.\dot{8}5714\dot{2}.\]
The decimal expansion of a rational number is either finite or periodic.

\begin{definition}{Real Numbers}
The set of {\bf real numbers} $\mathbb{R}$ is intuitively defined to be the set that contains all decimal numbers, which is not necessary periodic. 
\end{definition}

The set of real numbers contains the set of rational numbers $\mathbb{Q}$ as a subset. If a real number is not a rational number, we call it an {\bf irrational number}. The set of irrational numbers is $\mathbb{R}\setminus\mathbb{Q}$.

\begin{figure}[ht]
\centering
\includegraphics[scale=0.2]{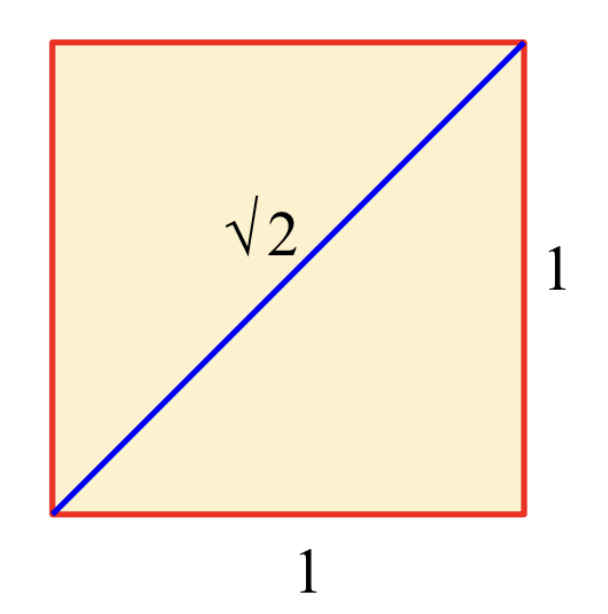}

\caption{The number $\sqrt{2}$.}\label{figure1}
\end{figure}

It has been long known that there are real numbers that are not   rational numbers. The best example is the number $\sqrt{2}$, which appears as the length of the diagonal of a unit square (see Figure \ref{figure1}).

The   addition and multiplication operations defined on the set of natural numbers can be extended to the set of real numbers consistently. 

If $a$ and $b$ are  real numbers, $a+b$ is the addition of  $a$ and $b$, and $ab$ is the multiplication of $a$ and $b$.

If $a$ and $b$ are positive real numbers, $a+b$ and $ab$ are also positive real numbers.

The set of real numbers with the addition and multuplication operations is a field, which you will learn in abstract algebra. These operations satisfy the following properties.

 \begin{highlight}{Properties of Real Numbers}

         \begin{enumerate} [1.]
\item
\textbf{Commutativity of Addition}
\[a+b=b+a\]
\item \textbf{Associativity of Addition}
\[(a+b)+c=a+(b+c)\]
\item\textbf{Additive Identity }
\[a+0=0+a=a\]
0 is called the additive identity.
\item\textbf{Additive Inverse}
\\For every real number $a$, the negative of $a$, denoted by $-a$, satisfies
\[a+(-a)=(-a)+a=0\]

\item
\textbf{Commutativity of Multiplication}
\[ab=ba\]
\item \textbf{Associativity of Multiplication}
\[(ab)c=a(bc)\]

\item\textbf{Multiplicative Identity }
\[a\cdot 1=1\cdot a=a\]
1 is called the multiplicative identity.
\item\textbf{Multiplicative Inverse}
\\For every nonzero real number $a$, the reciprocal of $a$, denoted by $1/a$, satisfies
\[a\cdot \frac{1}{a}=\frac{1}{a}\cdot a=1\]
\item \textbf{Distributivity}
\[a(b+c)=ab+ac\]
\end{enumerate}
 \end{highlight}

\vfill\pagebreak
 
The set of complex numbers $\mathbb{C}$ is the set that contains all numbers of the form $a+ib$, where $a$ and $b$ are real numbers, and $i$ is the purely imaginary number such that $i^2=-1$. It contains the set of real numbers $\mathbb{R}$ as a subset. Addition and multiplication can be extended to the set of complex numbers. These two operations on complex numbers also satisfy all the properties listed above. Nevertheless, we shall focus on the set of real numbers in this course.

There are special subsets of real numbers which are called {\bf intervals}.  There are nine types of intervals, four types are finite, five types are semi-infinite or infinite.  Their definitions are as follows. 

\begin{highlight}{Finite Intervals}

\begin{enumerate}
\item[1.] $(a,b)=\di\left\{x\in\mathbb{R}\,|\, a<x<b\right\}$
\item[2.] $[a,b)=\di\left\{x\in\mathbb{R}\,|\, a\leq x<b\right\}$
\item[3.] $(a,b]=\di\left\{x\in\mathbb{R}\,|\, a<x\leq b\right\}$
\item[4.] $[a,b]=\di\left\{x\in\mathbb{R}\,|\, a\leq x \leq b\right\}$
\end{enumerate}
\end{highlight}

For the intervals $(a, b)$, $[a, b)$, $(a, b]$, $[a, b]$, the points $a$ and $b$ are the {\bf end points} of the interval, while any point $x$ with $a<x<b$ is an {\bf interior point}.

\begin{highlight}{Semi-Infinite or Infinite Intervals}
\begin{enumerate}
\item[5.] $(a,\infty)=\di\left\{x\in\mathbb{R}\,|\, x>a\right\}$
\item[6.] $[a,\infty)=\di\left\{x\in\mathbb{R}\,|\, x\geq a\right\}$
\item[7.] $(-\infty, a)=\di\left\{x\in\mathbb{R}\,|\, x<a\right\}$
\item[8.] $(-\infty, a]=\di\left\{x\in\mathbb{R}\,|\, x\leq a\right\}$
\item[9.] $(-\infty, \infty)=\mathbb{R}$.
\end{enumerate}
\end{highlight}

  For  the intervals $(a, \infty)$, $[a, \infty)$, $(-\infty, a)$ and $(-\infty, a]$, $a$ is the {\bf end point} of the interval, while any other points in the interval besides $a$ is an {\bf interior point}.

The set of natural numbers is a well-ordered set. Every nonempty subset of positive integers has a smallest element.  This statement is equivalent to the principle of mathematical induction, which is one of the  important strategies in proving mathematical statements.

\begin{proposition}{Principle of Mathematical Induction}
Let $P(n)$ be a sequence of statements that are indexed by the set of positive integers $\mathbb{Z}^+$. Assume that the following two assertions are true.
\begin{enumerate}[1.]
\item The statement $P(1)$ is true.
\item For every positive integer $n$, if the statement $P(n)$ is true, the statement $P(n+1)$ is also true.
\end{enumerate}Then we can conclude that for all positive integers $n$, the statement $P(n)$ is true.

\end{proposition}

Before ending this section, let us discuss the absolute value and some useful inequalities.
\begin{definition}{Absolute Value}
Given a real number $x$, the {\bf absolute value} of $x$, denoted by $|x|$, is defined to be the nonnegative number
\[|x|=\begin{cases}x,\quad &\text{if}\;x\geq 0,\\-x,\quad &\text{if}\;x<0.\end{cases}\]
In particular, $|-x|=|x|$.
\end{definition}

For example, $|2.7|=2.7$, $|-2.7|=2.7$. 

The absolute value  $|x|$ can be interpreted as the distance between the number $x$ and the number $0$ on the number line. For any two real numbers $x$ and $y$, $|x-y|$ is the distance between $x$ and $y$. Hence, the  absolute value can be used to express an interval.
\begin{highlight}{Intervals Defined by Absolute Values}
Let $a$ be a real number.
\begin{enumerate}[1.]
\item
If $r$ is a positive number, 
\[|x-a|<r\iff -r<x-a<r \iff x\in (a-r, a+r).\]
\item 
If $r$ is a nonnegative number, 
\[|x-a|\leq r\iff  -r\leq x-a\leq r \iff x\in [a-r, a+r].\]
\end{enumerate}
\end{highlight}

Absolute values behave well with respect to multiplication operation. 
\begin{proposition}{}
Given real numbers $x$ and $y$,
\[|xy|=|x||y|.\]
\end{proposition}

In general, $|x+y|$ is not equal to $|x|+|y|$. Instead, we have an inequality, known as the triangle inequality, which is very important in analysis.

\begin{proposition}{Triangle Inequality}
Given real numbers $x$ and $y$, 
\[|x+y|\leq |x|+|y|.\]
\end{proposition}
This is proved by discussing all four  possible cases where $x\geq 0$ or $x<0$, $y\geq 0$ or $y<0$. 

A common mistake students tend to make is to replace both plus signs in the triangle equality directly by   minus signs. This is totally assurd. The correct one is
\[|x-y|\leq |x|+|-y|=|x|+|y|.\]

For the inequality in the other direction, we have
\begin{proposition}{}
Given real numbers $x$ and $y$, 
\[|x-y|\geq \left||x|-|y|\right|.\]
\end{proposition}
\begin{myproof}{Proof}
Since $|x-y|\geq 0$, the statement is equivalent to
\[-|x-y|\leq |x|-|y|\leq |x-y|.\]
By triangle inequality,
\[|x-y|+|y|\geq |x-y+y|=|x|.\]
Hence,
\[ |x|-|y|\leq |x-y|.\]
By triangle inequality again,
\[|x-y|+|x|=|y-x|+|x|\geq |y-x+x|=|y|.\]
Hence,
\[-|x-y|\leq |x|-|y|.\]
This completes the proof.

\end{myproof}

\begin{example}{}
If $|x-5|\leq 2$, show that 
\[9\leq x^2\leq 49.\]

\end{example}
\begin{solution}
{Solution}
$|x-5|\leq 2$ implies $3\leq x\leq 7$.
This means that $x$ is positive. The inequality $x\geq 3$ then  implies that $x^2\geq 9$, and the inequality $x\leq 7$ implies that $x^2\leq 49$. Therefore,
\[9\leq x^2\leq 49.\]
\end{solution}

Finally, we have the useful Cauchy's inequality.
\begin{proposition}{Cauchy's Inequality}
For any real numbers $a$ and $b$,
\[ab\leq \frac{a^2+b^2}{2}.\]
\end{proposition}
\begin{myproof}{Proof}
This is just a consequence of $(a-b)^2\geq 0$.
\end{myproof}

An immediate consequence of Cauchy's inequality is the arithmetic mean-geometric mean inequality. For any nonnegative numbers $a$ and $b$, the geometric mean of $a$ and $b$ is 
$\sqrt{ab}$, and the arithmetic mean is $\di \frac{a+b}{2}$. 
\begin{proposition}{}
 If $a\geq 0$, $b\geq 0$, then
\[\sqrt{ab}\leq \frac{a+b}{2}.
\]
\end{proposition}

\vp
\noindent
{\bf \large Exercises  \thesection}
\setcounter{myquestion}{1}

 \begin{question}[label=Q23020502]{\themyquestion}
Use induction to show that for any positive integer $n$,
 \[n!\geq 2^{n-1}.\]
 \end{question}
\atc

\begin{question}[label=Q23020501]{\themyquestion:\;Bernoulli's Inequality}
 Given that $a>-1$, use induction to show that 
 \[(1+a)^n\geq 1+na\]for all positive integer $n$.
 \end{question}
 \atc
\begin{question}{\themyquestion} 
Let $n$ be a positive integer.  If $c_1, c_2, \ldots, c_n$ are numbers that lie in the interval $(0,1)$, show that
\[(1-c_1)(1-c_2)\ldots (1-c_n)\geq 1-c_1-c_2-\cdots-c_n.\]
\end{question}
  
\vp
\section{Bounded Sets and the Completeness Axiom}
\label{sec1.3}

In this section, we discuss a property of real numbers called completeness. The set of rational numbers does not have this property.

First, we introduce the concept of boundedness.
\begin{definition}{Boundedness}
Let $S$ be a subset of $\mathbb{R}$.
\begin{enumerate}[1.]
\item
We say that $S$ is {\bf bounded above} if there is a number $c$ such that 
\[ x\leq c\quad \text{for all}\; x\in S.\]
Such a $c$ is called an upper bound of $S$.
\item We say that $S$ is {\bf bounded below} if there is a number $b$ such that
\[x\geq b\quad \text{for all}\; x\in S.\]
Such a $b$ is called a lower bound of $S$.
\item We say that $S$ is {\bf bounded} if it is bounded above and bounded below. In this case, there is a number $M$ such that
\[|x|\leq M\quad \text{for all}\; x\in S.\]
\end{enumerate}
\end{definition}

Let us look at some examples.
\begin{example}[label=23020705]{}
Determine whether each of the following sets of real numbers is bounded above, whether it is bounded below, and whether it is bounded.
\begin{enumerate}[(a)]
\item 
$\di A=\left\{x\,|\, x<2\right\}$
\item  $\di B=\left\{x\,|\, x>-2\right\}$
\item $\di C=\left\{x\,|\, -2<x<2\right\}$. 
\end{enumerate}
\end{example}

\begin{solution}{Solution}
\begin{enumerate}[(a)]
\item
The set $A$ is bounded above since every element of $A$ is less than or equal to 2. It is not bounded below, and so it is not bounded.

\item The set $B$ is bounded below since every element of $B$ is larger than or equal to $-2$. It is not bounded above, and so it is not bounded.

\item The set $C$ is equal to $A\cap B$. So it is bounded above and bounded below. Therefore, it is bounded.
\end{enumerate}
\end{solution}

\begin{figure}[ht]
\centering
\includegraphics[scale=0.2]{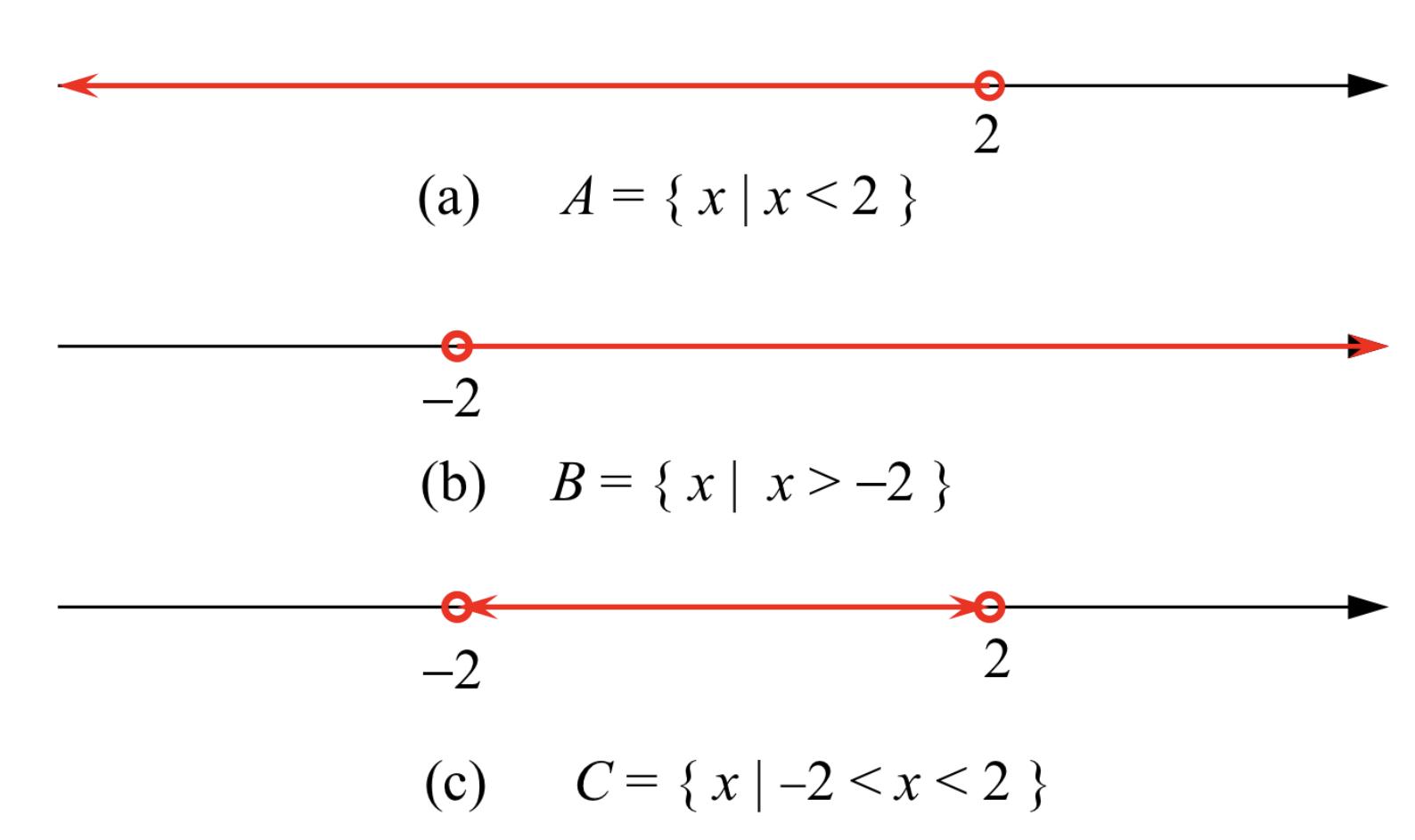}

\caption{  The sets $A$, $B$, $C$ in Example \ref{23020705}.}\label{figure5}
\end{figure}

If $S$ is a set of real numbers, the {\bf negative} of $S$, denoted by $-S$, is the set
\[-S=\left\{-x\,|\,x\in S\right\}.\]
For example, the set $B=\left\{x\,|\, x>-2\right\}$ is the negative of the set $\di A=\left\{x\,|\, x<2\right\}$, the set $\di C=\left\{x\,|\, -2<x<2\right\}$ is the negative of itself (see Figure \ref{figure5}). It is obvious that $S$ is bounded above if and only if $-S$ is bounded below.

Next, we recall the definition of maximum and minimum of a set.
\begin{definition}{Maximum and Minimum} Let $S$ be a nonempty subset of real numbers. 
\begin{enumerate}[1.]
\item  A number $c$ is called the {\bf largest element} or {\bf maximum} of $S$ if $c$ is an element of $S$ and
\[x\leq c\hspace{1cm}\text{for all}\;x\in S.\]
If the maximum of the set $S$ exists, we denote is by $\max S$.
\item A number $b$ is called the {\bf smallest element} or  {\bf minimum} of $S$ if $b$ is an element of $S$ and
\[x\geq b\hspace{1cm}\text{for all}\;x\in S.\]
If the minimum of the set $S$ exists, we denote it by $ \min S$.
\end{enumerate}
\end{definition}

Obviously, $b$ is the maximum of a set $S$ if and only if $-b$ is the minimum of the set $-S$.

\begin{example}{}
For the set $S_1=[-2,2]$, $-2$ is the minimum, and $2$ is the maximum.
\\
For the set $S_2=[-2, 2)$, $-2$ is the minimum, and there is no maximum.
\end{example}

This example shows that a bounded set does not necessarily have maximum or minimum.
However, a finite set always have a maximum and a minimum.

\begin{proposition}{}

If $S$ is a finite set, then $S$ has a maximum and a minimum.
\end{proposition}

Next, we introduce the concept of least upper bound.
\begin{definition}{Least Upper Bound}
Let $S$ be a nonempty subset of real numbers
 that is bounded above, and let $U_S$ be the set of upper bounds of $S$. Then $U_S$ is a nonempty set that is bounded below. If $U_S$ has a smallest element $u$, we say that $u$ is the {\bf least upper bound} or {\bf supremum} of $S$, and denote it by
\[u=\sup S.\]

\end{definition}
\begin{example}{}
For the sets $S_1=[-2,2]$ and  $S_2=[-2, 2)$, \[\sup S_1=\sup S_2=2.\]
\end{example}

Notice that $\sup S$, if exists, is not necessary an element of $S$.  The following proposition depicts the relation between the maximum of a set (if exists) and its least upper bound.
\begin{proposition}[label=23020606]{Supremum and Maximum}
Let $S$ be a nonempty subset of real numbers. Then $S$ has a maximum if and only if $S$ is bounded above and $\sup S$ is in $S$.
\end{proposition}

One natural question to ask is, if $S$ is a nonempty subset of real numbers that is bounded above, does $S$ necessarily have a least upper bound. The completeness axiom   asserts that this is true.
\begin{highlight}{Completeness Axiom}
If $S$  is a nonempty subset of real numbers that is bounded above, then $S$ has a least upper bound.
\end{highlight}

The reason this is formulated as an axiom is we cannot prove this   from our intuitive definition of real numbers. Therefore, we will assume this as a fact for the set of real numbers. A lots of theorems that we  are going to derive later is a consequence of this axiom. 

Actually, the set of real numbers can be constructed axiomatically, taken it to be a set that contains the set of rational numbers, satisfying all properties of addition and multiplication operations, as well as   the completeness axiom. However, this is a tedious construction and will drift us too far. 

 To show that the completeness axiom is not completely trivial, we  show in Example \ref{ex23020101} that if we only consider the set of rational numbers, we can find a subset of rational numbers $A$ that is bounded above but does not   have a least upper bound in the set of rational numbers. We look at the following example first.

\begin{example}[label=23020701]{}
Define the set of real numbers $S$ by
\[S=\left\{x\in \mathbb{R}\,|\, x^2<2\right\}.\]
Show that $S$ is nonempty and is bounded above. Conclude that the set 
\[A=\left\{x\in \mathbb{Q}\,|\, x^2<2\right\}\] is also nonempty and is bounded above by a rational number.
\end{example}

\begin{solution}{Solution}
 The number 1 is  in $S$, and so $S$ is nonempty.
For any $x\in S$, 
$x^2<2<4$, and hence
$x<2$.
This shows that $S$ is bounded above by 2.  
Since 1 and 2 are rational numbers, the same reasoning shows that the set $A$ is nonempty and is bounded above by a rational number.
\end{solution}
 
\begin{example}[label=ex23020101]{}
Consider the set 
\[A=\left\{x\in\mathbb{Q}\,|\, x^2<2\right\}.\]
By Example \ref{23020701},  $A$ is a nonempty subset of rational numbers that is bounded above by 2.
Let $U_A$ be the set of upper bounds of $A$ in $\mathbb{Q}$. Namely,
\[U_A=\left\{c\in\mathbb{Q}\,|\, x\leq c\;\text{for all}\;x\in A\right\}.\]
Show that $U_A$ does not have a smallest element.
\end{example}
\begin{solution}{Solution}
We use proof by contradiction. Assume that $U_A$ has a smallest element $c_1$, which is an upper bound of $A$ that is smaller than or equal to any upper bound of $A$. Then for any $x\in A$, 
\[x^2\leq c_1.\] 
 Since $1$ is in $A$,  $c_1$ is a positive rational number. Hence, there are poitive integers $p$ and $q$ such that
\[c_1=\frac{p}{q}.\]

Since there are no rational numbers whose square is 2, we must have either $c_1^2<2$ or $c_1^2>2$.  
 
Define the positive rational  number $c_2$ by
\[c_2=\frac{2p+2q}{p+2q}.\]

Notice that
\[c_1-c_2=\frac{p(p+2q)-q(2p+2q)}{q(p+2q)}=\frac{p^2-2q^2}{q(p+2q)},\]

\[c_1^2-2=\frac{p^2-2q^2}{q^2},\]
and 
\[c_2^2-2=\frac{4p^2+8pq+4q^2-2(p^2+4pq+4q^2)}{(p+2q)^2}=\frac{2(p^2-2q^2)}{(p+2q)^2}.\]

\textbf{Case 1: } $c_1^2<2$.\\
In this case, $p^2<2q^2$. It follows that $c_1<c_2$ and $c_2^2<2$. But then $c_1$ and $c_2$ are both in $A$, and $c_2$ is an element in $A$ that is larger than $c_1$, which contradicts to $c_1$ is an upper bound of $A$. Hence, we cannot have $c_1^2<2$.
\bs

\textbf{Case 2: } $c_1^2>2$.\\
In this case, $p^2>2q^2$. It follows that $c_1>c_2$ and $c_2^2>2$. 
Since $c_2^2>2$, we find that for any $x\in A$, \[x^2<2<c_2^2.\]Thus, 
\[-c_2<x<c_2.\] In particular, $c_2$ is also an upper bound of $A$. Namely, $c_2$ is in $U_A$.
But then $c_1$ and $c_2$ are both in $U_A$ and $c_1>c_2$. This contradicts to $c_1$ is the smallest element in $U_A$. Hence, we cannot have $c_1^2>2$.

Since both Case 1 and Case 2 lead to contradictions, we conclude that $U_A$ does not have a smallest element.
\end{solution}

In the solution above, the construction of the positive rational number $c_2$ seems a bit adhoc. In fact, we can define $c_2$ by
\[c_2=\frac{m p+2n q}{n p+m q}\]for any positive integers $m$ and $n$ with $m^2>2n^2$. Then the proof still works.

Now let us see how completeness axiom is used to guarantee that there is a real number whose square is 2. 

\begin{example}[label=23021011]{}
Use completeness axiom to show that there is a positive real number $c$ such that 
\[c^2=2.\]
\end{example}

\begin{solution}{Solution}
Define the set of real numbers $S$ by
\[S=\left\{x\in \mathbb{R}\,|\, x^2<2\right\}.\] Example \ref{23020701} asserts that $S$ is a nonempty subset of real numbers that is bounded above.
Completeness axiom  asserts that $S$ has a least upper bound $c$.\bs
Since $1$ is in $S$, $c\geq 1$.
We are going to prove that $c^2=2$   using proof by contradiction. If $c^2\neq 2$, then $c^2<2$ or $c^2>2$.

\textbf{Case 1:} $c^2<2$. \\
Let $d=2-c^2$. Then $0<d\leq 1$. Define the number $c_1$ by
\[c_1= c+\frac{d}{4c}.\]
Then $c_1>c$, and
\[c_1^2=c^2+\frac{d}{2}+\frac{d^2}{16c^2}\leq c^2+\frac{d}{2}+\frac{d}{16}<c^2+d=2.\]
This implies that $c_1$ is an element of $S$ that is larger than $c$, which contradicts to $c$  is an upper bound of $S$.

\textbf{Case 2:} $c^2>2$. \\
Let $d=c^2-2$. Then $d>0$. Define the number $c_1$ by
\[c_1= c-\frac{d}{2c}.\]
Then $c_1<c$, and
\[c_1^2=c^2-d+\frac{d^2}{4c^2}>c^2-d=2.\]
This implies that $c_1$ is an upper bound of $S$ that is smaller than $c$, which contradicts to $c$ is the least upper bound of $S$.

Since we obtain a contradiction if $c^2\neq 2$, we must have $c^2=2$.
\end{solution}
In fact, the completeness axiom can be used to show that for any positive real number $a$, there is a positive  real number $c$ such that 
\[c^2=a.\] We denote this number $c$ as $\sqrt{a}$, called the positive square root of $a$. The number $b=-\sqrt{a}$ is another real number such that $b^2=a$.

More generally, if $n$ is a positive integer, $a$ is a positive real number, then there is a positive real number $c$ such that $c^n=a$. We denote this number $c$ by
\[c=\sqrt[n]{a},\]called the positive $n^{\text{th}}$-root of $a$.

Using the interplay between a set and its negative, we can define the greatest lower bound of a set that is bounded below. 
 
\begin{definition}{Greatest Lower Bound}
Let $S$ be a nonempty subset of real numbers
 that is bounded below, and let $L_S$ be the set of lower bounds of $S$. Then $L_S$ is a nonempty set that is bounded above. If $L_S$ has a largest element $\ell$, we say that $\ell$ is the {\bf greatest lower bound} or {\bf infimum} of $S$, and denote it by
\[\ell=\inf S.\]

\end{definition}

From the completeness axiom, we have the following.
\begin{theorem}{}
If $S$ is a nonempty subset of real numbers
 that is bounded below, then $S$ has a greatest lower bound.
\end{theorem}

For a nonempty set $S$ that is bounded, it has a least upper bound $\sup S$ and a greatest lower bound $\inf S$. The following is quite obvious.

\begin{proposition}{}
If $S$ is a bounded nonempty subset of real numbers, it has a least upper bound $\sup S$ and a greatest lower bound $\inf S$. Moreover,
\[\inf S\leq \sup S,\]
and $\inf S=\sup S$ if and only if $S$ contains exactly one element.
\end{proposition}

Let us emphasize again the characterization of the least upper bound and greatest lower bound of a set.
\begin{highlight}{Characterization of Supremum and Infimum}
Let $S$ be a nonempty subset of real numbers, and let $a$ be a real number.
\begin{enumerate}[1.]
\item$a=\sup S$ if and only if the following two conditions are satisfied.
\begin{enumerate}[(i)]
\item For all $x\in S$, $x\leq a$.
\item If $b$ is a real number such that $x\leq b$ for all $x\in S$, then $a\leq b$.

\end{enumerate}
\item $a=\inf S$ if and only if the following two conditions are satisfied.
\begin{enumerate}[(i)]
\item For all $x\in S$, $x\geq a$.
\item If $b$ is a real number such that $x\geq b$ for all $x\in S$, then $a\geq b$.

\end{enumerate}
\end{enumerate}
\end{highlight}

\begin{example}{}
For each of the following set of real numbers, determine whether it has a least upper bound, and whether it has a greatest lower bound.  
\begin{enumerate}[(a)]
\item $A=\left\{x\in\mathbb{R}\,|\, x^3<2\right\}$
\item $B=\left\{x\in\mathbb{R}\,|\, x^2<10\right\}$.
\end{enumerate}
\end{example}
\begin{solution}{Solution}
\begin{enumerate}[(a)]
\item The set $A$ is bounded above, since if $x\in A$, then $x^3<2<2^3$, and so $x<2$. The set $A$ is not bounded below since it contains all negative numbers. Hence, $A$ has a least upper bound, but it does not have a greatest lower bound.
\item If $x^2<10$, then $x^2<16$, and so $-4<x<4$. This shows that $B$ is bounded. Hence, $B$ has a least upper bound, and a greatest lower bound.
\end{enumerate}
\end{solution}

Finally, we want to highlight again Proposition \ref{23020606} together with its lower bound versus infimum counterpart.
\begin{highlight}{Existence of Maximum and Minimum}
Let $S$ be a nonempty subset of real numbers. 
\begin{enumerate}[1.]
\item $S$ has a maximum if and only if $S$ is bounded above and $\sup S$ is in $S$.
\item $S$ has a minimum if and only if $S$ is bounded below and $\inf S$ is in $S$.
\end{enumerate}
\end{highlight}
\vp
 
\noindent
{\bf \large Exercises  \thesection}
\setcounter{myquestion}{1}

 \begin{question}{\themyquestion}
For each of the following sets of real numbers, find its least upper bound,   greatest lower bound,  maximum, and minimum if any of these  exists. If any of these does not exist, explain why.
\begin{enumerate}[(a)]
\item
$A=(-\infty, 20)$
\item $B=[-3, \infty)$
\item $C=[-10, -2)\cup (1, 12]$
\item $D=[-2, 5]\cap (-1, 7]$
\end{enumerate}
\end{question}

\atc
 \begin{question}{\themyquestion}
 Use completeness axiom to show that there is a positive real number $c$ such that 
\[c^2=5.\]
\end{question}
\atc
 \begin{question}{\themyquestion}
For each of the following set of real numbers, determine whether it has a least upper bound, and whether it has a greatest lower bound.  
\begin{enumerate}[(a)]
\item $A=\left\{x\in\mathbb{R}\,|\, x^3>10\right\}$
\item $B=\left\{x\in\mathbb{R}\,|\, x^2<2020\right\}$.
\end{enumerate}
\end{question}

\vp
\section{Distributions of Numbers }\label{sec1.4}
In this section, we consider additional properties of the set of integers, rational numbers and real numbers.

We start by a proposition about distribution of integers.
\begin{proposition}
{}
\begin{enumerate}[1.]
\item If $n$ is an integer, there is no integer in the interval $(n, n+1)$.
\item For any real number $c$, there is exactly one integer in the interval $[c, c+1)$, and there is exactly one integer in the interval $(c, c+1]$.

\end{enumerate}
\end{proposition}
These statements are quite obvious. For any real number $c$, the integer in the interval $[c, c+1)$ is  $\lceil c\rceil$, called the {\bf ceiling} of $c$. It is the smallest integer larger than or equal to $c$. For example $\lceil -2.5\rceil=-2$, $\lceil -3\rceil =-3$. The integer in the interval $(c, c+1]$ is $\lfloor c\rfloor+1$, where $\lfloor c\rfloor$ is the {\bf floor} of $c$.  It is the largest integer that is less than or equal to $c$. For example, $\lfloor -2.5\rfloor =-3$, $\lfloor -3\rfloor =-3$.

In Section \ref{sec1.3}, we have seen that a nonempty subset of real numbers that is bounded above does not necessary have a maximum. Example \ref{ex23020101} shows that a nonempty  subset of rational numbers that is bounded above also does not necessary have a maximum. However, for nonempty  subsets of integers, the same is not true.

\begin{proposition}{}
Let $S$ be a nonempty subset of integers.
\begin{enumerate}[1.]
\item If $S$ is bounded above, it has a maximum.
\item If $S$ is bounded below, it has a minimum.
\end{enumerate}
\end{proposition}

The two statements are equivalent, and the second statement is a generalization of the well-ordered principle for the set of positive integers. It can be proved using mathematical induction.

Next we discuss another important property called the Archimedean property. First let us show that the set of positive integers $\mathbb{Z}^+$ is not bounded above. 

\begin{theorem}
[label=thm23020202]{}The set of positive integers $\mathbb{Z}^+$ is not bounded above. 
\end{theorem}
\begin{myproof}{Proof}
Assume to the contrary that the set of positive integers $\mathbb{Z}^+$ is bounded above. By completeness axiom, it has a least upper bound $u$. 
Since $u-1<u$, $u-1$ is not an upper bound of $\mathbb{Z}^+$. Hence, there is a positive integer $n$ such that
\[n>u-1.\]
It follows that
\[n+1>u.\]
 Since $n+1$ is also a positive integer, this says that there is an element of $\mathbb{Z}^+$ that is larger than the least upper bound of $\mathbb{Z}^+$. This contradicts to the definition of least upper bound. Hence, $\mathbb{Z}^+$ cannot be bounded above.
\end{myproof}
The proof uses the key fact  that any number that is smaller than the least upper bound of a set is not an upper bound of the set. This is a standard technique in proofs.

\begin{theorem}{The Archimedean Property}
\begin{enumerate}[1.]
\item
For any positive number $M$, there is a positive integer $n$ such that $n> M$.
\item For any positive number $\varepsilon$, there is a positive integer $n$ such that $1/n<\varepsilon$.
\end{enumerate}
\end{theorem}
These two statements are equivalent, and the first statement is equivalent to the fact that the set of positive integers is not bounded above.
 
 In the following, we consider another property called denseness.
\begin{definition}{Denseness}
Let $S$ be a subset of real numbers. We say that $S$ is {\bf dense} in $\mathbb{R}$ if every open interval $(a, b)$ contains an element of $S$.
\end{definition}

A key fact we want to prove is that the set of rational numbers $\mathbb{Q}$ is dense in the set of real numbers. 

\begin{theorem}{Denseness of the Set of Rational Numbers}
The set of rational numbers $\mathbb{Q}$ is dense in the set of real numbers $\mathbb{R}$.

\end{theorem}
\begin{myproof}{Proof}
Let $(a, b)$ be an open interval. Then $\varepsilon=b-a>0$. By the Archimedean property, there is a positive integer $n$ such that
 $1/n <\varepsilon$.  Hence, \[nb-na=n\varepsilon>1,\]and so
\[ na+1<nb.\]
Consider the interval $(na, na+1]$. There is an integer $m$ that lies in this interval. In other words,
\[na<m\leq na+1<nb.\]
Dividing by $n$, we have
\[a<\frac{m}{n}<b.\]
This proves that the open interval $(a,b)$ contains the rational number $m/n$, and thus completes the proof that the set of rational numbers is dense in the set of real numbers.
\end{myproof}

Recall that a set $A$ is said to be   {\bf countably infinite} if there is a bijection $f:  \mathbb{Z}^+\to A$. A set that is either finite or countably infinite is said to be {\bf countable}.  We assume that students have seen the proofs of the following.

\begin{proposition}{}
The set of integers $\mathbb{Z}$ and the set of rational numbers $\mathbb{Q}$ are countable, while the set of real numbers $\mathbb{R}$ is not countable.
\end{proposition}

Since the union of countable sets is  countable, this proposition implies that the set of irrational numbers is uncountable. Therefore, there are far more irrational numbers than rational numbers. Hence, it should not be surprising that the set of irrational numbers is also dense in the set of real numbers. To prove this, 
let us recall the following facts.
\begin{highlight}{Rational Numbers and Irrational Numbers}
\begin{enumerate}[1.]
\item
If $a$ and $b$ are rational numbers, then $a+b$ and $ab$ are rational numbers.
\item If $a$ is a nonzero rational number, $b$ is an irrational number, then $ab$ is an irrational number.
\end{enumerate}
\end{highlight}

\begin{theorem}{Denseness of the Set of Irrational Numbers}
The set of irrational numbers $\mathbb{R}\setminus\mathbb{Q}$ is dense in the set of real numbers $\mathbb{R}$.

\end{theorem}
\begin{myproof}{Proof}
Let $(a, b)$ be an open interval.  Define
\[c=\frac{a}{\sqrt{2}}, \hspace{1cm}d=\frac{b}{\sqrt{2}}.\]
Then $c<d$. By the denseness of rational numbers, there is a rational number $u$ that lies in the interval $(c,d)$. Hence, 
\[\frac{a}{\sqrt{2}}=c<u<d=\frac{b}{\sqrt{2}}.\]
Let $v=\sqrt{2}u$. Then $v$ is an irrational number satisfying
\[a<v<b.\]
This proves that the open interval $(a, b)$ contains the irrational number $v$, and thus completes the proof that the set of irrational numbers is dense in the set of real numbers.
\end{myproof}

\begin{example}{}
Is the set of integers $\mathbb{Z}$ dense in $\mathbb{R}$? Justify your answer.
\end{example}
\begin{solution}{Solution}
$(0,1)$ is an open interval that does not contain any integers. Hence, the set of integers is not dense in $\mathbb{R}$.
\end{solution}
\vspace{0.8cm}
\hrule
\vspace{0.8cm}
\noindent
{\bf \large Exercises  \thesection}
\setcounter{myquestion}{1}

 \begin{question}{\themyquestion}
Let  $S=\mathbb{Q}\setminus\mathbb{Z}$. Is the set $S$ dense in $\mathbb{R}$? Justify your answer.
\end{question}

\vp
\section{The Convergence of Sequences}\label{sec1.5}

{\bf Infinite sequences} play important roles in analysis. We will consider infinite sequences that are indexed by the set of positive integers
\[a_1, a_2, \ldots, a_n, \ldots\]
This can be considered as a function $f:\mathbb{Z}^+\rightarrow\mathbb{R}$, where $a_n=f(n)$.  The general term in the sequence is denoted by $a_n$.  In some occasions, we may also want to consider sequences that start with $a_0$.

In the sequel, when we say a sequence, we always mean an infinite sequence that is indexed by the set of positive integers, unless otherwise specified.
A sequence can be denoted by $\{a_n\}$ or $\{a_n\}_{n=1}^{\infty}$. This should not be confused with the set $\{a_n\,|\, n\in\mathbb{Z}^+\}$ that contains all terms in the sequence. 

There are various ways to specify a sequence. One of the ways is to give an explicit formula for the general term $a_n$. For example $\{1/n\}$ is the sequence with $a_n=1/n$. More precisely, it is the sequence with first five terms given by
\[1, \frac{1}{2}, \frac{1}{3}, \frac{1}{4}, \frac{1}{5}, \ldots.\]
A sequence can also be defined recursively, such as the following example. 

\begin{example}[label=ex23020301]{}
Let $\{a_n\}$ be the sequence defined by $a_1=2$, and for $n\geq 2$,
\[a_n=a_{n-1}+3.\]
  Find the first 5 terms of the sequence.
\end{example}
\begin{solution}{Solution}We compute recursively.\\
$
a_1 =2 $\\
$a_2= a_1+3=5 $\\
$a_3 =a_2+3=8$\\
$a_4 =a_3+3=11$\\
$a_5 =a_4+3=14$
 
\end{solution}
The sequence $\{a_n\}$ in Example \ref{ex23020301} is  an example of an arithmetic sequence. One can prove by induction that
\[a_n=3n-1.\]

\begin{example}[label=ex23020302]{}
Let $\{s_n\}$ be the sequence defined by $s_1=\frac{1}{2}$, and for $n\geq 2$,
\[s_n=s_{n-1}+\frac{1}{2^n}.\]
  Find the first 5 terms of the sequence.
\end{example}
\begin{solution}{Solution}We compute recursively.
\begin{align*}
s_1&=\frac{1}{2} \\
s_2&=s_1+\frac{1}{2^2} =\frac{3}{4}  \hspace{8cm}\\
s_3&=s_2+\frac{1}{2^3} =\frac{7}{8}\\
s_4&=s_3+\frac{1}{2^4} =\frac{15}{16} \\
s_5&=s_4+\frac{1}{2^5} =\frac{31}{32} \hspace{8cm}
\end{align*}
\end{solution}
The sequence $\{s_n\}$ in Example \ref{ex23020302} is the partial sum of the geometric sequence $\di\left\{\frac{1}{2^n}\right\}$. One can prove by induction that
\[s_n=1-\frac{1}{2^n}.\]
\begin{example}[label=ex23020303]{}
Let $\{s_n\}$ be the sequence defined by 
\[s_n=1+\frac{1}{2}+\cdots+\frac{1}{n}.\]
  This sequence can also be defined recursively by $s_1=1$, and for $n\geq 2$,
\[s_n=s_{n-1}+\frac{1}{n}.\]
\end{example}For the sequence $\{s_n\}$ defined in Example \ref{ex23020303}, the general term $s_n$ cannot be expressed as an explicit elementary function of $n$.

\begin{example}[label=ex23020304]{}
Let $\{a_n\}$ be the sequence defined by $a_1=2$, and for $n\geq 1$,
\[a_{n+1}=\begin{cases} a_n+\frac{1}{n}\quad &\text{if}\;a_n<3,\\a_n-\frac{1}{n}\quad &\text{if}\;a_n\geq 3.\end{cases}\]
Find the first six terms of the sequence.
\end{example}
\begin{solution}{Solution}We compute recursively.
\begin{align*}
a_1&=2<3\\
a_2&= a_1+1=3\geq 3 \\
a_3&=a_2-\frac{1}{2}=\frac{5}{2}<3\\
a_4&=a_3+\frac{1}{3}=\frac{17}{6}<3\\
a_5&=a_4+\frac{1}{4}=\frac{37}{12}\geq 3\hspace{8cm}
\\
a_6&=a_5-\frac{1}{5}=\frac{173}{60}
\end{align*}
\end{solution}

From the examples above, we observe that some sequences are monotone. 

\begin{definition}{Increasing and Decreasing Sequences}
\begin{enumerate}[1.]
\item
We say that a sequence $\{a_n\}$ is {\bf increasing} if
\[a_n\leq a_{n+1}\hspace{1cm}\text{for all}\;n\in\mathbb{Z}^+.\]
\item We say that a sequence is {\bf decreasing} if 
\[a_n\geq a_{n+1}\hspace{1cm}\text{for all}\;n\in\mathbb{Z}^+.\]
\item We say that a sequence $\{a_n\}$ is {\bf monotone} if it is an increasing sequence or it is a decreasing sequence.
\end{enumerate}
\end{definition}

\begin{example}{}
\begin{enumerate}[1.]
\item The sequence  $\{a_n\}$ defined in Example \ref{ex23020301} is increasing.
\item The sequence $\di\left\{\frac{1}{n}\right\}$ is decreasing.
\item The sequence  $\{a_n\}$ defined in Example \ref{ex23020304} is neither increasing nor decreasing.

\end{enumerate}
\end{example}

In analysis, we are often led to consider the behavior of a sequence $\{a_n\}$ when $n$ gets larger than larger. We are interested to know whether the sequence would approach a fixed value. This leads to the idea of convergence.

\begin{definition}{Convergence of   Sequences}
A sequence $\{a_n\}$ is said to {\bf converge} to the number $a$ if for every positive number $\varepsilon$, there is a positive integer $N$ such that for all $n\geq N$, 
\[|a_n-a|<\varepsilon.\] 
\end{definition}
Here the positive number $\varepsilon$ is used to measure the distance from the term $a_n$ to the number $a$. Since $\varepsilon$ can be any positive number, the distance can get as small as possible. 

 \begin{figure}[ht]
\centering
\includegraphics[scale=0.2]{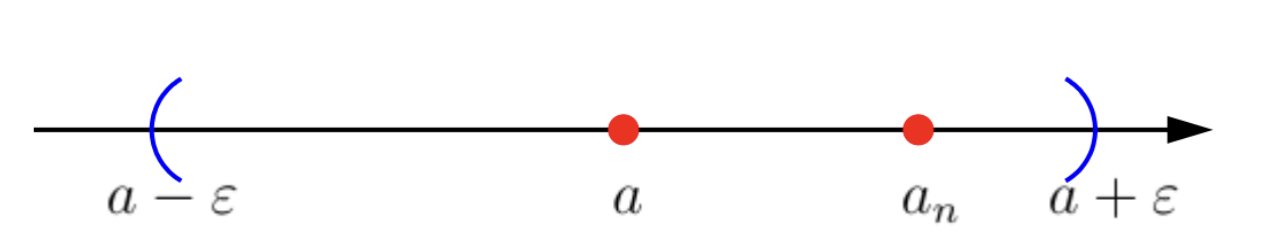}

\caption{  $|a_n-a|<\varepsilon$.}\label{figure2}
\end{figure}

One question that is natural to ask is whether a sequence $\{a_n\}$ can converge to two different numbers. This is   impossible.
\begin{theorem}[label=thm23020301]{}
A sequence cannot converge to two different numbers.
\end{theorem}
\begin{myproof}{Proof}
This is proved by contradiction. Assume that there is a sequence $\{a_n\}$ which converges to two different numbers $b$ and $c$.  Let
\[\varepsilon=\frac{|b-c|}{2}.\]
Since $b$ and $c$ are distinct, $|b-c|>0$ and so $\varepsilon>0$. By definition of convergence, there is a positive integer $N_1$ such that for all $n\geq N_1$,
\[|a_n-b|<\varepsilon.\]
Similarly, there is a positive integer $N_2$ such that for all $n\geq N_2$, 
\[|a_n-c|<\varepsilon.\]
If $N=\max\{N_1, N_2\}$, then $N\geq N_1$ and $N\geq N_2$. It follows that
\[|b-c|=|(a_N-c)-(a_N-b)|\leq |a_N-c|+|a_N-b|<\varepsilon+\varepsilon=|b-c|.\]
This gives $|b-c|<|b-c|$, which is a contradiction. Hence, we conclude that a sequence cannot converge to two different numbers.
\end{myproof}
\begin{figure}[ht]
\centering
\includegraphics[scale=0.2]{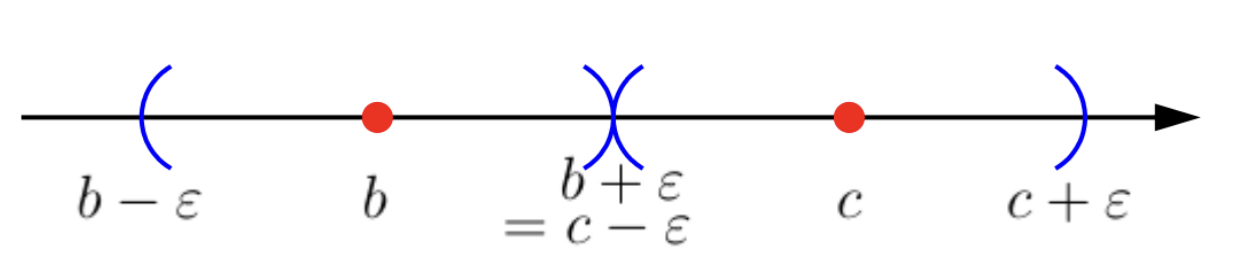}

\caption{ A sequence $\{a_n\}$ cannot converge to two different numbers $b$ and $c$.}\label{figure3}
\end{figure}

\begin{highlight}{Limit of a Sequence}
If a sequence $\{a_n\}$ converges to a number $a$, we say that the sequence is \textbf{convergent}. Otherwise, we say that it is \textbf{divergent}. Theorem \ref{thm23020301} says that for a convergent sequence $\{a_n\}$, the number $a$ that it converges to is unique. We call this unique number $a$ the \textbf{limit} of the convergent sequence $\{a_n\}$, and express the convergence of $\{a_n\}$ to $a$ as
\[\lim_{n\rightarrow\infty}a_n=a.\]

\end{highlight}

Using logical expression,     
\begin{highlight}{}
\[\lim_{n\rightarrow\infty}a_n=a\;\iff\; \forall \varepsilon>0, \;\exists N\in \mathbb{Z}^+, \;\forall n\geq N,\; |a_n-a|<\varepsilon.\]
 \end{highlight}
Let us look at a simple example of a constant sequence.

\begin{example}{}
Let $c$ be a real number and let $\{a_n\}$ be the sequence with $a_n=c$ for all $n\in\mathbb{Z}^+$. Then for any $\varepsilon>0$, we take $N=1$. For all $n\geq N=1$, we have
\[|a_n-c|=|c-c|=0<\varepsilon,\] which shows that
the limit of the constant sequence $\{a_n\}$   is $c$. Namely,
\[\lim_{n\rightarrow\infty}c=c.\]

\end{example}
Another simple example is the sequence $\{a_n\}$ with $\di a_n=\frac{1}{n}$.
\begin{example}[label=ex23020305]{}
Use the definition of convergence to show that 
\[\lim_{n\rightarrow\infty}\frac{1}{n}=0.\]
\end{example}
\begin{solution}{Solution}
Given $\varepsilon>0$, the Archimedean property asserts that there is a positive integer $N$ such that $1/N<\varepsilon$. If $n\geq N$, we have   
\[0<\frac{1}{n}\leq\frac{1}{N}<\varepsilon.\] 
This gives
\[\left|\frac{1}{n}-0\right|<\varepsilon\hspace{1cm}\text{for all}\;n\geq N.\]
By definition, we conclude that
\[\lim_{n\rightarrow\infty}\frac{1}{n}=0.\]
\end{solution}

Let $f:\mathbb{Z}^+\rightarrow\mathbb{Z}^+$ be a function satisfying
\[f(k)<f(k+1)\hspace{1cm}\text{for all}\;k\in \mathbb{Z}^+.\]
Then $f(\mathbb{Z}^+)$ is an infinite set of positive integers.   If we let
$n_k=f(k)$, then
\[n_1<n_2<n_3<\cdots.\]Namely, $n_1, n_2, n_3, \ldots$ is a strictly increasing sequence of positive integers. 
\begin{definition}{Subsequence}
Let $\{a_n\}$ be a sequence. A subsequence of $\{a_n\}$ is a sequence $\{a_{n_k}\}$ indexed by $k\in\mathbb{Z}^+$, where $n_k=f(k)$ is defined by a function $f:\mathbb{Z}^+\rightarrow\mathbb{Z}^+$  satisfying
\[f(k)<f(k+1)\hspace{1cm}\text{for all}\;k\in \mathbb{Z}^+.\]
\end{definition}

\begin{example}{}
The sequence $\{1/(2n-1)\}$ with first three terms given by 
\[1, \frac{1}{3}, \frac{1}{5},\]
is a subsequence of the sequence $\{1/n\}$ whose first five terms are
\[1, \frac{1}{2}, \frac{1}{3}, \frac{1}{4}, \frac{1}{5}.\]
\end{example}

If a sequence $\{a_n\}$ converges to $a$, what can we say about its subsequence? It is natural to expect any subsequence of $\{a_n\}$ also converges to $a$.

\begin{theorem}[label=thm23020305]{Subsequence of a Convergent Sequence}
If the sequence $\{a_n\}$ converges to $a$, then any of its subsequence also converges to $a$.
\end{theorem}
\begin{myproof}{Proof}
Let $\{a_{n_k}\}$ be a subsequence of $\{a_n\}$. Notice that for all $k\in\mathbb{Z}^+$, 
\[n_k\geq k.\]
Given $\varepsilon>0$, there is a positive integer $N$ such that for all $n\geq N$,
\[|a_n-a|<\varepsilon.\]
Take $K=N$. Then for all $k\geq K$, $n_k\geq n_K=n_N\geq  N$, and thus,
\[|a_{n_k}-a|<\varepsilon.\] This proves that $\{a_{n_k}\}$ indeed converges to $a$.
\end{myproof}

\begin{example}
{}
Find the limit 
\[\lim_{n\rightarrow\infty}\frac{1}{2^n}\] if it exists.
\end{example}
\begin{solution}{Solution}
Notice that $\di\left\{\frac{1}{2^n}\right\}$ is a subsequence of $\di\left\{\frac{1}{n}\right\}$ with $n_k=2^k$. By Example \ref{ex23020305}, 
\[\lim_{n\rightarrow\infty}\frac{1}{n}=0.\] We conclude from Theorem \ref{thm23020305} that \[\lim_{n\rightarrow\infty}\frac{1}{2^n}=0.\]
\end{solution}

\begin{example}{}
Show that the sequence $\{(-1)^n\}$ is divergent.
\end{example}
\begin{solution}{Solution}
Let $a_n=(-1)^n$. Then  for any positive integer $n$, $a_{2n-1}=-1$, and $a_{2n}=1$. The subsequence $\{a_{2n-1}\}$ of $\{a_n\}$ converges to $-1$, while the subsequence  $\{a_{2n}\}$ of $\{a_n\}$ converges to 1. Since there are two subsequences of $\{a_n\}$ that converge to two different limits, by Theorem \ref{thm23020305}, the sequence $\{a_n\}$ is not convergent.
\end{solution}

For the sequence $\{a_n\}$ defined in Example \ref{ex23020301}, we can see that the set $\{a_n\,|\, n\in \mathbb{Z}^+\}$ is not bounded above. Therefore, we would expect that the sequence does not converge to any number.

For simplicity, we say that a sequence $\{a_n\}$ is bounded above/bounded below/ bounded if the set $\{a_n\,|\, n\in \mathbb{Z}^+\}$  is bounded above/bounded below/bounded . If  the sequence $\{a_n\}$ is bounded above, we denote the supremum of the set $\{a_n\,|\, n\in \mathbb{Z}^+\}$  as $\sup\{a_n\}$. If  the sequence $\{a_n\}$ is bounded below, we denote the infimum of the set $\{a_n\,|\, n\in \mathbb{Z}^+\}$  as $\inf\{a_n\}$. 

We have the following theorem which guarantees that a convergent sequence must be bounded.

\begin{theorem}[label=thm23020304]{Boundedness of Convergent Sequence}
If a sequence $\{a_n\}$ is convergent, then it is bounded. Equivalently, if a sequence $\{a_n\}$ is not bounded, then it is not convergent.
\end{theorem}
\begin{myproof}{Proof}
Let $\{a_n\}$ be a convergent sequence that converges to the limit $a$. By definition of convergence with $\varepsilon=1$, there is a positive integer $N$ such that for all $n\geq N$,
\[|a_n-a|<1.\]
This implies that
\[|a_n|\leq |a_n-a|+|a|<1+|a|\hspace{1cm}\text{for all}\;n\geq N.\]
Define
\[M=\max\left\{|a_1|, |a_2|, \ldots, |a_{N-1}|, |a|+1\right\}.\]Then 
\[|a_n|\leq M\hspace{1cm}\text{for all}\;n\in \mathbb{Z}^+.\]This shows that the sequence $\{a_n\}$ is bounded.
\end{myproof}

\begin{example}{}
By Theorem  \ref{thm23020304},  the sequence $\{a_n\}$ defined in Example \ref{ex23020301} is not convergent.\end{example}

If the sequence $\{a_n\}$ is convergent, and $c$ is a constant, it is natural to expect that the sequence $\{ca_n\}$ is also convergent.  

\begin{proposition}[label=p23020401]{}
If the sequence $\{a_n\}$ converges to $a$, then the sequence $\{ca_n\}$ converges to $ca$.
\end{proposition}
\begin{myproof}{Proof}
  Given $\varepsilon>0$, the number\[\varepsilon_1=\di\frac{\varepsilon}{|c|+1}\] is also positive. Since $\{a_n\}$ converges to $a$, there is a positive integer $N$ such that for all $n\geq N$,
\[|a_n-a|<\varepsilon_1=\frac{\varepsilon}{|c|+1}.\]
It follows that for all $n\geq N$,
\[|ca_n-ca|=|c||a_n-a|<\frac{|c|}{|c|+1}\varepsilon<\varepsilon.\]
This proves that $\{ca_n\}$ converges to $ca$.
\end{myproof}

\begin{example}{}
By Proposition \ref{p23020401}, we find that for any constant $c$,
\[\lim_{n\rightarrow\infty}\frac{c}{n}=0.\]
\end{example}

 In the following, we establish a comparison theorem for limits.
\begin{theorem}[label=squeeze]{Squeeze Theorem}
Let $\{a_n\}$, $\{b_n\}$ and $\{c_n\}$ be three sequences. Assume that there is a positive integer $N_0$ such that for all $n\geq N_0$,
\[b_n\leq a_n\leq c_n.\]
If both the sequences $\{b_n\}$ and $\{c_n\}$ converge to $\ell$, then the sequence $\{a_n\}$ also converges to $\ell$.
\end{theorem}

\begin{myproof}{ Proof }
 
For a positive number $\varepsilon$, since the sequence $\{b_n\}$ converges to $\ell$, there is a positive integer $N_1$ such that for all $n\geq N_1$, 
\[|b_n-\ell|<\varepsilon.\]
This implies that for all $n\geq   N_1$, 
\[  b_n-\ell>-\varepsilon.\]
Similarly, since the sequence $\{c_n\}$ converges to $\ell$, there is a positive integer $N_2$ such that for all $n\geq N_2$,
\[|c_n-\ell|<\varepsilon.\]
This implies that for all $n\geq N_2$, 
\[c_n-\ell<\varepsilon.\]
Let $N=\max\{N_0, N_1, N_2\}$. If $n\geq N$, $n\geq N_0$, $n\geq N_1$ and $n\geq N_2$. Therefore, if $n\geq N$,
\[a_n-\ell\geq b_n-\ell>-\varepsilon,\]and
\[a_n-\ell\leq c_n-\ell<\varepsilon.\]
This proves that for all $n\geq N$,
\[|a_n-\ell|<\varepsilon.\]Therefore, the sequence $\{a_n\}$ converges to $\ell$.
\end{myproof}

When applying the squeeze theorem, we are interested in the limit of the sequence $\{a_n\}$. It is not enough to find two seqeunces $\{b_n\}$ and $\{c_n\}$ satisfying
\[b_n\leq a_n\leq c_n\] for all $n$ greater than or equal to a fixed $N_0$. The two sequences $\{b_n\}$ and $\{c_n\}$ must have the same limit. 
\begin{example}{}
For the sequence $\{a_n\}$ with
\[a_n=\frac{(-1)^n}{n},\]
we have
\[-\frac{1}{n}\leq a_n\leq \frac{1}{n}.\]
Since
\[\lim_{n\rightarrow\infty}\frac{1}{n}=0,\]we have
\[\lim_{n\rightarrow\infty}-\frac{1}{n}=0.\]
By squeeze theorem, 
\[\lim_{n\rightarrow\infty}\frac{(-1)^n}{n}=0.\]
\end{example}

More generally, we have the following.
\begin{theorem}[label=thm23020307]{}
The sequence $\{a_n\}$ converges to 0 if and only if the sequence $\{|a_n|\}$ converges to 0.
\end{theorem}
 
A word of caution. If the sequence $\{|a_n|\}$ is convergent, the sequence $\{a_n\}$ is not necessarily convergent. An example is the sequence $\{a_n\}$ with $a_n=(-1)^n$. 
Theorem \ref{thm23020307} asserts that if $\{|a_n|\}$ converges to $0$, then $\{a_n\}$ is convergent, and it converges to 0.
Nevertheless, if the sequence $\{a_n\}$ is convergent, the sequence $\{|a_n|\}$ is necessarily convergent (see Question \ref{absolute}).
 
\begin{myproof}{\linkt Proof of Theorem \ref{thm23020307}\linko }
 
First assume that the sequence $\{a_n\}$ converges to 0. Given $\varepsilon>0$, there is a positive integer $N$ such that for all $n\geq N$,
\[|a_n-0|<\varepsilon.\]\bp
Notice that \[\bigl|\,|a_n|-0\,\bigr|=|a_n|=|\,a_n-0\,|.\]
  Hence, for all $n\geq N$,
\[\bigl|\,|a_n|-0\,\bigr|<\varepsilon.\]This proves that the sequence $\{|a_n|\}$ converges to 0.

Next, we assume that the sequence $\{|a_n|\}$ converges to 0. Then the sequence $\{-|a_n|\}$ also converges to 0. Since
\[-|a_n|\leq a_n\leq |a_n|,\]
 squeeze theorem implies that the sequence $\{a_n\}$ converges to 0.
\end{myproof}
 
In the following, we discuss two useful results that can be deduced from specific information about a convergent sequence. They will be useful in the proofs of other theorems that we are going to discuss.

\begin{lemma}[label=23020405]{Sequence with Positive Limit}
If $\{a_n\}$ is a sequence that converge to a positive number $a$, there is a positive integer $N$ such that $a_n>a/2>0$ for all $n\geq N$. 
\end{lemma}
One can easily formulate a counterpart of this lemma for a sequence with negative limit.
\begin{myproof}{Proof}
Take $\varepsilon=a/2$. Then $\varepsilon>0$. Hence, there is a positive integer $N$ so that for all $n\geq N$, 
\[|a_n-a|<\frac{a}{2}.\]
This implies that for all $n\geq N$,
\[a_n-a>-\frac{a}{2}.\]
Thus, for all $n\geq N$,
\[a_n>\frac{a}{2}>0.\]
\end{myproof}

\begin{lemma}[label=23020406]{}
\begin{enumerate}[1.]
\item
Given that $\{a_n\}$ is a sequence that is bounded above by $c$. If  $\{a_n\}$ converges to $a$, then $a\leq c$. 
\item 
Given that $\{a_n\}$ is a sequence that is bounded below by $b$. If  $\{a_n\}$ converges to $a$, then $a\geq b$. 
\item Given that $\{a_n\}$ is a sequence satifying
\[b\leq a_n\leq c\hspace{1cm}\text{for all}\; n\in\mathbb{Z}^+.\]If  $\{a_n\}$ converges to $a$, then $b\leq a\leq c$.
\end{enumerate}
\end{lemma}
It is suffices to prove the first statement. The second statement follows by considering the negative of the sequence. The third statement follows by combining the results of the first two statements.
\begin{myproof}{Proof}

Given that \[\di\lim_{n\rightarrow\infty}a_n=a\quad\text{and}\quad a_n\leq c \quad \text{for all}\;n\in\mathbb{Z}^+,\] we want to show that $a\leq c$. Assume to the contrary that $a>c$. Take $\varepsilon=a-c$. Then $\varepsilon>0$. By definition of convergence, there is a positive integer $N$ such that for all $n\geq N$,
\[|a_n-a|<\varepsilon.\]
This implies that
\[a_n-a>-\varepsilon=c-a\hspace{1cm}\text{when}\;n\geq N.\]
Hence,
\[a_n>c\hspace{1cm}\text{when}\;n\geq N.\]This contradicts to $a_n\leq c$ for all $n\in\mathbb{Z}^+$.
Therefore, we must have $a\leq c$.

\end{myproof}

In Proposition \ref{p23020401}, we have seen what happens when   a convergent sequence is multiplied by a constant. In the following, we   inspect the behaviour of limits with respect to sums, products and quotients.
We start by sums.
\begin{theorem}[label=23020402]{Sums of Convergent Sequences}
If the sequences $\{a_n\}$ and $\{b_n\}$ converge to $a$ and $b$ respectively, the sequence $\{a_n+b_n\}$ converges to $a+b$. 
\end{theorem}
\begin{highlight}{Linearity of Limits of Sequences}
 Combining Proposition \ref{p23020401} and Theorem \ref{23020402}, we obtain  the following. If
\[\lim_{n\rightarrow \infty}a_n=a,\hspace{1cm}\lim_{n\rightarrow\infty}b_n=b,\]then for any constants $\alpha$ and $\beta$,
\[\lim_{n\rightarrow\infty}\left(\alpha a_n+\beta b_n\right)=\alpha a+\beta b.\]  

 \end{highlight}
\begin{myproof}{\linkt Proof of Theorem \ref{23020402}\linko}
 
Given a positive number $\varepsilon$, the number $\varepsilon/2$ is also positive. Since the sequence $\{a_n\}$ converges to $a$, there is a positive integer $N_1$ such that for all $n\geq N_1$,
\[|a_n-a|<\frac{\varepsilon}{2}.\]
Similarly, 
 there is a positive integer $N_2$ such that for all $n\geq N_2$,
\[|b_n-b|<\frac{\varepsilon}{2}.\]
Take $N=\max\{N_1, N_2\}$. Then $N$ is a positive integer and $N\geq N_1$, $N\geq N_2$. If $n\geq N$, triangle inequality implies that
\begin{align*}
\left|(a_n+b_n)-(a+b)\right|&=|(a_n-a)+(b_n-b)|\\
&\leq |a_n-a|+|b_n-b|\\
&<\frac{\varepsilon}{2}+\frac{\varepsilon}{2}\\
&=\varepsilon.
\end{align*}This proves that the sequence $\{a_n+b_n\}$ converges to $a+b$.
\end{myproof}

Now we consider products.
\begin{theorem}[label=23020403]{Products of Convergent Sequences}
If the sequences $\{a_n\}$ and $\{b_n\}$ converge to $a$ and $b$ respectively, the sequence $\{a_nb_n\}$ converges to $ab$. 
\end{theorem}

Notice that  Proposition \ref{p23020401} is actually a special case of this theorem when $\{b_n\}$ is a constant sequence.

\begin{myproof}{\linkt Proof of Theorem \ref{23020403}\linko}
 
Since $\{a_n\}$ and $\{b_n\}$ are  convergent sequences, Theorem \ref{thm23020304} says that each of them is bounded. We can choose a common positive number $M$ so that for all $n\in\mathbb{Z}^+$,
\[|a_n|\leq M,\hspace{1cm}|b_n|\leq M.\]By Lemma \ref{23020406}, 
\[|a|\leq M,\hspace{1cm}|b|\leq M.\]Now we want to show that the difference of $a_nb_n$ and $ab$ aproaches zero when $n$ gets large. This should be achieved by the fact that $|a_n-a|$ and $|b_n-b|$ both approach 0 when $n$ gets large. To compare $a_nb_n-ab$ to $a_n-a$ and $b_n-b$, we do some manipulations as follows.
\begin{align*}
a_nb_n-ab=(a_n-a)b_n+a(b_n-b).
\end{align*}It follows from triangle inequality that
\begin{equation}\label{eq230204_1}|a_nb_n-ab|\leq |a_n-a||b_n|+|a||b_n-b|\leq M\left(|a_n-a|+|b_n-b|\right).\end{equation}
Now we can show that $a_nb_n$ converges to $ab$. Given $\varepsilon>0$, since $\varepsilon/(2M)$ is also positive, there exists a positive integer $N_1$ such that
\[|a_n-a|<\frac{\varepsilon}{2M}\hspace{1cm}\text{when}\;n\geq N_1.\]Similarly, 
there exists a positive integer $N_2$ such that \bp
\[|b_n-b|<\frac{\varepsilon}{2M}\hspace{1cm}\text{when}\;n\geq N_2.\]
 Take $N=\max\{N_1, N_2\}$. When $n\geq N$, $n\geq N_1$ and $n\geq N_2$. It follows from 
\eqref{eq230204_1} that
\[|a_n-b_n|<M\left(\frac{\varepsilon}{2M}+\frac{\varepsilon}{2M}\right)=\varepsilon.\]
This completes the proof that the sequence $\{a_nb_n\}$ converges to $ab$.
\end{myproof}

For quotient of two  sequences, we notice that if $y\neq 0$,
\[\frac{x}{y}=x\times \frac{1}{y},\]which says that the quotient of $x$ by $y$ is a product of $x$ with the reciprocal of $y$. Hence, it is enough to consider the reciprocal of a nonzero sequence.

\begin{theorem}[label=23020407]
{Reciprocal of a Convergent Nonzero Sequence}
If $\{a_n\}$ is a nonzero sequence that converges to a nonzero limit $a$, the reciprocal sequence $\{1/a_n\}$ converges to $1/a$.
\end{theorem}
\begin{myproof}{Proof}
Without loss of generality, assume that $a>0$. Lemma \ref{23020405} implies that there is a positive integer $N_1$ such that
\[a_n>\frac{a}{2}>0\hspace{1cm}\text{when}\;n\geq N_1.\]
Given $\varepsilon>0$, $a^2\varepsilon/2$ is also positive. By definition of convergence, there is a positive integer $N_2$ such that when $n\geq N_2$,
\[|a_n-a|<\frac{a^2\varepsilon}{2}.\]
Take $N=\max\{N_1, N_2\}$. If $n\geq N$, 
\[
\left|\frac{1}{a_n}-\frac{1}{a}\right|=\frac{|a_n-a|}{|a_n||a|}<\frac{2}{a^2}\times\frac{a^2\varepsilon}{2}= \varepsilon.\]
This proves that the  sequence $\{1/a_n\}$ converges to $1/a$.
 
\end{myproof}

\begin{remark}[label=r23020401]{Reciprocal of a Sequence That Converges to 0}
In the statement of Theorem \ref{23020407}, it is crucial that $a\neq 0$. To see this, consider the sequence $\{a_n\}$ with $a_n=1/n$. It converges to $a=0$. The sequence $\{1/a_n\}$ is the sequence of natural numbers $\{n\}$, which does not converge. In fact, since $\{a_n\}$ converges to $0$, the sequence $\{1/a_n\}$ is not bounded. Hence, the sequence $\{1/a_n\}$ does not converge.
\end{remark}

\begin{corollary}[label=23020408]{Quotients of Convergent Sequences}
Given that $\{a_n\}$ is a sequence that converges to $a$, $\{b_n\}$ is a nonzero sequence that converges to $b$. If $b\neq 0$, the sequence $\{a_n/b_n\}$ converges to $a/b$. 
\end{corollary}
The results about sums, products and quotients of convergent sequences can be summarized in the following.
\begin{highlight}{Operations on Convergent Sequences}
Given that
\[\lim_{n\rightarrow\infty}a_n=a\hspace{1cm}\text{and}\hspace{1cm}\lim_{n\rightarrow\infty}b_n=b.\]
\begin{enumerate}[1.]
\item For any constants $\alpha$ and $\beta$, 
$\di \lim_{n\rightarrow\infty}(\alpha a_n+\beta b_n)=\alpha a+\beta b$.

\item $\di\lim_{n\rightarrow \infty}a_nb_n=ab$.

\item If $b_n\neq 0$ for all $n\in\mathbb{Z}^+$ and $b\neq 0$, $\di \lim_{n\rightarrow\infty}\frac{a_n}{b_n}=\frac{a}{b}$.

\end{enumerate}
\end{highlight}
These will be used repeatedly in the future.
Let us now look at some examples how these properties are applied.

 \begin{example}{}
Let $m$ be a positive integer.  
Product rule of limits implies that
\[\lim_{n\rightarrow\infty}\frac{1}{n^m}=\underbrace{\lim_{n\rightarrow\infty}\frac{1}{n}\times\cdots\times\lim_{n\rightarrow\infty}\frac{1}{n}}_{m\;\text{terms}}=0.\]
\end{example}

\begin{example}{}
Determine whether the limit exists. If it exists, find the limit.
\begin{enumerate}
\item[(a)]
$\di\lim_{n\rightarrow \infty}\left(3+\frac{(-1)^n}{3n-2}\right)$
\item[(b)]
$\di \lim_{n\rightarrow\infty}\frac{2n^2+3n+4}{5-7n^2}$

\item[(c)] $\di \lim_{n\rightarrow\infty}\frac{n+1}{n^2+1}$
\item[(d)]
$\di \lim_{n\rightarrow\infty}\frac{n^2+1}{n+1}$
\end{enumerate}
\end{example}
\begin{solution}{Solution}
\begin{enumerate}
\item[(a)]
Since $\{1/(3n-2)\}$ is a subsequence of the sequence $\{1/n\}$, it converges to $0$. By Theorem \ref{thm23020307},
\[\lim_{n\rightarrow \infty}\frac{(-1)^n}{3n-2}=0.\] 
Hence,
\[\lim_{n\rightarrow \infty}\left(3+\frac{(-1)^n}{3n-2}\right)=\lim_{n\rightarrow \infty} 3+\lim_{n\rightarrow \infty}\frac{(-1)^n}{3n-2}=3+0=3.\]
\item[(b)] The sequence $\{2n^2+3n+4\}$ is not bounded. So it does not have a limit. We cannot apply quotient rule of limits directly. Instead,  we need to do some manipulations. Divide the numerator and the denominator by $n^2$ and then apply the rules for limits, we have
\[\lim_{n\rightarrow\infty}\frac{2n^2+3n+4}{5-7n^2}=\lim_{n\rightarrow\infty}\frac{2+\di\frac{3}{n}+\frac{4}{n^2}}{\di \frac{5}{n^2}-7}=\frac{2+0+0}{0-7}=-\frac{2}{7}.\]
\end{enumerate}
\end{solution}
\begin{solution}{ }
\begin{enumerate}
\item[(c)]
Divide the numerator and the denominator by $n^2$ and then apply the rules for limits, we have
\[\lim_{n\rightarrow\infty}\frac{n+1}{n^2+1}=\lim_{n\rightarrow\infty}\frac{\di \frac{1}{n}+\frac{1}{n^2}}{1+\di \frac{1}{n^2}}=\frac{0+0}{1+0}=0.\]
\item[(d)] Since the reciprocal of the sequence has limit 0 by part (c), we find that
\[ \lim_{n\rightarrow\infty}\frac{n^2+1}{n+1}\] does not exist.
\end{enumerate}
\end{solution}

We have seen in Section \ref{sec1.3} that the supremum or infimum of a set is not necessarily an element of the set. The supremum of a set is an element of the set if and only if the set has a maximum. Analogously,  the infimum of a set is an element of the set if and only if the set has a minimum.

Even though the supremum and infimum of a set might fail to be an element of the set,  they are always limits of sequences in that set.
\begin{lemma}[label=23020510]{Supremum and Infimum as Limits}
Let $S$ be a subset of real numbers.
\begin{enumerate}[1.]
\item
If $S$ is bounded above, there is a sequence $\{u_n\}$ in $S$ that converges to $u=\sup S$.
\item If $S$ is bounded below, there is a sequence $\{\ell_n\}$ in $S$ that converges to $\ell=\inf S$.
\end{enumerate}
\end{lemma}

\begin{example}[label=23020601]{}
Consider the set $S=(-\infty, \pi)$. It is bounded above with $\sup S=\pi$. The sequence $\{u_n\}$ with 
\[u_n=\pi -\frac{1}{n}\] is a sequence in $S$ that converges to $\pi=\sup S$.
\end{example}

To prove Lemma \ref{23020510}, it suffices for us to prove the first statement.

\begin{myproof}{\linkt Proof of Lemma \ref{23020510}\linko}
 
 Assume that $S$ is bounded above. Then the completeness axiom asserts that $u=\sup S$ exists. For any positive integer $n$, $u-1/n$ is smaller than $u$. Hence, $u-1/n$ is not an upper bound of $S$.   This implies that there is an element $u_n$ of $S$ such that
\[u_n>u-\frac{1}{n}.\]
Since $u_n$ is in $S$ and $u$ is an upper bound of $S$, we have $u_n\leq u$. In other words, we have
\[u-\frac{1}{n}<u_n\leq u\hspace{1cm}\text{for all}\;n\in\mathbb{Z}^+.\]
Since
\[\lim_{n\rightarrow \infty}\left(u-\frac{1}{n}\right)=\lim_{n\rightarrow \infty} u=u,\]
  squeeze theorem implies that
\[\lim_{n\rightarrow \infty}u_n=u.\]This means $\{u_n\}$ is  a sequence in $S$ that converges to $u$.
\end{myproof}

\vp
\noindent
{\bf \large Exercises  \thesection}
\setcounter{myquestion}{1}
 \begin{question}{\themyquestion}
Let $a$ be a positive integer that is larger than 1. Show that
$\di \lim_{n\rightarrow\infty}\frac{1}{a^n}=0$.
\end{question}

\atc
\begin{question}{\themyquestion}
If $\{a_n\}$ is a sequence that converge to a negative number $a$,   show that there is a positive integer $N$ such that $a_n<a/2<0$ for all $n\geq N$. 
\end{question}
\atc
 \begin{question}{\themyquestion}
Determine whether the limit exists. If it exists, find the limit.
\begin{enumerate}[(a)]
\item
$\di\lim_{n\rightarrow \infty} \frac{3n+(-1)^n}{n+2} $
\item
$\di \lim_{n\rightarrow\infty}\frac{4n+2}{7n^2+3n}$

\item $\di \lim_{n\rightarrow\infty}\frac{n}{2n^2+n+5}$
\item
$\di \lim_{n\rightarrow\infty}\frac{n^2+4n}{n+3}$
\end{enumerate}
\end{question}

\atc

 \begin{question}[label=absolute]{\themyquestion}
If $\{a_n\}$ is a sequence that converges to   $a$, use the definition of convergence to show that the sequence $\{|a_n|\}$ converges to $|a|$. 
\end{question}

\atc

 \begin{question}{\themyquestion:\; Last Statement in Remark \ref{r23020401}}
Given that $\{a_n\}$ is a nonzero sequence that converges to 0.
\begin{enumerate}[(a)]
\item Show that $\{1/a_n\}$ is not bounded.
\item Conclude that the sequence $\{1/a_n\}$ is divergent. 

\end{enumerate}
\end{question}

\atc 
\begin{question}{\themyquestion}
Let $\{a_n\}$ and $\{b_n\}$ be  sequences. Assume that there is a real number $a$ such that
\[|a_n-a|\leq  b_n \hspace{1cm}\text{for all}\;n\in\mathbb{Z}^+.\]
If $\di\lim_{n\rightarrow\infty}b_n=0$, show that
\[\lim_{n\rightarrow\infty}a_n=a.\]
\end{question}

\atc

\begin{question}{\themyquestion:\; The Convergence of the Sequence in Example \ref{ex23020304}}
Consider the sequence $\{a_n\}$  defined in Example by \ref{ex23020304}. It is defined recursively by $a_1=2$, and for $n\geq 1$,
\[a_{n+1}=\begin{cases} a_n+\frac{1}{n}\quad &\text{if}\;a_n<3,\\a_n-\frac{1}{n}\quad &\text{if}\;a_n\geq 3.\end{cases}\]
\begin{enumerate}
[(a)]
\item
Show that $|a_{n+1}-3|\leq\di\frac{1}{n}$ for all $n\in\mathbb{Z}^+$.\\
$[$Hint: Use induction.$]$
\item Show that the sequence $\{a_n\}$ is convergent and find its limit. 
\end{enumerate}
\end{question}
\vp
\section{Closed Sets and Limit Points }\label{sec1.6}

When we study convergence of sequences, we measure the closeness between points  by a positive number $\varepsilon$. A point $x$ is within $\varepsilon$ from the point $a$ if $x$ is in the open interval $(a-\varepsilon, a+\varepsilon)$. More generally, we define a neighbourhood of the point $a$ as follows.
\begin{definition}{Neighbourhood}
Given $a$ is a point in $\mathbb{R}$, a neighbourhood of $a$ is an open interval $(b, c)$ that contains $a$.
\end{definition}
The concept of neighbourhood is closely related to the concept of interior point. 
\begin{definition}{Interior Point}
 If $S$ is a set of real numbers, and there is a neighbourhood of the point $a$ that is contained in $S$, we call $a$ an interior point of $S$.
\end{definition}
In this section, we use sequences to define and  study some properties of subsets of real numbers. Given a subset $S$ of real numbers, we say that a sequence $\{a_n\}$ is in $S$ if each of the terms $a_n$  is  a point in $S$. In other words, the sequence $\{a_n\}$ is in $S$ means that the set $\{a_n\,|\,n\in\mathbb{Z}^+\}$ is a subset of $S$. We will abuse notation and write this as $\{a_n\}\subset S$ when there is no confusion.
We start with a simple but useful lemma.
\begin{lemma}[label=23020505]{}
Let $S$ be a subset of real numbers. If $\{a_n\}$ is a sequence in $S$ that converges to $a$, then every neighbourhood of $a$ contains a point of $S$.
\end{lemma}
\begin{myproof}{Proof}
Let $(b, c)$ be a neighbourhood of $a$. Since $a$ is in $(b,c)$,  $b<a<c$, and hence the number
\[\varepsilon=\min\{a-b, c-a\}\] is positive. 
By definition, $a-b\geq\varepsilon$, $c-a\geq\varepsilon$. 
Since $\{a_n\}$ converges to $a$, there is a positive integer $N$ such that for all $n\geq N$, 
\[|a_n-a|<\varepsilon.\]
In particular,
\[b\leq a-\varepsilon<a_N<a+\varepsilon\leq c.\]This shows that $a_N$ is a point in $S$ that is in the neighbourhood $(b, c)$ of $a$. 
\end{myproof}
 \begin{figure}[ht]
\centering
\includegraphics[scale=0.2]{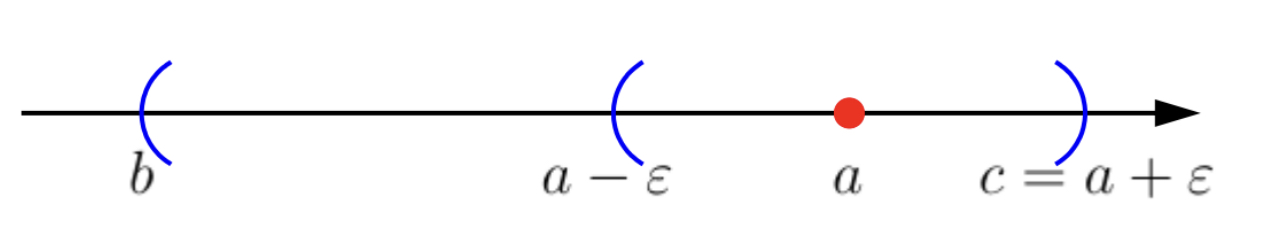}

\caption{  $b<a<c$ and $\varepsilon=c-a\leq a-b$.}\label{figure4}
\end{figure}
Next, we revisit the concept of denseness.
\begin{theorem}[label=23020409]{}
Let $S$ be a subset of real numbers. Then $S$ is dense in $\mathbb{R}$  if and only if every real number $x$ is the limit of a sequence in $S$.
\end{theorem}

Since we have proved that each of the set of rational numbers and the set of irrational numbers  is dense in the set of real numbers, we immediately obtain the following.
\begin{corollary}{}
Let $x$ be a real number.
\begin{enumerate}[1.]
\item
There is a sequence of rational numbers $\{p_n\}$ that converges to $x$.
\item There is a sequence of irrational numbers $\{q_n\}$ that converges to $x$.
\end{enumerate}
\end{corollary}

\begin{myproof}{\linkt Proof of Theorem \ref{23020409}}

First we assume that the set $S$ is dense in $\mathbb{R}$. Given a real number $x$, we want to show that there is a sequence in $S$ that converges to $x$. For each positive integer $n$, since $S$ is dense in $\mathbb{R}$, there is an element of $S$ in the open interval $(x-1/n, x)$. Choose one of these elements and denote it by $a_n$. Then $\{a_n\}$ is a sequence in $S$ satisfying
\[x-\frac{1}{n}<a_n<x\hspace{1cm}\text{for all}\;n\in\mathbb{Z}^+.\]
By squeeze theorem, the sequence $\{a_n\}$ converges to $x$.

Conversely, assume that every real number $x$ is the limit of a sequence in $S$. We want to show that $S$ is dense in $\mathbb{R}$. 
Let $(a,b)$ be an open interval. Take any point $x$ in the interval $(a,b)$.
By assumption, there is a sequence $\{c_n\}$ in $S$ which converges to $x$. By Lemma \ref{23020505}, the interval $(a,b)$ contains a point of $S$. Thus we have shown that every open interval $(a, b)$ contains a point of $S$. This proves that $S$ is dense in $\mathbb{R}$.
\end{myproof}

\begin{example}[label=23020501]{}
Let $x=\sqrt{2}$, and define the sequences $\{p_n\}$ and $\{q_n\}$ by
\[p_n=\frac{\lfloor 10^n\sqrt{2}\rfloor}{10^n}, \hspace{1cm} q_n=\sqrt{2}.\]
Here $\lfloor a\rfloor$ is the floor of $a$. By definition,
\[10^n\sqrt{2}-1<\lfloor 10^n\sqrt{2}\rfloor\leq 10^n\sqrt{2}.\]
Therefore,
\[\sqrt{2}-\frac{1}{10^n}<p_n\leq \sqrt{2}.\]
By squeeze theorem, 
$\{p_n\}$ converges to $\sqrt{2}$. Since $\lfloor 10^n\sqrt{2}\rfloor$ is an integer, $p_n$ is a rational number. Hence, $\{p_n\}$ is a sequence of rational numbers that converges to $x=\sqrt{2}$. Obviously, $\{q_n\}$ is a sequence of irrational numbers that converges to $x=\sqrt{2}$.
\end{example}
The number $p_n$ is the rational number obtained by truncating the decimal expansion of $\sqrt{2}$ to give a number with $n$ decimal places.
The first 7 terms of the sequence $\{p_n\}$ are
\[1.4, \;1.41, \;1.414, \;1.4142, \;1.41421,\; 1.414213,\; 1.4142135.\]

Now we introduce the concept  of closed sets.
\begin{definition}{Closed Set}
Let $S$ be a subset of $\mathbb{R}$. We say that $S$ is closed in $\mathbb{R}$ provided that if $\{a_n\}$ is a sequence of points in $S$ that converges to the limit $a$, the point $a$ is also in $S$. 
\end{definition}

\begin{example}{}
The three statements in Lemma \ref{23020406} imply that intervals of the form $(-\infty, a]$, $[a, \infty)$ and $[a,b]$ are closed subsets of $\mathbb{R}$. In particular, we call them closed intervals, and $[a,b]$ is a closed and bounded interval.
\end{example}

\begin{remark}{}
\begin{enumerate}[1.]
\item
By definition, $\mathbb{R}$ is closed in $\mathbb{R}$.
\item $\emptyset$ is closed in $\mathbb{R}$ because the statement that defines a closed set is a statement of the form $p\to q$, where $p$ is always false for an empty set. Hence, for an empty set, this statement $p\to q$ that defines a closed set is vacuously true.\end{enumerate}
\end{remark}

\begin{example}[label=23020502]{}
Is the interval $(0,2)$  closed in $\mathbb{R}$?
\end{example}
\begin{solution}{Solution}
The sequence $\{1/n\}$ is a sequence in the interval $(0, 2)$ that converges to the point $0$ that is not in $(0,2)$. Hence, the interval $(0,2)$ is not closed in $\mathbb{R}$.
\end{solution}

\begin{remark}{}
One can prove that if $S$ is an interval of the form $(a, b)$, or $(a, b]$, or $[a, b)$, or $(-\infty, a)$, or $(a, \infty)$, then $S$ is not closed in $\mathbb{R}$.
\end{remark}

\begin{example}{}
Is the set of rational numbers $\mathbb{Q}$ closed in $\mathbb{R}$? 
\end{example}
\begin{solution}{Solution}
We have seen in 
Example \ref{23020501} that there is a sequence in the set $\mathbb{Q}$ that converges to $\sqrt{2}$, which is not in $\mathbb{Q}$. Hence, $\mathbb{Q}$ is not closed in $\mathbb{R}$.
\end{solution}

The concept of closed sets is defined in terms of limits of sequences. This leads us to the concept of limit points.
\begin{definition}{Limit Points}
Let $S$ be a subset of real numbers. A point $x$ in $\mathbb{R}$ is called a limit point of the set $S$ if there is a sequence of points in $S\setminus \{x\}$ that converges to $x$.
\end{definition}
 
Notice that  $\{a_n\}$ is a sequence in $S\setminus\{x\}$ if and only if it is a sequence in $S$ with none of the terms $a_n$  equal to $x$. 
\begin{example}[label=23020503]{}
In the solution of Example \ref{23020502}, we have seen that the sequence $\{1/n\}$ in $(0,2)$  converges to the point $0$. Since none of the $a_n$ is 0,  0 is a limit point of the set $(0,2)$. 
\end{example}
\begin{highlight}{Limits and Limit Points}
Although the concepts of limits and limit points are closely related, one should not get confused. The limit of a convergent sequence $\{a_n\}$ is not necessarily the limit point of the set $\{a_n\,|\,n\in\mathbb{Z}^+\}$. For example, $c$ is the limit of the constant sequence $\{a_n\}$ with $a_n=c$ for all $n\in\mathbb{Z}^+$, but $c$ is not a limit point of the set $\{a_n\,|\,n\in\mathbb{Z}^+\}=\{c\}$.\end{highlight}

\begin{example}{}
Determine the set of limit  points of the set $(0,2)$.
\end{example}
\begin{solution}{Solution}
We claim that every point in $[0,2]$ is a limit point of the  set $(0,2)$. 

Example \ref{23020503} shows that 0 is a limit point of $(0,2)$. The sequence
$\{2-1/n\}$ is a sequence in $(0,2)$ that converges to 2. Hence, 2 is also a limit point of $(0,2)$.

For any $c\in (0,2)$, $c>0$. Let $m$ be a positive intger such that $1/m<c$. Then $\{c-1/(n+m)\}$ is  a sequence in $(0,2)$ that converges to $c$. Hence,  $c$ is a limit point of $(0,2)$.

This completes the proof that the set of limit points of $(0,2)$  is $[0,2]$.
\end{solution}

\begin{remark}{}
\begin{enumerate}[1.]
\item For intervals of the form $(a,b)$, $(a, b]$, $[a, b)$ or $[a,b]$, the set of limit points is $[a,b]$.
\item For intervals of the form $(-\infty, a)$ or $(-\infty, a]$, the set of limit points is $(-\infty, a]$.
\item For intervals of the form $(a, \infty)$ or $[a, \infty)$, the set of limit points is $[a, \infty)$.
\end{enumerate}

\end{remark}

\begin{example}[label=230322_1]{}
Show that the set $\mathbb{Z}$ does not have limit points.
\end{example}

\begin{solution}{Solution}

If $n$ is an integer, $x$ is contained in the open interval $(n-1,n+1)$ that does not contain any integer other than $n$ itself. Hence, there is no sequence in $\mathbb{Z}\setminus \{n\}$ that converges to $n$. Therefore, an integer $n$ is not a limit point of $\mathbb{Z}$.\bs

If $x$ is not an integer, it is contained in the interval $(\lfloor x\rfloor, \lceil x\rceil)$ that does not contain any integers. By Lemma \ref{23020505}, $x$ is not a limit of a sequence in $\mathbb{Z}$. Therefore, $x$ is not a limit point of $\mathbb{Z}$.
\end{solution}

\begin{definition}{Isolated Points}
Let $S$ be a subset of real numbers. We say that $x$ is an isolated point of $S$ if 
\begin{enumerate}[(a)]
\item $x$ is in $S$;
\item $x$ is not a limit point of $S$.
\end{enumerate}
\end{definition}

By definition, we have the following.
\begin{highlight}{Isolated Points vs Limit Points}
A point in a set $S$ is either a limit point or an isolated point of the set.
\end{highlight}

\begin{example}{}
By Example \ref{230322_1}, every point in the set of integers $\mathbb{Z}$ is an isolated point of the set.
\end{example}

The following is quite obvious from the definition of isolated points and Lemma \ref{23020505}.
\begin{theorem}[label=23020810]{}
Let $S$ be a subset of real numbers. A point $x$ in $S$ is an isolated  point if and only if there is a neighbourhood  $(a,b)$ of  $x$ that intersects the set $S$ only at the point $x$.
\end{theorem}

We have seen that a limit point of a set is not necessarily a point of that set. The following gives a characterization of closed sets in terms of limit points.
\begin{theorem}{}
Let $S$ be a subset of real numbers. The set $S$ is closed in $\mathbb{R}$ if and only if it contains all its limit points.

\end{theorem}
To prove a statement of the form $p\iff q$, we can prove $p\implies q$ and $\neg p \implies \neg q$.
\begin{myproof}{Proof}
Assume first $S$ is closed in $\mathbb{R}$. Let $x$ be a limit point of $S$. Then there is a sequence $\{a_n\}$ in $S\setminus\{x\}$ that converges to $x$. In particular, $\{a_n\}$ is a sequence in $S$ that converges to $x$. Since $S$ is closed in $\mathbb{R}$, $x$ is in $S$. This proves that $S$ contains all its limit points.

Now assume that $S$ is not closed in $\mathbb{R}$. Then there is a sequence $\{a_n\}$ in $S$ that converges to a point $x$, but $x$ is not in $S$. Since $x$ is not in $S$, none of the terms in the sequence $\{a_n\}$ is in $S$. Therefore, $x$ is a limit point of $S$. This shows that $S$ does not contain all its limit points.
\end{myproof}

\vp
\noindent
{\bf \large Exercises  \thesection}
\setcounter{myquestion}{1}

 \begin{question}{\themyquestion} 
 Show  that every real number is a limit point of the set of rational numbers.

\end{question}

\atc
 \begin{question}{\themyquestion} 
Let $S$ be the set
\[S=\left\{\left.\frac{1}{n}\,\right|\,n\in\mathbb{Z}^+\right\}\;\cup\;\{0\}.\]
\begin{enumerate}[(a)]
\item Find the set of limit points and the set of isolated points of $S$.
\item Is $S$ a closed set?
\end{enumerate}
\end{question}

\atc

 \begin{question}{\themyquestion} 
 Determine whether each of the following is a closed set.
\begin{enumerate}[(a)]\item $A=[2, 3]\;\cup\;[4, 7]$
\item $B=(-\infty, 2]\cup [3, 5]$
\item $C=\mathbb{R}\setminus (-1,1)$
\item $D=[1, 2)\cup [2, 4]$
\item $E=(1, 2)\cup (3, 4]$
\end{enumerate}
\end{question}

\vp
\section{The Monotone Convergence Theorem}\label{sec1.7}

Recall that a  sequence $\{a_n\}$ is monotone if it is increasing or it is decreasing. Obviously,  an increasing sequence is bounded below, and a decreasing sequence is  bounded above. However, a monotone sequence is not necessary convergent. A simple example is the sequence of natural numbers $\{n\}$. In the following, we give a characterization for a monotone sequence to be convergent. 
\begin{theorem}{The Monotone Convergence Theorem}
Let $\{a_n\}$ be a monotone sequence.
\begin{enumerate}[1.]
\item If $\{a_n\}$ is increasing, then $\{a_n\}$ is convergent if and only if it is bounded above. In this case, 
\[\lim_{n\rightarrow\infty}a_n=\sup \{a_n\}.\]
\item If $\{a_n\}$ is decreasing, then $\{a_n\}$ is convergent if and only if it is bounded below. In this case, 
\[\lim_{n\rightarrow\infty}a_n=\inf \{a_n\}.\]
\end{enumerate}
\end{theorem}
\begin{highlight}{Convergence Criteria for Monotone Sequences}
The monotone convergence theorem says that a montonone sequence is convergent if and only if it is bounded.
\end{highlight}
It is suffices to prove the case where $\{a_n\}$ is an increasing sequence. 
\begin{myproof}{Proof}
First suppose that $\{a_n\}$ is an increasing sequence that is convergent.  Then $\{a_n\}$ is bounded. So it is bounded above.

Conversely, suppose that $\{a_n\}$ is an increasing sequence that is bounded above. Then $a=\sup  \{a_n\}$ exists.  Now we use the same argument as in the proof of Lemma \ref{23020510}. Given $\varepsilon>0$, since $a-\varepsilon$ is less than $a$, it is not an upper bound of the set $S=\{a_n\,|\,n\in\mathbb{Z}^+\}$. Hence, there is a positive integer $N$ such that
\[a_N>a-\varepsilon.\]\bp
It follows that
\[a_n\geq a_N>a-\varepsilon\hspace{1cm}\text{for all}\;n\geq N.\]
Since $a$ is an upper bound of $S$, we also have $a_n\leq a$ for all $n$.  Thus,
\[|a_n-a|<\varepsilon\hspace{1cm}\text{for all}\;n\geq N.\]
This shows that the sequence $\{a_n\}$ converges to $a$.
\end{myproof}

The monotone convergence theorem is very useful because we can conclude the convergence of a sequence without apriori knowing  the limit of the sequence.  It is a consequence of the completeness axiom which asserts that any set that is bounded above has a supremum.  

\begin{example}[label=23020511]{}Let $a$ be a number in the interval $(0,1)$. Show that
\[\lim_{n\rightarrow\infty}a^n=0.\]

\end{example}

\begin{remark}{}

It follows from Theorem \ref{thm23020307} that for any $a$ in the interval $(-1,1)$, 
\[\lim_{n\rightarrow\infty}a^n=0.\]
\end{remark}

\begin{solution}{\linkt Solution to Example \ref{23020511}\linko}
Since $0<a<1$, for any positive integer $N$,
\[a^{n+1}=a^n\times a<a^n.\]
Hence, the sequence $\{a^n\}$ is decreasing.  On the other hand, $a^n>0$ for all $n\in\mathbb{Z}^+$. Hence, $\{a^n\}$ is a   decreasing sequence that is bounded below. By the monotone convergence theorem, $\{a^n\}$ converges to a number $\ell$. 
\bs
Since $\{a^{n+1}\}$ is a subsequence of $\{a^n\}$, it also converges to $\ell$.
Applying limit law to 
\[a^{n+1}=a\times a^n,\]
we have
\[\ell=\lim_{n\rightarrow\infty}a^{n+1}=a\lim_{n\rightarrow\infty}a^n=a\ell.\]
Since $a\neq 1$, we must have $\ell=0$.
\end{solution}

\begin{example}[label=23020512]{}
Define the sequence $\{a_n\}$ inductively by $a_1=1$ and for all $n\geq 1$,
\[a_{n+1}=\frac{2a_n+2}{a_n+2}.\]
Show that $\{a_{n}\}$ is convergent and find its limit.
\end{example}

\begin{solution}
{Solution} 
First notice that $a_n>0$ for all $n\in\mathbb{Z}^+$. 
When $n\geq 2$,
\[
a_{n+1}-a_n =\frac{2a_n+2}{a_n+2}-\frac{2a_{n-1}+2}{a_{n-1}+2}
 =\frac{2(a_n-a_{n-1})}{(a_n+2)(a_{n-1}+2)}.
\] Now, $a_2=4/3>a_1$. Hence, we deduce that $a_{n+1}-a_n>0$ for all $n\in\mathbb{Z}^+$.  In other words, $\{a_{n}\}$ is an increasing sequence.
 For all $n\geq 1$,
\[a_{n+1}=2-\frac{2}{a_n+2}<2.\]
Hence, $\{a_n\}$ is bounded above by 2.
Since $\{a_n\}$ is an increasing sequence that is bounded above, by monotone convergence theorem, it converges to a limit $u=\sup\{a_n\}$. 
Since $\{a_{n+1}\}$ is a subsequence of $\{a_n\}$,  it also converges to $u$. Apply the limit laws to \[a_{n+1}=\frac{2a_n+2}{a_n+2},\]we find that
\[u=\frac{2u+2}{u+2}.\]\bs
This implies that 
\[u^2=2.\]
Since $a_n>0$, we must have $u\geq 0$.  Hence, $u=\sqrt{2}$.

\end{solution}
Notice that Example \ref{23020512} is closely related to Example \ref{ex23020101}. The sequence $\{a_n\}$ defined in Example \ref{23020512}  is another sequence of rational numbers which converges to $\sqrt{2}$.

The next example is a classical one.

\begin{example}[label=23020507]{}

Show  that the limit
\[\lim_{n\rightarrow\infty} \left(1+\frac{1}{n}\right)^n\] exists. 
 
\end{example}

\begin{solution}{Solution}Let \[a_n=\left(1+\frac{1}{n}\right)^n.\]
  Given a positive integer $n$, notice that
\[\frac{ a_{n+1}}{a_n}=\frac{n+2}{n+1}\times \left(\frac{(n+2)n}{(n+1)^2}\right)^n.\]
By Bernoulli's inequality (see Question \ref{Q23020501}), 
\begin{align*}
 \left(\frac{(n+2)n}{(n+1)^2}\right)^n=\left(1-\frac{1}{(n+1)^2}\right)^n\geq 1-\frac{n}{(n+1)^2}=\frac{n^2+n+1}{(n+1)^2}.
\end{align*}
It follows that
\begin{align*}
\frac{a_{n+1}}{a_n}\geq \frac{(n+2)(n^2+n+1)}{(n+1)^3}=\frac{n^3+3n^2+3n+2}{n^3+3n^2+3n+1}>1.
\end{align*}
This shows that 
\[a_{n+1}>a_n\hspace{1cm}\text{for all}\;n\in\mathbb{Z}^+.\]\bs
Hence, $\{a_n\}$ is monotonically increasing.
 Using binomial expansion, we have
\[a_n=\sum_{k=0}^{n}\binom{n}{k}\frac{1}{n^k}.\]
For $k\geq 1$, Question \ref{Q23020502} shows that
\[\binom{n}{k}\frac{1}{n^k}=\frac{1}{k!}\frac{n(n-1)\cdots (n-k+1)}{n^k}\leq \frac{1}{k\,!}\leq \frac{1}{2^{k-1}}.\]
Therefore,
\[a_n\leq 1+1+\frac{1}{2}+\cdots+\frac{1}{2^{n-1}}=3-\frac{1}{2^{n-1}}\leq 3.\]
This proves that $\{a_n\}$ is bounded above by 3.
  Since $\{a_n\}$ is an increasing sequence that is bounded above, the monotone convergence theorem asserts that the limit
\[\lim_{n\rightarrow\infty}a_n=\lim_{n\rightarrow\infty} \left(1+\frac{1}{n}\right)^n\] exists. 
 
\end{solution}

\begin{highlight}{The Number $\pmb{e}$}
The number $e$ is defined as  
\[e=\lim_{n\rightarrow\infty} \left(1+\frac{1}{n}\right)^n.\] Correct to 15 decimal places, its numerical value is
\[e=2.718281828459046\]

One can show that the sequence $\{b_n\}_{n=0}^{\infty}$ defined by $b_0=1$, 
\[b_n=b_{n-1}+\frac{1}{n!}\hspace{1cm}\text{for all}\; n\geq 1,\]
also converges to $e$. In series notation, 
\[e=1+\frac{1}{1!}+\frac{1}{2!}+\frac{1}{3!}+\cdots+\frac{1}{n!}+\cdots.\]
\end{highlight}

\vp
\noindent
{\bf \large Exercises  \thesection}
\setcounter{myquestion}{1}

 \begin{question}{\themyquestion} 
Given that the sequence $\{a_n\}$ is defined by $a_1=2$, and for all $n\geq 1$,
\[a_{n+1}=\frac{3a_n+1}{a_n+2}.\]
Show that $\{a_n\}$ is convergent and find its limit.
\end{question}
\atc

 \begin{question}{\themyquestion} 
 For $n\geq 1$, let
\[a_n=\left(1+\frac{1}{n}\right)^n.\]
Define the sequence $\{b_n\}_{n=0}^{\infty}$ by $b_0=1$, and for all $n\geq 1$,
\[b_n=b_{n-1}+\frac{1}{n!}.\]

\begin{enumerate}[(a)]
\item Show that the sequence $\{b_n\}$ is convergent.
\item For a positive integer $n$, use the binomial expansion of $a_n$ to show that 
$a_n\leq b_n$ and 
\[b_n-a_n\leq\frac{3}{2n}.\]
\item Conclude that the sequence $\{b_n\}$ converges to $e$.
\end{enumerate}
 
\end{question}
\vp
\section{Sequential Compactness}\label{sec1.8}

Let us first look at an example.
 \begin{example}{}
 Let $\{a_n\}$ be the sequence defined by
 \[a_n=(-1)^{n-1}\frac{n}{n+1}.\]
 Obviously,
 \[|a_n|=\frac{n}{n+1}\leq 1.\] 
 Hence, the sequence $\{a_n\}$ is bounded. 
 Now,
 \[a_{2n-1}=\frac{2n-1}{2n},\hspace{1cm}a_{2n}=-\frac{2n}{2n+1}.\]
 The subsequence $\{a_{2n-1}\}$ converges to 1, whereas the subsequence $\{a_{2n}\}$ converges to $-1$. Since there are two subsequences that converge to two different limits,
 the sequence $\{a_n\}$ is not convergent.
 
 \end{example}
In this example, we find that although the sequence $\{a_n\}$ is not convergent, it has convergent subsequences.
 In this section, we are going to  prove that every bounded sequence has a convergent subsequence.
By monotone convergence theorem, it is sufficient to prove that every sequence has a monotone subsequence. It can be achieved via a concept called peak index.
 
 \begin{definition}{Peak Index}
 Let $\{a_n\}$ be a sequence of real numbers. A positive integer $m$ is called a {\bf peak index} of the sequence if 
 \[a_m\geq a_n\hspace{1cm}\text{for all}\;n\geq m.\]In other words,   there is no term after the $m^{\text{th}}$ term that is larger than $a_m$.
 \end{definition}
 
 If $\{a_n\}$ is a decreasing sequence, every positive integer is a peak index of the sequence. If $\{a_n\}$ is an increasing sequence, $m$ is a peak index if and only if $a_n=a_m$ for all $n\geq m$, which means $\{a_n\}$ is a constant from the $m^{\text{th}}$ term on. We can use the concept of peak indices to prove the following.
 
 \begin{theorem}{}
 Every sequence has a monotone subsequence.
 \end{theorem}
 \begin{myproof}{Proof}
 Given a sequence $\{a_n\}$, let $S$ be the set of its peak indices. It is a subset of positive integers.  We discuss the cases where $S$ is infinite and $S$ is finite.
 
 \textbf{Case 1:} $S$ is infinite.\\
Let
 $n_1, n_2, n_3, \ldots$ be the elements of $S$  arranged in increasing order, namely,
 \[n_1<n_2<n_3<\cdots.\] This is a subsequence of $\{n\}$. For any positive integer $k$, since $n_{k+1}>n_k$ and $n_k$ is a peak index, we have
 \[a_{n_{k+1}}\leq a_{n_k}.\]
 This shows that $\{a_{n_k}\}$ is a decreasing subsequence of $\{a_n\}$.
 
 \textbf{Case 2:} $S$ is finite.\\If $S$ is an empty set, let $n_1=1$. If $S$ is not empty,   it has a largest element $n_{\max}$. Let $n_1=n_{\max}+1$. Then for any integer $n$ such that $n\geq n_1$, $n$ is not a peak index of the sequence. Since  $n_1$ is not a peak index, there is an $n_2>n_1$ such that
$a_{n_2}>a_{n_1}$.
 Suppose that we have chosen the positive integers $n_1, n_2, \ldots, n_k$ such that $n_1<n_2<\cdots<n_k$ and
 \[a_{n_1}<a_{n_2}<\cdots<a_{n_k}.\]

 Now $n_k$ is   not a peak index implies that there is a positive integer $n_{k+1}$ larger than $n_k$ such that
 \[a_{n_{k+1}}>a_{n_{k}}.\] 
 This procedure constructs the increasing subsequence $\{a_{n_k}\}$ inductively.
 
 In both cases, we have shown that $\{a_n\}$ has a monotone subsequence.
 \end{myproof}
 
 Obviously, a subsequence of a bounded sequence is bounded. It follows from the monotone convergence theorem the following important assertion.
 \begin{theorem}{Bolzano-Weierstrass Theorem}
 Every bounded sequence has a convergent subsequence.
 \end{theorem}

Now we want to introduce a concept called Cauchy sequence, which is closely related to completeness axiom.
\begin{definition}{Cauchy Sequence}
A sequence $\{a_n\}$ is called a {\bf Cauchy sequence} provided that for any $\varepsilon>0$, the is a positive integer $N$ such that for all $m\geq n\geq N$,
\[|a_m-a_n|<\varepsilon.\]

\end{definition}
\begin{example}{}
For the sequence $\{a_n\}$ with $a_n=\di\frac{n+1}{n}$,  it is easy to check that it is a Cauchy sequence. Notice that if $m\geq n$,
\[|a_m-a_n|=\left|\frac{1}{n}-\frac{1}{m}\right|=\frac{1}{n}-\frac{1}{m}<\frac{1}{n}.\]
Given $\varepsilon>0$, the Archimedean property says that there is a positive integer $N$ such that $1/N<\varepsilon$. Hence, if $m\geq n\geq N$,
\[|a_m-a_n|<\frac{1}{n}\leq\frac{1}{N}<\varepsilon.\]
\end{example}

There is a similarity between the definition of a Cauchy sequence  and the definition of convergence of a sequence.  We can show that a linear combination of Cauchy sequences is a Cauchy sequence, and a product of Cauchy sequences is a Cauchy sequence. For the quotient, some care need to be taken.   We leave it to the students to formulate the precise statement.

In the definition of a Cauchy sequence, we do not need to know whether the sequence is convergent, or what is the limit of the sequence if it is convergent.
Nevertheless,   a convergent sequence is a Cauchy sequence.

\begin{theorem}[label=23020602]{}
If a sequence $\{a_n\}$ is convergent, then it is a Cauchy sequence.
\end{theorem}
\begin{myproof}{Proof}
Let $a$ be the limit of the convergent sequence $\{a_n\}$. Given $\varepsilon>0$, there is a positive integer $N$ such that for all $n\geq N$,
\[ |a_n-a|<\frac{\varepsilon}{2}.\]
It follows from triangle inequality that if $m\geq n\geq N$,
\[|a_m-a_n|\leq |a_m-a|+|a_n-a|<\frac{\varepsilon}{2}+\frac{\varepsilon}{2}=\varepsilon.\]
Hence, $\{a_n\}$ is a Cauchy sequence.
\end{myproof}

The converse is also true \emph{in the set of real numbers}. It is proved using the fact that every bounded sequence has a convergent subsequence.

\begin{theorem}[label=23020603]{Cauchy Criterion for Convergent Sequennce}
If $\{a_n\}$ is a Cauchy sequence of real numbers, then it converges to a real number.
\end{theorem}
\begin{myproof}
{Proof}First we prove that   $\{a_n\}$ is a Cauchy sequence implies that it is bounded. The proof is almost identical to the proof that a convergent sequence is bounded.
Take $\varepsilon=1$. There is  a positive integer $N_0$ such that for all $m\geq n\geq N_0$, 
\[|a_m-a_n|<1.\]
This implies that 
\[|a_m|\leq |a_{N_0}|+1\hspace{1cm}\text{for all}\;m\geq {N_0}.\]\bp
Let
\[M=\max\{|a_1|, \ldots, |a_{N_0-1}|, |a_{N_0}|+1.\}\] Then $|a_n|\leq M$ for all $n\in\mathbb{Z}^+$, proving that it is bounded.
Since $\{a_n\}$ is a bounded sequence, it has a convergent subsequence $\{a_{n_k}\}$ which converges to a limit $a$. We want to prove that the sequence $\{a_n\}$ also converges to $a$.
Given $\varepsilon>0$, there is a positive integer $N$ such that for all $m\geq n\geq N$,
\[|a_{m}-a_n|<\frac{\varepsilon}{2}.\]
There is a positive integer $K$ such that for all $k\geq K$,
\[|a_{n_k}-a|<\frac{\varepsilon}{2}.\]
 Now let $n$ be an integer such that $n\geq N$. Since $\{n_k\}$ is a subsequence of $\{n\}$, there is an integer $k$ such that $k\geq K$ and $n_k\geq n$. Then
\[|a_n-a|\leq |a_{n_k}-a_n|+|a_{n_k}-a|<\frac{\varepsilon}{2}+\frac{\varepsilon}{2}=\varepsilon.\]

 This proves that for all $n\geq N$,
\[|a_n-a|<\varepsilon.\] Hence, the sequence $\{a_n\}$ indeed converges to $a$.
\end{myproof}
 Theorem \ref{23020603} is proved using the fact that every bounded sequence has a convergent subsequence. The latter is a consequence of the monotone convergence theorem, whose validity relies on the completeness axiom for real numbers.
 Hence, the fact that every Cauchy sequence of real numbers is convergent is a consequence of the completeness axiom.
 
   If we consider the set of rational numbers, the assertion is not true. For example, we have shown that there is a sequence of rational numbers $\{a_n\}$ that converges to $\sqrt{2}$. Therefore, the sequence $\{a_n\}$ is a Cauchy sequence that does not converge in the set of rational numbers. 

The following combines the results of Theorem \ref{23020602} and Theorem \ref{23020603}.
\begin{highlight}{Cauchy Criterion for Convergent Sequennce}
A sequence of real numbers $\{a_n\}$ is convergent if and only if it is a Cauchy sequence.
\end{highlight}

As the monotone convergence theorem, the Cauchy criterion can be used to conclude the convergence of a sequence without apriori knowing the limit of the sequence. It has wide applications as we are going to see in latter chapters. 
\begin{example}
{}For a positive integer $n$, let
\[s_n=1+\frac{1}{2}+\cdots+\frac{1}{n}.\]
Show that the sequence $\{s_n\}$ is divergent.

\end{example}
\begin{solution}{Solution}
We prove that $\{s_n\}$ is not a Cauchy sequence, by showing that for $\varepsilon=1/2$, for any positive integer $N$, there are integers $m$ and $n$ with $m\geq n\geq N$ such that 
\[|s_m-s_n|\geq \frac{1}{2}.\]
For a given positive integer $N$, let $n=N$ and $m=2N$. Then $m\geq n\geq N$ and $m-n=N$.   Notice that
\begin{align*}s_m-s_n&=\frac{1}{N+1}+\frac{1}{N+2}+\cdots+\frac{1}{2N}\\& \geq \underbrace{\frac{1}{2N}+\frac{1}{2N}+\cdots+\frac{1}{2N}}_{m-n=N\;\text{terms}}\\&=\frac{1}{2}.\end{align*}This shows that $\{s_n\}$ is not a Cauchy sequence. Hence, it is not convergent.
\end{solution}

We have studied the convergence   of sequences, and the interplay between sequences and sets.
Now we define another property of sets called sequential compactness.
 \begin{definition}{Sequential Compactness}
 Let $S$ be a subset of real numbers. We say that $S$ is {\bf sequentially compact}  provided that   every sequence in $S$ has a subsequence that converges to a point in $S$.
 \end{definition}
 Using logic, we find that a set $S$ is not sequentially compact if there is a sequence in $S$ that do not have a convergent subsequence with limit in $S$.
 
 From the theories that we have developed in this chapter, it is not difficult to prove the following.
 \begin{theorem}[label=23020604]{}
 If $S$ is a closed and bounded subset of real numbers, then it is sequentially compact.
 \end{theorem}
 \begin{myproof}{Proof}
 Let $S$ be a subset of $\mathbb{R}$ that is closed and bounded.  Given a sequence $\{a_n\}$ in $S$, since $S$ is bounded,   the sequence $\{a_n\}$ is bounded. Therefore, there is a subsequence $\{a_{n_k}\}$ that converges to a number $a$. Since $\{a_{n_k}\}$ is a sequence in the  set $S$ that converges to  $a$, and $S$ is closed,  the limit $a$ must be in $S$. In other words, we have shown that the sequence $\{a_n\}$ in $S$ has a subsequence $\{a_{n_k}\}$ that converges to  a point $a$ that is in $S$. This proves that $S$ is sequentially compact.
 \end{myproof}
 
 \begin{example}{}
 Since an interval of the form $[a,b]$ is closed and bounded, it is sequentially compact.
 \end{example}
 
 The converse to Theorem \ref{23020604} is also true.
 \begin{theorem}[label=23020707]{}
Let $S$ be a subset of $\mathbb{R}$. If $S$ is sequentially compact, then it is closed and bounded.
 \end{theorem}
 This is a statement of the form $p\to q\wedge r$. It is equivalent to $(p\to q)\wedge (p\to r)$, which in turn is equivalent to $(\neg q\to \neg p)\wedge (\neg r\to \neg p)$. Hence, we will prove the following two statements: if $S$ is not closed, it is not sequentially compact; and if $S$ is not bounded, it is not sequentially compact.
 \begin{myproof}{Proof}
 First, we prove that if $S$ is not closed, it is not sequentially compact. If $S$ is not closed, there is a sequence $\{a_n\}$ in $S$ which converges to a point $a$ but $a$ is not in $S$. For this sequence, every subsequence is convergent with limit $a$. Hence, this sequence does not have a convergent subsequence with limit in $S$. This proves that $S$ is not sequentially compact.
 
 Next, we prove that if $S$ is not bounded, it is not sequentially compact. If $S$ is not bounded, for each integer $n$, there is a point $a_n$ in $S$ such that
 \[|a_n|\geq n.\]Consider the sequence $\{a_n\}$. If $\{a_{n_k}\}$ is a subsequence of $\{a_n\}$,
 \[|a_{n_k}|\geq n_k.\] Hence, the sequence
  $\{a_{n_k}\}$ is not bounded, and thus it is not convergent. This shows that the sequence $\{a_n\}$ does not have any convergent subsequence. Therefore, $S$ is not sequentially compact.

 \end{myproof}
 Combining Theorem \ref{23020604} and Theorem \ref{23020707}, we have the following.
 \begin{highlight}{Characterization of Sequentially Compact Sets}
 A subset of real numbers is sequentially compact if and only if it is closed and bounded.
 \end{highlight}
 
 Notice that the only type of intervals that is both closed and bounded is the type $[a, b]$. Hence, this is the only type of  intervals that are sequentially compact.
 \begin{example}{}
 Determine whether each of the following sets is sequentially compact.
 \begin{enumerate} [(a)]
  \item $\mathbb{Z}$
  \item 
 $A=[2, 5]\setminus\{3\}$

 \item  $B=(0, 6]\cap [4, 7]$.
 \end{enumerate}
 \end{example}
 \begin{solution}{Solution}
  \begin{enumerate}[(a)]
  \item The set $\mathbb{Z}$ is not bounded. Hence, it is not  sequentially compact.
  \item  $3$ is a limit point of the set $A$ but it is not in $A$. Hence, $A$ is not closed, and so it is not  sequentially compact. 
  \item  $B=[4, 6]$ is closed and bounded. Hence, $B$ is sequentially compact.
  \end{enumerate}
 \end{solution}
 
It might be wondered why there is a need to introduce the concept of sequential compactness if it is equivalent to closed and bounded.    We will see that for a subset of real numbers that is closed and bounded, every sequence in that set has a subsequence that converges to a point in that set  is   a very important characteristic. By introducing the concept of sequential compactness, we can avoid repeatedly proving this property for a set that is closed and bounded.

The next theorem gives an important feature of a sequentially compact set.
\begin{theorem}[label=23020908]{}Let $S$ be a subset of real numbers. If $S$ is closed and bounded, then it has a maximum and a minimum. Equivalently, if $S$ is sequentially compact, then it has a maximum and a minimum.
\end{theorem}
\begin{myproof}{Proof}
 Since $S$ is bounded, $S$ has a least upper bound $u$ and a greatest lower bound $\ell$. By Lemma \ref{23020510}, there are sequences $\{u_n\}$ and $\{\ell_n\}$ in $S$ that converge to $u$ and $\ell$ respectively. Since $S$ is closed, $u$ and $\ell$ are in $S$. Since $u=\sup S$ is in $S$, $S$ has a maximum. Since $\ell=\inf S$ is in $S$, $S$ has a minimum.
\end{myproof}
\vp
\noindent
{\bf \large Exercises  \thesection}
\setcounter{myquestion}{1}
 \begin{question}{\themyquestion} 
Given that the sequence $\{a_n\}$ is defined by
\[a_n=1+\frac{1}{3}+\ldots+\frac{1}{2n-1}.\]
Show that $\{a_n\}$ is not a Cauchy sequence. Then conclude that the sequence $\{a_n\}$ is divergent.
\end{question}
 
 \atc

 \begin{question}{\themyquestion} 
Determine whether each of the following sequence is a Cauchy sequence.
\begin{enumerate}[(a)]
\item The sequence $\{a_n\}$ with $a_n=\di \frac{n+(-1)^n}{n-(-1)^n}$
\item The sequence $\{b_n\}$ with $b_n=\di \frac{1+  n}{1-(-1)^n n}$

\end{enumerate}
\end{question}
\atc
 \begin{question}{\themyquestion} 
 Determine whether each of the following sets is sequentially compact.
 \begin{enumerate} [(a)]
  \item $A=\{1, 2, \cdots, 100\}$
  \item 
 $B=[4, 7]\cap (6, 8]$

 \end{enumerate}
\end{question}
 
\atc
 \begin{question}{\themyquestion} 
Show that the union of two sequentially compact sets is sequentially compact. 
\end{question}

\chapter{Limits of   Functions and Continuity}\label{ch2}

 In this chapter, we study functions $f:D\rightarrow \mathbb{R}$ defined on a subset of real numbers $D$, and taking values in the set of real numbers $\mathbb{R}$. 
 Polynomials and rational functions are special examples.  When we do not specify the domain of a function, we will take its domain $D$ to be the largest subset of real numbers where the function can be defined.
 
 \begin{definition}{Polynomials and Rational Functions}
 A polynomial is a function $p:\mathbb{R}\rightarrow\mathbb{R}$ of the form
 \[p(x)=a_nx^n+a_{n-1}x^{n-1}+\ldots+a_1x+a_0,\]
 where $a_0, a_1, \ldots, a_n$ are constants. We call $p(x)$ a polynomial of degree $n$ if $a_n\neq 0$. 
 A rational function is a function of the form 
 \[f(x)=\frac{p(x)}{q(x)},\]
 where $p(x)$ and $q(x)$ are polynomials, and $q(x)$ is not the zero polynomial. The domain of this function is the set $D=\mathbb{R}\setminus S$, where $S$ is the finite point set containing all $x$ for which $q(x)=0$.
 \end{definition}
 
 For example, the domain of the rational function \[ f(x)=\frac{x^2+1}{x+2}\] is the set $D=\mathbb{R}\setminus\{-2\}$.

 To be able to apply tools in analysis, we are   interested in functions that are continuous. 
 Continuity can be defined in two different  ways that are equivalent. One is using positive numbers $\delta$ and $\varepsilon$ to measure distances of points in the domain and range, while the other is using limits of sequences. 
 
 The limit of a function $f:D\rightarrow \mathbb{R}$ when the variable $x$ approaches a limit point $x_0$ of the domain $D$ is an important concept in defining derivatives.  This concept can be defined for any function $f:D\rightarrow \mathbb{R}$ whose domain $D$ contains limit points. 
 There is a close relation between the limit of a function $f(x)$ when $x$ approaches a limit point $x_0$, and the continuity of the function at $x_0$. 
 
 Although the continuity of a function can be defined independently of limits of functions, we choose to consider limits of functions first. 
 
\section{Limits of Functions }\label{sec2.1 }

 In Section \ref{sec1.6}, we have defined the concept of limit points of a set $D$. The point $x_0$ is a limit point of the set $D$ if there is a sequence of points in $D\setminus\{x_0\}$ that converges to $x_0$. A limit point of a set is not necessarily in that set. A set that contains all its limit points is a closed set. If a point $x_0$ is in a set $D$ but is not a limit point of $D$, it is called an isolated point of $D$. If $x_0$ is an isolated point of $D$, there is a neighbourhood $(a,b)$ of $x_0$ which intersects the set $D$ only at the point $x_0$.
 
 Limits of functions can be defined using the $\varepsilon-\delta$ language or using limits of sequences. We will define the concept using limits of sequences first, and then show that it is equivalent to the $\varepsilon-\delta$  definition.
 
 \begin{definition} 
 {Limits of   Functions}
 Let $D$ be a subset of real numbers and let $x_0$ be a limit point of $D$. Given a function $f:D\rightarrow \mathbb{R}$, we say that {\it the limit of $f(x)$ as $x$ approaches $x_0$ is $\ell$}, provided that whenever $\{x_n\}$ is a sequence of points in $D\setminus\{x_0\}$ that converges to $x_0$, the sequence $\{f(x_n)\}$ converges to $\ell$. 
 
If the limit of $f:D\rightarrow \mathbb{R}$ as $x$ approaches $x_0$ is $\ell$, we write
 \[\lim_{x\rightarrow x_0}f(x)=\ell.\]
 \end{definition}
 Notice that we do not define $\di\lim_{x\rightarrow x_0}f(x)$ if $x_0$ is not a limit point of the domain where the function is defined. 
  
 \begin{highlight}{Logical Expression for Definition of  Limits}
 ~\vspace{-0.5cm}
 \begin{align*}
 & \lim_{x\rightarrow x_0}f(x)=\ell
 \iff \\& \forall  \{x_n\} \subset   D\setminus\{x_0\},\;\lim_{n\rightarrow \infty}x_n=x_0 \,\implies\,\lim_{n\rightarrow\infty}f(x_n)=\ell.\end{align*}
 \end{highlight}
 
 Let us first look at some examples. 
 \begin{example}{}
 Find the  limit if it exists.
 \begin{enumerate}[(a)]
 \item
 $\di\lim_{x\rightarrow 1}\frac{2x+3}{x^2+1}$
 
 \item $\di \lim_{x\rightarrow 1}\frac{x^2-1}{x-1}$
 
 \item $\di \lim_{x\rightarrow 1}\frac{x^2+1}{x-1}$
 \end{enumerate}
 \end{example}
 \begin{solution}{Solution}
 \begin{enumerate}[(a)]
 \item The function \[f(x)= \frac{2x+3}{x^2+1}\] is defined on $D=\mathbb{R}$. If $\{x_n\}$ is a sequence in $\mathbb{R}\setminus \{1\}$ that converges to 1, limit laws imply that
 \[\lim_{n\rightarrow\infty}f(x_n)=\lim_{n\rightarrow \infty}\frac{2x_n+3}{x_n^2+1}=\frac{2\times 1+3}{1^2+1}=\frac{5}{2}.\]Hence,
 \[\lim_{x\rightarrow 1}\frac{2x+3}{x^2+1}=\frac{5}{2}.\]\end{enumerate} \bs
 \begin{enumerate}[(b)] 
  \item The function \[f(x)= \frac{x^2-1}{x-1}\] is defined on $D=\mathbb{R}\setminus\{1\}$. If $\{x_n\}$ is a sequence in $\mathbb{R}\setminus \{1\}$ that converges to 1, limit laws imply that
 \[\lim_{n\rightarrow\infty}f(x_n)=\lim_{n\rightarrow \infty}\frac{x_n^2-1}{x_n-1}=\lim_{n\rightarrow \infty}(x_n+1)=2.\] Hence,
 \[\lim_{x\rightarrow 1}\frac{x^2-1}{x-1}=2.\]\end{enumerate}
 \begin{enumerate}[(c)]
 \item  The function \[f(x)= \frac{x^2+1}{x-1}\] is defined on $D=\mathbb{R}\setminus\{1\}$. Consider the sequence $\{x_n\}$ with 
 \[x_n=1+\frac{1}{n}.\]
 We find that \[f(x_n)=\frac{\di \frac{1}{n^2}+\frac{2}{n}+2}{\di \frac{1}{n}}=2n+2+\frac{1}{n}.\]The sequence $\{f(x_n)\}$ is not bounded, and so it is divergent. Hence, the limit
  \[\lim_{x\rightarrow 1}\frac{x^2+1}{x-1} \]does not exist.
 \end{enumerate}
 \end{solution}
 In part $(b)$, we have used the fact that $x_n\neq 1$  to simplify $f(x_n)$ to $x_n+1$. 
 
 Using   laws for limits of sequences, it is immediate to  see that limits of functions respect taking linear combinations and multiplications. It also respects taking quotients provided that the function on the denominator does not approach 0.
 
 \begin{proposition}[label=23020813]{Limit Laws for Functions}
Let $D$ be a subset of real numbers. Given that $f:D\rightarrow\mathbb{R}$ and $g:D\rightarrow\mathbb{R}$ are functions defined on $D$,   $x_0$ is a limit point of $D$,  and
\[\lim_{x\rightarrow x_0}f(x)=\ell_1,\hspace{1cm}\lim_{x\rightarrow x_0}g(x)=\ell_2.\]
\begin{enumerate}[1.]
\item
For any constants $\alpha$ and $\beta$, $\di \lim_{x\rightarrow x_0}(\alpha f+\beta g)(x)=\alpha\ell_1+\beta \ell_2$.
\item $\di \lim_{x\rightarrow x_0}(  f  g)(x)= \ell_1  \ell_2$.
\item If  $g(x)\neq 0$ for all $x\in D$, and $\ell_2\neq 0$, then
\[\lim_{x\rightarrow x_0}\frac{f(x)}{g(x)}=\frac{\ell_1}{\ell_2}.\]
\end{enumerate}
 \end{proposition}
 
 From this proposition, it follows that we can take limits of a rational function easily at a point which is not a zero of the polynomial in the denominator.
 \begin{proposition}[label=23020807]{}
 Let $p(x)$ and $q(x)$ be polynomials. If $x_0$ is a real number such that
 $q(x_0)\neq 0$, then
 \[\lim_{x\rightarrow x_0}\frac{p(x)}{q(x)}=\frac{p(x_0)}{q(x_0)}.\]
 \end{proposition}

 Let us now look at an example that involves the absolute values.
 \begin{example}[label=23020803]{}
 Show that for any real number $x_0$, 
 \[\lim_{x\rightarrow x_0}|x|=|x_0|.\]

\end{example}
\begin{solution}{Solution}
Let $\{x_n\}$ be a sequence in $\mathbb{R}\setminus\{x_0\}$ that converges to $x_0$. By Question \ref{absolute}, the sequence $\{|x_n|\}$ converges to $|x_0|$. This proves that
\[\lim_{x\rightarrow x_0}|x|=|x_0|.\]
\end{solution}

 Now we want to formulate an equivalent definition for limits.
 \begin{theorem}[label=23020801]{Equivalent Definitions for Limits}
  Let $D$ be a subset of real numbers, and let $x_0$ be a limit point of $D$. Given a function $f:D\rightarrow \mathbb{R}$, 
  the following are two equivalent definitions for 
  \[\lim_{x\rightarrow x_0}f(x)=\ell.\]  
  \begin{enumerate}[(i)]
  \item 
  Whenever $\{x_n\}$ is a sequence of points in $D\setminus\{x_0\}$ that converges to $x_0$, the sequence $\{f(x_n)\}$ converges to $\ell$. 
  \item For any $\varepsilon>0$, there is a $\delta>0$ such that if the point $x$ is in $D$ and $0<|x-x_0|<\delta$, then $|f(x)-\ell|<\varepsilon$.
  \end{enumerate} 
 \end{theorem}
 
 In logical notation, we can express (ii) as follows.
  \begin{highlight}{Logical Expression for Second Definition of  Limits}
 ~\vspace{-0.5cm}
 \begin{align*}
 & \lim_{x\rightarrow x_0}f(x)=\ell 
 \iff \\& \forall  \varepsilon>0, \;\exists \delta>0,\;\forall x\; (x\in D)\,\wedge (0<|x-x_0|<\delta)\;\implies \;|f(x)-\ell|<\varepsilon.\end{align*}
 \end{highlight}
 Here $\delta$ is a measure of the closeness of the point $x\in D$ to the point $x_0$, and $\varepsilon$ is a measure of the closeness of the function value $f(x)$ to the number $\ell$. The condition $|x-x_0|>0$ is to stress that we only consider those points $x$ that is not $x_0$. From the definitions, we can see that the limit of a function $f(x)$ when $x$ approaches $x_0$ does not depend on how the function $f $ is defined at $x_0$, and  $f$   does  not need to be defined at $x_0$ for the limit to be defined.

 To prove Theorem \ref{23020801}, we need to show that (i) $\iff$ (ii). This is equivalent to (ii) $\implies$ (i) and $\neg$ (ii) $\implies$ $\neg$ (i). 
 \begin{myproof}{\linkt Proof of Theorem \ref{23020801}}
 We start by  showing that if (ii) holds, then (i) holds. Assume that (ii) holds. To prove (i), we take a sequence $\{x_n\}$ in $D\setminus\{x_0\}$ that converges to the point $x_0$. We want to show that the sequence $\{f(x_n)\}$ converges to $\ell$. We prove this using the definition of convergence of sequences. 
  Given $\varepsilon>0$, our assumption that (ii) holds implies that there is a $\delta>0$ such that for all $x$ that is in $D$ with $0<|x-x_0|<\delta$, we have $|f(x)-\ell|<\varepsilon$. Since $\{x_n\}$ converges to $x_0$, there is a positive integer $N$ such that for all $n\geq N$, $|x_n-x_0|<\delta$. By our definition, $x_n$ are points in $D\setminus \{x_0\}$. Hence, for all $n\geq N$,  
 \[|f(x_n)-\ell|<\varepsilon.\]Since we have shown that for all $\varepsilon>0$, there is a positive integer $N$ such that for all $n\geq N$, $|f(x_n)-\ell|<\varepsilon$, we conclude that the sequence $\{f(x_n)\}$ converges to $\ell$. This proves (i) holds.
 
Now assume that (ii) is false.  In logical notation, this means
\[\exists  \varepsilon>0, \;\forall \delta>0,\;\exists x\; (x\in D)\,\wedge (0<|x-x_0|<\delta)\;\wedge\;|f(x)-\ell|\geq \varepsilon.\] 
Namely, there is an $\varepsilon>0$ such that for any $\delta>0$, there is a point $x$ in $D\setminus\{x_0\}$ with $|x-x_0|<\delta$ but $|f(x)-\ell|\geq \varepsilon$.
For this $\varepsilon>0$, we  construct a sequence $\{x_n\}$ in $D\setminus \{x_0\}$ in the following way. For each positive integer $n$, there is a point $x_n$ in $D\setminus\{x_0\}$ such that
$|x_n-x_0|<1/n$ but $|f(x_n)-\ell|\geq\varepsilon$.
Then $\{x_n\}$ is a sequence in $D\setminus\{x_0\}$ that satisfies
\[|x_n-x_0|<\frac{1}{n}.\]
Since $\di\lim_{n\rightarrow \infty}1/n=0$, we find that the sequence $\{x_n\}$ converges to $x_0$. Since $|f(x_n)-\ell|\geq \varepsilon$ for all $n\in\mathbb{Z}^+$, the sequence $\{f(x_n)\}$ cannot converge to $\ell$. Hence, we have shown that there is a sequence $\{x_n\}$ in $D\setminus \{x_0\}$ that converges to $x_0$ but $\{f(x_n)\}$ does not converge to $\ell$. This proves that (i) does not hold.
 \end{myproof}
 
 \begin{example}[label=23020808]
 {The Heaviside Function}The Heaviside function  $H:\mathbb{R}\rightarrow\mathbb{R}$ is  defined by
 \begin{align*}
 H(x)=\begin{cases}1,\quad &\text{if}\; x\geq 0;
 \\0,\quad & \text{if}\;x<0.\end{cases}
 \end{align*}
 For any  real number $x_0$, determine whether the limit $\di\lim_{x\rightarrow x_0}H(x)$ exists.
 \end{example}
 
 \begin{figure}[ht]
\centering
\includegraphics[scale=0.2]{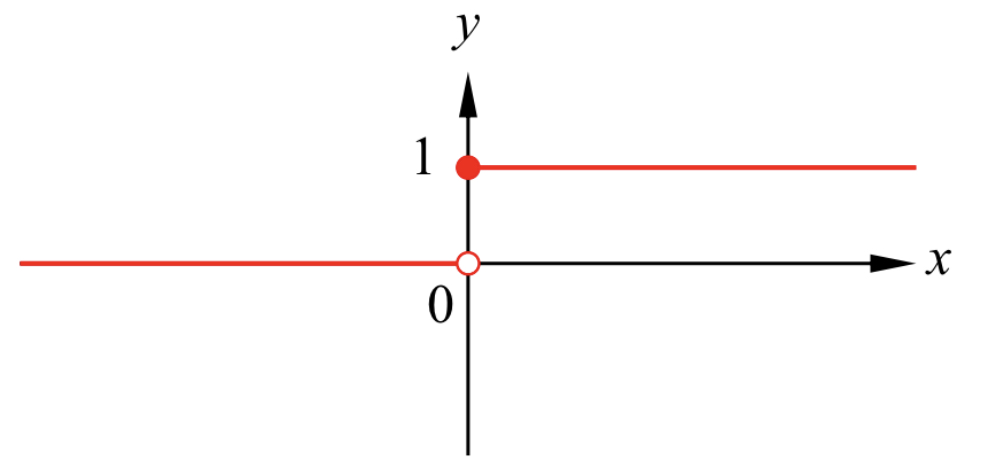}
\caption{  The Heaviside function $H(x)$.}\label{figure6}
\end{figure}

 \begin{solution}{Solution}
 We consider the cases where $x_0>0$, $x_0<0$ and $x_0=0$.
 
  \textbf{Case 1:} $x_0>0$.
 In this case, we claim that $\di\lim_{x\rightarrow x_0}f(x)=1$. \\Given $\varepsilon>0$, take $\delta=x_0$. Then $\delta>0$. If $x$ is in $\mathbb{R}$ and $0<|x-x_0|<\delta=x_0$, we have
 $x-x_0>-x_0$ and hence $x>0$. Thus $f(x)=1$ and 
 \[|f(x)-1|=0<\varepsilon.\]
 This proves that  $\di\lim_{x\rightarrow x_0}f(x)=1$.
 
 \textbf{Case 2:} $x_0<0$.
 In this case, we claim that $\di\lim_{x\rightarrow x_0}f(x)=0$. \\Given $\varepsilon>0$, take $\delta=-x_0$. Then $\delta>0$. If $x$ is in $\mathbb{R}$ and $0<|x-x_0|<\delta=-x_0$, we have
 $x-x_0<-x_0$ and hence $x<0$. Thus $f(x)=0$ and 
 \[|f(x)-0|=0<\varepsilon.\]
 This proves that  $\di\lim_{x\rightarrow x_0}f(x)=0$. 
 \bs
  \textbf{Case 3:} $x_0=0$.
 In this case, we claim that $\di\lim_{x\rightarrow x_0}f(x)$ does not exist. \
 Let $\{u_n\}$ and $\{v_n\}$ be the sequences $\{1/n\}$ and $\{-1/n\}$ respectively. They are both sequences in $\mathbb{R}\setminus\{0\}$ that converge to 0.
 \[f(u_n)=1 \quad\text{and}\quad f(v_n)=0 \hspace{1cm}\text{for all}\;n\in\mathbb{Z}^+.\]
 Therefore,
 \[\lim_{n\rightarrow\infty}f(u_n)=1\hspace{1cm}\text{and}\hspace{1cm}\lim_{n\rightarrow\infty}f(v_n)=0.\]
 Since $\{f(x_n)\}$ has different limits when we consider two different sequences $\{x_n\}$ in $\mathbb{R}\setminus\{0\}$ that converge to 0, we conclude that   $\di\lim_{x\rightarrow x_0}f(x)$ does not exist. 
 \end{solution}
 
 In this example, we can also use the $\varepsilon-\delta$ definition to show that  $\di\lim_{x\rightarrow  0}f(x)$ does not exist.  Assume that $\di\lim_{x\rightarrow  0}f(x)$ exists and is equal to $\ell$. Take $\varepsilon=1/2$. There exists $\delta>0$ such that for any $x\in \mathbb{R}$, if $0<|x-0|<\delta$, then
 \[|f(x)-\ell|<\varepsilon.\]Now the points $x=x_1=-\delta/2$ and $x=x_2=\delta/2$ both satisfy 
$0<|x-0|<\delta$. We have $f(x_1)=0$ and $f(x_2)=1$.
By triangle inequality,
\[|f(x_1)-f(x_2)|\leq |f(x_1)-\ell|+|f(x_2)-\ell|<2\varepsilon=1.\]
This gives
\[1=|f(x_1)-f(x_2)|<1,\]which is a contradiction. Hence, $\di\lim_{x\rightarrow 0}f(x)$ does not exist.

  In calculus, we have defined  the concepts of  left limits and right limits to deal with functions like the Heaviside function, which is defined by cases.
  Given a subset of real numbers $D$ and a point $x_0$, define
  \[D_{x_0, -}=\left\{x\in D\,|\,x< x_0\right\},\hspace{1cm} D_{x_0, +}=\left\{x\in D\,|\,x> x_0\right\}.\]
  For example, consider $D=[0, 2)$. If  $x_0=1$, then $D_{1,-}=[0,1)$ and $D_{1,+}=(1,2)$. If $x_0=0$, then $D_{0,-}=\emptyset$ and $D_{0,+}=(0,2)$. If $x_0=2$, then $D_{2,-}=[0,2)$ and $D_{2,+}=\emptyset$.
  
  Notice that even though $x_0$ is a limit point of $D$, it might not be a limit point of $D_{x_0,-}$ or $D_{x_0,+}$. 
  We define the left limit and right limit of a function $f:D\rightarrow\mathbb{R}$ when $x$ approaches $x_0$ in the following way.
  
  \begin{definition}{Left Limits and Right Limits}
   Let   $D$ be a subset of real numbers and let $f:D\rightarrow \mathbb{R}$ be a function defined on   $D$.  
   \begin{enumerate}[1.]
   \item If $x_0$ is a limit point of $D_{x_0,-}$, $D_{x_0,-}$ is not an empty set.  We say that the   limit of the function $f:D\rightarrow \mathbb{R}$ as $x$ approaches $x_0$ from the left exists provided that the limit of the function $f:D_{x_0,-}\rightarrow\mathbb{R}$ as $x$ approaches $x_0$ exists. If the left limit exists, it is denoted by
   \[\lim_{x\rightarrow x_0, x<x_0}f(x)\hspace{1cm}\text{or simply as}\hspace{1cm}\lim_{x\rightarrow x_0^-}f(x).\]
      \item If $x_0$ is a limit point of $D_{x_0,+}$, $D_{x_0,+}$ is not an empty set.  We say that the   limit of the function $f:D\rightarrow \mathbb{R}$ as $x$ approaches $x_0$ from the right exists provided that the limit of the function $f:D_{x_0,+}\rightarrow\mathbb{R}$ as $x$ approaches $x_0$ exists. If the right limit exists, it is denoted by
   \[\lim_{x\rightarrow x_0, x>x_0}f(x)\hspace{1cm}\text{or simply as}\hspace{1cm}\lim_{x\rightarrow x_0^+}f(x).\]
   \end{enumerate}
  
  \end{definition}
  
\begin{highlight}{Left Limits, Right Limits, and Limits}
\begin{enumerate}[1.]
\item If $x_0$ is a limit point of both $D_{x_0,+}$ and $D_{x_0,-}$, then $\di\lim_{x\rightarrow x_0}f(x)$ exists if and only if both $\di\lim_{x\rightarrow x_0^-}f(x)$ and $\di\lim_{x\rightarrow x_0^+}f(x)$ exist and they are equal. 
\item 
If $x_0$ is a limit point of $D_{x_0,-}$ but is not a limit point of $D_{x_0,+}$, then  $\di\lim_{x\rightarrow x_0}f(x)$ exists if and only if   $\di\lim_{x\rightarrow x_0^-}f(x)$ exists.
\item 
If $x_0$ is a limit point of $D_{x_0,+}$ but is not a limit point of $D_{x_0,-}$, then  $\di\lim_{x\rightarrow x_0}f(x)$ exists if and only if   $\di\lim_{x\rightarrow x_0^+}f(x)$ exists.
\end{enumerate}\end{highlight}
\begin{example}{}
For the Heaviside function, we have
  \[\lim_{x\rightarrow 0^-}H(x)=0\hspace{1cm}\text{and}\hspace{1cm}\lim_{x\rightarrow 0^+}H(x)=1.\]Since the left and right limits are not equal, $\di \lim_{x\to 0}H(x)$ dos not exist.\end{example}

 \begin{example}[label=23020809]
 {The Dirichlet's Function} The Dirichlet's function is the function $f:\mathbb{R}\rightarrow\mathbb{R}$ defined by
 \begin{align*}
 f(x)=\begin{cases}1,\quad &\text{if}\; x\;\text{is rational},
 \\0,\quad &\text{if}\; x\;\text{is irrational}.\end{cases}
 \end{align*}
 For any  real number $x_0$, determine whether the limit $\di\lim_{x\rightarrow x_0}f(x)$ exists.
 \end{example}This is a classical example of a function which we cannot visualize the graph.
 \begin{solution}{Solution}
 Fixed a real number $x_0$. For any positive integer $n$, there is a rational number $p_n$ and an irrational number $q_n$ in the open interval $(x_0-1/n, x_0)$. The sequences $\{p_n\}$ and $\{q_n\}$ are in $\mathbb{R}\setminus\{x_0\}$ and converge to $x_0$. Since
 \[f(p_n)=1\quad\text{and}\quad f(q_n)=0\hspace{1cm}\text{for all}\;n\in\mathbb{Z}^+,\]
 we find that 
 \[\lim_{n\rightarrow\infty}f(p_n)=1\hspace{1cm}\text{and}\hspace{1cm}\lim_{n\rightarrow\infty}f(q_n)=0.\]
 Since the sequence $\{f(x_n)\}$ has different limits when we consider two different sequences $\{x_n\}$ in $\mathbb{R}\setminus\{x_0\}$ that converge to $x_0$, we conclude that   $\di\lim_{x\rightarrow x_0}f(x)$ does not exist. 
 \end{solution}
 For this example, if one wants to use the $\varepsilon-\delta$ definition of limits, one can proceed in the following way. For fixed $x_0$ in $\mathbb{R}$, assume that  $\di\lim_{x\rightarrow x_0}f(x)=\ell$. When $\varepsilon=1/2$, there is a $\delta>0$ such that for any $x$ with $0<|x-x_0|<\delta$, $|f(x)-\ell|<\varepsilon$. The open interval $(x_0-\delta, x_0)$ contains a rational number $x_1$  and an irrational number $x_2$. Notice that $f(x_1)=1$ and $f(x_2)=0$.  Both $x=x_1$ and $x=x_2$ satisfy $0<|x-x_0|<\delta$. By triangle inequality,
\[|f(x_1)-f(x_2)|\leq |f(x_1)-\ell|+|f(x_2)-\ell|<2\varepsilon=1.\]
This gives
\[1=|f(x_1)-f(x_2)|<1,\]which is a contradiction. Hence, $\di\lim_{x\rightarrow x_0}f(x)$ does not exist.
 
Next, we consider composite functions. 
\begin{proposition}[label=23020815]{}
Given the two functions $f:D\rightarrow \mathbb{R}$ and $g: U\rightarrow\mathbb{R}$, if $f(D)\subset U$, we can define the composite function $h=g\circ f:D\rightarrow \mathbb{R}$ by
$h(x)=g(f(x))$. If $x_0$ is a limit point of $D$, $y_0$ is a limit point of $U$, $f(D\setminus\{x_0\})\subset U\setminus\{y_0\}$,
\[\lim_{x\rightarrow x_0}f(x)=y_0,\hspace{1cm}\lim_{y\rightarrow y_0}g(y)=\ell,\] then
\[\lim_{x\rightarrow x_0}h(x)=\lim_{x\rightarrow x_0}(g\circ f)(x)=\ell.\]
\end{proposition}
\begin{myproof}{Proof}
Let $\{x_n\}$ be a sequence in $D\setminus\{x_0\}$ that converges to $x_0$, and let $y_n=f(x_n)$ for all $n\in\mathbb{Z}^+$. By assumption,   $\{y_n\}$ is a sequence in $U\setminus\{y_0\}$.  Since $\di \lim_{x\rightarrow x_0}f(x)=y_0$, the sequence $\{f(x_n)\}$ converges to $y_0$. Since $\di\lim_{y\rightarrow y_0}g(y)=\ell$, the sequence
$\{g(y_n)\}$ converges to $\ell$. In other words, the sequence $\{(g\circ f)(x_n)\}$ converges to $\ell$.

Since we have proved that whenever  $\{x_n\}$ is a sequence in $D\setminus\{x_0\}$ that converges to $x_0$, the sequence $\{(g\circ f)(x_n)\}$ converges to $\ell$, we conclude that
\[\lim_{x\rightarrow x_0}(g\circ f)(x)=\ell.\]
\end{myproof}

Using the result of Example \ref{23020803}, we obtain the following.
\begin{corollary}{}
Let $D$ be a subset of real numbers. Given a function $f:D\rightarrow\mathbb{R}$, if $x_0$ is a limit point of $D$ and 
$\di \lim_{x\rightarrow x_0}f(x)=\ell$,
then
\[\lim_{x\rightarrow x_0}|f(x)|=|\ell|.\]
\end{corollary}

%Now let us consider the square-root function $f(x)=\sqrt{x}$, $x\geq 0$.
\begin{example}[label=23020901]{}
For any $x_0\geq 0$, show that
\[\lim_{x\rightarrow x_0}\sqrt{x}=\sqrt{x_0}.\]
\end{example}
\begin{solution}{Solution}
Let us use the $\varepsilon-\delta$ definition of limits. 
Consider the case $x_0=0$ first. Given $\varepsilon>0$, take $\delta=\varepsilon^2$. Then $\delta>0$. If $x\geq 0$ is such that $0<|x-0|<\delta=\varepsilon^2$, we have $0<x<\varepsilon^2$, which implies that
$\di 0< \sqrt{x}<\varepsilon$. Hence,  if $x\geq 0$ and $0<|x-0|<\delta$,
\[|\sqrt{x}-\sqrt{0}|<\varepsilon.\] This proves  that
\[\lim_{x\rightarrow 0}\sqrt{x}=0=\sqrt{x_0}.\]

Now consider the case $x_0>0$. 
Notice that
\[\sqrt{x}-\sqrt{x_0}=\frac{x-x_0}{\sqrt{x}+\sqrt{x_0}}.\]
 If $x>x_0/4$, then $\sqrt{x}>\sqrt{x_0}/2$ and
 \[\frac{1}{\sqrt{x}+\sqrt{x_0}}<\frac{2}{3\sqrt{x_0}}.\]
Given $\varepsilon>0$, let
$\di \delta=\min\left\{\frac{3}{4}x_0,\; \frac{3}{2}\varepsilon\sqrt{x_0} \right\}$.
Then $\delta>0$. If $x\geq 0$ and $0<|x-x_0|<\delta$, then $\di |x-x_0|<\frac{3}{4}x_0$ and so $\di x>\frac{1}{4}x_0$. 
Therefore,\bs
\[\left|\sqrt{x}-\sqrt{x_0}\right|=\frac{|x-x_0|}{\sqrt{x}+\sqrt{x_0}}<\delta\times \frac{2}{3\sqrt{x_0}}\leq\varepsilon.\]
This proves that 
$\di \lim_{x\rightarrow x_0}\sqrt{x} =\sqrt{x_0}$.
\end{solution}

\begin{figure}[ht]
\centering
\includegraphics[scale=0.2]{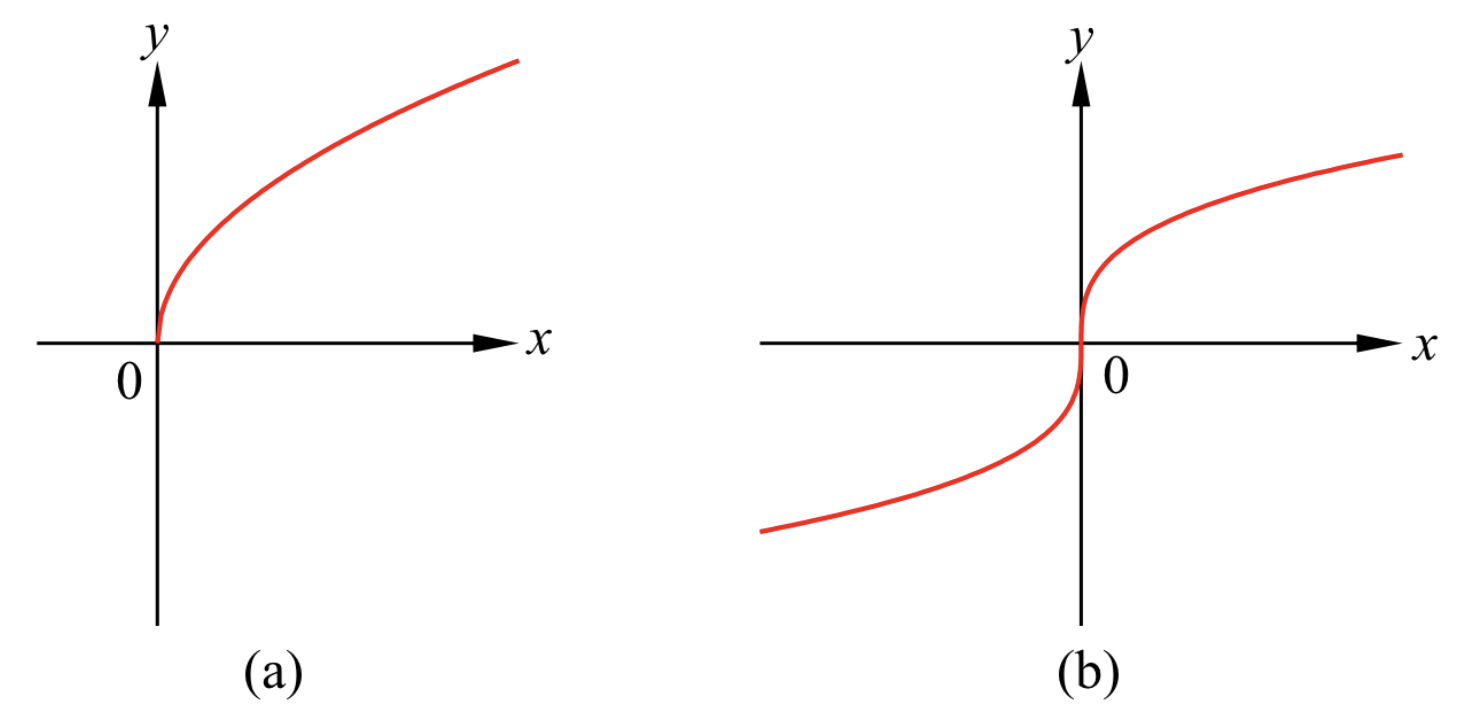}
\caption{  (a) The function $f(x)=\sqrt{x}$. (b) The function $f(x)=\sqrt[3]{x}$.}\label{figure7}
\end{figure}

Using similar methods, one can prove that  if $n$ is an   integer,   then for any $x_0$ in the domain of the function $f(x)=\sqrt[n]{x}$,
\[\lim_{x\rightarrow x_0}\sqrt[n]{x}= \sqrt[n]{x_0}.\]

Now we want to give a brief discussion about limits that involve infinities.  
\begin{definition}{Infinity as Limits of Sequences}
Given that $\{a_n\}$ is a sequence of real numbers.
\begin{enumerate}[1.]
\item
We say that the sequence $\{a_n\}$ diverges to $\infty$, written as
$\di \lim_{n\rightarrow \infty}a_n=\infty$,
if for every positive number $M$, there is a positive integer $N$ such that for all $n\geq N$,
$a_n\geq M$.
\item
We say that the sequence $\{a_n\}$ diverges to $-\infty$, written as
$\di \lim_{n\rightarrow \infty}a_n=-\infty$,
if for every positive number $M$, there is a positive integer $N$ such that for all $n\geq N$,
$a_n\leq -M$.
\end{enumerate}\end{definition}
 
\begin{example}{}
\begin{enumerate}[(a)]
\item
The sequence $\{n^2\}$ diverges to $\infty$.
\item The sequence $\{-n^2\}$ diverges to $-\infty$.
\item The sequence $\{(-1)^nn^2\}$ neither diverges to $\infty$ nor to $-\infty$.
\end{enumerate}
\end{example}
%Notice that the sequence $\{a_n\}$ diverges to $-\infty$ if and only if the sequence $\{-a_n\}$ diverges to $\infty$.

Given that $\{x_n\}$ is a sequence of  real numbers. If $\{x_n\}$ diverges to $\infty$ or $-\infty$, there is a positive integer $N$ such that $x_n\neq 0$ for all $n\geq N$. Hence, for sequences that diverge to $\infty$ and $-\infty$, we can assume none of the terms is zero.

The following is another characterization of boundedness for a set in terms of sequences that diverge to infinity.
\begin{highlight}{}

\begin{enumerate}[1.]
\item A set $D$ is not bounded above if and only if there is a sequence $\{x_n\}$ in $D$ that diverges to $\infty$. 
\item A set $D$ is not bounded below if and only if there is a sequence $\{x_n\}$ in $D$ that diverges to $-\infty$. 
\end{enumerate}
\end{highlight}

Using these, we can make the following definitions.
  
\begin{definition}{Limits of Functions at Infinity}
Let $D$ be a subset of real numbers that is not bounded above. Given  that $\ell$ is a real number and $f:D\rightarrow\mathbb{R}$  is a function  defined on the set $D$.  The following  two   definitions for
\[\lim_{x\rightarrow\infty}f(x)=\ell \]are equivalent.
\begin{enumerate}[(i)]
\item Whenever $\{x_n\}$ is a sequence of points in $D$ that diverges to $\infty$, the sequence $\{f(x_n)\}$ converges to  $\ell$. 
\item For any $\varepsilon>0$, there is a positive number $M  $ such that if the point $x$ is in $D$ and $x>M$, then
\[|f(x)-\ell|<\varepsilon.\]
\end{enumerate}
\end{definition}

\begin{definition}{Limits of Functions at Negative Infinity}
Let $D$ be a subset of real numbers that is not bounded below. Given  that $\ell$ is a real number and $f:D\rightarrow\mathbb{R}$  is a function  defined on the set $D$.  The following  are two equivalent definitions for
\[\lim_{x\rightarrow-\infty}f(x)=\ell.\]
\begin{enumerate}[(i)]
\item Whenever $\{x_n\}$ is a sequence of points in $D$ that diverges to $-\infty$, the sequence $\{f(x_n)\}$ converges to  $\ell$. 
\item For any $\varepsilon>0$, there is a positive number $M  $ such that if the point $x$ is in $D$ and $x<-M$, then
\[|f(x)-\ell|<\varepsilon.\]
\end{enumerate}
\end{definition}

Now let us look at a simple example.
\begin{example}[label=23020903]{}
Show that $\di\lim_{x\rightarrow\infty}\frac{1}{x}=0$.
\end{example}
\begin{solution}{Solution}
We use both definitions to prove the statement.

Using the sequence definition, let $\{x_n\}$ be a sequence of nonzero real numbers that diverges to $\infty$. We want to show that the sequence $\{1/x_n\}$ converges to 0. Given $\varepsilon>0$, the number $M=1/\varepsilon$ is also positive. Since the sequence $\{x_n\}$ diverges to $\infty$, there is a positive integer $N$ such that for all $n\geq N$, 
\[x_n>M=\frac{1}{\varepsilon}.\] In particular, for all $n\geq N$, $x_n>0$ and
$\di 0<\frac{1}{x_n}<\varepsilon$. This proves that the sequence $\{1/x_n\}$ converges to 0. Therefore, $\di \lim_{x\rightarrow\infty}\frac{1}{x}=0$.

Now consider the definition in terms of $\varepsilon$. Given $\varepsilon>0$, let $M=1/\varepsilon$. Then $M$ is a positive number. If $x$ in $\mathbb{R}\setminus\{0\}$ is such that $x>M$, then
\[0<\frac{1}{x}<\frac{1}{M}=\varepsilon.\]
This proves that $\di\lim_{x\rightarrow\infty}\frac{1}{x}=0$.
\end{solution}

This example demonstrates that working with the definition in terms of $\varepsilon$ is sometimes easier.

 \begin{figure}[ht]
\centering
\includegraphics[scale=0.2]{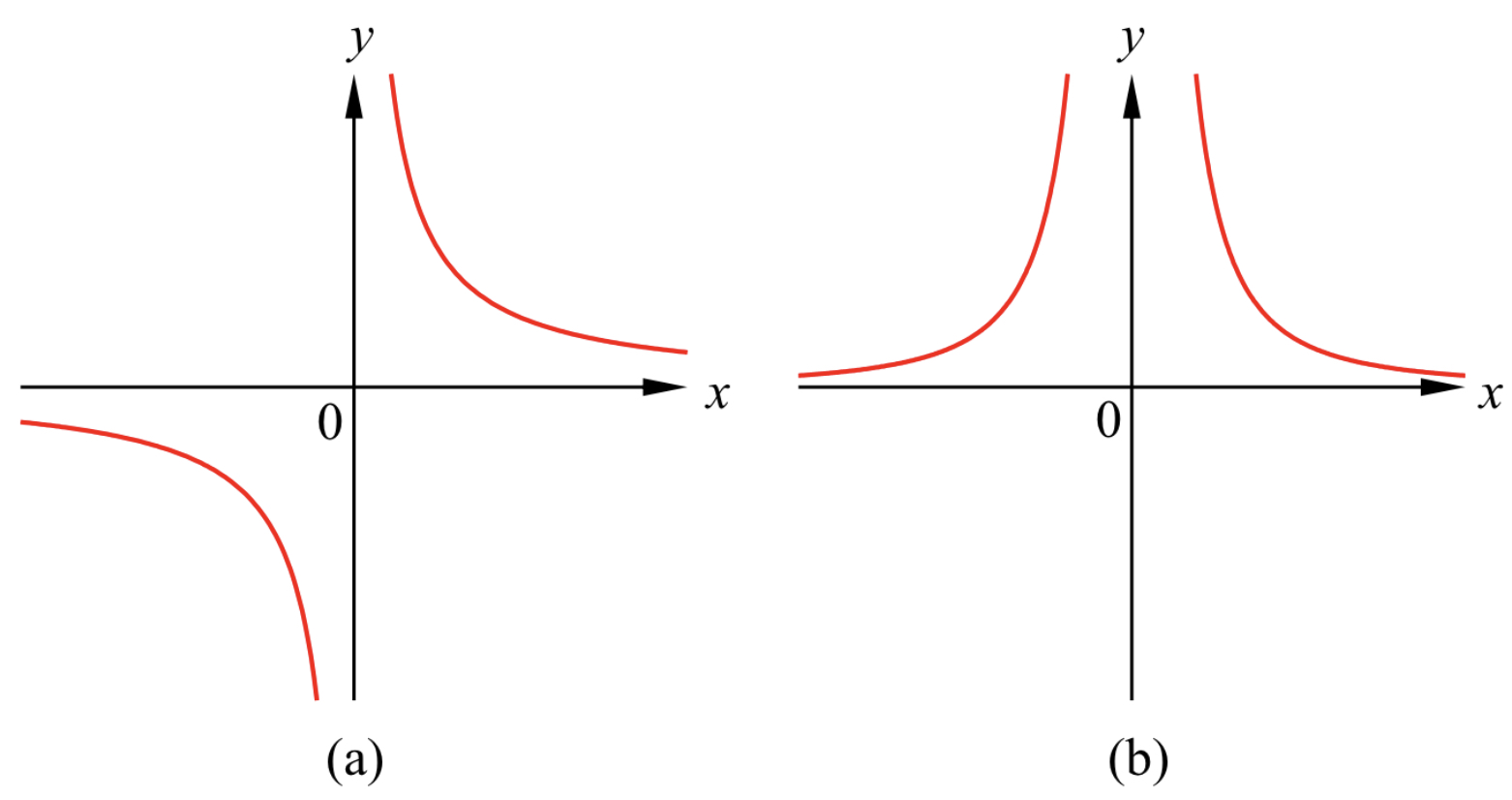}
\caption{  (a) The function $f(x)=1/x$. (b) The function $f(x)=1/x^2$.}\label{figure8}
\end{figure}

It is easy to see that the limit laws given in Proposition \ref{23020813} and Proposition \ref{23020815} also hold for the case where $x\rightarrow \infty$ or $x\rightarrow -\infty$. We will skip the formulation and use it directly. For example, we have the following.

\begin{example}{}
For any positive integer $n$, 
$\di \lim_{x\rightarrow\infty}\frac{1}{x^n}=0$.
\end{example}

Now let us look at some more examples.
\begin{example}{}
Determine whether the limit 
\[ \lim_{x\rightarrow\infty}\frac{2x^2+3x+4}{x^2+7}\]exists. If it  exists, find the limit.
 
\end{example}
\begin{solution}{Solution}
  Divide the numerator and denominator by $x^2$, we have
\[\frac{2x^2+3x+4}{x^2+7}=\frac{2+\di\frac{3}{x}+\frac{4}{x^2}}{1+\di\frac{7}{x^2}}.\]
Using limit laws and the fact that $\di\lim_{x\rightarrow\infty}1/x=0$, we find that
\[\lim_{x\rightarrow\infty}\frac{2x^2+3x+4}{x^2+7}=\frac{2+0+0}{1+0}=2.\]
 
\end{solution}

\begin{example}{}
Determine whether the limit 
\[  \lim_{x\rightarrow -\infty}\frac{x}{\sqrt{x^2+1}}\]exists. If it  exists, find the limit.

\end{example}

\begin{solution}{Solution}
  Notice that $\sqrt{x^2}=|x|$. Hence,
\[\frac{x}{\sqrt{x^2+1}}=\frac{x}{|x|\sqrt{1+\di\frac{1}{x^2}}}.\] 
When $x<0$, 
\[\frac{x}{|x|}=-1.\] 
Therefore,
\[\lim_{x\rightarrow-\infty}\frac{x}{|x|}=-1.\]On the other hand, since 
\[\lim_{x\rightarrow\infty}1+\di\frac{1}{x^2}=1\hspace{1cm}\text{and}\hspace{1cm}\lim_{y\rightarrow 1}\sqrt{y}=1,\] we find that
\[\lim_{x\rightarrow\infty}\frac{1}{\sqrt{1+\di\frac{1}{x^2}}}=1.\]Hence,
\[\lim_{x\rightarrow -\infty}\frac{x}{\sqrt{x^2+1}}=-1.\]

\end{solution}

 \begin{figure}[ht]
\centering
\includegraphics[scale=0.2]{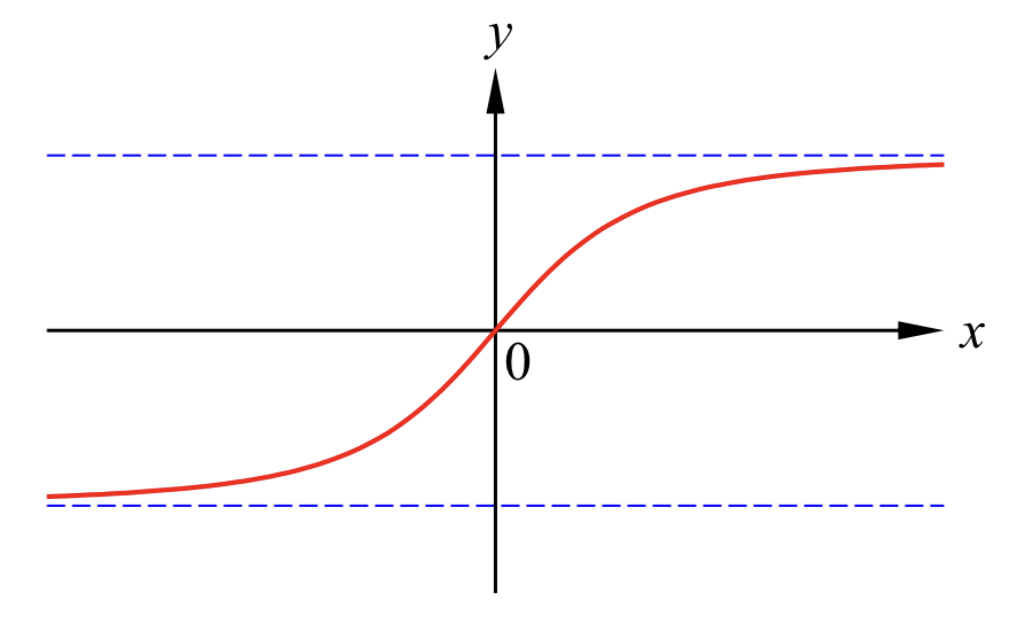}
\caption{  The function $\di f(x)=\frac{x}{\sqrt{x^2+1}}$.}\label{figure9}
\end{figure}

\begin{remark}{}
Using similar ideas, one can formulate analogous definitions for the following limits.
\[
 \lim_{x\rightarrow x_0}f(x)=\infty,\hspace{1cm}
 \lim_{x\rightarrow x_0}f(x)=-\infty,\]
\[\lim_{x\rightarrow \infty}f(x)=\infty,\hspace{1cm}
 \lim_{x\rightarrow \infty}f(x)=-\infty,\]
 \[
 \lim_{x\rightarrow -\infty}f(x)=\infty,\hspace{1cm}
 \lim_{x\rightarrow -\infty}f(x)=-\infty.\]
 
\end{remark}

Finally, we would like to mention that there is an analogue of the squeeze theorem for functions, whose proof is straightforward.
\begin{theorem}{Squeeze Theorem}
Let $D$ be a subset of real numbers. Given that   $f:D\rightarrow\mathbb{R}$, $g:D\rightarrow\mathbb{R}$, $h:D\rightarrow\mathbb{R}$ are functions defined on $D$ and
\[g(x)\leq f(x)\leq h(x)\hspace{1cm}\text{for all}\;x\in D.\]
If $x_0$ is a limit point of $D$ and
\[\lim_{x\rightarrow x_0}g(x)=\lim_{x\rightarrow x_0}h(x)=\ell,\]then
\[\lim_{x\rightarrow x_0}f(x)=\ell.\]
\end{theorem}

\sbr
\noindent
{\bf \large Exercises  \thesection}
\setcounter{myquestion}{1}

 \begin{question}{\themyquestion}
  Find the limit if it exists.
 \begin{enumerate}[(a)]
 \item
 $\di\lim_{x\rightarrow -1}\frac{x^2+3x+2}{x^2+1}$
 
 \item $\di \lim_{x\rightarrow -1}\frac{x^2-3x-4}{x+1}$
 
 \item $\di \lim_{x\rightarrow -1}\frac{x^2+1}{x+1}$
 
 \item  $\di \lim_{x\rightarrow -1}\left|\frac{x^2-3x-4}{x+1}\right|$
 \end{enumerate}
\end{question}
\atc
\begin{question}{\themyquestion}
  Define the function $f:\mathbb{R}\rightarrow\mathbb{R}$   by
 \begin{align*}
 f(x)=\begin{cases}x+3,\quad & \text{if}\; x\geq 1,
 \\5-x,\quad & \text{if}\;x<1.\end{cases}
 \end{align*}
 For any  real number $x_0$, determine whether the limit $\di\lim_{x\rightarrow x_0}f(x)$ exists.
\end{question}

\atc
\begin{question}{\themyquestion}
 Define the function $f:\mathbb{R}\rightarrow\mathbb{R}$   by
 \begin{align*}
 f(x)=\begin{cases}x,\quad & \text{if}\; x\;\text{is rational};
 \\-x,\quad & \text{if}\;x\;\text{is irrational}.\end{cases}
 \end{align*}
 \begin{enumerate}
 [(a)]
 \item Use squeeze theorem to show that  $\di\lim_{x\rightarrow 0}f(x)$ exists and find the limit.
 \item If $x_0\neq 0$, show that  $\di\lim_{x\rightarrow x_0}f(x)$ does not exist.
 \end{enumerate}
\end{question}

\atc

\begin{question}{\themyquestion}
Determine whether the limit exists. If it  exists, find the limit.
\begin{enumerate}[(a)]
\item[(a)]
$\di \lim_{x\rightarrow -\infty}\frac{2x^2+x+4}{5x^2+2}$
\item [(b)]$\di \lim_{x\rightarrow \infty}\frac{2x+3}{\sqrt{4x^2+1}}$
\end{enumerate}
\end{question}

\atc
\begin{question}{\themyquestion}
For any $x_0\geq 0$, show that
\[\lim_{x\rightarrow x_0}\sqrt[4]{x}=\sqrt[4]{x_0}.\]
\end{question}

\vp

\section{Continuity of Functions}\label{sec2.2}

In this section, we introduce the concept of continuity of functions.

\begin{definition}{Continuity}
Let $D$ be a subset of real numbers that contains the point $x_0$, and let $f:D\rightarrow\mathbb{R}$ be a function defined on $D$. We say that the function $f$ is {\bf continuous at } $x_0$ provided that whenever $\{x_n\}$ is a sequence of points in $D$ that converges to $x_0$, the sequence $\{f(x_n)\}$ converges to $f(x_0)$. 

We say that $f:D\rightarrow \mathbb{R}$ is a \textbf{continuous function} if it is continuous at every point of its domain $D$.
\end{definition}

The  definitions of limit and continuity are very similar. However, there is a slight difference. To define continuity at a point $x_0$, $x_0$ must be a point in the domain of the function $D$. To define limit, $x_0$ does not need to be a point in the domain $D$ but has to be a limit point of $D$.  When the point  $x_0$ is in $D$ and is also a limit point of $D$, the relation between limit and continuity is as follows.
\begin{proposition}
{Relation Between Limit and Continuity}
Let $D$ be a subset of real numbers that contains the point $x_0$. If $x_0$ is a limit point of $D$, then $f$ is continuous at $x_0$ if and only if
\[\lim_{x\rightarrow x_0}f(x)=f(x_0).\]
\end{proposition}

In other words, it says that if $x_0$ is a limit point of the domain $D$, then $f$ is continuous at $x_0$ if and only if 
\[\lim_{x\rightarrow x_0}(f(x)-f(x_0))=0,\]if and only if 
\[\lim_{x\rightarrow x_0}|f(x)-f(x_0)|=0.\]
The following fact is  quite obvious.
\begin{proposition}{}
Let $D$ be a subset of real numbers and let $f:D\rightarrow\mathbb{R}$ be a function defined on $D$. If $f:D\rightarrow \mathbb{R}$ is continuous, then for any subset $A$ of $D$, the function $f:A\rightarrow\mathbb{R}$, which is the restriction of $f$ to $A$, is also continuous.
\end{proposition}

\begin{example}{}
Proposition \ref{23020807} says that a rational function is continuous.
\end{example}

\begin{example}{}
Example \ref{23020808} says that the Heaviside function $H(x)$ is continuous at $x$ if $x\neq 0$. It is not continuous at $x=0$. 
\end{example}

\begin{example}{}
Example \ref{23020809} says that the Dirichlet's function is nowhere continuous.
\end{example}

\begin{example}{}
Example \ref{23020901} says that the function $f(x)=\sqrt{x}$ is continuous. 
In general, for any positive integer $n$, the function $f(x)=\sqrt[n]{x}$ is continuous.
\end{example}

A natural question to ask is the continuity of a function at an isolated point of its domain.  Let us first prove the following.
\begin{lemma}
[label=23020811]{}
Let $D$ be a subset of real numbers and let $x_0$ be an isolated point of $D$. If $\{x_n\}$ is a sequence of points in $D$ that converges to $x_0$, then there is a positive integer $n$ such that $x_n=x_0$ for all $n\geq N$.
\end{lemma}
\begin{myproof}{Proof}
By Theorem \ref{23020810},  there is a neighbourhood $(a, b)$ of $x_0$ which intersects $D$   at $x_0$ only.  
 Let $\varepsilon=\min\{x_0-a, b-x_0\}$. Then $\varepsilon>0$ and $ (x_0-\varepsilon, x_0+ \varepsilon)\subset (a,b)$. Since the sequence $\{x_n\}$ converges to $x_0$, there is a positive integer $N$ such that for all $n\geq N$, $|x_n-x_0|<\varepsilon$. Hence, for all $n\geq N$, $x_n\in (x_0-\varepsilon, x_0+ \varepsilon)\subset (a,b)$. Since $(a, b)\cap D=\{x_0\}$, we find that  $x_n=x_0$ for all $n\geq N$.
 
\end{myproof}

Using this lemma, it is easy to prove the continuity of a function at an isolated point of its domain. 
\begin{proposition}
{Continuity at an Isolated Point}
Let $D$ be a subset of real numbers that contains the point $x_0$. If $x_0$ is an isolated point of  $D$, then $f$ is continuous at $x_0$.
\end{proposition}
\begin{myproof}{Proof}
If $\{x_n\}$ is a sequence in $D$ that converge to $x_0$, Lemma \ref{23020811} says that there is a positive integer $N$ such that $x_n=x_0$ for all $n\geq N$.  Therefore, $f(x_n)=f(x_0)$ for all $n\geq N$.  This implies that the sequence $\{f(x_n)\}$ converges to $f(x_0)$. By the definition of continuity, $f$ is continuous at $x_0$.
\end{myproof}

\begin{example}{}
Since every point of the set of positive integers $\mathbb{Z}^+$ is an isolated point, any function $f:\mathbb{Z}^+\rightarrow \mathbb{R}$ defined on the set of positive integers is continuous.
\end{example}
This conclusion might be a bit counter intuitive for students that see it for the first time. One can think about it naively in the following way. For an isolated point, it has no close neighbours to be compared to. Hence, the limit operation does not work, and thus the function is continuous by default.

Let us summarize again the continuity of a function at a point.
\begin{highlight}{Continuity of a Function at a Point}
Let $D$ be a subset of real numbers and let $x_0$ be a point in $D$. Given that $f:D\rightarrow \mathbb{R}$ is a function defined on $D$.
\begin{enumerate}[1.]
\item
If $x_0$ is an isolated point of $D$, then $f$ is continuous at $x_0$.
\item If $x_0$  is a limit point of $D$, then $f$ is continuous at $x_0$ if and only if 
\[\lim_{x\rightarrow x_0}f(x)=f(x_0).\]
\end{enumerate}
\end{highlight}

Similar to limits, we also have an equivalent definition for continuity in terms of $\delta$ and $\varepsilon$.
 \begin{theorem}[label=23020812]{Equivalent Definitions for Continuity}
  Let $D$ be a subset of real numbers and let $x_0$ be   point in $D$. Given a function $f:D\rightarrow \mathbb{R}$, 
  the following two definitions for 
   $f$ to be continuous at $x_0$ are equivalent.
  \begin{enumerate}[(i)]
  \item 
  Whenever $\{x_n\}$ is a sequence of points in $D $ that converges to $x_0$, the sequence $\{f(x_n)\}$ converges to $f(x_0)$. 
  \item For any $\varepsilon>0$, there is a $\delta>0$ such that if the point $x$ is in $D$ and $|x-x_0|<\delta$, then $|f(x)-f(x_0)|<\varepsilon$.
  \end{enumerate} 
 \end{theorem}
 The proof of Theorem \ref{23020812} is almost identical to the proof of Theorem \ref{23020801}.

\begin{example}{}
Use the $\varepsilon-\delta$ definition to show that the function $f:\mathbb{R}\setminus\{0\}\rightarrow \mathbb{R}$ defined by $f(x)=\di \frac{1}{x}$ is continuous.
\end{example}
 \begin{solution}{Solution}
The domain of the function is $D=\mathbb{R}\setminus\{0\}$. Let $x_0$ be a point in $D$. Then $x_0\neq 0$.  
Notice that \bs
\begin{equation}\label{eq230209_1}\left|f(x)-f(x_0)\right|=\left|\frac{1}{x}-\frac{1}{x_0}\right|=\frac{|x-x_0|}{|x||x_0|}.\end{equation}
If \[|x-x_0|<\frac{|x_0|}{2},\]
then
\[|x|> \frac{|x_0|}{2}>0.\]
Given $\varepsilon>0$, let
\[\delta=\min\left\{ \frac{|x_0|}{2}, \frac{|x_0|^2}{2}\varepsilon\right\}.\]If $x$ in $D$ is such that $|x-x_0|<\delta$, then $\di |x-x_0|<\frac{|x_0|}{2}$ and so $\di |x|> \frac{|x_0|}{2}$. It follows from \eqref{eq230209_1} that
\[\left|f(x)-f(x_0)\right|<\delta\times \frac{2}{|x_0|^2}\leq\varepsilon.\]
This proves that $f$ is continuous at $x_0$.
\end{solution}
 
 From Proposition \ref{23020813}, it follows immediately that continuity is preserved when we perform certain operations on functions.
  \begin{proposition}[label=23020814]{}
Let $D$ be a subset of real numbers that contains the point $x_0$. Given that the functions $f:D\rightarrow\mathbb{R}$ and $g:D\rightarrow\mathbb{R}$ are  continuous at $x_0$.
\begin{enumerate}[1.]
\item
For any constants $\alpha$ and $\beta$, the function $\alpha f+\beta g: D\rightarrow\mathbb{R}$ is continuous at $x_0$.
\item The function $(f  g):D\rightarrow\mathbb{R}$ is continuous at $x_0$. 
\item If  $g(x)\neq 0$ for all $x\in D$,  
then the function $(f/g):D\rightarrow\mathbb{R}$ is continuous at $x_0$.
\end{enumerate}
 \end{proposition}
 
 For composition of functions, we have the following which is a counterpart of Proposition \ref{23020815}.
 \begin{proposition}[label=23020816]{}
Given the two functions $f:D\rightarrow \mathbb{R}$ and $g: U\rightarrow\mathbb{R}$ with $f(D)\subset U$. If $x_0$ is a  point of $D$,  $f$ is continuous at $x_0$, $g$ is continuous at $y_0=f(x_0)$, then
the composite function $(g\circ f):D\rightarrow \mathbb{R}$ is continuous at $x_0$.
\end{proposition}
This proposition can be proved easily using definition of continuity in terms of convergent sequences.

\begin{corollary}{}
Let $D$ be a subset of real numbers that contains the point $x_0$. If the function $f:D\rightarrow \mathbb{R}$ is continuous at $x_0$, then the function $|f|:D\rightarrow \mathbb{R}$ is also continuous at $x_0$.
\end{corollary}

Let us now look at an example of a piecewise function.
\begin{example}[label=23020902]{}
Let $f:[-2,3]\rightarrow\mathbb{R}$ be the function defined by
\[f(x)=\begin{cases} 2x^2-3,\quad &\text{if}\;-2\leq x\leq 1,\\
cx+2,\quad &\text{if}\;\;\;1<x\leq 3.\end{cases}\]
Show that there is a value of $c$ for which $f$ is a continuous function.
\end{example}

\begin{solution}{Solution}
The domain of the function $f$ is $D=[-2, 3]$.
First we show that if $x_0\in D\setminus \{1\}$, then $f$ is continuous at $x_0$. 

If $x_0\in [-2, 1)$, then $x_0<1$. If $\{x_n\}$ is a sequence in $D\setminus\{x_0\}$ that converges to $x_0$, then there is a positive integer $N$ such that $x_n<1$ for all $n\geq N$. This implies that for all $n\geq N$, $f(x_n)=2x_n^2-3$. Hence, the sequence $\{f(x_n)\}$ converges to $f(x_0)=2x_0^2-3$. This proves that $f$ is continuous at $x_0$.

Using similar arguments, we can show that if $x_0\in (1, 3]$,  $f$ is continuous at $x_0$. \bs

Now, by definitions of left limits and right limits,
\[\lim_{x\rightarrow 1^-}f(x)=\lim_{x\rightarrow 1^-}(2x^2-3)=-1;\]
\[\lim_{x\rightarrow 1^+}f(x)=\lim_{x\rightarrow 1^+}(cx+2)=c+2.\]
For $f$ to be continuous at $x_0=1$, $\di\lim_{x\rightarrow 1}f(x)$ must exist. So we must have
\[\lim_{x\rightarrow 1^-}f(x)=\lim_{x\rightarrow 1^+}f(x).\]
This gives $c=-3$. 
In fact, 
when $c=-3$, 
\[\lim_{x\rightarrow 1}f(x)=-1=f(1),\]and hence $f$ is continuous at $x=1$.
\end{solution}

 \begin{figure}[ht]
\centering
\includegraphics[scale=0.2]{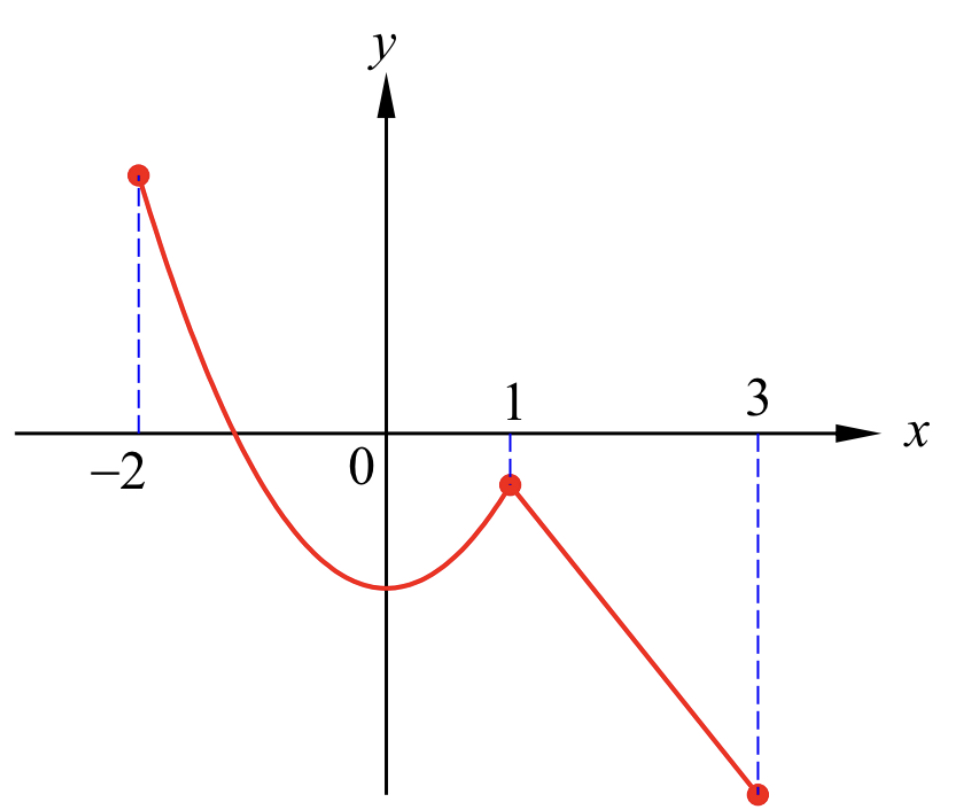}
\caption{  The function defined in Example \ref{23020902}.}\label{figure10}
\end{figure}

\begin{remark}{}
We can formulate a general theorem from Example \ref{23020902} as follows. 

Let $D$ be a subset of real numbers that contains the point $x_0$, and let $D_{x_0,-}$ and $D_{x_0,+}$ be the intersection of $D$ with the sets $\{x\,|\,x<x_0\}$ and $\{x\,|\,x>x_0\}$ respectively. Suppose that $x_0$ is a limit point of both $D_{x_0,-}$ and $D_{x_0, +}$, and $f:D\rightarrow\mathbb{R}$ is a function such that  its restrictions to $D_{x_0,-}$ and $D_{x_0, +}$ are continuous. If
\[\lim_{x\rightarrow x_0^-}f(x)=\lim_{x\rightarrow  x_0^+}f(x)=f(x_0),\] then $f:D\rightarrow\mathbb{R}$ is a continuous function.
\end{remark}

Finally, we define a special class of continuous functions called the Lipschitz function.

\begin{definition}{Lipschitz Function}
Let $D$ be a subset of real numbers. A function $f:D\rightarrow \mathbb{R}$ is said to be a Lipschitz function if there is a  constant $c$ such that
\[|f(x_1)-f(x_2)|\leq c|x_1-x_2|\hspace{1cm}\text{for all}\;x_1, x_2\in D.\]
The constant $c$ is called a Lipschitz constant of the function.
\end{definition}

Notice that a Lipschitz constant is nonnegative. The only  Lipschitz functions with 0 Lipschitz constant  are the constant functions. If $c_0$ is a Lipschitz constant of a Lipschitz function $f:D\rightarrow \mathbb{R}$, any number $c$ that is larger than $c_0$ is also a Lipschitz constant of $f$.

\begin{example}{}
Let $f:\mathbb{R}\rightarrow\mathbb{R}$ be the function given by $f(x)=ax+b$. Then $f$ is a Lipschitz function with Lipschitz constant $|a|$.
\end{example}

\begin{example}[label=23021007]{}
Let $f:[-10, 8]\rightarrow \mathbb{R}$ be the function defined by $f(x)=x^2$. Show that $f$ is Lipschitz.
\end{example}
\begin{solution}{Solution}
For any $x_1$ and $x_2$ in $[-10, 8]$, 
\[|f(x_1)-f(x_2)|=|x_1^2-x_2^2|=|x_1+x_2||x_1-x_2|.\]
Triangle inequality implies that
\[|x_1+x_2|\leq |x_1|+|x_2|\leq 10+10=20.\]
Hence,
\[|f(x_1)-f(x_2)|\leq 20|x_1-x_2|.\]
This shows that $f$ is a Lipschitz function with Lipschitz constant $20$.
\end{solution}

\begin{example}{}
Let $f:\mathbb{R}\rightarrow \mathbb{R}$ be the function defined by $f(x)=x^2$. Is $f$ a Lipschitz function?
\end{example}
\begin{solution}{Solution}
If $f$ is a Lipschitz function, there is a positive constant $c$ such that
\[|f(x_1)-f(x_2)|\leq c|x_1-x_2|\] for all real numbers $x_1$ and $x_2$. 
 Take   $x_1=c+1$ and $x_2=0$. We find that
\[(c+1)^2=|f(x_1)-f(x_2)|\leq c|x_1-x_2|= c(c+1),\] which implies that $c+1\leq 0$, a contradiction.
Hence, $f$ is not a Lipschitz function.
\end{solution}

Here we see that whether a function is Lipschitz or not depends on the domain.
Finally we prove that a Lipschitz function is continuous. 

\begin{theorem}[label=23021005]{}
Let $D$ be a subset of real numbers. If $f:D\rightarrow\mathbb{R}$ is a Lipschitz function, then it is continuous.
\end{theorem}
\begin{myproof}{Proof}
Since $f:D\rightarrow\mathbb{R}$ is Lipschitz,  there is a positive constant $c$ such that for any $x_1$ and $x_2$ in $D$,
\[|f(x_1)-f(x_2)|\leq c|x_1-x_2|.\]Let $x_0$ be a point in $D$. Given $\varepsilon>0$, take $\delta =\varepsilon/c$. Then $\delta>0$ and for any $x\in D$, if $|x-x_0|<\delta$,
\[|f(x)-f(x_0)|\leq c|x-x_0|<c\delta=\varepsilon.\] This proves that $f$ is continuous at $x_0$. Hence, $f:D\rightarrow\mathbb{R}$ is a continuous function.
\end{myproof}
\vp

\noindent
{\bf \large Exercises  \thesection}
\setcounter{myquestion}{1}

 \begin{question}{\themyquestion}
 Consider the function $f:\mathbb{R}\rightarrow\mathbb{R}$ defined by
 \[f(x)=\begin{cases} 2,\quad &\text{if}\;x>2,\\
 x,\quad & \text{if}\;x\leq 2\end{cases}.\]
 Show that $f$ is a continuous function.
\end{question}
\atc
 \begin{question}{\themyquestion}
  Consider the function $f:\mathbb{R}\rightarrow\mathbb{R}$ defined by
 \[f(x)=\begin{cases} x^2,\quad &\text{if $x$ is rational},\\
 -x^2,\quad & \text{if $x$ is irrational}.\end{cases}\]Show that $f$ is continuous at $x=0$.
 
\end{question}
\atc
 \begin{question}{\themyquestion}
  Consider the function $f:\mathbb{R}\rightarrow\mathbb{R}$ defined by
 \[f(x)=\begin{cases} 2x+5,\quad &\text{if}\;x<-1,\\
 ax^2+x,\quad & \text{if}\;x\geq -1\end{cases}.\]
 Show that there is a value of $a$ for which $f$ is a continuous function.
\end{question}
\atc
 \begin{question}{\themyquestion}
 Let $f:[-7, 5]\rightarrow\mathbb{R}$ be the function defined by $f(x)=2x^2+3x$. Show that $f$ is  a Lipschitz function.
\end{question}
\atc
 
 \begin{question}{\themyquestion}
Let $f:[1, \infty)\rightarrow \mathbb{R}$ be the function defined by $f(x)=\sqrt{x}$. Show that $f$ is a Lipschitz function.
\end{question}
\vp

\section{The Extreme Value Theorem}\label{sec2.3}
For a real-valued function $f:D\rightarrow\mathbb{R}$, the maximum value   is the largest  value the function can assume; while the minimum value  is the smallest value the function can assume.
\begin{definition}{Maximium and Minimum Values of a Function}
Let $D$  be a subset of real numbers. Given that $f:D\rightarrow\mathbb{R}$ is a  real-valued function  defined on $D$.
\begin{enumerate}[1.]
\item  $f$  has a maximum value if and only if there is a point $x_0$ in $D$ such that
\[f(x)\leq f(x_0)\hspace{1cm}\text{for all}\;x\in D.\]
Such a  $x_0$ is called a {\bf maximizer} of the function $f$.
\item  $f$  has a minimum value if and only if there is a point $x_0$ in $D$ such that
\[f(x)\geq f(x_0)\hspace{1cm}\text{for all}\;x\in D.\]
Such a  $x_0$ is called a {\bf minimizer} of the function $f$.
\end{enumerate}
\end{definition}

\begin{highlight}{Extreme Values}The maximum value of a function $f:D\rightarrow \mathbb{R}$ is the maximum of the set $f(D)$; while the minimum value is the minimum of the set $f(D)$.

A maximum value or a minimum value of a function  is called an {\bf extreme value} of the function. 
\end{highlight}

\begin{example}
[label=23020905]{}
\begin{enumerate}[(a)]
\item For the function $f:[-1, 2]\rightarrow \mathbb{R}$, $f(x)=2x$, $D=[-1,2]$ and $f(D)=[-2, 4]$. Thus, $f$ has minimum value $-2$ and maximum value 4.
\item For the function $g:[-1, 2)\rightarrow \mathbb{R}$, $g(x)=2x$, $D=[-1, 2)$ and $g(D)=[-2, 4)$. Thus, $g$ has minimum value $-2$, but it does not have maximum value.\end{enumerate}\be\begin{enumerate}[(a)]
\item[(c)] For the function $h:(-1, 2]\rightarrow \mathbb{R}$, $h(x)=2x$, $D=(-1, 2]$ and $h(D)=(-2, 4]$. Thus, $h$ has  maximum value 4, but it does not have minimum value.
\end{enumerate}
\end{example2}

 \begin{figure}[ht]
\centering
\includegraphics[scale=0.2]{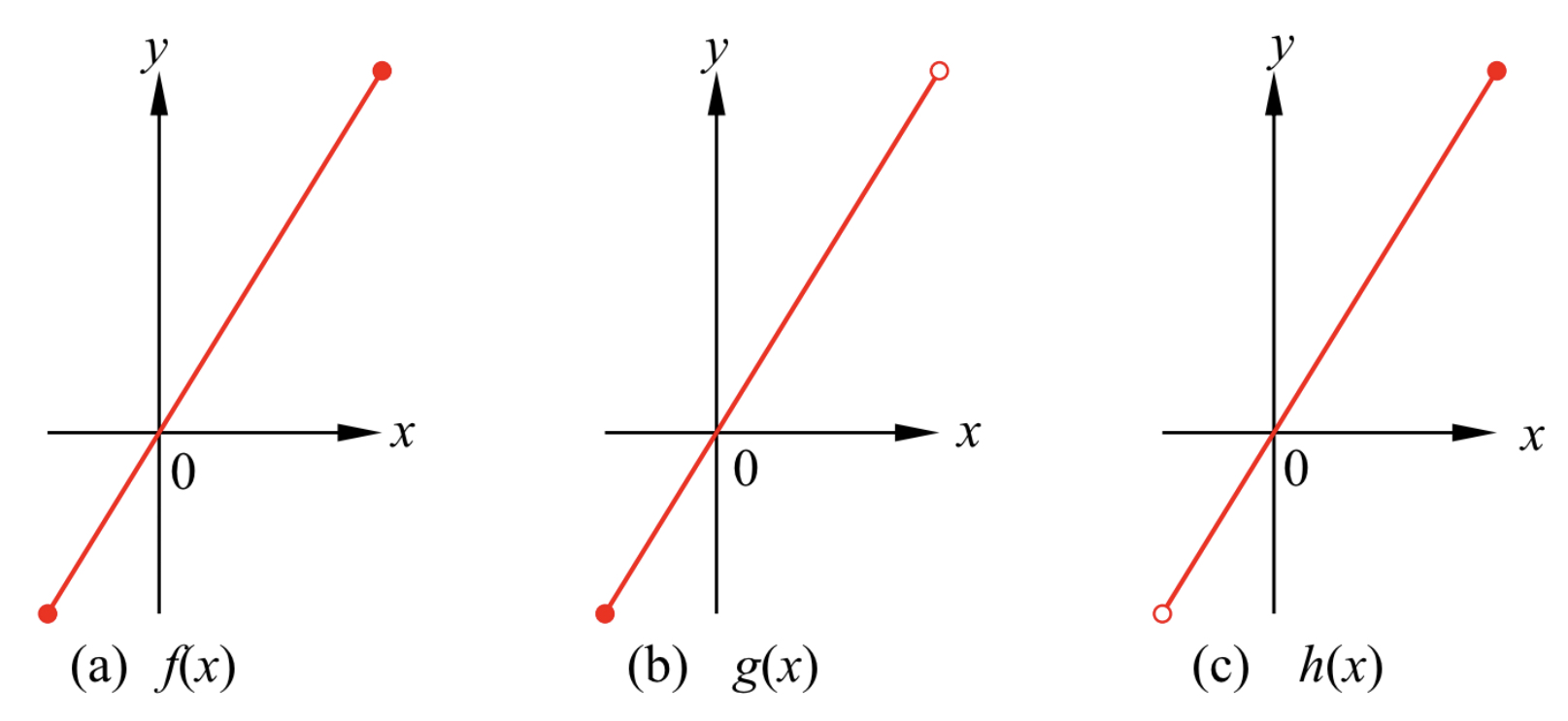}
\caption{  The functions $f(x)$, $g(x)$ and $h(x)$ defined in Example \ref{23020905}.}\label{figure11}
\end{figure}

Example \ref{23020905} shows that the existence of extreme values  depends on the domain of the function.

For a set to have maximum and minimum values, it is necessary (but not sufficient) that the set is bounded. Let us first define what it means for a function to be bounded.

\begin{definition}{Bounded Functions}
We say that a real-valued function $f:D\rightarrow\mathbb{R}$ is bounded if its range $f(D)$ is bounded. In other words, a function $f:D\rightarrow\mathbb{R}$ is bounded  if and only if there is a positive constant $M$ such that 
\[|f(x)|\leq M\hspace{1cm}\text{for all}\;x\in D.\]

\end{definition}

\begin{example}{}
All the three functions defined in Example \ref{23020905} are  bounded.
\end{example}

We are interested in a sufficient condition for a continuous function to have maximum and minimum values. Before we proceed, let us look at two examples.

  \begin{example}{}
 Consider the function $f:(0,1)\rightarrow\mathbb{R}$ defined by $\di f(x)=\frac{1}{x}$. Although the domain of the function $D=(0,1)$ is bounded, the range of the function $f(D)=(1, \infty)$ is not bounded.
 \end{example}
 This example shows that a continuous function does not necessarily map a bounded set to a bounded set.

  \begin{example}{}
 Consider the function $f:[1,\infty)\rightarrow\mathbb{R}$ defined by $\di f(x)=\frac{1}{x}$. Although the domain of the function $D=[1, \infty)$ is closed, the range of the function $f(D)=(0, 1]$ is not closed.
 \end{example}
 This example shows that a continuous function does not necessarily map a closed set to a closed set. 
  
   The situation changes   when we combine closed and bounded.
 Recall that we have defined the concept of sequential compactness in Chapter \ref{ch1}, Section \ref{sec1.8}. A set $D$ is sequentially compact if    every sequence in $D$ has a subsequence that converges to a point in $D$. We have proved that a subset of real numbers is sequentially compact if and only if it is closed and bounded.

  The following theorem says that a continuous function maps a  closed and bounded set to a closed and bounded set.
 
  \begin{theorem}[label=23020906]{}
  Let $D$ be a closed and bounded subset of $\mathbb{R}$. If $f:D\rightarrow \mathbb{R}$ is  a continuous function, then the set $f(D)$ is   closed and bounded.
  \end{theorem}
  Using the fact that a subset of real numbers is sequentially compact if and only if it is closed and bounded, Theorem \ref{23020906} is equivalent to the following.
  \begin{theorem}
  [label=23020607]{}
  Let $D$ be a sequentially compact subset of $\mathbb{R}$. If $f:D\rightarrow \mathbb{R}$ is  a continuous function, then the set $f(D)$ is sequentially compact.
  \end{theorem}
  \begin{myproof}{Proof}
  We use the definition of sequential compactness to prove this theorem.  Let $\{y_n\}$ be a sequence in $f(D)$. We need to prove that there is a subsequence of $\{y_n\}$ that converges to a point in $f(D)$. 
  
  For each positive integer $n$, since $y_n$ is in $f(D)$, there is an $x_n$ in $D$ such that $f(x_n)=y_n$. This gives a sequence $\{x_n\}$ in $D$. Since $D$ is sequentially compact, there is a subsequence $\{x_{n_k}\}$ of $\{x_n\}$ that converges to a point $x_0$ in $D$. Since $f$ is continuous at $x_0$, the sequence $\{f(x_{n_k})\}$ converges to $f(x_0)$. In other words, we have shown that the subsequence $\{y_{n_k}\}$ of $\{y_n\}$ converges to the point $f(x_0)$ in $f(D)$.

  \end{myproof}
  
Proving Theorem \ref{23020906}   without using  sequential compactness is tedious, and it essentially goes through some of the arguments used to prove that a subset of real numbers is sequentially compact if and only if it is closed and bounded. From here, we can see the usefulness of the concept of sequential compactness.
  
  In Theorem \ref{23020908}, we have seen that a set that is closed and bounded must have a maximum and a minimum. Hence, we obtain immediately the following   theorem.
  
  \begin{theorem}[label=23020909]{Extreme Value Theorem}
 
  Let $D$ be a closed and bounded subset of $\mathbb{R}$. If $f:D\rightarrow \mathbb{R}$ is  a continuous function, then $f$ has a maximum value and a minimum value. 
    \end{theorem}
    \begin{corollary}{}
  If $f:[a,b]\rightarrow\mathbb{R}$ is a continuous function defined on a closed and bounded interval, then $f$ is bounded, and it has a maximum value and a minimum value. \end{corollary}

  Extreme value theorem is used to guarantee the existence of a maximum value and a minimum value before we proceed to find these values, so that the attempt to look for extreme values is not futile. In some circumstances, knowing the existence of such extreme values is sufficient. 
  
  \begin{example}{}
  Show that the function $f:\mathbb{R}\rightarrow\mathbb{R}$ defined by 
  \[f(x)=|x-1|+|x-2|+|x-3|+|x-4|+|x-5|\] has a minimum value. 
  \end{example}
  \begin{solution}{Solution}
  In this example, the domain of the function is not closed and bounded. We cannot apply the extreme value theorem directly. However, we can proceed in the following way.
  First, we justify that the function  $f:\mathbb{R}\rightarrow\mathbb{R}$  is continuous. A function of the form $g(x)=x-a$ is continuous since it is a polynomial. Absolute value of a continuous function is continuous. Hence, a function of the form $h(x)=|x-a|$ is continuous. Being a sum of continuous functions, $f(x)$ is a continuous function.

 To prove the existence of a minimum value,  we notice that for $x\geq 5$,
  \[f(x)=x-1+x-2+x-3+x-4+x-5=5x-15\geq  10.\]
  For $x\leq 1$,  
  \[f(x)=1-x+2-x+3-x+4-x+5-x=15-5x\geq 10.\]
  Now restrict the domain to $[1,5]$, the function $f:[1,5]\rightarrow\mathbb{R}$ is continuous. Hence, it has a minimum value at some $x_0\in [1,5]$. It follows that
  \[f(x_0)\leq f(x)\hspace{1cm}\text{for all}\;x\in [1,5].\]
  In particular,
  \[f(x_0)\leq f(1)=10.\] 
  This proves that for all $x\in\mathbb{R}$, 
$f(x)\geq f(x_0)$.
  Hence, the function $f:\mathbb{R}\rightarrow\mathbb{R}$ has a minimum value.
  \end{solution}
  
\vp
\noindent
{\bf \large Exercises  \thesection}
\setcounter{myquestion}{1}
 \begin{question}{\themyquestion}
 Determine whether the function is bounded.
 \begin{enumerate}[(a)]
 \item $f:\mathbb{R}\rightarrow\mathbb{R}$, $\di f(x)=\frac{x}{\sqrt{x^2+4}}$.
 \item $f:(0,1)\rightarrow\mathbb{R}$, $f(x)=x+\di\frac{1}{x}$.
 \end{enumerate}
\end{question}
\atc
 \begin{question}{\themyquestion}
 If a function $f:D\rightarrow\mathbb{R}$ is continuous and bounded, does it necessarily have a maximum value and a minimum value? Justfiy your answer.
\end{question}
\atc
 \begin{question}{\themyquestion}
 Let $f:[-4,4]\rightarrow\mathbb{R}$ be the function defined by 
 \[f(x)=\frac{x^2+x+1}{\sqrt{4x^2+9}}.\]Show that it has a maximum value and a minimum value.
\end{question}
\vp

\section{The Intermediate Value Theorem}\label{sec2.4}

In this section, we are going to discuss  the intermediate value theorem, which is an  important theorem for continuous functions. It is essentially a theorem about existence of solutions for equations defined by continuous functions. 

\begin{theorem}{Intermediate Value Theorem}
Given that $f:[a,b]\rightarrow \mathbb{R}$ is a continuous function. For any real number $w$ that is between $f(a)$ and $f(b)$, there exists a point $c$ in $[a, b]$ such that
\[f(c)=w.\]

\end{theorem}

\begin{myproof}{Proof}
The proof is using bisection method, which provides a constructive way to find the point $c$. 

Without loss of generality, assume that $f(a)<w<f(b)$.

We construct two sequences $\{a_n\}$ and $\{b_n\}$ recursively. Define $a_1=a$, $b_1=b$, and let
\[m_1=\frac{a_1+b_1}{2}\] be the midpoint of $a_1$ and $b_1$. The interval $[a,b]=[a_1, b_1]$ is bisected into two  subintervals $[a_1, m_1]$ and $[m_1, b_1]$ by the point $m_1$. 

 We want to define the interval $[a_2, b_2]$ to be one of these, based on the value of $f(m_1)$. 
\begin{enumerate}[$\bullet$\;\;]\item
If $f(m_1)<w$, define $a_2=m_1$ and $b_2=b_1$.
 \item If $f(m_1)\geq w$, define $a_2=a_1$ and $b_2=m_1$.
 
 \end{enumerate}By  definition,
 \[a_1\leq a_2<b_2\leq b_1,\]
 \[f(a_2)<w\leq f(b_2),\] and  the length of the interval $[a_2, b_2]$ is half the length of the interval $[a_1, b_1]$. 
 \bp
 Suppose that we have defined $a_1, \ldots, a_n$, $b_1, \ldots, b_n$, such that
 \[a_1\leq   \ldots \leq a_{n-1}\leq a_n<b_n\leq b_{n-1}\leq \ldots \leq b_1,\]
 \[f(a_k)<w\leq f(b_k)\hspace{1cm}\;\text{for all}\;1\leq k\leq n,\]
 and
 \[b_k-a_k=\frac{b_{k-1}-a_{k-1}}{2} \hspace{1cm} \text{for all}\;2\leq k\leq n.\]
 Let
 \[m_n=\frac{a_n+b_n}{2}\]  be the midpoint of $a_n$ and $b_n$.
 \begin{enumerate}[$\bullet$\;\;]\item
If $f(m_n)<w$, define $a_{n+1}=m_n$ and $b_{n+1}=b_n$.
 \item If $f(m_n)\geq w$, define $a_{n+1}=a_n$ and $b_{n+1}=m_n$.
 
 \end{enumerate}  By   definition,
 \[a_n\leq a_{n+1}<b_{n+1}\leq b_n,\]
 \[f(a_{n+1})<w\leq f(b_{n+1}),\] and  
 \[b_{n+1}-a_{n+1}=\frac{b_n-a_n}{2}.\]
This constructs the sequences $\{a_n\}$ and $\{b_n\}$. Notice that $\{a_n\}$ is an increasing sequence that is bounded above by $b$, while $\{b_n\}$ is a decreasing sequence that is bounded below by $a$.
 
By monotone convergence theorem, the sequence $\{a_n\}$ converges to a number $c_1=\sup\{a_n\}$ and the sequence $\{b_n\}$ converges to a number $c_2=\inf\{b_n\}$. By induction, we find that
\[b_n-a_n=\frac{b-a}{2^{n-1}}.\]Taking $n\rightarrow \infty$ limits, we conclude that
\[c_2-c_1=0.\]\bp
It follows that the number   $c=c_1=c_2$   satisfies
\begin{equation}\label{eq230210_1}a_n\leq c\leq b_n\hspace{1cm}\text{for all}\;n\in\mathbb{Z}^+,\end{equation}and 
\begin{equation}\label{eq230210_2}\lim_{n\rightarrow\infty}a_n=c=\lim_{n\rightarrow\infty}b_n.\end{equation}
Eq. \eqref{eq230210_1} shows that $c$ is in $[a,b]$. The continuity of the function $f$ and \eqref{eq230210_2} implies that
\[f(c)=\lim_{n\rightarrow \infty}f(a_n)=\lim_{n\rightarrow\infty}f(b_n).\]
Since \[f(a_n)<w \quad\text{and}\quad f(b_n)\geq w\hspace{1cm}\text{for all}\; n\in\mathbb{Z}^+,\]
we find that
\[f(c)\leq w\hspace{1cm}\text{and}\hspace{1cm}f(c)\geq w.\]
This proves that $f(c)=w$, and hence completes the proof of the theorem.

\end{myproof}

 \begin{figure}[ht]
\centering
\includegraphics[scale=0.2]{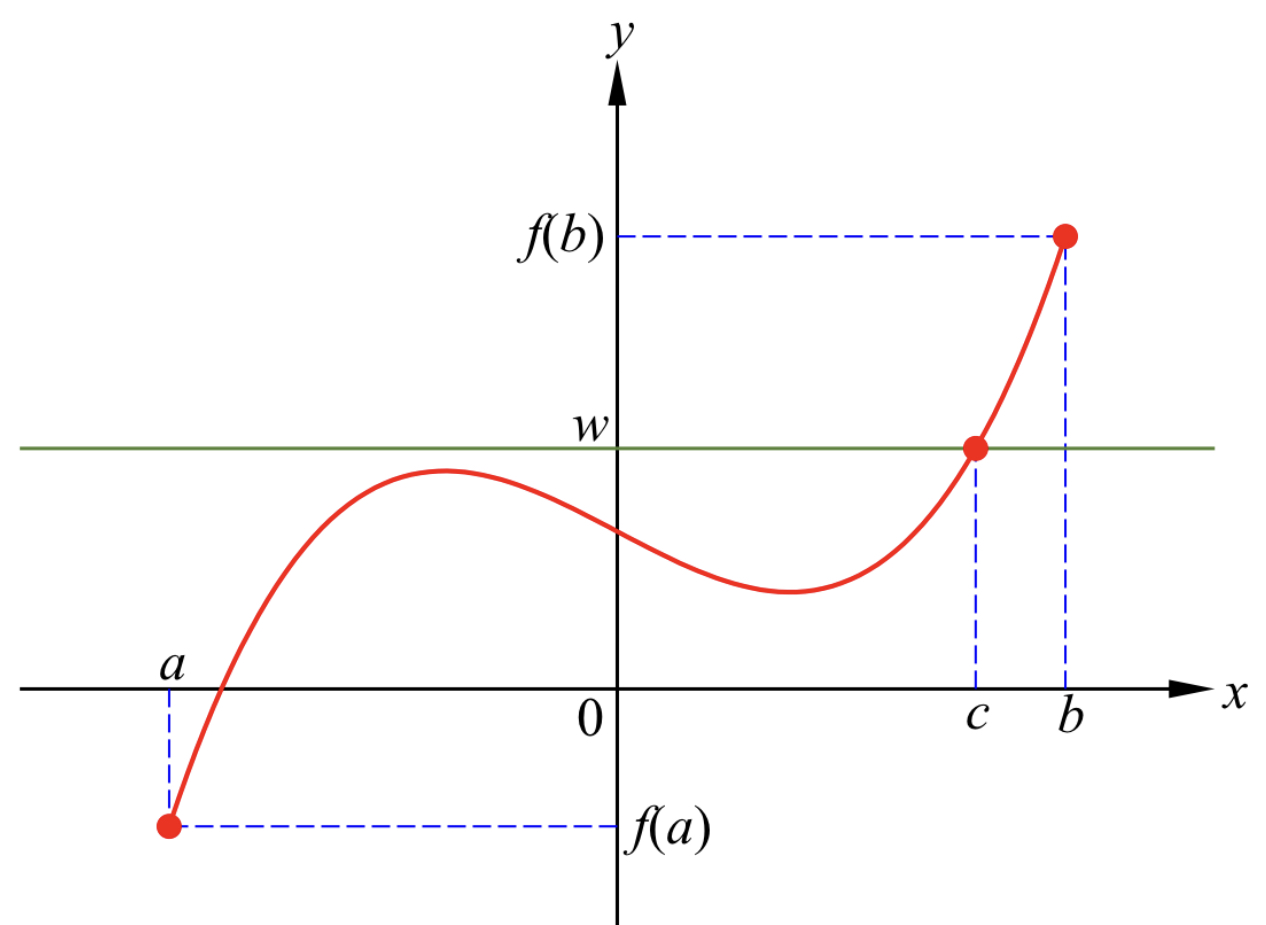}
\caption{The  intermediate value theorem.}\label{figure12}
\end{figure}
The following is an example which we use the intermediate value theorem to justify that an equation has a solution.
\begin{example}[label=ex230215_1]{}
Show that the equation 
\[x^6+6x+1=0\]has a real root.
\end{example}
\begin{solution}{Solution}
Let $f(x)=x^6+6x+1$.
Since $f(x)$ is a polynomial, it is a continuous function. Notice that
\[f(0)=1,\hspace{1cm}f(-1)=-4.\]
Hence, $f(-1)<0<f(0)$. Namely, $0$ is a value between $f(-1)$ and $f(0)$. By intermediate value theorem, there is a point $c$ in the interval $(-1,0)$ such that $f(c)=0$. Then $x=c$ is a root of the equation\[x^6+6x+1=0.\]
\end{solution}In this solution, the choice of $a=-1$ and $b=0$ are by trial and error. In practice, one can use a computer to sample some $x$ values and calculate the corresponding values of $f(x)$. The goal is to find   $a$ and   $b$ such that $f(a)$ and $f(b)$ have oppositive signs. To calculate the root $c$, one can implement the bisection method numerically.

\begin{example}{}
Let $n$ be a positive integer, and let $c$ be a positive number. Use the intermediate value theorem to show that there is a positive real number $x$ such that
\[x^n=c.\]
\end{example}

\begin{solution}{Solution}Take $a=0$, $b=c+1$, and
consider the function $f:[a, b]\rightarrow\mathbb{R}$ defined by $f(x)=x^n$. Then, \[f(a)=f(0)=0, \]
\[  f(b)=f(c+1)=(1+c)^n\geq 1+nc>c.\]
Hence, \[f(a)<c<f(b).\]
Since $f$ is a continuous function, intermediate value theorem asserts that there is a number $x$ in the interval $[0, c+1]$ such that $f(x)=x^n=c$.
\end{solution}
In Chapter \ref{ch1} Example \ref{23021011}, we use completeness axiom to solve this problem when $n=2$ and $c=2$. Here we use the intermediate value theorem to tackle the general problem. The tedious part has been settled in the proof of the intermediate value theorem.

In the following, we want to formulate a precise relation between intervals and the intermediate value theorem. We first introduce a concept called convexity.

\begin{definition}{Convex Sets}
Let $S$ be a subset of real numbers. We say that $S$ is convex if for any $u$ and $v$ in $S$, $(1-t)u+tv$ is in $S$ for all $t\in [0,1]$. \\Equivalently, $S$ is convex provided that whenever $u$ and $v$ are in $S$ and $u<v$, then any $w$ in the interval $[u,v]$ is also in $S$.
\end{definition}
The equivalence of the two definitions is seen by observing  that when $t$ changes from $0$ to $1$, $(1-t)u+tv$ goes through all the points in the interval $[u, v]$. 

Obviously, an interval is a convex set. The converse is also true.
\begin{theorem}{}
Let $S$ be a subset of real numbers. If $S$ is a convex set, then $S$ is an interval. 
\end{theorem}
\begin{myproof}{Sketch of Proof}
 If $S$ is bounded below, let $a=\inf S$. Otherwise, set $a=-\infty$. If $S$ is bounded above, let $b=\sup S$. Otherwise, set $b=\infty$. 

If $c$ is a point in $(a, b)$, then $a<c<b$. In particular, since $c>a$, it is not a lower bound of $S$. Hence, there is a point $u$ in $S$ such that $a\leq u<c$. Since $c<b$, $c$ is not an upper bound of $S$. Hence, there is a point $v$ in $S$ such that $c<v\leq b$. Since $u$ and $v$ are in $S$ and $S$ is convex, all points in the interval $[u,v]$ are in $S$. By construction, $u<c<v$. Hence, $c$ is in $S$. This proves that all the points in $(a, b)$ are in $S$.

 Finally, we just need to consider   whether $S$ contains $a$, and whether it contains $b$. 
\end{myproof}

\begin{highlight}{Convex Sets and Intervals}
Let $S$ be a convex set. If $S$ is bounded below, let $a=\inf S$.  If $S$ is bounded above, let $b=\sup S$.  
\begin{enumerate}[1.]
\item If $S$ is bounded, $S$ does not contain $\inf S$ and $\sup S$, then $S=(a,b)$.
\item If $S$ is bounded, $S$  contains $\inf S$ but does not contain $\sup S$, then $S=[a,b)$.
\item If $S$ is bounded, $S$  contains $\sup S$ but does not contain $\inf S$, then $S=(a,b]$.
\item If $S$ is bounded, and $S$  contains both $\inf S$ and $\sup S$, then $S=[a,b]$.
\item If $S$ is bounded below but not bounded above, and $S$ does not contain $\inf S$, then $S=(a,\infty)$.
\item If $S$ is bounded below but not bounded above, and $S$   contains $\inf S$, then $S=[a,\infty)$.
\end{enumerate}\begin{enumerate}[7.]
\item If $S$ is bounded above but not bounded below, and $S$ does not contain $\sup S$, then $S=(-\infty, b)$.
\end{enumerate}\end{highlight}\begin{highlight}{}\begin{enumerate}[8.]
\item If $S$ is bounded above but not bounded below, and $S$   contains $\sup S$, then $S=(-\infty, b]$.
\end{enumerate}\begin{enumerate}[9.]
\item If $S$ is not bounded above nor bounded below, then $S=(-\infty, \infty)=\mathbb{R}$.
\end{enumerate}
\end{highlight}

The following is a   reformulation of the intermediate value theorem.
\begin{theorem}{Intermediate Value Theorem}
Let $I$ be an interval. If the function $f:I\rightarrow\mathbb{R}$ is continuous, then $f(I)$ is an interval.
\end{theorem}
\begin{myproof}{Proof}
To show that $f(I)$ is an interval, take two distinct points $u$ and $v$ in $f(I)$. We need to show that any $w$ in between $u$ and $v$ is in $f(I)$. Since $u$ and $v$ are in $f(I)$, there exist $a$ and $b$ in $I$ such that $u=f(a)$ and $v=f(b)$. Without loss of generality, assume that $a<b$. Since $I$ is an interval, it contains the interval $[a,b]$. Since $f$ is continuous on $[a,b]$, and $w$ is in between $f(a)$ and $f(b)$,    the version of the intermediate value theorem that we have proved implies that there is a point $c$ in the interval $(a, b)$ such that $f(c)=w$. This shows that $w$ is also in $f(I)$.

\end{myproof}

\vp
\noindent
{\bf \large Exercises  \thesection}
\setcounter{myquestion}{1}
 \begin{question}{\themyquestion}
 Show that the equation 
$\di 2x+\sqrt{x^2+1}=0$ has a real solution.
\end{question}

\atc
\begin{question}{\themyquestion}
Given that $f:[-2, 10]\to [-2, 10]$ is a continuous function.
Show that there is a point $x$ in the interval $[-2, 10]$ such that $f(x)=x$.
\end{question}

\atc
\begin{question}{\themyquestion}
Suppose that $f:\mathbb{R}\rightarrow \mathbb{R}$ is a bounded continuous function.  Show that there is a real number $x$  such that $f(x)=x$.
\end{question}
\atc
\begin{question}{\themyquestion}
Let $n$ be an odd positive integer, and let \[p(x)=a_nx^n+a_{n-1}x^{n-1}+\cdots+a_1x+a_0\] be a polynomial of degree $n$. 
Show that $p(x)=0$ has a real root.
\end{question}
\vp

\section{Uniform Continuity}\label{sec2.5}
In Section \ref{sec2.2}, we have defined the concept of continuity at a point. This is a local property which only depends on the function value in a neighbourhood of a point. In this section, we want to define a concept called uniform continuity, which depends on the behaviour of the function on the whole domain. Such a property is called a global property.

\begin{definition}{Uniform Continuity}
Let $D$ be a subset of real numbers. A function  $f:D\rightarrow\mathbb{R}$  defined on $D$ is uniformly continuous provided that for any $\varepsilon>0$, there exists $\delta>0$ such that if $x_1$ and $x_2$ are in $D$ and $|x_1-x_2|<\delta$, then 
\[|f(x_1)-f(x_2)|<\varepsilon.\]
 
\end{definition}

\begin{theorem}[label=230404_1]{Equivalent Definition of Uniform Continuity}
Let $D$ be a subset of real numbers. A function  $f:D\rightarrow\mathbb{R}$  defined on $D$ is uniformly continuous if and only if whenever $\{u_n\}$ and $\{v_n\}$ are sequences in $D$ such that
\[\lim_{n\rightarrow\infty}(u_n-v_n)=0,\]
$\{f(u_n)\}$ and $\{f(v_n)\}$ are sequences in $f(D)$ such that
\[\lim_{n\rightarrow\infty}\Bigl(f(u_n)-f(v_n)\Bigr)=0.\] 
\end{theorem}
Notice that   we only require the sequence $\{u_n-v_n\}$ to converge to 0. We do not require the sequence $\{u_n\}$ nor the sequence $\{v_n\}$ to be convergent.

Theorem \ref{230404_1} can be proved in the same way as we prove the equivalence of two  definitions for limits of functions.

The following is quite obvious.
 \begin{theorem}[label=t23021006]{}
Let $D$ be a subset of real numbers. If $f:D\rightarrow\mathbb{R}$ is a  uniformly continuous function, it is continuous.
\end{theorem}

  Let us compare the definitions of continuity and uniform continuity using the definitions in terms of $\varepsilon-\delta$.

\begin{highlight}{Continuity versus Uniform Continuity}
\begin{enumerate}[$\bullet$]\item
A function $f:D\rightarrow \mathbb{R}$ is continuous if
\[\forall x_0\in D,\,\forall \varepsilon>0, \,\exists \delta>0,\, \forall x \in D, \, |x-x_0|<\delta\implies |f(x)-f(x_0)|<\varepsilon.\]
\item A function $f:D\rightarrow \mathbb{R}$  is uniformly continuous if
\[\forall \varepsilon>0, \,\exists \delta>0,\, \forall x_0\in D,\,\forall x \in D, \, |x-x_0|<\delta\implies |f(x)-f(x_0)|<\varepsilon.\]\end{enumerate}
\end{highlight}
The difference is in the order of the quantifiers. For a function to be continuous, it must be continuous at each point in the domain. For each point  $x_0$ in the domain, there should exist a positive $\delta$ for each positive $\varepsilon$. This number $\delta$ not only depends on $\varepsilon$, but also on the point $x_0$. For uniform continuity, one needs to be able to find a $\delta$ which only depends on $\varepsilon$ but not on the point in the domain. This is where the uniformity lies.

Let us look at some examples of functions that are continuous but not uniformly continuous.
\begin{example}{}
Show that the function $f:(0,1)\rightarrow\mathbb{R}$, $f(x)=\di\frac{1}{x}$ is not uniformly continuous. 
\end{example}
\begin{solution}{Solution}
For a positive integer $n$, let
\[u_n=\frac{1}{n+1},\hspace{1cm}v_n=\frac{1}{n+2}.\]
Then $\{u_n\}$ and $\{v_n\}$ are sequences in the domain $D=(0,1)$, and
\[\lim_{n\rightarrow\infty}(u_n-v_n)=\lim_{n\to\infty}\frac{1}{n+1}-\lim_{n\to\infty}\frac{1}{n+2}=0.\]\bs
Since $f(u_n)=n+1$ and $f(v_n)=n+2$, we find that
\[\lim_{n\rightarrow\infty}\Bigl(f(u_n)-f(v_n)\Bigr)=-1\neq 0.\]Hence, $f$ is not uniformly continuous.
\end{solution}

\begin{example}{}
Show that the function $f:(0,\infty)\rightarrow\mathbb{R}$, $f(x)=x^2$ is not uniformly continuous. 
\end{example}
\begin{solution}{Solution}
For a positive integer $n$, let
\[u_n=n+\frac{1}{n},\hspace{1cm}v_n=n.\]
Then $\{u_n\}$ and $\{v_n\}$ are sequences in the domain $D=(0,\infty)$, and
\[\lim_{n\rightarrow\infty}(u_n-v_n)=\lim_{n\rightarrow\infty}\frac{1}{n}=0.\] 
Since 
\[f(u_n)-f(v_n)=\left(n+\frac{1}{n}\right)^2-n^2=2+\frac{1}{n^2},\]we find that
\[\lim_{n\rightarrow\infty}\Bigl(f(u_n)-f(v_n)\Bigr)=2\neq 0.\]Hence, $f$ is not uniformly continuous.
\end{solution}

If we change the domain of the function, the conclusion is different.
\begin{example}[label=e23021008]{}
Let $f:[-10, 8]\rightarrow \mathbb{R}$ be the function defined by $f(x)=x^2$. Show that $f$ is uniformly continuous.
\end{example}
\begin{solution}{Solution}In the solution of Example \ref{23021007}, we have shown that
for any $x_1$ and $x_2$ in the domain $D=[-10, 8]$, 
\[|f(x_1)-f(x_2)|\leq 20|x_1-x_2|.\]
Given $\varepsilon>0$, take $\delta=\varepsilon/20$. Then $\delta>0$. If $x_1$ and $x_2$ are in $D$ and $|x_1-x_2|<\delta$, then
\[|f(x_1)-f(x_2)|<20\delta=\varepsilon.\]This proves that $f$ is uniformly continuous.
\end{solution}

Example \ref{e23021008} is a function that is Lipschitz. In fact, 
the proof of Theorem \ref{23021005} can be easily modified to prove that a Lipschitz function is uniformly continuous.
 \begin{theorem}[label=23021006]{}
Let $D$ be a subset of real numbers. If $f:D\rightarrow\mathbb{R}$ is a Lipschitz function, then it is uniformly continuous.
\end{theorem}

The converse is not true. For example, the function $f:[0,1]\rightarrow\mathbb{R}$, $f(x)=\sqrt{x}$ is not Lipschitz, but it is uniformly continuous. We leave this to the exercise.

In the following, we give a sufficient condition for a function to be uniformly continuous. 
\begin{theorem}[label=23021008]{}
Let $D$ be a closed and bounded subset of real numbers. If $f:D\rightarrow\mathbb{R}$ is a continuous function, then it is uniformly continuous.
\end{theorem}
To prove this theorem, we start with a technical lemma.
\begin{lemma}[label=23021009]{}
Let $S$ be a sequentially compact set in $\mathbb{R}$, and let $\{a_n\}$ and $\{b_n\}$ be two sequences in $S$. There is strictly increasing sequence of positive integers $\{n_1, n_2, n_3,  \ldots\}$ such that each of 
 the subsequences  $\{a_{n_1}, a_{n_2}, a_{n_3}, \ldots\}$     and 
  $\{b_{n_1}, b_{n_2}, b_{n_3}, \ldots\}$    converges to a point in $S$.
\end{lemma}
Using sequential compactness, we can guarantee that $\{a_n\}$ has a subsequence  that converges to a point in $S$, and $\{b_n\}$ also has a subsequence that converges to a point in $S$. However, the indices of these two subsequences might not be related. We need to choose the subsequences carefully to make sure that the indices $\{n_k\}$ are the same. 
\begin{myproof}{Proof}
First, the sequentially compactness of $S$ guarantees that there is a strictly increasing sequence of positive integers  $\{k_1, k_2, k_3, \ldots\}$  so that the subsequence $\{a_{k_1}, a_{k_2}, a_{k_3}, \ldots\}$  converges to a point $a$ in $S$. For a positive integer $j$, let   
\[c_j=b_{k_j}.\] Consider the sequence $\{c_j\}$ indexed by $j\in\mathbb{Z}^+$. It is a subsequence of $\{b_n\}$, and it is also a sequence in $S$.  Since $S$ is sequentially compact, there is a strictly increasing sequence of positive  integers $\{j_1, j_2, j_3, \ldots\}$  such that the subsequence $\{c_{j_1}, c_{j_2}, c_{j_3}, \ldots\}=\{b_{k_{j_1}}, b_{k_{j_2}}, b_{k_{j_3}}, \ldots\}$ converges to a point $b$ in $S$. For a positive integer $m$, let
$\di n_m=k_{j_m}.$  Then $\{n_1, n_2, n_3, \ldots\}$ is a strictly increasing sequence of positive integers. The sequence $\{a_{n_1}, a_{n_2}, a_{n_3}, \ldots\}$ is a subsequence of $\{a_{k_1}, a_{k_2}, a_{k_3}, \ldots\}$. Hence, it converges to $a$. The sequence $\{b_{n_1}, b_{n_2}, b_{n_3}, \ldots\}$ is the sequence  $\{c_{j_1}, c_{j_2}, c_{j_3}, \ldots\}$ which converges to $b$.
\end{myproof}
Now we return to the proof of Theorem \ref{23021008}.
\begin{myproof}{\linkt Proof of Theorem \ref{23021008}}
We use proof by contradiction. If $f$ is not uniformly continuous,
there is an $\varepsilon>0$ such that for all $\delta>0$, there are points $u$ and $v$ is $D$ with $|u-v|<\delta$ and $|f(u)-f(v)|\geq\varepsilon$.
This implies that for each $n\in\mathbb{Z}^+$, there are points $u_n$ and $v_n$ in $D$ with 
\[\left|u_n-v_n\right|<\frac{1}{n}\quad \text{and}\quad\left|f(u_n)-f(v_n)\right|\geq\varepsilon.\]
Thus, $\{u_n\}$ and $\{v_n\}$ are two sequences in $D$ such that 
\begin{equation}\label{eq230210_5}\lim_{n\rightarrow\infty}(u_n-v_n)=0.\end{equation}
Since $D$ is closed and bounded, it is sequentially compact. By Lemma \ref{23021009}, we find that there is a strictly increasing sequence of positive integers $\{n_1, n_2, n_3,\ldots\}$ so that the subsequence $\{u_{n_k}\}$ converges to a point $u_0$ in $D$; and the subsequence $\{v_{n_k}\}$ converges to a point $v_0$ in $D$.  
 Since $f:D\to\mathbb{R}$ is continuous, the sequence $\{f(u_{n_k})\}$ converges to $f(u_0)$, and the sequence $\{f(v_{n_k})\}$  converges to $f(v_0)$.
By construction,
 \[|f(u_{n_k})-f(v_{n_k})|\geq \varepsilon\hspace{1cm}\text{for all}\;k\in\mathbb{Z}^+.\]
 This implies that
 \[|f(u_0)-f(v_0)|\geq \varepsilon.\]Since $\{u_{n_k}-v_{n_k}\}$ is a subsequence of $\{u_n-v_n\}$, \eqref{eq230210_5} implies that
\[u_0-v_0=\lim_{k\to\infty}\left(u_{n_k}-v_{n_k}\right)=0.\]
 In other words, $u_0=v_0$. Then we should have $f(u_0)=f(v_0)$, which contradicts to $|f(u_0)-f(v_0)|\geq\varepsilon$. 
We conclude that $f$ must be uniformly continuous.
\end{myproof}

\begin{example}{}
Show that the function $f:(0, 100)\rightarrow\mathbb{R}$, $f(x)=\sqrt{x}$ is uniformly continuous.
\end{example}
Using the definition of uniform continuity to solve this problem is tedious.
\begin{solution}{Solution}
The domain of the function $D_f=(0,100)$ is not closed and bounded. We cannot apply Theorem \ref{23021008} directly. Consider the function $g:[0,100]\rightarrow \mathbb{R}$ defined by $g(x)=\sqrt{x}$. It is a continuous function.
Since the domain $D_g=[0, 100]$ is closed and bounded, $g$ is uniformly continuous.

Since $f:(0, 100)\rightarrow\mathbb{R}$ is the restriction of the function $g$ to $D_f$, it is also uniformly continuous.
\end{solution}
\vspace{0.5cm}
\hrule
\vspace{0.5cm}
\noindent
{\bf \large Exercises  \thesection}
\setcounter{myquestion}{1}

\begin{question}{\themyquestion}
Determine whether the function $f:[0, \infty)\rightarrow \mathbb{R}$, $f(x)=2x^2+3x$ is uniformly continuous.
\end{question}
\atc
\begin{question}{\themyquestion}
Show that the function $f:(0, 20)\rightarrow \mathbb{R}$, $f(x)=\di\frac{x}{\sqrt{x+1}}$ is uniformly continuous.
\end{question}
\atc
 \begin{question}{\themyquestion}
Let $f:[0,1]\rightarrow\mathbb{R}$ be the function defined by $f(x)=\sqrt{x}$.
Show  that $f$ is not Lipschitz, but it is uniformly continuous.
\end{question}
\atc
\begin{question}{\themyquestion}
Determine whether the function $f:(0, 1)\rightarrow \mathbb{R}$, $f(x)=\di\frac{1}{\sqrt{x}}$ is uniformly continuous.
\end{question}

\vp

\section{Monotonic Functions and Inverses of  Functions}\label{sec2.6}

In this section, we study monotonic functions.
\begin{definition}{Monotonic    and Strictly Monotonic Functions}
Let $D$ be a subset of real functions and let $f:D\rightarrow\mathbb{R}$ be a function defined on $D$.
\begin{enumerate}[1.]
\item
We say that $f:D\to\mathbb{R}$ is an increasing function if for any $x_1$ and $x_2$ in $D$,
\[x_1\leq x_2\;\implies\;f(x_1)\leq f(x_2).\]
 
\item
We say that $f:D\to\mathbb{R}$ is a strictly increasing function if for any $x_1$ and $x_2$ in $D$,
\[x_1<x_2\;\implies\;f(x_1)< f(x_2).\]
 
\item
We say that $f:D\to\mathbb{R}$ is a decreasing function if for any $x_1$ and $x_2$ in $D$,
\[x_1\leq x_2\;\implies\;f(x_1)\geq f(x_2).\]
 \item
We say that $f:D\to\mathbb{R}$ is a strictly decreasing function if for any $x_1$ and $x_2$ in $D$,
\[x_1<x_2\;\implies\;f(x_1)> f(x_2).\]
\item We say that $f:D\to\mathbb{R}$ is monotonic if it is increasing or it is decreasing.

\item We say that $f:D\to\mathbb{R}$ is strictly monotonic if it is strictly increasing or it is strictly decreasing.
\end{enumerate}
\end{definition}
The following is obvious from the definitions.
\begin{proposition}{}
Let $D$ be a subset of real functions and let $f:D\rightarrow\mathbb{R}$ be a function defined on $D$. If $f$ is strictly monotonic, then it is one-to-one.
\end{proposition}
\begin{example}[label=23021101]{}\begin{enumerate}[(a)]
\item
Let $f:[-3,4]\rightarrow\mathbb{R}$  be the function defined by 
\[f(x)=\begin{cases} -1,\quad &\text{if}\;-3\leq x\leq 0,\\
x-1,\quad &\text{if}\hspace{0.5cm}0<x\leq 4.\end{cases}\]It is an increasing function.
\item Let   $g:[-3,4]\rightarrow\mathbb{R}$ be the function defined by 
\[g(x)=\begin{cases} 1,\quad &\text{if}\;-3\leq x\leq 0,\\
1-x,\quad &\text{if}\hspace{0.5cm}0< x\leq 4.\end{cases}\]It is a decreasing function. 
 \end{enumerate} Neither $f$ nor $g$ is strictly monotonic.
\end{example}

\begin{figure}[ht]
\centering
\includegraphics[scale=0.2]{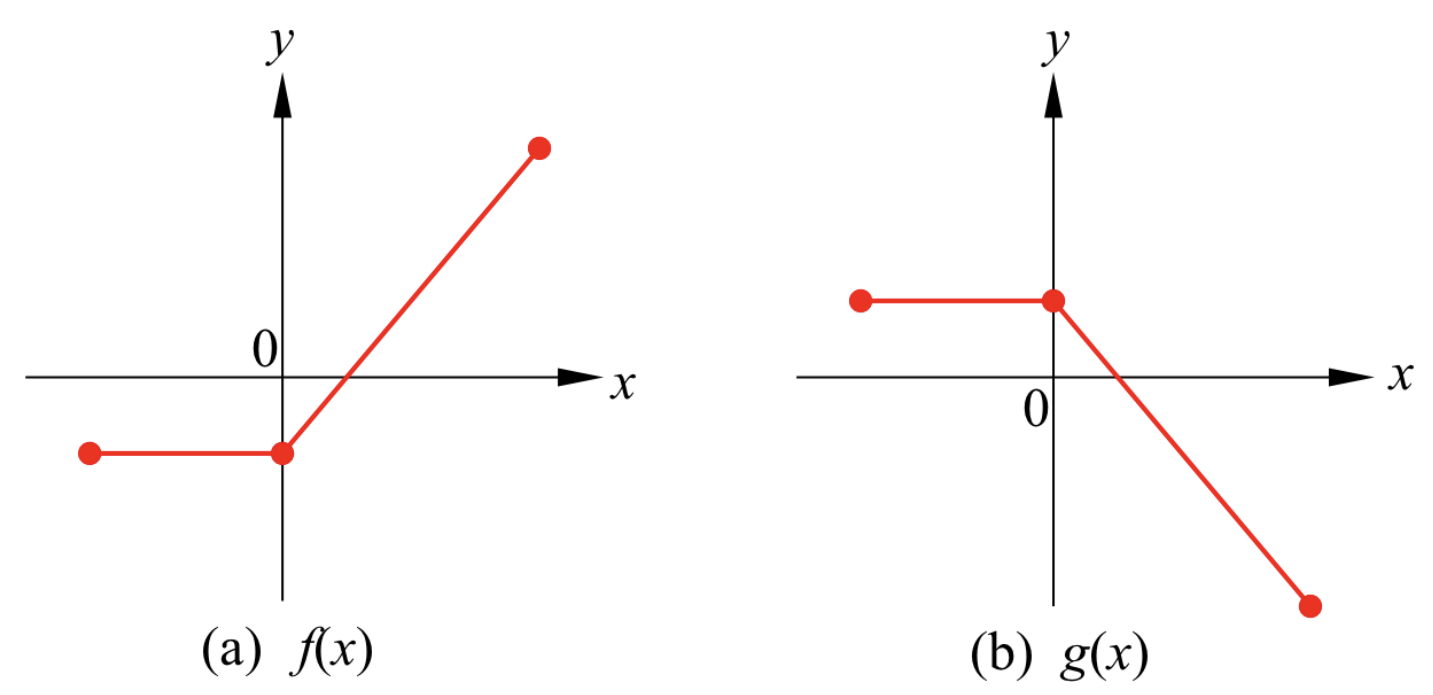}
\caption{  The functions $f(x)$ and $g(x)$ defined in Example \ref{23021101}.}\label{figure13}
\end{figure}

\begin{example}[label=23021102]{}
\begin{enumerate}[(a)]
\item Let $f:(-\infty, 0]\rightarrow\mathbb{R}$ be the function defined by $f(x)=x^2$. Then $f$ is strictly decreasing.

\item Let $g:[0, \infty)\rightarrow\mathbb{R}$ be the function defined by $g(x)=x^2$. Then $g$ is strictly increasing.\end{enumerate}
\end{example}

\begin{figure}[ht]
\centering
\includegraphics[scale=0.2]{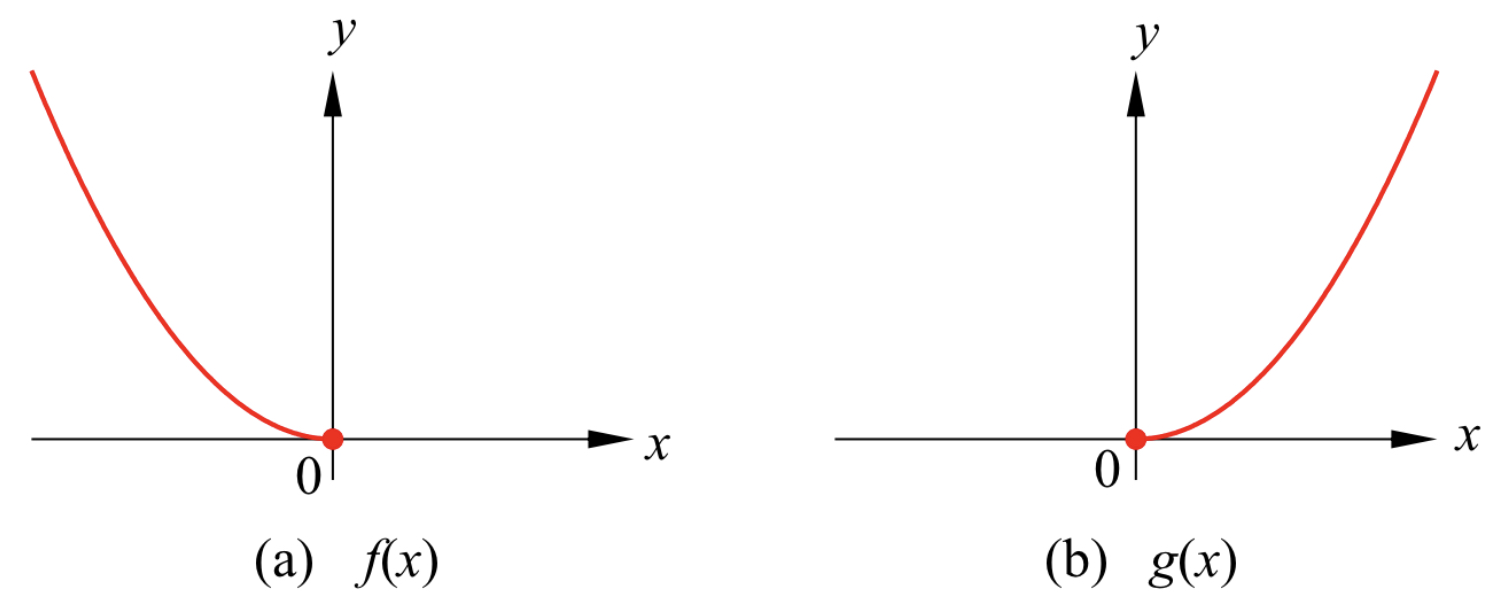}
\caption{  The functions $f(x)$ and $g(x)$ defined in Example \ref{23021102}.}\label{figure14}
\end{figure}

\begin{figure}[ht]
\centering
\includegraphics[scale=0.2]{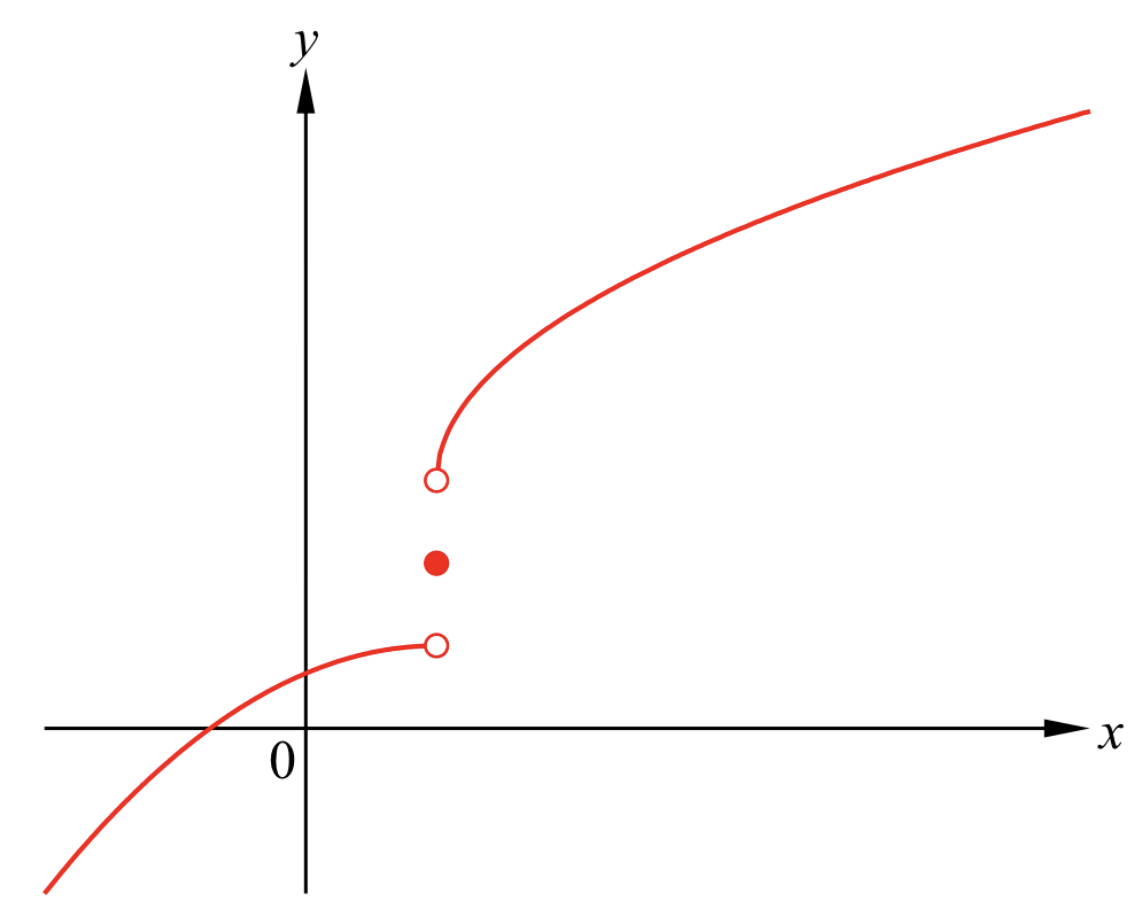}
\caption{  An increasing function with a jump discontinuity.}\label{figure15}
\end{figure}

The following is a characterization of the discontinuities of a  monotonic function.
\begin{theorem}[label=23021103]{}
Let $f:[a,b]\rightarrow\mathbb{R}$ be a  monotonic function. For any $x_0$ in $(a, b]$, 
the left limit \[f_-(x_0)=\di \lim_{x\rightarrow x_0^-}f(x)\] exists. For any $x_0$ in $[a, b)$,    the right limit \[f_+(x_0)=\di\lim_{x\rightarrow x_0^+}f(x)\] exists. 
Define
\[f_-(a)=f(a)\quad\text{and}\quad f_+(b)=f(b).\] Then the  function $f:[a,b]\rightarrow\mathbb{R}$ is continuous at the point $x_0$ in $[a,b]$ if and only if 
\[f_-(x_0)=f(x_0)=f_+(x_0).\]
Otherwise, $f$ has a jump discontinuity at $x_0$ with jump \[|f_+(x_0)-f_-(x_0)|.\]
\end{theorem}

\begin{myproof}{Proof}
If $f$ is  decreasing, then $-f$ is  increasing. Hence, we only need to consider the case where  $f:[a,b]\rightarrow\mathbb{R}$  is  increasing.
Fixed $x_0$ in $(a, b]$. Define the nonempty set $S_-$ by
\[S_-=\left\{f(x)\,|\, a\leq x<x_0\right\}.\]
Since $f$ is  increasing,  $f(x)\leq f(x_0)$ for any $x$ in $[a, x_0)$. Therefore the set $S_-$ is bounded above by $f(x_0)$. Let $u=\sup S_-$. Then $u\leq f(x_0)$. 
We claim that
\[u=\di \lim_{x\rightarrow x_0^-}f(x)=f_-(x_0).\]
Given $\varepsilon>0$,  $u-\varepsilon <u$ and thus it is not an upper bound of $S_-$. Hence, there is a point $x_1$ in $[a, x_0)$ such that
\[f(x_1)>u- \varepsilon.\]\bp 
Let $\delta=x_0-x_1$. Then $\delta>0$. If $x$ is a point in $[a, x_0)$  such that $|x-x_0|<\delta$, then $x_1<x<x_0$, and thus
\[u-\varepsilon<f(x_1)\leq f(x)\leq u.\]From this, we have
\[|f(x)-u|<\varepsilon.\]This proves that
\[f_-(x_0)= \lim_{x\rightarrow x_0^-}f(x)=u.\]
Using similar argument, we find that for any $x_0$ in $[a, b)$, the right limit $\di\lim_{x\rightarrow x_0^+}f(x)$ exists, and
\[f_+(x_0)=\lim_{x\rightarrow x_0^+}f(x)=\inf \left\{f(x)\,|\, x_0<x\leq b\right\}.\]By definition of continuity, the function $f$ is continuous at $x_0$ if and only if $f_-(x_0)=f_+(x_0)$.
The statement about the jump is obvious.

\end{myproof}
 
\begin{corollary}[label=230222_1]{}
Let $I$ be an interval. If $f:I\rightarrow \mathbb{R}$ is monotonic, then $f$ is continuous if and only if $f(I)$ is an interval.
\end{corollary}

\begin{myproof}{Proof}
If $f:I\rightarrow \mathbb{R}$ is continuous, intermediate value theorem implies that $f(I)$ is an interval.

If $f:I\rightarrow \mathbb{R}$ is not continuous, Theorem \ref{23021103} implies that there is a point $x_0$ in the interval $I$ for which either $f_{-}(x_0)\neq f(x_0)$ or $f_+(x_0)\neq f(x_0)$. In any case, $f(I)$ cannot be an interval.
\end{myproof}

For a function $f:I\rightarrow\mathbb{R}$  defined on an interval $I$, we have seen that if $f$ is strictly monotonic, it is one-to-one. It is true even if the function is not continuous. If we assume that the function is continuous, the converse is also true.
It is a consequence of the intermediate value theorem.
\begin{theorem}[label=23021108]{}
Let $f:I\to \mathbb{R}$ be a function defined on an interval $I$. If $f$ is continuous and one-to-one, then $f$ is strictly monotonic. 
\end{theorem}
\begin{myproof}{Proof}
If $f$  fails to be strictly monotonic, there   exist three points $a, x_0, b$ in $I$ such that $a<x_0<b$ and one of the following holds.
\begin{enumerate}[(i)]
\item  $f(a)<f(b)<f(x_0)$
\item $f(x_0)<f(a)<f(b)$
\item $f(a)>f(b)>f(x_0)$
\item $f(x_0)>f(a)>f(b)$

\end{enumerate}Consider case (i) where $f(a)<f(b)<f(x_0)$. Since $w=f(b)$ is a value between $f(a)$ and $f(x_0)$, intermediate value theorem implies that there is a point $c$ in the interval $(a, x_0)$ for which $f(c)=w$. But then $c\neq b$, but $f(c)=f(b)$. This contradicts to $f$ is one-to-one.

Using the same argument, we will reach a contradiction for the other three cases. 
This proves that $f$ must be strictly monotonic.
\end{myproof}

Now we consider invertibility of functions. We only consider functions that are defined on intervals.
\begin{definition}{Invertibility of a Function}
Let  $I$ be an interval and let $f:I\rightarrow\mathbb{R}$ be a function defined on $I$. The function $f:I\rightarrow\mathbb{R}$ is invertible if and only if it is injective. If $f:I\rightarrow\mathbb{R}$ is injective,   its inverse is the function $f^{-1}:f(I)\rightarrow\mathbb{R}$ defined in such a way so that
\[f^{-1}(y)=x\iff f(x)=y\hspace{1cm}\text{for all}\;y\in f(I).\]
\end{definition}

\begin{example}[label=23021105]{}
Consider the functions $f$ and $g$ that are defined in Example \ref{23021102}. 
\begin{enumerate}[(a)]
\item The inverse of the function
  $f:(-\infty, 0]\rightarrow\mathbb{R}$, $f(x)=x^2$, is the function $f^{-1}:[0, \infty)\to (-\infty, 0]$, $f^{-1}(x)=-\sqrt{x}$.
\item The inverse of the function
$g:[0, \infty)\rightarrow\mathbb{R}$, $g(x)=x^2$, is the function $g^{-1}:[0, \infty)\to [0, \infty)$, $g^{-1}(x)=\sqrt{x}$.\end{enumerate}
\end{example}
\begin{figure}[ht]
\centering
\includegraphics[scale=0.2]{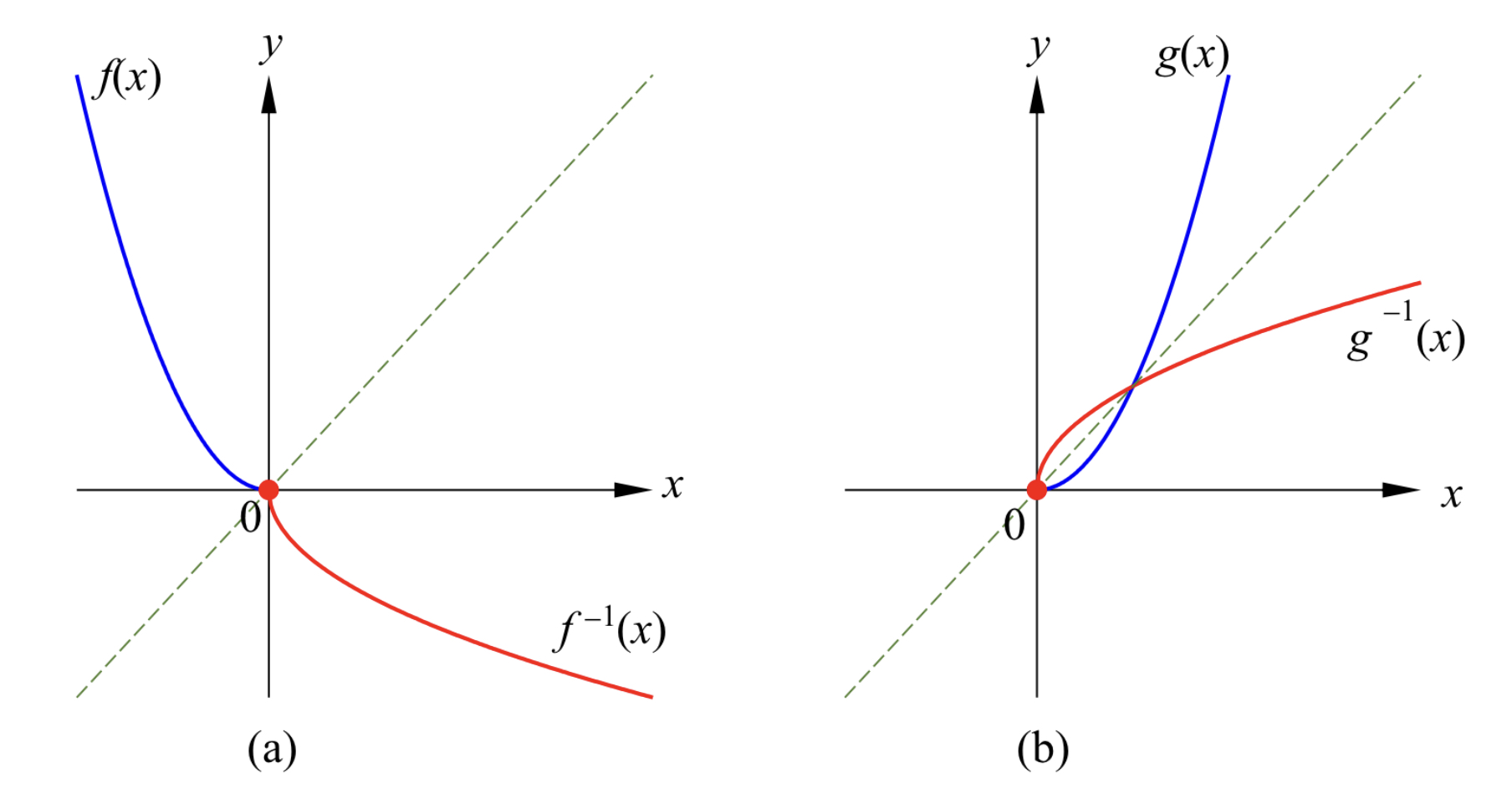}
\caption{  The functions  in Example \ref{23021105}.}\label{figure16}
\end{figure}
Notice that the inverse of a strictly increasing function is strictly increasing. The inverse of a strictly decreasing function is strictly decreasing.

In the next theorem, we prove that the inverse of a continuous function   is continuous. 
\begin{theorem}[label=23021106]{}
Let $I$ be an interval and let $f:I\to\mathbb{R}$ be a continuous function defined on $I$. If $f:I\rightarrow\mathbb{R}$ is one-to-one, then $f^{-1}:f(I)\to\mathbb{R}$ exists, and it is continuous.

\end{theorem}
\begin{myproof}
{Proof}
By Theorem \ref{23021108}, $f$ is strictly monotonic. Without loss of generality, we assume that $f$ is strictly increasing. 

 Given  a point $y_0$ in the interval $f(I)$, let $x_0$ be the unique point in $I$ such that $f(x_0)=y_0$. Given $\varepsilon>0$, we need to prove that there is a $\delta>0$ such that if $y$ is a point in $f(I)$ with $|y-y_0|<\delta$, then $|f^{-1}(y)-f^{-1}(y_0)|<\varepsilon$.

For simplicity, assume that $x_0$ is an interior point of $I$. Then there is a $r>0$ such that $[x_0-r,x_0+r]$ is in $I$. 
Take \[\varepsilon_1=\min\{\varepsilon, r\}.\] Then $\varepsilon_1>0$, $\varepsilon_1\leq \varepsilon$ and $[x_0-\varepsilon_1, x_0+\varepsilon_1]\subset [x_0-r, x_0+r]\subset I$.
Since $f$ is strictly increasing,
\[f(x_0-\varepsilon_1)<f(x_0)<f(x_0+\varepsilon_1),\]and the interval $[f(x_0-\varepsilon_1), f(x_0+\varepsilon_1)]$ is in $f(I)$. 
Let
\[\delta=\min\{f(x_0)-f(x_0-\varepsilon_1), f(x_0+\varepsilon_1)-f(x_0)\}.\]
Then $\delta>0$. If $y$ is a point in $I$ such that $|y-y_0|<\delta$, then \[f(x_0-\varepsilon_1)\leq y_0-\delta<y<y_0+\delta\leq f(x_0+\varepsilon_1).\]This implies that $y$ is also in $f(I)$.  Since $f^{-1}$ is  strictly increasing, we have
\[x_0-\varepsilon\leq x_0-\varepsilon_1<f^{-1}(y)<x_0+\varepsilon_1\leq x_0+\varepsilon.\]This implies that
\[|f^{-1}(y)-f^{-1}(y_0)|<\varepsilon,\]which completes the proof that $f^{-1}$ is continuous at $y_0$.

If $f^{-1}(y_0)=x_0$ is an endpoint of $I$, we need to modify the proof a bit to show that $f^{-1}$ is continuous at $y_0$. The details are left to the students.

\end{myproof}

\begin{remark}{}
If $I=(a, b)$ is an open interval and the function $f:(a, b)\to \mathbb{R}$ is continuous and one-to-one, we have seen in Theorem \ref{23021108} that $f$ is strictly monotonic. In fact, one can prove that $f(I)$ is also an open interval. 

Without loss of generality, assume that $f$ is strictly increasing.   Since $f$ is continuous, $f(I)$ is an interval. If $f(I)$ is not an open interval, either $\inf f(I)$ or $\sup f(I)$ is in $f(I)$. If $c=\inf f(I)$ is in $f(I)$, there is a point $u$ in $(a, b)$ such that $f(u)=c$. But then $u>a$ and so $u_1=\di \frac{u+a}{2}$ is also a point in $(a,b)$. Since $u_1<u$, $f(u_1)<f(u)=c$. This contradicts to $c=\inf f(I)$. In the same way, one can show that $\sup f(I)$ is not in $f(I)$. Hence, $f(I)$ must be an open interval.
\end{remark}

Although we can use limits to show that when $n$ is a positive integer, the  function $f(x)=\sqrt[n]{x}$ is continuous, it is tedious. Using Theorem \ref{23021106} is much more succint.
\begin{example}{}
Let $n$ be a positive integer. 
\begin{enumerate}[1.]
\item If $n$ is odd, the function $f:\mathbb{R}\to\mathbb{R}$, $f(x)=x^n$ is continuous and  one-to-one. Hence, its inverse
$f^{-1}:\mathbb{R}\to\mathbb{R}$, $f^{-1}(x)=\sqrt[n]{x}$ is a continuous function.
\item If $n$ is even, the function $f:[0,\infty)\to[0,\infty)$, $f(x)=x^n$ is continuous and  one-to-one. Hence, its inverse
$f^{-1}:[0,\infty)\to[0,\infty)$, $f^{-1}(x)=\sqrt[n]{x}$ is a continuous function.
\end{enumerate}
\end{example}

Recall that a rational number $r$ can be written as $r=p/q$, where $p$ is an integer and $q$ is a positive integer. For a positive real number $x$, we define
$x^r$ by
\[x^r=\sqrt[q]{x^p}=\left(\sqrt[q]{x}\right)^p.\]
It is easy to check that the two expressions for $x^r$ are equal.
Using the fact that composition of continuous functions is continuous, we obtain the following.

\begin{theorem}{}
Let $r$ be a rational number. 
\begin{enumerate}[1.]
\item If $r>0$, $f:[0,\infty)\to [0,\infty)$, $f(x)=x^r$ is a strictly increasing continuous function.
\item If $r<0$, $f:(0,\infty)\to (0,\infty)$, $f(x)=x^r$ is a strictly decreasing continuous function.
\end{enumerate}
\end{theorem}

\vspace{0.8cm}\hrule\vspace{0.8cm}
\noindent
{\bf \large Exercises  \thesection}
\setcounter{myquestion}{1}
 \begin{question}{\themyquestion}
 Show that the function $f:(-1,\infty)\rightarrow \mathbb{R}$, 
$\di f(x)=\frac{2x+1}{x+1}$ is strictly monotonic, and find  the inverse function $f^{-1}$.
\end{question}

\chapter{Differentiating Functions of a Single Variable}\label{ch3}
The simplest function is the constant function $f(x)=c$, whose function value does not vary with the input. The next  class  of functions that are relatively easy to study is a polynomial of degree one $f(x)=ax+b$, where $a\neq 0$. Sometimes we also call  any function of the form $f(x)=ax+b$ as a linear function, as its graph $y=ax+b$ is a straight line in the $xy$-plane. However, this should not be confused with a linear function that are considered in linear algebra, which in the single variable case, refers to a function of the form $f(x)=ax$. 

\begin{definition}{Graph of a Function}
If $f:D\rightarrow\mathbb{R}$ is a function defined on a subset $D$ of real numbers, its {\bf graph} is the subset $G_f$ in $\mathbb{R}^2$ defined as
\[G_f=\left\{(x,y)\,|\, x\in D, y=f(x)\right\}.\]
\end{definition}

 \begin{figure}[ht]
\centering
\includegraphics[scale=0.25]{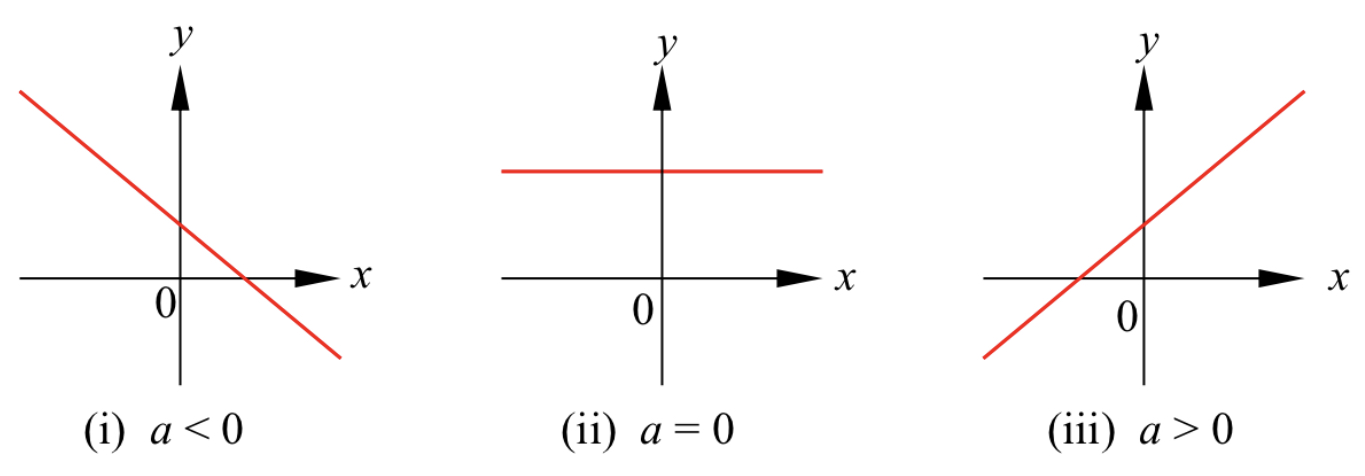}
\caption{  The graph of the function $f(x)=ax+b$ when (i) $a<0$, (ii) $a=0$ and (iii) $a>0$.\fa}\label{figure17}
\end{figure}

  For the function $y=f(x)=ax+b$,  we find that for any two distinct points $x_1$ and $x_2$,
\[\frac{f(x_2)-f(x_1)}{x_2-x_1}=\frac{a(x_2-x_1)}{x_2-x_1}=a.\]
In other words, the change in the $y$ values,
\[\Delta y=y_2-y_1=f(x_2)-f(x_1)\]
is proportional to the change in the $x$-values
\[\Delta x=x_2-x_1,\] with propotionality constant $a$. This constant $a$ is called the rate of change of the function, and it is the slope of the line $y=ax+b$. Its magnitude $|a|$ characterizes how fast $y$ is changing with respect to $x$, and its sign determine the way $y$ changes. When $a>0$,  $y$ increases as   $x$ increases. When $a<0$, $y$ decreases as $x$ increases. 

For a function that is more complicated, such as a quadratic function $y=f(x)=x^2$, we find that
\[\Delta y=f(x_2)-f(x_1)=x_2^2-x_1^2=(x_2+x_1)(x_2-x_1)=(x_2+x_1)\Delta x.\]
In this case, $\Delta y/\Delta x$ is not a constant. It depends on the points $x_1$ and $x_2$. 

In real-life scenario, functions are used to describe the dependence of a variable on the other. For example, if one wants to record the distance that has been travelled by an object,  the independent variable is the time $t$, while the dependent variable is the distance $s$. In this case, one obtains a function $s=s(t)$. The average speed the object is travelling between the time $t=t_1$ and the time $t=t_2$ is
\[\frac{s(t_2)-s(t_1)}{t_2-t_1}.\]In general, one cannot expect that this speed is a constant. If we are interested in the instantaneous speed that the object is travelling at time $t=t_1$, one can fix $t_1$ and  take  $ t_2$ to be closer and closer to $t_1$, and study the behaviour of the average speed. This leads to the concept of derivatives.

\section{Derivatives }\label{sec3.1}

The derivative of a function $y=f(x)$ is a measure of the rate of change of the $y$-values with respect to the change in the $x$- values. To be able to measure this rate of change at a particular point $x_0$, the function has to be defined in a neighbourhood of the point $x_0$. Henceforth, when we define derivatives, we will assume that the function is defined on an open interval $(a, b)$. This includes the case where $a$ is $-\infty$ or $b$ is $\infty$.

\begin{definition}{Derivatives}
  Given a function $f:(a,b)\to \mathbb{R}$ and a point $x_0$ in the interval $(a, b)$, the derivative of $f$ at $x_0$ is defined to be the limit
\[\lim_{x\rightarrow x_0} \frac{f(x)-f(x_0)}{x-x_0}\] if it exists. If the limit exists, we say that $f$ is differentiable at $x_0$, and its derivative a $x_0$ is denoted by $f'(x_0)$. Namely,
\[f'(x_0)=\lim_{x\rightarrow x_0} \frac{f(x)-f(x_0)}{x-x_0}.\]

\end{definition}

Notice that in defining the derivative of a function $f(x)$ at a point $x_0$, the function that we are taking limit of is the function
\[g(x)=\frac{f(x)-f(x_0)}{x-x_0},\]
which is defined on the set $D=(a,b)\setminus\{x_0\}$. It is easy to check that $x_0$ is indeed a limit point of the set  $D$. The function $g(x)$ is the quotient of the function $p(x)=f(x)-f(x_0)$ and the function $q(x)=x-x_0$. It is not defined at $x=x_0$ since $q(x_0)=0$. Moreover, since  $\di\lim_{x\rightarrow x_0}q(x)=q(x_0)=0$,  a necessary condition for $f$ to be differentiable at the point $x_0$ is $\di\lim_{x\rightarrow x_0}p(x)=0$, which says that  the function $f(x)$ is continuous at $x_0$.

\begin{theorem}[label=23021302]{Differentiability Implies Continuity}
Let $x_0$ be a point in the open interval $(a, b)$. If the function $f:(a,b)\to \mathbb{R}$ is differentiable at $x_0$, it is continuous  at $x_0$.
\end{theorem}
\begin{myproof}{Proof}
If $f$ is differentiable at $x_0$, the limit
\[f'(x_0)=\lim_{x\rightarrow x_0} \frac{f(x)-f(x_0)}{x-x_0}\] exists. By limit laws,
\[\lim_{x\to x_0}(f(x)-f(x_0))=\lim_{x\rightarrow x_0} \frac{f(x)-f(x_0)}{x-x_0}\lim_{x\to x_0}(x-x_0)=f'(x_0)\times 0=0.\]
Hence,
\[\lim_{x\to x_0}f(x)=f(x_0),\] which proves that $f$ is continuous at $x_0$.
\end{myproof}

\begin{highlight}{}
By writing $x=x_0+h$, where $h=x-x_0$ is the change in the $x$-values, we can write the derivative of a function $f:(a,b)\to\mathbb{R}$ at a point $x_0$ as
\[f'(x_0)=\lim_{h\to 0}\frac{f(x_0+h)-f(x_0)}{h}.\]
The continuity of the function $f$ at the point $x_0$ is then equivalent to
\[\lim_{h\to 0}f(x_0+h)=f(x_0).\]
\end{highlight}

\begin{definition}{Differentiable Functions}
We say that a function $f:(a,b)\to \mathbb{R}$ is {\bf differentiable} if it is differentiable at all points in $(a,b)$. In this case, the derivative of $f$ is the function $f':(a,b)\to\mathbb{R}$, where
\[f'(x)=\lim_{h\to 0}\frac{f(x+h)-f(x)}{h}.\]
\end{definition}
Let us look at the simplest example where $f(x)=ax+b$.
\begin{example}
{}
Let $f:\mathbb{R}\to\mathbb{R}$ be the function $f(x)=ax+b$. For any $x$ and $x_0$ where $x\neq x_0$, we have
\[\frac{f(x)-f(x_0)}{x-x_0}=a.\]
This implies that
\[f'(x_0)=\lim_{x\to x_0}\frac{f(x)-f(x_0)}{x-x_0}=a.\] Hence, $f$ is a differentiable function and its derivative is
\[f'(x)=a\hspace{1cm}\text{for all}\;x\in\mathbb{R}.\]
\end{example}

Now let us look at a quadratic function.
\begin{example}[label=23021301]{}
Let $f:\mathbb{R}\rightarrow \mathbb{R}$ be the function $f(x)=x^2$. Show that $f$ is differentiable and find its derivative.
\end{example}
\begin{solution}{Solution}
For any real number $x$,
\begin{align*}
f'(x)&=\lim_{h\to 0}\frac{f(x+h)-f(x)}{h}\\
&=\lim_{h\to 0}\frac{(x+h)^2-x^2}{h}\\
&=\lim_{h\to 0}\frac{2xh+h^2}{h}\\
&=\lim_{h\to 0}\;(2x+h)\\
&=2x.
\end{align*}This shows that $f$ is differentiable and its derivative is $f'(x)=2x$.
\end{solution}

 Finding derivative is finding the limit of 
\[m_{x;x_0}=\frac{f(x)-f(x_0)}{x-x_0}\] as $x$ approaches $x_0$. Since $m_{x;x_0}$ is the slope of the secant line joining the two points $(x_0, f(x_0))$ and $(x,f(x))$ on the graph of the function, in the limit $x\to x_0$, we obtain a straight line that only touches the graph in a neighbourhood of the point $(x_0, f(x_0))$ at this point. This line is called the tangent line of the curve $y=f(x)$ at the point $(x_0, f(x_0))$.

 \begin{figure}[ht]
\centering
\includegraphics[scale=0.2]{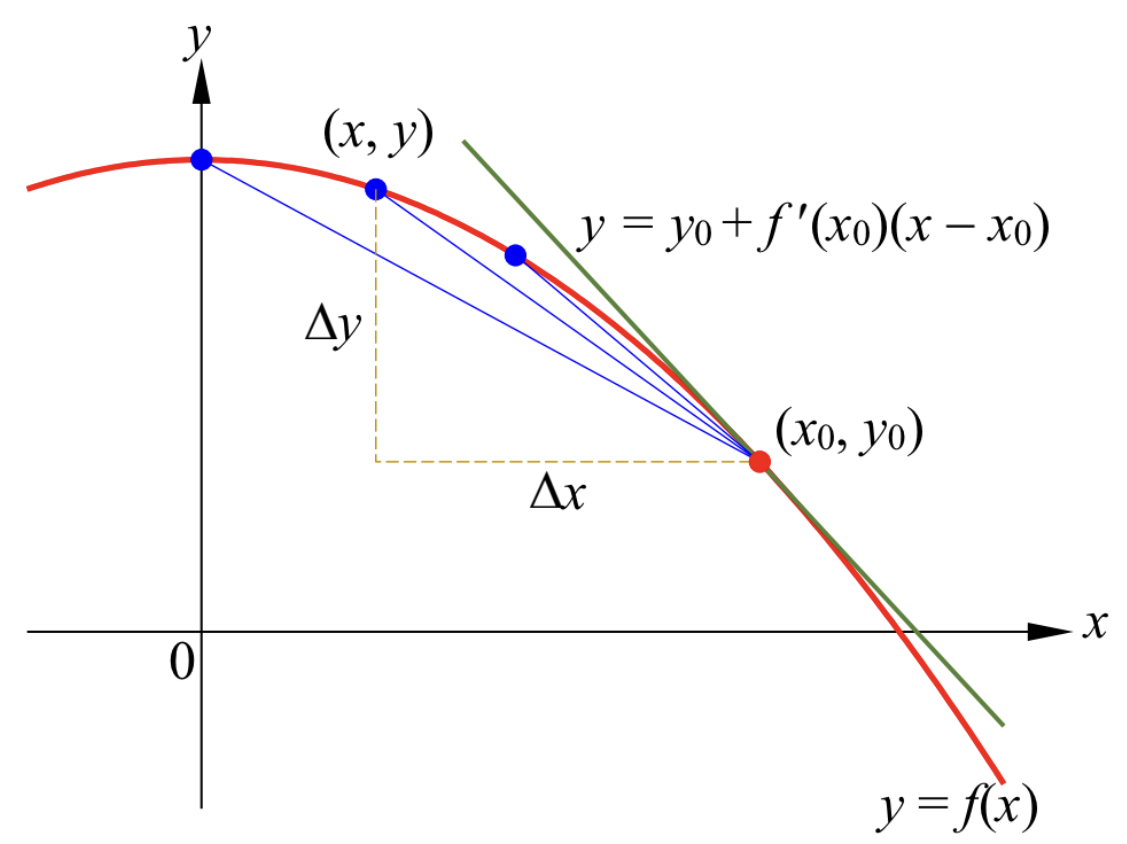}
\caption{  Derivative as slope of tangent line. \fa}\label{figure18}
\end{figure}

\begin{definition}{Tangent Line}
Let $x_0$ be a point in the interval $(a,b)$. If the function $f:(a,b)\to\mathbb{R}$ is  differentiable at $x_0$, then the tangent line to the curve $y=f(x)$ at the point $(x_0, f(x_0))$ is
\[y=y_0+f'(x_0)(x-x_0).\]
\end{definition}

\begin{example}{}
We have found in Example \ref{23021301} that the derivative of the function $f(x)=x^2$ is $f'(x)=2x$. At the point $x=3$, $f(3)=9$ and $f'(3)=6$. Hence, the equation of the tangent line to the curve $y=x^2$ at the point $(3,9)$ is
\[y=9+6(x-3)=6x-9.\]
\end{example}

Theorem \ref{23021302} says that if a function $f:(a,b)\to \mathbb{R}$ is differentiable at $x_0$, then it is continuous at $x_0$. A natural question to ask is whether the converse is true. The answer is no, as shown by the following classical example.

\begin{example}[label=23021303]{}
Consider the function $f:\mathbb{R}\to\mathbb{R}$, $f(x)=|x|$. We have seen in Chapter \ref{ch2} that this is a continuous function. Let $x_0$ be a point in $\mathbb{R}$.

\textbf{Case 1:} If $x_0>0$, then for any $x$ in the neighbourhood $(0, \infty)$ of $x_0$, $f(x)=|x|=x$.  It follows that
\[\lim_{x\to x_0}\frac{f(x)-f(x_0)}{x-x_0}=\lim_{x\to x_0}\frac{x-x_0}{x-x_0}=1.\]
This implies that $f$ is differentiable at $x_0$ and $f'(x_0)=1$.

\textbf{Case 2:} If $x_0<0$, then for any $x$ in the neighbourhood $(-\infty, 0)$ of $x_0$, $f(x)=|x|=-x$.  It follows that 
\[\lim_{x\to x_0}\frac{f(x)-f(x_0)}{x-x_0}=\lim_{x\to x_0}\frac{-x-(-x_0)}{x-x_0}=-\frac{x-x_0}{x-x_0}=-1.\]
This implies that $f$ is differentiable at $x_0$ and $f'(x_0)=-1$.

\textbf{Case 3:} When $x_0=0$, we find that
\[\lim_{x\to 0^+}\frac{f(x)-f(0)}{x-0}=\lim_{x\to 0^+}\frac{x}{x}=1,\]
\[\lim_{x\to 0^-}\frac{f(x)-f(0)}{x-0}=\lim_{x\to 0^-}\frac{-x}{x}=-1.\]
This implies that the limit 
\[\lim_{x\to 0}\frac{f(x)-f(0)}{x-0}\] does not exist. Hence, $f$ is not differentiable at $x=0$.
\end{example}

Graphically, we find that the curve $y=|x|$ has a "sharp turn" at the point $(0,0)$, and there is no well-defined tangent there (see Figure \ref{figure19}).

 \begin{figure}[ht]
\centering
\includegraphics[scale=0.18]{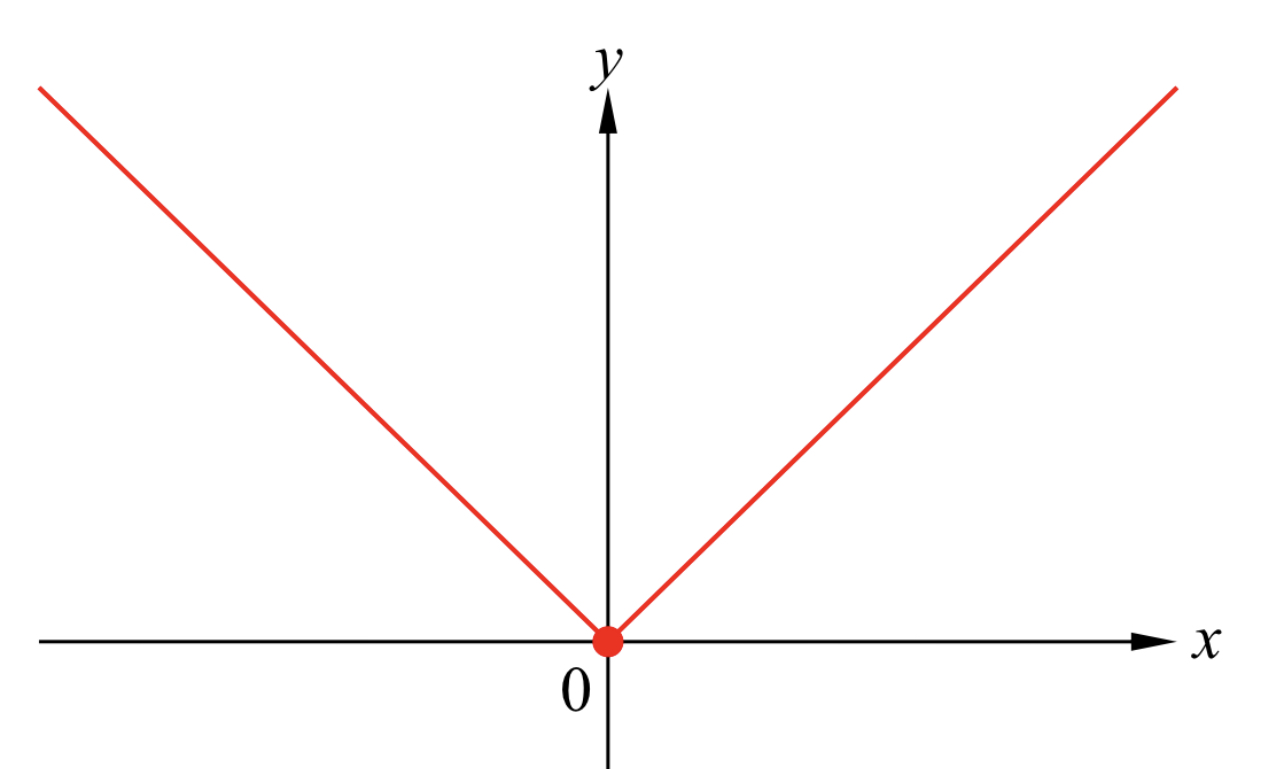}
\caption{  The graph of the function $f(x)=|x|$ has a "sharp turn" at $x=0$. \fa}\label{figure19}
\end{figure}

\begin{remark}{Left Derivatives and Right Derivatives}

We can use left limits and right limits to define left derivatives and right derivatives.  
Let $f:[a,b]\rightarrow \mathbb{R}$ be a function defined on the closed interval $[a,b]$.

\begin{enumerate}[1.]
\item 
For any $x_0\in (a, b]$, we say that the function is left-differentiable at $x_0$ provided that the left derivative at $x_0$, $f_-'(x_0)$, defined as the left limit
\[f_-'(x_0)=\lim_{x\to x_0^-}\frac{f(x)-f(x_0)}{x-x_0},\] exists.  

\item 
For any $x_0\in [a, b)$, we say that the function is right-differentiable at $x_0$ provided that the right derivative at $x_0$, $f_+'(x_0)$, defined as the right limit
\[f_+'(x_0)=\lim_{x\to x_0^+}\frac{f(x)-f(x_0)}{x-x_0},\] exists.  

\item For any $x_0\in (a, b)$, $f$ is differentiable at $x_0$ if and only if it is both left differentiable and right differentiable at $x_0$.

\item We say that the function  $f:[a,b]\rightarrow \mathbb{R}$ is differentiable if it is differentiable at all $x_0\in (a, b)$, right differentiable at $a$ and left differentiable at $b$.
\end{enumerate}

\end{remark}
 In the following, we will mainly discuss derivatives of functions defined on open intervals. The extension to closed intervals is straighforward by considering the one-sided derivatives at the end points.

\begin{example}{\linkt Example \ref{23021303} Revisited}
The function $f(x)=|x|$ is left differentiable and right differentiable at $x=0$, with 
\[f_-'(0)=-1,\hspace{1cm}f_+'(0)=1.\]
It is not differentiable at $x=0$ since $f_-'(0)\neq f_+'(0)$.
\end{example}

\begin{highlight}{Leibniz Notation for Derivatives}
For the function $y=f(x)$, its derivative $f'(x)$ is also denoted by $\di \frac{dy}{dx}$ or $\di \frac{d}{dx}f(x)$.
\end{highlight}

For example, in Example \ref{23021301}, we have shown that
\[\frac{d}{dx} x^2=2x.\]

In the following, we are going to derive derivative formulas. 
The simplest   derivative formula is the one for the function $f(x)=x^n$, where $n$ is a positive integer.
\begin{proposition}{}
Let $n$ be a positive integer. The function $f(x)=x^n$ is differentiable with derivative 
\[\frac{d}{dx}x^n=nx^{n-1}.\]
\end{proposition}
\begin{myproof}{Proof}
We use the formula
\[x^n-x_0^n=(x-x_0)(x^{n-1}+x^{n-2}x_0+\cdots+xx_0^{n-2}+x_0^{n-1}).\]
Let $f(x)=x^n$. Then\bp
\begin{align*}
 f'(x_0)&=\lim_{x\to x_0}\frac{f(x)-f(x_0)}{x-x_0}\\
&=\lim_{x\to x_0}\frac{x^n-x_0^n}{x-x_0}\\
&=\lim_{x\to x_0}\left(x^{n-1}+x^{n-2}x_0+\cdots+xx_0^{n-2}+x_0^{n-1}\right)\\
&=\underbrace{x_0^{n-1}+x_0^{n-1}+\cdots+x_0^{n-1}+x_0^{n-1}}_{n \;\text{terms}}\\
&=nx_0^{n-1}.
\end{align*}
\end{myproof}

Now let us look at the square root function.
\begin{example}{}
Determine the points where the function $f:[0, \infty)\to\mathbb{R}$, $f(x)=\sqrt{x}$ is differentiable, and find the derivatives at those points.
\end{example}
\begin{solution}{Solution}
First we consider the case where $x_0>0$. When $x>0$ and $x\neq x_0$,
\begin{equation}\label{eq230213_4}
\frac{f(x)-f(x_0)}{x-x_0}=\frac{\sqrt{x}-\sqrt{x_0}}{x-x_0}=\frac{1}{\sqrt{x}+\sqrt{x_0}}.
\end{equation}
Hence, 
\[f'(x_0)=\lim_{x\to x_0}\frac{f(x)-f(x_0)}{x-x_0}=\lim_{x\to x_0}\frac{1}{\sqrt{x}+\sqrt{x_0}}=\frac{1}{2\sqrt{x_0}}.\]
This shows that $f$ is differentiable at $x_0$ with derivative $f'(x_0)=\di \frac{1}{2\sqrt{x_0}}$.

For $x_0=0$, we can only consider the right derivative. When $x>0$,
the formula \eqref{eq230213_4} still holds. However, the limit
\[\lim_{x\rightarrow 0^+}\frac{f(x)-f(0)}{x-0}=\lim_{x\to 0^+}\frac{1}{\sqrt{x}}\] does not exist. Hence, the function $f(x)=\sqrt{x}$ is not differentiable at $x=0$.

\end{solution}

 \begin{figure}[ht]
\centering
\includegraphics[scale=0.18]{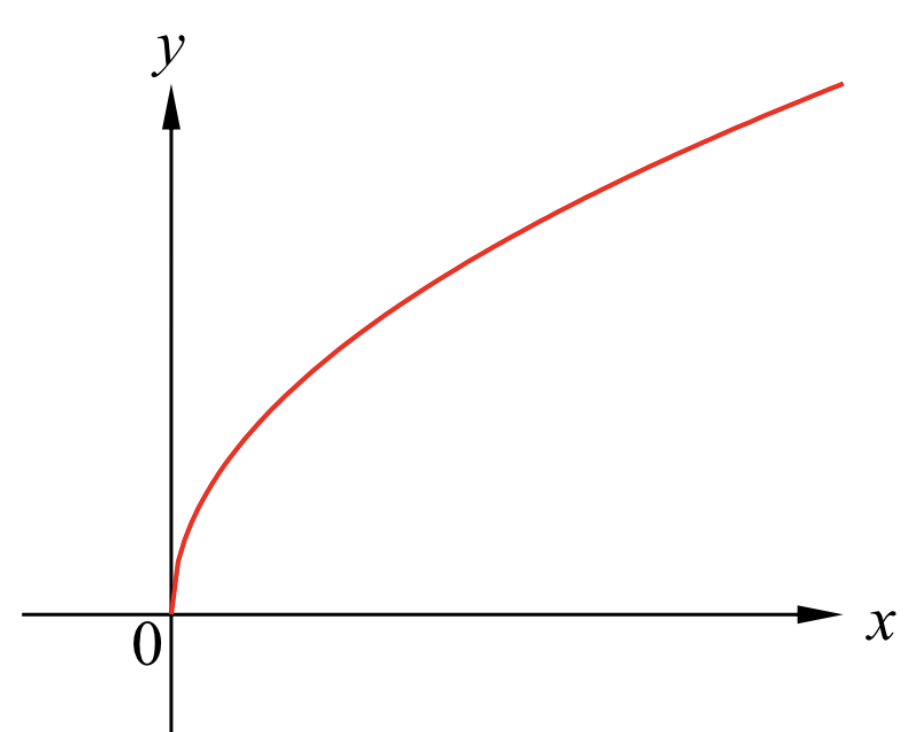}
\caption{  The function $f(x)=\sqrt{x}$ is not differentiable at $x=0$. \fa}\label{figure20}
\end{figure}

Using limit laws, one can find derivatives of linear combinations, products and quotients of  functions.
\begin{proposition}[label=23021305]{Linearity of Derivatives}
Let $x_0$ be a point in $(a, b)$. Given that the functions $f:(a,b)\to \mathbb{R}$ and $g:(a,b)\to \mathbb{R}$ are differentiable at $x_0$. For any constants $\alpha$ and $\beta$, the function $\alpha f+\beta g:(a,b)\rightarrow\mathbb{R}$ is also differentiable at $x_0$ and
\[(\alpha f+\beta g)'(x_0)=\alpha f'(x_0)+\beta g'(x_0).\]

\end{proposition}

\begin{myproof}{Proof}
This is straightforward derivation from the limit laws. By assumption, we have
\[f'(x_0)=\lim_{x\to x_0}\frac{f(x)-f(x_0)}{x-x_0}\quad\text{and}\quad g'(x_0)=\lim_{x\to x_0}\frac{g(x)-g(x_0)}{x-x_0}.\]
It follows that
\begin{align*}
(\alpha f+\beta g)'(x_0)&=\lim_{x\to x_0}\frac{(\alpha f+\beta g)(x)-(\alpha f+\beta g)(x_0)}{x-x_0}\\
&=\lim_{x\to x_0}\left(\alpha\frac{f(x)-f(x_0)}{x-x_0}+\beta\frac{g(x)-g(x_0)}{x-x_0}\right)\\
&=\alpha \lim_{x\to x_0}\frac{f(x)-f(x_0)}{x-x_0}+\beta \lim_{x\to x_0}\frac{g(x)-g(x_0)}{x-x_0}\\
&=\alpha f'(x_0)+\beta g'(x_0).
\end{align*}
\end{myproof}

This formula can be extended to $k$ functions for any positive integer $k$. If $f_1,  \ldots, f_k$ are functions defined on $(a,b)$ and differentiable at the point $x_0$, then for any constants $c_1, \ldots, c_k$, the function $c_1 f_1+\ldots+c_k f_k$ is also differentiable at $x_0$, and
\[(c_1 f_1+\ldots+c_k f_k)'(x_0)=c_1f_1'(x_0)+\ldots+c_kf_k'(x_0).\]

\begin{proposition}[label=23021306]{Product  Rule for Derivatives}
Let $x_0$ be a point in $(a, b)$. Given that the functions $f:(a,b)\to \mathbb{R}$ and $g:(a,b)\to \mathbb{R}$ are differentiable at $x_0$,   the function $(fg) :(a,b)\rightarrow\mathbb{R}$ is also differentiable at $x_0$ and
\[(fg)'(x_0)=  f'(x_0)g(x_0)+f(x_0) g'(x_0).\]

\end{proposition}

\begin{myproof}{Proof}
Again, we are given that
\[f'(x_0)=\lim_{x\to x_0}\frac{f(x)-f(x_0)}{x-x_0}\quad\text{and}\quad g'(x_0)=\lim_{x\to x_0}\frac{g(x)-g(x_0)}{x-x_0}.\] Since differentiability implies continuity, we have
\[\lim_{x\to x_0}f(x)=f(x_0) \quad\text{and}\quad \lim_{x\to x_0}g(x)=g(x_0).\]
Just like the proof of the product rule for limits, we need to do some manipulations.
\[f(x)g(x)-f(x_0)g(x_0)=(f(x)-f(x_0))g(x_0)+f(x)(g(x)-g(x_0)).\]
It follows that
\begin{align*}
(fg)'(x_0)&=\lim_{x\to x_0}\frac{f(x)g(x)-f(x_0)g(x_0)}{x-x_0}\\
&=\lim_{x\to x_0}\left(\frac{f(x)-f(x_0)}{x-x_0}g(x_0)+f(x)\frac{g(x)-g(x_0)}{x-x_0}\right)\\
&=\lim_{x\to x_0} \frac{f(x)-f(x_0)}{x-x_0}\lim_{x\to x_0} g(x_0)+\lim_{x\to x_0} f(x)\lim_{x\to x_0}\frac{g(x)-g(x_0)}{x-x_0}\\
&=f'(x_0)g(x_0)+f(x_0)g'(x_0).
\end{align*}
\end{myproof}
For a different perspective, we denote $f(x_0)$ and $g(x_0)$ by $u$ and $v$ respectively, and let
\[\Delta u =f(x)-f(x_0),\hspace{1cm} \Delta v=g(x)-g(x_0).\]
Then
\begin{align*}
f(x)g(x)-f(x_0)g(x_0)&=(u+\Delta u)(v+\Delta v)-uv\\
&=v\Delta u+u\Delta v +\Delta u\Delta v.
\end{align*}After didiving by $\Delta x=x-x_0$, the term $\Delta u\Delta v/\Delta x$ vanishes in the limit $\Delta x\to 0$, and we obtain the product rule.

\begin{highlight}{General Product Rule}
The product rule for derivatives can be expressed as
\[\frac{d}{dx}(uv)=v\frac{du}{dx}+u\frac{dv}{dx}.\] 
In general, when we have $k$ functions $u_1, u_2, \ldots, u_k$,
the product rule says that
\[\frac{d}{dx}(u_1u_2\cdots u_k)=(u_2\cdots u_k)\frac{d u_1}{dx}
+(u_1u_3\cdots u_k)\frac{du_2}{dx}+\cdots+(u_1\cdots u_{k-1})\frac{du_k}{dx}.\]
This can be proved by induction on $k$.
\end{highlight}
Notice that the formula
\[\frac{d}{dx}x^n=nx^{n-1}, \quad\text{where}\;n\in\mathbb{Z}^+ \] follows from the general product rule and $\di\frac{d}{dx}x=1$.

Finally we turn to the quotient rule.
\begin{proposition}[label=23021307]{Quotient Rule for Derivatives} 
Let $x_0$ be a point in $(a, b)$. Given that the functions $f:(a,b)\to \mathbb{R}$ and $g:(a,b)\to \mathbb{R}$ are differentiable at $x_0$, and $g(x)\neq 0$ for all $x$ in $(a,b) $. Then the function $(f/g):(a,b)\to\mathbb{R}$ is differentiable at $x_0$ and
\[\left(\frac{f}{g}\right)'(x_0)=\frac{f'(x_0)g(x_0)-f(x_0)g'(x_0)}{g (x_0)^2}.\]
\end{proposition}The assumption $g(x)\neq 0$ for all $x$ in $ (a,b)$ is to make sure that the function $f/g$ is well-defined on $(a,b)$. In practice, we only need $g(x_0)\neq 0$ and $g$ is differentiable at $x_0$. For  then we find that $g$ is continuous at $x_0$. The assumption $g(x_0)\neq 0$ will imply that $g(x)\neq 0$ in a neighbourhood of $x_0$.
\begin{myproof}{Proof}First, notice that
\begin{align*}
\frac{f(x)}{g(x)}-\frac{f(x_0)}{g(x_0)}&=\frac{f(x)g(x_0)-f(x_0)g(x)}{g(x)g(x_0)}\\&=\frac{(f(x)-f(x_0))g(x_0)-f(x_0)(g(x)-g(x_0))}{g(x)g(x_0)}.
\end{align*}Using the same reasoning as in the proof of the product rule, we obtain
\begin{align*}
 \left(\frac{f}{g}\right)'(x_0) &=\lim_{x\to x_0}\frac{1}{g(x)g(x_0)}\\&\quad \times\left\{g(x_0)\lim_{x\to x_0}\frac{f(x)-f(x_0)}{x-x_0}-f(x_0)\lim_{x\to x_0}\frac{g(x)-g(x_0)}{x-x_0}\right\}\\
&=\frac{f'(x_0)g(x_0)-f(x_0)g'(x_0)}{g (x_0)^2}.
\end{align*}
\end{myproof}
Again, using the $u$ and $v$ notations, we have
\[\frac{u+\Delta u}{v+\Delta v}-\frac{u}{v}=\frac{(u+\Delta u)v-(v+\Delta v)u}{v(v+\Delta v)}=\frac{v\Delta u-u\Delta v}{v(v+\Delta v)}.\] This gives a different perspective on the quotient rule.

Let us use the quotient rule to derive the derivative for $f(x)=x^n$, when $n$ is a negative integer. 

\begin{proposition}[label=prop230215_1]{}
For any integer $n$,
\begin{equation}\label{eq230213_8}\frac{d}{dx}x^n =nx^{n-1}.\end{equation}

\end{proposition}
\begin{myproof}{Proof}
We have proved the formula   \eqref{eq230213_8} when $n\geq 0$. When $n<0$, let $m=-n$. Then $m$ is a positive integer.  
By quotient rule, we have
\[\frac{d}{dx}x^n=\frac{d}{dx}\frac{1}{x^m}=\frac{x^m \di \frac{d}{dx} 1- \frac{d}{dx}x^m}{x^{2m}}=-\frac{mx^{m-1}}{x^{2m}}=-\frac{m}{x^{m+1}}=nx^{n-1}.\]
Hence, the formula \eqref{eq230213_8} also holds when $n$ is a negative integer.  
\end{myproof}

\begin{definition}{Higher Order Derivatives}
If the function $f:(a,b)\to\mathbb{R}$ is differentiable, its derivative $f':(a,b)\to\mathbb{R}$ is also a function  defined on $(a, b)$. We can investigate whether $f'$ is differentiable. If $f'$ is differentiable at a point $x_0$ in $(a, b)$, we denote its derivative by $f''(x_0)$, called the second (order) derivative of the function $f$ at $x_0$. 

In the same way, we can define the $n^{\text{th}}$-order derivative of the function $y=f(x)$ at  a point $x_0$ for any positive integer $n$. We use the notation \[f^{(n)}(x)\quad\text{or}\quad\frac{d^n y}{dx^n}\] to denote the $n^{\text{th}}$-derivative of the function $y=f(x)$. It is defined recursively by
\[f^{(n)}(x)=\lim_{h\to 0}\frac{f^{(n-1)}(x+h)-f^{(n-1)}(x)}{h},\]
where by default, $f^{(0)}(x)=f(x)$.

We say that a function $f:(a,b)\rightarrow\mathbb{R}$ is $n$ times differentiable if $f^{(n)}(x)$ exists for all $x$ in $(a,b)$. A function is infinitely differentiable if it is $n$ times differentiable for any positive integer $n$.
\end{definition}

\begin{example}{}
Polynomial functions are infinitely differentiable. Moreover, if the degree of a polynomial $p(x)$ is  $n$, then
$p^{(k)}(x)=0$ for all $k\geq n+1$.
\end{example}

\begin{example}[label=23021401]{}
Define the function $f:\mathbb{R}\to\mathbb{R}$ by
\[f(x)=\begin{cases} ax^2,\quad &\text{if}\;x< 1,\\
x+\di \frac{b}{x},\quad &\text{if}\; x\geq 1.\end{cases}
\]Find the values of $a$ and $b$ so that $f$ is differentiable.

\end{example}
\begin{solution}{Solution}
The function $f$ is differentiable at any point  $x_0$ in the interval $(-\infty, 1)$ or the interval $(1, \infty)$.

For $f$  to be  differentiable, $f$ has to be continuous and differentiable at $x=1$. For $f$ to be continuous at $x=1$, we must have
\[\lim_{x\rightarrow 1^-}f(x)=\lim_{x\rightarrow 1^+}f(x).\]
This gives
\[a=1+b.\]
For $f$ to be differentiable at $x=1$, we must have
\[\lim_{x\rightarrow 1^-}\frac{f(x)-f(1)}{x-1}=\lim_{x\rightarrow 1^+}\frac{f(x)-f(1)}{x-1}.\]
Notice that
\[\lim_{x\rightarrow 1^-}\frac{f(x)-f(1)}{x-1}=\left.\frac{d}{dx}\right|_{x=1}ax^2=2a,\]
\[\lim_{x\rightarrow 1^+}\frac{f(x)-f(1)}{x-1}=\left.\frac{d}{dx}\right|_{x=1}\left(x+\frac{b}{x}\right)=1-b.\]
Hence, we must have
\[2a=1-b.\]
Solving for $a$ and $b$, we have 
\[a= \frac{2}{3},\;b=-\frac{1}{3}.\]
\end{solution}

\vp
 \begin{figure}[ht]
\centering
\includegraphics[scale=0.2]{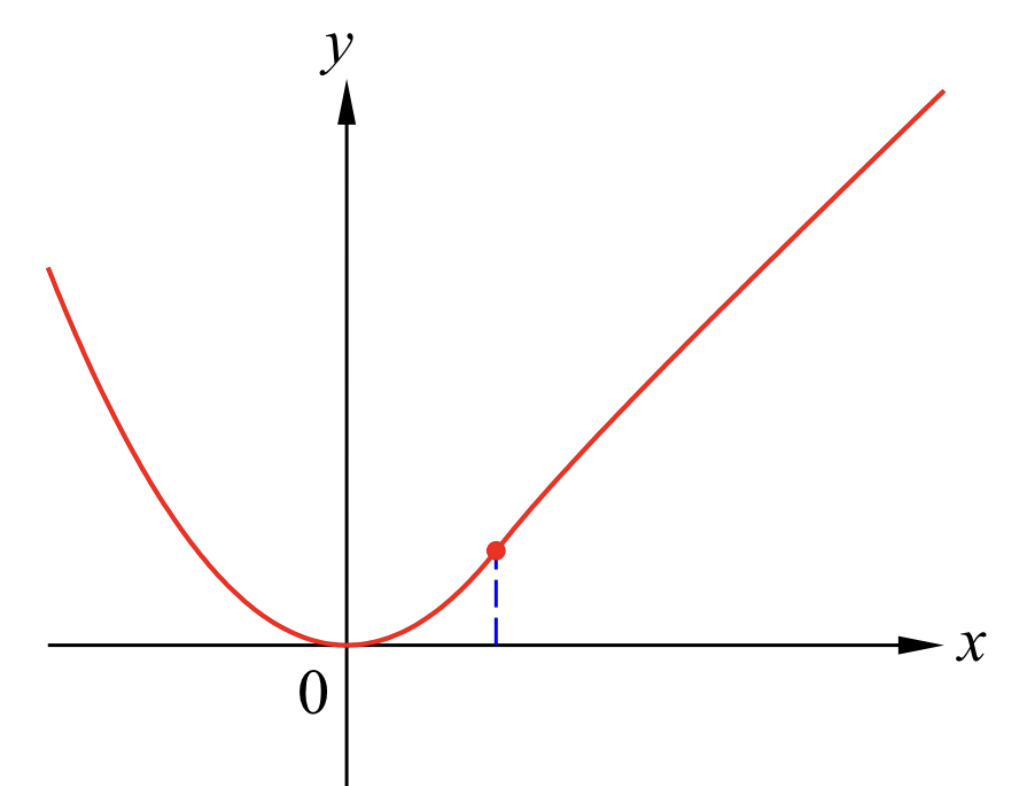}
\caption{  The function $f(x)$ defined in Example \ref{23021401}. \fa}\label{figure21}
\end{figure}

\vp
 
\noindent
{\bf \large Exercises  \thesection}
\setcounter{myquestion}{1}
\begin{question}{\themyquestion} 
Define the function $f:\mathbb{R}\to\mathbb{R}$ by
\[f(x)=\begin{cases} ax^2+ x,\quad &\text{if}\;x< 1,\\
bx+\di \frac{3 }{x^2},\quad &\text{if}\;x\geq 1.\end{cases}
\]Find the values of $a$ and $b$ so that $f$ is differentiable.
\end{question}
\atc
\begin{question}{\themyquestion} 
Let $x_0$ be a point in $(a, b)$. Given that $f:[a,b]\to\mathbb{R}$ is a continuous function defined on $[a, b]$ and differentiable at $x_0$. Let
$g:[a,b]\to\mathbb{R}$ be the function defined by
\[g(x)=\begin{cases} \di\frac{f(x)-f(x_0)}{x-x_0},\quad& \text{if}\;x\in [a,b]\setminus\{ x_0\}\\f'(x_0),\quad & \text{if}\;x=x_0.\end{cases}\]Show that $g:[a,b]\to\mathbb{R}$ is a continuous function.
\end{question}
 
\atc
\begin{question}{\themyquestion} 
Let $x_0$ be a point in $(a, b)$ and let   $f:(a,b)\to\mathbb{R}$ be a function defined on $(a, b)$. 
\begin{enumerate}[(a)]
\item If $f:(a,b)\to\mathbb{R}$ is differentiable at $x_0$, show that
\[\lim_{h\to 0}\frac{f(x_0+h)-f(x_0-h)}{2h}=f'(x_0).\]
\item If the limit \[\lim_{h\to 0}\frac{f(x_0+h)-f(x_0-h)}{2h}\] exists, is $f$ necessarily differentiable at $x_0$?
\end{enumerate}
\end{question}

\vp

\section{Chain Rule and Derivatives of Inverse Functions}\label{sec3.2}

In this section, we are going to derive derivative formulas for composite functions and inverse functions. First we discuss a different perspective for differentiability of a function at a point.

\begin{highlight}{Differentiability}
Let $x_0$ be a point in the interval $(a, b)$ and let $f:(a,b)\rightarrow \mathbb{R}$ be a function defined on $(a, b)$. 
If $f$ is differentiable at $x_0$, then
\[f'(x_0)=\lim_{h\to 0}\frac{f(x_0+h)-f(x_0)}{h}.\]
This implies that
\[\lim_{h\to 0}\frac{f(x_0+h)-f(x_0)-f'(x_0)h}{h}=0.\]
Conversely, if there is a number $c$ such that
\[\lim_{h\to 0}\frac{f(x_0+h)-f(x_0)-ch}{h}=0,\]
limit laws imply that
\[c=\lim_{h\to 0}\frac{f(x_0+h)-f(x_0)}{h}.\]
This implies that $f$ is differentiable at $x_0$ and $f'(x_0)=c$.
In other words, the function $f$ is differentiable at $x_0$ if and only if there is a number $c$ such that
\[\lim_{h\to 0}\frac{f(x_0+h)-f(x_0)-ch}{h}=0.\]

Since $x_0\in (a, b)$, there is an $r>0$ such that $(x_0-r, x_0+r)\subset (a,b)$.  For a given real number $c$, let
$\varepsilon:(-r,r)\to\mathbb{R}$ be the function defined by
\[\varepsilon(h)=\frac{f(x_0+h)-f(x_0)-ch}{h}.\] \end{highlight}\begin{highlight}{} Then the differentiability of $f$ at $x_0$ is equivalent to $\di\lim_{h\rightarrow 0}\varepsilon(h)=0$.
Hence, $f:(a,b)\to\mathbb{R}$ is differentiable at $x_0$ if and only if there is a number $c$ and a function $\varepsilon(h)$ such that
\[f(x_0+h)=f(x_0)+ch+h\varepsilon(h),\]
and
\[\varepsilon(h)\to 0\quad\text{when}\;h\to 0.\]
\end{highlight}

\begin{theorem}{Chain Rule}
Given that  $f:(a,b)\to \mathbb{R}$ and $g:(c,d)\to\mathbb{R}$ are  functions such that $f(a,b)\subset (c,d)$. If $x_0$ is a point in $(a,b)$, $f$ is differentiable at $x_0$, $g$ is differentiable at $f(x_0)$, then the composite function $ (g\circ f):(a, b)\to \mathbb{R}$ is differentiable at $x_0$ and 
\[(g\circ f)'(x_0)=g'(f(x_0))f'(x_0).\]

\end{theorem}
\begin{myproof}{Proof}
Let  $y_0=f(x_0)$, and define the functions $\varepsilon_1(h)$ and $\varepsilon_2(k)$ by
\begin{align*}
\varepsilon_1(h)&=\frac{f(x_0+h)-f(x_0)-f'(x_0)h}{h},\\
\varepsilon_2(k)&=\frac{g(y_0+k)-g(y_0)-g'(y_0)k}{k}.
\end{align*}Since $f$ is differentiable at $x_0$ and $g$ is differentiable at $y_0$, we have
$\di\lim_{h\to 0}\varepsilon_1(h)=0$ and $\di \lim_{k\to 0}\varepsilon_2(k)=0$.
Let
\[k(h)=f(x_0+h)-f(x_0)=f'(x_0)h+\varepsilon_1(h)h.\] Then by the definitions of $\varepsilon_1(h)$ and $\varepsilon_2(k)$,
\begin{align*}
(g\circ f)(x_0+h)-(g\circ f)(x_0)&=g(y_0+k(h))-g(y_0)\\
&=g'(y_0)k(h)+\varepsilon_2(k(h))k(h)\\
&=g'(y_0)f'(x_0)h+\varepsilon_3(h)h,
\end{align*}\bp
where 
\[\varepsilon_3(h)=g'(y_0)\varepsilon_1(h)+\varepsilon_2(k(h))\frac{k(h)}{h}.\]
Since $f$ is differentiable at $x_0$, \[\lim_{h\to 0}\frac{k(h)}{h}=\lim_{h\to 0}\frac{f(x_0+h)-f(x_0)}{h}=f'(x_0).\] This implies that $\di\lim_{h\rightarrow 0}k(h) =0$. By limit law  for composite functions,
\[\lim_{h\to 0}\varepsilon_2(k(h))=\lim_{k\to 0}\varepsilon_2(k)=0.\]
Limit laws then imply that
\[\lim_{h\to 0}\frac{(g\circ f)(x_0+h)-(g\circ f)(x_0)-g'(y_0)f'(x_0)h}{h}=\lim_{h\to 0}\varepsilon_3(h)=0.\]
This proves that the function $g\circ f$ is differentiable at $x_0$ and 
\[(g\circ f)'(x_0)=g'(y_0)f'(x_0)=g'(f(x_0))f'(x_0).\]
\end{myproof}
Heuristically, if we let $u=f(x)$ and $y=g(u)=g(f(x))$,   chain rule says that
\[\frac{dy}{dx}=\frac{dy}{du}\times\frac{du}{dx},\]
which is the limit of 
\[\frac{\Delta y}{\Delta x}=\frac{\Delta y}{\Delta u}\times\frac{\Delta u}{\Delta x}\]when $\Delta x\to 0$.
The rigorous proof we give above do not use this because we might face the problem that $\Delta u=f(x)-f(x_0)$ can be zero even when $x\neq x_0$. 

\begin{example}
{}
Let $f:\mathbb{R}\to\mathbb{R}$ be a differentiable function and let $a$ be a constant. Show that the function $g:\mathbb{R}\to\mathbb{R}$ defined by $g(x)=f(ax)$ is differentiable, and $g'(x)=af'(ax)$.

\end{example}
\begin{solution}{Solution}
The function $u:\mathbb{R}\to\mathbb{R}$, $u(x)=ax$ is differentiable with $u'(x)=a$. By  chain rule, the function $g(x)=(f\circ u)(x)$ is also differentiable and
\[g'(x)=f'(u(x))u'(x)=af'(ax).\]
\end{solution}

\begin{example}{}
Given that the function $f:(0,2)\to \mathbb{R}$ is differentiable at $x=1$ and $f'(1)=a$, find the value of
\[\lim_{x\rightarrow 1}\frac{f(x^3)-f(1)}{x-1}\] in terms of $a$.
\end{example}
\begin{solution}{Solution}
Let $g(x)=x^3$. Then $g(1)=1$ and $g$ is differentiable at $x=1$ with $g'(1)=3$.
\[\lim_{x\rightarrow 1}\frac{f(x^3)-f(1)}{x-1}=\lim_{x\rightarrow 1}\frac{(f\circ g)(x )-(f\circ g)(1)}{x-1}.\]
Since $g$ is differentiable at $x=1$ and  $f$ is differentiable at $g(1)$,   chain rule implies that
\[\lim_{x\rightarrow 1}\frac{f(x^3)-f(1)}{x-1}=(f\circ g)'(1)=f'(g(1))g'(1)=3f'(1)=3a.\]
\end{solution}

Recall that we have proved in Section \ref{sec2.6} that if $I$ is an interval, $f:I\to\mathbb{R}$ is strictly monotonic and continuous, then $f$ is invertible and $f^{-1}:f(I)\to \mathbb{R}$ is also continuous. The strictly monotonicity is a necessary and sufficient condition for a continuous function to be one-to-one. If $x_0$ is a point in the interior of $I$, and $f$ is differentiable at $x_0$, we can ask whether the inverse function $f^{-1}$ is differentiable at the point $y_0=f(x_0)$. 
Since $(f^{-1}\circ f)(x)=x$ for all $x\in I$, if $f^{-1}$ is differentiable at $y_0$,  chain rule implies that
\[(f^{-1})'(y_0)f'(x_0)=(f^{-1})'(f(x_0))f(x_0)=1.\]
Therefore, a necessary condition for $f^{-1}$ to be differentiable at $y_0$ is $f'(x_0)$ cannot be zero. In the following theorem, we show that this condition is also sufficient. 

\begin{theorem}[label=thm230218_9]{Derivative for Inverse Function}
Let $I$ be an open interval containing   $x_0$, and let $f:I\rightarrow\mathbb{R}$ be a function that is strictly monotonic  and continuous. If $f$ is differentiable at $x_0$ and $f'(x_0)\neq 0$, the inverse function $f^{-1}:f(I)\to \mathbb{R}$ is differentiable at $y_0=f(x_0)$, and 
\[(f^{-1})'(y_0)=\frac{1}{f'(x_0)}.\]
\end{theorem}
The formula for $(f^{-1})'(y_0)$ would follow from the chain rule if we know apriori that $f^{-1}$ is   differentiable at $y_0$. The gist of this theorem is to state that $f^{-1}$ is indeed differentiable at $y_0$.

\begin{myproof}{Proof}Without loss of generality, assume that $f$ is strictly increasing. 
By Theorem \ref{23021106}, $f^{-1}:f(I)\to\mathbb{R}$ is also continuous.
There is a  $\delta>0$ so that $[x_0-\delta, x_0+\delta]\subset I$. Then $(f(x_0-\delta), f(x_0+\delta))$ is an open interval in $f(I)$ containing the point $y_0$.  This implies that there is an $r>0$ so that $(y_0-r, y_0+r)\subset f(I)$. For any $k\in (-r, r)$, let 
\[h(k)=f^{-1}(y_0+k)-f^{-1}(y_0).\]Then $h$ is a strictly increasing continuous function of $k$ and  $\di\lim_{k\to 0}h(k)=0$. Notice that
\[y_0+k=f(x_0+h(k)).\]
Therefore,
\[\frac{f^{-1}(y_0+k)-f^{-1}(y_0)}{k}=\frac{h(k)}{f(x_0+h(k))-f(x_0)}.\]
\bp
Hence, by limit laws for quotients and composite functions,
we find that
\begin{align*}
\lim_{k\to 0}\frac{f^{-1}(y_0+k)-f^{-1}(y_0)}{k}&=\frac{1}{\di \lim_{k\rightarrow 0}\frac{f(x_0+h(k))-f(x_0)}{h(k)}}\\&=\frac{1}{\di \lim_{h\rightarrow 0}\frac{f(x_0+h)-f(x_0)}{h}}\\&=\frac{1}{f'(x_0)}.\end{align*}
This proves that $f^{-1}$ is differentiable at $y_0$ and
\[(f^{-1})'(y_0)=\frac{1}{f'(x_0)}.\]

 \end{myproof}

As a corollary, we have the following.
\begin{corollary}{}
Let $I$ be an open interval, and let $f:I\to\mathbb{R}$ be a strictly monotonic differentiable function. If $f'(x)\neq 0 $ for all $x\in I$, then the inverse function $f^{-1}:f(I)\to\mathbb{R}$ is also a strictly monotonic differentiable function with
\[(f^{-1})'(x)=\frac{1}{f'(f^{-1}(x))}.\]
\end{corollary}
 
\begin{example}{}
Let $r$ be a rational number, and let $f:(0,\infty)\to \mathbb{R}$ be the function $f(x)=x^r$. Show that $f$ is differentiable and 
\[f'(x)=rx^{r-1}.\]

\end{example}
\begin{solution}{Solution}
First we consider the case  $r=1/n$, where $n$ is a positive integer. The function $f(x)=x^{1/n}$ is the inverse of the function $g(x)=x^n$, which is differentiable and strictly increasing. Hence, $f(x)=x^{1/n}$ is differentiable and strictly increasing. Moreover, since $g'(x)=nx^{n-1}$, we have
\[f'(x)=\frac{1}{g'(f(x))}=\frac{1}{g'(x^{1/n})}=\frac{1}{n(x^{1/n})^{n-1}}=\frac{1}{n}x^{\frac{1}{n}-1}.\]
Now for a general rational number $r$, there is an integer $p$ and a positive integer $q$ such that
$r=p/q$. It follows that
\[f(x)=(x^p)^{1/q}=(g\circ h)(x),\]
where
\[h(x)=x^p,\hspace{1cm} g(x)=x^{1/q}.\]
By Proposition \ref{prop230215_1}, $h'(x)=px^{p-1}$. We have just shown that $\di g'(x)=\frac{1}{q}x^{\frac{1}{q}-1}$. By chain rule,
\[f'(x)=g'(h(x))h'(x)=\frac{1}{q}\left(x^p\right)^{\frac{1}{q}-1} \times px^{p-1}=\frac{p}{q}x^{\frac{p}{q}-1}=rx^{r-1}.\]
\end{solution}

\vp
\noindent 
{\bf \large Exercises  \thesection}
\setcounter{myquestion}{1}
\begin{question}{\themyquestion}
Given that the function $f:(0,\infty)\to \mathbb{R}$ is defined by
\[f(x)=\frac{1}{\sqrt{4+x^2}}.\]
\begin{enumerate}[(a)]
\item Show that $f$ is one-to-one.
\item Show that $f$ is differentiable.
\item Show that $f^{-1}$ exists and is differentiable.
\item Find $f^{-1}(x)$ and $(f^{-1})'(x)$.
\end{enumerate}
\end{question}

\atc

\begin{question}{\themyquestion}
Let $a$ be a positive number. 
Recall that a function $f:(-a, a)\to\mathbb{R}$ is even if and only if
\[f(-x)=f(x)\hspace{1cm}\text{for all}\;x\in (-a,a);\]
and a function $f:(-a, a)\to\mathbb{R}$ is odd if and only if
\[f(-x)=-f(x)\hspace{1cm}\text{for all}\;x\in (-a,a).\]
Let $f:(-a,a)\to\mathbb{R}$ be a differentiable function.
\begin{enumerate}[(a)]
\item If $f$ is even, show that $f'$ is odd.
\item If $f$ is odd, show that $f'$ is even.

\end{enumerate}
\end{question}
\vp

\section{The Mean Value Theorem  and Local Extrema}\label{sec3.3}

The mean value theorem is one of the most important theorems in analysis. 
We will first prove a special case of the mean value theorem called Rolle's theorem. To prove this, we need the extreme value theorem, which asserts the existence of global maximum and global minimum for a continuous function defined on a closed and bounded interval. As a matter of fact, what we actually need is a local extremum, which we define as follows.

\begin{definition}{Local Maximum and Local Minimum}
Let $D$ be a subset of real numbers that contains the point $x_0$, and let $f:D\to\mathbb{R}$ be a function defined on $D$. 
\begin{enumerate}[1.]
\item The point  $x_0$ is a local maximizer of $f$
provided that there is a $\delta>0$ such that for all $x$ in $D$ with $|x-x_0|<\delta$,   we have
\[f(x)\leq f(x_0).\]  The value $f(x_0)$ is then  a local maximum value of $f$.
\item  The point  $x_0$ is a local minimizer of $f$
provided that there is a $\delta>0$ such that for all $x$ in $D$ with $|x-x_0|<\delta$,   we have
\[f(x)\geq f(x_0).\]  The value $f(x_0)$ is then  a local minimum value of $f$.
\item The point $x_0$ is a local extremizer if it is a local maximizer or a local minimizer.
 The value $f(x_0)$ is a local extreme value if it is a local maximum value or a local minimum value.
\end{enumerate}
\end{definition}
The definition of local extremum that we give here is quite general. We do not impose conditions on the set $D$, nor require $x_0$ to be an interior point of $D$. Other mathematicians might define it differently. Under our definition, a global extremum of a function is also a local extremum of the function.

 \begin{figure}[ht]
\centering
\includegraphics[scale=0.2]{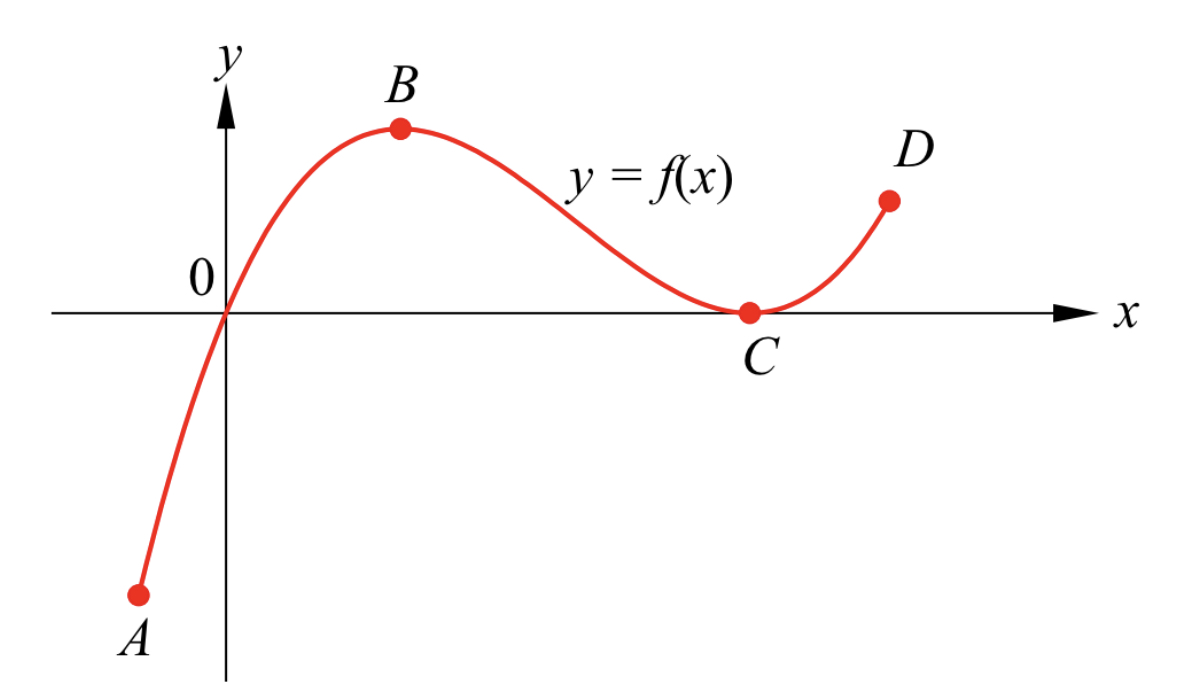}
\caption{  The function $y=f(x)$ has local maxima at the points $B$ and $D$, and local minima at the points $A$ and $C$. The point $A$ is also where global minimum appears; while the point $B$ is where the global maximum appears.\fa}\label{figure22}
\end{figure}

Derivative is an useful tool in the search for local extrema. 
When a local extremizer of a function is an interior point of the domain, and $f$ is differentiable at that point, the derivative of the function can only be zero at that point.
\begin{theorem}[label=thm230215_2]{}
Let $(a, b)$ be a neighbourhood of the point $x_0$, and let $f:(a,b)\to\mathbb{R}$ be a function defined on $(a,b)$. If $x_0$ is a local extremizer of $f$, and $f$ is differentiable at $x_0$, then $f'(x_0)=0$.
\end{theorem}
\begin{myproof}{Proof}Without loss of generality, assume that $x_0$ is a local maximizer. Then there is a $\delta>0$ such that $(x_0-\delta, x_0+\delta)\subset (a,b)$, and for all $x$ in $(x_0-\delta, x_0+\delta)$, $f(x)\leq f(x_0)$. 
Since $f$ is differentiable at $x_0$, the limit
\[\lim_{x\to x_0}\frac{f(x)-f(x_0)}{x-x_0}\] exists and is equal to $f'(x_0)$. This implies that the left limit and the right limit both exist and both equal to $f'(x_0)$.
Namely,
\[f'(x_0)=\lim_{x\to x_0^-}\frac{f(x)-f(x_0)}{x-x_0}=\lim_{x\to x_0^+}\frac{f(x)-f(x_0)}{x-x_0}.\]\bp
For the left limit, when $x$ is in $(x_0-\delta, x_0)$, $x-x_0<0$ and $f(x)-f(x_0)\leq 0$. Therefore,
\[\frac{f(x)-f(x_0)}{x-x_0}\geq 0\hspace{1cm}\text{when}\;x\in (x_0-\delta, x_0).\]
Taking the $x\to x_0^-$ limit, we find that $f'(x_0)\geq 0$. For the right limit, when $x$ is in $(x_0, x_0+\delta)$, $x-x_0>0$ but $f(x)-f(x_0)\leq 0$.  Therefore,
\[\frac{f(x)-f(x_0)}{x-x_0}\leq 0\hspace{1cm}\text{when}\;x\in (x_0, x_0+\delta).\]
Taking the $x\to x_0^+$ limit, we find that $f'(x_0)\leq 0$. Since the left limit shows that $f'(x_0)\geq 0$ while the right limit shows that $f'(x_0)\leq 0$,  we conclude that $f'(x_0)=0$.
\end{myproof}
This theorem gives a necessary condition for a function $f:(a,b)\rightarrow\mathbb{R}$ to have a local extremum at a point where it is differentiable. Notice that it cannot be applied if the local extremizer is not an interior point of the domain. 

\begin{definition}{Stationary Points}
Let $D$ be a subset of real numbers and let $f:D\to\mathbb{R}$ be a function defined on $D$. If $x_0$ is an interior point of $D$,   $f$ is differentiable at $x_0$ and $f'(x_0)=0$, we call $x_0$ a stationary point of the function $f$. 
\end{definition}

Hence, Theorem \ref{thm230215_2} says that if $x_0$ is an interior point of $D$, and the function $f:D\to\mathbb{R}$ is differentiable at $x_0$, a necessary condition for $x_0$ to be a local extremum of the function $f$ is $x_0$ must be a stationary point.
Nevertheless, this condition is not sufficient. For example, the function $f(x)=x^3$ has a stationary point at $x=0$, but $x=0$ is not a local extremizer of the funnction.

Now let us return to the mean value theorem. 
As a motivation, let us consider the distance $s$ travelled by an object as a function of time $t$. We have discussed in Section \ref{sec3.1} that to find the instantaneous speed of the object at a particular time $t_0$, we first find the average speed over the time interval from $t_0$ to $t_0+\Delta t$, and take the limit  $\Delta t\to 0$. Namely, the instantaneous speed at time $t_0$ is
\[\lim_{\Delta t\to 0}\frac{s(t_0+\Delta t)-s(t_0)}{\Delta t},\] which is precisely $s'(t_0)$, the derivative of $s(t)$ at $t=t_0$. The mean value theorem asserts that the average speed of the object in a time interval $[t_1, t_2]$ must equal to the instantaneous speed $s'(t_0)$ for some $t_0$ in that interval. Intuitively, this is something one would expect to be true.

Now let us prove a special case of the mean value theorem.
\begin{theorem}[label=thm_Rolles]{Rolle's Theorem}
Let $f:[a,b]\to\mathbb{R}$ be a function that satisfies the following conditions.
\begin{enumerate}[(i)]
\item $f:[a,b]\to\mathbb{R}$ is continuous.
\item $f:(a, b)\to\mathbb{R}$ is differentiable.
\item $f(a)=f(b)$.
\end{enumerate}
Then there is a point $x_0$ in $(a, b)$ such that $f'(x_0)=0$. 
\end{theorem}
\begin{myproof}{Proof}
Since $f:[a,b]\to\mathbb{R}$ is a continuous function defined on a closed and bounded interval, the extreme value theorem says that it must have minimum value and maximum value. In other words, there are two points $x_1$ and $x_2$ in  $[a,b]$ such that
\[f(x_1)\leq f(x)\leq f(x_2)\hspace{1cm}\text{for all}\;x\in [a,b].\] Notice that $x_1$ and $x_2$ are also local extremizers of the function $f:[a,b]\to\mathbb{R}$.
If $f(x_1)=f(x_2)$, then $f$ is a constant function. In this case,  $f'(x_0)=0$ for all $x_0$ in $(a, b)$.
If $f(x_1)\neq f(x_2)$, then $f(x_1)<f(x_2)$. Since $f(a)=f(b)$, either $x_1$ or $x_2$ must be in the open interval $(a, b)$. In other words, there is a local extremizer $x_0$ in the interval $(a,b)$. Since $f$ is differentiable at $x_0$, Theorem \ref{thm230215_2}  says that we must have $f'(x_0)=0$.
In either case,   there is an $x_0$ in $(a, b)$ satisfying $f'(x_0)=0$.
\end{myproof}

 \begin{figure}[ht]
\centering
\includegraphics[scale=0.2]{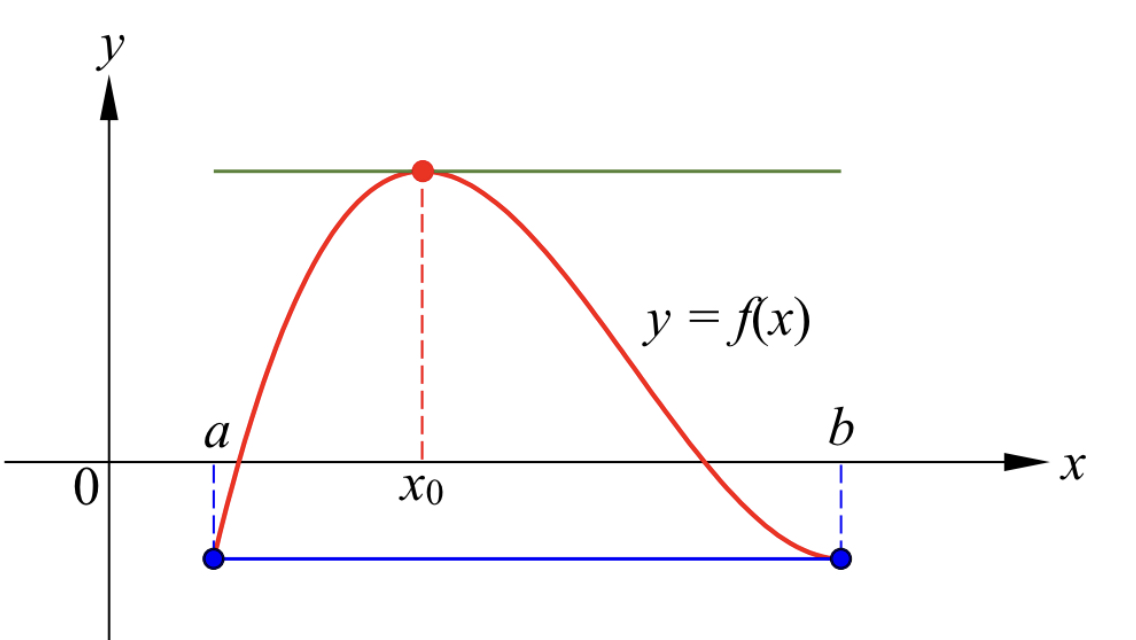}
\caption{  The Rolle's theorem.\fa}\label{figure23}
\end{figure}

Now we can prove the mean value theorem.
\begin{theorem}[label=thm_mvt]{Mean Value Theorem}
Let $f:[a,b]\to\mathbb{R}$ be a function that satisfies the following conditions.
\begin{enumerate}[(i)]
\item $f:[a,b]\to\mathbb{R}$ is continuous.
\item $f:(a, b)\to\mathbb{R}$ is differentiable.
 
\end{enumerate}
Then there is a point $x_0$ in $(a, b)$ such that \[f'(x_0)=\frac{f(b)-f(a)}{b-a}.\]
\end{theorem}
The mean value theorem stated in Theorem \ref{thm_mvt} is also referred to as Lagrange's mean value theorem.
Notice that Rolle's theorem is a special case of the mean value theorem where $f(a)=f(b)$. The quantity 
\[\frac{f(b)-f(a)}{b-a}\] gives the average rate of change of the function $f(x)$ over the interval $[a, b]$, and the mean value theorem says that this average rate of change is equal to the rate of change at a particular point. To prove the mean value theorem, we apply a transformation to the function $f(x)$ to get a function $g(x)$ that satisfies the conditions in the Rolle's theorem.  
\begin{myproof}{Proof}
Let $g:[a,b]\to\mathbb{R}$ be the function defined by \[g(x)= f(x)-mx,\]
where the constant $m$ is determined by $g(a)=g(b)$.
This gives
\[f(a)-ma=f(b)-mb,\]
and so
\[m=\frac{f(b)-f(a)}{b-a}.\]Notice that the function $g:[a,b]\to\mathbb{R}$ is continuous, and $g:(a, b)\to \mathbb{R}$ is differentiable with
\[g'(x)=f'(x)-m=f'(x)-\frac{f(b)-f(a)}{b-a}.\]
By construction, $g(a)=g(b)$. Hence, we can apply Rolle's theorem to the function $g$ and conclude that there is a point $x_0$ in $(a,b)$ such that $g'(x_0)=0$. For this point $x_0$,
\[f'(x_0)=\frac{f(b)-f(a)}{b-a}.\]This proves the mean value theorem.
\end{myproof}

 \begin{figure}[ht]
\centering
\includegraphics[scale=0.2]{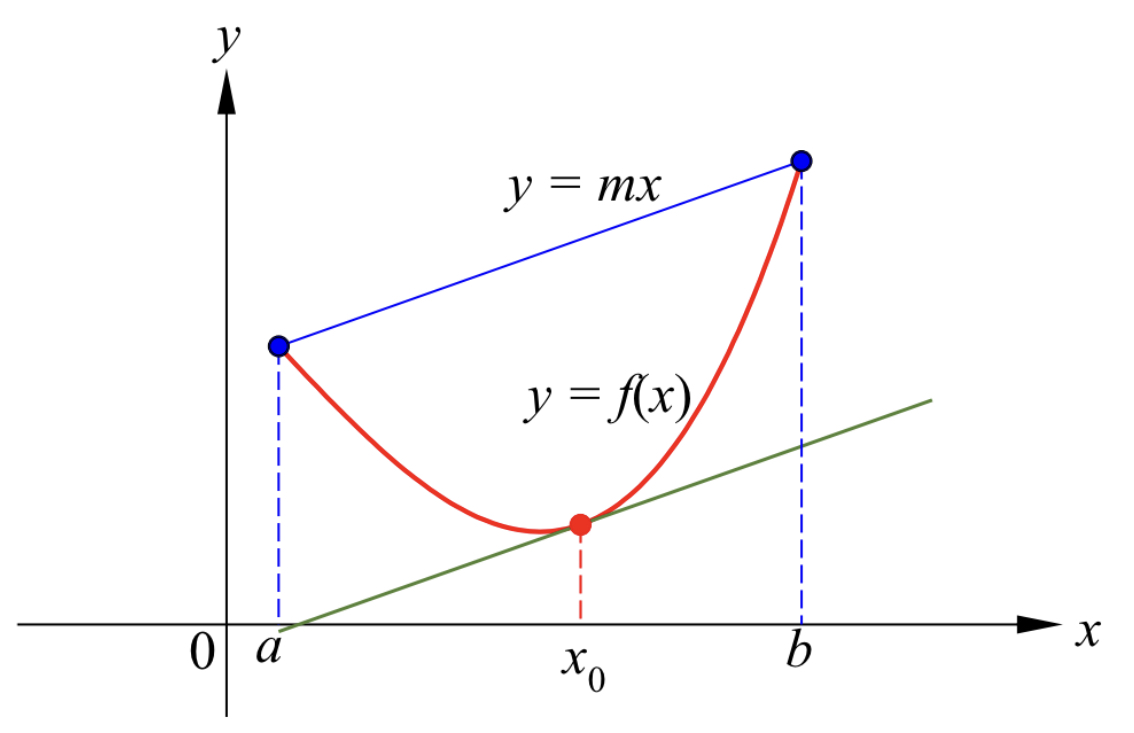}
\caption{  The mean value theorem.\fa}\label{figure24}
\end{figure}
Notice that for the mean value theorem to hold, the function $f:[a,b]\to\mathbb{R}$ do not need to be differentiable at the end points of the interval $[a,b]$, and the point $x_0$ is guaranteed to be a point in the interior of the interval. 

The mean value theorem has very wide applications. We will discuss a few in this section.

Recall that the derivative of a constant function is 0. The converse is not obvious, but it is an easy consequence of the mean value theorem.
\begin{lemma}[label=lemma230215_2]{}
If the function $f:[a,b]\to\mathbb{R}$ is continuous on $[a, b]$, differentiable on $(a,b)$, and $f'(x)=0$ for all $x\in (a,b)$, then $f$ is a constant function.
\end{lemma}
\begin{myproof}{}
Take any $x\in (a, b]$. Then $f$ is continuous on $[a, x]$, differentiable an $(a, x)$. Therefore, we can apply mean value theorem to conclude that there is a point $c$ in $(a, x)$  such that
\[\frac{f(x)-f(a)}{x-a}=f'(c)=0.\]
This proves that $f(x)=f(a)$. Therefore, the function $f$ is a constant.
\end{myproof}

From this, we immediately obtain the following.
\begin{theorem}[label=thm230215_3]{ }
Assume that the functions $f:[a,b]\to\mathbb{R}$ and $g:[a,b]\to\mathbb{R}$ are continuous on $[a,b]$, differentiable on $(a,b)$, and 
\[f'(x)=g'(x)\hspace{1cm}\text{for all}\;x\in (a,b).\]
Then there is a constant $C$ such that
\[f(x)=g(x)+C\hspace{1cm}\text{for all}\; x\in [a,b].\]
\end{theorem}
\begin{myproof}{Proof}
Define the function $h:[a,b]\to\mathbb{R}$ by
\[h(x)=f(x)-g(x).\] Then the function $h$ is continuous on $[a, b]$, differentiable on $(a,b)$, and $h'(x)=0$ for all $x\in (a,b)$. By Lemma \ref{lemma230215_2}, $h$ is a cosntant function. Namely, there is a constant $C$ so that $h(x)=C$ for all $x\in [a, b]$. Therefore,
\[f(x)=g(x)+C\hspace{1cm}\text{for all}\; x\in [a,b].\]
\end{myproof}

\begin{highlight}{}
Theorem \ref{thm230215_3} implies   \emph{the identity criterion}, which says that if two  functions are differentiable in an open interval, their derivatives are the same, and their values at a single point in the interval coincide, then these two functions must be identical.
\end{highlight}

We have seen that if $p(x)$ is a polynomial of degree $n$, and $k$ is an integer larger than $n$, then the $k^{\text{th}}$-order derivative of $p(x)$ is identically zero. Using Lemma \ref{lemma230215_2}, we can prove that the converse is also true.

\begin{example}{}
Let $n$ be a nonnegative integer. Assume that the function $p:\mathbb{R}\to\mathbb{R}$ is $(n+1)$ times differentiable and $p^{(n+1)}(x)=0$ for all real numbers $x$. Then $p(x)$ is a polynomial of degree at most $n$.
\end{example}
\begin{myproof}{Proof}
We prove this by induction on $n$. When $n=0$, the statement says that if $p:\mathbb{R}\to\mathbb{R}$ is a differentiable function and $p'(x)=0$ for all $x\in\mathbb{R}$, then $p(x)$ is a polynomial of degree  0. Since a polynomial of degree 0 is a constant, this statement is true by Lemma \ref{lemma230215_2}.

Now let $n\geq 1$, and assume that
 we have proved that for any $k<n$, if $q:\mathbb{R}\to\mathbb{R}$ is a function that is $(k+1)$ times differentiable and $q^{(k+1)}(x)=0$ for all real numbers $x$, then $q(x)$ is a polynomial of degree at most $k$. \bp Let $p:\mathbb{R}\to\mathbb{R}$ be a function that is $(n+1)$ times differentiable and $p^{(n+1)}(x)=0$ for all real numbers $x$. Lemma \ref{lemma230215_2} says that there is a constant $C$ such that
\[p^{(n)}(x)=C.\]
Consider the function $q:\mathbb{R}\to\mathbb{R}$ defined by
\[q(x)=p(x)-\frac{C}{n!}x^n.\]
It is $n$ times differentiable and
\[q^{(n)}(x)=p^{(n)}(x)-C=0\hspace{1cm}\text{for all}\;x\in\mathbb{R}.\]
By inductive hypothesis,
$q(x)$ is a polynomial of degree at most $n-1$. Namely, there are constants $a_0$, $a_1$, $\ldots$, $a_{n-1}$ such that
\[q(x)=a_{n-1}x^{n-1}+\cdots+a_1 x+a_0.\]
This implies that
\[p(x)=a_nx^n+a_{n-1}x^{n-1}+\cdots+a_1x+a_0,\]where $a_n=C/n!$. Hence, $p(x)$ is a polynomial of degree at most $n$.
\end{myproof}

Mean value theorem can be used to estimate the magnitude of a function provided that we know the derivative.
\begin{example}{}
Given that the function $f:[0,10]\to\mathbb{R}$ is continuous on $[0,10]$, differentiable on $(0,10)$, and $-3<f'(x)<8$ for all $x$ in $ (0,10)$. If $f(0)=-2$, find a range for the values of $f(x)$.
\end{example}
\begin{solution}{Solution}
Let $x$ be point in  $(0, 10]$. By mean value theorem, there is a $c\in (0, x)$ such that
\[\frac{f(x)-f(0)}{x-0}=f'(c).\]
Since $-3<f'(c)<8$, we find that
\[-3x<f(x)+2<8x.\]
This implies that
\[ -32<-3x-2<f(x)<8x-2<78.\]Therefore, a range for the values of $f(x)$ is $(-32, 78)$.
\end{solution}

The next example shows that the mean value theorem can be used to determine the number of solutions of an equation.
\begin{example}{}
Recall that in Example \ref{ex230215_1}, we have shown that the equation 
\[x^6+6x+1=0\]has a real root. Determine the exact number of real roots of this equation.

\end{example}
\begin{solution}{Solution}
Let $f:\mathbb{R}\to\mathbb{R}$ be the function $f(x)=x^6+6x+1$. This is a  differentiable function with\[f'(x)=6x^5+6.\]From this, we find that $f'(x)=0$ if only if $x^5=-1$, if and only if $x=-1$.   

If $x_1$ and $x_2$ are two points such that $x_1<x_2$ and $f(x_1)=f(x_2)=0$,  Rolle's theorem says that there is a point  $u$ in $(x_1, x_2)$  such that $f'(u)=0$. 

\bs
If $f(x)=0$ has three distinct real roots, we can assume that these real roots are $x_1$, $x_2$ and $x_3$ with $x_1<x_2<x_3$. Then there is a $u_1$ in $(x_1, x_2)$, and a $u_2$ in $(x_2, x_3)$ such that  
$f'(u_1)=f'(u_2)=0$. In other words, $f'(x)=0$ has two distinct real  roots $u_1$ and $u_2$. But we have shown that there is only one $x$ such that $f'(x)=0$. Therefore, $f(x)=0$ can have at most two real solutions.

Since $f(0)=1$, we have $f(-1)<0<f(1)$. By intermediate value theorem, there is a $c_1\in (-1, 0)$ such that $f(c_1)=0$.

Since $f(-2)=53$, we have $f(-1)<0<f(-2)$. By intermediate value theorem, there is a $c_2\in (-2, -1)$ such that $f(c_2)=0$. 

We conclude that $f(x)=0$ has exactly two real roots.

\end{solution}
Another important application of the mean value theorem is to determine the increasing or decreasing patterns of functions.

\begin{theorem}[label=thm230215_4]{}
Given that $f:[a,b]\to\mathbb{R}$ is a function continuous on $[a,b]$, and differentiable on $(a, b)$.
\begin{enumerate}[1.]
\item If $f'(x)>0$ for all $x\in (a, b)$, then $f:[a,b]\to\mathbb{R}$ is a strictly increasing function.
\item If $f'(x)<0$ for all $x\in (a, b)$, then $f:[a,b]\to\mathbb{R}$ is a strictly decreasing function.
\end{enumerate}
\end{theorem}
Notice that we only assume that $f'$ is positive or negative on the open interval $(a, b)$. If $f'$ exists at the end points, it can be 0 there, and the conclusion about the strict monotonicity still holds for the entire closed interval $[a,b]$.
\begin{myproof}{Proof}
It suffices for us to prove the first statement. Given any two points $x_1$ and $x_2$ in the closed interval $[a,b]$ with $x_1<x_2$, the function $f$ is continuous on $[x_1, x_2]$, differentiable on $(x_1, x_2)$, and $f'(x)>0$ for any $x\in (x_1, x_2)$. By mean value theorem, there is a point $c$ in $(x_1, x_2)$ such that \bp
\[\frac{f(x_2)-f(x_1)}{x_2-x_1}=f'(c).\]  Since $f'(c)>0$ and $x_2-x_1>0$, we conclude that
\[f(x_2)>f(x_1).\]
This proves that $f:[a,b]\to\mathbb{R}$ is  strictly increasing.
\end{myproof}

We look at a simple example.
\begin{example}[label=ex230216_8]{}
Consider the function $f:\mathbb{R}\to\mathbb{R}$, $f(x)=x^3$. Notice that $f$ is differentiable and $f'(x)=3x^2$. Hence, $f'(x)>0$ for $x\neq 0$, but $f'(0)=0$. Therefore, we cannot apply Theorem \ref{thm230215_4} directly to conclude that $f:\mathbb{R}\to\mathbb{R}$, $f(x)=x^3$ is a strictly increasing function. However, we can proceed in the following way. Since $f'(x)>0$ on the open interval $(-\infty, 0)$, Theorem \ref{thm230215_4} implies that $f$ is strictly increasing on the closed interval $(-\infty, 0]$. Since  $f'(x)>0$ on the open interval $(0, \infty)$, Theorem \ref{thm230215_4} again implies that $f$ is strictly increasing on the closed interval $[0, \infty)$. Combining together, we conclude that  $f:\mathbb{R}\to\mathbb{R}$, $f(x)=x^3$ is strictly increasing.
\end{example}

\begin{remark}{}
Let   $f:[a, b]\to\mathbb{R}$ be a function defined on $[a, b]$, and let $x_1, \ldots, x_n$ be points in $(a, b)$ such that the following conditions are satisfied.
\begin{enumerate}[(i)]
\item $f$ is  continuous on $[a,b]$, differentiable on $(a, b)$.
\item   $f'(x_k)=0$ for $1\leq k\leq n$.
\item $f'(x)>0$ for any $x\in (a, b)\setminus\{x_1, \ldots, x_n\}$.
\end{enumerate}
Using the same reasoning as in Example \ref{ex230216_8},  one can prove that $f$ is strictly increasing on $[a, b]$.
\end{remark}

Example \ref{ex230216_8} shows that if a function $f:(a, b)\to\mathbb{R}$ is differentiable and strictly increasing, it is not necessary that $f'(x)>0$ for all $x\in (a, b)$. If we relax the strict monotonicity to  monotonicity,  we will find that $f'(x)\geq 0$ for all $x\in (a, b)$ is sufficient and necessary for $f$ to be increasing.

\begin{theorem}[label=thm230216_9]{}
Given that the function $f:[a,b]\to\mathbb{R}$ is  continuous on $[a,b]$, and differentiable on $(a, b)$.
\begin{enumerate}[1.]
\item $f:[a,b]\to\mathbb{R}$ is an  increasing function if and only if   $f'(x)\geq 0$ for all $x\in (a, b)$.
\item $f:[a,b]\to\mathbb{R}$ is a decreasing function if and only if   $f'(x)\leq 0$ for all $x\in (a, b)$.
\end{enumerate}
\end{theorem}
\begin{myproof}{Proof} Again, let us consider the first statement. 
If $f'(x)\geq 0$ for all $x\in (a,b)$,  the proof that $f$ is increasing is almost verbatim the proof in Theorem \ref{thm230215_4}, with $>$ replaced by $\geq $. For the converse, if $f$ is increasing on $[a, b]$, we want to show that $f'(x_0)\geq 0$ for any $x_0$ in $(a, b)$. This follows from the fact that
\[\frac{f(x)-f(x_0)}{x-x_0}\geq 0\] for any $x$ in $(a, b)\setminus\{0\}$ since $f$ is increasing. Taking limit gives $f'(x_0)\geq 0$.
\end{myproof}

For a   function $f(x)$ that is differentiable,  the condition $f'(x_0)=0$ is   necessary  for an interior point $x_0$ to be a local extremizer, but not sufficient.
Theorem  \ref{thm230215_4} provides the tool  for determining whether such point is a local extremizer. It is called the \emph{first derivative test}. We would not go into the general formulation. Instead, we will apply Theorem  \ref{thm230215_4} or Theorem \ref{thm230216_9} directly to solve such problems.

\begin{example}[label=ex230216_1]{}
Consider the function $f:\mathbb{R}\to\mathbb{R}$ defined by
\[f(x)=\frac{x}{x^2+1}.\]
Find the local maximum value and the local minimum value of $f$, and find the range of the function $f$. 
 
\end{example}
\begin{solution}{Solution}
Since $f$ is a rational function, it is differentiable, and
\[f'(x)=\frac{(x^2+1)-x(2x)}{(x^2+1)^2}=\frac{1-x^2}{(x^2+1)^2}=-\frac{(x+1)(x-1)}{(x^2+1)^2}.\]Since $f$ is differentiable everywhere, the only candidates for the local maximizer and the local minimizer are those points $x$ where $f'(x)=0$, which are the poins $x=-1$ and $x=1$.
\begin{enumerate}[$\bullet$\;\;]
\item
 When $x\in (-\infty, -1)$, $f'(x)<0$,  and so $f$ is   decreasing on $(-\infty, 1]$. 
\item When $x\in (-1, 1)$, $f'(x)>0$,  and so $f$ is   increasing on $[-1, 1]$. 
\item When $x\in (1, \infty)$, $f'(x)<0$,  and so $f$ is   decreasing on $[1, \infty)$. 
\end{enumerate}These imply that $x=-1$ is a local minimizer, and   $x=1$ is a local maximizer. 
The local maximum value of $f$ is $f(1)=\frac{1}{2}$, and the local minimum value is $f(-1)=-\frac{1}{2}$.
Notice that
\[\lim_{x\to -\infty}f(x)=\lim_{x\to\infty} f(x)=0.\]Since $f$ is   decreasing on $(-\infty, -1]$,   for any $x$ in $(-\infty, -1]$, \[-\frac{1}{2}=f(-1)\leq f(x)<0.\]
Since $f$ is   increasing on $[-1, 1]$, for any $x$ in $[-1, 1]$, 
\[-\frac{1}{2}=f(-1)\leq f(x)\leq f(1)=\frac{1}{2}.\]
\bs Since $f$ is decreasing   on $[1, \infty)$,  for any $x$ in $[1, \infty)$, \[0<f(x)\leq f(1)=\frac{1}{2}.\]
Combining together, we conclude that the range of $f$ is $[-\frac{1}{2}, \frac{1}{2}]$.
 
\end{solution}

 \begin{figure}[ht]
\centering
\includegraphics[scale=0.2]{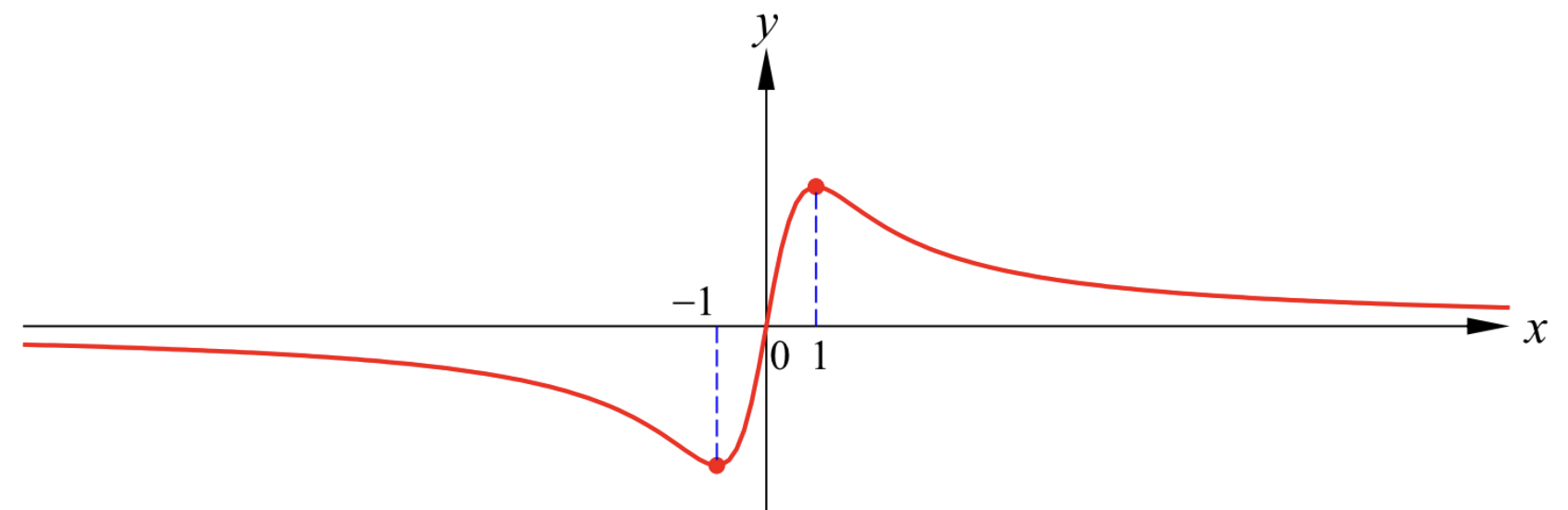}
\caption{  The function $\di f(x)=\frac{x}{x^2+1}$.\fa}\label{figure25}
\end{figure}

There is also a \emph{second derivative test} for determining whether a stationary point is a local minimizer or a local maximizer.

\begin{theorem}[label=thm230215_6]{Second Derivative Test}
Let $(a,b)$ be an interval that contains the point $x_0$, and let $f:(a, b)\to\mathbb{R}$ be a  differentiable function. Assume that $f'(x_0)=0$, and $f''(x_0)$ exists.
\begin{enumerate}[1.]
\item
If $f''(x_0)>0$, then $x_0$ is a local minimizer of $f$.
\item If $f''(x_0)<0$, then $x_0$ is a local maximizer of $f$.
\end{enumerate}
\end{theorem}

\begin{highlight}{}The second derivative test is inconclusive  if $f''(x_0)=0$, as can be shown by considering the three functions $f_1(x)=x^4$, $f_2(x)=-x^4$ and $f_3(x)=x^3$. All these three functions have $x=0$ as  a stationary point. Their second derivatives are all equal to zero at $x=0$. However, $x=0$ is a local minimizer of $f_1(x)=x^4$, it is a local maximizer of the function $f_2(x)=-x^4$, and it is not a local extremizer for the function $f_3(x)=x^3$.\end{highlight}

\begin{myproof}{\linkt Proof of Theorem \ref{thm230215_6}}
We will give a proof of the first statement. The proof of the second statement is similar.
For the first statement, we are given that $f'(x_0)=0$ and $f''(x_0)>0$. By definition, 
\[f''(x_0)=\lim_{x\to x_0}\frac{f'(x)-f'(x_0)}{x-x_0}=\lim_{x\to x_0}\frac{f'(x)}{x-x_0}.\]
Take $\varepsilon$ to be the positive number $f''(x_0)/2$. The definition of limit implies that there is a $\delta>0$ such that $(x_0-\delta, x_0+\delta)\subset (a, b)$, and for all the points $x$ in $(x_0-\delta, x_0)\cup (x_0, x_0+\delta)$, 
\[\left|\frac{f'(x)}{x-x_0}-f''(x_0)\right|<\frac{f''(x_0)}{2}.\]
This implies that for all $x\in (x_0-\delta, x_0)\cup (x_0, x_0+\delta)$,
\begin{equation}\label{eq230215_7}\frac{f'(x)}{x-x_0}>f''(x_0)-\frac{f''(x_0)}{2}=\frac{f''(x_0)}{2}>0.\end{equation}

If $x\in (x_0-\delta, x_0)$, $x-x_0<0$. Equation \eqref{eq230215_7} implies that $f'(x)<0$. Therefore, $f$ is   decreasing on $(x_0-\delta, x_0]$. This implies that
\[f(x)\geq f(x_0)\hspace{1cm} \text{for all}\;x\in (x_0-\delta, x_0).\]
If $x\in (x_0, x_0+\delta)$, $x-x_0>0$. Equation \eqref{eq230215_7} implies that $f'(x)>0$. Therefore, $f$ is   increasing on $[x_0, x_0+\delta)$. This implies that
\[f(x)\geq f(x_0)\hspace{1cm} \text{for all}\;x\in (x_0, x_0+\delta).\]
Combining together, we find that $f(x)\geq f(x_0)$ for all $x$ in $(x_0-\delta, x_0+\delta)$. This proves that $x_0$ is a local minimizer of $f$.

\end{myproof}
 
\begin{example}{}
For the function $\di f(x)$ considered in Example \ref{ex230216_1}, we have shown that the stationary points are $x=-1$ and $x=1$. A tedious computation gives
\[f''(x)=\frac{2x(x^2-3)}{(x^2+1)^3}.\]
Hence,
\[f''(-1)=\frac{1}{2},\hspace{1cm}f''(1)=-\frac{1}{2}.\]The second derivative test can then be used to conclude that
$x=-1$ is a local minimizer, and $x=1$ is a local maximizer.
\end{example}
Although applying the second derivative test seems straightforward, an analysis using the first derivative test is more conclusive. Finding the second derivative can also be tedious, as shown in the example above.

At the end of this section, we want to prove an analogue of intermediate value theorem for derivatives.
\begin{theorem}{Darboux's Theorem}
Let $f:[a,b]\to\mathbb{R}$ be a differentiable function. If $w$ is a value strictly between $f_+'(a)$ and $f_-'(b)$, then there is a point $c$ in $(a, b)$ such that
$f'(c)=w$.
\end{theorem}If the function $g'$ is continuous, then Darboux's theorem follows immediately from the intermediate value theorem. The strength of Darboux's theorem lies in the fact that it does not     assume the continuity of $g'$. 

\begin{myproof}{Proof} The proof   is an again an application of the extreme value theorem.

 Without loss of generality,   assume that $f_+'(a)<w<f_-'(b)$. 
Since $f:[a,b]\to\mathbb{R}$ is a differentiable function, it is continuous. 
Define the function $g:[a,b]\to\mathbb{R}$ by
\[g(x)=f(x)-wx.\]\bp
Then $g$ is differentiable and
\[g'(x)=f'(x)-w.\]
Notice that $g:[a,b]\to\mathbb{R}$ is also continuous. By extreme value theorem, $g:[a,b]\to\mathbb{R}$ has a  minimum value. Now,
\[g'_+(a)=f_+'(a)-w<0,\hspace{1cm} g_-'(b)=f_-'(b)-w>0.\]
 By definition,
\[g_+'(a)=\lim_{x\to a^+}\frac{g(x)-g(a)}{x-a}.\] Taking $\varepsilon$ to be the positive number $-g_+'(a)/2$, we find that there is a $\delta>0$ such that $\delta\leq b-a$, and for all $x\in (a, a+\delta)$,  
\[\frac{g(x)-g(a)}{x-a}<g_+'(a)+\varepsilon= \frac{g_+'(a)}{2}<0.\]In particular, for all $x\in (a, a+\delta)$, $g(x)<g(a)$, and thus $g(a)$ is not a minimum value of the function $g$. Similarly, since $g_-'(b)>0$, we find that $g(b)$ is not a minimum value of the function $g$. In other words, the minimizer of $g$ must be a point $c$ inside $(a, b)$. This is then also a local minimizer. Since $g$ is differentiable, we must have $g'(c)=0$. This implies that $f'(c)=w$.

\end{myproof}

\begin{remark}[label=remark230216_1]{}
As a consequence of the Darboux's theorem, we find that if a function $f:(a,b)\to\mathbb{R}$ is differentiable and $f'(x)\neq 0$ for any $x\in (a, b)$, then either $f'(x)>0$ for all $x\in (a, b)$, or $f'(x)<0$ for all $x\in (a, b)$. In any case, this means that such  a function must be strictly monotonic. 
\end{remark}

Before closing this section, let us define a terminology.
\begin{definition}{$\pmb{C^k}$ functions}
Let $k$ be a nonnegative integer.
We say that a function $f:(a, b)\to\mathbb{R}$ is a $C^k$-function it is has $k$ times derivatives and the $k^{\text{th}}$-derivative $f^{(k)}:(a,b)\to\mathbb{R}$ is also continuous. 
\end{definition}
It s easy to see that if  $f:(a, b)\to\mathbb{R}$ is a $C^k$-function, then for any $0\leq j<k$, $f^{(j)}:(a,b)\to\mathbb{R}$ is   continuous. 

A $C^0$ function is just a continuous function.
A $C^1$ function is   called a continuously differentiable function.
In general, a $C^k$ function is called a $k$-times continuously differentiable function.

The definition of $C^k$ functions can be extended to the case where the function $f$ is defined on a closed interval $[a, b]$.
\vp
\noindent
{\bf \large Exercises  \thesection}
\setcounter{myquestion}{1}

\begin{question}{\themyquestion}
Given that the function $f:[-5,8]\to\mathbb{R}$ is continuous on $[-5,8]$, differentiable on $(-5,8)$, and $-4<f'(x)<4$ for all $x$ in $ (-5,8)$. If $f(0)=2$, find a range for the values of $f(x)$.
\end{question}

\atc
\begin{question}{\themyquestion}
Show that the function $f:\mathbb{R}\to\mathbb{R}$, 
\[f(x)=\frac{x^3}{x^2+1}\]
is strictly increasing, and find the range of the function.

\end{question}
\atc
\begin{question}{\themyquestion}
Show that the equation
\[ x^5+x+32=0\] has exactly one real solution.
\end{question}
\atc
\begin{question}{\themyquestion}
Find the number of real solutions of the equation
\[ \frac{32x}{x^4+16}=1.\]
\end{question}
\atc
\begin{question}{\themyquestion}
Let $n$ be a nonnegative integer, and let $f:(a, b)\to\mathbb{R}$ be a differentiable function.  If the equation $f'(x)=0$ has   $n$ distinct real roots in the interval $(a,b)$, show that the equation $f(x)=0$ has at most $(n+1)$ distinct real roots in the interval $(a, b)$.
\end{question}
 
\atc
\begin{question}{\themyquestion}
Consider the function $f:\mathbb{R}\to\mathbb{R}$ defined by
\[f(x)=\frac{x+1}{x^2+15}.\]
\begin{enumerate}[(a)]
\item 
Find the local maximum value and the local minimum value of $f$.
\item Find the range of the function $f$. 
\end{enumerate}
\end{question}

\atc
\begin{question}{\themyquestion}
Let $f:[a, b]\to \mathbb{R}$ be a function such that the limit 
\[L=\lim_{x\to b^-}\frac{f(x)-f(b)}{x-b} \]
 exists. If $L>0$, show that there is a $\delta>0$ such that $\delta\leq b-a$ and for all $x\in (b-\delta, b)$,
\[f(x)<f(b).\]

\end{question}

\atc
\begin{question}{\themyquestion}
Let $f:(a, b)\to \mathbb{R}$ be a differentiable function. Suppose that $f':(a,b)\to\mathbb{R}$ is monotonic, show that $f':(a,b)\to\mathbb{R}$ is continuous.

\end{question}
\vp
\section{The Cauchy Mean Value Theorem}\label{sec3.4}
In previous section, we have seen that the mean value theorem is very useful in analysing the behavior of a differentiable function. For future applications, we will often quote it in the following form.
\begin{highlight}{Alternative Form of Mean Value Theorem} If $f:(a, b)\to \mathbb{R}$ is a differentiable function, $x_0$ is a point in $(a, b)$, $h$ is such that $x_0+h$ is also in $(a, b)$, then there is a number $c\in (0,1)$ such that 
\begin{equation}\label{eq230216_3}f(x_0+h)-f(x_0)=f'(x_0+ch)h.\end{equation}
\end{highlight}
To see this, let $x_1=x_0+h$.  If $h=0$, \eqref{eq230216_3} is obviously true for any $c$ in $(0, 1)$. If $h\neq 0$, then when $c$ runs through all values from 0 to 1, $x_0+ch$ runs through all points in the open interval  $I$ with $x_0$ and $x_1$ as endpoints. Thus, \eqref{eq230216_3} says that
\[ f(x_1)-f(x_0)=f'(u)(x_1-x_0)\] for some $u$ in the open intefval $I$, which is precisely the statement of the mean value theorem.

When finding limits of functions, we often encounter situations like
\[\lim_{x\to x_0}\frac{f(x)}{g(x)}\] where both $\di\lim_{x\to x_0}f(x)$ and $\di\lim_{x\to x_0}g(x)$ are   zero. For example, let  \[f(x)=x^{20}+2x^9-3\hspace{1cm}\text{and} \hspace{1cm} g(x)=x^7-1.\] Then \[\lim_{x\to 1}f(x)=f(1)=0\hspace{1cm}\text{and} \hspace{1cm}\lim_{x\to 1}g(x)=g(1)=0.\] Hence, we cannot apply limit quotient law to evaluate 
\[\lim_{x\to 1}\frac{f(x)}{g(x)}=\lim_{x\to 1}\frac{x^{20}+2x^9-3}{x^7-1}.\]
Observe that 
\[\lim_{x\to 1}\frac{f(x)}{g(x)}=\lim_{h\to 0}\frac{f(1+h)}{g(1+h)}.\] Since we are only interested in the limit when $x$ approaches 1, we are prone to use the mean value theorem in the form \eqref{eq230216_3} and conclude that there are $c_1 $ and $c_2 $ in $(0,1)$ such that
\begin{equation}\label{eq230217_1}\lim_{x\to 1}\frac{f(x)}{g(x)}=\lim_{h\to 0}\frac{f(1+h)-f(1)}{g(1+h)-g(1)} =\lim_{h\to 0}\frac{f'(1+c_1h)}{g'(1+c_2h)}.\end{equation}
For the functions $f$ and $g$ that we consider above, $f'$ and $g'$ are both continuous at $x=1$ and $g'(1)\neq 0$. Hence, we find that
\[\lim_{h\to 0}\frac{f'(1+c_1h)}{g'(1+c_2h)}=\frac{f'(1)}{g'(1)}.\]  For general differentiable functions $f(x)$ and $g(x)$ with $f(1)=g(1)=0$, if we do not assume that $f'$ and $g'$ are continuous, we cannot conclude the limit from \eqref{eq230217_1}  
 since $c_1$ and $c_2$ are in general different functions of $h$.

In this section, we are going to prove a generalization of the mean value theorem, called the Cauchy mean value theorem,  which ensures that we can have the same value for $c_1$ and $c_2$.

\begin{theorem}[label=thm230216_2]{Cauchy Mean Value Theorem}
Let $f:[a,b]\to\mathbb{R}$ and $g:[a,b]\to\mathbb{R}$ be two functions that satisfy the following conditions.
\begin{enumerate}[(i)]
\item $f:[a,b]\to\mathbb{R}$ and $g:[a,b]\to\mathbb{R}$  are continuous.
\item $f:(a, b)\to\mathbb{R}$ and $g:(a, b)\to\mathbb{R}$ are differentiable.
\item $g'(x)\neq 0$ for all $x\in (a, b)$. 
 
\end{enumerate}
Then there is a point $x_0$ in $(a, b)$ such that \[\frac{f'(x_0)}{g'(x_0)}=\frac{f(b)-f(a)}{g(b)-g(a)}.\]
\end{theorem}
Notice that when $g(x)=x$,   we have the Lagrange's mean value theorem.
\begin{myproof}{Proof}
The proof uses the same idea as the proof of Lagrange's mean value theorem, with the function $g(x)=x$ replaced by a general $g(x)$.
By Remark \ref{remark230216_1}, the condition $g'(x)\neq 0$ for all $x\in (a, b)$ implies that $g$ is strictly monotonic. Hence, $g(a)\neq g(b)$.

 Define the function $h:[a,b]\to\mathbb{R}$ by
\[h(x)=f(x)-mg(x),\] where the number $m$ is determined by $h(a)=h(b)$. This means
\[f(a)-mg(a)=f(b)-mg(b),\] which gives
\[m=\frac{f(b)-f(a)}{b-a}.\]
Again, the function $h:[a,b]\to\mathbb{R}$ is continuous on $[a,b]$, differentiable on $(a, b)$, and satisfies $h(a)=h(b)$. By Rolle's theorem, there is a point $x_0$ in $(a, b)$ such that
$h'(x_0)=0$. For this $x_0$,
\[f'(x_0)-mg'(x_0)=0.\] Since $g'(x_0)\neq 0$ by assumption, we find that
\[\frac{f'(x_0)}{g'(x_0)}=m=\frac{f(b)-f(a)}{b-a}.\]

\end{myproof}

\begin{example}[label=ex230216_10]{}
Consider the functions $f:[1, 7]\to\mathbb{R}$, $f(x)=x^2$ and $g:[1,7]\to \mathbb{R}$, $g(x)=x^3-9x^2$. 
By Lagrange's mean value theorem, there are points $c_1$ and $c_2$ in $(1, 7)$ such that
\[2c_1=f'(c_1)=\frac{f(7)-f(1)}{7-1}=8,\]
and
\[3c_2^2-18c_2=g'(c_2)=\frac{g(7)-g(1)}{7-1}=-15.\]\end{example}

\begin{example2}{}
Solving for $c_1$ and $c_2$, we have $c_1=4$ and  and $c_2=5$. 

By Cauchy mean value theorem, there is a point $c$ in $(1, 7)$ such that
\[\frac{2c}{3c^2-18c}=\frac{f'(c)}{g'(c)}=\frac{f(7)-f(1)}{g(7)-g(1)}=-\frac{8}{15}.\]Solving this equation gives
\[c=\frac{19}{4}.\]
\end{example2}

An important application of the Cauchy mean value theorem is the following.
\begin{theorem}[label=thm230216_11]{}
Let $n$ be a positive integer, and let $(a, b)$ be an open interval that contains the point $x_0$. If the function $f:(a, b)\to\mathbb{R}$ is $n$ times differentiable, and
\[f(x_0)=f'(x_0)=\cdots=f^{(n-1)}(x_0)=0,\]
then for any $x$ in $ (a, b)$, there is a $c\in (0,1)$ such that
\begin{equation}\label{eq230216_12}f(x)=\frac{h^n}{n!}f^{(n)}(x_0+ch),\hspace{1cm}\text{where}\;h=x-x_0.\end{equation}
\end{theorem}
\begin{myproof}{Proof}
We apply the Cauchy mean value theorem $n$ times to the given  function $f:(a, b)\to\mathbb{R}$ and the function $g:(a,b)\to\mathbb{R}$ defined by $g(x)=(x-x_0)^n$.
Notice that $g$ is also $n$ times differentiable,
\[g(x_0)=g'(x_0)=\cdots=g^{(n-1)}(x_0)=0,\] and 
\[g^{(n)}(x)=n!\hspace{1cm}\text{for all}\;x\in (a, b).\]
Since $f(x_0)=0$, eq. \eqref{eq230216_12} obviously holds for $x=x_0$ with any $c\in (0,1)$. So we only need to consider a point $x =x_1$ in $(a, b)\setminus\{ x_0\}$. First assume that $ x_1>x_0$. For any $1\leq k\leq n$, $g^{(k)}(x)\neq 0$ for any $x\in (x_0, x_1)$. Thus we can apply the Cauchy mean value theorem for the pairs $(f, g)$, $(f', g')$, $\ldots$, $(f^{(n-1)}, g^{(n-1)})$ over the interval $(x_0, x_1)$. \bp
Since $f(x_0)=g(x_0)=0$, Cauchy mean value theorem implies that there exists a point $u_1$ in $(x_0, x_1)$ such that
\[\frac{f(x_1)}{g(x_1)}=\frac{f(x_1)-f(x_0)}{g(x_1)-g(x_0)}=\frac{f'(u_1)}{g'(u_1)}.\]If $n=1$, we are done. If $n\geq 2$, then $f'(x_0)=g'(x_0)=0$. Apply Cauchy mean value theorem again, we find that there is a $u_2$ in $(x_0, u_1)$ such that
\[\frac{f(x_1)}{g(x_1)}=\frac{f'(u_1)}{g'(u_1)}=\frac{f'(u_1)-f'(x_0)}{g'(u_1)-g'(x_0)}=\frac{f''(u_2)}{g''(u_2)}.\]
Continue with this $n$ times, we find that there are points $u_1, \ldots, u_n$ such that $x_0<u_n<u_{n-1}<\cdots<u_1<x_1$, and
\begin{equation}\label{eq230216_13}\frac{f(x_1)}{g(x_1)}=\frac{f'(u_1)}{g'(u_1)}=\cdots=\frac{f^{(n)}(u_n)}{g^{(n)}(u_n)}.\end{equation}
Since $u_n\in (x_0, x_1)$, there is $c\in (0,1)$ such that $u_n=x_0+ch$, where $h=x_1-x_0$. Eq. \eqref{eq230216_13} then  implies that
\[f(x_1)=\frac{h^n}{n!} f^{(n)}(x_0+ch),\hspace{1cm}\text{where}\;h=x_1-x_0.\]
This completes the proof if $x_1>x_0$. The proof for $x_1<x_0$ is similar.
\end{myproof}

Let us look at a  classical example.
\begin{example}[label=ex230216_17]{}
Let $x_0$ be a point in the interval $(a, b)$, and assume that the function $f:(a,b)\to \mathbb{R}$ is twice continuously differentiable. Prove that
\[\lim_{h\to 0}\frac{f(x_0+h)+f(x_0-h)-2f(x_0)}{h^2}=f''(x_0).\]
\end{example}
\begin{solution}{Solution}
Let $r=\min\{x_0-a, b-x_0\}$. Then $r>0$ and $(x_0-r,x_0+r)\subset (a, b)$. Define the function $g:(-r,r)\to \mathbb{R}$   by
\[g(h)=f(x_0+h)+f(x_0-h)-2f(x_0).\]
Then $g$ is twice continuously diferentiable, and
\[g'(h)=f'(x_0+h)-f'(x_0-h),\hspace{1cm} 
g''(h)=f''(x_0+h)+f''(x_0-h).\]
It is easy to check that
\[g(0)=g'(0)=0.\]By Theorem \ref{thm230216_11}, for any $h\in (-r,r)$, there is a $c(h)\in (0,1)$ such that
\[g(h)=\frac{h^2}{2}g''(c(h)h).\]
Hence,
\[\frac{g(h)}{h^2}=\frac{f''(x_0+c(h)h)+f''(x_0-c(h)h)}{2}.\]
Now since $c(h)\in (0,1)$,
\[|c(h)h|\leq |h|.\]
Therefore, 
\[\lim_{h\to 0}c(h)h=0.\]
Since $f''$ is continuous,
\[\lim_{k\to 0}f''(x_0+k)=f''(x_0).\]
By limit law for composite functions, we find that
\[\lim_{h\to 0}\frac{f''(x_0+c(h)h)+f''(x_0-c(h)h)}{2}=\frac{f''(x_0)+f''(x_0)}{2}=f''(x_0).\]
Therefore,
\[\lim_{h\to 0}\frac{f(x_0+h)+f(x_0-h)-2f(x_0)}{h^2}=\lim_{h\to 0}\frac{g(h)}{h^2}=f''(x_0).\]
\end{solution}
\vp
\noindent
{\bf \large Exercises  \thesection}
\setcounter{myquestion}{1}
\begin{question}{\themyquestion}
Given that $p(x)$ is a polynomial of degree at most 5, and
\[p(1)=p^{(1)}(1)=p^{(2)}(1)=p^{(3)}(1)=p^{(4)}(1)=0,\hspace{1cm}p^{(5)}(1)=1200.\]
Find the polynomial $p(x)$.
\end{question}

\atc
\begin{question}{\themyquestion}
Let $(a, b)$ be an interval that contains the point $x_0$. Given that the function $f:(a,b)\to\mathbb{R}$ is three times continuously differentiable, find the limit
\[\lim_{h\to 0}\frac{f(x_0+2h)-2f(x_0+h)+2f(x_0-h)- f(x_0-2h)}{h^3}.\]
\end{question}

\vp
\section{Transcendental Functions}\label{sec3.5}

Up to now we have only dealt with  algebraic functions, which are functions that can be obtained by performing algebraic operations of addition,  multiplication, division and taking roots on polynomials. 
In this section, we introduce other useful elementary functions -- the class of transcendental functions which  includes exponential, logarithmic and trigonometric functions.  These functions have been introduced in  a pre-calculus course, but not rigorously.  

In this section, we are going to define these functions  and derive their properties using calculus. Everything would be done rigorously using the analytic tools that we have developed so far, except for an existence theorem that we are going to prove in Chapter \ref{ch4}.

Let us first state this existence theorem.
\begin{theorem}[label=thm230217_3]{Existence and Uniqueness Theorem}
Let $(a,b)$ be an open interval that contains the point $x_0$, and let $y_0$ be any real number.  Given that $f:(a,b)\to\mathbb{R}$ is a continuous function, there exists a unique differentiable function $F:(a,b)\to\mathbb{R}$ such that \[F'(x)=f(x)\quad \text{for all}\;x\in (a, b), \hspace{1cm}F(x_0)=y_0.\]
\end{theorem}
The function $F(x)$ that satisfies $F'(x)=f(x)$ is called an \emph{antiderative} of $f(x)$.
\begin{definition}{Antiderative}
Let $I$ be an interval. If $f:I\to\mathbb{R}$ and $F:I\to\mathbb{R}$ are functions on $I$ such that 
$F$ is differentiable and 
\[F'(x)=f(x)\hspace{1cm}\text{for all}\;x\in I,\]
then $F(x)$ is called an {\bf antiderivative} of $f(x)$.
\end{definition}
Theorem \ref{thm230217_3} asserts that  a continuous function has an antiderivative. One way to construct an antiderivative is to use integrals, a topic we are going to discuss in Chapter \ref{ch4}.   Theorem \ref{thm230215_3} says that any two  antiderivatives of a given function differ by a constant. The \emph{initial condition} $F(x_0)=y_0$ fixes the constant. Hence, only the existence part of Theorem \ref{thm230217_3} is pending a proof. The uniqueness follows from what we have discussed.

\subsection{The Logarithmic Function }
It is easy to check that for any integer $n$ that is not equal to $-1$, an antiderivative of the function $f(x)=x^n$ is the function
\[F(x)=\frac{x^{n+1}}{n+1}.\]
So far we haven't seen any algebraic function whose antiderivative is equal to $f(x)=\di\frac{1}{x}$.
We define one such function and call it the natural logarithm function.

\begin{definition}{The Natural Logarithm Function}
The natural logarithm function $f:(0, \infty)\to \mathbb{R}$, $f(x)=\ln x$ is defined to be the unique differentiable function satisfying
\[f'(x)=\di\frac{1}{x},\hspace{1cm}f(1)=0.\]
\end{definition}
Since $g:(0, \infty)\to \mathbb{R}$, $g(x)=\di\frac{1}{x}$ is a continuous function,
the existence and uniqueness of the function $f(x)=\ln x$ is guaranteed by Theorem \ref{thm230217_3}.

\begin{highlight}{The Natural Logarithm Function}
By definition,
\[\frac{d}{dx}\ln x=\frac{1}{x},\hspace{1cm} x>0.\]
Since $1/x>0$ for all $x>0$, we find that $f(x)=\ln x$ is a strictly increasing function. Moreover, since $\ln 1=0$,
\begin{enumerate}[$\bullet$\;\;]
\item when $0<x<1$, $\ln x<0$;
\item when $x>1$, $\ln x>0$.
\end{enumerate}
\end{highlight}

\begin{figure}[ht]
\centering
\includegraphics[scale=0.2]{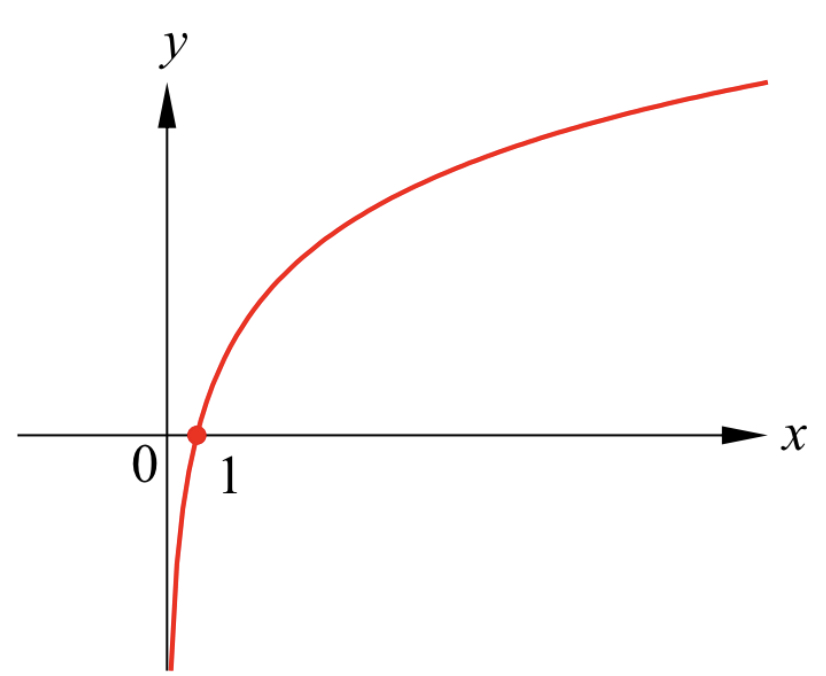}
\caption{  The function $y=\ln x$.\fa}\label{figure26}
\end{figure}
The following gives some useful properties of the natural logarithmic function.
\begin{proposition}[label=prop230217_5]{Properties of the Natural Logarithm Function}
Let $x$ and $y$ be any positive numbers, and let $r$ be a rational number. We have the following.
\begin{enumerate}[(a)]
\item $\ln (xy)=\ln x+\ln y$
\item $\ln\di \frac{x}{y}=\ln x-\ln y$
\item $\ln x^r=r\ln x$

\end{enumerate}
\end{proposition}
In part (c), we require $r$ to be a rational number since we have not defined $x^r$ when $r$ is an irrational number.
\begin{myproof}{Proof}
To prove (a), we fixed $y>0$ and define the function $f:(0,\infty)\to\mathbb{R}$ by
\[f(x)=\ln(xy)-\ln y.\] 
Then $f(1)=\ln y-\ln y=0$, and
\[f'(x)=\frac{y}{xy}=\frac{1}{x}.\]
By the uniquesness asserted in Theorem \ref{thm230217_3} and the definition of the natural logarithm function, we conclude that $f(x)=\ln x$. This proves (a).\bp

To prove (b), we notice that part (a) gives
\[ \ln\left(\frac{x}{y}\right)+\ln y=\ln\left(\frac{x}{y}\times y\right)=\ln x.\] 

For (c), notice that it is obvious if $r=0$. If $r\neq 0$,   define the function $f:(0,\infty)\to\mathbb{R}$ by
\[f(x)= \frac{1}{r}\ln x^r.\]   Then $f(1)=0$, and
\[f'(x)=\frac{1}{r}\times\frac{rx^{r-1}}{x^r}=\frac{1}{x}.\]
This allows us to conclude that $f(x)=\ln x$, and (c) is thus proved.

\end{myproof}

From part (b) of Proposition \ref{prop230217_5}, we find that for any $x>0$,
\[\ln \frac{1}{x}= \ln x^{-1}=-\ln x;\]
and if $n$ is a positive integer, 
\[\ln x^n=n\ln x.\]
In particular, we find that
\[\ln 2^n=n\ln 2,\]
\[\ln \frac{1}{2^n}=-n\ln 2.\]
Since $\ln 2>0$, we conclude the following.
\begin{proposition}{}
$f:(0, \infty)\to\mathbb{R}$, $f(x)=\ln x$ is a strictly increasing function with
\[\lim_{x\rightarrow 0^+}\ln x=-\infty,\hspace{1cm}\lim_{x\to \infty}\ln x=\infty.\]
Hence, the range of $f(x)=\ln x$ is $\mathbb{R}$. 
\end{proposition}

\subsection{The Exponential Functions}
Since the function $f:(0,\infty)\to\mathbb{R}$, $f(x)=\ln x$ is continuous and strictly increasing, its inverse function exists. We define this inverse function as the exponential function $\exp(x)$. The domain of $\exp(x)$ is the range of $\ln x$, which is $\mathbb{R}$. The range of $\exp(x)$ is the domain of $\ln x$, which is $(0, \infty)$.
\begin{definition}{The Natural Exponential Function}
The natural exponential function $\exp :\mathbb{R}\to \mathbb{R}$ is defined to be the inverse of the function $ \ln x$. It satisfies
\[\ln \exp(x)=x \quad\text{for any}\;x\in\mathbb{R},\hspace{1cm}
 \exp(\ln x)=x\quad\text{for any}\;x>0.\]
\end{definition}

\begin{figure}[ht]
\centering
\includegraphics[scale=0.2]{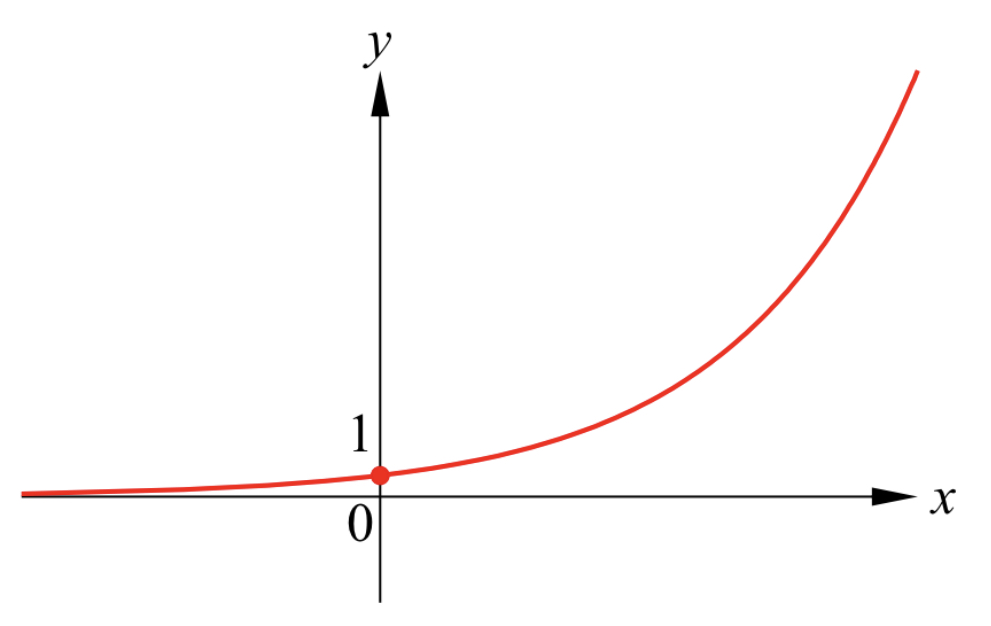}
\caption{  The function $y=\exp(x)$.\fa}\label{figure27}
\end{figure}
We can deduce the following properties.
\begin{proposition}{Properties of the Natural Exponential Function I}
The exponential function $\exp(x)$ is a strictly increasing differentiable function defined on the set of real numbers. It has the following properties.
\begin{enumerate}[(a)]
\item $\exp(x)>0$ for all $x\in\mathbb{R}$ and $\exp(0)=1$.
\item $\di\lim_{x\to-\infty}\exp(x)=0$, $\di\lim_{x\to \infty}\exp(x)=\infty$.
\item $\di\frac{d}{dx}\exp(x)=\exp(x)$.

\end{enumerate}
\end{proposition}
\begin{myproof}{Proof}
(a) and (b) are obvious from the corresponding properties of $\ln x$.
 For part (c), we employ the derivative formula for inverse function. To make it less confusing, let $y=\exp(x)$. Then
\[\ln y=x.\]
Differentiating both sides with respect to $x$, we find that
\[ \frac{1}{y}\frac{dy}{dx}=1.\]Therefore,
\[\frac{d}{dx}\exp(x)=\frac{dy}{dx}=y=\exp(x).\]

\end{myproof}

From the properties of the natural logarithm stated in Proposition \ref{prop230217_5}, we have the following.
\begin{proposition}{Properties of the Natural Exponential Function II}
Let $x$ and $y$ be any real numbers, and let $r$ be a rational number. We have the following.
\begin{enumerate}[(a)]
\item $\exp(x+y)=\exp(x)\exp(y)$.
\item $\exp(x-y)=\di \frac{\exp(x)}{\exp(y)}$.
\item $\exp(x)^r=\exp(rx)$.\end{enumerate}
\end{proposition}
\begin{myproof}{Proof}
Let $u=\exp(x)$ and $v=\exp(y)$. Then $u$ and $v$ are positive numbers and
\[x=\ln u,\hspace{1cm}y=\ln v.\]
By Proposition \ref{prop230217_5},
\[\ln(uv)=\ln u+\ln v=x+y.\]\bp
Therefore,
\[\exp(x)\exp(y)=uv=\exp(x+y).\]
Part (b) is proved in the same way. For part (c),
Proposition \ref{prop230217_5} implies that
\[\ln(u^r)=r\ln u=rx.\]
Therefore,
\[\exp(x)^r=u^r=\exp(rx).\]
\end{myproof}

Notice that part (c) says that for any postive number $u$, and any rational number $r$, 
\[u^r=\exp(r\ln u).\]We can use this to define power functions with    irrational powers.
\begin{definition}{Power Functions}
For any real number $r$, the power function $f(x)=x^r$ is the function defined on $(0, \infty)$ by the formula
\[x^r=\exp(r\ln x).\] When $r>0$, we can extend the definition to the point $x=0$ by definining $f(0)=0$.
\end{definition}

We have seen that this definition coincides with the old definition when $r$ is a rational number. 
For $r>0$, since $\ln x\to -\infty$ as $x\to 0^+$, $r\ln x\to -\infty$ as $x\to 0^+$. Since $\exp(x)\to 0$ as $x\to -\infty$, we conclude that $x^r\to 0$ as $x\to 0^+$. Therefore, the definition $f(0)=0$ makes the function $f(x)=x^r$ continuous.
Using the fact that $\exp(x)$ and $\ln x$ are inverses of each other, we have the following.
\begin{highlight}{}
For any positive number $x$ and any real number $r$,
\[\ln (x^r)=r\ln x.\]
\end{highlight}

\begin{figure}[ht]
\centering
\includegraphics[scale=0.2]{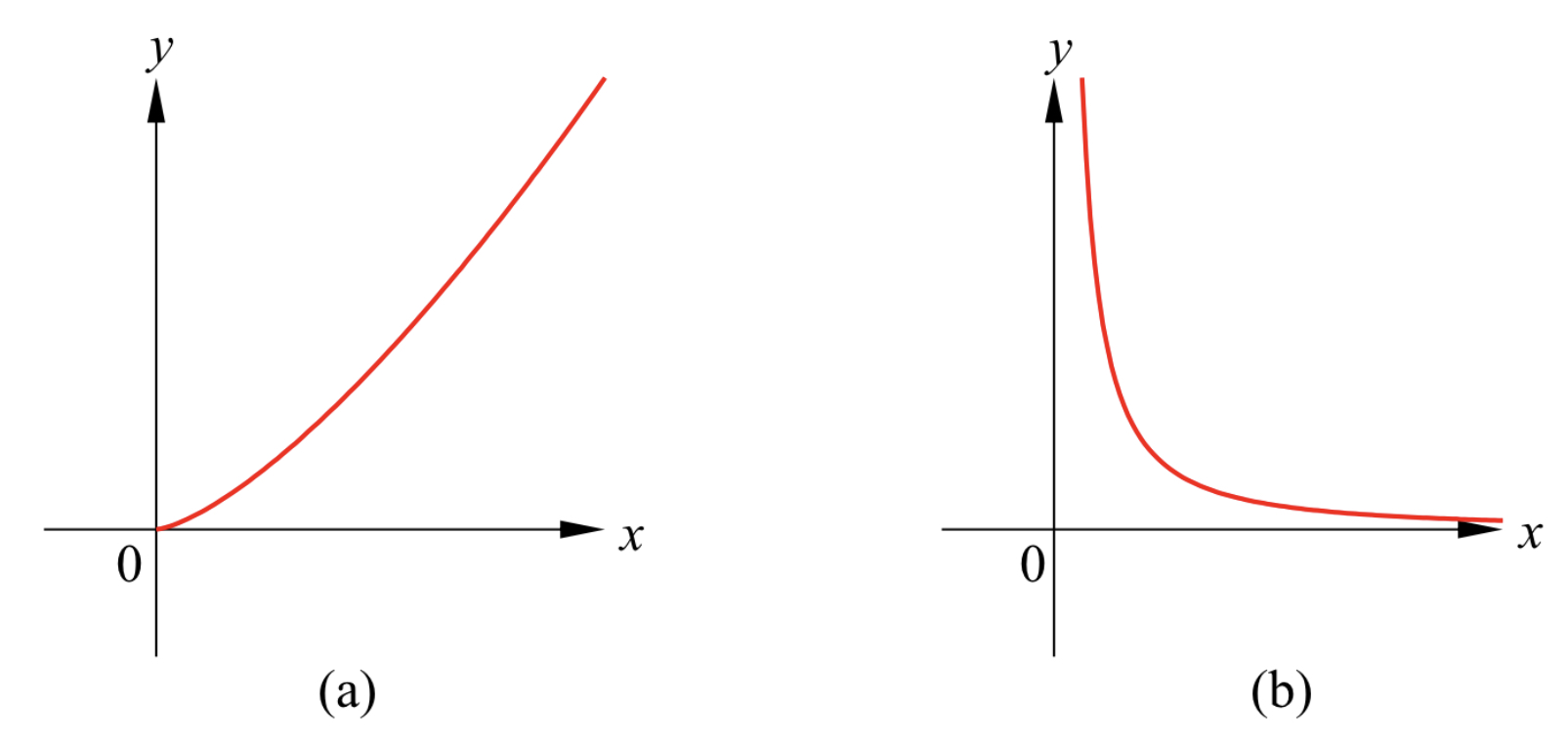}
\caption{ (a) The function $y=x^{\sqrt{2}}$.   (b) The function $y=x^{-\sqrt{2}}$.\fa}\label{figure28}
\end{figure}

Since both $\ln x$ and $\exp(x)$ are strictly increasing functions, it is easy to deduce the following.
\begin{highlight}{Monotonicity of Power Functions}
\begin{enumerate}[1.]
\item When $r>0$, the function $f:[0,\infty)\to\mathbb{R}$, $f(x)=x^r$ is strictly increasing.
\item When $r<0$, the function $f:(0,\infty)\to\mathbb{R}$, $f(x)=x^r$ is strictly decreasing.
\end{enumerate}
\end{highlight}

The following gives the properties of power functions.
\begin{proposition}{}
For any positive numbers $x$ and $y$, and any real numbers $r$ and $s$,
\begin{enumerate}[(a)]
\item $\di (xy)^r=x^ry^r$
\item $\di \left(\frac{x}{y}\right)^r=\frac{x^r}{y^r}$
\item $x^{r+s}=x^rx^s$
\item $\di x^{r-s}=\frac{x^r}{x^s}$
\item $(x^r)^s=x^{rs}$

\end{enumerate}
\end{proposition}
\begin{myproof}{Proof}For part (a), we have
\begin{align*}
(xy)^r&=\exp\left(r\ln(xy)\right)=\exp(r\ln x+r\ln y)\\&=\exp(r\ln x)\exp(r\ln y)=x^ry^r.
\end{align*}Part (b) is proved in the same way. For part (c),
\[x^{r+s}=\exp\left((r+s)\ln x\right)=\exp(r\ln x)\exp(s\ln x)=x^rx^s.\]
Part (d) is proved in the same way.  For part (e),
\[(x^r)^s=\exp\left(s\ln (x^r)\right)=\exp\left(rs\ln x\right)=x^{rs}.\]
\end{myproof}

Using chain rule, we find that $f(x)=x^r$ is a differentiable function.
\begin{proposition}{}
For any real number $r$ and any positive number $x$,
\[\frac{d}{dx}x^r=rx^{r-1}.\]
\end{proposition}
\begin{myproof}{Proof}
This follows from straightforward computation.
\[\frac{d}{dx}x^r=\frac{d}{dx}\exp\left(r\ln x\right)= \exp(r\ln x)\frac{d}{dx}(r\ln x)=x^r\times\frac{r}{x}=rx^{r-1}.\]
\end{myproof}

Now we want to show that $\exp(1)=e$, where $e$ is the number we defined as
\[e=\lim_{n\to \infty}\left(1+\frac{1}{n}\right)^n\] in Chapter \ref{ch1}.

\begin{theorem}[label=thm230218_1]{} We have
\[\exp(1) =\lim_{n\to\infty}\left(1+\frac{1}{n}\right)^n=e.\]This imples that
\[\ln e=1.\]
\end{theorem}
\begin{myproof}{Proof}
We consider the differentiable function $g(x)=\ln (1+x)$, $x>-1$, whose derivative is
\[g'(x)=\frac{1}{1+x}.\]   By definition of derivative, 
\[\lim_{x\to 0}\frac{\ln(1+x)}{x}=\lim_{x\to 0}\frac{g(x)-g(0)}{x-0}=g'(0)=1.\]
Since $\exp(x)$ is a continuous function, we find that
\begin{equation}\label{eq230217_7}\lim_{x\to 0 }\exp\left(\frac{\ln(1+x)}{x}\right)=\exp\left(\lim_{x\to 0 }\frac{\ln(1+x)}{x}\right)= \exp(1).\end{equation} 
Notice that $\{1/n\}$ is a sequence of positive numbers that converges to 0. Therefore, eq. \eqref{eq230217_7} implies that
\[\lim_{n\to \infty}\exp\left(n\ln\left(1+\frac{1}{n}\right)\right)=\exp(1).\] 
By definition,\[\exp\left(n\ln\left(1+\frac{1}{n}\right)\right)=\left(1+\frac{1}{n}\right)^n.\] Thus, we have shown that
\[\exp(1) =\lim_{n\to\infty}\left(1+\frac{1}{n}\right)^n=e.\]
Since $\exp(x)$ and $\ln x$ are inverses of each other, we find that $\ln e=1$.
\end{myproof}
 
\begin{definition}{General Exponential Functions}
Let $a$ be a positive real number such that $a\neq 1$. The exponential function $f(x)=a^x$ is defined by
\[a^x=\exp\left(x\ln a\right),\hspace{1cm}x\in\mathbb{R}.\]
When $a=e$, $f(x)=e^x$ is the natural exponential function 
\[e^x=\exp(x).\]
\end{definition}
Henceforth, we will also use $e^x$ to denote the natural exponential function $\exp(x)$.

The following properties of the general exponential functions can be easily derived from the  corresponding properties of the $\exp(x)$ function.  
\begin{proposition}{}Let $a$ be a positive number.
\begin{enumerate}[1.]
\item When $0<a<1$, $f(x)=a^x$ is a strictly decreasing function.
\item When $a>1$, $f(x)=a^x$ is a strictly increasing function.

\end{enumerate}
\end{proposition}
\begin{proposition}{}
Let $a$ be a positive number such that $a\neq 1$. The function $f(x)=a^x$ is differentiable, and 
\[\frac{d}{dx}a^x=a^x\ln a.\]
\end{proposition}

\begin{proposition}{}Let $a$ be a positive number such that $a\neq 1$. For any real numbers $x$ and $y$, 
\begin{enumerate}[1.]
\item  $a^{x+y}=a^xa^y$
\item $a^{x-y}=\di\frac{a^x}{a^y}$
\item $(a^x)^y=a^{xy}$
\end{enumerate}
\end{proposition}

\subsection{The Trigonometric Functions}
Now we consider the trigonometric functions.  
Recall that an angle is usually measured in degrees, so that the angle of a full circle is $360^{\circ}$. But for analysis, we need to make a change of units to radians.

\begin{figure}[ht]
\centering
\includegraphics[scale=0.2]{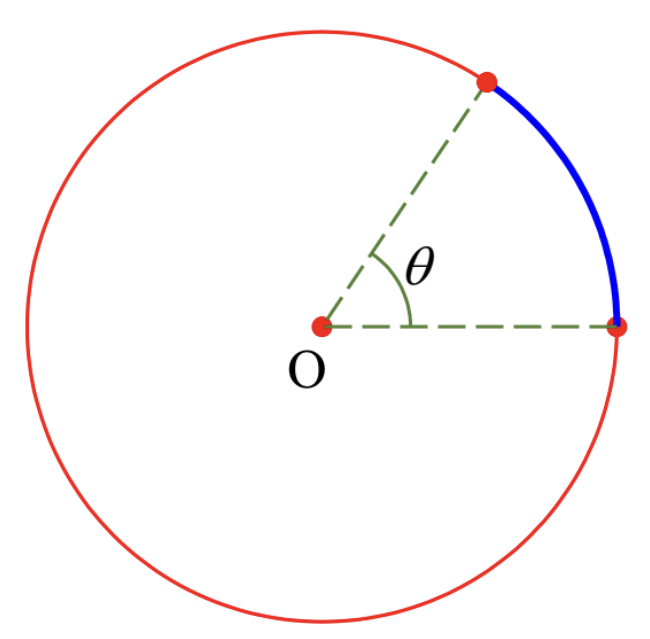}
\caption{An arc with central angle $\theta$.\fa}\label{figure29}
\end{figure}
The number $\pi$ is defined as the ratio of the circumsference of a circle to its diameter. Hence, a circle of radius 1 would have circumsference $2\pi$. This  number $\pi$ can be shown to be an irrational number. The radian measurement of an angle is so   that an arc with central angle $\theta$ radians on a circle of radius $r$ has length $r\theta$, so that the circumsference of the circle is $2\pi r$. Hence, the conversion between degrees and radians is
\[\theta^{\circ} =\frac{\pi}{180}\theta \,\text{rad}.\]

Historically, sine and cosine are defined using right-angled triangles, as shown in Figure \ref{figure30}.
\begin{figure}[ht]
\centering
\includegraphics[scale=0.2]{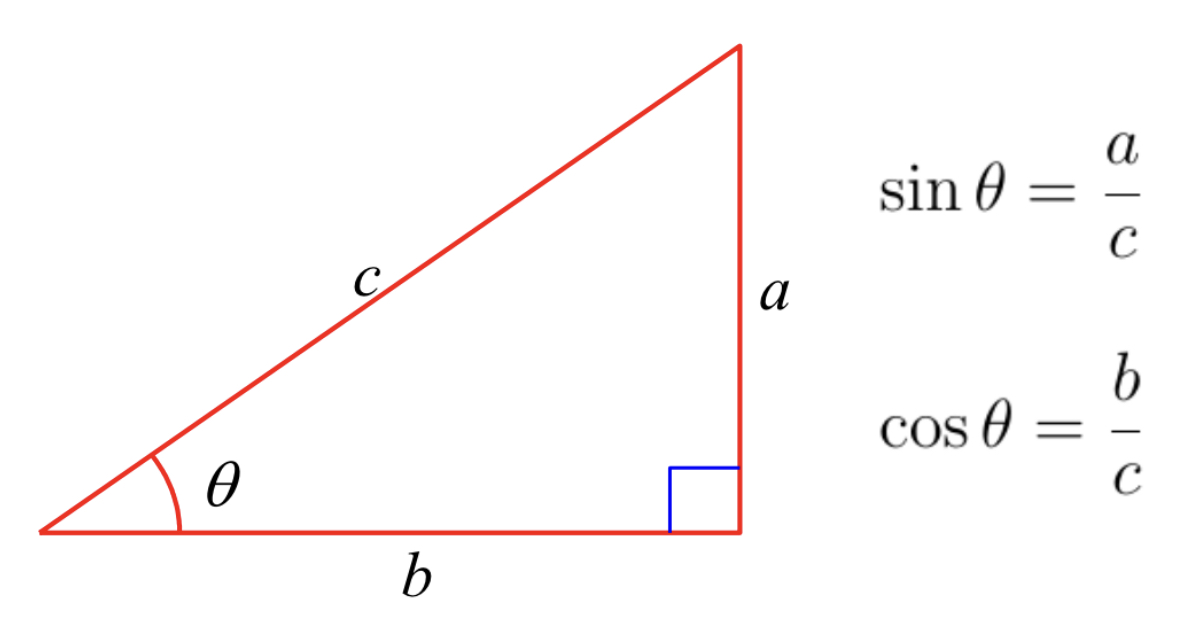}
\caption{Classical definitions of sine and cosine functions.\fa}\label{figure30}
\end{figure}

To extend the definitions of $\sin\theta$ and $\cos\theta$ so that $\theta$ can be any real numbers, we use the unit circle $x^2+y^2=1$. 
The angle measurement starts from the positive $x$-axis and we take the counter-clockwie direction as positive direction. 
For any real number $\theta$, 
find a point $P(x, y)$ on the unit circle such that the line segment between the origin $O$ and the point $P$ makes an angle $\theta$ radians with the positive $x$-axis (see Figure \ref{figure31}). Then we define $\cos\theta$ and $\sin\theta$ to be the $x$ and $y$ coordinates of $P$:
\[x=\cos\theta, \hspace{1cm} y=\sin\theta.\]

\begin{figure}[ht]
\centering
\includegraphics[scale=0.2]{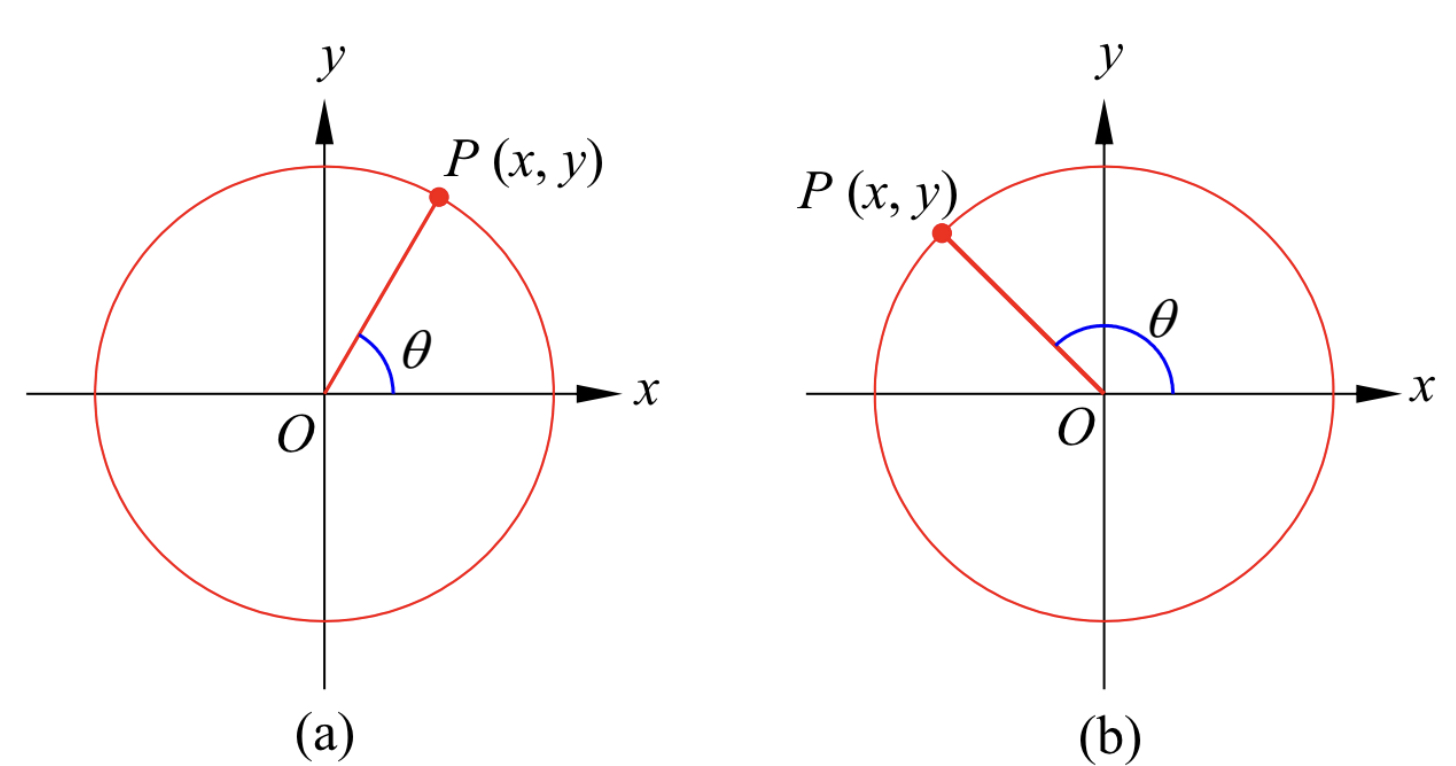}
\caption{The definitions of $\sin\theta$ and $\cos\theta$ for   (a)  $\di \theta=\frac{7\pi}{3}$ and  (b)  $\di\theta =\frac{3\pi}{4}$.\fa}\label{figure31}
\end{figure}
In this way, the function $\sin\theta$ and $\cos\theta$ are defined rigorously, and when $\theta$ is an acute angle, it coincides with the definition using right-angled triangles. From the definitions, it is obvious that $\sin\theta$ and $\cos\theta$ are periodic functions of periodic $2\pi$.

\begin{definition}{Periodic Functions}
A   function $f:\mathbb{R}\to\mathbb{R}$ is said to be periodic if there is a positive number $L$ so that
\[f(x+L)=f(x)\hspace{1cm}\text{for all}\; x\in \mathbb{R}.\]
Such a number $L$ is called a period of the function $f$. If $L$ is a period of $f$, then for any positive integer $n$, $nL$ is also a period of $f$.
\end{definition}

From the definitions, it is quite obvious that $\sin\theta$ and $\cos\theta$ are continuous functions. A rigorous proof is tedious. To show that these two functions are differentiable is also possible, but complicated. Two crucial formulas are
\begin{subequations}\label{eq230218_6}
\begin{align}
\sin(\theta_1+\theta_2)&=\sin\theta_1\cos\theta_2+\cos\theta_1\sin\theta_2,\label{eq230218_6a}\\
\cos(\theta_1+\theta_2)&=\cos\theta_1\cos\theta_2-\sin\theta_1\sin\theta_2.\label{eq230218_6b}
\end{align}\end{subequations}The proofs of these two formulas by elementary means are tedious.

In this section, we are going to define the sine and cosine functions using a different approach. We will show that the functions thus defined agree with the old definitions.

First, we present an existence and uniquess theorem.
\begin{theorem}[label=thm230218_3]{Existence and Uniqueness Theorem}
Let $\alpha$ and $\beta$ be any two real numbers. There exists a unique twice  differentiable function $f:\mathbb{R}\to\mathbb{R}$ satisfying
\[f''(x)+f(x)=0,\hspace{1cm}f(0)=\alpha,\;f'(0)=\beta.\]
\end{theorem}
Again, the proof of the existence requires knowledge from later chapters. We will prove uniqueness here. We begin by a lemma that will be useful later.
\begin{lemma}[label=lemma230218_5]{}Let  
 $f:\mathbb{R}\to\mathbb{R}$ be a twice differentiable  function that satisfies
\[f''(x)+f(x)=0.\]The following holds.
\begin{enumerate}[1.]
\item
  $f$ is infinitely differentiable. 
\item For any positive integer $n$, the $n^{\text{th}}$ derivative of $f$, $g(x)=f^{(n)}(x)$,   satisfies 
\[g''(x)+g(x)=0.\]
\item The function 
$ f(x)^2+f'(x)^2$ is  a constant.\end{enumerate}
\end{lemma}
\begin{myproof}{Proof}
Since $f$ is twice differentiable, $f$ is  continuous and differentiable. Since
$f''(x)=-f(x)$,   $f''$ is   continuous and differentiable.  This implies that $f$ is three times differentiable and $f'''=-f'$. Continue arguing in this way, we find that $f$ is infinitely differentiable, and for any nonengative integer $n$, \[f^{(n+2)}{x}=-f^{(n)}(x).\] The latter says that if $g=f^{(n)}$, then
\[g''(x)+g(x)=0.\]
These prove the first and second statements. For the third statement, 
we notice that
\begin{align*}\frac{d}{dx}\left(f'(x)^2+f(x)^2\right)&=f'(x)f''(x)+f(x)f'(x)\\&=2f'(x)\left(f''(x)+f(x)\right)=0.\end{align*}
This implies that $f(x)^2+f'(x)^2$ is  a constant.
\end{myproof}
Now we return to Theorem  \ref{thm230218_3}.
\begin{myproof}{\linkt   Proof of Theorem \ref{thm230218_3}}
If $f_1$ and $f_2$ are two functions that satisfy the given conditions, then the function $f=(f_1-f_2):\mathbb{R}\to\mathbb{R}$ is a  twice   differentiable function
satifying
\[f''(x)+f(x)=0,\hspace{1cm}f(0)=0,\; f'(0)=0.\]
To prove uniqueness, we only need to show that this function $f$ must be identically zero. 
By Lemma \ref{lemma230218_5}, there is a constant $C$ such that 
\[f'(x)^2+f(x)^2=C.\]
Setting $x=0$, we find that $C=0$. Hence,
\[f'(x)^2+f(x)^2=0.\]\bp
Since the square of a nonzero number is always positive, we must have
\[f(x)=f'(x)=0\hspace{1cm}\text{for all}\;x\in \mathbb{R}.\]
This completes the proof that $f$ is identically zero.

\end{myproof}
Notice that for a function $f:\mathbb{R}\to\mathbb{R}$ that satisfies $f''(x)+f(x)=0$, we have
\[f^{(4)}(x)=-f''(x)=f(x).\]
This implies that for all positive integers $n$,
\[f^{(4n)}(x)=f(x),\quad f^{(4n+1)}(x)=f'(x),\]
\[f^{(4n+2)}(x)=f''(x),\quad f^{(4n+3)}(x)=f'''(x).\]If $f$ is the unique solution to \[f''(x)+f(x)=0,\hspace{1cm}f(0)=\alpha,\; f'(0)=\beta,\] then its derivative $g=f'$ is the unique solution to
\[g''(x)+g(x)=0,\hspace{1cm}g(0)=\beta,\; g'(0)=-\alpha.\]

\begin{definition}{The Sine and Cosine functions}
The sine function $S(x)=\sin x$ is defined to be   the unique twice   differentiable function satisfying
\[S''(x)+S(x)=0,\hspace{1cm}S(0)=0,\;S'(0)=1.\]The cosine function $C(x)=\cos x$ is defined as the derivative of $S(x)$. Namely, 
$C(x)=S'(x)$. It is the unique twice differentiable function satisfying
\[C''(x)+C(x)=0, \hspace{1cm}C(0)=1,\;C'(0)=0.\]
\end{definition}

Notice that once we prove the existence of the function $S(x)=\sin x$, then the function $C(x)=\cos x$ exists. One can then check that the function
\[f(x)=\alpha C(x)+\beta S(x)\] is a twice differentiable function satifying \[f''(x)+f(x)=0,\hspace{1cm}f(0)=\alpha,\;f'(0)=\beta.\]In other words, to prove the existence part in Theorem 
\ref{thm230218_3}, we only need to establish the existence of the function $S(x)=\sin x$.

In the following, we establish the properties of the functions $S(x) $ and $C(x) $.

\begin{theorem}[label=thm230218_7]{}
The functions $S(x)$ and $C(x)$ are infinitely differentiable functions that satisfy  the following.
\begin{enumerate}[(a)]
\item $S'(x)=C(x)$ and $C'(x)=-S(x)$ for all $x\in\mathbb{R}$.
\item $S(x)$ is an odd function, $C(x)$ is an even function.
\item $S(x)^2+C(x)^2=1$ for all $x\in\mathbb{R}$.
\item For any real numbers $x$ and $y$, $S(x+y)=S(x)C(y)+C(x)S(y)$.
\item  For any real numbers $x$ and $y$, $C(x+y)=C(x)C(y)-S(x)S(y)$.

\end{enumerate}
\end{theorem}
\begin{myproof}{Proof} 
$S'(x)=C(x)$ is  by the definition of $C(x)$. Differentiating gives $C'(x)=S''(x)=-S(x)$. To prove (b), one check that the function $f(x)=-S(-x)$ satisfies
$f''(x)+f(x)=0$, $f(0)=0$ and $f'(0)=1$. By uniquess of the function $S(x)$, we have $f(x)=S(x)$, which proves that $S(x)$ is an odd function. Since $C(x)=S'(x)$, $C(x)$ is an even function.
Lemma \ref{lemma230218_5} says that $S(x)^2+S'(x)^2$ is a constant. Hence, there is a constant $A$ such that
\[S(x)^2+C(x)^2=A\hspace{1cm}\text{for all}\;x\in\mathbb{R}.\]
Setting $x=0$ gives $A=1$. This proves part (c). For part (d), fixed a real number $y$ and consider the function
\[f(x)=S(x+y).\]\bp
We find that \[f'(x)=S'(x+y)=C(x+y),\]  
\[f''(x)+f(x)=S''(x+y)+S(x+y)=0,\]
and
\[f(0)=S(y),\hspace{1cm}f'(0)=C(y).\]
Since the function $g(x)=S(y)C(x)+C(y)S(x)$  satisfies
\[g''(x)+g(x)=0,\hspace{1cm}g(0)=S(y),\;g'(0)=C(y),\]
by uniquesness, we find that $f(x)=g(x)$ for all $x\in \mathbb{R}$. Therefore,
\[S(x+y)=S(x)C(y)+C(x)S(y).\]
Differentiate with respect to $x$ gives
\[C(x+y)=C(x)C(y)-S(x)S(y).\]
\end{myproof}

Here we have used advanced analytic tools to prove the identities \eqref{eq230218_6} in a simple way.  Part (c) in Theorem \ref{thm230218_7} says that
\[\sin^2 x+\cos^2 x=1\hspace{1cm}\text{for all}\;x\in\mathbb{R}.\]
This implies that
\[|\sin x|\leq 1,\hspace{1cm} |\cos x|\leq 1\hspace{1cm}\text{for all}\;x\in\mathbb{R}.\]By definition, $\sin 0=S(0)=0$ and $\cos 0=C(0)=1$. What is not obvious is that  0 is in the range of $C(x)$. 
\begin{theorem}[label=thm230218_8]
{}
There is a smallest positive number $u$ such that $C(u)=0$. 
\end{theorem}
\begin{myproof}{Proof}
Since $S(x)$ is differentiable, we can apply mean value theorem to conclude that there is a point $v$ in $(0, 2)$ such that
\[\frac{S(2)-S(0)}{2-0}=S'(v)=C(v).\]
This gives
\[|C(v)|=\frac{1}{2}|S(2)|\leq\frac{1}{2}.\]By part (e) and part (c) in Theorem \ref{thm230218_7},
\[C(2v)= C(v)^2-S(v)^2= 2C(v)^2-1\leq \frac{1}{2}-1<0.\]
Since $C(2v)<0<C(0)$, and $C(x)$ is a continuous function, intermediate value theorem implies that there is a point $w$ in $(0, 2v)$ such that $C(w)=0$. Let 
\[A=\left\{w>0\,|\, C(w)=0.\right\}.\] We have just shown that $A$ is a nonempty set. By definition, $A$ is bounded below by 0. Hence, $u=\inf A$ exists. By Lemma \ref{23020510}, there is a sequence $\{w_n\}$ in $A$ that converges to $u$. Since $C(x)$ is continuous, the sequence $\{C(w_n)\}$ converges to $C(u)$. But $C(w_n)=0$ for all $n$. Hence, $C(u)=0$. Since $C(0)=1$, $u\neq 0$. Hence, $u>0$. This proves that 
$u$ is the smallest positive number such that $C(u)=0$. \end{myproof}

Let $u$ be the smallest positive number such that $C(u)=0$. Then we must have $C(x)>0$ for all $x\in [0, u)$. Since $S'(x)=C(x)$, $S(x)$ is strictly increasing on $[0, u]$. Thus, $S(x)>0$ for all $x\in (0, u]$. This, and $S(u)^2+C(u)^2=1$, implies that $S(u)=1$.
From part (d) and part (e) in Theorem \ref{thm230218_7}, we find that
\begin{subequations}
\begin{align*}
S(x+u)&=S(x)C(u)+C(x)S(u)=C(x),\\
C(x+u)&=C(x)C(u)-S(x)S(u)=-S(x).
\end{align*}
\end{subequations}It follows that
\begin{subequations}
\begin{align*}
S(x+2u)&=C(x+u)=-S(x),\\
C(x+2u)&= -S(x+u)=-C(x).
\\
S(x+3u)&=C(x+2u)=-C(x),\\
C(x+3u)&= -S(x+2u)=S(x).
\\
S(x+4u)&=C(x+3u)=S(x),\\
C(x+4u)&= -S(x+3u)=C(x).
\end{align*}
\end{subequations}The last pair of equations show that $S(x)$ and $C(x)$ are periodic functions of period $4u$. Since $S(x)>0$ and $C(x)>0$ for $x\in (0,u)$, we have the following.
\begin{enumerate}[$\bullet$\;\;]
\item For $x\in (0, u)$, $C(x)>0$, $S(x)>0$.
\item For $x\in (u, 2u)$, $C(x)<0$, $S(x)>0$.
\item For $x\in (2u, 3u)$, $C(x)<0$, $S(x)<0$.
\item For $x\in (3u, 4u)$, $C(x)>0$, $S(x)<0$.
\end{enumerate}
Together with $S(0)=0$, $C(0)=1$,  $S(u)=1$, $C(u)=0$, we find that $S(2u)=0$, $C(2u)=-1$, $S(3u)=-1$, $C(3u)=0$. These imply that for every $P(x, y)$ on the unit circle $x^2+y^2=1$, there is a unique $\theta\in [0, 4u)$ such that
\[x=C(\theta), \quad y=S(\theta).\]
 What is not obvious is that this $\theta$ is exactly the radian of the angle that the line segment $OP$ makes with the positive $x$-axis. To show this, we can argue in the following way. Assume that an object is travelling on the circle $x^2+y^2=1$, and its position at time $t$ is $(x(t), y(t))$, where
\[x=C(t), \quad y=S(t).\]
It follows that the velocity of the object at time $t$ is $(x'(t), y'(t))$, where
\[x'(t)=-S(t),\quad y'(t)=C(t).\]
This implies that the speed is
\[\sqrt{x'(t)^2+y'(t)^2}=\sqrt{S(t)^2+C(t)^2}=1.\]
Hence, the object is travelling at a constant speed 1. The distance travelled up to time $t$ is then $t$. This proves that the arclength of the arc from $(1,0)$ to the point $P(C(t), S(t))$ is $t$. Then $t$ must be the radian of the angle $OP$ makes with the positive $x$ axis. Hence, the functions $C(t)$ and $S(t)$ coincide with the classical $\cos t$ and $\sin t$ functions. Having proved this, by the definition of $\pi$, we have
\[2u=\pi.\]
Hence, we can summarize the facts above as follows.
\begin{highlight}{Properties of the Sine and Cosine Functions}
The functions $S(x)=\sin x$ and $C(x)=\cos x$ are $2\pi$ periodic infinitely differentiable functions. 
\[\frac{d}{dx}\sin x=\cos x,\hspace{1cm}\frac{d}{dx}\cos x=-\sin x.\]

Moreover, they have the following properties.
\begin{enumerate}[1.]
\item $\sin\left(x+\frac{\pi}{2}\right)=\cos x$, $\cos\left(x+\frac{\pi}{2}\right)=-\sin x$.
\item $\sin\left(x+\pi\right)=-\sin x$, $\cos\left(x+\pi \right)=-\cos x$.
\item $\sin\left(x+\frac{3\pi}{2}\right)=-\cos x$, $\cos\left(x+\frac{3\pi}{2}\right)=\sin x$.
\item $\sin x$ is an odd function, $\cos x$ is an even function.
\item $\sin x=0$ if and  only if $x=n\pi$, where $n$ is an integer.
\item $\cos x=0$ if and only if $x=\left(n+\frac{1}{2}\right)\pi$,  where $n$ is an integer.
\end{enumerate}
\end{highlight}

\begin{figure}[ht]
\centering
\includegraphics[scale=0.2]{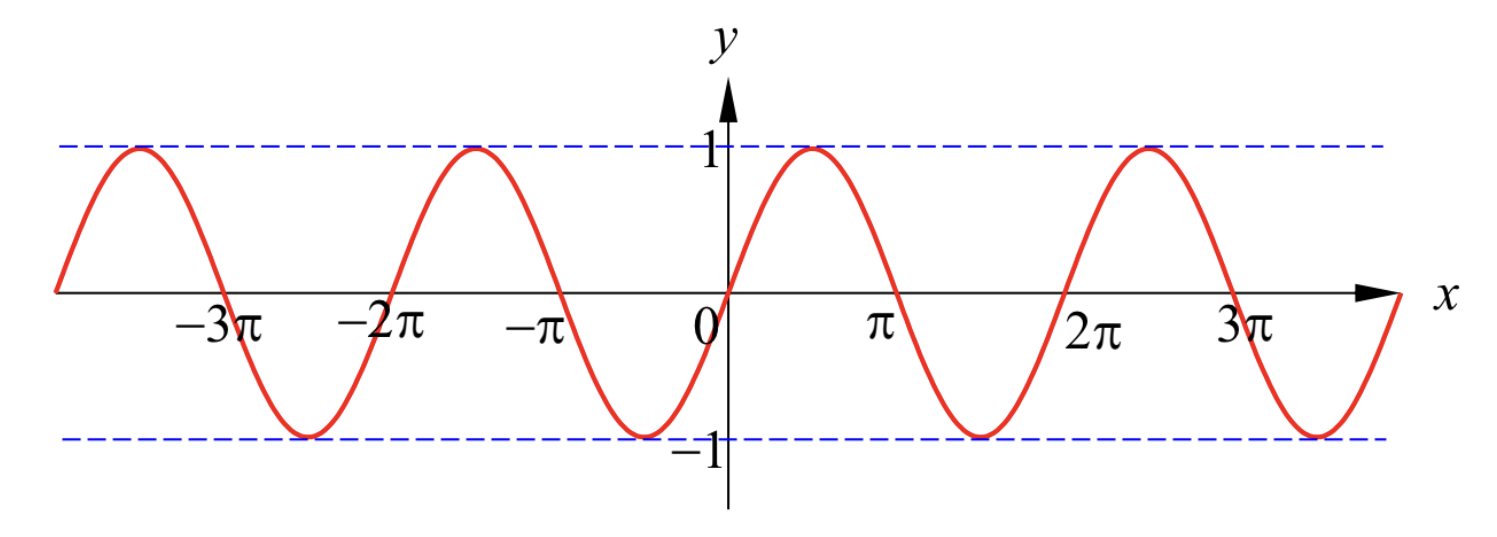}
\caption{The sine function $S(x)=\sin x$.\fa}\label{figure32}
\end{figure}

\begin{figure}[ht]
\centering
\includegraphics[scale=0.2]{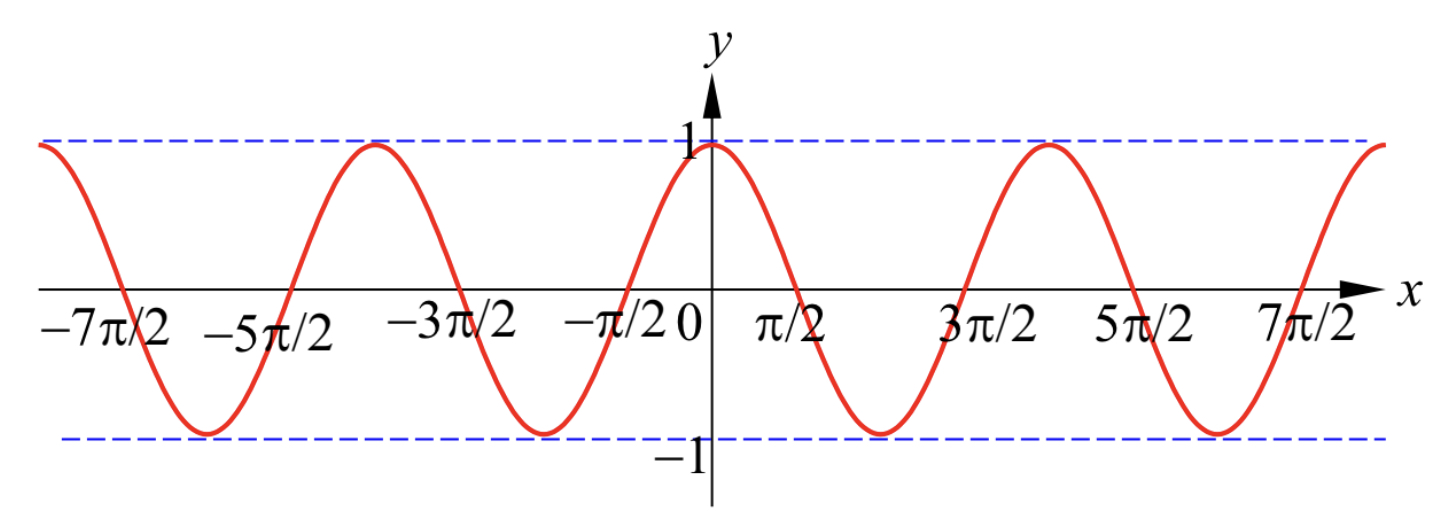}
\caption{The cosine function $C(x)=\cos x$.\fa}\label{figure33}
\end{figure}

There are four other trigonometric functions. They are defined in terms of $\sin x$ and $\cos x$ in the usual way.
\begin{definition}{Trigonmetric Functions}
The tangent, cotangent, secant and cosecant functions are defined as
\begin{align*}
\tan x&=\frac{\sin x}{\cos x},\hspace{1cm}\sec x=\frac{1}{\cos x},\hspace{1cm} x\neq \left(n+\frac{1}{2}\right)\pi, n\in\mathbb{Z};\\
\cot x&=\frac{\cos x}{\sin x},\hspace{1cm}\csc x=\frac{1}{\sin x},\hspace{1cm} x\neq n\pi, n\in\mathbb{Z}.\end{align*}
\end{definition}

\begin{figure}[ht]
\centering
\includegraphics[scale=0.2]{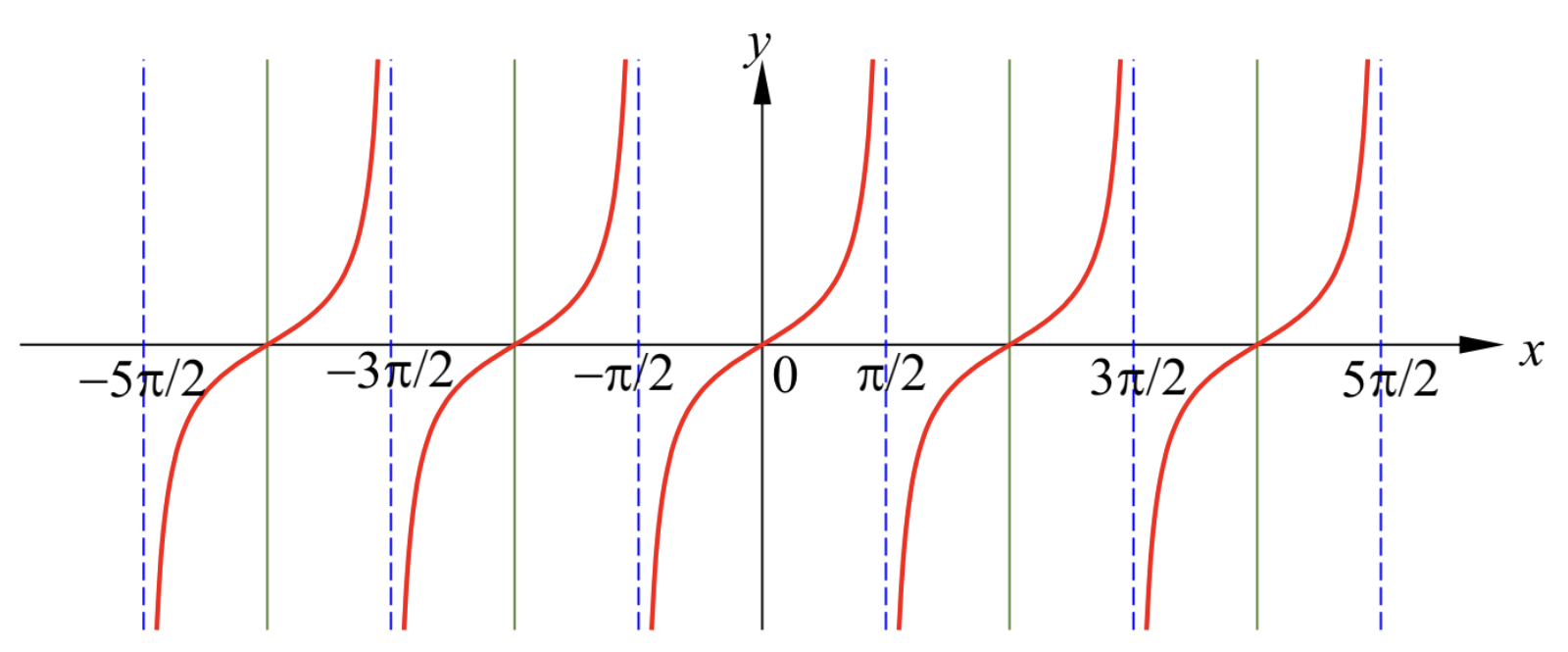}
\caption{The   function $f(x)=\tan x$.\fa}\label{figure34}
\end{figure}
The  following are easy to derive.
\begin{proposition}{}
$\tan x$, $\cot x$, $\sec x$ and $\csc x$ are infinitely differentiable functions with
\begin{align*}
\frac{d}{dx}\tan x&=\sec^2 x,\hspace{1.8cm} \frac{d}{dx}\cot x=-\csc^2x,\\
\frac{d}{dx}\sec x&=\sec x\tan x,\hspace{1cm}\frac{d}{dx}\csc x=-\csc x\cot x.
\end{align*}
\end{proposition}

Before closing this subsection, we want to prove some important limits and  inequalities for the function $\sin x$.
\begin{theorem}[label=230307_8]{}
\begin{enumerate}[1.]
\item
For any real number $x$, $|\sin x|\leq |x|$.
\item $\di \lim_{x\to 0}\frac{\sin x}{x}=1$.
\item For any $x\in [0, \pi/2]$, 
\[\frac{2}{\pi}x\leq \sin x \leq x.\]
\end{enumerate}
\end{theorem}
\begin{myproof}{Proof}
 
When $x=0$, $\sin x=0$ and $|\sin x|\leq |x|$ is obviously true. If $x\neq 0$, mean value theorem implies that there is a number $c$ in $(0,1)$ such that
\begin{equation}\label{eq230218_9}\frac{\sin x}{x}=\frac{\sin x-\sin 0}{x-0}=\cos (cx).\end{equation}
Hence, \[\left|\frac{\sin x}{x}\right|=|\cos (cx)|\leq 1,\]which implies that $|\sin x|\leq |x|$. This proves the first statement.\bp
For the second statement, the definition of derivative implies that
\begin{equation*}\lim_{x\to 0}\frac{\sin x}{x}=\lim_{x\to 0}\frac{\sin x-\sin 0}{x-0}=\left.\frac{d}{dx}\right|_{x=0}\sin x=\cos 0=1.\end{equation*}
  For the third statement, define the function $g:[0,\frac{\pi}{2}]\to\mathbb{R}$ by
\[g(x)=\begin{cases}\di\frac{\sin x}{x},\quad &\text{if}\;0<x\leq\di\frac{\pi}{2},\\1,\quad &\text{if}\;x=0.\end{cases}\]Then $g$ is continuous on $[0, \frac{\pi}{2}]$, and differentiable on $(0, \frac{\pi}{2})$, with
\[g'(x)=\frac{x\cos x-\sin x}{x^2}=\frac{1}{x}\left(\cos x-\frac{\sin x}{x}\right)\hspace{1cm}\text{when}\;0<x<\frac{\pi}{2}.\]
As before, for each $x\in (0, \frac{\pi}{2}]$, mean value theorem implies that there is a $u\in (0, x)$ such that
\[\frac{\sin x}{x}=\cos u.\]
Since $0<u<x$ and the cosine function is strictly decreasing on $(0, \frac{\pi}{2})$, we find that
\[g'(x)=\frac{1}{x}\left(\cos x-\frac{\sin x}{x}\right)=\frac{1}{x}\left(\cos x-\cos u\right)<0.\]
This shows that $g:[0,\frac{\pi}{2}]\to\mathbb{R}$ is a strictly decreasing function. Since $g(0)=1$ and $g(\frac{\pi}{2})=\frac{2}{\pi}$, we find that for all $x\in (0, \frac{\pi}{2}]$,
\[\frac{2}{\pi}\leq\frac{\sin x}{x}\leq 1.\]
Thus, for all $x\in [0, \pi/2]$, 
\[\frac{2}{\pi}x\leq \sin x \leq x.\]

\end{myproof}

\subsection{The Inverse Trigonometric Functions}
In this section, we are going to define inverse functions for $\sin x$, $\cos x$ and $\tan x$.
Since trigonometric functions are periodic functions, they are not one-to-one. Hence, we cannot find their inverses over the whole domain of their definitions. However, we can restrict each of  their domains to  an interval on which each of  them is one-to-one to define the inverse. Such interval should contain the interval $(0, \pi/2)$ which is where these functions are classically defined. 

\begin{enumerate}[$\bullet$\;\;]
\item The largest interval that contains  the interval $(0, \pi/2)$ and on which $\sin x$ is one-to-one is $[-\frac{\pi}{2}, \frac{\pi}{2}]$.
\item The largest interval that contains  the interval $(0, \pi/2)$ and on which $\cos x$ is one-to-one is $[0, \pi]$.
\item The largest interval that contains  the interval $(0, \pi/2)$ and on which $\tan x$ is one-to-one is $(-\frac{\pi}{2}, \frac{\pi}{2})$.
\end{enumerate}

\begin{definition}{Inverse Sine Function }
  The function  $\sin^{-1}x$ is a function defined on $[-1, 1]$ and with range $[-\frac{\pi}{2}, \frac{\pi}{2}]$ such that
\begin{align*}\sin(\sin^{-1}x)&=x\hspace{1cm}\text{for all}\;x\in[-1,1];\\
\sin^{-1}(\sin x)&=x\hspace{1cm}\text{for all}\; x\in\left[-\frac{\pi}{2}, \frac{\pi}{2}\right].\end{align*}

\end{definition}
\begin{definition}{Inverse Cosine Function }
  The function  $\cos^{-1}x$ is a function defined on $[-1, 1]$ and with range $[0, \pi]$ such that
\begin{align*}\cos(\cos^{-1}x)&=x\hspace{1cm}\text{for all}\;x\in[-1,1];\\
\cos^{-1}(\cos x)&=x\hspace{1cm}\text{for all}\; x\in [0, \pi].\end{align*}

\end{definition}
\begin{definition}{Inverse Tangent Function }
  The function  $\tan^{-1}x$ is a function defined on $\mathbb{R}$ and with range $(-\frac{\pi}{2}, \frac{\pi}{2})$ such that
\begin{align*}\tan(\tan^{-1}x)&=x\hspace{1cm}\text{for all}\;x\in\mathbb{R};\\
\tan^{-1}(\tan x)&=x\hspace{1cm}\text{for all}\; x\in \left(-\frac{\pi}{2}, \frac{\pi}{2}\right).\end{align*}

\end{definition}
The differentiability of the inverse trigonometric functions and their derivative formulas follow  immediately from Theorem \ref{thm230218_9}.

\begin{theorem}{}
$\sin^{-1}:(-1,1)\to \mathbb{R}$ is a differentiable function with
\[\frac{d}{dx}\sin^{-1}x=\frac{1}{\sqrt{1-x^2}}.\]
\end{theorem}

\begin{theorem}{}
$\cos^{-1}:(-1,1)\to \mathbb{R}$ is a differentiable function with
\[\frac{d}{dx}\cos^{-1}x=-\frac{1}{\sqrt{1-x^2}}.\]
\end{theorem}

\begin{theorem}{}
$\tan^{-1}:\mathbb{R}\to \mathbb{R}$ is a differentiable function with
\[\frac{d}{dx}\tan^{-1}x=\frac{1}{1+x^2}.\]
\end{theorem}
\vp

\noindent
{\bf \large Exercises  \thesection}
\setcounter{myquestion}{1}
\begin{question}{\themyquestion}
Determine the following limits.
\begin{enumerate}[(a)]
\item
$\di\lim_{n\to \infty}\left(1-\frac{1}{n}\right)^n$
\item
$\di\lim_{n\to \infty}\left(1+\frac{2}{n}\right)^n$
\item
$\di\lim_{n\to \infty}\left(1-\frac{2}{n}\right)^n$
 
\end{enumerate}
\end{question}
\atc

\begin{question}{\themyquestion}
For any $x\in [-1,1]$, show that $\sin^{-1}x+\cos^{-1}x$ is a constant and find this constant.
\end{question}
\atc
\begin{question}{\themyquestion}
Determine the following limits.
\begin{enumerate}[(a)]
\item
$\di\lim_{x\to \frac{\pi}{2}^-}\tan x$
\item
$\di\lim_{x\to -\frac{\pi}{2}^+}\tan x$
\item
$\di\lim_{x\to -\infty}\tan^{-1} x$
\item
$\di\lim_{x\to \infty}\tan^{-1} x$
\end{enumerate}
\end{question}

\atc

\begin{question}{\themyquestion}
Consider the function $f:\mathbb{R}\to\mathbb{R}$ defined by
\[f(x)=\begin{cases} \di \sin\left(\frac{1}{x}\right),\quad &\text{if}\;x\neq 0,\\
0,\quad &\text{if}\;x= 0.\end{cases}\]
Determine whether $f$ is a continuous function. If not, find the points where the function $f$ is not continuous.
\end{question}

\atc

\begin{question}{\themyquestion}
Consider the function $f:\mathbb{R}\to\mathbb{R}$ defined by
\[f(x)=\begin{cases} \di x\sin\left(\frac{1}{x}\right),\quad &\text{if}\;x\neq 0,\\
0,\quad &\text{if}\;x= 0.\end{cases}\]
Show that $f$ is a continuous function.
\end{question}

\atc

\begin{question}{\themyquestion}
Consider the function $f:\mathbb{R}\to\mathbb{R}$ defined by
\[f(x)=\begin{cases} \di x^2\sin\left(\frac{1}{x}\right),\quad &\text{if}\;x\neq 0,\\
0,\quad &\text{if}\;x= 0.\end{cases}\]Let $g:\mathbb{R}\to\mathbb{R}$ be the function defined by
\[g(x)=x+f(x).\]
\begin{enumerate}[(a)]
\item Show that $f$ is a differentiable function.
\item Show that $f':\mathbb{R}\to\mathbb{R}$ is not continuous.
\item Show that $g'(0)=1$, but for any neighbourhood $(a,b)$ of 0,   $g:(a,b)\to\mathbb{R}$ is not increasing.
\end{enumerate}
\end{question}

\vp
\section{L' H$\hat{\text{o}}$pital's Rules}\label{sec3.6}
In this section, we will apply the Cauchy mean value theorem to prove the l' H$\hat{\text{o}}$pital's rules. The latter are useful rules for finding limits of the form
\[\lim_{x\to x_0}\frac{f(x)}{g(x)},\]
when we have one of the following two indeterminate forms.
\begin{enumerate}[1.]
\item Type $0/0$, where
$\di\lim_{x\to x_0}f(x)=0$ and $\di\lim_{x\to x_0}g(x)=0$.
\item Type $\infty/\infty$, where $\di\lim_{x\to x_0}f(x)=\infty$ and $\di\lim_{x\to x_0}g(x)=\infty$.
\end{enumerate}
Here $x_0$ can be $\infty$ or $-\infty$.

Let us first prove the following special case.
\begin{theorem}[label=thm230216_14]{}
Let $f:(a,b)\to\mathbb{R}$ and $g:(a,b)\to\mathbb{R}$ be differentiable functions that satisfy the following conditions.
\begin{enumerate}[(i)]
\item $\di\lim_{x\rightarrow a^+}f(x)=\lim_{x\rightarrow a^+}g(x)=0$.
\item $\di\lim_{x\rightarrow a^+}\frac{f'(x)}{g'(x)}=L$.

\end{enumerate}
Then $\di \lim_{x\rightarrow a^+}\frac{f(x)}{g(x)}=L$.
\end{theorem}
\begin{myproof}{Proof}
The condition (i) implies that we can extend $f$ and $g$ to be continuous functions on $[a, b)$ by defining $f(a)=g(a)=0$. Then by Cauchy mean value theorem, for any $x\in (a, b)$, there is a point $u(x)\in (a, x)$ such that
\begin{equation}\label{eq230216_15}\frac{f(x)}{g(x)}=\frac{f(x)-f(a)}{g(x)-g(a)}=\frac{f'(u(x))}{g'(u(x))}.\end{equation}\bp
 Since 
\[a<u(x)<x,\]
squeeze theorem implies that
\[\lim_{x\to a^+}u(x)=a.\] 
By limit law for composite functions, we find that
\[\lim_{x\to a^+}\frac{f'(u(x))}{g'(u(x))}=\lim_{u\to a^+}\frac{f'(u)}{g'(u)}=L.\]
By \eqref{eq230216_15}, this proves that
\[\lim_{x\rightarrow a^+}\frac{f(x)}{g(x)}=L.\]
\end{myproof}

It is easy to see that an analogue of Theorem \ref{thm230216_14} holds for left limits.  Combine the left limit and the right limit, we have the following.

\begin{theorem}[label=thm230216_16]{l' H$\hat{\text{o}}$pital's Rule I}
Let $x_0$ be a point in the open interval $(a, b)$, and let   $D=(a,b)\setminus\{x_0\}$. Given that $f:D\to\mathbb{R}$ and $g:D\to\mathbb{R}$ are diferentiable  functions that satisfy the following conditions.
\begin{enumerate}[(i)]
\item $\di\lim_{x\rightarrow x_0}f(x)=\lim_{x\rightarrow x_0}g(x)=0$.
\item  $\di\lim_{x\rightarrow x_0}\frac{f'(x)}{g'(x)}=L$.

\end{enumerate}
Then  we have $\di \lim_{x\rightarrow x_0}\frac{f(x)}{g(x)}=L$.
\end{theorem}
 We return to a problem that we discussed earlier.
\begin{example}{}Determine the limit \[\lim_{x\to 1}\frac{x^{20}+2x^9-3}{x^7-1}\] if it exists.
\end{example}
\begin{solution}{Solution}
Let $f(x)=x^{20}+2x^9-3$ and $g(x)=x^7-1$. Then
\[\lim_{x\to 1}f(x)=f(1)=0\hspace{1cm} \text{and}\hspace{1cm} \lim_{x\to 1}g(x)=g(1)=0.\]
$f$ and $g$ are continuously differentiable functions with
\[f'(x)=20x^{19}+18x^8\hspace{1cm} \text{and}\hspace{1cm}  g'(x)=7x^6.\]
Since\[\lim_{x\to 1}\frac{f'(x)}{g'(x)}=\lim_{x\to 1}\frac{20x^{19}+18x^8}{7x^6}=\frac{38}{7},\]
l' H$\hat{\text{o}}$pital's rule implies that
\[\lim_{x\to 1}\frac{x^{20}+2x^9-3}{x^7-1}=\frac{38}{7}.\]
\end{solution}
Let us look at some other examples.
\begin{example}{}
Determine whether the limit exists. If it exists, find the limit.
\begin{enumerate}[(a)]
\item
$\di \lim_{x\to 0}\frac{e^x-1-x}{x^2}$
\item $\di \lim_{x\to 0}\frac{\sin 2x}{3x}$
\item $\di \lim_{x\to 0}\frac{\cos 2x-1}{x^2}$
\end{enumerate}
\end{example}
\begin{solution}{Solution}
\begin{enumerate}[(a)]
\item
This is a limit of the form $0/0$. Applying l' H$\hat{\text{o}}$pital's rule, we have
\[\lim_{x\to 0}\frac{e^x-1-x}{x^2}=\lim_{x\to 0}\frac{e^x-1}{2x}.\]
Again, we have a limit of the form $0/0$. Applying l' H$\hat{\text{o}}$pital's rule again, we have
\[\lim_{x\to 0}\frac{e^x-1-x}{x^2}=\lim_{x\to 0}\frac{e^x-1}{2x}=\lim_{x\to 0}\frac{e^x}{2}=\frac{1}{2}.\]
\item This is a limit of the form $0/0$.  Apply l' H$\hat{\text{o}}$pital's rule, we have
\[\lim_{x\to 0}\frac{\sin 2x}{3x}=\lim_{x\to 0}\frac{2\cos 2x}{3}=\frac{2}{3}.\]
\item This is a limit of the form $0/0$. Applying l' H$\hat{\text{o}}$pital's rule twice, we have
\[\lim_{x\to 0}\frac{\cos 2x-1}{x^2}=\lim_{x\to 0}\frac{-2\sin 2x}{2x}=\lim_{x\to 0}\frac{-4\cos 2x}{2}=-2.\]
\end{enumerate}
\end{solution}
 
Using l' H$\hat{\text{o}}$pital's rule, we can give a second solution to Example \ref{ex230216_17}.
\begin{example}{}
 Since $f$ is continuous, we have
\[\lim_{h\to 0}\left(f(x_0+h)+f(x_0-h)-2f(x_0)\right)=0.\]Since we also have $\di\lim_{h\to 0}h^2=0$, we can apply l' H$\hat{\text{o}}$pital's rule to get
\[\lim_{h\to 0}\frac{f(x_0+h)+f(x_0-h)-2f(x_0)}{h^2} =\lim_{h\to 0}\frac{f'(x_0+h)-f'(x_0-h) }{2h}.\] 
 Since $f'$ is  continuous, \[\lim_{h\to 0}\left(f'(x_0+h)-f'(x_0-h) \right)=0.\]\be Since we also have $\di\lim_{h\to 0}(2h)=0$, applying l' H$\hat{\text{o}}$pital's rule again give
\[\lim_{h\to 0}\frac{f'(x_0+h)-f'(x_0-h) }{2h}=\lim_{h\to 0}\frac{f''(x_0+h)+f''(x_0-h)}{2}.\]
It follows from the continuity of $f''$ that
\[ \lim_{h\to 0}\frac{f''(x_0+h)+f''(x_0-h)}{2}=f''(x_0).\]These prove that
\[\lim_{h\to 0}\frac{f(x_0+h)+f(x_0-h)-2f(x_0)}{h^2}  =f''(x_0).\]

\end{example2}

In the future, we are going to see that Taylor's approximation is an alternative to l' H$\hat{\text{o}}$pital's rule when the point $x_0$ is finite and the indeterminate form if of the type $0/0$. However, when $x_0$ is infinite or the indeterminate form is of type $\infty/\infty$,  l' H$\hat{\text{o}}$pital's rule becomes useful.

The following is for the case where $x_0$ is infinite, and the limit is of the form $0/0$.
\begin{theorem}[label=thm230216_18]{l' H$\hat{\text{o}}$pital's Rule II}
Let $a$ be a positive number. Given that $f:(a,\infty)\to\mathbb{R}$ and $g:(a,\infty)\to\mathbb{R}$ are diferentiable  functions that satisfy the following conditions.
\begin{enumerate}[(i)]
\item $\di\lim_{x\rightarrow \infty}f(x)=\lim_{x\rightarrow \infty}g(x)=0$.
\item  $\di\lim_{x\rightarrow \infty}\frac{f'(x)}{g'(x)}=L$.

\end{enumerate}
Then  we have $\di \lim_{x\rightarrow \infty}\frac{f(x)}{g(x)}=L$.
\end{theorem}
\begin{myproof}{Proof}
Let $b=1/a$, and define the functions $f_1:(0, b)\to\mathbb{R}$ and $g_1:(0,b)\to\mathbb{R}$ by
\[f_1(x)=f\left(\frac{1}{x}\right),\hspace{1cm}g_1(x)=g\left(\frac{1}{x}\right).\]  Then $f_1$ and $g_1$ are differentiable functions and
\[f_1'(x)=-\frac{1}{x^2}f'\left(\frac{1}{x}\right),\hspace{1cm}g_1'(x)=-\frac{1}{x^2}g'\left(\frac{1}{x}\right).\]
Moreover,
\[\lim_{x\to 0^+}f_1(x)=\lim_{x\to\infty}f(x)=0,\hspace{1cm}\lim_{x\to 0^+}g_1(x)=\lim_{x\to\infty}g(x)=0,\]
and
\[\lim_{x\to 0^+}\frac{f_1'(x)}{g_1'(x)}=\lim_{x\to 0^+}\frac{f'\left(\frac{1}{x}\right)}{g'\left(\frac{1}{x}\right)}=\lim_{x\to\infty}\frac{f'(x)}{g'(x)}=L.\]
By Theorem \ref{thm230216_14}, 
\[\lim_{x\to 0^+}\frac{f_1(x)}{g_1 (x)}=L.\]This implies that 
\[ \lim_{x\rightarrow \infty}\frac{f(x)}{g(x)}=L.\]
\end{myproof}

 Let us look at the following example.
\begin{example}{}
Determine whether the limit $\di\lim_{x\to \infty}\left(\frac{x}{x+2}\right)^{x+1}$ exists. If it exists, find the limit.
\end{example}
This is not of the type $0/0$. But the logarithm of it can be turned into that form.
\begin{solution}{Solution}
Consider the function 
\[g(x)=\ln \left(\frac{x}{x+2}\right)^{x+1}=(x+1)\ln \left(\frac{x}{x+2}\right).\]
When $x\to \infty$, we have something of the form $\infty \,\cdot\,0$. We turn it to the form $0/0$ by
\[g(x)=\frac{\di \ln\left(\frac{x}{x+2}\right)}{\di \frac{1}{x+1}}=\frac{\ln x-\ln(x+2)}{\di \frac{1}{x+1}}.\]
 l' H$\hat{\text{o}}$pital's rule implies that
\begin{align*}\lim_{x\to\infty}g(x)&=\lim_{x\to\infty}\frac{\di\frac{1}{x}-\frac{1}{x+2}}{-\di \frac{1}{(x+1)^2}}\\
&=-2\lim_{x\to \infty}\frac{x^2+2x+1}{x^2+2x}\\
&=-2.
\end{align*}
By continuity of the exponential function, we have
\[\lim_{x\to \infty}\left(\frac{x}{x+2}\right)^{x+1}=\lim_{x\to \infty}e^{g(x)}=\exp\left(\lim_{x\to \infty}g(x)\right)=e^{-2}.\]
\end{solution}

Suppose we want to find the limit
\begin{equation}\label{eq230218_2}\lim_{x\rightarrow \infty}\frac{x}{e^x}.\end{equation}
This is a limit of the form $\infty/\infty$. One may say that we can turn it to a limit of the form $0/0$ by writing
\[\frac{x}{e^x}=\frac{e^{-x}}{x^{-1}}.\]
Then l' H$\hat{\text{o}}$pital's rule says that if the limit
\begin{equation}\label{eq230218_3}\lim_{x\to\infty}\frac{e^{-x}}{- x^{-2}}\end{equation} exists and is equal to $L$, the limit \eqref{eq230218_2} also exists and is equal to $L$. However, the limit \eqref{eq230218_3} is   more complicated than the limit \eqref{eq230218_2}. So this strategy is useless.
Hence, there is a need for us to consider the $\infty/\infty$ indeterminate case. We only  prove the theorem in the case $x_0$ is finite. The case where $x_0$ is infinite can be dealt with in the same way as in the proof of Theorem \ref{thm230216_18}.

\begin{theorem}[label=thm230216_19]{l' H$\hat{\text{o}}$pital's Rule III}
Let $x_0$ be a point in the open interval $(a, b)$, and let  $D$ be the set $D=(a,b)\setminus\{x_0\}$. Given that $f:D\to\mathbb{R}$ and $g:D\to\mathbb{R}$ are diferentiable  functions that satisfy the following conditions.
\begin{enumerate}[(i)]
\item $\di\lim_{x\rightarrow x_0}f(x)=\lim_{x\rightarrow x_0}g(x)=\infty$.
\item  $\di\lim_{x\rightarrow x_0}\frac{f'(x)}{g'(x)}=L$.

\end{enumerate}
Then  we have $\di \lim_{x\rightarrow x_0}\frac{f(x)}{g(x)}=L$.
\end{theorem}
 The proof of this theorem is technical because of the infinite limits. The strategegy to rewrite this as 
\[\lim_{x\to x_0}\frac{1/g(x)}{1/f(x)}\] is not useful, as have been demonstrated in our discussion before this theorem.
\begin{myproof}{Proof}
We will prove that  the right limit $\di \lim_{x\rightarrow x_0^+}\frac{f(x)}{g(x)}$ is equal to $L$. The proof that the left limit is equal to $L$ is similar.  
Observe that if we fix a point $u$ in $(x_0, b)$, then for any $x$ in $(x_0, u)$, Cauchy mean value theorem asserts that there is a $c_x$ in $(x, u)$ such that 
\[\frac{f(x)-f(u)}{g(x)-g(u)}=\frac{f'(c_x)}{g'(c_x)}.\]\bp
This implies that
\[\frac{(f(x)-Lg(x))-(f(u)-Lg(u))}{g(x)-g(u)}=\frac{f'(c_x)}{g'(c_x)}-L.\]Thus,
\begin{equation}\label{eq230217_2}\frac{f(x)}{g(x)}-L=\frac{g(x)-g(u)}{g(x)}\left(\frac{f'(c_x)}{g'(c_x)}-L\right)+\frac{f(u)-Lg(u)}{g(x)}.\end{equation}

Fixed $\varepsilon>0$. By assumption of
\[\lim_{x\rightarrow x_0^+}\frac{f'(x)}{g'(x)}=L,\]
  there exsits a $\delta_1>0$ such that $(x_0, x_0+\delta_1)\subset (a,b)$, and for any $x\in (x_0, x_0+\delta_1)$, 
\[\left|\frac{f'(x)}{g'(x)}-L\right|<\frac{\varepsilon}{3}.\]Take $u=x_0+\delta_1/2$. 
Since $\di \lim_{x\rightarrow x_0}g(x)=\infty$, we find that
\[\lim_{x\to x_0^+}\frac{g(x)-g(u)}{g(x)}=1\hspace{1cm}\lim_{x\to x_0^+}\frac{f(u)-Lg(u)}{g(x)}=0.\]
Therefore, there exists a number $\delta$ such that $0<\delta\leq\delta_1/2$, and for all $x\in (x_0, x_0+\delta)$,
\[\left|\frac{g(x)-g(u)}{g(x)}\right|<2,\hspace{1cm}\left|\frac{f(u)-Lg(u)}{g(x)}\right|<\frac{\varepsilon}{3}.\]
If $x$ is in $(x_0, x_0+\delta)$, $x_0<x<u$ and hence $x_0<c_x<u<x_0+\delta_1$. This implies that
\[\left|\frac{f'(c_x)}{g'(c_x)}-L\right|<\frac{\varepsilon}{3}.\] 
Eq. \eqref{eq230217_2} then implies that for all $x\in (x_0, x_0+\delta)$,
\begin{align*}
\left|\frac{f(x)}{g(x)}-L\right|&\leq \left|\frac{g(x)-g(u)}{g(x)}\right|\left| \frac{f'(c_x)}{g'(c_x)}-L\right|+\left|\frac{f(u)-Lg(u)}{g(x)}\right|\\
&<2\times\frac{\varepsilon}{3}+\frac{\varepsilon}{3}=\varepsilon.
\end{align*}\bp
This proves that 
\[ \lim_{x\rightarrow x_0^+}\frac{f(x)}{g(x)}=L.\]

\end{myproof} Notice that in the proof, we do not use the assumption that $\di\lim_{x\to x_0}f(x)=\infty$. Hence, this can be ommited from the conditions in the theorem. If $f(x)$ is bounded in a neighbourhood of $x_0$, there is no need to apply l' H$\hat{\text{o}}$pital's rule.

Let us now look at some examples.
\begin{example}[label=230307_11]{}
Let $r$ be a positive number.
Prove that
\[\lim_{x\to \infty}\frac{\ln x}{x^{r}}=0.\]
Deduce that for any positive number $s$,
\[\lim_{x\to \infty}x^se^{-x}=0.\]
\end{example}
\begin{solution}{Solution}
The limit \[\lim_{x\to \infty}\frac{\ln x}{x^{r}}\] is of the form $\infty/\infty$. Apply l' H$\hat{\text{o}}$pital's rule, we have
\[\lim_{x\to \infty}\frac{\ln x}{x^{r}}=\lim_{x\to\infty}\frac{\di \frac{1}{x}}{\di rx^{r-1}}=\frac{1}{r}\lim_{x\to \infty}\frac{1}{x^r}=0.\]  Since $\di\lim_{x\to \infty}e^x=\infty$, and the function $f(x)=x^s$ is a continuous function, we find that
\begin{align*}
\lim_{x\to\infty}x^se^{-x}&=\lim_{u\to \infty} \frac{(\ln u)^s}{u}\\
&=\left(\lim_{u\to \infty} \frac{\ln u}{u^{1/s}}\right)^s\\
&=0^s=0.
\end{align*}
\end{solution}
\begin{highlight}{}
The result of this example shows that when $x$ becomes large,\begin{enumerate}[$\bullet$\;\;]
\item $\ln x$ goes to infinity slower than any positive powers of $x$;
\item any positive powers of $x$ goes to infinity slower than $e^x$.\end{enumerate}
\end{highlight}

\begin{example}{}
Show that there exists a number $c$  so that the function
\[f(x)=\begin{cases} x^x,\quad &\text{if}\;x>0,\\
c,\quad &\text{if}\;x=0,\end{cases}\] is continuous.
\end{example}\begin{solution}{Solution}
Since
\[f(x)=x^x=e^{x\ln x}\hspace{1cm}\text{when}\;x>0,\] $f(x)$ is continuous on $(0, \infty)$.  To make $f$ continuous, $f$ must be continuous at $x=0$. This means
\[c=f(0)=\lim_{x\to 0^+}f(x)=\lim_{x\to 0^+}e^{x\ln x}.\]
Let us look at the limit $\di\lim_{x\to 0^+}x\ln x$. It is of the form $0\,\cdot\,\infty$. We turn it to the form $\infty/\infty$, and use l' H$\hat{\text{o}}$pital's rule.  
\[\lim_{x\to 0^+}x\ln x=\lim_{x\to 0^+}\frac{\ln x}{\di\frac{1}{x}}=\lim_{x\to 0^+}\frac{\di \frac{1}{x}}{-\di\frac{1}{x^2}}=-\lim_{x\to 0^+}x=0.\]Therefore, when
\[c=\exp\left(\lim_{x\to 0^+}x\ln x\right)=e^0=1,\]the function $f:[0, \infty)\to\mathbb{R}$ is continuous.
\end{solution}

\vp
\noindent
{\bf \large Exercises  \thesection}
\setcounter{myquestion}{1}

\begin{question}{\themyquestion}
Determine whether the limit exists. If it exists, find the limit.
\begin{enumerate}[(a)]
\item
$\di\lim_{x\to 1}\frac{x^{100}+x^{50}-2}{3x^{101}+4x^{55}-7x}$
\item $\di \lim_{x\to 0}\frac{2e^{-x}-2+2x}{x^2+3x^3}$
\item $\di\lim_{x\to 0}\frac{\tan^{-1}x}{x}$
\item $\di\lim_{x\to 0}\frac{\tan x-x}{x^3}$
\end{enumerate}
\end{question}
\atc
\begin{question}{\themyquestion}
Find the limit
\[\lim_{x\to \infty}\left(\frac{3x-1}{3x+1}\right)^{2x+1}.\]
\end{question}
\atc
\begin{question}{\themyquestion}
Let $r$ be a positive number. Prove that
\[\lim_{x\to 0^+}x^r\ln x=0.\]
\end{question}

\atc
\begin{question}{\themyquestion}
Determine whether the limit
 exists. If it exists, find the limit.
\begin{enumerate}[(a)]
\item
$\di \lim_{x\to 1}\frac{x^x-1}{x-1}$
\item $\di \lim_{x\to 0^+}\frac{x^x-1}{x\ln x}$
\item $\di \lim_{x\to 0^+}\frac{x^x\ln (1+x)}{x}$
\end{enumerate}
\end{question}

\vp

\section{Concavity of Functions}\label{sec3.7}

In this section, we study concavity of functions. If a function is twice differentiable, its concavity is determined by the second derivative.

Recall that if $x_1$ and $x_2$ are two points in $\mathbb{R}$, then as $t$ runs through all numbers in the interval $[0,1]$, 
\[x_1+t(x_2-x_1)=(1-t)x_1+tx_2\]runs through all points in the interval $[x_1, x_2]$. We say that a subset $S$ of $\mathbb{R}$ is convex if and only if for any two points $x_1$ and $x_2$ in $S$, and for any $t$ in $[0, 1]$, the point $(1-t)x_1+tx_2$ is also in $S$. We have proved that a subset of $\mathbb{R}$ is convex if and only if it is an interval.
\begin{definition}{Concavity of Functions}
Let $I$ be an interval.
\begin{enumerate}[1.]
\item A function $f:I\rightarrow \mathbb{R}$ is   concave up (or convex) provided that for any two points $x_1$ and $x_2$ in $I$, and for any $t\in [0,1]$,
    \[f((1-t)x_1+tx_2)\leq (1-t)f(x_1)+tf(x_2).\]
 \item A function $f:I\rightarrow \mathbb{R}$ is   concave down provided that for any  two points $x_1$ and $x_2$ in $I$, and for any $t\in [0,1]$,
   \[f((1-t)x_1+tx_2)\geq (1-t)f(x_1)+tf(x_2).\]
\item A function $f:I\rightarrow \mathbb{R}$ is strictly  concave up (or strictly convex)  provided that for any two distinct points $x_1$ and $x_2$ in $I$, and for any $t\in (0,1)$,
    \[f((1-t)x_1+tx_2)<(1-t)f(x_1)+tf(x_2).\]
 \item A function $f:I\rightarrow \mathbb{R}$ is strictly concave  down  provided that for any two distinct points $x_1$ and $x_2$ in $I$, and for any $t\in (0,1)$,
   \[f((1-t)x_1+tx_2)> (1-t)f(x_1)+tf(x_2).\]
\end{enumerate}
\end{definition}

Notice that a function $f:I\to\mathbb{R}$ is concave up if and only if the function $-f:I\to\mathbb{R}$ is concave down. Same for the strict concavity.

Geometrically, we draw a line   $L$ passing through the points $P_1(x_1, f(x_1))$  and $P_2(x_2, f(x_2))$ on the graph $y=f(x)$. If the equation of this line $L$ is $y=g(x)$, and $x_0=(1-t)x_1+tx_2$,  then 
\[g(x_0)= (1-t)f(x_1)+tf(x_2).\] Hence, $(x_0, g(x_0))$ is  point on the line $L$. Therefore, a function $y=f(x)$ is strictly concave up if its graph is always below a line segment joining two points on the graph; and it is strictly concave down if its graph is always above a line segment joining two points on the graph.

\begin{figure}[ht]
\centering
\includegraphics[scale=0.2]{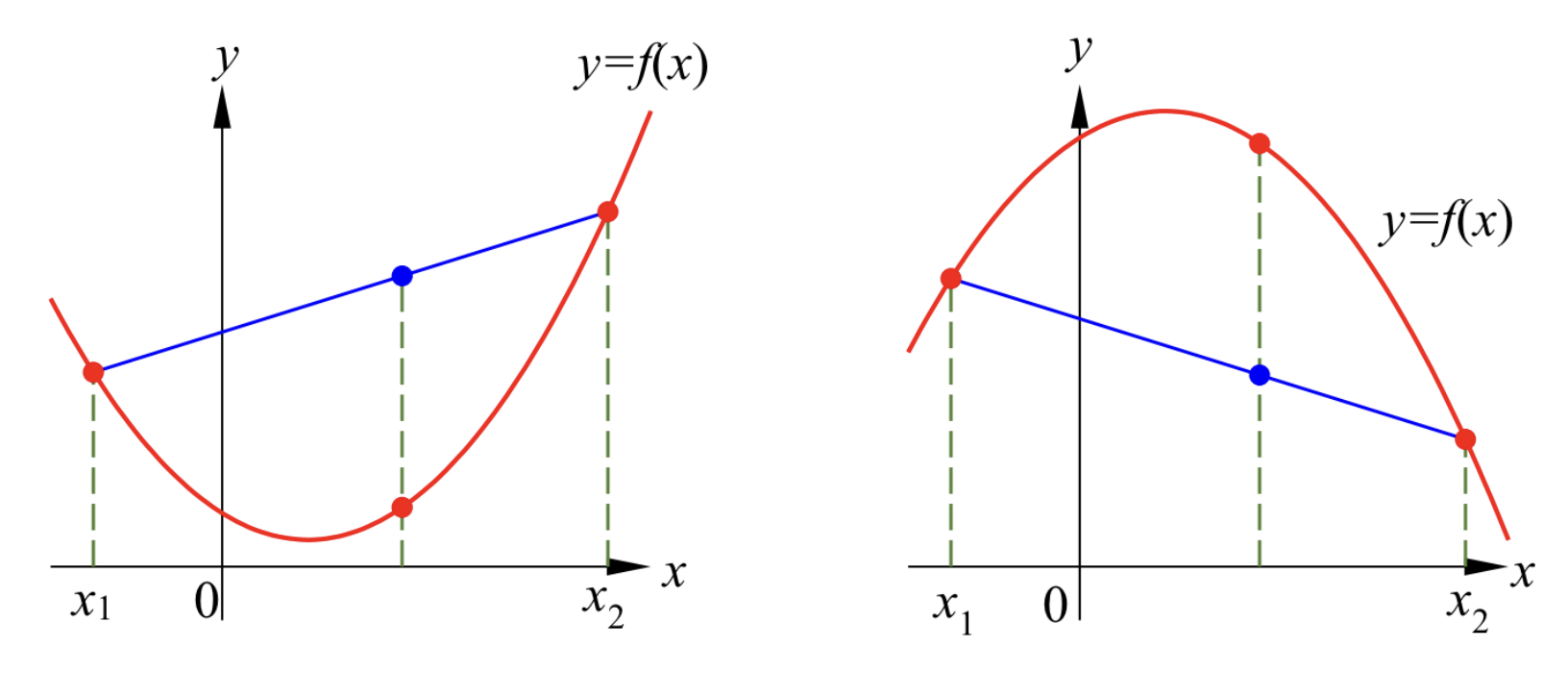}
\caption{(a) A strictly concave up function. (b)  A strictly concave down function.\fa}\label{figure35}
\end{figure}

\begin{example}{}
For any constants $m$ and $c$, the function $f(x)=mx+c$ is concave up and concave down. It is neither strictly concave up nor strictly concave down.
\end{example}

\begin{example}
{}
Show that the function $f:\mathbb{R}\to\mathbb{R}$, $f(x)=x^2$ is strictly concave up.
\end{example}
\begin{solution}{Solution}
Let $x_1$ and $x_2$ be any two distinct real numbers, and let $t$ be a number in the interval $(0,1)$. Then
\begin{align*}
&f((1-t)x_1+tx_2)-(1-t)f(x_1)-tf(x_2)\\&=(1-t)^2x_1^2+2t(1-t)x_1x_2+t^2x_2^2-(1-t)x_1^2-tx_2^2\\
&=-t(1-t)x_1^2+2t(1-t)x_1x_2-t(1-t)x_2^2\\
&=-t(1-t)(x_1-x_2)^2.
\end{align*}
Since   $(x_1-x_2)^2>0$, $t>0$ and $1-t>0$, we find that
\[f((1-t)x_1+tx_2)-(1-t)f(x_1)-tf(x_2)<0.\]
This proves that $f$ is strictly concave up.
\end{solution}

In the definition of concavity, we do not assume any regularity about the function. If a function is differentiable, we can characterize the concavity of the function in terms of its tangent lines.

For a point $x_0\in (a, b)$, the equation of the tangent line to the curve $y=f(x)$ at $x=x_0$ is
\[y=f(x_0)+f'(x_0)(x-x_0).\] 
We say that the graph of $f$ is   above the tangent line at $x=x_0$ provided that 
\[f(x)\geq f(x_0)+f'(x_0)(x-x_0)\hspace{1cm}\text{for all}\; x\in [a,b].\]
We say that the graph of $f$ is  strictly  above the tangent line at $x=x_0$ except at the tangential point provided that 
\[f(x)> f(x_0)+f'(x_0)(x-x_0)\hspace{1cm}\text{for all}\; x\in [a,b]\setminus \{x_0\}.\] 
Similarly, one can define what it means for the graph of $f$ to be below a tangent line, or strictly below.
\begin{theorem}[label=thm230219_1]{}
 Let  $f:[a, b]\rightarrow\mathbb{R}$  be a   function that is continuous on $[a,b]$, and   differentiable on $(a, b)$. The following three conditions are equivalent.

\begin{enumerate}[(a)]

\item $f'$ is strictly increasing on $(a, b)$.
\item The graph of $f$ is strictly  above every tangent line except at the tangential point.
\item $f$ is strictly concave up.
\end{enumerate}
 
\end{theorem}

\begin{myproof}{Proof}
First we prove (a) $\implies$ (b). Take any $x_0\in (a, b)$. The equation of the tangent line at $x=x_0$ is
\[y=g(x)=f(x_0)+f'(x_0)(x-x_0).\] If $x\in [a, b]$, 
\[f(x)-g(x)=f(x)-f(x_0)+f'(x_0)(x-x_0).\]
When $x\neq x_0$,  mean value theorem implies that there exists $u$ strictly between $x_0$ and $x$ such that
\[f(x)-f(x_0)=f'(u)(x-x_0).\]
Therefore
\[f(x)-g(x)=(x-x_0)(f'(u)-f'(x_0)).\]If $a\leq x<x_0$, $u<x_0$ and so $f'(u)<f'(x_0)$. This implies that $f(x)-g(x)>0$. If $x>x_0$, $u>x_0$ and so $f'(u)>f'(x_0)$. Then we also have $f(x)-g(x)>0$.  In other words, we have proved that for any $x_0$ in $(a, b)$, for any $x\in [a,b]\setminus\{x_0\}$,
\[f(x)>f(x_0)+f'(x_0)(x-x_0).\]\bp
This proves that the graph of $f$ is strictly  above every tangent line except at the tangential point.

Next, we prove (b) $\implies$ (c). Given $x_1$ and $x_2$ in $[a,b]$ with $x_1<x_2$, and $t\in (0,1)$, let $x_0=(1-t)x_1+tx_2$. Then $x_1<x_0<x_2$, and 
\[x_1-x_0=-t(x_2-x_1),\hspace{1cm}x_2-x_0=(1-t)(x_2-x_1).\]
By assumption,
\[f(x_1)>f(x_0)+f'(x_0)(x_1-x_0)=f(x_0)-tf'(x_0)(x_2-x_1),\]
\[f(x_2)>f(x_0)+f'(x_0)(x_2-x_0)=f(x_0)+(1-t)f'(x_0)(x_2-x_1).\]
Therefore,
\[(1-t)f(x_1)+tf(x_2)>f(x_0)=f((1-t)x_1+tx_2).\]
This proves that $f$ is strictly concave up.

Finally, we prove (c) $\implies $ (a). First we will prove that $f'$ is increasing on $(a, b)$. Given $x_1$ and $x_2$ in $(a, b)$ with $x_1<x_2$, we want to show that  $f'(x_1)\leq f'(x_2)$. For any $x\in (x_1, x_2)$, there exists $t\in (0,1)$ such that $x=(1-t)x_1+tx_2$. Since $f$ is strictly concave up, we have
\[f(x)=f((1-t)x_1+tx_2)<(1-t)f(x_1)+tf(x_2).\]
This implies that
\[(1-t)(f(x)-f(x_1))<t(f(x_2)-f(x)).\]
Since $x-x_1=t(x_2-x_1)$ and $x_2-x=(1-t)(x_2-x_1)$, we have
\begin{equation}\label{eq230219_1}\frac{f(x)-f(x_1)}{x-x_1}<\frac{f(x_2)-f(x)}{x_2-x}.\end{equation} Letting $x\to x_1^+$ in \eqref{eq230219_1}, we find that
\[f'(x_1)\leq \frac{f(x_2)-f(x_1)}{x_2-x_1}.\]\bp
Letting $x\to x_2^-$ in \eqref{eq230219_1}, we find that
\[f'(x_2)\geq \frac{f(x_2)-f(x_1)}{x_2-x_1}.\]These prove that
\[f'(x_2)\geq f'(x_1).\]
In other words, we have proved that $f'$ is increasing on $(a,b)$. If $f'$ is not strictly increasing on $(a, b)$, there exist $x_1$ and $x_2$ in $(a, b)$ with $x_1<x_2$ but $f'(x_1)=f'(x_2)$. Since $f'$ is increasing, we will have $f'(x)=f'(x_1)=f'(x_2)=m$ for all $x\in (x_1, x_2)$. This implies that 
\[f(x)=mx+c \hspace{1cm}\text{for all}\; x\in [x_1, x_2].\]
But then for any $t\in (0,1)$,
\[f((1-t)x_1+tx_2)=m\left[(1-t)x_1+tx_2\right]+b=(1-t)f(x_1)+tf(x_2),\]
which constradicts to the strict concavity of $f$. Therefore, $f'$ must be strictly increasing on $(a,b)$.
\end{myproof}

By replacing $f$ by $-f$ in Theorem \ref{thm230219_1}, we obtain the following immediately.
\begin{theorem}[label=thm230219_2]{}
 Let  $f:[a, b]\rightarrow\mathbb{R}$  be a   function that is continuous on $[a,b]$, and   differentiable on $(a, b)$. The following three conditions are equivalent.
\begin{enumerate}[(a)]

\item $f'$ is strictly decreasing on $(a, b)$.
\item The graph of $f$ is strictly below every tangent line except at the tangential point.
\item $f$ is strictly concave down.
\end{enumerate}
 
\end{theorem}

In Theorem \ref{thm230219_1}, if we relax the strictness, the proofs are actually easier.  
\begin{theorem}[label=thm230219_3]{}
 Let  $f:[a, b]\rightarrow\mathbb{R}$  be a   function that is continuous on $[a,b]$, and   differentiable on $(a, b)$. The following three conditions are equivalent.
\begin{enumerate}[(a)]

\item $f'$ is   increasing on $(a, b)$.
\item The graph of $f$ is   above every tangent line.
\item $f$ is  concave up.
\end{enumerate}
 
\end{theorem}

\begin{theorem}[label=thm230219_4]{}
 Let  $f:[a, b]\rightarrow\mathbb{R}$  be a   function that is continuous on $[a,b]$, and   differentiable on $(a, b)$. The following three conditions are equivalent.
\begin{enumerate}[(a)]

\item $f'$ is  decreasing on $(a, b)$.
\item The graph of $f$ is  below every tangent line.
\item $f$ is  concave down.
\end{enumerate}
 
\end{theorem}

If a function is twice differentiable, we can characterize concavity using second derivatives.
\begin{theorem}{}
 Let  $f:[a, b]\rightarrow\mathbb{R}$ be a function that is continuous on $[a,b]$, and twice differentiable on $(a, b)$.
 
\begin{enumerate}[1.]
\item $f(x)$ is   concave up if and only if $f''(x)\geq 0$ for all $x\in (a, b)$.
\item $f(x)$ is   concave down if and only if $f''(x)\leq 0$ for all $x\in (a, b)$.
\item If $f''(x)>0$ for all $x\in (a, b)$,  then $f$ is strictly concave up.

\item  If $f''(x)<0$ for all $x\in (a, b)$, then  $f$ is strictly concave down.

\end{enumerate}
\end{theorem}
\begin{myproof}{Proof}
For (a) and (b),  this is just the fact that $f''(x)\geq 0$ for all $x\in (a, b)$ if and only if $f'$ is increasing;  $f''(x)\leq 0$ for all $x\in (a, b)$ if and only if $f'$ is decreasing.

For (c) and (d), we note that $f''(x)>0$ for all $x\in (a,b)$ implies that $f'$ is strictly increasing on $(a,b)$; while  $f''(x)<0$ for all $x\in (a,b)$ implies that $f'$ is strictly decreasing on $(a,b)$
\end{myproof}
We have seen that if a differentiable function $g:(a,b)\to\mathbb{R}$ is strictly increasing, it is not necessary that $g'(x)>0$ for all $x\in (a,b)$. This is why for $f:(a,b)\to\mathbb{R}$ to be strictly concave up, it is not necessary that $f''(x)>0$ for all $x\in (a,b)$. The function $f:\mathbb{R}\to\mathbb{R}$, $f(x)=x^4$ gives an example of a function that is strictly concave up, but it is not true that $f''(x)>0$ for all $x\in\mathbb{R}$. 
\begin{example}{}
Show that the function $f:[0, \pi]\to\mathbb{R}$, $f(x)=\sin x$ is strictly concave down.
\end{example}
\begin{solution}{Solution}
Since $f''(x)=-\sin x<0$ for all $x\in (0, \pi)$, we find that $f:[0, \pi]\to\mathbb{R}$, $f(x)=\sin x$ is strictly concave down.
\end{solution}
\begin{example}
{}Consider the power function $f(x)=x^r$. Since $f''(x)=r(r-1)x^{r-2}$, we have the following.
\begin{enumerate}[1.]
\item When $r<0$,  the function $f:(0,\infty)\to\mathbb{R}$, $f(x)=x^r$ is strictly concave up. 
\item When $0<r<1$, the function $f:[0,\infty)\to\mathbb{R}$, $f(x)=x^r$ is strictly concave down. 
\item When $r>1$, the  function $f:[0,\infty)\to\mathbb{R}$, $f(x)=x^r$ is strictly concave up.\end{enumerate}
\end{example}

Next we look at a classical example where the concavity of a function can be used to prove inequalities.
\begin{example}{Young's Inequality}
Given that $p$ and $q$ are positive numbers such that
\[\frac{1}{p}+\frac{1}{q}=1.\]If $a$ and $b$ are positive numbers, show that 
\[ab\leq \frac{a^p}{p}+\frac{b^q}{q}.\]Equality holds if and only if $a^p=b^q$.
\end{example}
\begin{solution}{Solution}Notice that since $p$ and $q$ are positive, we have $1/p<1$ and $1/q<1$. This implies that $p>1$ and $q>1$.
Consider the function $f:(0,\infty)\to\mathbb{R}$, $f(x)=\ln x$. We find that
\[f'(x)=\frac{1}{x},\quad f''(x)=-\frac{1}{x^2}<0\hspace{1cm}\text{for all}\; x>0.\]  Hence, the function $f:(0,\infty)\to\mathbb{R}$, $f(x)=\ln x$ is strictly concave down. 
This implies that for any two   positive numbers $x_1$ and $x_2$, and any $t\in (0,1)$, 
\begin{equation}\label{eq230219_5}\ln \left((1-t)x_1+tx_2\right)\geq (1-t)\ln x_1+t\ln x_2.\end{equation} The equality can hold if and only if $x_1=x_2$.   Now, let
$t=1/q$.
Then $t\in (0,1)$ and $1-t=1/p$. Given the positive numbers $a$ and $b$, let
\[x_1=a^p, \quad x_2=b^q.\]
 Then $x_1$ and $x_2$ are positive numbers. Eq. \eqref{eq230219_5} implies that
\[\ln \left(\frac{a^p}{p}+\frac{b^q}{q}\right)\geq \frac{1}{p}\ln a^p+\frac{1}{q}\ln b^q=\ln (ab),\]\bs with equality holds if and only if $a^p=b^q$. Therefore,
\[ab\leq \frac{a^p}{p}+\frac{b^q}{q}.\]Equality holds if and only if $a^p=b^q$.
\end{solution}

\vp
\noindent
{\bf \large Exercises  \thesection}
\setcounter{myquestion}{1}
\begin{question}{\themyquestion}
\begin{enumerate}[(a)]
\item Show that the function $f:\mathbb{R}\to \mathbb{R}$, $f(x)=e^{-x}$ is strictly concave up.
\item Show that the function $f:(0,\pi)\to \mathbb{R}$, $f(x)=\di\frac{1}{\sin x}$ is strictly concave up.
\end{enumerate}
\end{question}

 \atc
\begin{question}{\themyquestion}
Let $f:(a,b)\to\mathbb{R}$ be a twice differentiable function. If $f:(a,b)\to\mathbb{R}$ is concave down, and $f(x)>0$ for all $x\in (a, b)$, prove that the function $g:(a,b)\to\mathbb{R}$,
\[g(x)=\frac{1}{f(x)}\] is concave up.
\end{question}

\atc
\begin{question}{\themyquestion}
Given that the function $f:[a,b]\to\mathbb{R}$ is concave down. Show that for any  $x_1, x_2, \ldots, x_n$ in $[a, b]$, if  $t_1, t_2, \ldots, t_n$ are nonnegative numbers satifying \[t_1+t_2+\ldots+t_n=1,\]then
\[f(t_1x_1+t_2x_2+\cdots+t_nx_n)\geq t_1f(x_1)+t_2f(x_2)+\cdots+t_nf(x_n).\]
 
\end{question}

\atc
\begin{question}{\themyquestion: Arithmetic Mean-Geometric Mean Inequality}
The arithmetic mean-geometric mean inequality states that if $a_1, a_2, \ldots, a_n$ are $n$ positive numbers, then
\[\frac{a_1+a_2+\cdots+a_n}{n}\geq \sqrt[n]{a_1a_2\cdots a_n}.\]
Use the concavity of the function $f:(0,\infty)\to\mathbb{R}$, $f(x)=\ln x$ and the result of the previous question to prove this inequality.
 
\end{question}

\chapter{Integrating Functions of a Single Variable}\label{ch4}

The concept of integrals arises naturally when one wants to compute the area bounded by a curve, such as the area of   a circle. Since the ancient time, our ancestors have found a good strategy to deal with such problems. For example, they used the area of polygons to approximate the area of a circle. The circle is partitioned into  sectors, and the area of each sector is approximated by the area of the inscribed triangle (see Figure \ref{figure36}). When the circle is partitioned into more sectors, better approximation is obtained.

 \begin{figure}[ht]
\centering
\includegraphics[scale=0.2]{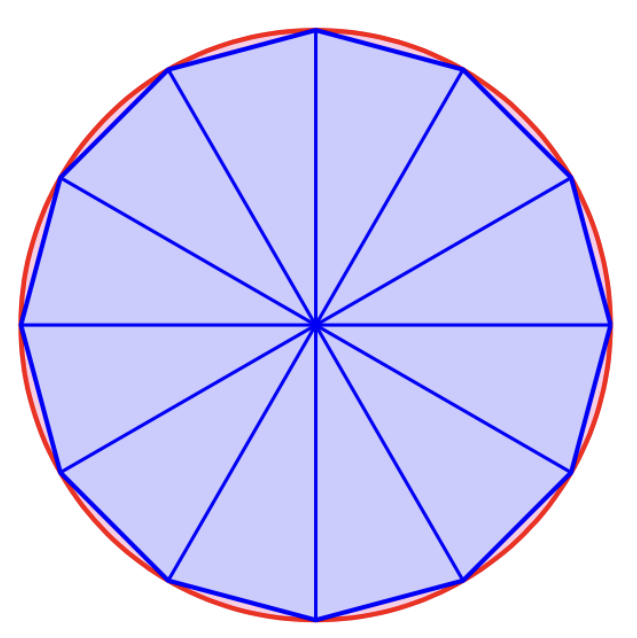}
\caption{ Approximating the area of a circle by the area of a polygon.\fa}\label{figure36}
\end{figure}

The same idea can be used to find the area enclosed by any curves. This motivated the definition of integrals. For curves defined by continuous functions, it is not difficult to formulate a well-defined definition for integrals. However, mathematicians soon discovered that we need to work with functions that are not continuous as well. The process to make integrals rigorously defined is long and tedious. We will follow the historical path and study the Riemann integrals in this course. This lays down the foundation for advanced theory of integration $\grave{\text{a}}$ la Lebesgue. For practical applications and computations, Riemann integrals are sufficient and easier to calculate.

\section{Riemann Integrals of Bounded Functions }\label{sec4.1}

In this section, we  define the Riemann integral for a function $f:[a,b]\to\mathbb{R}$ that is defined on a closed and bounded interval $[a,b]$. For this purpose, the function is necessarily bounded. In Section \ref{sec4.6}, we will discuss how to deal with functions that are not necessarily bounded, via some limiting processes. 

For  a closed and bounded interval $[a,b]$, we will always assume that $a<b$.
We start by a few definitions.
\begin{definition}{Partitions}
Let $[a, b]$ be a closed and bounded interval. A partition of $[a,b]$ is a finite sequence of points $x_0, x_1, x_2, \ldots, x_k$, where
\[a=x_0<x_1<\cdots<x_{k-1}<x_k=b.\]
 It is denoted by $P=\{x_0, x_1, \ldots, x_k\}$. For each $0\leq i\leq k$, $x_i$ is a partition point. These points partition the interval $[a, b]$ into $k$ subintervals $[x_0, x_1]$, $[x_1, x_2]$, $\ldots$, $[x_{k-1}, x_k]$. The $i^{\text{th}}$-subinterval is $[x_{i-1}, x_i]$.
\end{definition}We have slightly abused notation and used set notation for a partition.  
\begin{example}{}
$P=\{0, 2, 3, 5, 9, 10\}$ is a partition of the interval $[0,10]$ into 5 subintervals 
\[[0, 2], [2, 3], [3, 5], [5, 9]\;\text{and}\;[9,10].\]
\end{example}
We use the lengths of the subintervals to measure how fine a partition is. 
\begin{definition}{Gap of a Partition}
Let $P=\{x_0, x_1, \ldots, x_k\}$ be a partition of the interval $[a,b]$. The gap of the partition $P$, denoted by $|P|$ or $\text{gap}\,P$, is the length of the largest subinterval in the partition. Namely,
\[|P|=\text{gap}\, P=\max\left\{x_i-x_{i-1}\,|\, 1\leq i\leq k\right\}.\]
\end{definition}
\begin{example}{}For the partition $P=\{0, 2, 3, 5, 9, 10\}$ of $[0, 10]$,
\[|P|=\max\{2,1,2,4,1\}=4.\]
\end{example}
A partition where all subintervals have equal lengths is very useful.
\begin{definition}{Regular Partitions}
Let $[a, b]$ be a closed and bounded interval. A regular partition of $[a,b]$ into $k$ intervals is the partition 
$P=\{x_0, x_1, \ldots, x_k\}$, where 
\[|P|=x_1-x_0=x_2-x_1=\cdots=x_k-x_{k-1}=\frac{b-a}{k}.\]This implies that
$\di x_i=x_0+i\frac{b-a}{k}$, $1\leq i\leq k$.

\end{definition}
 \begin{example}{}
The  regular partition of the interval $[0,10]$ into 5 intervals is the partition 
\[P=\{0, 2, 4, 6, 8, 10\}.\]The gap of this partition is
$\di |P|=\frac{10-0}{5}=2$.
\end{example}

Next, we define the Riemann sums and Darboux sums.

\begin{definition}{Riemann Sums}
Let $f:[a,b]\to \mathbb{R}$ be a   function, and let $P=\{x_0, x_1, \ldots, x_k\}$ be a partition of $[a, b]$. For each $1\leq i\leq k$, choose an intermediate point $\xi_i$ in the $i^{\text{th}}$-subinterval $[x_{i-1}, x_i]$. Denote this sequence of points    $\{\xi_i\}_{i=1}^k$ by $A$. Then the Riemann sum of $f$ with respect to the partition $P$ and the intermediate points $A=\{\xi_i\}_{i=1}^k$ is the sum
\[R(f, P, A)=\sum_{i=1}^kf(\xi_i)(x_i-x_{i-1}).\]
\end{definition}

\begin{example}[label=ex230220_1]{}
Consider the function $f:[0, 6]\to\mathbb{R}$, $f(x)=6x-x^2$, and the partition  $P=\{0, 2, 3, 5, 6\}$ of $[0, 6]$. Let 
\[A=\left\{1, 3, 4, 5\right\}.\]Then 
\[R(f, P, A)=5\times 2+9\times 1+8\times 2+5\times   1=40.\]
\end{example}
As shown in Figure \ref{figure37}, the Riemann sum $R(f, P, A)$ is the sum of the areas of rectangles that are used to approximate the region bounded by the curve $y=6x-x^2$ and the $x$-axis. 

 \begin{figure}[ht]
\centering
\includegraphics[scale=0.2]{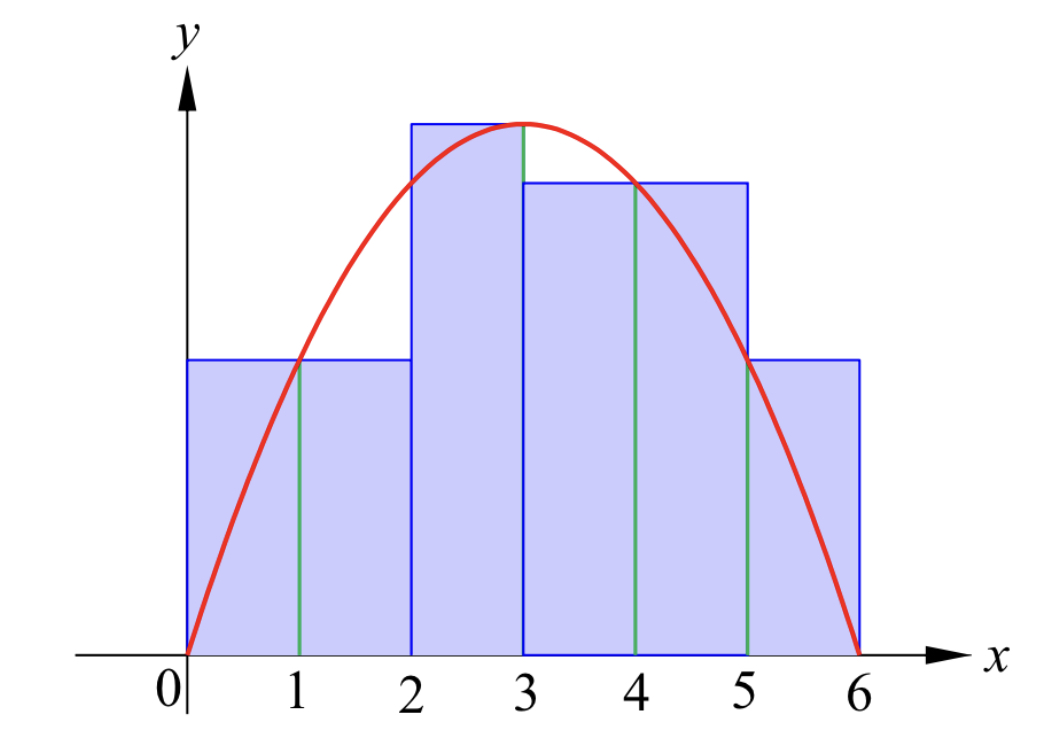}
\caption{Riemann sum is an approximation of area under a curve.\fa}\label{figure37}
\end{figure}

In general, if $f:[a,b]\to\mathbb{R}$ is a nonnegative function, then a Riemann sum $R(f, P, A)$ is an approximation to the area bounded by the curve $y=f(x)$, the $x$-axis, and the lines $x=a$ and $x=b$. In its definition, we do not need to assume that $f$ is a bounded function. 
Since Riemann sum involves an arbitrary choice of points in each subinterval, to give a bound to Riemann sums, we need the concept of Darboux sums, whose definition requires $f:[a,b]\to\mathbb{R}$ to be a bounded function.

\begin{definition}{Darboux Sums}
Let $f:[a,b]\to \mathbb{R}$ be a bounded function, and let $P=\{x_0, x_1, \ldots, x_k\}$ be a partition of $[a, b]$. For each $1\leq i\leq k$, let
\[m_i =\inf_{x_{i-1}\leq x\leq x_i}f(x),\hspace{1cm}
M_i =\sup_{x_{i-1}\leq x\leq x_i}f(x).\]The Darboux lower sum $L(f,P)$ and the Darboux upper sum $U(f,P)$ are defined by
\begin{align*}
L(f,P)&=\sum_{i=1}^km_i(x_i-x_{i-1}),\\
U(f,P)&=\sum_{i=1}^kM_i(x_i-x_{i-1}).
\end{align*}
\end{definition}
\begin{remark}{}
For convenience, we   denote \[ \inf \{f(x)\,|\,x_{i-1}\leq x\leq x_i\}\quad\text{and}\quad\sup\{f(x)\,|\,x_{i-1}\leq x\leq x_i\}\] by
$\di \inf_{x_{i-1}\leq x\leq x}f(x)$ and $\di \sup_{x_{i-1}\leq x\leq x_i}f(x)$ respectively. The assumption that $f$ is bounded is needed to ensure that $m_i$ and $M_i$ exist for all $1\leq i\leq k$. The reason we use infimum and supremum is obvious, as the function $f$ might not have minimum or maximum on an interval.
\end{remark}
 \begin{figure}[ht]
\centering
\includegraphics[scale=0.2]{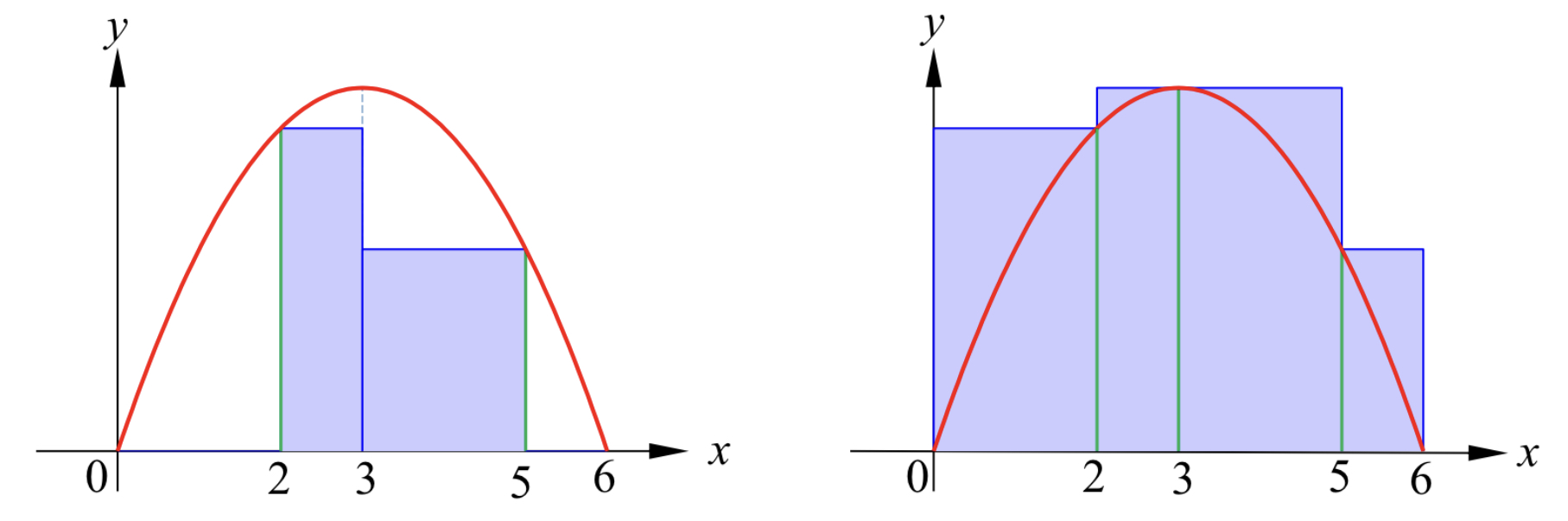}
\caption{Darboux lower sum and Darboux upper sum.\fa}\label{figure38}
\end{figure}
\begin{example}[label=ex230220_2]{}
For the   function $f:[0, 6]\to\mathbb{R}$, $f(x)=6x-x^2$  and the partition  $P=\{0, 2, 3, 5, 6\}$ of $[0, 6]$ considered in Example \ref{ex230220_1}, we find that

\vspace{0.5cm} ~\hspace{1.5cm}
\begin{tabular}{||c||c|c|c|c||}
\hline
\hline
$i$ &$\;$ interval$\;$ & $x_i-x_{i-1}$ & $m_i$ & $M_i$\\
\hline
\hline
$\;\;1\;\;$ & $[0,2]$ & 2&$\quad 0\quad $ & $\quad 8\quad $\\
\hline
2 & $[2,3]$ &1& 8 & 9\\
\hline
3 & $[3, 5]$ &2&  5&9\\
\hline 
4 & $[5, 6]$ &1& 0 & 5\\
\hline

\hline

\end{tabular}

\vspace{0.5cm}
Hence, the Darboux lower sum $L(f,P)$ and the Darboux upper sum $U(f,P)$ are
\begin{align*}
L(f,P)&=0\times 2+8\times 1+5\times 2+0\times   1=18,\\
U(f,P)&=8\times 2+9\times 1+9\times 2+5\times   1= 48.
\end{align*}
\end{example}

\begin{example}{}
If $f:[a,b]\to\mathbb{R}$ is the constant function $f(x)=c$, it is obvious that for any partition $P=\{x_i\}_{i=0}^k$ of $[a,b]$, and for any choices of intermediate points $A=\{\xi_i\}_{i=1}^k$,
\[L(f,P)=U(f,P)=R(f,P,A)=c(b-a).\]
\end{example}

The following can be easily deduced  from the definitions.
\begin{proposition}[label=230304_6]{}
Let $f:[a,b]\to\mathbb{R}$ be a bounded function such that
\[m\leq f(x)\leq M\hspace{1cm}\text{for all}\;a\leq x\leq b.\] For any partition $P=\{x_i\}_{i=0}^k$ of the interval $[a,b]$, and any  choice of intermediate points  $A=\{\xi_i\}_{i=1}^k$ for the partition $P$, we have
\[m(b-a)\leq L(f,P)\leq R(f,P,A)\leq U(f,P)\leq M(b-a).\]
\end{proposition}
\begin{myproof}{Proof}
For any $1\leq i\leq k$, let
\[m_i=\inf_{x_{i-1}\leq x\leq x}f(x),\hspace{1cm}
M_i =\sup_{x_{i-1}\leq x\leq x_i}f(x).\]Then 
\[m\leq m_i\leq f(\xi_i)\leq M_i\leq M.\] Therefore,
\begin{align*}
m(x_i-x_{i-1}) \leq m_i(x_i-x_{i-1})&\leq f(\xi_i)(x_i-x_{i-1})\\&\leq M_i(x_i-x_{i-1})\leq M(x_i-x_{i-1}).\end{align*}
Summing over $i$ from $i=1$ to $i=k$, we obtain
\[m(b-a)\leq L(f,P)\leq R(f,P,A)\leq U(f,P)\leq M(b-a).\]
\end{myproof}

For a  bounded nonnegative   function $f:[a,b]\to \mathbb{R}$,  if the region bounded by the $x$-axis, the curve $y=f(x)$, the lines $x=a$ and $x=b$ has an area, then a Darboux lower sum  is always less than or  equal to the area, and a Darboux upper sum   is always larger than or equal to the area. This leads to the fact that a Darboux lower sum is always less than or equal to a Darboux upper sum. To prove this for any bounded functions, we introduce the concept of refinement.

\begin{definition}{Refinement of a Partition}
Let $P$ and $P^*$ be partitions of the interval $[a,b]$. We say that $P^*$ is refinement of $P$ if every partition point of $P$ is also a partition point of $P^*$. In other words, the set of points in $P$ is a subset of the set of points in $P^*$.
\end{definition}

 \begin{example}[label=ex230220_3]{}
For the partition $P=\{0, 2, 3, 5, 9, 10\}$ of $[0,10]$, \[P^*=\{0, 1, 2, 3, 5, 6, 8, 9, 10\}\] is a refinement.
\end{example}
 \begin{figure}[ht]
\centering
\includegraphics[scale=0.2]{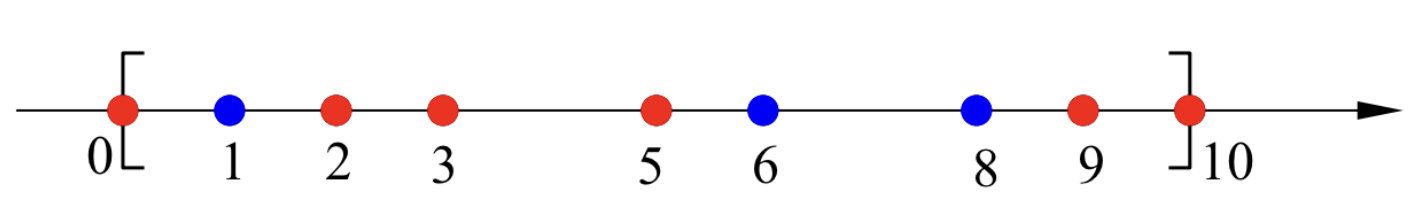}
\caption{A partition $P$ of $[0,10]$ and its refinement $P^*$.\fa}\label{figure39}
\end{figure}

If $P^*$ is a refinement of $P=\{x_i\}_{i=1}^k$, then for each $1\leq i\leq k$, $P^*$ induces a partition $P_i$ of the interval $[x_{i-1}, x_i]$. 

 \begin{example}[label=ex230220_4]{}
For the partition $P$ and $P^*$ in Example \ref{ex230220_3}, $P^*$ induces the partition $P_1=\{0,1,2\}$, $P_2=\{2,3\}$, $P_3=\{3,5\}$, $P_4=\{5,6,8,9\}$ and $P_5=\{9,10\}$ of the intervals $[0,2]$, $[2,3]$,  $[3,5]$, $[5,9]$ and $[9,10]$ respectively. 
\end{example}

Since the union of all the subintervals in the partition $P_i$, $1\leq i\leq k$ is the collection of all the subintervals in the partition $P^*$, the following is quite obvious.
\begin{proposition}[label=230220_5]{}
Let $f:[a,b]\to\mathbb{R}$ be a bounded function, and let $P=\{x_i\}_{i=0}^k$ be a partition of $[a,b]$. Given a refinement $P^*$ of $P$,  let $P_i$, $1\leq i\leq k$, be the partition that $P^*$ induces on the interval $[x_{i-1}, x_i]$. Then
\[\sum_{i=1}^kL(f,P_i)=L(f,P^*),\hspace{1.5cm} \sum_{i=1}^kU(f,P_i)=U(f,P^*).\]
\end{proposition}

From this, it is quite easy to obtain the following. 
\begin{theorem}[label=230220_6]{}
Let $f:[a,b]\to\mathbb{R}$ be a bounded function, and let $P$ and $P^*$ be   partitions of $[a,b]$. If $P^*$ is a refinement of $P$, then
\[L(f,P)\leq L(f,P^*)\leq U(f,P^*)\leq U(f,P).\]

\end{theorem}
\begin{myproof}{Proof}Let $P=\{x_i\}_{i=0}^k$. For each $1\leq i\leq k$, 
let
\[m_i=\inf_{x_{i-1}\leq x\leq x_i}f(x),\hspace{1cm}M_i=\sup_{x_{i-1}\leq x\leq x_i}f(x).\]
Since 
\[m_i\leq f(x)\leq M_i\hspace{1cm}\text{for all}\;x\in [x_{i-1}, x_i],\]
we find that
\[m_i(x_i-x_{i-1})\leq L(f, P_i)\leq U(f,P_i)\leq M_i(x_i-x_{i-1}).\]Summing over $i$ from $1$ to $k$ gives
\[\sum_{i=1}^{k}m_i(x_i-x_{i-1})\leq\sum_{i=1}^{k} L(f, P_i)\leq \sum_{i=1}^{k}U(f,P_i)\leq \sum_{i=1}^{k}M_i(x_i-x_{i-1}).\]By Proposition \ref{230220_5}, this gives
\[L(f,P)\leq L(f,P^*)\leq U(f,P^*)\leq U(f,P).\]
\end{myproof}

\begin{figure}[ht]
\centering
\includegraphics[scale=0.2]{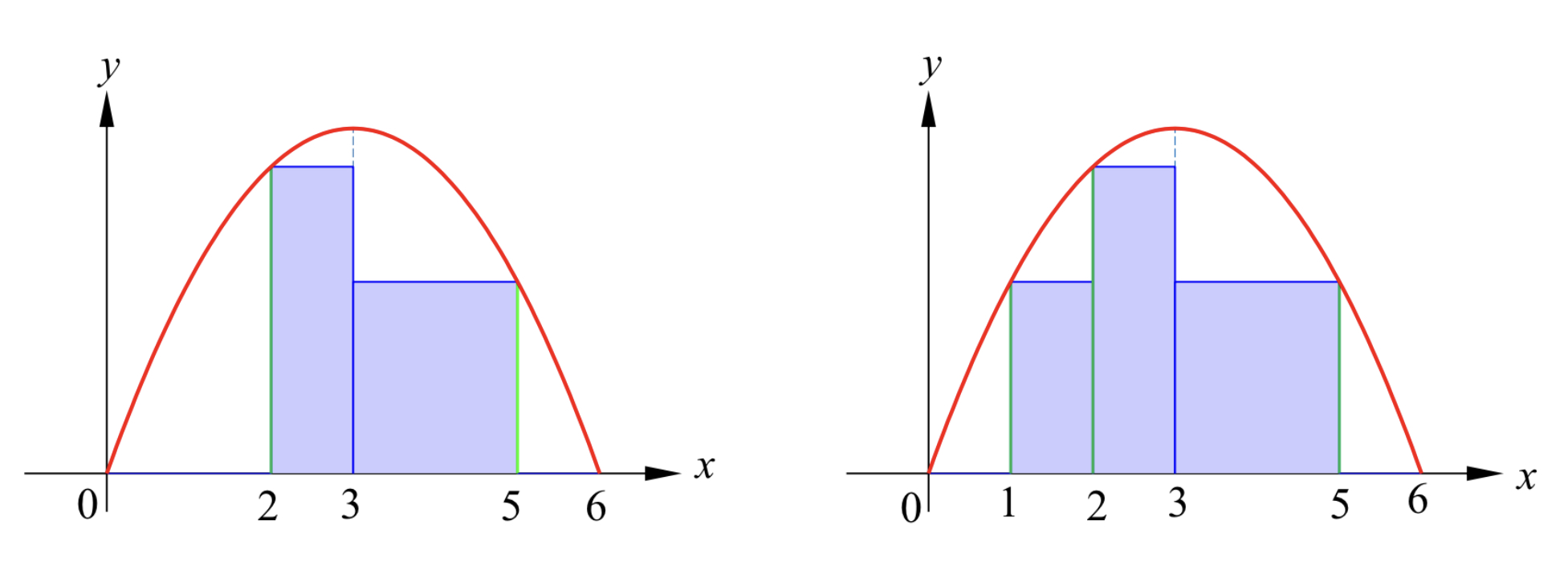}
\caption{When a partition is refined, Darboux lower sum gets larger.\fa}\label{figure40}
\end{figure}

\begin{figure}[ht]
\centering
\includegraphics[scale=0.2]{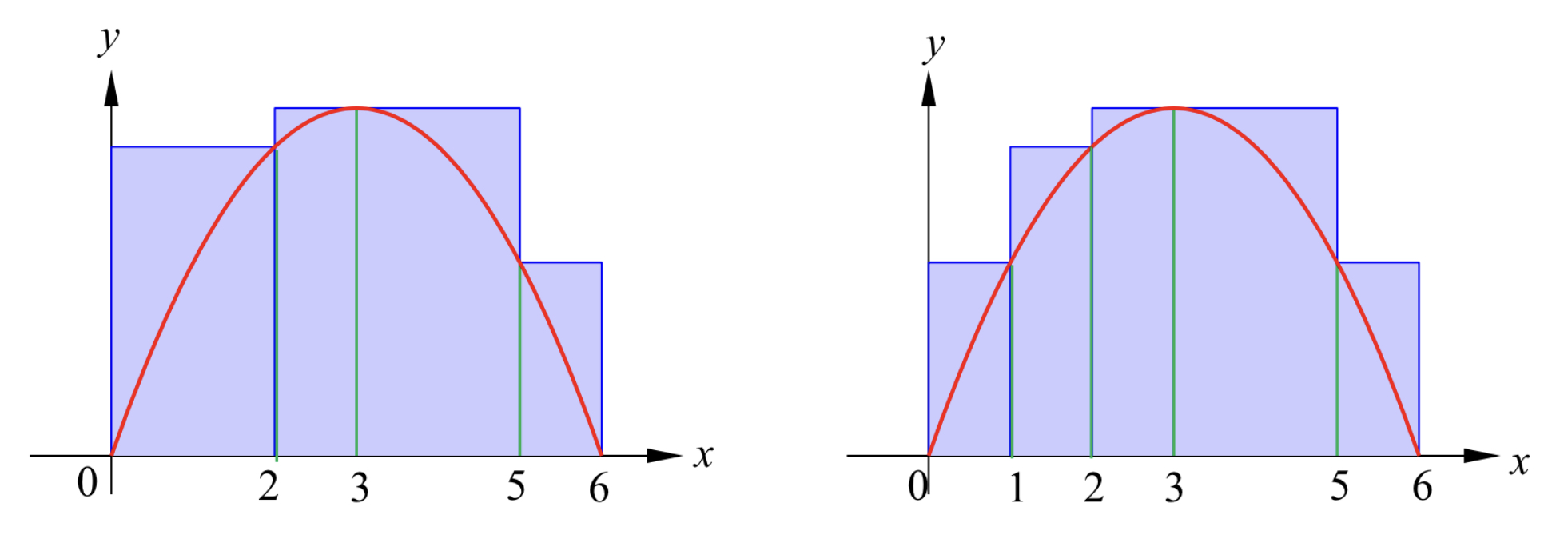}
\caption{When a partition is refined, Darboux upper sum gets smaller.\fa}\label{figure41}
\end{figure}

If $P_1$ and $P_2$ are partitions of $[a,b]$, a common refinement of $P_1$ and $P_2$ is a partition $P^*$ which contains all the partition points of $P_1$ and $P_2$. Such a common refinement always exists. The smallest one is the one whose set of points is the union of the set of points in $P_1$ and the set of points in $P_2$.
\begin{corollary}[label=230220_7]{}
Let $f:[a,b]\to\mathbb{R}$ be a bounded function, and let $P$ and $P_2$ be any  two partitions of $[a,b]$. Then
\[L(f, P_1)\leq U(f, P_2).\]
\end{corollary}
\begin{myproof}{Proof}
Take a common refinement $P^*$ of the partitions $P_1$ and $P_2$. By Theorem \ref{230220_6},
\[L(f,P_1)\leq L(f,P^*)\leq U(f,P^*)\leq U(f,P_2).\]
\end{myproof}

Given a bounded function $f:[a,b]\to\mathbb{R}$, we consider the set of Darboux lower sums and the set of Darboux upper sums of $f$.
\begin{align*}
S_L(f)&=\left\{L(f,P)\,|\, P\;\text{is a partition of}\;[a,b]\right\},\\
S_U(f)&=\left\{U(f,P)\,|\, P\;\text{is a partition of}\;[a,b]\right\}.
\end{align*}If $m$ and $M$ are lower and upper bounds of $f$, then 
\[m(b-a)\leq L(f,P)\leq U(f,P)\leq M(b-a).\]
This implies that the sets $S_L(f)$ and $S_U(f)$ are bounded. When we use the Darboux lower sums and upper sums to approximate areas, we are interested in the least upper bound of the lower sums and the greatest lower bound of the upper sums.
\begin{definition}{Lower Integrals and Upper Integrals}
Let $f:[a,b]\to\mathbb{R}$ be a bounded function.
\begin{enumerate}[1.]
\item
The lower integral of $f$, denoted by $\di\underline{\int_a^b}f$, is defined as the least upper bound of the Darboux lower sums.
\[\underline{\int_a^b}f=\sup S_L(f)=\sup \left\{L(f,P)\,|\, P\;\text{is a partition of}\;[a,b]\right\}.\]
\item
The upper integral of $f$, denoted by $\di\overline{\int_a^b}f$, is defined as the greatest lower bound of the Darboux upper sums.
\[\overline{\int_a^b}f=\inf S_U(f)=\inf \left\{U(f,P)\,|\, P\;\text{is a partition of}\;[a,b]\right\}.\]
\end{enumerate}
\end{definition}

\begin{example}{}
For the constant function $f:[a.,b]\to\mathbb{R}$, $f(x)=c$, 
\[L(f,P)=U(f,P)=c(b-a)\] for any partition $P$ of $[a,b]$. Thus,
$S_L(f)=S_U(f)=\{c(b-a)\}$, and
\[\underline{\int_a^b}f=\overline{\int_a^b}f=c(b-a).\]
\end{example}

By Corollary \ref{230220_7}, we have the following.
\begin{proposition}{}
Let $f:[a,b]\to\mathbb{R}$ be a bounded function. We have
\[\underline{\int_a^b}f\leq \overline{\int_a^b}f.\]
\end{proposition}
\begin{myproof}{Proof}
By definitions of infimum and supremum, for any positive integer $n$, there are partition $P_1$ and $P_2$ such that
\[L(f,P_1)>\underline{\int_a^b}f-\frac{1}{n},\hspace{1cm}U(f,P_2)< \overline{\int_a^b}f+\frac{1}{n}.\]
These, together with Corollary \ref{230220_7}, give the following.
\[ \overline{\int_a^b}f-\underline{\int_a^b}f>U(f,P_2)-L(f,P_1)-\frac{2}{n}\geq -\frac{2}{n}.\]
Taking the limit $n\to\infty$, we deduce that
\[\overline{\int_a^b}f-\underline{\int_a^b}f\geq 0.\]
\end{myproof}
\begin{highlight}{}
Let $f:[a,b]\to\mathbb{R}$ be a bounded function, and let  $P$ be a partition of $[a,b]$.
\begin{enumerate}[1.]
\item $\di L(f,P)\;\leq \;\underline{\int_a^b} f\;\leq \;\overline{\int_a^b }f \;\leq U(f,P)$.
\item $\di 0\;\leq \;\overline{\int_a^b }f \;-\;\underline{\int_a^b} f\;\leq\;U(f,P)\;-\;L(f,P)$.

\end{enumerate}
\end{highlight}
Notice that if $w$ is a number such that
\[\underline{\int_a^b}f\leq w\leq \overline{\int_a^b}f,\]then
for any partition $P$ of $[a,b]$,
\[L(f,P)\leq w\leq U(f,P).\]
As we mentioned above, for a nonnegative bounded function $f:[a,b]\to\mathbb{R}$, a Darboux lower sum $L(f,P)$ is less  than or equal to the area below the curve $y=f(x)$, while a Darboux upper sum is larger than or equal to the area, if such an area is well-defined. Intuitively, the area would be well-defined if there is a single number $A$ that is larger than or equal to all the Darboux lower sums, and less than or equal to all the Darboux upper sums. This is the case if and only if the lower and the upper integrals are the same. 
\begin{definition}{Riemann Integrability}
Let $f:[a,b]\to \mathbb{R}$ be a bounded function. We say that $f$ is Riemann integrable, or simply integrable, if 
\[\underline{\int_a^b}f= \overline{\int_a^b}f.\]In this case, we define the integral of $f$ over $[a,b]$ as
\[\int_a^b f=\underline{\int_a^b}f=\overline{\int_a^b}f.\]
It is the unique number  that is larger than or equal to all the Darboux lower sums, and less than or equal to all the Darboux upper sums.
\end{definition}

\begin{remark}{}
If $f:[a,b]\to \mathbb{R}$ is a continuous nonnegative function, we are going to prove that $f$ is Riemann integrable.  It follows from our discussions above that 
the integrable $\di\int_a^b f$ is the area bounded by the curve $y=f(x)$, the $x$-axis, and the lines $x=a$ and $x=b$.
\end{remark}

\begin{highlight}{Leibniz Notation}
In Leibniz notation, the integral of $f:[a,b]\to\mathbb{R}$ over $[a,b]$ is denoted by \[\int_a^b f(x)dx.\]
\end{highlight}
\begin{example}{}
A constant functon $f:[a,b]\to\mathbb{R}$, $f(x)=c$  is integrable and 
\[\int_a^b f =c(b-a).\]
\end{example}

Let us look at an example of a function that is not integrable.
 \begin{example}[label=230220_9]
 {Non-Integrability of the Dirichlet's Function} The Dirichlet's function is the function $f:[0,1]\rightarrow\mathbb{R}$ defined by
 \begin{align*}
 f(x)=\begin{cases}1,\quad &\text{if}\; x\;\text{is rational},
 \\0,\quad &\text{if}\; x\;\text{is irrational}.\end{cases}
 \end{align*}
Show that $f$ is not Riemann integrable.
 \end{example}
\begin{solution}{Solution}
Let $P=\{x_i\}_{i=0}^k$ be any partition of the interval $[0,1]$. For any $1\leq i\leq k$, by denseness of the set of rational numbers and the set of irrational numbers, there exist a rational number and an irrational number in the interval $[x_{i-1}, x_i]$.  This shows that \[m_i=0,\hspace{1cm}M_i=1,\hspace{1cm}\text{for all}\;1\leq i\leq k.\]
Hence, 
\begin{align*}L(f,P)&=\sum_{i=1}^k m_i(x_i-x_{i-1})=0,\\ U(f,P)&=\sum_{i=1}^k M_i(x_i-x_{i-1})=1.\end{align*}
This shows that  
\[S_L(f)=\{0\},\hspace{1cm}S_U(f)=\{1\}.\]
Therefore, the lower integral and the upper integral of $f$ are
\[\underline{\int_a^b}f=0,\hspace{1cm}\overline{\int_a^b}f=1\]respectively. Since they are not equal, $f$ is not Riemann integrable.
\end{solution}

An interesting question now is what functions are Riemann integrable. 
Let us first give alternative criteria for Riemann integrability.
\begin{lemma}[label=230220_10]
{}
Let $f:[a,b]\to\mathbb{R}$ be a bounded function. Then the following are equivalent.
\begin{enumerate}[(a)]
\item $f:[a,b]\to\mathbb{R}$  is Riemann integrable.
\item For any $\varepsilon>0$, there exists a partition $P$ of $[a,b]$ such that
\[U(f,P)-L(f,P)<\varepsilon.\]

\end{enumerate}
\end{lemma}
\begin{myproof}{Proof}
First, let us prove (a)  implies  (b). If $f$ is Riemann integrable, \[\int_a^bf=\underline{\int_a^b}f =\overline{\int_a^b} f.\] Given $\varepsilon>0$, by definitions of lower and upper integrals as supremums and infimums, there exist partitions $P_1$ and $P_2$ of $[a,b]$ such that
\[L(f,P_1)>\int_a^b f-\frac{\varepsilon}{2}\hspace{1cm}\text{and}\hspace{1cm}
U(f,P_2)<\int_a^b f+\frac{\varepsilon}{2}.\]
This gives 
\[U(f,P_2)-L(f,P_1)<\varepsilon.\]
Let $P$ be a common refinement of $P_1$ and $P_2$. Then
\[L(f,P_1)\leq L(f,P)\leq U(f,P)\leq U(f,P_2).\]
This implies that
\[
U(f,P)-L(f,P) \leq U(f,P_2)-L(f,P_1)< \varepsilon.
\]
Conversely, assume that (b) holds. Then for every positive integer $n$, there is a partition $P_n$ of $[a,b]$ such that\bp
\[0\leq \;U(f,P_n)\;-\;L(f,P_n)<\frac{1}{n}.\]
Therefore,
\[0\leq\overline{\int_a^b }f \;-\;\underline{\int_a^b} f <\frac{1}{n}.\]
Taking the  $n\to \infty$ limit, squeeze theorem implies that
\[\underline{\int_a^b} f\;=\;\overline{\int_a^b }f.\]
This shows that $f$ is Riemann integrable.
\end{myproof}
 
\begin{theorem}[label=230220_11]
{The Archimedes-Riemann Theorem}
Let $f:[a,b]\to\mathbb{R}$ be a bounded function. Then $f:[a,b]\to\mathbb{R}$  is Riemann integrable if and only if there is a sequence $\{P_n\}$ of partitions of $[a,b]$ such that
\begin{equation}\label{eq230220_12}\lim_{n\to \infty}(U(f,P_n)-L(f,P_n))=0.\end{equation}
In this case, the Riemann integral of $f$ over $[a,b]$ can be computed by
\begin{equation}\label{eq230220_13}\int_a^b f=\lim_{n\to\infty}L(f,P_n)=\lim_{n\to\infty}U(f,P_n).\end{equation}
\end{theorem}

This theorems says that the Riemann integrability of a function can be checked by the existence of a sequence of partitions satisfying \eqref{eq230220_12}. Such sequence of partitions can also be used to compute the   Riemann integral. Thus, we give such sequence a special name.
\begin{definition}{Archimedes Sequence of Partitions}
Let $f:[a,b]\to\mathbb{R}$ be a bounded function. A sequence $\{P_n\}$ of partitions of $[a,b]$ is called an Archimedes sequence of partitions for $f$ provided that
\[\lim_{n\to \infty}(U(f,P_n)-L(f,P_n))=0.\]
\end{definition}Hence, Theorem \ref{230220_11} says that $f:[a,b]\to\mathbb{R}$ is Riemann integrable if and only if it has an Archimedes sequence of partitions. 

\begin{myproof}{\linkt Proof of Theorem \ref{230220_11}}
If $f:[a,b]\to\mathbb{R}$ is Riemann integrable,  by Lemma \ref{230220_10}, for every positive integer $n$, there is a partition $P_n$ of $[a,b]$ such that
\[0\leq U(f,P_n)-L(f,P_n)<\frac{1}{n}.\]
 By squeeze theorem,
\[\lim_{n\to\infty}(U(f,P_n)-L(f,P_n))=0.\] Conversely, if there is a sequence $\{P_n\}$ of partitions of $[a,b]$ such that
\[\lim_{n\to \infty}(U(f,P_n)-L(f,P_n))=0,\]  the definition of limit of sequences implies that for every $\varepsilon>0$, there is a positive integer $N$ such that for all $n\geq N$,
\[|U(f,P_n)-L(f,P_n)|<\varepsilon.\]In particular, we find that $P_N$ is a partition of $[a,b]$ satisfying
\[U(f,P_N)-L(f,P_N)<\varepsilon.\]
By Lemma \ref{230220_10} again, we find that $f$ is Riemann integrable.
This means $\di\int_a^bf=\underline{\int_a^b}f =\overline{\int_a^b} f$. Therefore,
\[0\;\leq\; U(f,P_n)\;-\;\int_a^b f\;\leq U(f,P_n)-L(f,P_n).\]
By taking the $n\to \infty$ limit, we find that
\[\lim_{n\to \infty}U(f,P_n)=\int_a^bf.\]Since $\di \lim_{n\to \infty}(U(f,P_n)-L(f,P_n))=0$, we find that
  \[\lim_{n\to \infty}L(f,P_n)=\lim_{n\to \infty}U(f,P_n)=\int_a^bf.\]  
\end{myproof}

Let us look at an example how to apply the Archimedes-Riemann theorem to compute integrals.
\begin{example}[label=ex230221_4]{}
Let $f:[1, 4]\to\mathbb{R}$ be the function $f(x)=x^2 $. Show that $f$ is Riemann integrable and find the integral $\di \int_1^4f(x)dx$.
\end{example}
Here we need the formulas
\[\sum_{i=1}^ni=\frac{n(n+1)}{2},\hspace{1cm}\sum_{i=1}^ni^2=\frac{n(n+1)(2n+1)}{6}.\]
\begin{solution}{Solution}
 Let $n$ be a positive integer, and let $P_n=\{x_0, x_1, \ldots, x_n\}$ be the regular partition of $[1,4]$ into $n$  intervals. Then
\[x_i=1+\frac{3i}{n}.\]
Notice that  the function $f:[1, 4]\to\mathbb{R}$, $f(x)=x^2 $ is an increasing function. Therefore, on the interval $[x_{i-1}, x_i]$, 
\[m_i=f(x_{i-1})=\left(1+\frac{3i-3}{n}\right)^2,\hspace{1cm}M_i=f(x_i)=\left(1+\frac{3i}{n}\right)^2.\]
From this, we find that
\begin{align*}
L(f,P_n)&=\sum_{i=1}^nm_i(x_i-x_{i-1})=\frac{3}{n}\sum_{i=1}^n\left(1+\frac{6(i-1)}{n}+\frac{9(i-1)^2}{n^2}\right)\\
&=\frac{3}{n}\left(n+3(n-1)+\frac{3(n-1)(2n-1)}{2n}\right)\\
&=\frac{3}{2n^2}(14n^2-15n+3),
\end{align*}\bs
\begin{align*}
U(f,P_n)&=\sum_{i=1}^nM_i(x_i-x_{i-1})=\frac{3}{n}\sum_{i=1}^n\left(1+\frac{6i}{n}+\frac{9i^2}{n^2}\right)\\
&=\frac{3}{n}\left(n+3(n+1)+\frac{3(n+1)(2n+1)}{2n}\right)\\
&=\frac{3}{2n^2}(14n^2+15n+3).
\end{align*}Moreover,
\[U(f,P_n)-L(f,P_n)=\frac{45}{n}.\]
It follows that
\[\lim_{n\to\infty}\left(U(f,P_n)-L(f,P_n)\right)=0.\]This proves that $\{P_n\}$ is an Archimedes sequence of partitions for $f$. Therefore, $f$ is Riemann integrable, and
\[\int_1^4f(x)dx=\lim_{n\to\infty}U(f,P_n)=\lim_{n\to\infty} \frac{3}{2 }\left(14 +\frac{15}{n}+\frac{3}{n^2}\right)=21.\]
\end{solution}

The following gives an $\varepsilon-\delta$ characterization of Riemann integrability   in terms of Darboux sums.
\begin{theorem}[label=230221_2]{Equivalent Definitions of Riemann Integrability}
Let  $f:[a,b]\to\mathbb{R}$ be a bounded function. Then the following two statements are equivalent.
\begin{enumerate}[(i)]
\item $f:[a,b]\to\mathbb{R}$ is Riemann integrable, in the sense that $\di \underline{\int_a^b }f=\overline{\int_a^b} f$.
\item 
For any $\varepsilon>0$, there exists a $\delta>0$ so that   if $P=\{x_i\}_{i=0}^{k}$ is a partition of $[a,b]$ with $|P|<\delta$,  then
\[U(f,P)-L(f,P)<\varepsilon.\]
\end{enumerate}

\end{theorem}
 By Lemma  \ref{230220_10}, $f:[a,b]\to\mathbb{R}$ is Riemann integrable if and only if for every $\varepsilon>0$, there is a partition $P$ satisfying $U(f,P)-L(f,P)<\varepsilon$. The highly nontriviality of this theorem is   the existence of a single partition satisfying $U(f,P)-L(f,P)<\varepsilon$ is equivalent to the existence of a positive number $\delta$ such that \emph{all} partitions $P$ with gaps less than $\delta$   satisfy $U(f,P)-L(f,P)<\varepsilon$.

\begin{myproof}{Proof}
 (ii) implies (i) follows trivially from  Lemma  \ref{230220_10}.

Now assume that (i) holds.
Since $f:[a,b]\to\mathbb{R}$ is bounded, there exists a positive number $M$ such that
\[|f(x)|\leq M\hspace{1cm}\text{for all}\;x\in [a,b].\]
Given $\varepsilon>0$,  Lemma \ref{230220_10} implies that  there is a partition \[P_0=\{\widetilde{x}_0, \widetilde{x}_1, \ldots, \widetilde{x}_s\}\] of $[a,b]$ such that
\[U(f,P_0)-L(f,P_0)<\frac{\varepsilon}{2}.\]
Take
\[\delta=\frac{\varepsilon}{8sM}.\] Then $\delta>0$. If $P=\{x_0, x_1, \ldots, x_k\}$ is a partition of $[a,b]$ with $|P|<\delta$, we want to show that
\[U(f, P)-L(f,P)=\sum_{i=1}^k(M_i-m_i)(x_i-x_{i-1})<\varepsilon.\]
Here\[m_i=\inf_{x_{i-1}\leq x\leq x_i}f(x),\hspace{1cm}M_i=\sup_{x_{i-1}\leq x\leq x_i}f(x).\]   
 Let
\[E_1=\left\{1\leq i\leq k\,|\, \exists j,\;\widetilde{x}_j\in [x_{i-1}, x_i]\right\}\]
be the set that contains those indices $i$ where the interval $[x_{i-1}, x_i]$ contains a partition point of $P_0$, and
let $E_2=\{1, 2, \ldots, k\}\setminus E_1$ be the set of those indices that are not in $E_1$. 
\bp
  By definition, the point $\widetilde{x}_0$ can only be in $[x_0, x_1]$, and the point $\widetilde{x}_s$ can only be in $[x_{k-1}, x_k]$. For any $1\leq j\leq s-1$, $\widetilde{x}_j$ can be in at most two different subintervals of $P$. Hence, $E_1 $ contains at most $2s$ elements.

Splitting the sum over $i$ to a sum over $E_1$ and a sum over $E_2$, we have
\[U(f,P)-L(f,P)=\sum_{i\in E_1}(M_i-m_i)(x_i-x_{i-1})+\sum_{i\in E_2}(M_i-m_i)(x_i-x_{i-1}).\] 
First we estimate the sum over $E_1$. 
 Since \[-M\leq f(x)\leq M\quad\text{ for all }\;x\in [a,b],\] we find that for any $1\leq i\leq k$,
\[0\leq M_i-m_i\leq 2M.\]
Since $x_i-x_{i-1}\leq|P|<\delta$, we have 
\[\sum_{i\in E_1}(M_i-m_i)(x_i-x_{i-1})\leq \sum_{i\in E_1}2M(x_i-x_{i-1})\leq 2M\delta |E_1|\leq 4sM\delta \leq\frac{\varepsilon}{2}.\]Let $P^*$ be the common refinement of $P$ and $P_0$ obtained by taking the union of their partition points. By our definitions of $E_1$ and $E_2$, for each $i$ in $E_2$, $[x_{i-1}, x_i]$ is also a partition interval in $P^*$. 
 Therefore,
\begin{align*} \sum_{i\in E_2}\left(M_i-m_i\right)(x_i-x_{i-1}) &\leq U(f,P^*)-L(f,P^*)\\&\leq U(f,P_0)-L(f,P_0)<\frac{\varepsilon}{2}.\end{align*}
 The two estimates above imply that
\[U(f,P)-L(f,P)<  \varepsilon,\]
which completes the proof that (i)  implies (ii).
\end{myproof}

A disadvantage of working with Darboux sums is we need to figure out the infimum and supremum of a function over the partition intervals. Let us turn to Riemann sums.

\begin{lemma}[label=230221_3]{}
Let $f:[a,b]\to\mathbb{R}$ be a bounded function, and let $P=\{x_i\}_{i=0}^k$ be a partition of $[a,b]$.
 For every  $\varepsilon>0$, there exist   choices of intermediate points $A$ and  $B$ for the partition $P$ such that
\begin{align*}
0\leq R(f,P,A)-L(f,P)<\varepsilon,\quad
 0\leq U(f, P)- R(f,P,B) <\varepsilon.\end{align*}
 
\end{lemma}
\begin{myproof}{Proof}For $1\leq i\leq k$, let
$\di m_i=\inf_{x_{i-1}\leq x\leq x_i}f(x)$ and $\di M_i=\sup_{x_{i-1}\leq x\leq x_i}f(x)$.
By definitions of infimum and supremum, for each $1\leq i\leq k$, there are points $\xi_i$ and $\eta_i$ in $[x_{i-1}, x_i]$ such that
\begin{gather*}
m_i\leq f(\xi_i)<m_i+\frac{\varepsilon}{ (b-a)},\\  M_i-\frac{\varepsilon}{ (b-a)}<f(\eta_i)\leq M_i.\end{gather*}
Multiply  by $(x_i-x_{i-1})$ and sum  over $i$, we find that
\begin{equation*} \begin{split}
\sum_{i=1}^km_i(x_i-x_{i-1})\leq \sum_{i=1}^kf(\xi_i)(x_i-x_{i-1}) <\sum_{i=1}^km_i(x_i-x_{i-1})+ \varepsilon,\\
\sum_{i=1}^kM_i(x_i-x_{i-1})- \varepsilon<\sum_{i=1}^kf(\eta_i)(x_i-x_{i-1})\leq  \sum_{i=1}^kM_i(x_i-x_{i-1}).
\end{split}\end{equation*}
Let $A=\{\xi_i\}_{i=1}^k$ and  $B=\{\eta_i\}_{i=1}^k$. They are choices of intermediate points for the partition $P$.  The two inequalities above give
\begin{gather*}
L(f,P)\leq R(f,P,A)<L(f,P)+\varepsilon,\\ U(f,P)-\varepsilon<R(f,P,B)\leq U(f,P),\end{gather*}which are the desired results.
\end{myproof}

The following gives an $\varepsilon-\delta$ definition for Riemann integrability of a bounded function.
\begin{theorem}[label=230220_15]{Equivalent Definitions of Riemann Integrability}
Let $f:[a,b]\to\mathbb{R}$ be a bounded function. Consider  the following two definitions for $f$ to be Riemann integrable.
\begin{enumerate}[(i)]
\item $\di \underline{\int_a^b }f=\overline{\int_a^b} f$.

\item There is a number $I$ such that for any $\varepsilon>0$, there exists a $\delta>0$ so that   if $P=\{x_i\}_{i=0}^{k}$ is a partition of $[a,b]$ with $|P|<\delta$, $A=\{\xi_i\}_{i=1}^k$ is a choice of intermediate points for $P$, then
\[|R(f,P,A)-I|<\varepsilon.\]

\end{enumerate}
These two statements are equivalent, and in case $f$ is Riemann integrable, \[I=\underline{\int_a^b }f=\overline{\int_a^b} f=\int_a^b f.\]
\end{theorem}
Note that statement  (ii) can be  expressed as saying the limit of Riemann sums
\[I=\lim_{|P| \to 0}R(f, P, A)\]exists.

\begin{myproof}{Proof}
 
Assume that (i) holds. 
Let \[I=\underline{\int_a^b }f=\overline{\int_a^b} f.\]
Given $\varepsilon>0$, by Theorem \ref{230221_2}, there exists a $\delta>0$ such that if $P=\{x_i\}_{i=0}^k$ is a partition of $[a,b]$ with $|P|<\delta$, then 
\[U(f,P)-L(f,P)<\varepsilon.\]\bp
If $A=\{\xi_i\}_{i=1}^k$ is any choice of intermediate points for the partition $P$,  
\[L(f,P)\leq R(f, P, A)\leq U(f,P).\] Since we also have
\[L(f,P)\leq I\leq U(f,P),\]we find that
\[|R(f, P, A)-I|\leq  U(f,P)-L(f,P)<\varepsilon.\]
This proves that (i) implies (ii).

Conversely,   assume that (ii) holds. By Lemma \ref{230220_10}, to show that (i) holds, it suffices to prove that
  for any $\varepsilon>0$, there is a partition $P$ of $[a,b]$ so that
\[U(f,P)-L(f,P)<\varepsilon.\]
Given $\varepsilon>0$, (ii) implies   there is  a $\delta>0$ such that if $P=\{x_i\}_{i=0}^{k}$ is a partition of $[a,b]$ with $|P|<\delta$, $A=\{\xi_i\}_{i=1}^k$ is a choice of intermediate points for the partition $P$, then
\begin{equation}\label{eq230220_16}|R(f,P,A)-I|<\frac{\varepsilon}{4}.\end{equation}Here $I$ is the limit of Riemann sums implied by (ii).  Let $P=\{x_i\}_{i=0}^n$ be  a regular partition into $n$ intervals, where $n$ is large enough so that \[|P|=\frac{b-a}{n} <\delta.\] By Lemma \ref{230221_3}, there exist choices of intermediate points $A$ and $B$ for the partition $P$ which satisfy
\[U(f,P)<R(f,P,A)+\frac{\varepsilon}{4},\hspace{1cm} L(f,P)>R(f,P,B)-\frac{\varepsilon}{4}.\] 
These imply that
\[U(f,P)-L(f,P)<R(f,P,A)-R(f,P,B)+\frac{\varepsilon}{2}.\]
 \bp By \eqref{eq230220_16},
\[|R(f,P,A)-R(f,P,B)|\leq |R(f,P,A)-I|+|R(f,P,B)-I|<\frac{\varepsilon}{2}.\]
This proves that
\[U(f,P)-L(f,P)<\varepsilon,\]which completes the proof that  (ii)  implies  (i).

\end{myproof}
As a consequence of Theorem \ref{230221_2} and Theorem \ref{230220_15}, we have the following.
\begin{corollary}[label=230618_1]{}
Let $f:[a,b]\to\mathbb{R}$ be a bounded function that is Riemann integrable, and let $\{P_n\}$ be a sequence of partitions of $[a,b]$ such that
\[\lim_{n\to \infty}|P_n|=0.\]Then 
 \begin{enumerate}[(a)]
\item
$\di \int_a^b f=\lim_{n\to \infty}U(f,P_n)= \lim_{n\to \infty}L(f,P_n)$
\item $\di \int_a^b f=\lim_{n\to \infty}R(f,P_n, A_n)$, where for each $n\in\mathbb{Z}^+$, $A_n$ is a choice of intermediate points for the partition $P_n$.
\end{enumerate}
\end{corollary}
This corollary says that if we \emph{ know apriori } that $f:[a,b]\to\mathbb{R}$ is Riemann integrable, then we can evaluate the integral by a sequence of partitions whose gaps goes to 0, using either the Darboux upper sums, or the Darboux lower sums, or Riemann sums for any choice of intermediate points.

\begin{myproof}{Proof}
Let $I=\di\int_a^bf$. Given $\varepsilon>0$, Theorem \ref{230221_2} and Theorem \ref{230220_15} imply that there is a $\delta>$ such that for any partition $P$ with $|P|<\delta$, and any choice of intermdiates points $A$ for the partition $P$,  
\[U(f,P)-L(f,P)<\varepsilon \hspace{1cm}\text{and}\hspace{1cm}|R(f, P, A)-I|<\varepsilon.\]\bp
Since $\di\lim_{n\to \infty}|P_n|=0$, there is a positive integer $N$ so that for all $n\geq N$, $|P_n|<\delta$. This implies that
for all $n\geq N$, 
\[U(f,P_n)-L(f,P_n)<\varepsilon \quad\text{and}\quad|R(f, P_n, A_n)-I|<\varepsilon.\]
Since $L(f,P_n)\leq I\leq U(f,P_n)$,
we find that
\[|U(f,P_n)-I|<\varepsilon\quad\text{and}\quad|L(f,P_n)-I|<\varepsilon \quad\text{for all}\;n\geq N.\]
These prove that
\[I=\lim_{n\to \infty}U(f,P_n) =\lim_{n\to \infty}L(f,P_n)=\lim_{n\to \infty}R(f,P_n, A_n).\]
\end{myproof}
\begin{highlight}{}For every positive integer $n$, take $P_n$ to be the regular partition of $[a,b]$ into $n$ intervals.  This gives a sequence of partitions $\{P_n\}$ with \[\lim_{n\to \infty}|P_n|=\lim_{n\to\infty}\frac{b-a}{n}=0.\] For the choices of intermediate points $A_n$, one can take the left end point of each interval, or the right end point, or the midpoint. \end{highlight}
\begin{example}{}
We are going to prove in Section \ref{sec4.3} that a continuous function is integrable. The function $f:[0, 6]\to\mathbb{R}$, $f(x)=6x-x^2 $ is continuous. Use Riemann sums to evaluate the integral $\di \int_0^6f(x)dx$.
\end{example}
\begin{solution}{Solution}
For a positive integer $n$, let $P_n=\{x_i\}_{i=0}^n$ be the regular partition of $[0,6]$ into $n$ intervals. Then
\[x_i=\frac{6i}{n},\hspace{1cm} 0\leq i\leq n.\]\bs
Let $A_n=\{\xi_i\}_{i=1}^n$, where
\[\xi_i=x_i=\frac{6i}{n}, \quad 1\leq i\leq n.\]  Then
\begin{align*}
R(f,P_n,A_n)&=\sum_{i=1}^nf(\xi_i)(x_i-x_{i-1})\\&=\frac{6}{n}\sum_{i=1}^n\left(\frac{36i}{n}-\frac{36i^2}{n^2}\right)\\
&=\frac{216}{n^2}\left(\frac{n (n+1)}{2}-\frac{ (n+1)(2n+1)}{6}\right)\\
&=\frac{36(n^2-1)}{n^2}.
\end{align*}
Therefore,
\[\int_0^6f(x)dx=\lim_{n\to\infty}R(f,P_n,A_n)=36.\]
\end{solution}
\vp
\noindent
{\bf \large Exercises  \thesection}
\setcounter{myquestion}{1}

\begin{question}{\themyquestion}
Let $f:[0, 2]\to\mathbb{R}$ be the function $f(x)=4-x^2 $. Given a positive integer $n$, let $P_n$ be the regular partition of $[0,2]$ into $n$ subintervals.
\begin{enumerate}[(a)]
\item Compute $L(f,P_n)$ and $U(f, P_n)$.
\item Show directly that $\di \lim_{n\to\infty}\left(U(f,P_n)-L(f,P_n)\right)=0$.
\item Use part (b) to conclude that $f$ is Riemann integrable and find the integral $\di \int_0^2f(x)dx$.
\end{enumerate}
\end{question}
 \atc

\begin{question}{\themyquestion}Given that the functon $f:[0, 4]\to\mathbb{R}$, $f(x)=x^2-2x+3$ is Riemann integrable. Use Riemann sums to evaluate the integral $\di \int_0^4f(x)dx$.
\end{question}
\vp

\section{Properties of Riemann Integrals  }\label{sec4.2}
In this section, we derive some properties of the Riemann integrals.
First we show that integral of a nonnegative function is nonnegative.
\begin{theorem}[label=230221_5]{}
If $f:[a,b]\to\mathbb{R}$ is a bounded function that is Riemann integrable, and \[f(x)\geq 0\hspace{1cm} \text{for all}\;x\in [a, b],\] then \[\int_a^b f\geq 0.\]
\end{theorem}
\begin{myproof}{Proof}
Since $f$ is Riemann integrable,  
\[\int_a^b f=\lim_{n\to \infty}L(f,P_n),\]
where $P_n$ is the regular partition of $[a,b]$ into $n$ intervals. 
Since $f(x)\geq 0$ for all $x\in [a, b]$, we find that\[L(f,P_n)\geq 0\hspace{1cm}\text{for all }\;n\in\mathbb{Z}^+.\] Therefore,
\[\int_a^b f\geq 0.\]
\end{myproof}

Linearity is always an important property. 
\begin{theorem}[label=230221_6]{Linearity of Integrals}
Let $f:[a,b]\to\mathbb{R}$ and $g:[a,b]\to\mathbb{R}$ be bounded functions. If $f$ and $g$ are Riemann integrable, then for any constants $\alpha$ and $\beta$, $\alpha f+\beta g:[a,b]\to\mathbb{R}$ is also Riemann integrable, and
\[\int_a^b(\alpha f+\beta g)=\alpha\int_a^b f+\beta \int_a^b g.\]
\end{theorem}
\begin{myproof}{Proof}
Here we   use   the fact that a function $h:[a,b]\to\mathbb{R}$ is Riemann integrable if and only if the limit
$\di \lim_{|P|\to 0}R(h,P,A)$ exists.  Since  $f:[a,b]\to\mathbb{R}$ and $g:[a,b]\to\mathbb{R}$ are bounded, $\alpha f+\beta g$ is also bounded. Since  $f:[a,b]\to\mathbb{R}$ and $g:[a,b]\to\mathbb{R}$ are Riemann integrable, 
\[\int_a^b f=\lim_{|P|\to 0}R(f,P,A),\hspace{1cm}\int_a^bg=\lim_{|P|\to 0}R(g,P,A).\]
Notice that for any partition $P=\{x_i\}_{i=0}^k$ of $[a,b]$, and any choice of intermediate points $A=\{\xi_i\}_{i=1}^k$ for the partition $P$,
\begin{align*}R(\alpha f+\beta g, P, A)&=\sum_{i=1}^k\left(\alpha f(\xi_i)+\beta g(\xi_i)\right)(x_i-x_{i-1})\\
&=\alpha R(f,P,A)+\beta R(g,P,A).\end{align*}
Limit laws imply that
\begin{align*}\lim_{|P|\to 0}R(\alpha f+\beta g, P, A)&=\alpha\lim_{|P|\to 0}R(f,P,A)+\beta\lim_{|P|\to 0}R(g,P,A)\\&=\alpha \int_a^b f+\beta\int_a^b g.\end{align*}
This proves that $\alpha f+\beta g $ is  Riemann integrable and
\[\int_a^b(\alpha f+\beta g)=\alpha\int_a^b f+\beta \int_a^b g.\]
\end{myproof}

From the previous two theorems, we   obtain a comparison theorem for integrals.
\begin{theorem}{Monotonicity}
Let $f:[a,b]\to\mathbb{R}$ and $g:[a,b]\to\mathbb{R}$ be bounded functions. If $f$ and $g$ are Riemann integrable, and 
\[f(x)\geq g(x)\hspace{1cm}\text{for all}\;x\in [a,b],\]
then
\[\int_a^b f\geq \int_a^b g.\]
\end{theorem}
\begin{myproof}{Proof}
Define $h:[a,b]\to\mathbb{R}$ to be the function
$h(x)=f(x)-g(x)$.
Then 
$h(x)\geq 0$ for all $x\in [a,b]$. 
By Theorem \ref{230221_6}, $h$ is Riemann integrable and 
\[\int_a^b h=\int_a^b f-\int_a^b g.\]
By Theorem \ref{230221_5}, $\di\int_a^b h\geq 0$. Hence,
\[\int_a^b f\geq \int_a^b g.\]
\end{myproof}
We can apply the monotonicity theorem to obtain  bounds for  an integral  from the lower bound and the upper bound of the function.
\begin{example}{}
Let $f:[a,b]\to\mathbb{R}$ be a Riemann integrable funtion satisfying 
\[m\leq f(x)\leq M\hspace{1cm}\text{for all}\;x\in [a,b].\]
Then
\[m(b-a)\;\leq\;\int_a^b f\;\leq\; M(b-a).\]
\end{example}

When an interval is partitioned into a finite collection of intervals,
the integral over the whole interval is expected to equal to the sum of the intergrals over the subintervals. It is enough for us to consider two subintervals.
\begin{theorem}{Additivity}
Let $f:[a,b]\to\mathbb{R}$ be a bounded function, and let $c$ be a point in $(a,b)$.
\begin{enumerate}[(a)]
\item
If $f:[a,b]\to\mathbb{R}$ is Riemann integrable, then  $f:[a,c]\to\mathbb{R}$  and $f:[c,b]\to\mathbb{R}$  are Riemann integrable.
\item If $f:[a,c]\to\mathbb{R}$  and $f:[c,b]\to\mathbb{R}$  are Riemann integrable, then $f:[a,b]\to\mathbb{R}$ is Riemann integrable. 
\end{enumerate}In either case,
\[\int_a^b f=\int_a^c f+\int_c^b f.\]
\end{theorem}

\begin{myproof}{Proof}
We use  Lemma \ref{230220_10}.  First we prove (a). Given $\varepsilon>0$, since $f:[a,b]\to\mathbb{R}$ is Riemann integrable, there is a partition $P$ of $[a,b]$ such that
\[U(f,P)-L(f,P)<\varepsilon.\]
Let $P^*$ be the partition of $[a,b]$ that is obtained  by taking the union of the partition points in $P$ and $P_0=\{a, c, b\}$. If $P$ already contains $c$ as a partition point, then $P^*=P$. In any case, $P^*$ is a refinement of $P$.   Therefore,
\[U(f,P^*)-L(f,P^*)\leq U(f,P)-L(f,P)<\varepsilon.\]
Consider $P^*$ as a refinement of $P_0$, let $P_1$ be the partition of $[a,c]$   induced by $P^*$, and let $P_2$ be the partition of $[c, b]$ induced by $P^*$. Then
\[L(f,P^*)=L(f,P_1)+L(f,P_2), \hspace{1cm}U(f,P^*)=U(f,P_1)+U(f,P_2).\]
These imply that
\[(U(f,P_1)-L(f,P_1))+(U(f,P_2)-L(f,P_2))=U(f,P^*)-L(f,P^*)<\varepsilon.\]
\bp Since $U(f,P_1)-L(f,P_1)\geq 0$ and $U(f,P_2)-L(f,P_2)\geq 0$, we find that
\[U(f,P_1)-L(f,P_1)<\varepsilon\quad\text{and}\quad U(f,P_2)-L(f,P_2)<\varepsilon.\]
By Lemma \ref{230220_10},  we conclude that $f:[a,c]\to\mathbb{R}$  and $f:[c,b]\to\mathbb{R}$  are Riemann integrable.
 
Next, we prove (b). Given $\varepsilon>0$,  since $f:[a,c]\to\mathbb{R}$  and $f:[c,b]\to\mathbb{R}$  are Riemann integrable, there exists a partition $P_1$ of $[a,c]$, and a partition $P_2$ of $[c,b]$ such that 
\[U(f,P_1)-L(f,P_1)<\frac{\varepsilon}{2}\quad\text{and}\quad U(f,P_2)-L(f,P_2)<\frac{\varepsilon}{2}.\]
Let $P$ be the partition of $[a,b]$ obtained by taking the union of the partition points in $P_1$ and $P_2$. Then 
\[L(f,P )=L(f,P_1)+L(f,P_2), \hspace{1cm}U(f,P )=U(f,P_1)+U(f,P_2).\]Therefore,
\[U(f,P)-L(f,P)=(U(f,P_1)-L(f,P_1))+(U(f,P_2)-L(f,P_2))<\varepsilon.\]
This proves that $f:[a,b]\to\mathbb{R}$ is Riemann integrable.

Now we prove the last statement.
For any positive integer $n$,   let $P_{1,n}$ be the regular partition of $[a,c]$ into $n$ intervals, and let $P_{2,n}$ be the regular partition of $[c,b]$ into $n$ intervals. Then let $P_n$ be the partition of $[a,b]$ obtained by taking the union of the partition points in $P_{1,n}$ and $P_{2,n}$.  For the Darboux upper sums, we have 
\[U(f,P_n)=U(f,P_{1,n})+U(f,P_{2,n}).\] 
Taking the $n\to \infty$ limits on both sides, we conclude that
\[\int_a^b f=\int_a^c f+\int_c^b f.\]
\end{myproof}
\begin{highlight}{Extension of Definition of Integrals}
The additivity allows us to extend the definition of the integral $\di \int_a^b f$ to the case where $a\geq b$. We define
\[\int_a^a f=0.\] If $a>b$, define
\[\int_a^b f=-\int_b^a f.\]
Then one can check that 
as long as two of the three  integrals $\di\int_a^b f$, $\di \int_a^cf$, $\di \int_c^bf$ exist, the third one also exists, and we always have
\[\int_a^b f=\int_a^c f+\int_c^b f.\]

\end{highlight}

Using induction, we can extend the additivity theorem.
\begin{corollary}[label=230221_9]{General Additivity Theorem}
Let $f:[a,b]\to\mathbb{R}$ be a bounded function, and let $P_0=\{a_0, a_1, \ldots, a_k\}$ be a partition of $[a,b]$. Then $f:[a,b]\to\mathbb{R}$ is Riemann integrable if and only if for each $1\leq i\leq k$, $f:[a_{i-1},a_i]\to\mathbb{R}$ is Riemann integrable. In this case,
\[\int_a^b f(x)dx=\int_{a}^{a_1}f+\int_{a_1}^{a_2}f+\cdots+\int_{a_{k-2}}^{a_{k-1}}f+\int_{a_{k-1}}^bf.\]
\end{corollary}
\vp
\noindent
{\bf \large Exercises  \thesection}
\setcounter{myquestion}{1}

\begin{question}{\themyquestion}
Given that $f:[2, 7]\to\mathbb{R}$ is a function satisfying
\[-3\leq f(x)\leq 11\hspace{1cm} \text{for all}\;x\in [2, 7].\]
Find a lower bound and an upper bound for $\di \int_2^7f(x)dx$. 
\end{question}

\atc
\begin{question}{\themyquestion}
Given that $f:[a,b]\to\mathbb{R}$ and $g:[a,b]\to\mathbb{R}$ are bounded functions,  $P$ is a partition of $[a,b]$, and $c$ and $d$ are two points in $[a, b]$ with $c<d$. Prove the following.
\begin{enumerate}[(a)]
\item $\di\inf_{c\leq x\leq d}(f+g)(x)\geq \inf_{c\leq x\leq d}f(x)+\inf_{c\leq x\leq d}g(x)$.
\item  $\di\sup_{c\leq x\leq d}(f+g)(x)\leq \sup_{c\leq x\leq d}f(x)+\sup_{c\leq x\leq d}g(x)$.
\item $L(f+g, P)\geq L(f, P)+L(g, P)$.
\item $U(f+g,P)\leq U(f,P)+U(g,P)$.
\end{enumerate}
  Then use (c) and (d) to give a proof of the following statement: If  $f:[a,b]\to\mathbb{R}$ and $g:[a,b]\to\mathbb{R}$ are Riemann integrable, then $f+g:[a,b]\to\mathbb{R}$ is also Riemann integrable, and
\[\int_a^b(f+g)=\int_a^b f+\int_a^b g.\]
\end{question}

\vp

\section{Functions that are Riemann Integrable}\label{sec4.3}
In this section, we are going to derive Riemann integrability of a few classes of functions. The first class of functions that are of interest is the class of continuous  functions. If $f:[a,b]\to\mathbb{R}$ is a continuous function, then $f([a,b])$ is sequentially compact. In particular, $f([a,b])$ is bounded. Hence, if $f:[a,b]\to\mathbb{R}$ is a continuous function, it is bounded. 
The crucial property for a continuous function defined on a closed and bounded interval to be integrable is   uniform continuity.  

\begin{theorem}[label=230221_7]{}
Let $f:[a,b]\to\mathbb{R}$ be a continuous function. Then $f:[a,b]\to\mathbb{R}$ is Riemann integrable.
\end{theorem}
\begin{myproof}{Proof}
Since  $f:[a,b]\to\mathbb{R}$ is a continuous function defined on a closed and bounded interval, it is uniformly continuous. Given $\varepsilon>0$, there exists $\delta>0$ such that for any $u$ and $v$ in $[a,b]$ with $|u-v|<\delta$,
\[|f(u)-f(v)|< \frac{\varepsilon}{b-a}.\]
Let $P=\{x_i\}_{i=0}^k$ be a partition of $[a,b]$ with $|P|<\delta$. For any $1\leq i\leq k$, $f:[x_{i-1},x_i]\to\mathbb{R}$ is continuous. By extreme value theorem, there exists $u_i$ and $v_i$ in $[x_{i-1}, x_i]$ such that
\[m_i=\inf_{x_{i-1}\leq x\leq x_i}f(x)=f(u_i),\hspace{1cm}M_i=\sup_{x_{i-1}\leq x\leq x_i}f(x)=f(v_i).\]Then
\[|u_i-v_i|\leq x_i-x_{i-1}\leq |P|<\delta.\]
Therefore,
\[M_i-m_i=f(v_i)-f(u_i)< \frac{\varepsilon}{b-a}.\]\bp
Hence,
\begin{align*}
U(f,P)-L(f,P)&=\sum_{i=1}^k (M_i-m_i)(x_i-x_{i-1}) \\&<\frac{\varepsilon}{b-a}\sum_{i=1}^k(x_i-x_{i-1}) =\varepsilon.\end{align*}
This proves that $f:[a,b]\to\mathbb{R}$ is Riemann integrable.
\end{myproof}

\begin{highlight}{}
It follows from this theorem that all the following classes of functions are integrable on a  closed and bounded interval that is contained in their domains.
\begin{enumerate}[$\bullet$\;\;]
\item Polynomials
\item Rational Functions
\item Exponential Functions
\item Logarithmic Functions
\item Trigonometric Functions

\end{enumerate}
\end{highlight} 

\begin{figure}[ht]
\centering
\includegraphics[scale=0.2]{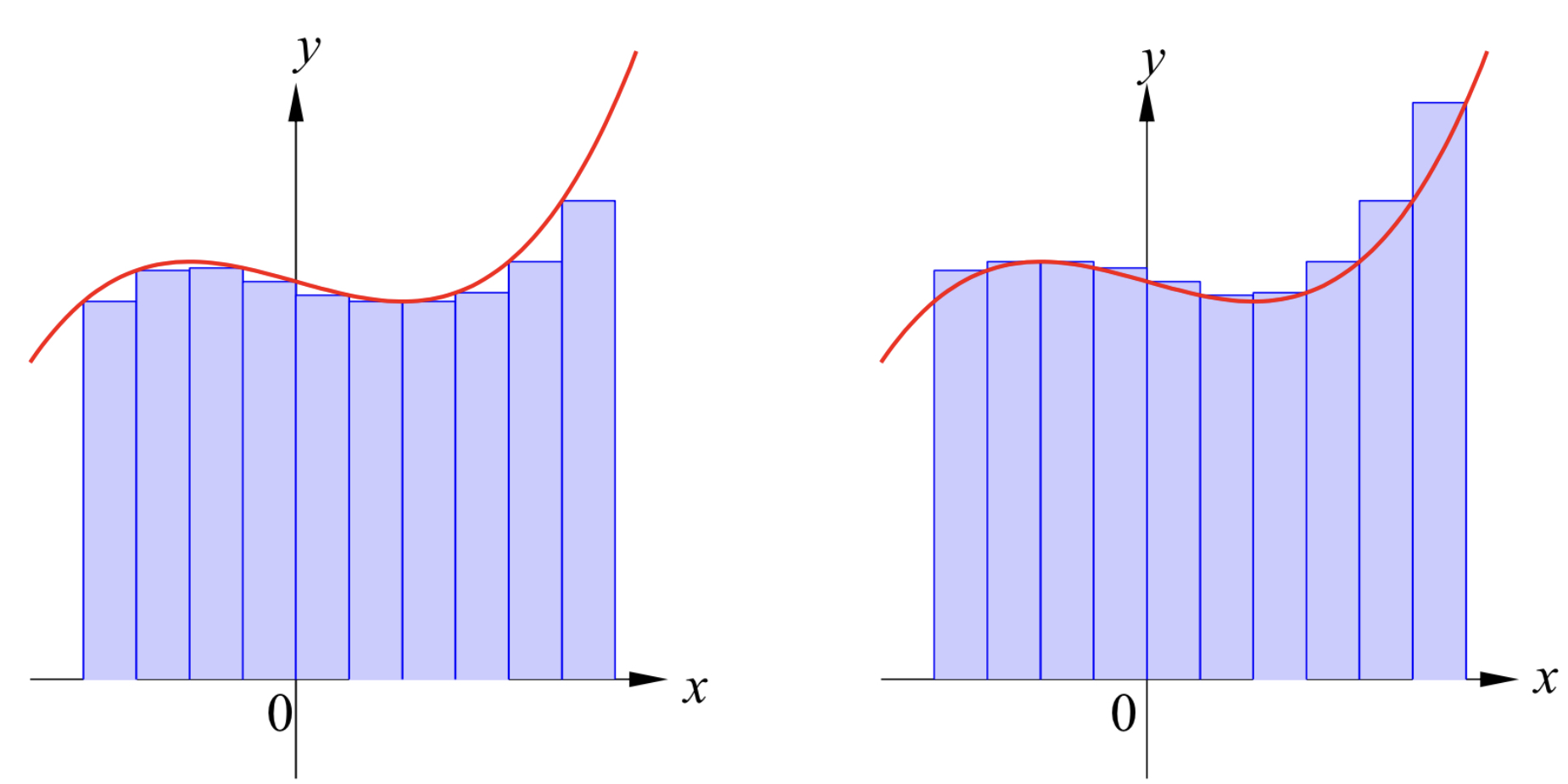}
\caption{Darboux lower sum underestimates area while Darboux upper sum overestimates area.  }\label{figure70}
\end{figure}
Let us revisit the concept of area, which is our original motivation to define integrals.
\begin{remark}{Area}
Let $f:[a,b]\to \mathbb{R}$ be a continuous function such that $f(x)\geq 0$ for all $x\in [a,b]$, and let $R$ be the region bounded between the $x$-axis, the lines $x=a$ and $x=b$, as well as the curve $y=f(x)$. 
\end{remark}\begin{highlight}{}Given $P$ a partition of $[a,b]$, the Darboux lower sum $L(f,P)$ is a sum of areas of rectangles that are inside $R$. The Darboux upper sum $U(f,P)$ is a sum of areas of rectangles whose union contains $R$. Therefore, if $R$ has a area $A$, $L(f,P)$ is less than or equal to  $A$, while $U(f,P)$ is larger than or equal to $A$. Since $f$ is continuous, the Riemann integral $\di I=\int_a^b f(x)dx$ exists. By definition, $I$ is the unique number such that
\[L(f,P)\leq I\leq U(f,P)\] for all partitions $P$ of $[a,b]$. Therefore, we define the area of $R$ to be this number $I$. Namely,
\[\text{Area of}\, R=\int_a^b f(x)dx.\]\end{highlight}

There are also other classes of functions that are Riemann integrable, which are useful. First, we relax the continuity condition slightly in the previous theorem.

\begin{theorem}[label=230221_8]{}
Let $f:[a,b]\to\mathbb{R}$ be a bounded  function that is continuous on $(a, b)$. Then $f:[a,b]\to\mathbb{R}$ is Riemann integrable.
\end{theorem}
Here we only assume $f$ is continuous on $(a,b)$. The function can take on any values on the boundary points $a$ and $b$.
\begin{myproof}{Proof}
Since  $f:[a,b]\to\mathbb{R}$ is bounded, there is a positive constant $M$ such that
\[|f(x)|\leq M\hspace{1cm}\text{for all}\;x\in [a,b].\]
Given $\varepsilon>0$, let 
\[r=\min\left\{\frac{\varepsilon}{8M}, \frac{b-a}{3}\right\}.\]Then $r>0$ and $a+r<b-r$. The function $f:[a+r, b-r]\to\mathbb{R}$ is continuous.  By Theorem \ref{230221_7}, $f:[a+r, b-r]\to\mathbb{R}$ is Riemann integrable.

 Therefore, there is a partition $P_1$ of $[a+r, b-r]$ such that
\[U(f,P_1)-L(f,P_1)<\frac{\varepsilon}{2}.\]
Let $P$ be the partition of $[a,b]$ obtained by adding the points $a$ and $b$ to $P_1$.  Then
\begin{align*}U(f,P)-L(f,P)&=U(f,P_1)-L(f,P_1)\\&\quad +r\left(\sup_{a\leq x\leq a+r}f(x)-\inf_{a\leq x\leq a+r}f(x)\right)\\&\quad +r\left(\sup_{b-r\leq x\leq b}f(x)-\inf_{b-r\leq x\leq b}f(x)\right)\\
&<\frac{\varepsilon}{2}+4Mr\\
&\leq \frac{\varepsilon}{2}+\frac{\varepsilon}{2}=\varepsilon.\end{align*}
This proves that $f:[a,b]\to\mathbb{R}$ is Riemann integrable.

\end{myproof}

As we can see in the proof above, the integral $\di\int_a^b f$ does not depend on the function value at the end points. In fact, this is true for  any finite number of points.
\begin{theorem}[label=230222_2]{}
Let $f:[a,b]\to\mathbb{R}$  be a bounded function that is Riemann integrable. Assume that $g:[a,b]\to\mathbb{R}$ is a function  and $S=\{ a_1, a_2, \ldots, a_k\}$ is a finite subset of $[a,b]$ such that
\[g(x)=f(x)\hspace{1cm}\text{for all}\; x\in [a,b]\setminus S.\]
Then  $g:[a,b]\to\mathbb{R}$ is Riemann integrable, and
\[\int_a^b g=\int_a^b f.\]
\end{theorem}

\begin{myproof}{Proof}

Let $h:[a,b]\to\mathbb{R}$ be the function $h(x)=g(x)-f(x)$. Then   $h(x)=0$ for $x\in  [a,b]\setminus S$. Since $S$ is a finite set, $h$ is bounded, and so there is a positive constant $M$ such that $|h(x)|\leq M$ for all $x\in [a,b]$. Given a positive integer $n$, let $P_n=\{x_0, x_1, \ldots, x_n\}$  be the regular partition of $[a,b]$ into $n$ intervals. There are at most $2k$ of the intervals $[x_{i-1}, x_i]$ that contains a point of $S$. In these intervals, 
\[-M\leq \inf_{x_{i-1}\leq x\leq x_i}h(x)\leq \sup_{x_{i-1}\leq x\leq x_i}h(x)\leq  M.\] If $[x_{i-1},x_i]$   does not contain any points of $S$, then 
\[\inf_{x_{i-1}\leq x\leq x_i}h(x)=\sup_{x_{i-1}\leq x\leq x_i}h(x)=0.\]
These imply that
\[U(h,P_n) =\sum_{i=1}^n \sup_{x_{i-1}\leq x\leq x_i}h(x) (x_{i}-x_{i-1})\leq  \frac{2Mk(b-a)}{n}.\]
\[L(h,P_n) =\sum_{i=1}^n \inf_{x_{i-1}\leq x\leq x_i}h(x) (x_{i}-x_{i-1})\geq   -\frac{2Mk(b-a)}{n}.\]
\bp
Therefore,
\[-\frac{2Mk(b-a)}{n}\leq L(h,P_n)\leq U(h,P_n)\leq \frac{2Mk(b-a)}{n}.\]
Taking $n\to \infty$ limits, we find that
\[\lim_{n\to\infty}  U(h,P_n)=\lim_{n\to\infty}L(h,P_n) =0.\]By the Archimedes-Riemann theorem, $h:[a,b]\to\mathbb{R}$
 is Riemann integrable and
$\di  \int_a^b h=0$.
Therefore, $g=h+f$ is also Riemann integrable and 
\[\int_a^b g=\int_a^bh+\int_a^bf=\int_a^b f.\]
\end{myproof}

\begin{remark}{}
If $f:(a,b)\to\mathbb{R}$ is a bounded function, we can extend the function to $[a,b]$ and discuss its integrability. By Theorem \ref{230222_2}, this is not affected by how we define the function at $x=a$ and $x=b$. In case the extension is Riemann integrable, we still denote the integral by $\di\int_a^b f$.
\end{remark}
\begin{definition}{Piecewise Continuous Functions}
We say that a function $f:[a,b]\to\mathbb{R}$ is piecewise continuous if there is a partition $P_0=\{a_0, a_1, \ldots, a_k\}$ of $[a,b]$ such that for each $1\leq i\leq k$, $f:(a_{i-1}, a_i)\to \mathbb{R}$ is continuous. 
\end{definition} 

Using the general additivity theorem (Corollary \ref{230221_9}) and Theorem \ref{230221_8}, we obtain the following immediately.
\begin{theorem}[label=230221_10]{}
Let $f:[a,b]\to\mathbb{R}$ be a function that is bounded and piecewise continuous. Then $f:[a,b]\to\mathbb{R}$ is Riemann integrable.
\end{theorem}

\begin{example}[label=ex230221_10]{}
The function $f:[-1, 2]\to\mathbb{R}$ defined by
\[f(x)=\begin{cases} 2-x,\quad &\text{if}\; -1\leq x<0,\\
x^2,\quad &\text{if}\; \quad 0\leq x\leq 2,\end{cases}\]is piecewise continuous and bounded. Hence,  $f:[-1, 2]\to\mathbb{R}$ is Riemann integrable.
\end{example}

 \begin{figure}[ht]
\centering
\includegraphics[scale=0.2]{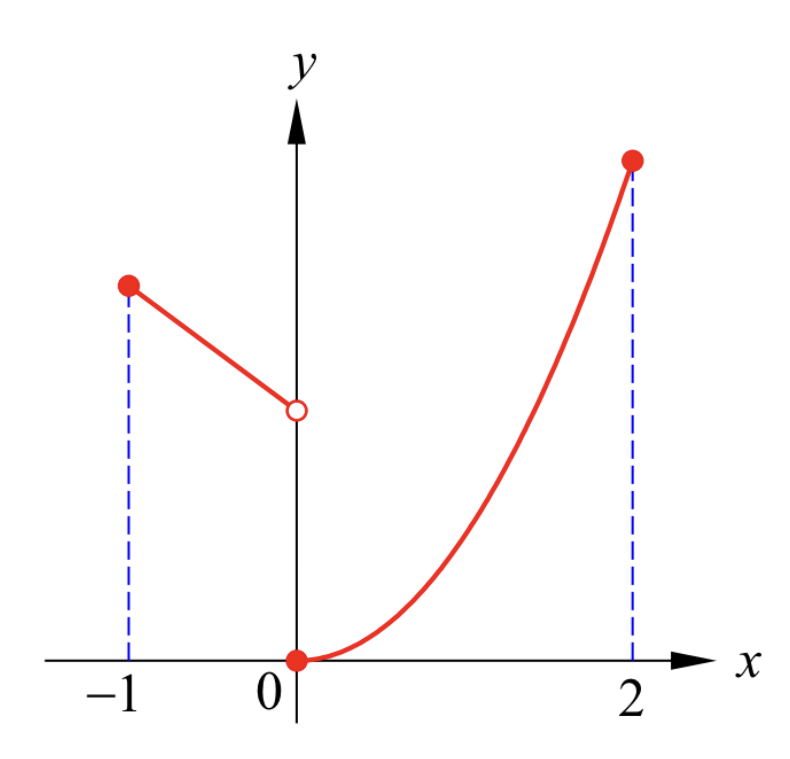}
\caption{The piecewise continuous function defined in Example \ref{ex230221_10}.\fa}\label{figure42}
\end{figure}
A special class of function that is bounded and piecewise continuous is the class of step functions. 
\begin{definition}{Step Functions}
We say that   $f:[a,b]\to\mathbb{R}$ is a step function if there is a partition $P_0=\{a_0, a_1, \ldots, a_k\}$ of $[a,b]$ such that for each $1\leq i\leq k$, $f:(a_{i-1}, a_i)\to\mathbb{R}$ is a constant function.
\end{definition}
By previous theorem, a step function is Riemann integrable. In fact, it is easy to compute its integral.
\begin{proposition}{}
Let $P_0=\{a_0, a_1, \ldots, a_k\}$ be a partition of $[a,b]$, and let $f:[a,b]\to\mathbb{R}$ be a step function such that for $1\leq i\leq k$,
\[f(x)=c_i,\quad\text{when}\quad a_{i-1}< x< a_i.\] Then $f:[a,b]\to\mathbb{R}$  is Riemann integrable  and
\[\int_a^b f=\sum_{i=1}^k c_i(a_i-a_{i-1}).\]
\end{proposition}

\begin{example}[label=ex230221_14]{}
Let $f:[0,5]\to\mathbb{R}$ be the function defined as 
\begin{align*}
f(x)=\begin{cases} 1,\quad &\text{if}\;0\leq x\leq \di 1,\\
\di\left\lfloor 5/x\right\rfloor,\quad &\text{if}\; \di 1< x\leq 5.\end{cases}
\end{align*}Show that $f$ is Riemann integrable and find $\di\int_0^5 f$. 
\end{example}

\begin{solution}{Solution}
The function $f$ is given explicitly by
\begin{align*}
f(x)=\begin{cases} 1,\quad &\text{if}\hspace{0.7cm} 0\leq x\leq \di 1,\\4,\quad &\text{if}\hspace{0.7cm} \di 1< x\leq 5/4,\\
3,\quad &\text{if}\;\; \di 5/4< x\leq 5/3,\\
2,\quad &\text{if}\; \;\di 5/3<x\leq 5/2,\\
1,\quad &\text{if}\; \;\di 5/2< x\leq 5.
 \end{cases}
\end{align*}
This is a step function. Hence, it is integrable, and
\[\int_0^5f=1\times 1+4\times \frac{1}{4}+3\times \frac{5}{12}+2\times\frac{5}{6}+1\times\frac{5}{2}=\frac{89}{12}.\]
\end{solution}
 \begin{figure}[ht]
\centering
\includegraphics[scale=0.2]{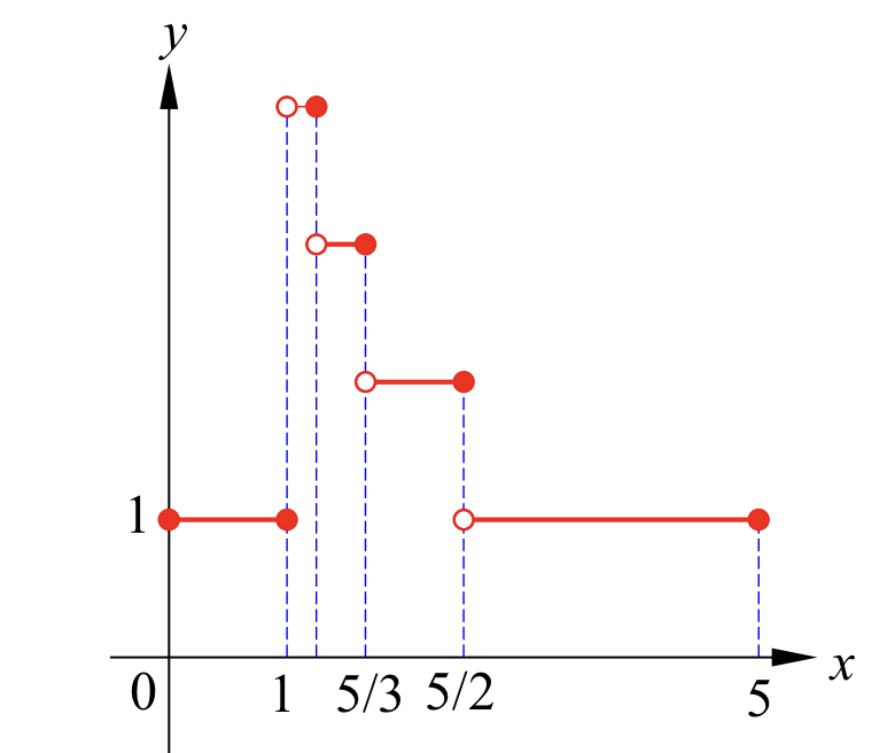}
\caption{The   function defined in Example \ref{ex230221_14}.\fa}\label{figure43}
\end{figure}

The following theorem shows that  monotonic functions are also Riemann integrable.
\begin{theorem}{}
If  $f:[a,b]\to\mathbb{R}$   is a monotonic function, then   it  is Riemann integrable.
\end{theorem}
\begin{myproof}{Proof} Without loss of generality, assume that $f:[a,b]\to\mathbb{R}$ is an increasing function. If $P=\{x_i\}_{i=0}^k$ is a partition of $[a,b]$, then for any $1\leq i\leq k$,
\[m_i=\inf_{x_{i-1}\leq x\leq x_i}f(x)=f(x_{i-1}),\hspace{1cm}M_i=\inf_{x_{i-1}\leq x\leq x_i}f(x)=f(x_i).\]
Therefore,
\[U(f,P)-L(f,P)=\sum_{i=1}^n\left(f(x_i)-f(x_{i-1}\right)(x_i-x_{i-1}).\]
For each positive integer $n$, let $P_n$ be the regular partition of $[a, b]$ into $n$ intervals. Then
\begin{align*}U(f, P_n)-L(f, P_n)&=\frac{b-a}{n}\sum_{i=1}^n\left(f(x_i)-f(x_{i-1})\right)\\&=\frac{(b-a)\left(f(b)-f(a)\right)}{n}.\end{align*}
\bp
This implies that
\[\lim_{n\to \infty}\left(U(f, P_n)-L(f, P_n)\right)=\lim_{n\to \infty}\frac{(b-a)\left(f(b)-f(a)\right)}{n}=0.\]
In other words, $\{P_n\}$ is an Archimedes sequence of partitions for $f$. By the Archimedes-Riemann theorem, this proves that $f$ is Riemann integrable.
\end{myproof}
\begin{example}{}
Let $f:[0,1]\to\mathbb{R}$ be the function defined by $f(0)=1$, and for each positive integer $n$, 
\[f(x)=1-\frac{1}{n+1},\hspace{1cm}\text{when}\; \frac{1}{n+1}<x\leq\frac{1}{n}.\]
One can verify that $f:[0,1]\to\mathbb{R}$ is a decreasing function. Hence, it is Riemann integrable.
However, $f$ is not a piecewise continuous function, since it has discontinuities at infinitely many points. 
\end{example}

The Riemann integrability of a function implies the Riemann integrability of its absolute value.
\begin{theorem}[label=230221_15]{}
Let $f:[a,b]\to\mathbb{R}$   be a bounded function. If  $f:[a,b]\to\mathbb{R}$ is Riemann integrable, then the function  $|f|:[a,b]\to\mathbb{R}$ is Riemann integrable.
\end{theorem}

\begin{myproof}{Proof}We will first prove the following: For any $c$ and $d$ in $[a,b]$ with $c<d$,
\begin{equation}\label{eq230221_19}\sup_{c\leq x\leq d}|f(x)|-\inf_{c\leq x\leq d}|f(x)|\leq  \sup_{c\leq x\leq d}f(x) -\inf_{c\leq x\leq d}f(x).\end{equation}
 There are two sequences of points $\{u_n\}$ and $\{v_n\}$ in $[c,d]$ such that
\[\lim_{n\to\infty}|f(u_n)|=\inf_{c\leq x\leq d}|f(x)|,\hspace{1cm}\lim_{n\to\infty}|f(v_n)|=\sup_{c\leq x\leq d}|f(x)|.\]\bp
Since $u_n$ and $v_n$ are points in $[c,d]$, we find that
\[|f(v_n)|-|f(u_n)|\leq|f(v_n)-f(u_n)|\leq \sup_{c\leq x\leq d}f(x) -\inf_{c\leq x\leq d}f(x).\]
Passing to the $n\to\infty$ limit, we obtain \eqref{eq230221_19}.

Now, given $\varepsilon>0$, since $f:[a,b]\to\mathbb{R}$ is Riemann integrable, there is a partition $P=\{x_i\}_{i=0}^k$ of $[a,b]$ such that 
\[U(f,P)-L(f,P) <\varepsilon.\]
But then \begin{align*}
&U(|f|,P)-L(|f|,P)\\&=\sum_{i=1}^n\left(\sup_{x_{i-1}\leq x\leq x_i}|f(x)|-\inf_{x_{i-1}\leq x\leq x_i}|f(x)|\right)(x_i-x_{i-1})\\
&\leq \sum_{i=1}^n\left(\sup_{x_{i-1}\leq x\leq x_i}f(x)-\inf_{x_{i-1}\leq x\leq x_i}f(x)\right)(x_i-x_{i-1})\\&=U(f,P)-L(f,P)<\varepsilon.\end{align*}This prove that $|f|:[a,b]\to\mathbb{R}$ is Riemann integrable.
\end{myproof}

\begin{figure}[ht]
\centering
\includegraphics[scale=0.2]{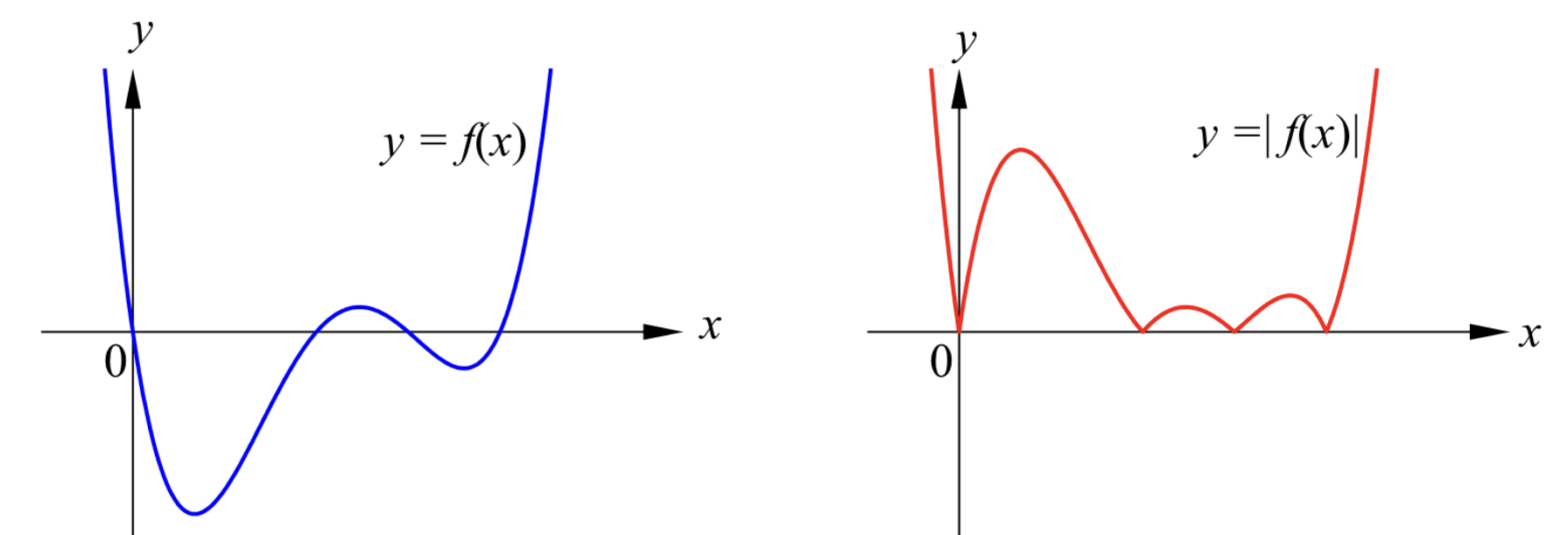}
\caption{A function $y=f(x)$ and its absolute value $y=|f(x)|$.\fa}\label{figure44}
\end{figure}
\begin{remark}{}
The converse of Theorem \ref{230221_15} is not true. Namely, if a function $f:[a,b]\to\mathbb{R}$ is bounded, $|f|:[a,b]\to\mathbb{R}$ is Riemann integrable does not imply that $f:[a,b]\to\mathbb{R}$ is Riemann integrable. For a counter example, consider the function $f:[0,1]\to\mathbb{R}$ defined as
\[f(x)=\begin{cases} 1,\quad &\text{if $x$ is rational},\\-1,\quad & \text{if $x$ is irrational}.\end{cases} \]
One can prove that $f:[0,1]\to\mathbb{R}$ is not integrable, exactly the same way as in Example \ref{230220_9}. On the other hand, since $|f|:[0,1]\to\mathbb{R}$ is a constant function, it is Riemann integrable.
\end{remark}

 The following theorem says that the product of Riemann integrable functions is Riemann integrable.
\begin{theorem}{}
Let $f:[a,b]\to\mathbb{R}$ and  $g:[a,b]\to\mathbb{R}$ be bounded functions. If  $f:[a,b]\to\mathbb{R}$ and  $g:[a,b]\to\mathbb{R}$ are Riemann integrable, then the function  $(fg):[a,b]\to\mathbb{R}$ is also Riemann integrable.
\end{theorem}
\begin{myproof}{Proof}We will apply Lemma \ref{230220_10} to prove the Riemann integrability of the function  $h=(fg):[a,b]\to\mathbb{R}$.
Since  $f:[a,b]\to\mathbb{R}$ and  $g:[a,b]\to\mathbb{R}$ are bounded functions, there is a positive number $M$ so that 
\[|f(x)|\leq M\quad\text{and}\quad|g(x)|\leq M\hspace{1cm}\text{for all}\;x\in [a,b].\]
We claim that for any $c$ and $d$ in $[a,b]$ with $c<d$,
\begin{equation}\label{230524_1}\begin{split}&\sup_{c\leq x\leq d}h(x)-\inf_{c\leq x\leq d}h(x)\\&\leq M\left(\sup_{c\leq x\leq d}f(x)-\inf_{c\leq x\leq d}f(x)+\sup_{c\leq x\leq d}g(x)-\inf_{c\leq x\leq d}g(x)\right).\end{split}\end{equation}\bp
 There are two sequences of points $\{u_n\}$ and $\{v_n\}$ in $[c,d]$ such that
\[\lim_{n\to\infty}h(u_n)=\inf_{c\leq x\leq d}h(x),\hspace{1cm}\lim_{n\to\infty}h(v_n)=\sup_{c\leq x\leq d}h(x).\] 
Notice that
\begin{align*}
|h(v_n)-h(u_n)|&=|g(v_n)(f(v_n)-f(u_n))+f(u_n)(g(v_n)-g(u_n))|\\
&\leq |g(v_n)||f(v_n)-f(u_n)|+|f(u_n)||g(v_n)-g(u_n)|\end{align*}  Since $u_n$ and $v_n$ are in $[c,d]$, we find that
\[|f(v_n)-f(u_n)|\leq \sup_{c\leq  x\leq d}f(x)-\inf_{c\leq x\leq d}f(x),\]
\[|g(v_n)-g(u_n)|\leq \sup_{c\leq  x\leq d}g(x)-\inf_{c\leq x\leq d}g(x).\]
Therefore,
\begin{align*}
&|h(v_n)-h(u_n)|\\&\leq M\left(\sup_{c\leq  x\leq d}f(x)-\inf_{c\leq x\leq d}f(x)+\sup_{c\leq x\leq d}g(x)-\inf_{c\leq x\leq d}g(x)\right).
\end{align*}Passing to the $n\to\infty$ limit gives \eqref{230524_1}.

Now given $\varepsilon>0$, there are partitions $P_1$ and $P_2$ of $[a,b]$ such that
\begin{gather*}U(f,P_1)-L(f,P_1)<\frac{\varepsilon}{2M},\\
U(g,P_2)-L(g,P_2)<\frac{\varepsilon}{2M}.\end{gather*}
Let $P^*=\{x_0, x_1, \ldots, x_k\}$ be  a common refinement of $P_1$ and $P_2$. Then
\begin{gather*}
U(f,P^*)-L(f,P^*)<\frac{\varepsilon}{2M},\\
U(g,P^*)-L(g,P^*)<\frac{\varepsilon}{2M}.\end{gather*}

\bp
 It follows that
\begin{align*}
&U(h,P^*)-L(h,P^*)\\&=\sum_{i=1}^k\left(\sup_{x_{i-1}\leq x\leq x_i}h(x)-\inf_{x_{i-1}\leq x\leq x_i}h(x)\right)(x_i-x_{i-1})\\
&\leq M\sum_{i=1}^k\left(\sup_{x_{i-1}\leq x\leq x_i}f(x)-\inf_{x_{i-1}\leq x\leq x_i}f(x)\right)(x_i-x_{i-1})\\
&\quad +M\sum_{i=1}^k\left(\sup_{x_{i-1}\leq x\leq x_i}g(x)-\inf_{x_{i-1}\leq x\leq x_i}g(x)\right)(x_i-x_{i-1})\\
&=M\left(U(f,P^*)-L(f,P^*)+U(g,P^*)-L(g,P^*)\right)\\
&<M\times\left(\frac{\varepsilon}{2M}+\frac{\varepsilon}{2M}\right)=\varepsilon.
\end{align*}This proves that $(fg):[a,b]\to\mathbb{R}$ is   Riemann integrable.
\end{myproof}
\vp
\noindent
{\bf \large Exercises  \thesection}
\setcounter{myquestion}{1}
 
\begin{question}{\themyquestion}
Given that $f:[-1,1]\to\mathbb{R}$ is the function defined by
\[f(x)=\begin{cases} \di \sin\left(\frac{4}{x}\right),\quad &\text{if}\;x\neq 0,\\
0,\quad &\text{if}\;x= 0.\end{cases}\]Explain why $f:[-1,1]\to\mathbb{R}$ is Riemann integrable.
\end{question}
\atc

\begin{question}{\themyquestion}
Let $f:[a,b]\to\mathbb{R}$ be a bounded function. If $f:[a,b]\to\mathbb{R}$ is Riemann integrable, show that
\[\left|\int_a^b f(x)dx\right|\leq \int_a^b |f(x)|dx.\]
\end{question}
\atc
\begin{question}{\themyquestion}
Let $f:[0,6]\to\mathbb{R}$ be the function defined as 
\begin{align*}
f(x)=\begin{cases} -2,\quad &\text{if}\;0\leq x<\di 1,\\
\di\left\lfloor 4/x\right\rfloor,\quad &\text{if}\; \di 1\leq x\leq 6.\end{cases}
\end{align*}Show that $f$ is Riemann integrable and find $\di\int_0^6 f$. 
\end{question}
\atc

\begin{question}{\themyquestion}
Let $f:\mathbb{R}\to\mathbb{R}$ be the function defined as
\[f(x)=\left(x-\lfloor x\rfloor\right)^2.\]
Explain why the function $f$ is Riemann integrable on any closed and bounded interval $[a,b]$.
\end{question}

\atc

\begin{question}[label=ex230224_7]{\themyquestion}
Let $f:[a,b]\to\mathbb{R}$ and $g:[a,b]\to\mathbb{R}$ be bounded functions. Define the function
$h:[a,b]\to\mathbb{R}$ by
\[h(x)=\max\{f(x), g(x)\}.\]
\begin{enumerate}[(a)]
\item
Show that
\[h=\frac{f+g+|f-g|}{2}.\]
\item If  $f:[a,b]\to\mathbb{R}$ and $g:[a,b]\to\mathbb{R}$ are Riemann integrable, show that $h:[a,b]\to\mathbb{R}$ is also Riemann integrable.
\end{enumerate}
\end{question}
\atc
 \begin{question}{\themyquestion\;\;[Cauchy Schwarz Inequality]}
Let $f:[a,b]\to\mathbb{R}$ and $g:[a,b]\to\mathbb{R}$ be bounded functions that are Riemann integrable. Prove that
\[\left(\int_a^b f(x)g(x)dx\right)^2\leq \left(\int_a^bf(x)^2dx\right)\left(\int_a^bg(x)^2dx\right).\]
\end{question}
\vp
\section{The Fundamental Theorem of Calculus}\label{sec4.4}
In this section, we prove the fundamental theorem of calculus, which gives a relation between integration and differentiation. It also provides a useful method to compute integrals of certain functions. 
We first prove a few results about integrals.

 Given that a function $f:[a,b]\to\mathbb{R}$ is bounded and Riemann integrable on $[a,b]$, it is Riemann integrable on any interval $[c,d]$ that is contained in the interval $[a,b]$. Thus,  we can define a new function $F:[a,b]\to\mathbb{R}$ by
\[F(x)=\int_a^x f(u)du.\]
By definition, $F(a)=0$. 
For any $c$ and $d$ in $[a, b]$, one can check that
\[F(d)-F(c)=\int_a^df(u)du-\int_a^cf(u)du=\int_c^df(u)du.\]
The followng theorem says that  $F:[a,b]\to\mathbb{R}$ is a continuous function.
\begin{theorem}[label=230222_11]{}
Let $f:[a,b]\to\mathbb{R}$ be a bounded function that is Riemann integrable, and let $F:[a,b]\to\mathbb{R}$ be the function defined by \[F(x)=\int_a^x f(u)du.\]Then $F:[a,b]\to\mathbb{R}$ is a Lipschitz function, and hence it is continuous.
\end{theorem}
\begin{myproof}{Proof}
It is sufficient to prove that $F:[a,b]\to\mathbb{R}$ is Lipschitz. The continuity follows. Since $f:[a,b]\to\mathbb{R}$ is bounded, there is a positive constant $M$ such that
\[|f(x)|\leq M\hspace{1cm}\text{for all}\;x\in [a,b].\]\bp

For any $x_1$ and $x_2$ in $[a,b]$ with $x_1<x_2$, we have
\[ F(x_2)-F(x_1)= \int_{x_1}^{x_2}f(u)du.\]
Therefore,
\begin{align*}
|F(x_2)-F(x_1)|&=\left|\int_{x_1}^{x_2}f(u)du\right|\leq \int_{x_1}^{x_2}|f(u)|du\\&\leq \int_{x_1}^{x_2}Mdu=M(x_2-x_1)=M|x_2-x_1|.
\end{align*}This proves that $F:[a,b]\to\mathbb{R}$ is a Lipschitz function wth Lipschitz constant $M$. 
\end{myproof}

The next is a mean value theorem for integrals.
\begin{theorem}{Mean Value Theorem for Integrals}
Let $f:[a,b]\to\mathbb{R}$ be a continuous function. Then there exists $c$ in $[a,b]$ such that
\[\frac{1}{b-a}\int_a^bf(x)dx=f(c).\]
\end{theorem}This theorem is known as the mean value theorem since 
\[\frac{1}{b-a}\int_a^bf(x)dx\]can be interpreted as the average of the values of $f$ over the interval $[a,b]$. 
\begin{myproof}{Proof}
Since $f:[a,b]\to\mathbb{R}$ is a continuous function, the extreme value theorem says that there are points $u$ and $v$ in $[a,b]$ such that
\[f(u)\leq f(x)\leq f(v)\hspace{1cm}\text{for all}\;x\in [a,b].\]This implies that
\[f(u)(b-a)\leq \int_a^b f(x)dx\leq f(v)(b-a).\]\bp

Therefore, the number
\[w=\frac{1}{b-a}\int_a^bf(x)dx\]satisfies
\[f(u)\leq w\leq f(v).\]
By intermediate value theorem, there is a point $c$ in $[a,b]$ such that $f(c)=w$. This gives
\[\frac{1}{b-a}\int_a^bf(x)dx=f(c).\]
\end{myproof}
In fact, one can argue that the number $c$ can be chosen to be in  $(a,b)$.  

\begin{example}[label=ex230222_5]{}
Let $f:[a,b]\to\mathbb{R}$ be a continuous function. If  $f(x)\geq m$ for all $x\in [a,b]$ and 
\[\int_a^b f(x)dx=m(b-a),\]prove that
\[f(x)=m\hspace{1cm}\text{for all}\;x\in [a,b].\]
\end{example}
\begin{solution}{Solution}
Let $g:[a,b]\to\mathbb{R}$ be the function $g(x)=f(x)-m$. Then $g:[a,b]\to\mathbb{R}$ is a continuous function and $g(x)\geq 0$ for all $x\in [a,b]$. Moreover,
\[\int_a^b g(x)dx=\int_a^b f(x)dx-\int_a^b mdx=0.\]
We want to show that $g(x)=0$ for all $x\in [a,b]$. Suppose to the contrary that there is a point $x_0$ in $[a,b]$ such that $g(x_0)\neq 0$. Then $g(x_0)>0$. Since $g$ is continuous, there exists a $\delta>0$ such that for all $x\in (x_0-\delta,x_0+\delta)\cap (a,b)$, 
\[g(x)>\frac{g(x_0)}{2}.\]\bs
Without loss of generality, assume that $\delta<\di\frac{b-a}{2}$. Then either $x_0-\delta>a$ or $x_0+\delta<b$. In any case, $(x_0-\delta,x_0+\delta)\cap (a,b)=(c,d)$ is an interval of length at least $\delta$. But then
\begin{align*}\int_a^b g(x)dx&=\int_{a}^cg(x)dx+\int_c^dg(x)dx+\int_d^bg(x)dx\\&\geq 0\times (c-a)+(d-c)\frac{g(x_0)}{2}+0\times (b-d)\\&=(d-c)\frac{g(x_0)}{2}>0,\end{align*}
which is a contradiction. Therefore, we must have $g(x)=0$ for all $x\in [a,b]$.
\end{solution}

\begin{remark}{}
In the mean value theorem for integrals, we can strengthen the theorem to have the point $c$  being a point in the open interval $(a,b)$. In the proof, we have shown that for $w=\di\int_a^b f(x)dx$,
 \[f(u)\leq w\leq f(v).\] Here $f(u)$ is the minimum value of $f:[a,b]\to\mathbb{R}$, and $f(v)$ is the maximum value. By Example \ref{ex230222_5}, if $w=f(u)$, then $f$ is a constant. In this case, we can take $c$ to be any point in $(a,b)$. If $w=f(v)$, the same reasoning as in Example \ref{ex230222_5} also shows that $f$ is a constant, and so $c$ also can be any point in $(a,b)$. If $w\neq f(u)$ and $w\neq f(v)$, then $c$ is a point strictly between $u$ and $v$, and thus it is strictly between $a$ and $b$.
\end{remark}

Now we turn to the fundamental theorem of calculus. Consider the case that an object is moving with speed $v(t)$ at time $t$. To find $s(t)$, the distance  travelled up to time $t$, we can partition the time interval $[0, t]$ into a finite number of subintervals $[0, t_1], [t_1, t_2], \ldots, [t_{k-1}, t_k]$, where $t_k=t$. For each time interval $[t_{i-1}, t_i]$, where $1\leq i\leq k$, take a point $t_i^*\in [t_{i-1}, t_i]$, and approximate the average speed over the time interval $[t_{i-1}, t_i]$ by the speed at time $t_i^*$, $v(t_i^*)$. Then the distance travelled up to time $t$ is approximately 
\[\sum_{i=1}^k v(t_i^*)(t_i-t_{i-1}).\]
We recognize that this is  a Riemann sum of the speed function $v(t)$. The distance travelled $s(t)$ should be calculated as the limit where the gap of the partition goes to zero. In other words,
\[s(t)=\int_0^t v(\tau)d\tau.\]In Chapter \ref{ch3}, we have motivated that $v(t)=s'(t)$. Hence, 
\[s(t)=\int_0^ts'(\tau)d\tau,\]which means that differentiation and integration are inverse processes of each other. The fundamental theorem of calculus gives a rigorous setting for this.
\begin{theorem}{Fundamental Theorem of Calculus I}
Let $f:[a,b]\to\mathbb{R}$ be a bounded function that is Riemann integrable, and let $F:[a,b]\to\mathbb{R}$ be the function defined by
\[F(x)=\int_a^xf(u)du.\]
If $x_0$ is a point in $ (a,b)$, and $f$ is continuous at $x_0$, then $F$ is differentiable at $x_0$ and
\[F'(x_0)=f(x_0).\]
\end{theorem}
 
\begin{myproof}{Proof of Fundamental Theorem of Calculus I}
We need to show that the limit
\[\lim_{h\to 0}\frac{F(x_0+h)-F(x_0)}{h}\] exists and is equal to $f(x_0)$. Given $\varepsilon>0$, since $f$ is continuous at $x_0$, there is a $\delta>0$ such that $(x_0-\delta,x_0+\delta)\subset (a,b)$ and 
\begin{equation}\label{eq230222_10}\left|f(x)-f(x_0)\right|<\frac{\varepsilon}{2}\hspace{1cm}\text{for all}\; x\in (x_0-\delta, x_0+\delta).\end{equation}
\bp
For $h\in (-\delta, \delta)$,
\begin{align*}
F(x_0+h)-F(x_0)-f(x_0)h&=\int_{x_0}^{x_0+h}f(u)du-f(x_0)h\\&=\int_{x_0}^{x_0+h}\left(f(u)-f(x_0)\right)du.\end{align*}
Eq. \eqref{eq230222_10} implies that
\[\left|F(x_0+h)-F(x_0)-f(x_0)h\right|\leq\frac{\varepsilon}{2}|h|.\]
Hence, if $h\in (-\delta, \delta)\setminus\{0\}$, 
\[\left|\frac{F(x_0+h)-F(x_0)}{h}-f(x_0) \right|\leq\frac{\varepsilon}{2}<\varepsilon.\]
This proves that
\[\lim_{h\to 0}\frac{F(x_0+h)-F(x_0)}{h}=f(x_0).\]
\end{myproof}
In fact, the point $x_0$ can be $a$ or $b$ if we consider one-sided derivatives.

\begin{example}[label=20230527]{}
Consider the piecewise continuous function $f:[0,2]\to\mathbb{R}$ given by 
\[f(x)=\begin{cases} -1,\quad&\text{if}\; 0\leq x<1,\\
2,\quad&\text{if}\; 1\leq x<2,\\
1,\quad &\text{if}\;\quad x=2.\end{cases}\]
We find that
\[F(x)=\int_0^x f(u)du=\begin{cases} -x,\quad&\text{if}\; 0\leq x<1,\\
2x-3,\quad&\text{if}\; 1\leq x\leq 2.\end{cases}\]\be
The function $F:[0,2]\to\mathbb{R}$ is   continuous.
For $x\in (0,1)$, $F$ is differentiable and $F'(x)=-1=f(x)$. For  $x\in (1,2)$, $F$ is differentiable and $F'(x)=2$. However, $F$ is not differentiable at $x=1$ since
\[\lim_{x\to 1^-}\frac{F(x)-F(1)}{x-1}=-1,\hspace{1cm}\lim_{x\to 1^+}\frac{F(x)-F(1)}{x-1}=2.\]
\end{example2}
\begin{figure}[ht]
\centering
\includegraphics[scale=0.2]{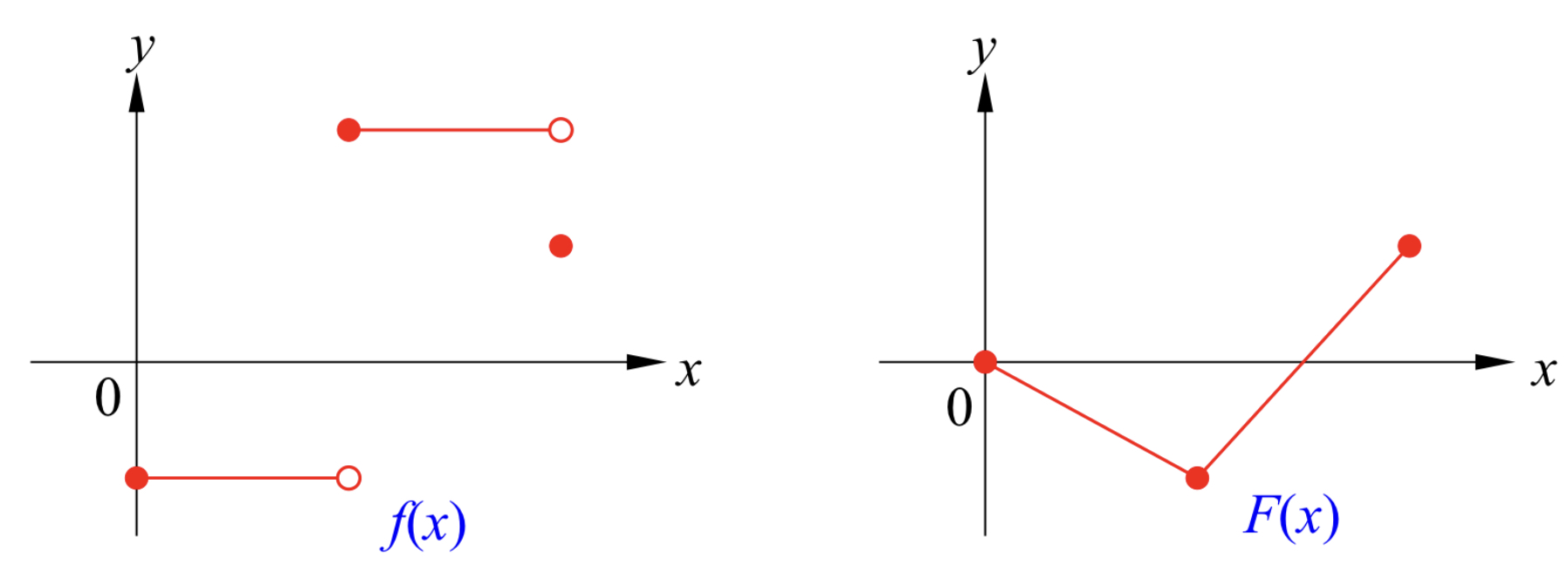}
\caption{The function  $f:[0,2]\to\mathbb{R}$  and  $F:[0,2]\to\mathbb{R}$ in Example \ref{20230527}.}\label{figure69}
\end{figure}
\begin{example}{}
Evaluate the following derivatives.
\begin{enumerate}[(a)]
\item
$\di\frac{d}{dx}\int_0^x\sin(u^2)du$
\item
$\di\frac{d}{dx}\int_x^{1}\sin(u^2)du$
\item
$\di\frac{d}{dx}\int_x^{x^2}\sin(u^2)du$
\end{enumerate}
\end{example}
\begin{solution}{Solution}
The function $f:\mathbb{R}\to\mathbb{R}$, $f(x)=\sin(x^2)$ is continuous. Hence, it is Riemann integrable over any closed and bounded intervals. Let
\[F(x)=\int_0^xf(u)du=\int_0^x\sin(u^2)du.\]\bs
 By fundamental theorem of calculus I, $F'(x)=f(x)=\sin(x^2)$.
\begin{enumerate}[(a)]
\item$\di \frac{d}{dx}\int_0^x\sin(u^2)du=\frac{d}{dx}F(x) =\sin (x^2)$.
\item $\di \frac{d}{dx}\int_x^{1}\sin(u^2)du=\frac{d}{dx}\left(F(1)-F(x)\right) =-F'(x)= -\sin(x^2)$.
\item $\di \frac{d}{dx}\int_x^{x^2}\sin(u^2)du=\frac{d}{dx}\left(F(x^2)-F(x)\right)=2xF'(x^2)-F'(x) $\\~$~\hspace{5cm}~$$=2x\sin(x^4)-\sin(x^2)$.
\end{enumerate}

\end{solution}

Now we turn to the second fundamental theorem of calculus, which provides a mean for   calculating integrals of continuous functions.
\begin{theorem}{Fundamental Theorem of Calculus II}
Let $F:[a,b]\to\mathbb{R}$ be a continuous function, and let $f:[a,b]\to\mathbb{R}$ be a bounded function that is continuous on $(a,b)$. If 
\[F'(x)=f(x)\hspace{1cm}\text{for all}\;x\in (a,b),\]
Then
\[\int_a^b f(x)dx=F(b)-F(a).\]
\end{theorem}
Recall that if the functions $F(x)$ and $f(x)$ are related by
\[F'(x)=f(x)\hspace{1cm}\text{for all}\;x\in (a,b),\]$F$ is called an antiderivative of $f$. Hence, the fundamental theorem of calculus II states that if the function $f:(a,b)\to\mathbb{R}$ is continuous and it has an antiderivative $F(x)$ which can be extended to a continuous function $F:[a,b]\to\mathbb{R}$, then 
\[\int_a^bf(x)dx=\left[ F(x)\right]_{ a}^{ b}=F(b)-F(a).\]
We will present two proofs of the fundamental theorem of calculus II. The first one uses fundamental theorem of calculus I.
\begin{myproof}{First  Proof of Fundamental Theorem of Calculus II}
Since $f:[a,b]\to\mathbb{R}$ is a bounded function that is continuous on $(a,b)$, it is integrable over any subinterval of $[a,b]$. Let $G:[a,b]\to\mathbb{R}$ be the function defined by
\[G(x)=\int_a^xf(u)du.\]By default, $G(a)=0$.  By Theorem \ref{230222_11}, $G$ is continuous on $[a,b]$.  By fundamental theorem of calculus I, 
\[G'(x)=f(x)\hspace{1cm} \text{for all}\;x\in (a,b).\]Hence, $F:[a,b]\to\mathbb{R}$ and $G:[a,b]\to\mathbb{R}$ are continuous functions satisfying
\[G'(x)=F'(x)\hspace{1cm} \text{for all}\;x\in (a,b).\] By Theorem \ref{thm230215_3}, there is a constant $C$ such that 
\[G(x)=F(x)+C.\]
Since $G(a)=0$, we find that 
$C=-F(a)$. Therefore, 
\[G(x)=F(x)-F(a),\]and so
\[\int_a^bf(x)dx=G(b)=F(b)-F(a).\]
\end{myproof}
In different textbooks, the ordering of the two fundamental theorem of calculus might be different. One can use one to deduce the other. This is why a   proof of the fundamental theorem of calculus II without using the fundamental theorem of calculus I is of interest.

\begin{myproof}{Second Proof of Fundamental Theorem of Calculus II}
Here we use the Lagrange mean value theorem. Since $f:[a,b]\to\mathbb{R}$ is a bounded function that is continuous on $(a,b)$, it is integrable. Now if $P=\{x_i\}_{i=0}^k$ is a partition of $[a,b]$, for each $1\leq i\leq k$, the mean value theorem implies that there is a $\xi_i\in (x_{i-1}, x_i)$ such that
\[F(x_i)-F(x_{i-1})=F'(\xi_i)(x_i-x_{i-1})=f(\xi_i)(x_i-x_{i-1}).\]
Summing over $i$ gives
\[F(b)-F(a)=\sum_{i=1}^k\left(F(x_i)-F(x_{i-1})\right)=\sum_{i=1}^k f(\xi_i)(x_i-x_{i-1})=R(f,P,A),\]where $A=\{\xi_i\}_{i=1}^k$. Since
\[L(f,P)\leq R(f,P,A)\leq U(f,P),\]we find that
\[L(f,P)\leq F(b)-F(a)\leq U(f,P).\]
Notice that this is true for any partition $P$ of $[a,b]$. By definitions of the lower integral and the upper integral, we find that
\[\underline{\int_a^b}f\;\leq\;F(b)-F(a)\;\leq\;\overline{\int_a^b}f.\]
Since $f:[a,b]\to\mathbb{R}$ is Riemann integrable, the lower integral and the upper integral are the same. Thus,
\[\int_a^b f=\underline{\int_a^b}f=\overline{\int_a^b}f=F(b)-F(a).\]

\end{myproof}

We can relax the conditions in the fundamental theorem of calculus II to let $f$ to be a piecewise continuous function.
\begin{corollary}{Generalized Fundamental Theorem of Calculus II}
 Let $S=\{a_0, a_1, \ldots, a_k\}$ be a finite subset of $[a,b]$ that contains $a$ and $b$, and let $f:[a,b]\to\mathbb{R}$ be a bounded function that is continuous on $[a,b]\setminus S$. If $F:[a,b]\to\mathbb{R}$ is a continuous function, differentiable on $[a,b]\setminus S$, and
\[F'(x)=f(x)\hspace{1cm}\text{for all}\;x\in [a,b]\setminus S,\] then
\[\int_a^b f(x)dx=F(b)-F(a).\]
\end{corollary}
\begin{myproof}{Proof}
We can assume that
\[a=a_0<a_1<\ldots<a_k=b.\]
Since $f:[a,b]\to\mathbb{R}$ is bounded and piecewise continuous, it is Riemann integrable. Moreover, by the generalized additivity theorem, we have
\[\int_a^bf(x)dx=\sum_{i=1}^k\int_{a_{i-1}}^{a_i}f(x)dx.\]
Applying the fundamental theorem of calculus II to each of the integrals $\di \int_{a_{i-1}}^{a_i}f(x)dx$, we find that
\[\int_a^b f(x)dx=\sum_{i=1}^k (F(a_i)-F(a_{i-1}))=F(b)-F(a).\]This completes the proof.
Note that it is crucial here that $F$ is continuous on $[a,b]$.
\end{myproof}

As is well known, the fundamental theorem of calculus provides a practical method for computing integrals of  functions that have antiderivatives. 
\begin{example}
{}
Compute the integral of the piecewise continuous function   $f:[-1, 2]\to\mathbb{R}$,
\[f(x)=\begin{cases} 2-x,\quad &\text{if}\; -1\leq x<0,\\
x^2,\quad &\text{if}\; \quad 0\leq x\leq 2,\end{cases}\] that is defined in Example \ref{ex230221_10}. 
\end{example}
\begin{solution}{Solution}
Using additivity,
\[\int_{-1}^2f(x)dx=\int_{-1}^0f(x)dx+\int_0^2f(x)dx.\]
Using fundamental theorem of calculus II, 
\[\int_{-1}^0f(x)dx=\int_{-1}^0(2-x)dx=\left[2x-\frac{x^2}{2}\right]_{-1}^0=0-\left(-\frac{5}{2}\right)=\frac{5}{2},\]
\[\int_0^2f(x)dx=\int_0^2x^2dx=\left[\frac{x^3}{3}\right]_0^2=\frac{8}{3}-0=\frac{8}{3}.\]
Hence,
\[\int_{-1}^2f(x)dx=\frac{5}{2}+\frac{8}{3}=\frac{31}{6}.\]

\end{solution}

\begin{remark}{Alternative Proof of Mean Value Theorem for Integrals}
Using the fundamental theorem of calculus, we can give an alternative proof of the mean value theorem for integrals as follows. Since the function $f:[a,b]\to\mathbb{R}$ is continuous, the function $F:[a,b]\to\mathbb{R}$ defined by
\[F(x)=\int_a^xf(u)du\] is continuous on $[a,b]$, differentiable on $(a,b)$, and $F'(x)=f(x)$ for all $x\in (a,b)$. By Lagrange mean value theorem, there is a $c\in (a,b)$ such that
\[\frac{1}{b-a}\int_a^bf(x)dx=\frac{F(b)-F(a)}{b-a}=F'(c)=f(c).\]
 
\end{remark}

Finally, we can prove the existence and uniqueness theorem mentioned in Chapter \ref{ch3}, Theorem \ref{thm230217_3}.
\begin{theorem}[label=thm230222_13]{Existence and Uniqueness Theorem}
Let $(a,b)$ be an open interval that contains the point $x_0$, and let $y_0$ be any real number.  Given that $f:(a,b)\to\mathbb{R}$ is a continuous function, there exists a unique differentiable function $F:(a,b)\to\mathbb{R}$ such that \[F'(x)=f(x)\quad \text{for all}\;x\in (a, b), \hspace{1cm}F(x_0)=y_0.\]
\end{theorem}
\begin{myproof}{Proof}
As we mentioned before, the uniqueness follows from the identity criterion. For the existence, notice that $f$ is continuous on any closed and bounded interval that is contained in $(a,b)$. Hence, we can define the function $F:(a,b)\to\mathbb{R}$ by
\[F(x)=\int_{x_0}^xf(u)du+y_0.\] 
Then $F(x_0)=y_0$ by default. By fundamental theorem of calculus, $F'(x)=f(x)$ for all $x\in (a,b)$.
\end{myproof}

Let us look at some other examples how integrals can be applied.  
\begin{example}{}
Find the limit 
\[\lim_{n\to\infty}\frac{1}{n}\sum_{k=1}^n\sin\left(\frac{\pi k}{n}\right).\]
\end{example}
\begin{solution}{Solution}
We try to identify \[\frac{1}{n}\sum_{k=1}^n\sin\left(\frac{\pi k}{n}\right)\] as a Riemann sum. For $1\leq k\leq n$, let $\xi_k=\di\frac{\pi k}{n}$.  These are equally spaced points in the interval $[0, \pi]$. This motivates us to define the function $f:[0,\pi]\to\mathbb{R}$, $f(x)=\sin x$. Since $f$ is a continuous function, it is Riemann integrable. Let $P_n$ be the regular partition of $[0,\pi]$ into $n$ intervals. Then with $A_n=\{\xi_k\}_{k=1}^n$, we have
\[R(f,P_n,A_n)=\sum_{k=1}^n\sin \left(\frac{\pi k}{n}\right)\frac{\pi }{n}.\]Since $f$ is Riemann integrable,
\[\lim_{n\to\infty}R(f, P_n, A_n)=\int_0^{\pi}f(x)dx.\]
By fundamental theorem of calculus,
\[\int_0^{\pi}f(x)dx=\int_0^{\pi}\sin xdx=\left[-\cos x\right]_0^{\pi}=2.\]
Therefore,
\[\lim_{n\to\infty}\frac{1}{n}\sum_{k=1}^n\sin\left(\frac{\pi k}{n}\right)=\frac{1}{\pi} \lim_{n\to\infty}R(f, P_n, A_n)=\frac{2}{\pi}.\]
\end{solution}
\vp
\noindent
{\bf \large Exercises  \thesection}
\setcounter{myquestion}{1}

 \begin{question}{\themyquestion}
Evaluate the following derivatives.
\begin{enumerate}[(a)]
\item
$\di\frac{d}{dx}\int_0^xe^{u^2}du$
\item
$\di\frac{d}{dx}\int_x^{1}\cos(u^2)du$
\item
$\di\frac{d}{dx}\int_x^{x^3}\sqrt{2+\sin u}\;du$
\end{enumerate}
\end{question}
\atc

 \begin{question}{\themyquestion}
Let $f:[-2, 6]\to\mathbb{R}$ be the function defined by
\[f(x)=\begin{cases} x^2-x,\quad & \text{if}\; -2\leq x<1,\\ \di x-\frac{1}{x},\quad &\text{if}\;\quad 1\leq x\leq 6.\end{cases} \]
Find a continuous function $F:[-2,6]\to\mathbb{R}$ such that $F$ is differentiable on $(-2,6)$, $F(0)=0$, and 
\[F'(x)=f(x)\hspace{1cm}\text{for all}\;x\in (-2,6).\]

\end{question}
\atc

 \begin{question}{\themyquestion}
Find the limit 
\[\lim_{n\to\infty}\frac{1^7+2^7+\cdots+n^7}{n^8}.\]
\end{question}
\atc
 \begin{question}{\themyquestion}
Find the limit 
\[\lim_{n\to\infty}\frac{1}{n}\sum_{k=1}^n\cos^2\left(\frac{2\pi k}{n}\right).\]
\end{question}
\vp

\section{Integration by Substitution and Integration by Parts}\label{sec4.5}

In this section, we prove the integration by substitution formula and integration by parts formula. We will only deal with the case where the function that we are integrating is continuous in the interior of the integration interval. For general case where the function is piecewise continuous, one can apply the additivity theorem.

\subsection{Integration by Substitution}
\begin{theorem}[label=230223_5]{Integration by Substitution}
Let  $g:[a,b]\to\mathbb{R}$ be a function that satisfies the following conditions:
\begin{enumerate}[(i)]
\item $g$ is continuous and one-to-one on $[a,b]$;
\item $g$  is continuously differentiable on $(a,b)$;
\item $g'(x)$ is bounded on $(a,b)$.
\end{enumerate}  Then $g  $ maps the interval $[a,b]$ onto  a closed and bounded interval $[c,d]$ with end points $g(a)$ and $g(b)$.   If $f:[c,d]\to\mathbb{R}$ is a   function that is   bounded and continuous on $(c,d)$, then the function $h:[a,b]\to \mathbb{R}$,
\[h(x)=f(g(x))g'(x)\]is  Riemann integrable and
\begin{equation}\label{eq230223_1}\int_a^bh(x)dx=\int_a^b f(g(x))g'(x)dx=\int_{g(a)}^{g(b)}f(u)du.\end{equation}
This is equivalent to
\begin{equation}\label{eq230223_2}\int_c^d f(u)du=\int_a^bf(g(x))|g'(x)|dx.\end{equation}
\end{theorem}
The function $g:[a,b]\to\mathbb{R}$ that satisfies all the three given conditions defines  a \emph{smooth} change of variables $u=g(x)$ from $x$ to $u$, in the sense that $g$ is continuously differentiable on $(a,b)$. 

\begin{myproof}{Proof}
Since $g$ is one-to-one, we have $g((a,b))\subset (c,d)$. Therefore, the function $h:[a,b]\to\mathbb{R}$, $h(x)=f(g(x))g'(x)$ is continuous and bounded on $(a,b)$, and hence, it is Riemann integrable.  For any $x\in [a,b]$, let
\begin{gather*}H_1(x) =\int_a^xh(u)du=\int_a^x f(g(u))g'(u)du,\\
F(x)=\int_c^xf(u)du,\\H_2(x) =\int_{g(a)}^{g(x)}f(u)du=F(g(x))-F(g(a)).\end{gather*}
Then $H_1(a)=H_2(a)=0$.
By fundamental theorem of calculus, $H_1$ and $H_2$ are differentiable on $(a,b)$, and for any $x\in (a,b)$,
\[H_1'(x)=h(x)=f(g(x))g'(x),\hspace{1cm} H_2'(x)=f(g(x))g'(x).\]
 Since $H_1'(x)=H_2'(x)$ for all $x\in (a,b)$, and $H_1(a)=H_2(a)$, we conclude that
$H_1(x)=H_2(x)$ for all $x\in [a,b]$. Namely,
\[\int_a^b f(g(x))g'(x)dx=\int_{g(a)}^{g(b)}f(u)du.\]
From this, we see that integration by substitution is just the inverse of the chain rule for differentiation.
To prove the equivalence of \eqref{eq230223_1} and \eqref{eq230223_2}, we consider two cases.

\textbf{Case I:}  $g$ is strictly increasing on $[a,b]$.\\In this case, $g'(x)\geq 0$, and $c=g(a)$, $d=g(b)$. So \eqref{eq230223_1} is equivalent to \eqref{eq230223_2}.

\textbf{Case II:} $g$ is strictly decreasing.\\In this case, $g'(x)\leq 0$, $g(a)=d$ and $g(b)=c$. Therefore,
\[\int_a^bf(g(x))|g'(x)|dx=-\int_a^bf(g(x))g'(x)dx\] and
\[\int_{g(a)}^{g(b)}f(u)du=\int_d^c f(u)du=-\int_c^df(u)du.\]Thus, \eqref{eq230223_1} and \eqref{eq230223_2} are equivalent.

\end{myproof}

If we impose the condition that $f$ is continuous at the boundary points $c$ and $d$, the condition that $g$ is one-to-one can be removed. The points $g(a)$ and $g(b)$ might not be the boundary points of the interval $J=g([a, b])$, but the proof still holds.
\begin{theorem}{General Integration by Substitution}
Let  $g:[a,b]\to\mathbb{R}$ be a function that satisfies the following conditions:
\begin{enumerate}[(i)]
\item $g$ is continuous  on $[a,b]$;
\item $g$  is continuously differentiable on $(a,b)$;
\item $g'(x)$ is bounded on $(a,b)$.
\end{enumerate}  Then $g$ maps $[a,b]$ to a closed and bounded  interval $J$.  If $f:J\to\mathbb{R}$ is a continuous  function, then the function $h:[a,b]\to \mathbb{R}$,
\[h(x)=f(g(x))g'(x)\]is  Riemann integrable and
\[\int_a^bh(x)dx=\int_a^b f(g(x))g'(x)dx=\int_{g(a)}^{g(b)}f(u)du.\]
 
\end{theorem}

\begin{example}{}
 Evaluate the integral $\di \int_{-2}^3x\sqrt{16+x^2}dx$. 
\end{example}\begin{solution}{Solution} Let $f(x)=\sqrt{x}$ and $g(x)=16+x^2$. The function $g$ is continuously differentiable, with $g'(x)=2x$, and it maps the interval $[-2, 3]$ onto the interval $[16, 25]$. However, it is not one-to-one. The function $f$ is continuous on $[16, 25]$, so we can  apply the integration by substitution. In practice, we will do substitution by letting $u=16+x^2$, and find that \[\frac{du}{dx}=2x.\]\bs
This implies that we can replace $xdx$ by $du/2$. When $x=-2$, $u=20$; when $x=3$, $u=25$. Thus,
\[\int_{-2}^3x\sqrt{16+x^2}dx=\frac{1}{2}\int_{20}^{25}\sqrt{u}du=\left[\frac{1}{3}u^{\frac{3}{2}}\right]_{20}^{25}=\frac{125-20\sqrt{20}}{3}.\]
Students are invited to split the integral into a sum of two integrals, one over the interval $[-2,0]$, and one over the interval $[0, 3]$. The function $g(x)$ is one-to-one on each of these two intervals. Check that the same answer is obtained.

\end{solution}
As we mentioned before, if the change of variables  is given by a one-to-one function $u=g(x)$, the function $f$ does not need to be continuous at the boundary points. Using addivitivity theorem, Theorem \ref{230223_5} still holds when $f$ is a bounded piecewise continuous function.
\begin{example}
{}Let $a$ be a positive number, and let $f:[0,a]\to\mathbb{R}$ be a piecewise continuous function that is bounded.
Show that 
\[\int_0^af(x)dx=\int_0^af(a-x)dx.\]
\end{example}
\begin{solution}{Solution}
We consider the change of variables $u=g(x)=a-x$. This is a strictly monotonic function with $g'(x)=-1$. Therefore, $du=-dx$. When $x=0$, $u=a$; when $x=a$, $u=0$. Hence,
\[\int_0^af(x)dx= \int_a^0f(a-u)(-du)=\int_0^af(a-x)dx.\]
\end{solution}
\begin{example}
{}Let $a$ be a positive number, and let $f:[-a,a]\to\mathbb{R}$ be a piecewise continuous function that is bounded.
\begin{enumerate}[(a)]
\item If $f$ is an even function, show that 
\[\int_{-a}^a f(x)dx=2\int_0^af(x)dx.\]
\item If $f$ is an odd function, show that 
\[\int_{-a}^a f(x)dx=0.\]

\end{enumerate}
\end{example}
\begin{solution}{Solution}Notice that
\[\int_{-a}^af(x)dx=\int_{-a}^0f(x)dx+\int_0^af(x)dx.\]For the   integral 
$\di \int_{-a}^0f(x)dx$, 
we consider the change of variables $u=g(x)= -x$. This is a strictly monotonic function with $g'(x)=-1$. Therefore, $du=-dx$. When $x=-a$, $u=a$; when $x=0$, $u=0$. Hence,
\[\int^0_{-a}f(x)dx= \int_a^0f(-u)(-du)=\int_0^af(-x)dx.\] 
\begin{enumerate}[(a)]
\item
When $f$ is an even function, $f(-x)=f(x)$ for all $x\in [0,a]$. Therefore,
\[\int_{-a}^af(x)dx= \int_0^af(x)dx+\int_0^af(x)dx=2\int_0^af(x)dx.\]
\item[(b)] When  $f$ is an odd function, $f(-x)=-f(x)$ for all $x\in [0,a]$. Therefore,
\[\int_{-a}^af(x)dx= -\int_0^af(x)dx+\int_0^af(x)dx=0.\]
\end{enumerate}
\end{solution}

\begin{figure}[ht]
\centering
\includegraphics[scale=0.2]{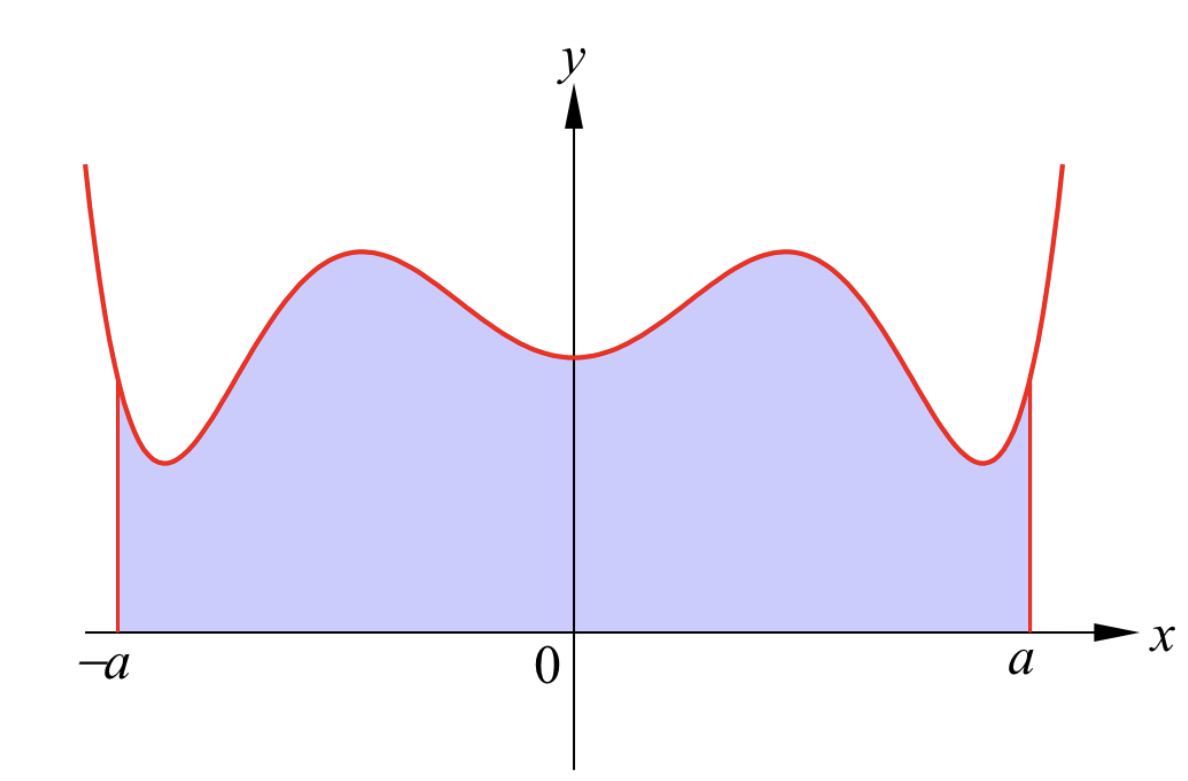}
\caption{An even function.\fa}\label{figure46}
\end{figure}

\begin{figure}[ht]
\centering
\includegraphics[scale=0.2]{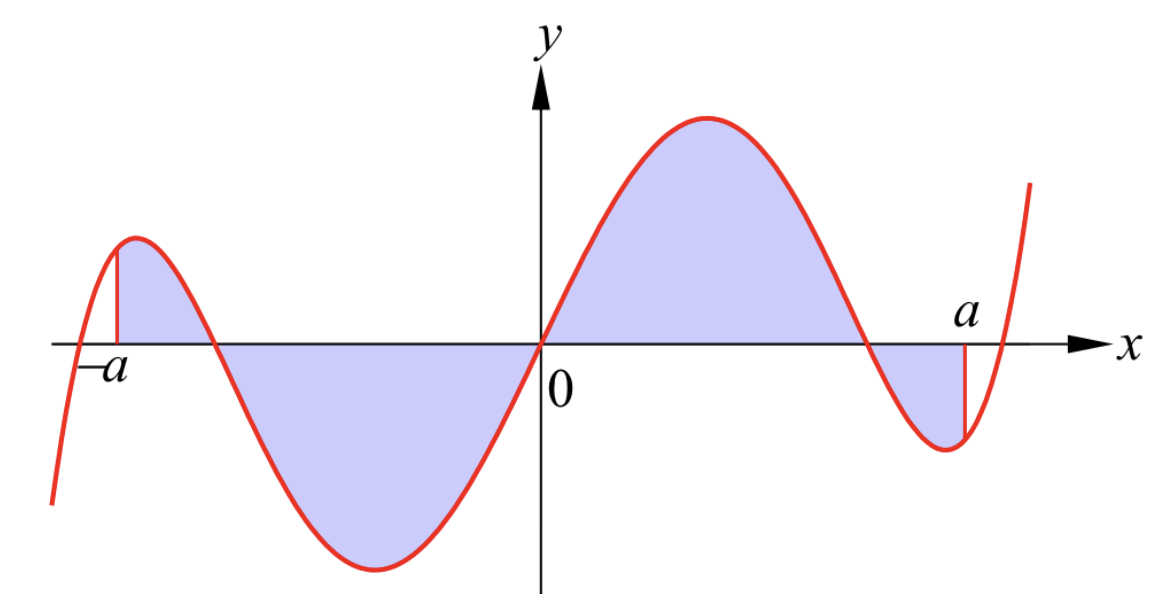}
\caption{An odd function.\fa}\label{figure47}
\end{figure}
\begin{example}{Area of a Circle}
 Find the area of a circle of radius $r$.

\begin{center}
\includegraphics[scale=0.2]{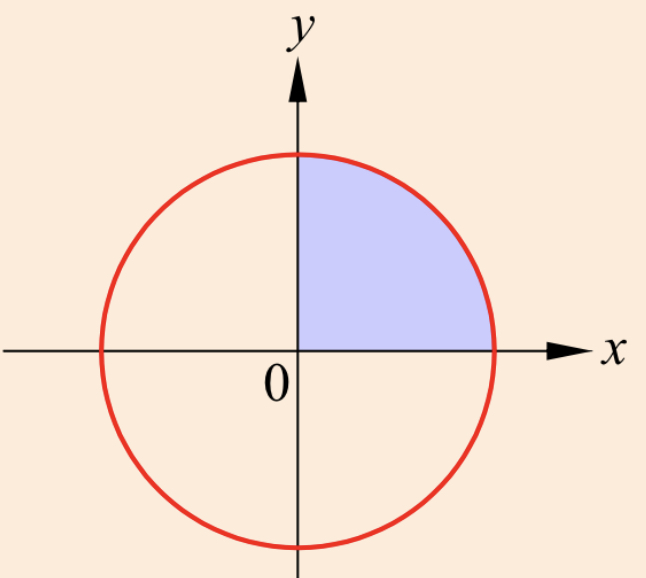}\end{center}

\end{example}
 
\begin{solution}{Solution}
A circle of radius $r$ with center at the origin has equation $x^2+y^2=r^2$. By symmetry, it is enough for us to find the area in the first quadrant, and then multiply by 4. The sector in the first quadrant is bounded by the curve $y=\sqrt{r^2-x^2}$, the lines $x=0$, $x=r$, and the $x$-axis.  Hence, the area of a circle of radius $r$ is
\[A=4\int_0^{r}\sqrt{r^2-x^2}dx.\]
Making a change of variables $x=r\sin\theta$, we find that 
\[\frac{dx}{d\theta}=r\cos\theta.\]When $x=0$, $\theta=0$; when $x=r$, $\theta=\pi/2$. Therefore,
\begin{align*}
A&=4\int_0^{\frac{\pi}{2}}\sqrt{r^2-r^2\sin^2\theta}\;r\cos\theta \,d\theta\\&=4r^2\int_0^{\frac{\pi}{2}}\cos^2\theta\,d\theta.\end{align*}
 Using the formula
\[\cos^2\theta=\frac{1+\cos 2\theta}{2},\]
we have
\begin{align*}A &=2r^2\int_0^{\frac{\pi}{2}}\left(1+\cos2\theta\right)d\theta \\&=2r^2\left[\theta+\frac{\sin2\theta}{2}\right]_0^{\frac{\pi}{2}} \\&=2r^2\times\frac{\pi}{2}\\&=\pi r^2.\end{align*}
\end{solution}
\subsection{Integration by Parts}

\begin{theorem}{Integration by Parts}
Let $f:[a,b]\to\mathbb{R}$ and $g:[a,b]\to\mathbb{R}$ be  functions that satisfy the following conditions:
\begin{enumerate}[(i)]
\item $f$ and $g$ are continuous on $[a,b]$;
\item $f$ and $g$ are  continuously differentiable on $(a,b)$;
\item $f'(x)$ and $g'(x)$ are bounded on $(a,b)$.
\end{enumerate}Then $fg'$ and $gf'$ are Riemann integrable on $[a,b]$, and
\[\int_a^b f(x)g'(x)dx=f(b)g(b)-f(a)g(a)-\int_a^bg(x) f'(x)dx.\]
\end{theorem}
\begin{myproof}{Proof}
Since  $f$ and $g$ are continuous on $[a,b]$, they are bounded. Since $f'(x)$ and $g'(x)$ are conitnuous and bounded on $(a,b)$, $fg'$ and $f'g$ are continuous and bounded on $(a,b)$. Therefore, $fg'$ and $gf'$ are Riemann integrable on $[a,b]$. By product rule, for any $x\in (a,b)$, 
\[(fg)'(x)=f(x)g'(x)+g(x)f'(x).\]
So $(fg)'$ is also bounded and continuous on $(a,b)$, and hence Riemann integrable on $[a,b]$. Since $fg$ is also continuous on $[a,b]$, we can apply  fundamental theorem of calculus, which gives
\[\int_a^b (fg)'(x)dx=(fg)(b)-(fg)(a).\]
Therefore,
\[\int_a^b f(x)g'(x)dx+\int_a^b g(x)f'(x)dx= f(b)g(b)-f(a)g(a).\]This proves the integration by parts formula.
\end{myproof}In a nutshell, the integration by parts formula is just the inverse of the product rule of differentiation. But it is a very useful integration technique.
\begin{highlight}{Integration by Parts}
The integration by parts formula is often expressed as
\[\int udv=uv-\int vdu.\]
In practice, we  identify which part should be $u$ and which part should be $dv$. The function $v$ is defined up to a constant.  One can verify directly that if $v$ is replaced by $v+C$, where $C$ is a constant, the right hand side of the formula is not changed. Hence, we can choose a $v$ that is most convenient.
\end{highlight}

\begin{example}{}
Let $n$ be a positive integer. Evaluate the integral
\[\int_1^e\frac{\ln x}{x^n}dx.\]

\end{example}
\begin{solution}{Solution}If $n=1$, we use integration by substitution with $u=\ln x$. Then
\[\frac{du}{dx}=\frac{1}{x}.\]
When $x=1$, $u=0$; when $x=e$, $u=1$. Therefore,
\[\int_1^e\frac{\ln x}{x}dx=\int_0^1udu=\left[\frac{u^2}{2}\right]_{0}^1=\frac{1}{2}.\]
If $n\geq 2$, we use integration by parts. Let 
\[u(x)=\ln x,\hspace{1cm}v'(x)=\frac{1}{x^{n}}.\]
Then
\[\frac{du}{dx}=\frac{1}{x},\quad v(x)=-\frac{1}{n-1}\times\frac{1}{x^{n-1}}.\]Both of $u(x)$ and $v(x)$ are continuously differentiable functions on $(0,\infty)$. 
\bs Therefore,
\begin{align*}\int_1^e\frac{\ln x}{x^n}dx&=\left[-\frac{1}{n-1}\times\frac{\ln x}{x^{n-1}}\right]_1^e+\frac{1}{n-1}\int_1^e\frac{1}{x^{n}}dx\\
&=-\frac{1}{n-1}\times\frac{1}{e^{n-1}}-\frac{1}{(n-1)^2}\left[\frac{1}{x^{n-1}}\right]_1^e\\
&=\frac{1}{(n-1)^2}-\frac{n}{(n-1)^2}\frac{1}{e^{n-1}}.
\end{align*}
\end{solution}
\begin{example}[label=230307_10]{}Let $I$ be an open interval that contains the point $x_0$, and 
let $f:I\to\mathbb{R}$ be a continuous function. Given  a positive integer $n$, define the function $F:I\to\mathbb{R}$ by
\[F(x)=\frac{1}{n!}\int_{x_0}^x(x-t)^nf(t)dt.\]
Prove that $F$ is $(n+1)$ times continuously differentiable, 
\[F(x_0)=F'(x_0)=\ldots=F^{(n)}(x_0)=0,\]
and
\[F^{(n+1)}(x)=f(x)\hspace{1cm}\text{for all}\;x\in I.\]
\end{example}
\begin{solution}
{Solution}
Define the function $g:\mathbb{R}\to\mathbb{R}$ by
\[g(x)=\int_{x_0}^xf(t)dt.\] Then $g(x_0)=0$, and 
by fundamental theorem of calculus,
\[g'(x)=f(x)\hspace{1cm}\text{for all}\;x\in I.\]

Now we prove the statement by induction on $n$. When $n=1$, 
\[F(x)=\int_{x_0}^x(x-t)f(t)dt.\]\bs
By definition, $F(x_0)=0$. 
For a fixed $x$, using integration by parts with $u(t)=x-t$ and $v'(t)=f(t)$, we find that  
\[\frac{du}{dt}=-1, \quad v(t)=g(t).\]  
It follows that
\[F(x)=\Bigl[(x-t)g(t)\Bigr]_{t=x_0}^{t=x}+\int_{x_0}^x g(t)dt=\int_{x_0}^x g(t)dt.\] Notice that $g(t)$ is continuously differentiable, and hence it is continuous. By fundamental theorem of calculus, 
\[F'(x)=g(x)\hspace{1cm}\text{for all}\;x\in I.\]
Therefore, $F'(x_0)=g(x_0)=0$, and 
\[F''(x)=g'(x)=f(x)\hspace{1cm}\text{for all}\;x\in I.\]
This proves that $F(x)$ is twice continuously differentiable. Since we have also shown that $F(x_0)=F'(x_0)=0$, and $F''(x)=f(x)$ for all $x\in I$.
 the statement is true when $n=1$.

Assume that we have proved the statement when $n=k-1$, where $k\geq 2$. When $n=k$, 
\[F(x)=\frac{1}{k!}\int_{x_0}^x(x-t)^kf(t)dt.\]
  For a fixed $x$, using integration by parts with $u(t)=(x-t)^k$ and $v'(t)=f(t)$, we find that 
\[\frac{du}{dt}=-k(x-t)^{k-1}, \quad v(t)=g(t).\] 
It follows that
\begin{align*}F(x)&=\frac{1}{k!}\left[(x-t)^kg(t)\right]_{t=x_0}^{t=x}+\frac{1}{(k-1)!}\int_{x_0}^x (x-t)^{k-1}g(t)dt\\&=\frac{1}{(k-1)!}\int_{x_0}^x (x-t)^{k-1}g(t)dt.\end{align*}
By inductive hypothesis, the function
$F(x)$ satisfies 
\[F(x_0)=F'(x_0)=\cdots=F^{(k-1)}(x_0)=0,\]\bs
and
\[F^{(k)}(x)=g(x)\hspace{1cm}\text{for all}\;x\in I.\] 
The latter implies that $F^{(k)}(x_0)=g(x_0)=0$, and $F(x)$ is $(k+1)$ times differentiable, with
\[F^{(k+1)}(x)=g'(x)=f(x) \] a continuous function. Therefore, when $n=k+1$, the statement also holds.

By principle of mathematical induction, the statement is true for all positive integers $n$.
\end{solution}

\vp
\noindent
{\bf \large Exercises  \thesection}
\setcounter{myquestion}{1}
\atc
\begin{question}{\themyquestion}Let $f:[a,b]\to\mathbb{R}$ be a bounded function that is Riemann integrable.
Show that for any real number $c$,
\[\int_a^bf(x)dx=\int_{a+c}^{b+c}f(x-c)dx.\]
\end{question}
\begin{question}{\themyquestion}
Explain why  
\[\int_{-1}^1 (x+1)e^{-x^2}dx=2\int_0^1e^{-x^2}dx.\]
\end{question}
 
\atc
\begin{question}{\themyquestion}
 Let $a$ be a positive number. Assume that the functions $f:[0,a]\to\mathbb{R}$ and $g:[0,a]\to\mathbb{R}$ are bounded and piecewise continuous, prove that
\[\int_0^af(x)g(a-x)dx=\int_0^af(a-x)g(x)dx.\]
\end{question}
 \atc
\begin{question}[label=ex230225_1]{\themyquestion}
 Let $m$ and $n$ be nonnegative integers. Show that
\[\int_0^1x^m(1-x)^ndx=\frac{m!\,n!}{(m+n+1)!}.\]
\end{question}

 \atc
\begin{question}{\themyquestion}
 Let $f:[a,b]\to\mathbb{R}$ be a continuous and strictly increasing function which maps the interval $[a,b]$ bijectively onto the interval $[c,d]$, where $c=f(a)$, and $d=f(b)$. Denote by $g:[c,d]\to\mathbb{R}$  the inverse function of $f$. Notice that  $f:[a,b]\to\mathbb{R}$ and $g:[c,d]\to\mathbb{R}$ are Riemann integrable. This question is regarding the proof of the formula
\begin{equation}\label{eq230224_1}\int_a^bf(x)dx=bf(b)-af(a)-\int_c^dg(x)dx.\end{equation}
 \begin{enumerate}[(a)]
\item If $a>0$ and $c>0$,  draw a figure to illustrate  the formula.
\item If $f$ is continuously differentiable on $(a,b)$,  
use integration by substitution with $u=g(x)$ to prove the formula \eqref{eq230224_1}.

 \item  Let $P_n=\{x_i\}_{i=0}^n$ be the regular partition of the interval $[a,b]$ into $n$ intervals. For $0\leq i\leq n$, let $y_i=f(x_i)$. Then $\widetilde{P}_n=\{y_i\}_{i=0}^n$ is a partition of $[c,d]$. 
\begin{enumerate}[(i)]
\item
Show that \[\sum_{i=1}^n f(x_{i-1})(x_i-x_{i-1})+\sum_{i=1}^n g(y_i)(y_i-y_{i-1})=bf(b)-af(a).\]
\item Show that $\di\lim_{n\to\infty}|P_n|=0$ and  $\di\lim_{n\to\infty}|\widetilde{P}_n|=0$. You might want to use uniform continuity.
\item Use part (i) and part (ii) to prove the formula
 \eqref{eq230224_1}.
\end{enumerate}
\end{enumerate}
\end{question}
\vp

\section{Improper Integrals}\label{sec4.6}
In this section, we want to discuss Riemann integrals for functions $f:I\to\mathbb{R}$ defined on an interval $I$, where either $I$ is not bounded, or $f$ is not bounded on $I$, or both. 

\begin{definition}{Improper Integral}
Let $I$ be an interval and let  $f:I\to\mathbb{R}$ be a function defined on $I$.  An integral of the form
\[\int_If\] is an improper integral if either $f$ is not bounded on $I$, or $I$ is an unbounded interval.
\end{definition}
This is not a rigorous definition. We will only be interested in the case where we can make sense of $\di\int_If$.

As an example,  Theorem  \ref{thm230222_13} says that there exists a differentiable function $g:(-1,1)\to\mathbb{R}$ satifying
\[g'(x)=\frac{1}{\sqrt{1-x^2}},\hspace{1cm}g(0)=0.\]It is given by
\[g(x)=\int_0^x\frac{du}{\sqrt{1-u^2}}du.\]
Notice that the function $\di f(u)=\di\frac{1}{\sqrt{1-u^2}}$ is bounded and continuous on the interval $[0,x]$ if $0< x<1$, and on $[x,0]$ if $-1<x<0$. Therefore, $g(x)$ is a well-defined Riemann integral when $-1<x<1$. We are interested to extend the definition of $g(x)$ to $x= 1$ and $x=-1$. But $f$ is not bounded on $(-1,1)$, so we cannot define the Riemann integral of $f$ on $[0, 1]$ or $[-1,0]$.  Our studies on the function $\sin x$ shows that $g(x)=\sin^{-1}x$ when $x\in(-1,1)$. Thus,
\[\lim_{x\to 1^-}g(x)=\sin^{-1}1=\frac{\pi}{2},\hspace{1cm}\lim_{x\to -1^+}g(x)=\sin^{-1}(-1)=-\frac{\pi}{2}.\]
Hence, it is reasonable to say that the improper integrals 
\[\int_0^{1}\frac{1}{\sqrt{1-u^2}}du\quad\text{and}\quad \int_0^{-1}\frac{1}{\sqrt{1-u^2}}du\] have values
\[\lim_{x\to 1^-}\int_0^{x}\frac{1}{\sqrt{1-u^2}}du=\frac{\pi}{2}\quad\text{and}\quad\lim_{x\to -1^+}\int_0^{x}\frac{1}{\sqrt{1-u^2}}du=-\frac{\pi}{2}\]respectively. This is how we are going to make sense of improper integrals.  
\begin{definition}
{Improper Integrals of Unbounded Functions}
\begin{enumerate}[1.]
\item
If the function $f:(a, b]\to\mathbb{R}$ is not bounded, but it is bounded and Riemann integrable on any interval $[c,b]$ with $a<c<b$, then we say that the improper integral $\di\int_a^bf(x)dx$ is convergent if  the limit
\[\lim_{c\to a^+}\int_c^bf(x)dx\]   exists. Otherwise, we say that the improper integral   is divergent. When the improper integral is convergent, we define its value as
\[ \int_a^bf(x)dx=\lim_{c\to a^+}\int_c^bf(x)dx.\] 
\item If the function $f:[a, b)\to\mathbb{R}$ is not bounded, but it is bounded and Riemann integrable on any interval $[a,c]$ with $a<c<b$, then we say that the improper integral $\di\int_a^bf(x)dx$ is convergent if  the limit
\[\lim_{c\to b^-}\int_a^cf(x)dx\]   exists. Otherwise, we say that the improper integral   is divergent. When the improper integral is convergent, we define its value as
\[ \int_a^bf(x)dx=\lim_{c\to b^-}\int_a^cf(x)dx.\] 

\end{enumerate}
\end{definition}

\begin{highlight}{Improper Integrals of Unbounded Functions}
Putting in another way, if the function $f:(a, b]\to\mathbb{R}$ is not bounded, but  it is bounded and Riemann integrable on any  intervals $[x,b]$  when $a<x<b$, we define the function $F:(a,b]\to\mathbb{R}$ by
\[F(x)=\int_x^{b}f(u)du.\]
Then $F$ is a continuous function. We say that the improper integral $\di\int_a^{b}f(x)dx$ is convergent if and only if 
\[\lim_{x\to a^+}F(x)\] exists.
Similarly, for a function $f:[a,b)\to\mathbb{R}$ that is not bounded, but is bounded and Riemann integral on any  intervals $[a,x]$ when $a<x<b$,  we say that the improper integral $\di\int_{a}^bf(x)dx$ is convergent if and only if the continuous function $F:[a,b)\to\mathbb{R}$ defined by
\[F(x)=\int_{a}^x f(u)du\] has a limit when $x\to b^-$.
\end{highlight}
\begin{example}{}
$\di\int_0^1\frac{1}{\sqrt{1-x^2}}dx$ is an improper integral as the function $f:[0,1)\to\mathbb{R}$, $f(x)=\di \frac{1}{\sqrt{1-x^2}}$ is not bounded. We have seen that this improper integral is convergent and has value $\di \frac{\pi}{2}$.
\end{example}

\begin{example}[label=ex230227_10]{}
Let $p$ be a positive number. Determine those values of $p$ for which the improper integral $\di\int_0^1\frac{1}{x^p}dx$ is convergent. Find the value of the improper integral when it is convergent.
\end{example}
\begin{solution}
{Solution}
For $p>0$, define the function $F:(0, 1]\to\mathbb{R}$ by 
\[F(x)=\int_x^1\frac{1}{u^p}du.\]Then
\begin{align*}
F(x)=\begin{cases}-\ln x,\quad &\text{if}\; p=1,\\
\di \frac{1-x^{1-p}}{1-p},\quad &\text{if}\; p\neq 1.\end{cases}
\end{align*} From this, we see that
$\di \lim_{x\to 0^+}F(x)$ exists if and only if $0<p<1$. Hence, the improper integral  $\di\int_0^1\frac{1}{x^p}dx$ is convergent if and only if $0<p<1$. In this case, 
\[\int_0^1\frac{1}{x^p}dx=\frac{1}{1-p},\quad 0<p<1.\]
\end{solution}
When $r\geq 0$, the integral $\di\int_0^1 x^rdx$ is just an ordinary integral. However, we will sometimes abuse terminology and say that the integral $\di\int_0^1x^rdx$ is convergent if and only if $r>-1$. 
\begin{highlight}{}
If $c$ is a point in $(a,b)$ and we have a function $f:[a,b]\setminus\{c\}\to\mathbb{R}$  that is not bounded, we will define the improper integral $\di\int_a^bf(x)dx$ as 
\[\int_a^cf(x)dx+\int_c^bf(x)dx.\]
We say that  the improper integral $\di\int_a^bf(x)dx$ is convergent provided that both improper integrals $\di \int_a^cf(x)dx$ and $\di \int_c^bf(x)dx$ are convergent. 
\end{highlight}
Next we consider improper integrals defined on unbounded intervals.
\begin{definition}{Improper integrals on Unbounded Intervals}
\begin{enumerate}[1.]
\item
If $f:[a,\infty)\to\mathbb{R}$ is a  function that is bounded and Riemann integrable on any bounded intervals $[a, b]$, we say that the improper integral $\di\int_a^{\infty}f(x)dx$ is convergent if the limit
\[\lim_{b\to\infty}\int_a^{b}f(x)dx\] exists. Otherwise, we say that the improper integral is divergent.
If the improper integral is convergent, we define its value as  
\[\int_a^{\infty}f(x)dx=\lim_{b\to\infty}\int_a^{b}f(x)dx.\]
\item
If $f:(-\infty,b]\to\mathbb{R}$ is a   function that is bounded and Riemann integrable on any bounded intervals $[a, b]$, we say that the improper integral $\di\int^b_{-\infty}f(x)dx$ is convergent if the limit
\[\lim_{a\to-\infty}\int_a^{b}f(x)dx\] exists. Otherwise, we say that the improper integral is divergent.
If the improper integral is convergent, we define its value as   
\[\int^b_{-\infty}f(x)dx=\lim_{a\to-\infty}\int_a^{b}f(x)dx.\]
\item If $f:\mathbb{R}\to\mathbb{R}$ is a function that is bounded and Riemann integrable on any bounded intervals $[a,b]$, we say that the improper integral $\di\int_{-\infty}^{\infty}f(x)dx$ is convergent if and only if for any  real number $c$, both the improper integrals 
\[\int_{-\infty}^{c}f(x)dx\hspace{1cm}\text{and}\hspace{1cm}\int_{c}^{\infty}f(x)dx\] are convergent. In such a case, we define the improper integral as
\begin{equation}\label{eq230224_2}\int_{-\infty}^{\infty}f(x)dx=\int_{-\infty}^cf(x)dx+\int_c^{\infty}f(x)dx.\end{equation}
\end{enumerate}
\end{definition}
\begin{remark}{}
 To make the integral $\di\int_{-\infty}^{\infty}f(x)dx$ well defined when it is convergent, we need to check that the right hand side of \eqref{eq230224_2} does not depend on the point $c$. 
In fact, we can show that if there is a real number $c_0$ so that both the improper integrals 
\[\int_{-\infty}^{c_0}f(x)dx\hspace{1cm}\text{and}\hspace{1cm}\int_{c_0}^{\infty}f(x)dx\] are convergent,  then for any other values of $c$, \[\int_{-\infty}^{c}f(x)dx\hspace{1cm}\text{and}\hspace{1cm}\int_{c}^{\infty}f(x)dx\] are convergent. This is just due to additivity, which says that
\begin{align*}\int_a^cf(x)dx&=\int_a^{c_0}f(x)dx+\int_{c_0}^cf(x)dx,\\
\int_c^bf(x)dx&=\int_{c}^{c_0}f(x)dx+\int_{c_0}^bf(x)dx.\end{align*}
Thus, $\di\lim_{a\to-\infty}\int_a^cf(x)dx$ exists if and only if $\di\lim_{a\to-\infty}\int_a^{c_0}f(x)dx$ exists, and 
$\di\lim_{b\to\infty}\int_c^bf(x)dx$ exists if and only if 
$\di\lim_{b\to \infty}\int_{c_0}^{b}f(x)dx$ exists. Moreover,
 \begin{align*}\lim_{a\to-\infty}\int_a^cf(x)dx&=\lim_{a\to-\infty}\int_a^{c_0}f(x)dx+\int_{c_0}^cf(x)dx,\\
\lim_{b\to\infty}\int_c^bf(x)dx&=\int_{c}^{c_0}f(x)dx+\lim_{b\to\infty}\int_{c_0}^bf(x)dx.\end{align*}
Since $\di \int_{c_0}^cf(x)dx=-\int_{c}^{c_0}f(x)dx$, we find that
\[\int_{-\infty}^{c_0}f(x)dx+\int_{c_0}^{\infty}f(x)dx=\int_{-\infty}^{c}f(x)dx+\int_{c}^{\infty}f(x)dx.\]
\end{remark}
\begin{highlight}{Improper Integrals on Unbounded Intervals}
Putting in another way, if $f:[a,\infty)\to\mathbb{R}$ is a function that is bounded and Riemann integrable on any bounded intervals, we define the function $F:[a,\infty)\to\mathbb{R}$ by
\[F(x)=\int_a^{x}f(u)du.\]
Then $F$ is a continuous function. We say that the improper integral $\di\int_a^{\infty}f(x)dx$ is convergent if and only if the limit
\[\lim_{x\to\infty}F(x)\] exists.
Similarly, for a function $f:(-\infty, b]\to\mathbb{R}$ that is bounded and Riemann integrable on any bounded intervals $[a,b]$, we say that the improper integral $\di\int_{-\infty}^bf(x)dx$ is convergent if and only if the continuous function $F:(-\infty, b]\to\mathbb{R}$ defined by
\[F(x)=\int_{x}^b f(u)du\] has a limit when $x\to-\infty$.
\end{highlight}

\begin{example}{}
Let $p$ be any real number. Determine those values of $p$ for which the improper integral $\di\int_1^{\infty}\frac{1}{x^p}dx$ is convergent. Find the value of the improper integral when it is convergent.
\end{example}
\begin{solution}{Solution}
For a fixed real number $p$, define the function $F:[1,\infty)\to\mathbb{R}$ by
\[F(x)=\int_1^x\frac{1}{u^p}dx.\]\bs
Then
\begin{align*}
F(x)=\begin{cases} \ln x,\quad &\text{if}\; p=1,\\
\di\frac{x^{1-p}-1}{1-p},\quad &\text{if}\;p\neq 1.\end{cases}
\end{align*}From this, we see that the limit $\di\lim_{x\to\infty}F(x)$ exists if and only if $p>1$. Hence, the improper integral  $\di\int_1^{\infty}\frac{1}{x^p}dx$ is convergent if and only if $p>1$, and
\[  \int_1^{\infty}\frac{1}{x^p}dx=\frac{1}{p-1},\quad p>1.\]
\end{solution}

\begin{example}{}
Determine whether the improper integral  is convergent. If yes, find the value of the integral.
\begin{enumerate}[(a)]
\item 
$\di \int_0^{\infty}\frac{1}{1+x^2}dx$
\item $\di \int_{-\infty}^0e^{x}dx$
\item $\di\int_{-\infty}^{\infty}\frac{x}{x^2+1}dx$
\end{enumerate}
\end{example}
\begin{solution}{Solution}
\begin{enumerate}[(a)]
\item
Since $\di\frac{d}{dx}\tan^{-1}x=\frac{1}{1+x^2}$, we find that
\[\int_0^b \frac{1}{1+x^2}dx=\tan^{-1}b-\tan^{-1}0=\tan^{-1}b.\]
Since
\[\lim_{b\to\infty}\tan^{-1}b=\frac{\pi}{2},\]the improper integral $\di \int_0^{\infty}\frac{1}{1+x^2}dx$ is convergent and its value is \end{enumerate}\bs\begin{enumerate}[(a)]\item[]
\[ \int_0^{\infty}\frac{1}{1+x^2}dx=\lim_{b\to\infty} \int_0^b \frac{1}{1+x^2}dx =\lim_{b\to\infty} \tan^{-1}b=\frac{\pi}{2}.\]
\item[(b)] Since $e^a\to 0$ as $a\to-\infty$, we have
\[ \int_{-\infty}^0e^{x}dx=\lim_{a\to -\infty}\int_a^{0}e^xdx=\lim_{a\to-\infty}\left(1-e^a\right)=1.\]The improper integral  $\di \int_{-\infty}^0e^{x}dx$ is convergent and is equal to 1.
\item[(c)] Here, we consider the improper integrals
\[\int_{-\infty}^0\frac{x}{x^2+1}dx\hspace{1cm}\text{and}\hspace{1cm}\int_0^{\infty} \frac{x}{x^2+1}dx.\]
Since
\[\frac{d}{dx}\ln(1+x^2)=\frac{2x}{1+x^2},\]we find that
\[\int_0^b\frac{x}{1+x^2}dx=\frac{1}{2}\ln(1+b^2).\]
But \[\lim_{b\to\infty}\ln(1+b^2)=\infty.\]
Hence,
 the improper integral $\di\int^{\infty}_0\frac{x}{x^2+1}dx$ is divergent. So, the improper integral  $\di\int_{-\infty}^{\infty}\frac{x}{x^2+1}dx$ is also divergent.
\end{enumerate}
\end{solution}
One is tempted to  
define the improper integral  
$\di\int_{-\infty}^{\infty}f(x)dx$  as
\[\lim_{a\to \infty}\int_{-a}^af(x)dx\] if it exists. For part (c) in the example above, $\di f(x)=\frac{x}{1+x^2}$ is an odd function. Thus, $\di \int_{-a}^a\frac{x}{1+x^2}dx=0$ for any $a$, and so
\[\lim_{a\to \infty}\int_{-a}^a\frac{x}{1+x^2}dx=0.\]
In fact, if $f:\mathbb{R}\to\mathbb{R}$ is an odd function, then we always have
\[\lim_{a\to \infty}\int_{-a}^af(x)dx=0.\]If we use   the limit \[\lim_{a\to \infty}\int_{-a}^af(x)dx\] as a definition for the improper integral $\di\int_{-\infty}^{\infty}f(x)dx$, it will lead to  undesirable results, such as that the integral $\di \int_{-\infty}^{\infty} xdx$ is convergent. Nevertheless,  the limit \[\lim_{a\to \infty}\int_{-a}^af(x)dx,\] if it exists, has some applications. It is called the Cauchy principal value of $\di\int_{-\infty}^{\infty}f(x)dx$.

\begin{definition}{Cauchy Principal Value}
If $f:\mathbb{R}\to\mathbb{R}$ is a function that is bounded and Riemann integrable on any symmetric bounded intervals $[-a,a]$, the Cauchy principal value of the improper integral $\di\int_{-\infty}^{\infty}f(x)dx$, denoted by $\text{P.V.}\, \di\int_{-\infty}^{\infty}f(x)dx$, is defined as
\[\text{P.V.}\,\int_{-\infty}^{\infty}f(x)dx=\lim_{a\to \infty}\int_{-a}^af(x)dx,\]
if the limit exists.
\end{definition}  
Thus, we find that if $f:\mathbb{R}\to\mathbb{R}$ is an odd function, then $\di\text{P.V.}\,\int_{-\infty}^{\infty}f(x)dx=0$. It is also easy to prove the following.

\begin{proposition}[label=230224_10]{} If the improper integral  $\di\int_{-\infty}^{\infty}f(x)dx$ is convergent, then 
its Cauchy principal value exists, and is equal to the improper integral. Namely, 
\[\text{P.V.}\,  \int_{-\infty}^{\infty}f(x)dx=  \int_{-\infty}^{\infty}f(x)dx.\]
\end{proposition}
\begin{myproof}{Proof}
If  the improper integral $\di \int_{-\infty}^{\infty}f(x)dx$ is convergent, then  the limits
\[\lim_{  c\to-\infty } \int_c^0f(x)dx\quad\text{and}\quad\lim_{b\to\infty}\int_0^bf(x)dx\]exists and
\[ \int_{-\infty}^{\infty}f(x)dx=\lim_{  c\to-\infty } \int_c^0f(x)dx+\lim_{b\to\infty}\int_0^bf(x)dx. \] This implies that 
\[ \lim_{a\to \infty} \int_{-a}^af(x)dx=
\lim_{a\to\infty}\int_{-a}^0f(x)dx+\lim_{a\to\infty}\int_0^af(x)dx=\int_{-\infty}^{\infty}f(x)dx.\] 
\end{myproof}

Consider the integral
\begin{equation}\label{eq230224_4}\int_0^{\infty}\frac{1}{\sqrt{x}(x+1)}dx.\end{equation}
The function $f:(0,\infty)\to\mathbb{R}$,
\[f(x)=\frac{1}{\sqrt{x}(x+1)}\] is not bounded on any interval $(0, b]$ when $b>0$. Hence, the integral is an improper integral of an unbounded function defined on an unbounded interval. Using the same principle, we will say that it is convergent if and only if for any $c>0$,  the improper integrals 
\[\int_0^cf(x)dx\hspace{1cm}\text{and}\hspace{1cm}\int_c^{\infty}f(x)dx\] are convergent. 

Another natural question to ask is whether one can determine whether an improper integral is convergent without explicitly computing the integral. There are some partial solutions to this. 

\begin{highlight}
{}If $J$ is an interval that is contained in the interval $I$, and the integral $\di\int_J f(x)dx$ is divergent, then the integral $\di\int_I f(x)dx$ is divergent. 
\end{highlight}
For instance, the integral $\di\int_0^{\infty}f(x)dx$ is divergent  if  the integral $\di\int_1^{\infty}f(x)dx$ is divergent.

The next proposition says that linear combination  of convergent integrals must be convergent.
\begin{proposition}{Linearity}
Let $I$ be an interval. If the improper integrals $\di\int_If(x)dx$ and $\di\int_Ig(x)dx$ are convergent, then for any constants $\alpha$ and $\beta$, the improper integral $\di\int_I(\alpha f+\beta g)$ is also convergent, and
\[\int_I(\alpha f+\beta g)=\alpha\int_If+\beta\int_I g.\]
\end{proposition}This follows easily from limit laws.
Now we want to prove some comparison theorems for improper integrals. We start with integrals of nonnegative functions. If a function $f $ is nonpositive, one just consider the function $-f$, which is then nonnegative. 
 
\begin{lemma}[label=230224_5]{} Let $I$ be an interval. Given that $f:I\to\mathbb{R}$  is a nonnegative function that is bounded and Riemann integrable on any closed and bounded intervals that are contained in $I$. Fixed $x_0$ in $I$ and define the function $F:I\to\mathbb{R}$ by
\[F(x)=\int_{x_0}^xf(u)du.\]
\begin{enumerate}[1.]
\item If $I=(a, b]$ or $I=(-\infty,b]$, then the integral $\di\int_I f(x)dx$ is convergent if and only if the function $F(x)$ is bounded below. 
\item If $I=[a,b)$ or $I=[a,\infty)$, then the integral $\di\int_If(x)dx$ is convergent if and only if the function $F(x)$ is bounded above. 
 \end{enumerate}
\end{lemma}

\begin{myproof}{Proof}
Notice that since $f(u)\geq 0$ for all $u\in I$, for any $x_1$ and $x_2$ in $I$, if $x_1<x_2$, then 
\[F(x_2)-F(x_1)=\int_{x_1}^{x_2}f(u)du\geq 0.\]
This implies that $F:I\to\mathbb{R}$ is an increasing function. \begin{enumerate}[1.]
\item If $I=(a,b]$ or $I=(-\infty, b]$, the limit $\di\lim_{x\to a^+}F(x)$ or the limit $\di\lim_{x\to-\infty}F(x)$ exists if and only if $F(x)$ is bounded below. \item If $I=[a, b)$ or $I=[a, \infty)$, the  limit $\di\lim_{x\to b^-}F(x)$ or the limit $\di\lim_{x\to \infty}F(x)$ exists if and only if $F(x)$ is bounded above. \end{enumerate}
\end{myproof}

In Proposition \ref{230224_10}, we have stated that if the improper integral $\di\int_{-\infty}^{\infty}f(x)dx$ is convergent, then the Cauchy principal value $\di \text{P.V.} \di\int_{-\infty}^{\infty}f(x)dx$  exists. The converse is true if the function $f:\mathbb{R}\to\mathbb{R}$ is nonnegative.
\begin{theorem}{}
Let $f:\mathbb{R}\to\mathbb{R}$ be a nonnegative function that is bounded and Riemann integrable on any closed and bounded intervals. The improper integral $\di\int_{-\infty}^{\infty}f(x)dx$ is convergent if and only if the Cauchy principal value $\di \text{P.V.} \di\int_{-\infty}^{\infty}f(x)dx$  exists. Moreover, 
\[\int_{-\infty}^{\infty}f(x)dx=\text{P.V.} \di\int_{-\infty}^{\infty}f(x)dx.\]

\end{theorem}
\begin{myproof}{Proof}We just need to show that if the Cauchy principal value $\di \text{P.V.} \di\int_{-\infty}^{\infty}f(x)dx$  exists, then the improper integral $\di\int_{-\infty}^{\infty}f(x)dx$ is convergent.\bp
Assume that the Cauchy principal value $\di \text{P.V.} \di\int_{-\infty}^{\infty}f(x)dx$  exists and is equal to $I$. As in the proof of Lemma \ref{230224_6}, the function 
\[F(x)=\int_{0}^xf(u)du\]  is an increasing function.  For any real numbers $b$ and $c$ with $b\leq c$, there is a positive number $a$ such that
\[-a\leq b\leq c\leq a.\]Hence,
\[F(c) -F(b)=\int_b^cf(x)dx\leq\int_{-a}^af(x)dx \leq I.\]
This proves that $-I\leq F(x)\leq I$ for all $x\in \mathbb{R}$. In other words, the function $F:\mathbb{R}\to\mathbb{R}$ is bounded. Therefore, the improper integrals $\di\int_0^{\infty}f(x)dx$ and  $\di\int_{-\infty}^0f(x)dx$ are convergent, and so the improper integral $\di\int_{-\infty}^{\infty}f(x)dx$ is convergent.
\end{myproof}

Now, we can present the comparison theorem for improper integrals.
\begin{theorem}[label=230224_6]{Comparison Theorem}
 Let $I$ be an interval. Given that $f:I\to\mathbb{R}$ and $g:I\to\mathbb{R}$ are nonnegative functions that are bounded and Riemann integrable on any closed and bounded intervals that are contained in $I$. Assume that
\[0\leq f(x)\leq g(x)\hspace{1cm}\text{for all}\;x\in I.\]
\begin{enumerate}[1.]\item If the integral $\di\int_Ig(x)dx$ is convergent, then the integral $\di\int_I f(x)dx$ is convergent. 
\item If the integral $\di\int_If(x)dx$ is divergent, then the integral $\di\int_I g(x)dx$ is divergent. 
\end{enumerate}
\end{theorem}
\begin{myproof}{Proof}Notice that the second statement is the contrapositive of the first statement. Hence, we only need to prove the first statement. 
Fixed $x_0$ in the interval $I$, and define
\[F(x)=\int_{x_0}^x f(u)du,\hspace{1cm} G(x)=\int_{x_0}^x g(u)du.\]
If $x>x_0$,
\[0\leq F(x)\leq G(x).\]
Therefore, $G$ is bounded above implies $F$ is bounded above.
  If $x<x_0$,
\[F(x)=-\int_{x}^{x_0}f(u)du,\hspace{1cm}G(x)=-\int_{x}^{x_0}g(u)du.\] Since
\[0\leq \int_{x}^{x_0}f(u)du\leq \int_{x}^{x_0}g(u)du,\] we find that
\[0\geq F(x)\geq G(x).\]  Therefore, $G$ is bounded below implies that $F$ is bounded below. The assertions about the convergence of the integrals then follow  from Lemma \ref{230224_5}.
\end{myproof}
\begin{example}
{} We can show that the integral $\di\int_0^{\infty}\frac{x}{x^2+1}dx$ is divergent without explicitly computing the integral. Notice that for $x\geq 1$,
\[0\leq \frac{1}{2x}\leq\frac{x}{x^2+1}.\]
Since the integral $\di\int_1^{\infty}\frac{1}{x}dx$ is divergent, the integral $\di \int_1^{\infty}\frac{x}{x^2+1}dx$ is also divergent. Hence,  the integral $\di\int_0^{\infty}\frac{x}{x^2+1}dx$ is divergent.
\end{example}
  \begin{example}{}
Determine whether the improper integral $\di\int_0^{\infty}\frac{1}{\sqrt{x}(x+1)}dx$ is convergent.
\end{example}
\begin{solution}{Solution}
We determine  the convergence of the two improper integrals
\[\int_0^1\frac{1}{\sqrt{x}(x+1)}dx\hspace{1cm}\text{and}\hspace{1cm}\int_1^{\infty}\frac{1}{\sqrt{x}(x+1)}dx\] separately.
For $0<x\leq 1$,
\[0\leq \frac{1}{\sqrt{x}(x+1)} \leq \frac{1}{\sqrt{x}}.\] Since the integral $\di\int_0^1\frac{1}{\sqrt{x}}dx$ is convergent, the integral $\di  \int_0^1\frac{1}{\sqrt{x}(x+1)} dx$ is convergent. 
For $x\geq 1$,
\[0\leq \frac{1}{\sqrt{x}(x+1)} \leq \frac{1}{x\sqrt{x}}=\frac{1}{x^{3/2}}.\]   Since the integral $\di\int_1^{\infty}\frac{1}{x^{3/2}}dx$ is convergent, the integral $\di\int_1^{\infty}  \frac{1}{\sqrt{x}(x+1)} dx$ is convergent. From these, we conclude that the integral  $\di \int_0^{\infty}\frac{1}{\sqrt{x}(x+1)}dx$ is convergent.
\end{solution}

Since the integral $\di\int_0^1x^{-p}dx$ is convergent when $p<1$, while the integral $\di\int_1^{\infty}x^{-p}dx$ is convergent if $p>1$, $\di\int_0^{\infty}x^{-p}dx$ is not convergent for any values of $p$. Hence, to determine the convergence of the integral in the example above, we need to split the integral into two parts and compare to different $g(x)=x^{-p}$. For $x\to 0^+$, we ignore the part $1/(x+1)$ which has a finite limit. For $x\to \infty$, the leading term of $1/(x+1)$ is $1/x$. This is how we identify the correct values of $p$ to compare to.

 Theorem \ref{230224_6}  provides a useful strategy to determine the convergence of an integral in the case that the function is nonnegative. For a function that can take  both positive and negative values, we need other strategies.
 
\begin{theorem}[label=230224_8]{} Let $I$ be an interval.
Assume that $f:I\to\mathbb{R}$ is a function that is bounded and Riemann integrable on any closed and bounded intervals that are contained in $I$. If the improper integral $\di\int_I|f(x)|dx$ is convergent, then the improper integral $\di\int_If(x)dx$ is convergent.
\end{theorem}This theorem can be interpreted as absolute convergence implies convergence.
\begin{myproof}{Proof}
 Define the functions $f_+:I\to\mathbb{R}$ and $f_-:I\to\mathbb{R}$ by
\[f_+(x)=\max\{f(x), 0\},\hspace{1cm}f_-(x)=\max\{-f(x),0\}.\]
In other words,  \[f_+(x)=\begin{cases} f(x),\quad &\text{if}\; f(x)\geq 0,\\0,\quad &\text{if}\;f(x)<0,\end{cases}\]
\[f_-(x)=\begin{cases} -f(x),\quad &\text{if}\; f(x)\leq 0,\\0,\quad &\text{if}\;f(x)>0.\end{cases}\]  Notice that $f_+$ and $f_-$ are nonnegative functions, and
\[f(x)=f_+(x)-f_-(x),\hspace{1cm}|f(x)|=f_+(x)+f_-(x).\]
The second equality implies that
\[0\leq f_+(x)\leq |f(x)|,\hspace{1cm}0\leq f_-(x)\leq |f(x)|\hspace{1cm}\text{for all}\;x\in I.\]Theorem \ref{230221_15} says that the function $|f|:I\to\mathbb{R}$ is Riemann integrable on any closed and bounded intervals that are contained in $I$. \bp Question \ref{ex230224_7} says that the functions $f+:I\to\mathbb{R}$ and $f_-:I\to\mathbb{R}$ are also Riemann integrable on any closed and bounded intervals that are contained in $I$.  By Theorem  \ref{230224_6}, the improper integrals $\di \int_If_+(x)dx$ and $\di \int_If_-(x)dx$ are convergent. By linearity, the improper integral $\di\int_I f(x)dx$ is also convergent.

\end{myproof}

Combining Theorem \ref{230224_6} and Theorem \ref{230224_8}, we have the following.
\begin{theorem}[label=230224_9]{General Comparison Theorem}
 Let $I$ be an interval. Given that $f:I\to\mathbb{R}$ and $g:I\to\mathbb{R}$ are  functions that are bounded and Riemann integrable on any closed and bounded  intervals that are contained in $I$. If
\[|f(x)|\leq g(x)\hspace{1cm}\text{for all}\;x\in I,\]
and the integral $\di\int_Ig(x)dx$ is convergent, then the integral $\di\int_I f(x)dx$ is convergent. 
\end{theorem}

\begin{example}
{}Show that the improper integral $\di\int_1^{\infty}\frac{\sin x}{x^2}dx$ is convergent. 
\end{example}
\begin{solution}{Solution}
For any $x\geq 1$,
\[\left|\frac{\sin x}{x^2}\right|\leq\frac{1}{x^2}.\]
Since the integral $\di\int_1^{\infty}\frac{1}{x^2}dx$ is convergent, the integral  $\di\int_1^{\infty}\frac{\sin x}{x^2}dx$ is convergent. 
\end{solution}

There are some important special functions in mathematics and physics which are defined in terms of improper integrals. One such function is the gamma function, which students have probably seen in probability theory. In fact, gamma function is ubiquitous in mathematics. 

\begin{example}{}
Let $s$ be a real number. Show that the improper integral $\di\int_0^{\infty}t^{s-1}e^{-t}dt$ is convergent if and only if $s>0$.
\end{example}
\begin{solution}{Solution}
We split the integral into the two integrals $\di\int_0^{1}t^{s-1}e^{-t}dt$ and $\di\int_1^{\infty}t^{s-1}e^{-t}dt$. Notice that
\[0\leq t^{s-1}e^{-1}\leq t^{s-1}e^{-t}\leq t^{s-1}\hspace{1cm}\text{for all}\; t\in (0,1].\]
Since $\di\int_0^1t^{s-1}dt$ is convergent if and only if $s>0$, $\di\int_0^1 t^{s-1}e^{-t}dt$ is convergent if and only if $s>0$. For the integral $\di \int_1^{\infty}t^{s-1}e^{-t}dt$, notice that
\[\lim_{t\to \infty}t^{s-1}e^{-t/2}=\lim_{t\to\infty} \frac{t^{s-1}}{e^{t/2}}=0.\]
Therefore, there is a number $t_0>1$ such that for all $t\geq t_0$, $t^{s-1}e^{-t/2}\leq 1$. Now the function
\[g(t)=t^{s-1}e^{-t/2}\] is continuous on the interval $[0, t_0]$. Hence, it is bounded on $[0, t_0]$. These imply that there is a number $M\geq 1$ such that
\[t^{s-1}e^{-t/2}\leq M\hspace{1cm}\text{for all}\;t\geq 1.\]
Hence,
\[0\leq t^{s-1}e^{-t}\leq Me^{-t/2}\hspace{1cm}\text{for all}\;t\geq 1.\]
Since the integral $\di \int_1^{\infty}e^{-t/2}dt$ is convergent, the integral  $\di\int_1^{\infty}t^{s-1}e^{-t}dt$ is convergent.

Hence, the integral  $\di\int_0^{\infty}t^{s-1}e^{-t}dt$ is convergent if and only if $s>0$.

\end{solution}
\begin{highlight}{The Gamma Function}
The gamma function $\Gamma(s)$ is defined as the improper integral
\[\Gamma(s)=\int_0^{\infty}t^{s-1}e^{-t}dt\]when $s>0$. It is easy to find that
\[\Gamma(1)=\int_0^{\infty}e^{-t}dt=1.\]
When $s>0$, using integration by parts with $u(t)=t^s$ and $v(t)=-e^{-t}$, we have
\begin{align*}
\Gamma(s+1)&=\lim_{\substack{a\to 0^+\\b\to\infty}}\int_a^bt^se^{-t}dt\\
&=\lim_{\substack{a\to 0^+\\b\to\infty}}\left\{ \left[-t^se^{-t}\right]_a^b +s\int_a^b t^{s-1}e^{-t}dt\right\}\\
&=\lim_{\substack{a\to 0^+\\b\to\infty}}\left\{a^se^{-a}- b^s e^{-b} \right\}+s\Gamma(s)\\
&=s\Gamma(s).
\end{align*}
This gives the formula
\[\Gamma(s+1)=s\Gamma(s).\]
By induction, one can show that
\[\Gamma(n+1)=n!.\]
Hence, the gamma function is a function that interpolates the factorials. Another special value is
\[\Gamma\left(\frac{1}{2}\right)=\int_0^{\infty}t^{-1/2}e^{-t}dt.\]
Students have probably seen in multivariable calculus or probability that 
\[\int_{-\infty}^{\infty}e^{-x^2}dx=\sqrt{\pi}.\]\end{highlight}
\begin{highlight}{}Making a change of variables $t=u^2$, we find that
\begin{align*}
\int_0^{\infty}t^{-1/2}e^{-t}dt&=\lim_{\substack{a\to 0^+\\b\to\infty}}\int_a^b t^{-1/2}e^{-t}dt\\
&=\lim_{\substack{a\to 0^+\\b\to\infty}}2\int_{\sqrt{a}}^{\sqrt{b}}e^{-u^2}du\\
&=\int_{-\infty}^{\infty}e^{-x^2}dx.
\end{align*}

Hence,
\[\Gamma\left(\frac{1}{2}\right)=\int_0^{\infty}t^{-1/2}e^{-t}dt=\sqrt{\pi}.\]
In the future, we are going to explore more about the gamma function. For example, we will prove the useful formula for the beta integral, 
which says that if $\alpha>0$, $\beta>0$,
\[\int_0^1t^{\alpha-1}(1-t)^{\beta-1}dt=\frac{\Gamma(\alpha)\Gamma(\beta)}{\Gamma(\alpha+\beta)}.\] A lots of other proper or improper integrals can be transformed to this. When $\alpha$ and $\beta$ are positive integers, this formula can be proved by induction. See Question \ref{ex230225_1}.

\end{highlight}
 
\vp
\noindent
{\bf \large Exercises  \thesection}
\setcounter{myquestion}{1}
 \begin{question}{\themyquestion}
Let $a$ be a positive real number. Show that the integral $\di\int_0^{\infty} e^{-ax}dx$ is convergent and find its value.
\end{question}

 \atc

\begin{question}{\themyquestion}
Let $n$ be a positive integer. Find the value of the integral $\di\int_0^{\infty} x^ne^{-x^2}dx$. 
\end{question}
\atc
\begin{question}{\themyquestion}
Explain why the given integral is an improper integral, and determine whether it is convergent. If yes, find the value of the integral.
\begin{enumerate}[(a)]
\item $\di\int_{-3}^0\frac{x}{\sqrt{9-x^2}}dx$
\item $\di\int_0^1\sqrt{x}\ln xdx$
\item $\di \int_0^2\frac{dx}{ (x-1)^2 }$
\end{enumerate}
\end{question}
\atc

\begin{question}{\themyquestion}
Determine whether the   improper integral is convergent. If yes, find its value.
\begin{enumerate}[(a)]
\item $\di\int_0^{\infty}\frac{\sqrt{x}}{x+1}dx$
\item $\di \int_1^{\infty}\frac{\ln x}{x^2}dx$
\end{enumerate}
\end{question}

\atc
\begin{question}{\themyquestion}
Determine whether the   improper integral is convergent.  
\begin{enumerate}[(a)]
\item $\di\int_0^{\infty}\frac{1}{(\sqrt{x}+1)^2}dx$
\item $\di \int_0^1\frac{e^x}{\sqrt{x}}dx$
\item $\di\int_0^{2\pi}\frac{\sin x}{x^{3/2}}dx$
\item $\di\int_{-\infty}^{\infty}\frac{x^3}{(x^2+x+1)^2}dx$
\end{enumerate}
\end{question}

\chapter{Infinite Series of Numbers and Infinite Products }\label{ch5}
In this chapter, we discuss infinite series of numbers and infinite products.

\section{Limit Superior and Limit Inferior}\label{sec5.1}
In Chapter \ref{ch1}, we have seen that a bounded sequence might not be convergent. In this section, we will   discuss the concepts called limit inferiors and limit superiors, which characterize the limits of subsequences of a sequence. 

First we extend the definitions of supremum and infimum as follows. 
\begin{highlight}{Extensions of Infmum and Supremum}
\begin{enumerate}[1.]
\item  If a nonempty set $S$ is not bounded below, we write $\inf S=-\infty$.
\item If a nonempty set $S$ is not bounded above, we write $\sup S=\infty$.

\item $\inf\{-\infty\}=-\infty, \inf\{\infty\}=\infty$.
\item $\sup\{-\infty\}=-\infty, \sup\{\infty\}=\infty$.
\end{enumerate}
\end{highlight}
The definition of limits are also extended to include $-\infty$ and $\infty$ as limits.
 
\begin{theorem}[label=230603_2]{}
Let $\{a_n\}$ be a sequence of real numbers.
\begin{enumerate}[1.]
\item The sequence $\{a_n\}$ is  not bounded above if and only if there is a strictly increasing subsequence $\{a_{n_k}\} $ such that $\di\lim_{k\to\infty}a_{n_k}=\infty$.
\item The sequence $\{a_n\}$ is  not bounded below if and only if there is a strictly decreasing subsequence $\{a_{n_k}\} $ such that $\di\lim_{k\to\infty}a_{n_k}=-\infty$.

\end{enumerate}
\end{theorem}
\begin{myproof}
{Proof}
It is sufficient to prove the first statement. If there is a subsequence $\{a_{n_k}\}$ of $\{a_n\}$ such that $\di\lim_{k\to\infty}a_{n_k}=\infty$, it is obvious that $\{a_n\}$ is not bounded above. 

Conversely,
 given that $\{a_n\}$ is not bounded above, we want to construct a strictly increasing subsequence $\{a_{n_k}\}$  such that $\di\lim_{k\to\infty}a_{n_k}=\infty$.  Let $n_1=1$. 
Since $\{a_n\}$ is not bounded above, there is a $n_2 >1$ such that $a_{n_2}\geq a_{n_1}+1$. Assume we have found $n_1, n_2, \ldots, n_{k-1}$, such that
\[n_1<n_2<\cdots<n_{k-1},\]and
\[a_{n_{j+1}}\geq a_{n_{j}}+1, \hspace{1cm}\text{for all}\;1\leq j\leq k-2.\]Since $\{a_n\}$ is not bounded above, there is an $n_k>n_{k-1}$ such that $a_{n_k}\geq a_{n_{k-1}}+1$. This constructs a strictly increasing sequence $\{a_{n_k}\}$ inductively which satisfies
\[ a_{n_{k+1}}\geq a_{n_{k}}+1, \hspace{1cm}\text{for all}\;k\in\mathbb{Z}^+.\] From this, we find that
\[a_{n_k}\geq a_{n_1}+k-1.\]  Therefore,  $\di\lim_{k\to\infty}a_{n_k}=\infty$.

\end{myproof}
Associated with a given sequence $\{a_n\}$, we can define two sequences $\{b_n\}$ and $\{c_n\}$.
\begin{definition}{}
Given a sequence $\{a_n\}$, we can define two sequences $\{b_n\}$ and $\{c_n\}$ as follows. For each positive integer $n$,  
\[b_n=\inf_{k\geq n}a_k=\inf\{a_k\,|\,k\geq n\},\hspace{1cm}c_n=\sup_{k\geq n}a_k=\sup\{a_k\,|\,k\geq n\}.\]

\end{definition}

\begin{example}[label=ex230226_1]{}
For the sequence $\{a_n\}$ with $a_n=\di\frac{1}{n}$,
\[b_n=0,\hspace{1cm} c_n=\frac{1}{n}\hspace{1cm}\text{for all}\;n\geq 1.\]

\end{example}

\begin{example}[label=ex230226_2]{}
For the sequence $\{a_n\}$ with $a_n=n$,
\[b_n=n,\hspace{1cm} c_n=\infty\hspace{1cm}\text{for all}\;n\geq 1.\]

\end{example}

\begin{example}[label=ex230226_3]{}
For the sequence $\{a_n\}$ with $a_n=(-1)^n$,
\[b_n=-1,\hspace{1cm} c_n=1\hspace{1cm}\text{for all}\;n\geq 1.\]

\end{example}
The following are obvious from the definitions and Theorem \ref{230603_2}.
\begin{proposition}[label=230603_1]{}
Given that $\{a_n\}$ is a sequence of real numbers, for each $n\in\mathbb{Z}^+$, let $\di b_n=\inf_{k\geq n}a_k$ and $\di c_n=\sup_{k\geq n}a_k$.  
\begin{enumerate}[1.]
\item For any positive integer $n$,  $b_n\leq a_n\leq c_n$.
\item  For any positive integer $n$,  $b_n$ cannot be $\infty$, $c_n$ cannot be $-\infty$.
\item $\{a_n\}$ is not bounded below if and only    $b_n=-\infty$ for all $n\geq 1$.
\item $\{a_n\}$ is not bounded above if and only if   $c_n=\infty$ for all $n\geq 1$.
 
\item $\{b_n\}$ is an increasing sequence.
\item $\{c_n\}$ is a decreasing sequence.   
\end{enumerate}
\end{proposition}

Since $\{b_n\}$ is an increasing sequence, $\di\lim_{n\to \infty}b_n=\sup\{b_n\}$ in the general sense.  Similarly, $\di\lim_{n\to \infty}c_n=\inf\{c_n\}$. 
\begin{definition}{Limit Inferior and Limit Superior}
Let $\{a_n\}$ be a sequence of real numbers.
\begin{enumerate}
[1.]\item The limit inferior or limit infimum of $\{a_n\}$, denoted by $\di\liminf_{n\to\infty}a_n$ or $\di\varliminf_{n\to\infty}a_n$, is defined as
\[\liminf_{n\to\infty}a_n=\lim_{n\to \infty}b_n=\sup_{n\geq 1}\inf_{k\geq n}a_k.\]
\item The limit superior or limit supremum of $\{a_n\}$, denoted by $\di\limsup_{n\to\infty}a_n$ or $\di\varlimsup_{n\to\infty}a_n$, is defined as
\[\limsup_{n\to\infty}a_n=\lim_{n\to \infty} c_n=\inf_{n\geq 1}\sup_{k\geq n}a_k.\]
\end{enumerate}
\end{definition}

Notice that using extended definitions of infimum and supremum, the limit infimum and limit supremum of a sequence always exist, either as a finite number, or $\pm \infty$.

\begin{example}{}
\begin{enumerate}[1.]
\item
For the sequence $\{a_n\}$ with $a_n=\di\frac{1}{n}$ defined in Example \ref{ex230226_1},
\[\liminf_{n\to\infty}a_n=0,\hspace{1cm} \limsup_{n\to\infty}a_n=0.\]
\item
For the sequence $\{a_n\}$ with $a_n=n$ defined in Example \ref{ex230226_2},
\[\liminf_{n\to\infty}a_n=\infty, \hspace{1cm} \limsup_{n\to\infty}a_n=\infty.\]
\item
For the sequence $\{a_n\}$ with $a_n=(-1)^n$ defined in Example \ref{ex230226_3},
\[\liminf_{n\to\infty}a_n=-1,\hspace{1cm} \limsup_{n\to\infty}a_n=1.\]
\end{enumerate}
\end{example}
The following are obvious.
\begin{highlight}{}
\begin{enumerate}[1.]\item
$\di \liminf_{n\to\infty}(-a_n)=-\limsup_{n\to\infty}a_n$
\item $\di\limsup_{n\to\infty}(-a_n)=-\liminf_{n\to\infty}a_n$.
\end{enumerate}
\end{highlight}
Since \[b_n\leq c_n\hspace{1cm}\text{for all}\;n\in\mathbb{Z}^+,\]we obtain the following immediately.
\begin{proposition}[label=230227_21]{}
For any sequence $\{a_n\}$,
\[\liminf_{n\to \infty}a_n\leq \limsup_{n\to\infty}a_n.\]
\end{proposition}

We also have the following comparison theorem.
\begin{proposition}{}
Let $\{u_n\}$ and $\{v_n\}$ be sequences of real numbers. If $u_n\leq v_n$ for all positive integers $n$, then
\[\liminf_{n\to\infty}u_n\leq\liminf_{n\to \infty}v_n\hspace{1cm}\limsup_{n\to\infty}u_n\leq\limsup_{n\to \infty}v_n.\]
\end{proposition}

\begin{example}{}Find $\di\liminf_{n\to\infty}a_n$ and $\di \limsup_{n\to \infty}a_n$ for the sequence $\{a_n\}$   defined by
\[a_n=(-1)^n\left(1+\frac{1}{n}\right).\]

\end{example}
\begin{solution}{Solution}
Notice that for any $n\geq 1$,
\[
a_{2n-1}=-1-\frac{1}{2n-1},\hspace{1cm}a_{2n}=1+\frac{ 1}{2n}.
\] \bs We observe that
\[a_1<a_3<\cdots<-1<1<\cdots<a_4<a_2.\]The sequence $\di\{a_{2n-1}\}$ increases to $-1$, while the sequence $\{a_{2n}\}$ decreases to 1. 
Therefore,
\[b_{2n}=b_{2n+1}=-1-\frac{1}{2n+1},\hspace{1cm}c_{2n-1}=c_{2n}=1+\frac{1}{2n}.\]
It follows that
\[\liminf_{n\to\infty}a_n=\lim_{n\to\infty}b_n=-1,\hspace{1cm}
\limsup_{n\to\infty}a_n=\lim_{n\to\infty} c_n= 1.\]
\end{solution}

\begin{example}{}Let $\{a_n\}$ be the sequence defined by
$a_n=(-1)^nn$.
Find $\di\liminf_{n\to\infty}a_n$ and $\di \limsup_{n\to \infty}a_n$.
\end{example}
\begin{solution}{Solution}
For any $n\geq 1$,
\[
a_{2n-1}=-(2n-1),\hspace{1cm}a_{2n}=2n.
\]  The sequence
$\{a_n\}$ is not bounded below nor bounded above. Therefore,
\[b_n=-\infty,\hspace{1cm}c_n=\infty.\]
It follows that
\[\liminf_{n\to\infty}a_n= -\infty,\hspace{1cm}
 \limsup_{n\to\infty}a_n=\infty.\]
\end{solution}

For a monotoic sequence, it is easy to find its limit inferior and limit superior.
\begin{theorem}{}
Let $\{a_n\}$ be a monotonic sequence.
\begin{enumerate}[1.]
\item If $\{a_n\}$ is increasing, then $b_n=a_n$ and $c_n=\sup\{a_n\}$. Therefore,
\[\liminf_{n\to\infty}a_n=\limsup_{n\to\infty}a_n=\lim_{n\to\infty}a_n=\sup\{a_n\}.\]
\item  If $\{a_n\}$ is decreasing, then $b_n=\inf\{a_n\}$ and $c_n=a_n$. Therefore,
\[\liminf_{n\to\infty}a_n=\limsup_{n\to\infty}a_n=\lim_{n\to\infty}a_n=\inf\{a_n\}.\]

\end{enumerate}
\end{theorem}
In other words, for monotonic sequence, the limit inferior, limit superior, and the limit are all the same.

In fact, if a sequence $\{a_n\}$ has a finite limit, then its limit inferior, limit superior and limit are all the same.
\begin{theorem}{}
Let $\{a_n\}$ be a sequence, and let $a$ be a finite number. Then the following two statements are equivalent.
\begin{enumerate}[(a)]
\item   $\di\lim_{n\to\infty}a_n=a$.
\item $\di \liminf_{n\to \infty}a_n=\limsup_{n\to \infty}a_n=a$.
 \end{enumerate}
\end{theorem}
\begin{myproof}
{Proof}
For a positive integer $n$, define
$\di b_n=\inf_{k\geq n}a_k$, $c_n=\sup_{k\geq n}a_n$.
Then
\[\liminf_{n\to \infty}a_n=\lim_{n\to\infty}b_n,\hspace{1cm}\limsup_{n\to \infty}a_n=\lim_{n\to\infty}c_n.\]
By definition,
\[b_n\leq a_n\leq c_n.\]
Hence, (b) implies (a) follows from squeeze theorem.
\bp
Now we prove that (a) implies (b). Given $\varepsilon>0$, since $\di\lim_{n\to\infty}a_n=a$, there is a positive integer $N$ such that for all $n\geq N$, 
\[a-\frac{\varepsilon}{2}<a_n<a+\frac{\varepsilon}{2}.\]
Hence, for all $n\geq N$,
\[ a-\frac{\varepsilon}{2}\leq b_n\leq a+\frac{\varepsilon}{2}\hspace{1cm}\text{and}\hspace{1cm}a-\frac{\varepsilon}{2}\leq c_n\leq a+\frac{\varepsilon}{2}.\]
These prove that for all $n\geq N$,
\[|b_n-a|<\varepsilon\hspace{1cm}\text{and}\hspace{1cm}|c_n-a|<\varepsilon.\]
Thus,
\[\liminf_{n\to \infty}a_n=\limsup_{n\to \infty}a_n=a.\]
\end{myproof}

Hence, we are left to consider sequences which does not have a finite limit.
 Let us first characterize when  the limit inferior and the limit superior of a sequence  can be $-\infty$ or $\infty$. 
\begin{theorem}[label=230226_1]{}
Let $\{a_n\}$ be a sequence of real numbers. Then the following three statements are equivalent.
\begin{enumerate}[(a)]
\item $\di\limsup_{n\to\infty}a_n=\infty$.
\item $\{a_n\}$ is not bounded above.
\item  There is a strictly increasing subsequence $\{a_{n_k}\}$ such that $\di\lim_{k\to\infty}a_{n_k}=\infty$.

\end{enumerate}
\end{theorem}

\begin{myproof}
{Proof}
Notice that $\di\limsup_{n\to\infty}a_n= \lim_{n\to \infty}c_n$, where $\di c_n=\sup_{k\geq n}a_k$. Since $\{c_n\}$ is a decreasing sequence, $\di\lim_{n\to \infty}c_n=\infty$ if and only if $c_n=\infty$ for all $n\geq 1$. Hence, this theorem  follows  from Theorem \ref{230603_2} and Proposition \ref{230603_1}.

\end{myproof}
 
The limit inferior  version of Theorem \ref{230226_1} is straightforward.
\begin{theorem}[label=230226_2]{}
Let $\{a_n\}$ be a sequence of real numbers. Then the following three statements are equivalent.
\begin{enumerate}[(a)]
\item $\di\liminf_{n\to\infty}a_n=-\infty$.
\item $\{a_n\}$ is not bounded below.
\item  There is a strictly decreasing subsequence $\{a_{n_k}\}$ such that $\di\lim_{k\to\infty}a_{n_k}=-\infty$.

\end{enumerate}
\end{theorem}

What is more nontrivial is when  limit superior is $-\infty$, or limit inferior is $\infty$.

\begin{theorem}[label=230226_3]{}
Let $\{a_n\}$ be a sequence of real numbers.
\begin{enumerate}[1.]
\item
$\di\limsup_{n\to\infty}a_n=-\infty$ if and only if   $\di\lim_{n\to\infty}a_n=-\infty$.
\item $\di\liminf_{n\to\infty}a_n= \infty$ if and only if  $\di\lim_{n\to\infty}a_n= \infty$.
\end{enumerate}
\end{theorem}

\begin{myproof}{Proof}Notice that if $\di\limsup_{n\to\infty}a_n=-\infty$, we must have $\di\liminf_{n\to\infty}a_n=-\infty$. Similarly, if $\di\liminf_{n\to\infty}a_n= \infty$, it is necessary that $\di\limsup_{n\to\infty}a_n= \infty$.

It is enough for us to prove the second statement. 
Let $b_n=\di\inf_{k\geq n }a_k$. Notice that $b_n\leq a_n$. If 
 $\di\liminf_{n\to\infty}a_n= \lim_{n\to\infty}b_n=\infty$, we must have
$\di\lim_{n\to\infty}a_n=\infty$.  Conversely, assume that $\di\lim_{n\to\infty}a_n=\infty$. Given $M>0$, there is a positive integer $N$ such that 
\[a_n\geq M\hspace{1cm}\text{for all}\;n\geq M.\]
This implies that
\[b_n\geq M\hspace{1cm}\text{for all}\;n\geq M.\] Hence, $\di\liminf_{n\to\infty} a_n=\lim_{n\to\infty}b_n=\infty$.

\end{myproof}

\begin{example}{}
Consider the sequence $\{a_n\}$ with $a_n=n+(-1)^{n}$. The first few terms are given by
$0, 3, 2, 5, 4, 7, \ldots$.
This sequence is neither increasing nor decreasing. For any $n\geq 1$,
\[b_{2n-1}=b_{2n}=a_{2n-1}=2n-2,\]
while $c_n=\infty$ for all $n\in\mathbb{Z}^+$. Therefore, \[\liminf_{n\to\infty}a_n=\lim_{n\to\infty}b_n=\infty.\]
\end{example}

Combining Theorem \ref{230226_1}, Theorem \ref{230226_2} and Theorem \ref{230226_3}, we can summarize the cases where the limit inferior or the limit superior is $-\infty$ or $\infty$.
\begin{highlight}{Infinities as Limit Superior or Limit Inferior}
Let $\{a_n\}$ be a sequence of real numbers.
\begin{enumerate}[1.]
\item $\di\liminf_{n\to\infty}a_n=\limsup_{n\to\infty}a_n=-\infty$ if and only if $\di\lim_{n\to\infty}a_n=-\infty$. In this case, $\{a_n\}$ is bounded above, not bounded below.
\item $\di\liminf_{n\to\infty}a_n=\limsup_{n\to\infty}a_n= \infty$ if and only if $\di\lim_{n\to\infty}a_n= \infty$. In this case, $\{a_n\}$ is bounded below, not bounded above.
\item $\di\liminf_{n\to\infty}a_n=-\infty$ and $\di\limsup_{n\to\infty}a_n= \infty$ if and only if $\{a_n\}$ is not bounded above nor bounded below.
\item  $-\infty<\di\liminf_{n\to\infty}a_n<\infty$ and $\di\limsup_{n\to\infty}a_n= \infty$ if and only if $\{a_n\}$ is bounded below but not bounded above.
\item $\di\liminf_{n\to\infty}a_n=-\infty$ and $\di -\infty<\limsup_{n\to\infty}a_n< \infty$ if and only if $\{a_n\}$ is bounded above but not bounded below.
\end{enumerate}
\end{highlight}
 
The following gives a relation of limit inferior and limit superior with limits of subsequences.
\begin{theorem}[label=230226_8]{}
Let $\{a_n\}$ be sequence with \[ \liminf_{n\to\infty}a_n=b\quad \text{and}\quad \limsup_{n\to\infty}a_n=c.\] If $\{a_{n_k}\}$ is a subsequence of $\{a_n\}$ that converges to a number $\ell$, then
\[b\leq \ell \leq c.\] Here $b$ and $c$ can be $\pm\infty$.
\end{theorem}
\begin{myproof}{Proof}For every positive integer $n$, let 
\[b_n=\inf_{k\geq n}a_k,\hspace{1cm}c_n=\sup_{k\geq n}a_n.\]
Then
\[b=\lim_{n\to\infty}b_n,\hspace{1cm}c =\lim_{n\to\infty} c_n.\]Given a subsequence $\{a_{n_k}\}$, 
$\{b_{n_k}\}$ and $\{c_{n_k}\}$ are subsequences of the monotonic sequences $\{b_n\}$ and $\{c_n\}$. Hence, we  have
\[\lim_{k\to\infty}b_{n_k}=b,\hspace{1cm}\lim_{k\to\infty}c_{n_k}=c,\]which still holds in the extended sense.
Since \[b_{n_k}\leq a_{n_k}\leq c_{n_k}\hspace{1cm}\text{for all}\;k\in\mathbb{Z}^+,\]we find that
\[b\leq\ell\leq c.\]
\end{myproof}

Now we turn to the case of finite limit superior and finite limit inferior.
We have the following equivalence.
\begin{theorem}[label=230226_5]{}
Let $\{a_n\}$ be a sequence of real numbers. Then  the following two statements are equivalent.
\begin{enumerate}[(a)]
\item $\di\limsup_{n\to\infty}a_n=c$ is finite.
\item Given $\varepsilon>0$, \begin{enumerate}[(i)]\item
there exists a positive integer $N$ such that for all $n\geq N$, $a_n<c+\varepsilon$; and 
\item for every 
  positive integer $N$, there exists an integer $n\geq N$, such that $a_n>c-\varepsilon$.
 \end{enumerate}\end{enumerate}
\end{theorem}

\begin{myproof}{Proof}
For a positive integer $n$, define $\di c_n=\sup_{k\geq n}a_n$, so that $\di\limsup_{n\to\infty}a_n=\lim_{n\to\infty}c_n$.

Let us first prove that (a) implies (b). Given $\varepsilon>0$, since $\di\lim_{n\to\infty}c_n=c$,  there is a positive integer $N$ such that for all $n\geq N$, $
\di |c_n-c|<\varepsilon$.
If $n\geq N$, we find that
\[a_n\leq\sup_{k\geq n}a_k=c_n<c+\varepsilon.\]This proves (b)(i). Now given $\varepsilon>0$,  there is a positive integer $N_0$ such that $c-\varepsilon<c_n<c+\varepsilon$ for all $n\geq N_0$.
For any positive integer $N$, let $N'=\max\{N, N_0\}$. Then $N'\geq N_0$. Hence,
\[ \sup_{k\geq N'}a_k=c_{N'}>c-\varepsilon.\]This implies that $c-\varepsilon$ is not an upper bound of $\{a_k\,|\,k\geq N'\}$. Therefore, there is an $n\geq N'\geq N$ such that $a_n>c-\varepsilon$. This proves (b)(ii).

Now we prove (b) implies (a).  Given $\varepsilon>0$, (b)(i) implies that there is a positive integer $N$ such that for all $n\geq N$, 
\[a_n<c+\frac{\varepsilon}{2}.\]\bp
This implies that for all $n\geq N$,
\[c_n\leq c+\frac{\varepsilon}{2}<c+\varepsilon.\] 
For any $n\geq N$, (b)(ii) implies that there is a $k\geq n$ such that $a_k>c-\varepsilon$. Therefore, $c_n\geq a_k>c-\varepsilon$. This shows that for all $n\geq N$, 
\[c-\varepsilon<c_n<c+\varepsilon.\]
Hence, $\di\limsup_{n\to\infty}a_n=\lim_{n\to\infty}c_n=c$.

\end{myproof}

The limit inferior counterpart of Theorem \ref{230226_5} is the following.
\begin{theorem}[label=230226_11]{}
Let $\{a_n\}$ be a sequence of real numbers. Then  the following two statements are equivalent.
\begin{enumerate}[(a)]
\item $\di\liminf_{n\to\infty}a_n=b$ is finite.
\item Given $\varepsilon>0$, \begin{enumerate}[(i)]\item
there exists a positive integer $N$ such that for all $n\geq N$, $a_n>b-\varepsilon$; and 
\item for every 
  positive integer $N$, there exists an integer $n\geq N$, such that $a_n<b+\varepsilon$.
 \end{enumerate}\end{enumerate}
\end{theorem}

 The following theorem says that the limit inferior and the limit superior are limits of subsequences.

\begin{theorem}[label=230226_9]{}Let $\{a_n\}$ be a sequence.
\begin{enumerate}[1.]

\item If $c=\di\limsup_{n\to\infty}a_n$ is finite, it is the limit of a subsequene of $\{a_n\}$.
\item If $b=\di\liminf_{n\to\infty}a_n$ is finite,  it is the limit of a subsequene of $\{a_n\}$.
\end{enumerate}
\end{theorem}
\begin{myproof}{Proof}
It is sufficient for us to prove the first statement.
We use (b)(ii) of Theorem \ref{230226_5}. Take $\varepsilon=1$ and $N=1$. There is an integer $n_1\geq 1$ such that $a_{n_1}>c-1$.  
Suppose we have chosen $n_1, n_2, \ldots, n_{k-1}$ such that 
$n_1<n_2<\ldots<n_{k-1}$ and 
\[ a_{n_j}>c-\frac{1}{j}\hspace{1cm}\text{for all}\;1\leq j\leq k-1.\]Take $\varepsilon=1/k$ and $N=n_{k-1}+1$. There is an $n_{k}\geq N>n_{k-1}$ such that  $a_{n_k}>c-1/k$.  
 By induction, we have constructed the subsequence $\{a_{n_k}\}$ with 
\[ a_{n_k}>c-\frac{1}{k}\hspace{1cm}\text{for all}\;k\in\mathbb{Z}^+.\]Notice that we also have $a_{n_k}\leq c_{n_k}$. Therefore,
\[c-\frac{1}{k}<a_{n_k}\leq c_{n_k}\hspace{1cm}\text{for all}\;k\in\mathbb{Z}^+.\] Being a subsequence of $\{c_n\}$, $\di\lim_{k\to\infty}c_{n_k}=c$.
By squeeze theorem,
\[\lim_{k\to\infty}a_{n_k}=c.\]

\end{myproof}
\begin{highlight}{}
Theorem \ref{230226_1}, Theorem \ref{230226_2}, Theorem \ref{230226_8} and Theorem \ref{230226_9} give characterization of  limit superior and limit inferior as follows.
 
\begin{enumerate}[1.]
\item[1.] The sequence $\{a_n\}$ is not bounded above  if and only if $\di\limsup_{n\to\infty}a_n=\infty$, if and only if there is a strictly increasing subsequence $\{a_{n_k}\}$ such that $\di\lim_{k\to\infty}a_{n_k}=\infty$.
\item[2.] The sequence$\{a_n\}$ is not bounded below  if and only if $\di\liminf_{n\to\infty}a_n=-\infty$, if and only if there is a strictly decreasing subsequence $\{a_{n_k}\}$ such that $\di\lim_{k\to\infty}a_{n_k}=-\infty$.\end{enumerate}\end{highlight}\begin{highlight}{}\begin{enumerate}[1.]
\item[3.] If the sequence $\{a_n\}$ is   bounded above, there is a subsequence $\{a_{n_k}\}$ such that $\di\lim_{k\to\infty}a_{n_k}=\limsup_{n\to\infty}a_n$.
\item[4.] If the sequence $\{a_n\}$ is   bounded below, there is a subsequence $\{a_{n_k}\}$ such that $\di\lim_{k\to\infty}a_{n_k}=\liminf_{n\to\infty}a_n$.
\item[5.] The limit of any subsequence $\{a_{n_k}\}$ must be between $\di\liminf_{n\to\infty}a_n$ and $\di\limsup_{n\to\infty}a_n$. 
\end{enumerate}
\end{highlight}
Let us look at the following example.
 \begin{example}{}
Find $\di\liminf_{n\to\infty}a_n$ and $\di \limsup_{n\to \infty}a_n$ for the sequence $\{a_n\}$   defined by
\[ a_n=\sin \frac{2\pi n}{5}.\]
\end{example}
\begin{solution}{Solution}
For any $n\geq 1$,  
\[a_{5n}=0,\quad a_{5n+1}=\sin\frac{2\pi}{5}, \quad a_{5n+2}=\sin\frac{4\pi}{5}, \]\[a_{5n+3}=\sin\frac{6\pi}{5}, \quad a_{5n+4}=\sin\frac{8\pi}{5}.\]
Hence, the limit of a convergent subsequence of $\{a_n\}$ can and can only be
\[0, \,\sin\frac{2\pi}{5}, \;\sin\frac{4\pi}{5},\;\sin\frac{6\pi}{5},\;\sin\frac{8\pi}{5}.\] 
Since
\[\sin\frac{8\pi}{5}<\sin\frac{6\pi}{5}<0<\sin\frac{4\pi}{5}<\sin\frac{2\pi}{5},\]we find that
\[\liminf_{n\to\infty}a_n=\sin\frac{8\pi}{5},\hspace{1cm}\limsup_{n\to\infty}a_n=\sin\frac{2\pi}{5}.\] 

\end{solution}
\vp
\noindent
{\bf \large Exercises  \thesection}
\setcounter{myquestion}{1}
 \begin{question}{\themyquestion}
Find $\di\liminf_{n\to\infty}a_n$ and $\di \limsup_{n\to \infty}a_n$ for the sequence $\{a_n\}$, where
\[ a_n=\frac{2n}{n+1}.\]
\end{question}

\atc
 
 \begin{question}{\themyquestion}
Find $\di\liminf_{n\to\infty}a_n$ and $\di \limsup_{n\to \infty}a_n$ for the sequence $\{a_n\}$.
\begin{enumerate}[(a)]
\item $\di a_n=(-1)^n\frac{2n}{2n+1}$
\item $\di a_n=(-1)^n\frac{2n+1}{2n}$
\end{enumerate}
\end{question}

\atc
 
 \begin{question}{\themyquestion}
Find $\di\liminf_{n\to\infty}a_n$ and $\di \limsup_{n\to \infty}a_n$ for the sequence $\{a_n\}$.
\begin{enumerate}[(a)]
\item $\di a_n= -2n+(-1)^n n$
\item $\di a_n=n-2(-1)^n n$
\end{enumerate}
\end{question}

\atc
 \begin{question}{\themyquestion}
Find $\di\liminf_{n\to\infty}a_n$ and $\di \limsup_{n\to \infty}a_n$ for the sequence $\{a_n\}$   defined by
\[a_n=\cos\frac{2\pi n}{9}.\]
\end{question}

\atc
 \begin{question}{\themyquestion}
Prove or disprove: Given two sequences $\{a_n\}$ and $\{b_n\}$, 
\[\limsup_{n\to\infty}(a_n+ b_n)=\limsup_{n\to\infty}a_n+\limsup_{n\to\infty}b_n.\]
\end{question}
\vp

\section{Convergence of Series}\label{sec5.2}

In this section, we consider infinite series and its convergence. A series is a   sum of the form
\[\sum_{n=1}^{\infty}a_n=a_{1}+a_{2}+\cdots+a_n+\cdots,\]where $\{a_n\}_{n=1}^{\infty}$ is an infinite sequence.  Sometimes a series might start with the $n=0$ term.
Since a series is an infinite sum, we need to study whether the sum makes sense. The natural thing to do is to define it using limits.

For convenience, we will deal with series that starts with the $n=1$ term in this chapter. When necessary, we will explain what changes   need to be made if the series starts with the $n=0$ term.

\begin{definition}{Convergence of Series}
Given an infinite series $\di\sum_{n=1}^{\infty}a_n$, we define its $n^{\text{th}}$  {\bf partial sum} $s_n$ by
\[s_n=\sum_{k=1}^n a_k.\]We say that the series is {\bf convergent} or has a finite sum, if the sequence $\{s_n\}$ has a finite limit. Otherwise, we say that the series is {\bf divergent}.
If the series is convergent, we define its sum by
\[\sum_{n=1}^{\infty}a_n=s=\lim_{n\to\infty}s_n=\lim_{n\to\infty}\sum_{k=1}^na_k.\]
 
\end{definition}If the infinite   series   $\di\sum_{n=0}^{\infty}a_n$   starts with the $n=0$ term, we still define its $n^{\text{th}}$  partial sum  by
\[s_n=\sum_{k=0}^n a_k\hspace{1cm}\text{when}\;n\geq 0.\] 

\begin{highlight}{}
The convergence of a series $\di\sum_{n=1}^{\infty}a_n$ is not affected by a finite number of terms in the series. Given a positive integer $n_0$, the series $\di\sum_{n=n_0}^{\infty}a_n$ is convergent if and only if the series $\di\sum_{n=1}^{\infty}a_n$ is convergent. In case they are convergent, 
\[\sum_{n=1}^{\infty}a_n-\sum_{n=n_0}^{\infty}a_n=\sum_{n=1}^{n_0-1}a_n.\]
\end{highlight}

\begin{example}{}
For the series $\di\sum_{n=1}^{\infty}\frac{1}{2^n}$, $a_n=\di\frac{1}{2^n}$, and the $n^{\text{th}}$ partial sum is 
\[s_n=\frac{1}{2}+\frac{1}{2^2}+\cdots+\frac{1}{2^n}=1-\frac{1}{2^n}.\]
Since $\di\lim_{n\to\infty}\frac{1}{2^n}=0$, we find that $\di\lim_{n\to\infty}s_n=1$. Hence, the series  $\di\sum_{n=1}^{\infty}\frac{1}{2^n}$ is convergent and 
\[  \sum_{n=1}^{\infty}\frac{1}{2^n}=1.\]
\end{example}

\begin{example}[label=230227_1]{Harmonic Series}
Determine whether the series $\di\sum_{n=1}^{\infty}\frac{1}{n}$ is convergent.
\end{example}
\begin{solution}{Solution}
The $n^{\text{th}}$ partial sum of the series is 
\[s_n=1+\frac{1}{2}+\cdots+\frac{1}{n}.\]\bs
For any positive integer $n$,
\[s_{2n}-s_n=\frac{1}{n+1}+\cdots+\frac{1}{2n}\geq \frac{n}{2n}=\frac{1}{2}.\]
Therefore, for any positive integer $k$,
\[s_{2^k}=s_{2^k}-s_{2^{k-1}}+s_{2^{k-1}}-s_{2^{k-2}}+\cdots+s_2-s_1+s_1\geq 1+\frac{k}{2}.\] This shows that the sequence $\{s_n\}$ is not bounded above. Hence, it is not convergent. Therefore, the series $\di\sum_{n=1}^{\infty}\frac{1}{n}$ is  divergent.
\end{solution}
 Example \ref{230227_1} is a typical example where we determine the convergence of a series without compute the exact value of its partial sum. In this section, we are going to learn various strategies that can be used to do so.

From linearity of limits, we immediately deduce the following.
\begin{proposition}{Linearity}
Let $\di\sum_{n=1}^{\infty}a_n$ and $\di\sum_{n=1}^{\infty}b_n$ be convergent series. Then for any constants $\alpha$ and $\beta$, the series $\di\sum_{n=1}^{\infty}(\alpha a_n+\beta b_n)$ is also convergent and
\[\sum_{n=1}^{\infty}(\alpha a_n+\beta b_n)=\alpha\sum_{n=1}^{\infty}a_n+\beta\sum_{n=1}^{\infty}b_n.\]
\end{proposition}
 
Let us  look at a simple criteria that can be used to conclude that a series is divergent.
\begin{theorem}{}
 
If a series $\di\sum_{n=1}^{\infty}a_n$ is convergent, then $\di\lim_{n\to\infty}a_n=0$. 
Equivalently, if $\di\lim_{n\to\infty}a_n\neq 0$, then the series  $\di\sum_{n=1}^{\infty}a_n$ is  divergent. 
\end{theorem} 
\begin{myproof}{Proof}
We just need to prove the first statement. If $\di\sum_{n=1}^{\infty}a_n$ is convergent, the sequence of partial sums $\{s_n\}$ converges to a number $s$. Notice that $a_n=s_n-s_{n-1}$, where $s_0=0$ by default. Therefore,
\[\lim_{n\to\infty}a_n=\lim_{n\to\infty}s_n-\lim_{n\to\infty}s_{n-1}=s-s=0.\]
\end{myproof}

When one is determining the convergence of a series $\di\sum_{n=1}^{\infty}a_n$, it is always good to start with checking whether the limit  $\di\lim_{n\to\infty}a_n $ is zero.
\begin{example}{}
The series $\di\sum_{n=1}^{\infty}(-1)^n$ is divergent since the limit $\di\lim_{n\to\infty}(-1)^n$ does not exist.
\end{example}

\begin{example}{}
Determine the convergence of the series 
\[\sum_{n=1}^{\infty}\left(1+\frac{1}{n}\right)^n.\]
\end{example}
\begin{solution}{Solution}
Since 
\[\lim_{n\to\infty}\left(1+\frac{1}{n}\right)^n=e\neq 0,\]
the series 
\[\sum_{n=1}^{\infty}\left(1+\frac{1}{n}\right)^n\] is divergent.
\end{solution}

The geometric series is a series which we can find the partial sums explicitly. It is   useful for comparisons.
\begin{theorem}[label=230227_2]{Geometric Series}
The geometric series $\di\sum_{n=0}^{\infty}r^n$ is convergent if and only if $|r|<1$. Moreover,
\[\sum_{n=0}^{\infty}r^n=\frac{1}{1-r}\hspace{1cm}\text{when}\;|r|<1.\]
\end{theorem}
\begin{myproof}{Proof}
When $|r|\geq 1$, the limit $\di\lim_{n\to\infty}r^n$ is not 0. Thus, the series $\di\sum_{n=0}^{\infty}r^n$ is divergent. 

When $|r|<1$, the $n^{\text{th}}$ partial sum is 
\[s_n=1+r+r^2+\cdots+r^n=\frac{1-r^{n+1}}{1-r}.\]
In this case,
$\di \lim_{n\to \infty}r^{n+1}=0$, and so
\[\lim_{n\to\infty}s_n=\lim_{n\to\infty} \frac{1-r^{n+1}}{1-r}=\frac{1}{1-r}.\]
Hence, the series $\di\sum_{n=0}^{\infty}r^n$ is convergent when $|r|<1$, and 
\[\sum_{n=0}^{\infty}r^n=\frac{1}{1-r}\hspace{1cm}\text{when}\;|r|<1.\]
\end{myproof}

\begin{highlight}{}
A general geometric series is a series of the form $\di \sum_{n=1}^{\infty} \left(ar^{n-1}\right)$, where $a$ is the first term of the series. It is convergent if and only if $|r|<1$. 
\end{highlight}

If all the terms $a_n$ in the series $\di \sum_{n=1}^{\infty}a_n$ are nonnegative, we notice that the partial sums $\{s_n\}$ form an increasing sequence. For an increasing sequence, we have the monotone convergence theorem. Applying to the sequence of partial sums, we have the following.
\begin{theorem}[label=230227_3]{}
If $a_n\geq 0$ for all $n\geq 1$, then the series $\di \sum_{n=1}^{\infty}a_n$ is convergent if and only if the sequence of partial sums $\{s_n\}$ is bounded above.
\end{theorem} 
In practice, it is sufficient that there is a positive integer $N$ so that $a_n\geq 0$ for all $n\geq N$.
In Example \ref{230227_1}, we have used this criterion to show that the harmonic series $\di\sum_{n=1}^{\infty}\frac{1}{n}$ is divergent.

Besides the geometric series, a series that is useful for comparisons is the $p$-series $\di\sum_{n=1}^{\infty}\frac{1}{n^p}$. When $p\leq 0$, this series is not convergent since $\di\lim_{n\to\infty}\frac{1}{n^p}\neq 0$. Hence, we will concentrate on the case where $p>0$. To determine the convergence of this series, a convenient tool is the integral test.
\begin{theorem}{Integral Test}
Suppose that $f:[1, \infty)\to\mathbb{R}$ is a function that satisfies the following conditions.
\begin{enumerate}[(i)]
\item $f$ is continuous.
\item $f$ decreases to 0 monotonically. 

\end{enumerate} For $n\geq 1$, let $a_n=f(n)$. Then the series $\di\sum_{n=1}^{\infty}a_n$ is convergent if and only if the improper integral $\di\int_1^{\infty}f(x)dx$ is convergent.
\end{theorem}

\begin{myproof}{Proof}
Since $f(x)$ decreases to 0 monotonically, $f(x)\geq 0$ for all $x\geq 1$,  and so $\{a_n\}$ is a nonnegative decreasing sequence with $\di\lim_{n\to\infty}a_n=0$.   Let $s_n=a_1+a_2+\cdots+a_n$ be the $n^{\text{th}}$ partial sum of the series, and let $\di F(x)=\int_1^x f(u)du$ when $x\geq 1$.\bp

Given a positive integer $n$, since $f$ is decreasing, we find that
\[f(n+1)\leq f(x)\leq f(n)\hspace{1cm}\text{for}\;n\leq x\leq n+1.\]
This implies that
\[a_{n+1}\leq\int_n^{n+1}f(x)dx\leq a_{n}.\]
Therefore, when $n\geq 2$,
\begin{align*}&\int_1^2f(x)dx+\cdots+\int_{n-1}^nf(x)dx+\int_{n}^{n+1}f(x)dx\\&\leq s_n\leq a_1+\int_1^2f(x)dx+\cdots+\int_{n-1}^nf(x)dx.\end{align*}
This gives
\begin{equation}\label{eq230227_8}F(n+1)\leq s_n\leq a_1+F(n).\end{equation}
If the improper integral $\di\int_1^{\infty}f(x)dx$ is convergent, $\{F(n)\}$ is bounded above by a number $M$. Therefore,
\[s_n\leq a_1+M\hspace{1cm}\text{for all}\;n\geq 1,\] and so the sequence $\{s_n\}$ is bounded above. Therefore, the series $\di\sum_{n=1}^{\infty}a_n$ is convergent.

If the improper integral $\di\int_1^{\infty}f(x)dx$  is divergent,  $\di\lim_{n\to\infty}F(n+1)=\infty$. By \eqref{eq230227_8},  we find  that  the sequence $\{s_n\}$ is not bounded above. Thus, the series $\di\sum_{n=1}^{\infty}a_n$ is divergent.
\end{myproof}

\begin{figure}[ht]
\centering
\includegraphics[scale=0.2]{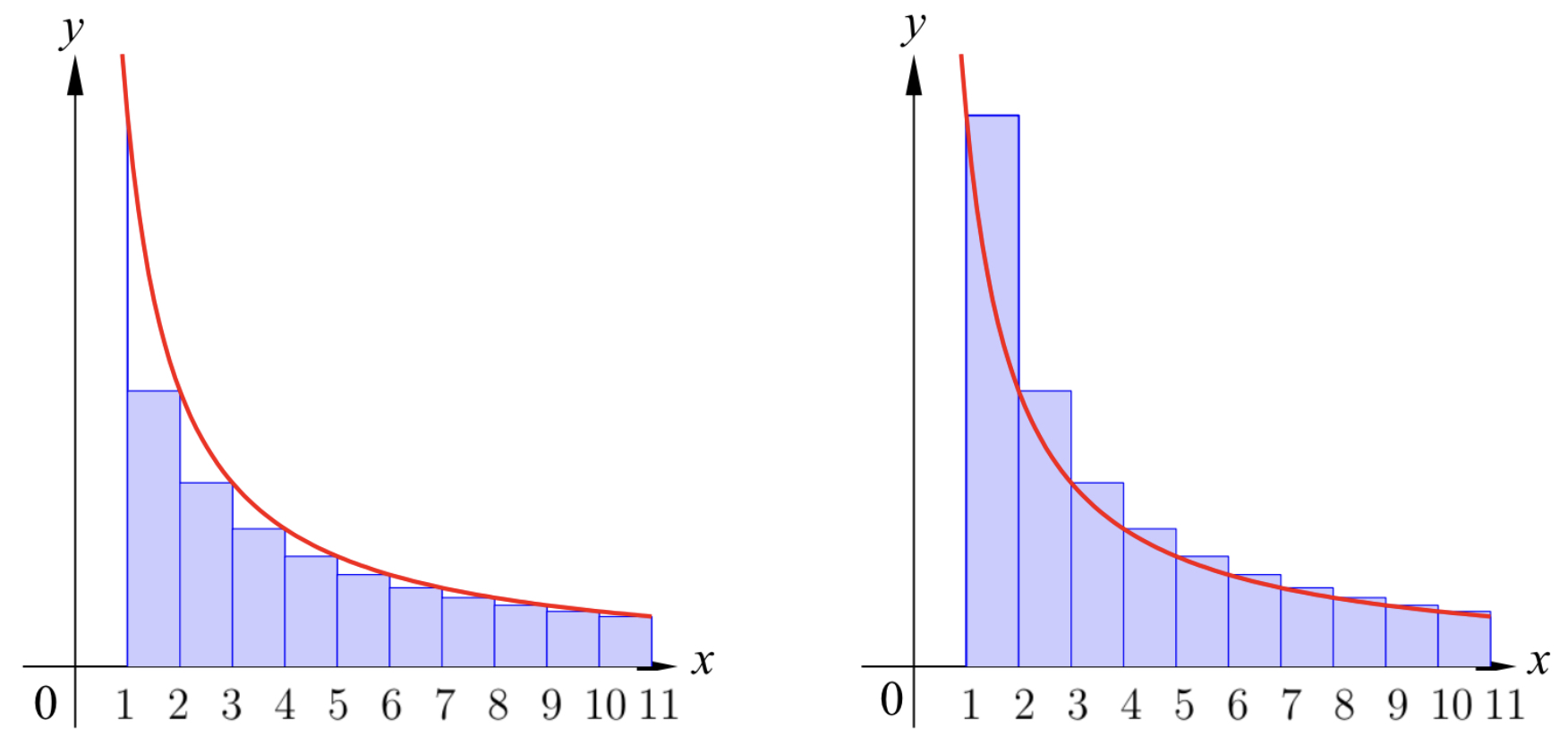}
\caption{The integral test.\fa}\label{figure48}
\end{figure}

Let us now use the integral test to determine the convergence of the $p$-series.
\begin{theorem}{p-Series}
Let $p$ be a positive number. The $p$-series $\di\sum_{n=1}^{\infty}\frac{1}{n^p}$ is convergent if and only if $p>1$.

\end{theorem}
\begin{myproof}{Proof}
Define the function $f:[1,\infty)\to\mathbb{R}$ by 
\[f(x)=\frac{1}{x^p}.\] Then $f$ is a continuous function that decreases monotonically to 0. By Example \ref{ex230227_10}, the improper integral $\di \int_1^{\infty}\frac{1}{x^p}dx$ is convergent if and only if $p>1$. By integral test, the series $\di\sum_{n=1}^{\infty}\frac{1}{n^p}$ is convergent if and only if $p>1$.
\end{myproof}

\begin{example}{}
The series $\di\sum_{n=1}^{\infty}\frac{1}{\sqrt{n}}$ is divergent since it is a $p$-series with $p=\frac{1}{2}\leq 1$.
\end{example}
\begin{remark}
{Integral Approximation to Partial Sums}
Given that $f:[1,\infty)\to\mathbb{R}$ is a continuous function that monotonically decreases to 0, let 
\[s_n=\sum_{k=1}^n f(k),\hspace{1cm} t_n=\int_1^nf(x)dx.\] 
From  the proof of the integral test, we have
\[t_n+f(n+1)\leq s_n\leq a_1+t_n.\]This implies that
\[f(n+1)\leq s_n-t_n\leq a_1,\]which gives bounds for the error when   the partial sum $s_n=\di\sum_{k=1}^na_k$ is approximated by the integral $\di\int_1^nf(x)dx$.
When the improper integral $\di\int_1^{\infty}f(x)dx$ is convergent, the sum $\di\sum_{n=1}^{\infty}a_n$ is also convergent. In this case, the sum of the infinite series satisfies
\[\int_1^{\infty}f(x)dx\leq \sum_{n=1}^{\infty}a_n\leq \int_1^{\infty}f(x)dx+a_1.\]
 If we use $s_n$ to approximate the sum $s=\di\sum_{n=1}^{\infty}a_n$, the error is
\[s-s_n=\sum_{k=n+1}^{\infty}a_k.\]
The same reasoning shows that if $n\geq 1$,
\[\int_{n+1}^{\infty}f(x)dx\leq s-s_n\leq\int_n^{\infty}f(x)dx.\]
\end{remark}

\begin{example}[label=ex230228_7]{Euler's Constant}
We can prove that the limit
\[\lim_{n\to\infty}\left(1+ \frac{1}{2}+\cdots+\frac{1}{n}-\ln n\right)\]exists as follows. Let
\[c_n=1+ \frac{1}{2}+\cdots+\frac{1}{n}-\ln n.\] Then
\[\ln n=\int_1^n\frac{1}{x}dx\leq \int_1^2\frac{1}{x}dx+\cdots+\int_{n-1}^{n}\frac{1}{x}dx
\leq 1+ \frac{1}{2}+\cdots+\frac{1}{n-1}.\]
Therefore, $c_n\geq 0$ for all $n\geq 1$. On the other hand,
\[c_{n+1}-c_n=\frac{1}{n+1}-\ln (n+1)+\ln n=\frac{1}{n+1}-\int_n^{n+1}\frac{1}{x}dx\leq 0.\]Hence, $\{c_n\}$ is a decreasing sequence that is bounded below by 0. By monotone convergence theorem, $\{c_n\}$ converges to a limit $\gamma$. This number \[\gamma=\lim_{n\to\infty}\left(1+ \frac{1}{2}+\cdots+\frac{1}{n}-\ln n\right)\] is called the Euler-Mascheroni constant, or simply as Euler's constant. It is an important constant in mathematics. Numerically, it is equal to
\[0.577215664901532\]correct to 15 decimal places.

\end{example}

Now we return to the comparison test. 
Using Theorem \ref{230227_3}, we obtain the following   test for nonnegative series.
\begin{theorem}{Comparison Test}
Let $\di\sum_{n=1}^{\infty}a_n$ and $\di\sum_{n=1}^{\infty}b_n$ be two series satisfying
\[0\leq a_n\leq b_n\hspace{1cm}\text{for all}\;n\geq 1.\]
\begin{enumerate}[1.]
\item
If $\di\sum_{n=1}^{\infty}b_n$ is convergent, $\di\sum_{n=1}^{\infty}a_n$ is convergent.
\item
If $\di\sum_{n=1}^{\infty}a_n$ is divergent, $\di\sum_{n=1}^{\infty}b_n$ is divergent.
\end{enumerate}
\end{theorem}
\begin{myproof}{Proof}
Let
$s_n=a_1+\ldots+a_n$ and $t_n=b_1+\ldots+b_n$  be respectively the $n^{\text{th}}$ partial sums of the series $\di\sum_{n=1}^{\infty}a_n$ and $\di\sum_{n=1}^{\infty}b_n$.   Then $\{s_n\}$ and $\{t_n\}$ are increasing sequences and \[s_n\leq t_n.\]
\begin{enumerate}[1.]\item 
If $\di\sum_{n=1}^{\infty}b_n$ is convergent, the sequence $\{t_n\}$ is bounded above. Then the sequence $\{s_n\}$  is also bounded above. Hence, $\di\sum_{n=1}^{\infty}a_n$ is convergent.
\item 
If $\di\sum_{n=1}^{\infty}a_n$ is divergent, the sequence $\{s_n\}$ is not bounded above. Then the sequence $\{t_n\}$  is also not bounded above. Hence, $\di\sum_{n=1}^{\infty}b_n$ is divergent.
\end{enumerate}
\end{myproof}

\begin{example}[label=ex230227_6]{}
Determine the convergence of the series 
\[\sum_{n=1}^{\infty}\frac{2^n}{3^n-1}.\]

\end{example}
\begin{solution}{Solution}
For $n\geq 1$, 
\[3^n-1\geq \frac{1}{2}\times 3^n.\]
Therefore,
\[\frac{2^n}{3^n-1}\leq   \frac{2^{n+1 }}{3^n}.\]Since the series \[\sum_{n=1}^{\infty}\frac{2^{n+1}}{3^n}=2\sum_{n=1}^{\infty}\frac{2^{n}}{3^n}\] is a geometric series with $r=2/3$, it is convergent. By comparison test, the series $\di\sum_{n=1}^{\infty}\frac{2^n}{3^n-1}$ is convergent.
\end{solution}

\begin{example}[label=ex230227_7]{}
Determine the convergence of the series 
\[\sum_{n=1}^{\infty}\frac{n}{n\sqrt{n}+1}.\]

\end{example}
\begin{solution}{Solution}For $n\geq 1$,\[
\frac{n}{n\sqrt{n}+1}\geq\frac{n}{n\sqrt{n}+n\sqrt{n}}=\frac{1}{2 \sqrt{n}}.\]
Since the series $\di\sum_{n=1}^{\infty}\frac{1}{\sqrt{n}}$ is a $p$-series with $p=1/2\leq 1$, it is divergent. So   the series $\di\sum_{n=1}^{\infty}\frac{1}{2\sqrt{n}}$ is also divergent. By comparison test, the series
$\di \sum_{n=1}^{\infty}\frac{n}{n\sqrt{n}+1}$ is divergent.
\end{solution}

In applying the comparison test, we need to identify the correct series to compare to, and prove some strict inequalities. In Example \ref{ex230227_6}, we compare  $a_n=\di \frac{2^n}{3^n-1}$ to $b_n=\di \frac{2^n}{3^n}$, since $2^n$ and $3^n$ are the leading terms of the numerator and the denominator of $a_n$ when $n$ is large. Since we know that the series $\di\sum_{n=1}^{\infty}\frac{2^n}{3^n}$ is convergent, we need to prove that $a_n $ is up to a constant, less than or equal to $b_n$, in order to use the comparison test to conclude that $\di\sum_{n=1}^{\infty}a_n$ is convergent.

Simiarly, for Example \ref{ex230227_7}, we compare  $a_n=\di\frac{n}{n\sqrt{n}+1}$ to $\di b_n=\frac{n}{n\sqrt{n}}$ since $n$ and $n\sqrt{n}$ are respectively the leading terms of the numerator and denominator of $a_n$. Since $\di\sum_{n=1}^{\infty}b_n$ is divergent, so we want to conclude that $\di\sum_{n=1}^{\infty}a_n$ is divergent. For this, we need to show that $a_n$ is larger than a constant times $b_n$.

Proving strict inequalities is tedious, and we see that it might not be necessary. In fact, we obtain the series to compare to by investigating the leading terms. This is somehow a limit. Hence, we can replace the comparison test by limit comparison test.

\begin{theorem}{Limit Comparison Test}
Given the two series  $\di\sum_{n=1}^{\infty}a_n$ and $\di\sum_{n=1}^{\infty}b_n$ that satisfy the following conditions.
 
\begin{enumerate}[(i)]
\item $a_n\geq 0$ and $b_n>0$ for all $n\in\mathbb{Z}^+$.
\item The limit $\di L=\lim_{n\to\infty} \frac{a_n}{b_n} $ exists and is finite. 
\end{enumerate}Since $\di \frac{a_n}{b_n}\geq 0$, we must have $L\geq 0$.
\begin{enumerate}[1.]
\item If $L=0$, and the  series $\di\sum_{n=1}^{\infty}b_n$ is convergent, then the series $\di\sum_{n=1}^{\infty}a_n$ is convergent.
\item If  $L>0$,   the  series $\di\sum_{n=1}^{\infty}a_n$ is convergent if and only if the series $\di\sum_{n=1}^{\infty}b_n$ is convergent.
\end{enumerate}

\end{theorem} The   condition  (ii)  says that when $n$ is large, $a_n$ is smaller than or equal to a multiple of   $b_n$.  
\begin{myproof}{Proof}
First consider the case where $L=0$.
By definition of limit with $\varepsilon=1$, there is a positive integer $N$ such that for all $n\geq N$, 
\[\left| \frac{a_n}{b_n} \right|< 1.\]
Therefore,
\[0\leq a_n\leq b_n\hspace{1cm}\text{for all}\; n\geq N.\]
Since the series $\di\sum_{n=1}^{\infty}b_n$ is convergent,  the series $\di\sum_{n=N}^{\infty}  b_n$ is convergent. Then comparison test implies that the series $\di \sum_{n=N}^{\infty} a_n $ is convergent. Therefore, the series $\di \sum_{n=1}^{\infty} a_n $ is convergent.  

Now for the case $L>0$, take $\varepsilon=L/2$. There is a positive integer $N$ such that for all $n\geq N$, 
\[\left| \frac{a_n}{b_n} -L\right|<\frac{L}{2}.\]
This implies that
\[0\leq \frac{L}{2}b_n \leq a_n\leq\frac{3L}{2} b_n\hspace{1cm}\text{for all}\; n\geq N.\]
Comparison test then shows that  the  series $\di\sum_{n=1}^{\infty}a_n$ is convergent if and only if the series $\di\sum_{n=1}^{\infty}b_n$ is convergent.
\end{myproof}

\begin{example}{\linkt Example \ref{ex230227_6} Revisited}
For the series $\di \sum_{n=1}^{\infty}\frac{2^n}{3^n-1}$ considered in Example \ref{ex230227_6}, we take $a_n=\di\frac{2^n}{3^n-1}$ and $\di b_n=\frac{2^n}{3^n }$.\be Then
\[\lim_{n\to\infty}\frac{a_n}{b_n}=\lim_{n\to\infty} \frac{1}{1-\di\frac{1}{3^n}}=1.\] 
 Since the series $\di\sum_{n=1}^{\infty}\frac{2^n}{3^n}$  is convergent, the series $\di\sum_{n=1}^{\infty}\frac{2^n}{3^n-1}$ is convergent.
\end{example2}

\begin{example}{\linkt Example \ref{ex230227_7} Revisited}
For the series $\di \sum_{n=1}^{\infty}\frac{n}{n\sqrt{n}+1}$ considered in Example \ref{ex230227_7}, we take $a_n=\di\frac{n}{n\sqrt{n}+1}$ and $\di b_n=\frac{1}{\sqrt{n}}$.  Then
\[\lim_{n\to\infty}\frac{a_n}{b_n}=\lim_{n\to\infty} \frac{1}{1+\di \frac{1}{n\sqrt{n}}}=1.\] 
 Since the series $\di\sum_{n=1}^{\infty}\frac{1}{\sqrt{n}}$  is divergent, the series $\di\sum_{n=1}^{\infty}\frac{n}{n\sqrt{n}+1}$ is divergent.
\end{example}

Let us now turn to series that can have negative terms. First we formulate a Cauchy criterion for convergence of series. Recall that a sequence $\{s_n\}$ is a Cauchy sequence if for every $\varepsilon>0$, there is a positive integer $N$ so that for all $m\geq n\geq N$, 
\[|s_m-s_n|<\varepsilon.\]
Applying the Cauchy criterion for convergence of sequences (see Theorem \ref{23020602}), and the fact that if $m\geq  n>1$,
\[s_m-s_{n-1}=a_{n}+a_{n+1}+\cdots+a_m,\] we obtain the following Cauchy criterion for convergence of infinite series.
\begin{theorem}{Cauchy Criterion for Infinite Series}
An infinite series $\di\sum_{n=1}^{\infty}a_n$ is convergent if and only if for every $\varepsilon>0$, there is a positive integer $N$ such that for all $m\geq n\geq N$,
\[|a_n+a_{n+1}+\cdots+a_m|<\varepsilon.\]
\end{theorem}

Using this, we can prove the following.

\begin{theorem}[label=230227_11]{}
If the series $\di\sum_{n=1}^{\infty}|a_n|$ is convergent, then the series $\di\sum_{n=1}^{\infty}a_n$ is convergent.
\end{theorem}
\begin{myproof}{Proof}
Given $\varepsilon>0$, since the series $\di\sum_{n=1}^{\infty}|a_n|$ is convergent, there is a positive integer $N$ such that for all $m\geq n\geq N$,
\[|a_n|+|a_{n+1}|+\cdots+|a_m|=\left||a_n|+|a_{n+1}|+\cdots+|a_m|\right|<\varepsilon.\] 
By triangle inequaltiy, we find that for $m\geq n\geq N$, 
\[|a_n+a_{n+1}+\cdots+a_m|<|a_n|+|a_{n+1}|+\cdots+|a_m|<\varepsilon.\]
Using Cauchy criterion, we conclude that the series $\di\sum_{n=1}^{\infty}a_n$ is convergent.
\end{myproof}

The converse of Theorem \ref{230227_11} is not true. Namely, there exists series $\di\sum_{n=1}^{\infty}a_n$ which is  convergent but the corresponding absolute series $\di\sum_{n=1}^{\infty}|a_n|$ is not convergent. Therefore, let us make the following definitions.
\begin{definition}{Absolute Convergence and Conditional Convergence}
Given that the series $\di\sum_{n=1}^{\infty}a_n$ is convergent.
\begin{enumerate}[1.]
\item
We say that the series $\di\sum_{n=1}^{\infty}a_n$ {\bf converges absolutely} if the series $\di\sum_{n=1}^{\infty}|a_n|$ is convergent.
\item We say that the series $\di\sum_{n=1}^{\infty}a_n$ {\bf converges conditionally} if the series $\di\sum_{n=1}^{\infty}|a_n|$ is divergent.
\end{enumerate}
\end{definition}

\begin{example}{}
If $p>1$, the series $\di\sum_{n=1}^{\infty}\frac{(-1)^n}{n^p}$ converges absolutely.
\end{example}

From the limit comparison test, we have the following. 
\begin{theorem}{Limit Comparison Test II}
Given the  series  $\di\sum_{n=1}^{\infty}a_n$, assume that there is a series $\di\sum_{n=1}^{\infty}b_n$ such that $b_n>0$ for all $n\in\mathbb{Z}^+$, and the limit $\di L=\lim_{n\to\infty} \frac{|a_n|}{b_n} $ exists and is finite. 
 If  the series $\di\sum_{n=1}^{\infty}b_n$ is convergent, then the  series $\di\sum_{n=1}^{\infty}a_n$ converges absolutely.
 
\end{theorem}

\begin{example}{}
Show that the series
\[\sum_{n=1}^{\infty} \frac{2^n+(-1)^n3^n}{5^n+1}\] is convergent.
\end{example}
\begin{solution}{Solution}
Let
\[a_n=\frac{2^n+(-1)^n3^n}{5^n+1},\hspace{1cm}b_n=\frac{3^n}{5^n}.\]
Then $b_n>0$ for all $n\in\mathbb{Z}^+$ and $\di\sum_{n=1}^{\infty}b_n$ is convergent. Now,
\[\lim_{n\to\infty}\frac{|a_n|}{b_n}=\lim_{n\to\infty}\frac{\di 1+(-1)^n\left(\frac{2}{3}\right)^n}{1+\di \frac{1}{5^n}}=1.\]
Therefore, the series $\di \sum_{n=1}^{\infty} \frac{2^n+(-1)^n3^n}{5^n+1}$ converges absolutely, and thus is convergent.
\end{solution}

 To give  an example of series that converges conditionally, let us   discuss a convergence test for a special class of series called alternating series.

\begin{definition}{Alternating Series}
A series of the form \[\sum_{n=1}^{\infty}(-1)^{n-1}b_n=b_1-b_2+b_3-b_4+\cdots+b_{2n-1}-b_{2n}+\cdots,\]where $b_n\geq 0$ for all $n\geq 1$, is called an alternating series.
\end{definition}

\begin{example}[label=ex230227_12]{}
The series 
\[1-\frac{1}{2}+\frac{1}{3}-\frac{1}{4}+\cdots\] is an alternating series.
\end{example} 
A necessary condition for an alternating series $\di\sum_{n=1}^{\infty}(-1)^{n-1}b_n$ to be convergent is $\di\lim_{n\to\infty}b_n=0$. The following theorem says that if $\{b_n\}$ is also decreasing, then the alternating series is convergent. 

\begin{theorem}[label=230227_15]{Alternating Series Test}
If  $\{b_n\}$ is a monotonically decreasing sequence with $\di\lim_{n\to \infty}b_n=0$, the alternating series $\di\sum_{n=1}^{\infty}(-1)^{n-1}b_n$ is convergent.
\end{theorem}
\begin{myproof}{Proof}Since $\{b_n\}$ decreases monotonically to 0, $b_n\geq 0$ for all $n\in\mathbb{Z}^+$. 
Let $a_n=(-1)^{n-1}b_n$ be the $n^{\text{th}}$ term of the series $\di\sum_{n=1}^{\infty}(-1)^{n-1}b_n$ , and let $s_n=a_1+a_2+\cdots+a_n$ be the $n^{\text{th}}$ partial sum. We are given that 
\[b_1\geq b_2\geq \cdots\geq b_n\geq b_{n+1}\cdots.\]
Therefore,
\begin{align*}
s_{2n+1}&= s_{2n-1}+a_{2n}+a_{2n+1}=s_{2n-1}-(b_{2n}-b_{2n+1})\leq s_{2n-1},\\
s_{2n+2}&=s_{2n}+a_{2n+1}+a_{2n+2} =s_{2n}+(b_{2n+1}-b_{2n+2})\geq s_{2n}. 
\end{align*}This shows that $\{s_{2n-1}\}$ is a decreasing sequence and $\{s_{2n }\}$ is an increasing sequence.
Since 
\[s_{2n }=s_{2n-1}-b_{2n},\]
we find that 
\[s_2\leq s_{2n}\leq s_{2n-1}\leq s_1.\]
Namely, the sequence $\{s_{2n-1}\}$ is bounded below by $s_2$, while the sequence $\{s_{2n}\}$ is bounded above by $s_1$. By the monotone convergence theorem, the limits
\[s_o=\lim_{n\to\infty}s_{2n-1}\hspace{1cm}\text{and}\hspace{1cm}s_e=\lim_{n\to\infty}s_{2n}\] exist.  Since
\[-b_{2n}=a_{2n}=s_{2n}-s_{2n-1},\] 
taking the $n\to \infty$ limits give
\[s_o=s_e.\]\bp
This proves that the sequence $\{s_n\}$ has a limit $s=s_0=s_e$, and thus the alternating series $\di\sum_{n=1}^{\infty}(-1)^{n-1}b_n$ is convergent.

\end{myproof}
Notice that the sum of the alternating series $s=\di\sum_{n=1}^{\infty}(-1)^{n-1}b_n$ is the least upper bound of $\{s_{2n}\}$, and the greatest lower bound of $\{s_{2n-1}\}$.

\begin{remark}{Approximating   the Sum of An Alternating Series}
If $\{b_n\}$ is a sequence that decreases monotonically to 0, the alternating series $\di\sum_{n=1}^{\infty}(-1)^{n-1}b_n$ converges to a sum $s$. If 
\[s_n=\sum_{k=1}^n(-1)^{k-1}b_k\] is the $n^{\text{th}}$ partial sum, then the error in approximating $s$ by $s_n$ is
\[s-s_n=\sum_{k=n+1}^{\infty}(-1)^{k-1}b_k,\] which is also an alternating series. From the proof of Theorem \ref{230227_15}, we obtain a simple estimate
\[|s-s_n|\leq |b_{n+1}|.\]
\end{remark}

\begin{example}[label=ex230227_13]{}
For the alternating   series 
\[\sum_{n=1}^{\infty}\frac{(-1)^{n-1}}{n}=1-\frac{1}{2}+\frac{1}{3}-\frac{1}{4}+\cdots\] in Example \ref{ex230227_12}, $b_n=\di\frac{1}{n}$.  Since $\{b_n\}$ decreases monotonically to 0, by the alternating series test, the series \be \[ \sum_{n=1}^{\infty}\frac{(-1)^{n-1}}{n}=1-\frac{1}{2}+\frac{1}{3}-\frac{1}{4}+\cdots\]  is convergent.   
  Since the harmonic series $\di\sum_{n=1}^{\infty}\frac{1}{n}$ is divergent, the series $\di \sum_{n=1}^{\infty}\frac{(-1)^{n-1}}{n}$ converges conditionally.
\end{example2}

\begin{example}[label=ex230227_14]{}
For any $0<p\leq 1$, the sequence $\{1/n^p\}$ decreases to 0 montonically. Hence, the alternating series $\di \sum_{n=1}^{\infty}\frac{(-1)^{n-1}}{n^p}$ is convergent. Since the series $\di\sum_{n=1}^{\infty}\frac{1}{n^p}$ is divergent,  the series   $\di \sum_{n=1}^{\infty}\frac{(-1)^{n-1}}{n^p}$ converges conditionally.
\end{example}

Now we turn to two useful tests that are used for testing convergence of power series. They both based on comparisons with geometric series. We first prove the following.

\begin{theorem}[label=230227_22]{}
Let $\{a_n\}$ be a sequence of positive numbers. Then
\[\liminf_{n\to\infty}\frac{a_{n+1}}{a_n}\leq\liminf_{n\to\infty}\sqrt[n]{a_n}\leq \limsup_{n\to\infty}\sqrt[n]{a_n}\leq \limsup_{n\to\infty}\frac{a_{n+1}}{a_n}.\]
Hence, if the limit $\di\lim_{n\to\infty}\frac{a_{n+1}}{a_n}$ exists, the limit $\di\lim_{n\to\infty}\sqrt[n]{a_n}$ also exists, and the two limits are equal.
\end{theorem}Since $a_n>0$ for all $n\in\mathbb{Z}^+$, all the four limits in the theorem are nonnegative. \begin{myproof}{Proof}

If $\{c_n\}$ is a sequence of postive numbers,  it is easy to verify that
\[\sup \left\{\frac{1}{c_n}\right\}=\frac{1}{\inf\{c_n\}},\hspace{1cm}\inf \left\{\frac{1}{c_n}\right\}=\frac{1}{\sup\{c_n\}}.\] \bp Therefore, if we prove that
\begin{equation}\label{eq230227_20}\limsup_{n\to\infty}\sqrt[n]{a_n}\leq \limsup_{n\to\infty}\frac{a_{n+1}}{a_n},\end{equation} 
then 
\[\liminf_{n\to\infty}\frac{a_{n+1}}{a_n}\leq\liminf_{n\to\infty}\sqrt[n]{a_n}\] follows by applying \eqref{eq230227_20} to the reciprocal sequence $\{1/a_n\}$. 
From Proposition \ref{230227_11}, we have the inequality $\di \liminf_{n\to\infty}\sqrt[n]{a_n}\leq \limsup_{n\to\infty}\sqrt[n]{a_n}$. 
Hence, we only need to prove  \eqref{eq230227_20}.

 If $\di  \limsup_{n\to\infty}\frac{a_{n+1}}{a_n}=\infty$, there is nothing to prove. Hence, we consider the case
\[u= \limsup_{n\to\infty}\frac{a_{n+1}}{a_n}\] is finite. Given $\varepsilon>0$, there is a positive integer $N$ such that
\[\frac{a_{n+1}}{a_n}<u+ \varepsilon \hspace{1cm}\text{for all}\;n\geq N.\]
By induction, we find that
\[a_n\leq a_N\left(u+ \varepsilon\right)^{n-N}\hspace{1cm}\text{for all}\;n\geq N.\]Let 
$c=a_N(u+\varepsilon)^{-N}$. 
Then 
\[\sqrt[n]{a_n}\leq c^{1/n}(u+\varepsilon)\hspace{1cm}\text{for all}\;n\geq N.\]
This implies that
\[\limsup_{n\to\infty}\sqrt[n]{a_n}\leq \limsup_{n\to\infty}c^{1/n}(u+\varepsilon)=(u+\varepsilon)\lim_{n\to\infty}c^{1/n}=(u+\varepsilon).\]
Since $\varepsilon>0$ is arbitrary, we conclude that
\[\limsup_{n\to\infty}\sqrt[n]{a_n}\leq u= \limsup_{n\to\infty}\frac{a_{n+1}}{a_n}.\]This completes the proof of the theorem.
\end{myproof}

Now we come to the proof of the root test.
\begin{theorem}[label=230227_23]{Root Test}
Given a series $\di\sum_{n=1}^{\infty}a_n$, let
\[\rho=\limsup_{n\to\infty}\sqrt[n]{|a_n|}.\]
\begin{enumerate}[1.]
\item
If $\rho<1$, the series $\di\sum_{n=1}^{\infty}a_n$ converges absolutely.
\item If $\rho>1$,  the series $\di\sum_{n=1}^{\infty}a_n$ is divergent.
\item If $\rho=1$, the test is inconclusive.
\end{enumerate}
\end{theorem}\begin{myproof}{Proof}
If $\rho<1$, take $\di\varepsilon=\frac{1-\rho}{2}$ in (b)(i) of Theorem \ref{230226_5}. There is a positive integer $N$ such that 
\[\sqrt[n]{|a_n|}<\rho+\varepsilon=\rho_1\hspace{1cm}\text{for all}\;n\geq N.\]
Thus, we have
\[|a_n|<\rho_1^n\hspace{1cm}\text{for all}\;n\geq N.\]Notice that
\[\rho_1=\frac{1+\rho}{2}<1.\]Therefore, the geometric series $\di\sum_{n=N}\rho_1^n$ is convergent. By comparison test, the series $\di\sum_{n=1}^{\infty}|a_n|$ is convergent. Thus, the series $\di\sum_{n=1}^{\infty}a_n$ converges absolutely.

If $\rho>1$,  take $\di\varepsilon=\frac{ \rho-1}{2}$ in (b)(ii) of Theorem \ref{230226_5}. There are positive  integers $n_1, n_2, \ldots$ such that $1\leq n_1<n_2<\ldots$ and
\[\sqrt[\leftroot{-3}\uproot{ 5} n_k]{|a_{n_k}|}>\rho-\varepsilon=\rho_2\hspace{1cm}\text{for all}\; k\in\mathbb{Z}^+.\]\bp
Thus, we have
\begin{equation}\label{eq230228_1}|a_{n_k}|>\rho_2^{n_k}\hspace{1cm}\text{for all}\;k\in\mathbb{Z}^+.\end{equation}
Since
\[\rho_2=\frac{1+\rho}{2}>1,\]and $n_k\to\infty$ as $k\to\infty$, we find that 
$\di\lim_{k\to\infty}\rho_2^{n_k}=\infty$. In other words, the sequence $\{\rho_2^{n_k}\}$ is not bounded above. Eq. \eqref{eq230228_1} then implies that $\{|a_{n_k}|\}$ is also not bounded above. Therefore, the limit $\di\lim_{n\to\infty}a_n$ is not zero. Hence, the series $\di\sum_{n=1}^{\infty}a_n$ is divergent.

Now, let us look at some examples where $\rho=1$. First, notice that
\begin{equation}\label{eq230305_7}\lim_{n\to\infty}\sqrt[n]{n}=\lim_{n\to\infty}\exp\left(\frac{\ln n}{n}\right)=\exp\left(\lim_{x\to\infty}\frac{\ln x}{x}\right)=e^0=1.\end{equation}
For the $p$-series $\di\sum_{n=1}^{\infty}\frac{1}{n^p}$, $a_n=n^{-p}$. Thus,
\[\rho=\lim_{n\to\infty}\sqrt[n]{n^{-p}}=\left(\lim_{n\to\infty}\sqrt[n]{n}\right)^{-p}=1.\]
But we have seen that the $p$-series is divergent if $p\leq 1$, and it is convergent when $p>1$. This shows that the root test is conclusive when $\rho=1$.
\end{myproof} 

\begin{example}{}
Determine the convergence of the  series 
$\di\sum_{n=1}^{\infty}\left(\frac{1-n}{2n+1}\right)^n$.
 
\end{example}
\begin{solution}{Solution}
Applying root test, 
\[\rho=\limsup_{n\to\infty}\sqrt[\leftroot{-5}\uproot{15} n]{\left|\left(\frac{1-n}{2n+1}\right)^n\right|}=\lim_{n\to\infty}\frac{n-1}{2n+1}=\frac{1}{2}.\]
Since $\rho<1$, we find that the series   is convergent.
\end{solution}

Finally, we have the ratio test.
\begin{theorem}[label=230305_4]{Ratio Test}
Given a series $\di\sum_{n=1}^{\infty}a_n$ with $a_n\neq 0$ for all $n\in\mathbb{Z}^+$, let
\[r=\liminf_{n\to\infty}\left|\frac{a_{n+1}}{a_n}\right|,\hspace{1cm}R=\limsup_{n\to\infty}\left|\frac{a_{n+1}}{a_n}\right|.\]
\begin{enumerate}[1.]
\item
If $R<1$, the series $\di\sum_{n=1}^{\infty}a_n$ converges absolutely.
\item If $r>1$,  the series $\di\sum_{n=1}^{\infty}a_n$ is divergent.
\item If $r\leq 1\leq R$, the test is inconclusive.
\end{enumerate}
\end{theorem}\begin{myproof}{Proof}
 
If $R<1$, Theorem \ref{230227_22} implies that $\rho=\di\limsup_{n\to\infty}\sqrt[n]{|a_n|}<1$. Theorem \ref{230227_23} implies that $\di\sum_{n=1}^{\infty}a_n$ converges absolutely.
 
If $r>1$, Theorem \ref{230227_22} implies that $\rho=\di\limsup_{n\to\infty}\sqrt[n]{|a_n|}>1$. Theorem \ref{230227_23} implies that $\di\sum_{n=1}^{\infty}a_n$  is divergent.
 
The $p$-series $\di\sum_{n=1}^{\infty}\frac{1}{n^p}$ provides examples of $r=R=1$, but the series is convergent if $p>1$, divergent when $p\leq1 $. Hence, ratio test is also inconclusive when $r\leq 1\leq R$.

 \end{myproof}
Ratio test is useful to determine the convergence of power series. We are going to study this in Chapter 6.

\begin{example}
{}
Determine whether the series is convergent.
\begin{enumerate}[(a)]
\item $\di \sum_{n=1}^{\infty}(-1)^{n-1}\frac{2^n}{n+1}$
\item  $\di \sum_{n=1}^{\infty}(-1)^{n-1}\frac{n+1}{2^n}$

\end{enumerate}
\end{example}
\begin{solution}{Solution}
\begin{enumerate}[(a)]
\item Using  ratio test with $\di a_n=(-1)^{n-1}\frac{2^n}{n+1}$, we find that
\[r=R=\lim_{n\to\infty}\left|\frac{a_{n+1}}{a_n}\right|=\lim_{n\to\infty}\frac{2^{n+1}}{n+2}\times\frac{n+1}{2^n}=2\lim_{n\to\infty}\frac{n+1}{n+2}=2.\]
Therefore, the series $\di \sum_{n=1}^{\infty}(-1)^{n-1}\frac{2^n}{n+1}$ is divergent.
\item Using  ratio test with $\di a_n=(-1)^{n-1}\frac{n+1}{2^n}$, we find that
\[r=R=\lim_{n\to\infty}\left|\frac{a_{n+1}}{a_n}\right|=\lim_{n\to\infty}\frac{n+2}{2^{n+1}}\times\frac{2^n}{n+1}=\frac{1}{2}\lim_{n\to\infty}\frac{n+2}{n+1}=\frac{1}{2}.\]
Therefore, the series $\di \sum_{n=1}^{\infty}(-1)^{n-1}\frac{n+1}{2^n}$ is convergent.

\end{enumerate}
\end{solution}

\begin{highlight}{Convergence Tests}
In this section, we have explored various strategies to determine the convergence of a series $\di\sum_{n=1}^{\infty}a_n$. We make a summary as follows. This is a useful manual for beginners, but it is not binding.

\begin{enumerate}[1.]
\item   Check whether $\di\lim_{n\to\infty}a_n$ is 0. If not, the series $\di\sum_{n=1}^{\infty}a_n$ is divergent.\end{enumerate}\end{highlight}\begin{highlight}{}\begin{enumerate}[1.] 
\item[2.] Check whether it is a geometric series $\di\sum_{n=1}^{\infty}ar^{n-1}$ or a $p$-series $\di\sum_{n=1}^{\infty}\frac{1}{n^p}$. A geometric series $\di\sum_{n=1}^{\infty}ar^{n-1}$ is convergent if and only if $|r|<1$. A $p$-series $\di\sum_{n=1}^{\infty}\frac{1}{n^p}$ is convergent if and only if $p>1$.
\item[3.] If $a_n$ contains powers of $n$ and functions such as $\ln n$, use integral test.
\item[4.] If  $a_n$ involves only expressions of the form $r^n$ for more than one $r$, do limit comparison test to compare with a geometric series.
\item[5.] If  $a_n$ is a rational function of powers of  $n$, do limit comparison test with a $p$-series.
\item[6.]  For alternating series which does not converge absolutely, check whether alternating series test can be applied.
\item[7.] If $a_n=b_n^n$ for each $n\in\mathbb{Z}^+$, and $\di\limsup_{n\to\infty}b_n$ exists, use root test.
\item[8.]  If $a_n$ is a product of a rational function of powers of $n$ and expressions of the form $r^n$, use ratio test.
\end{enumerate}
\end{highlight}

Finally, we want to prove the following useful fact.
\begin{theorem}[label=230306_1]{}
Let $r$ be a real number with $|r|<1$. For any real number $\alpha$, \[\lim_{n\to\infty}n^{\alpha}r^n=0.\]
\end{theorem}\begin{myproof}{Proof}
If $r=0$, the limit is trivial. Hence, we consider the case $|r|<1$ and $r\neq 0$.

If $\alpha\leq 0$, the statement is also easy to prove since $\di\lim_{n\to\infty}r^n=0$ and \bp\[\di\lim_{n\to\infty}n^{\alpha}=\begin{cases}0,\quad&\text{if}\;\alpha<0, \\1, \quad & \text{if}\;\alpha=0.\end{cases}\]

The highly nontrivial case is when $\alpha>0$. In this case, $\di\lim_{n\to\infty}n^{\alpha}=\infty$. Since
\[n^{\alpha}|r|^n=n^{\alpha}e^{n\ln |r|},\]
and $\ln |r|<0$, we can deduce that $\di \lim_{n\to\infty}n^{\alpha}|r|^n=0$ from $\di\lim_{x\to\infty}\frac{x^{\alpha}}{e^x}=0$. Nevertheless, let us present an alternative argument here which is interesting by its own. 

Consider the series $\di\sum_{n=1}^{\infty}n^{\alpha}r^n$ with $a_n=n^{\alpha}r^n$. When $|r|<1$ and $r\neq 0$,
\[\lim_{n\to\infty}\left|\frac{a_{n+1}}{a_n}\right|=\lim_{n\to\infty}|r|\left(\frac{n+1}{n}\right)^{\alpha}=|r|\left(\lim_{n\to\infty}\frac{n+1}{n}\right)^{\alpha}=|r|<1.\]
By ratio test, the  series $\di\sum_{n=1}^{\infty}n^{\alpha}r^n$ is convergent. Therefore,
\[\lim_{n\to\infty}n^{\alpha}r^n=\lim_{n\to\infty}a_n=0.\]
\end{myproof}
\vp
\noindent
{\bf \large Exercises  \thesection}
\setcounter{myquestion}{1}

\begin{question}{\themyquestion}
Determine whether the series $\di\sum_{n=1}^{\infty}\frac{n^2+1}{3n^2+n+1}$ is convergent.
\end{question}
\atc

\begin{question}{\themyquestion}
Let $p$ be a positive number. Show that the series $\di \sum_{n=1}^{\infty}\frac{\ln n}{n^p}$ is convergent if and only if $p>1$.
\end{question}
\atc

\begin{question}{\themyquestion}
Let $p$ be a positive number. Show that the series \[\sum_{n=1}^{\infty}(-1)^{n-1}\frac{\ln n}{n^p}\] is convergent.
\end{question}
\atc

\begin{question}{\themyquestion}
Determine whether the series is convergent.
\begin{enumerate}[(a)]
\item $\di \sum_{n=1}^{\infty}\frac{3^n+(-1)^n4^n}{5^n+2^n}$
\item  $\di \sum_{n=1}^{\infty}\frac{2^n-5^n}{4^n+3^n+1}$
\end{enumerate}
\end{question}
\atc

\begin{question}{\themyquestion}
Determine whether the series $\di\sum_{n=1}^{\infty}(-1)^{n-1}\frac{\sqrt{n}}{n+1}$ is convergent.
\end{question}
\atc
\begin{question}{\themyquestion}
Determine whether the series is convergent.
\begin{enumerate}[(a)]
\item $\di \sum_{n=1}^{\infty}\frac{2n\sqrt{n}+3}{5n^2-2}$
\item  $\di \sum_{n=1}^{\infty}\frac{4n^2-7}{6n^3\sqrt{n}+1}$
\end{enumerate}
\end{question}
\atc
\begin{question}{\themyquestion}
Use Theorem \ref{230227_22} to determine
\[\lim_{n\to\infty}\sqrt[n]{n!}.\]
\end{question}
\atc

\begin{question}{\themyquestion}
Determine whether the series is convergent.
\begin{enumerate}[(a)]
\item $\di \sum_{n=1}^{\infty}\left(\frac{2 \sqrt{n}-1}{\sqrt{n}+1}\right)^n$
\item  $\di \sum_{n=1}^{\infty}\left(\frac{2 \sqrt{n}-1}{3\sqrt{n}+1}\right)^n$
\end{enumerate}
\end{question}
\atc

\begin{question}{\themyquestion}
Determine whether the series is convergent.
\begin{enumerate}[(a)]
\item $\di \sum_{n=1}^{\infty}(-1)^{n-1}\frac{\sqrt{n}2^n}{3^n}$
\item  $\di \sum_{n=1}^{\infty}(-1)^{n-1}\frac{4^n}{3^nn^2}$
\end{enumerate}
\end{question}

\vp

\section{Rearrangement of Series}\label{sec5.3}
In this section, we want to explore more about the difference between a series that converges absolutely and one that converges conditionally.

Given a series $\di\sum_{n=1}^{\infty}a_n$ with terms $\{a_n\}$, define 
\begin{align*}
p_n&=\frac{|a_n|+a_n}{2}=\begin{cases}a_n,\hspace{0.7cm} &\text{if}\;a_n\geq 0,\\0,\quad &\text{if}\;a_n<0;\end{cases}\\q_n&=\frac{|a_n|-a_n}{2}=\begin{cases}-a_n,\quad &\text{if}\;a_n\leq 0,\\0,\quad &\text{if}\;a_n>0.\end{cases}.\end{align*}
Then $0\leq p_n\leq |a_n|$, $0\leq q_n\leq |a_n|$, and
\[|a_n|=p_n+q_n,\hspace{1cm}a_n=p_n-q_n.\]

\begin{example}[label=230228_4]{}
For the series $\di\sum_{n=1}^{\infty}\frac{(-1)^{n-1}}{n}=1-\frac{1}{2}+\frac{1}{3}-\frac{1}{4}+\cdots$,
\[p_{2n-1}=\frac{1}{2n-1},\quad p_{2n}=0;\hspace{1cm}q_{2n-1}=0,\hspace{1cm}q_{2n}=\frac{1}{2n}.\]
\end{example}

\begin{theorem}[label=230228_10]{}Let $\di\sum_{n=1}^{\infty}a_n$ be a convergent series.
\begin{enumerate}[1.]
\item 
If the series $\di\sum_{n=1}^{\infty}a_n$ converges absolutely, then the series $\di\sum_{n=1}^{\infty}p_n$ and the series $\di\sum_{n=1}^{\infty}q_n$ are convergent.
\item 
If the series $\di\sum_{n=1}^{\infty}a_n$ converges conditionally, then the series $\di\sum_{n=1}^{\infty}p_n$ and the series $\di\sum_{n=1}^{\infty}q_n$ are divergent.
 
\end{enumerate}
\end{theorem}
\begin{myproof}{Proof}
First we show that  the two series $\di\sum_{n=1}^{\infty}p_n$ and $\di\sum_{n=1}^{\infty}q_n$ can only be both convergent or both divergent.
We have
\[p_n=a_n+q_n,\hspace{1cm} q_n=p_n-a_n,\]
and we are given that the series $\di\sum_{n=1}^{\infty}a_n$ is convergent. Therefore, the series $\di\sum_{n=1}^{\infty}q_n$ is convergent implies that the series $\di\sum_{n=1}^{\infty}p_n$ is convergent. 
Similarly, the series $\di\sum_{n=1}^{\infty}p_n$ is convergent implies that the series $\di\sum_{n=1}^{\infty}q_n$ is convergent.

If the series  $\di\sum_{n=1}^{\infty}a_n$ converges absolutely, the series  $\di\sum_{n=1}^{\infty}|a_n|$ is convergent.
Since
\[0\leq p_n\leq|a_n|,\hspace{1cm}0\leq q_n\leq |a_n|,\]
comparison test implies that the series $\di\sum_{n=1}^{\infty}p_n$ and $\di\sum_{n=1}^{\infty}q_n$ are convergent.

Conversely, if  the series $\di\sum_{n=1}^{\infty}p_n$ and $\di\sum_{n=1}^{\infty}q_n$ are convergent, since
\[|a_n|=p_n+q_n,\]
the series $\di\sum_{n=1}^{\infty}|a_n|$ must be convergent. Therefore, if the series $\di\sum_{n=1}^{\infty}a_n$ converges conditionally, which means the series $\di\sum_{n=1}^{\infty}|a_n|$ is divergent, then the series $\di\sum_{n=1}^{\infty}p_n$ and the series $\di\sum_{n=1}^{\infty}q_n$ must be both divergent.
\end{myproof}
 \begin{example} {}For the series $\di\sum_{n=1}^{\infty}\frac{(-1)^{n-1}}{n}=1-\frac{1}{2}+\frac{1}{3}-\frac{1}{4}+\cdots$ in Example \ref{230228_4},
the series $\di\sum_{n=1}^{\infty}p_{n}=\sum_{n=1}^{\infty}\frac{1}{2n-1}$ and the series  $\di\sum_{n=1}^{\infty}q_{n}=\sum_{n=1}^{\infty}\frac{1}{2n}$ are divergent. 
\end{example}

\begin{definition}{Rearrangement of a Series}
A rearrangement of a series $\di\sum_{n=1}^{\infty}a_n$ is the series $\di\sum_{n=1}^{\infty}a_{\pi(n)}$, where $\pi:\mathbb{Z}^+\to\mathbb{Z}^+$ is a bijective correspondence.
\end{definition}

\begin{example}[label=ex230228_6]{}
Let $\pi:\mathbb{Z}^+\to\mathbb{Z}$ be the bijective correspondence
\[\pi(1)=1,\;\pi(2)=3, \;\pi(3)=2,\;\pi(4)=5,\;\pi(5)=7,\;\pi(6)=4, \;\ldots.\]Namely,
\[\pi(n)=\begin{cases}4k-3,\quad &\text{if}\;n=3k-2,\\
4k-1,\quad &\text{if}\;n=3k-1,\\
2k,\quad &\text{if}\;n=3k.
\end{cases}\]
The rearrangement of the series $\di\sum_{n=1}^{\infty}\frac{(-1)^{n-1}}{n}=1-\frac{1}{2}+\frac{1}{3}-\frac{1}{4}+\cdots$
induced by $\pi$ is
\[1+\frac{1}{3}-\frac{1}{2}+\frac{1}{5}+\frac{1}{7}-\frac{1}{4}+\cdots.\]
\end{example}
The main thing we want to discuss in this section is whether rearrangment will affect the convergence of a series. Consider the rearrangment discussed in Example \ref{ex230228_6}, we know that original series\[\sum_{n=1}^{\infty}a_n=\sum_{n=1}^{\infty}\frac{(-1)^{n-1}}{n}=1-\frac{1}{2}+\frac{1}{3}-\frac{1}{4}+\cdots\] is convergent. 
We can find its sum in the following way. Let $s_n=a_1+a_2+\cdots+a_n$ be its $n^{\text{th}}$ partial sum. 
Then
\begin{align*}
s_{2n}&=1-\frac{1}{2}+\frac{1}{3}-\frac{1}{4}+\cdots+\frac{1}{2n-1}-\frac{1}{2n}
\\&=\left(1+\frac{1}{2}+\frac{1}{3}+\frac{1}{4}+\cdots+\frac{1}{2n-1}+\frac{1}{2n}\right)-2 \left( \frac{1}{2}+\frac{1}{4}+\cdots+\frac{1}{2n}\right)
\\&=\left(1+\frac{1}{2}+\frac{1}{3}+\frac{1}{4}+\cdots+\frac{1}{2n-1}+\frac{1}{2n}\right)-  \left( 1+\frac{1}{2}+\cdots+\frac{1}{n}\right)\\
&=\sum_{k=1}^{2n}\frac{1}{k}-\sum_{k=1}^{n}\frac{1}{k}.
\end{align*}Let
\[c_n=1+\frac{1}{2}+\ldots+\frac{1}{n}-\ln n=\sum_{k=1}^{n}\frac{1}{k}-\ln n.\]
By Example \ref{ex230228_7}, $\di\lim_{n\to\infty}c_n=\gamma$ is the Euler's constant.
We can write $s_{2n}$ as 
\[s_{2n}=c_{2n}+\ln(2n)-\left(c_n+\ln n\right)=c_{2n}-c_n+\ln 2.\]
Then we find that
\[\lim_{n\to\infty} s_{2n}=\lim_{n\to\infty}(c_{2n}-c_n+\ln 2)=\gamma-\gamma+\ln 2=\ln 2.\]
This shows that  \[\sum_{n=1}^{\infty}a_n=\sum_{n=1}^{\infty}\frac{(-1)^{n-1}}{n}=\ln 2.\]
For the rearranged series $\di\sum_{n=1}^{\infty}b_n=\di\sum_{n=1}^{\infty}a_{\pi(n)}$,
\[b_{3k-2}=\frac{1}{4k-3},\quad b_{3k-1}=\frac{1}{4k-1},\quad b_{3k}=-\frac{1}{2k}\hspace{1cm}\text{for all}\;k\in\mathbb{Z}^+.\]
Let $t_n=b_1+b_2+\cdots+b_n$ be the $n^{\text{th}}$ partial sum of the series $\di\sum_{n=1}^{\infty}b_n$. Now \[t_{3n}=\sum_{k=1}^n\left(\frac{1}{4k-3}+\frac{1}{4k-1}\right) - \sum_{k=1}^n\frac{1}{2k}.\]As $k$ runs from 1 to $n$, $4k-3$ and $4k-1$ run through all positive odd integers between 1 and $4n$. Therefore,
\[t_{3n}=\sum_{k=1}^{2n}\frac{1}{2k-1}- \sum_{k=1}^n\frac{1}{2k}= \sum_{k=1}^{4n}\frac{1}{k}-\sum_{k=1}^{2n}\frac{1}{2k}- \sum_{k=1}^n\frac{1}{2k}.\]
Using $c_{n}$, we can rewrite this as
\[t_{3n}=c_{4n}+\ln(4n)-\frac{1}{2}\left(c_{2n}+\ln(2n)\right)-\frac{1}{2}\left(c_n+\ln n\right)=c_{4n}-\frac{1}{2}c_{2n}-\frac{1}{2}c_n+\frac{3}{2}\ln 2.\]This allows us to conclude that
\[\lim_{n\to\infty}t_{3n}=\frac{3}{2}\ln 2.\]Since
\[t_{3n+1}=t_{3n}+b_{3n+1},\hspace{1cm}t_{3n+2}=t_{3n+1}+b_{3n+1}+b_{3n+2},\]and $\di\lim_{n\to\infty}b_n=0$, we find that 
\[\lim_{n\to\infty}t_{3n+1}=\lim_{n\to\infty}t_{3n+2}=\lim_{n\to\infty}t_{3n}.\]
This proves that the series $\di\sum_{n=1}^{\infty} b_n$ is  convergent, and it converges to $\di\lim_{n\to\infty}t_{3n}=\frac{3}{2}\ln 2$.
 
Hence, although the series $\di\sum_{n=1}^{\infty}b_n$ is a rearrangement  of the series $\di\sum_{n=1}^{\infty}a_n$, it has a different sum.

In the following, we prove that rearrangement of a nonnegative series would not lead to different sums.
\begin{lemma}[label=230228_9]{}
If $a_n\geq 0$ for all $n\in\mathbb{Z}^+$ and the series $\di\sum_{n=1}^{\infty}a_n$ is convergent, then any rearrangement of the series has the same sum. Namely, for any bijecion $\pi:\mathbb{Z}^+\to\mathbb{Z}^+$, 
\[\sum_{n=1}^{\infty}a_{\pi(n)}=\sum_{n=1}^{\infty}a_n.\]
\end{lemma}\begin{myproof}{Proof}
In a nutshell, this is just the fact that a nonnegative series is convergent if and only if the sequence of partial sums is bounded above, and the sum of the series is the least upper bound of the sequence of partial sums. 

For a rigorous argument, define $s_n=a_1+\cdots+a_n$ to be the $n^{\text{th}}$ partial sum of the series  $\di\sum_{n=1}^{\infty}a_n$, and $t_n=a_{\pi(1)}+\cdots+a_{\pi(n)}$ to be the 
 $n^{\text{th}}$ partial sum of the series  $\di\sum_{n=1}^{\infty}a_{\pi(n)}$. Notice that both $\{s_n\}$ and $\{t_n\}$ are increasing sequences. We are given that $s=\sup\{s_n\}$ exists. For any positive integer $n$, the set $\{\pi(1), \pi(2), \ldots, \pi(n)\}$ has a maximum $N_n$. This means that the set $\{\pi(1), \pi(2), \ldots, \pi(n)\}$ is contained in the set $\{1, 2, \ldots, N_n\}$. Therefore,
\[t_{n}\leq s_{N_n}\leq s.\]

This shows that the increasing sequence $\{t_n\}$ is bounded above by $s$. Hence, $\di t=\lim_{n\to \infty}t_n=\sup\{t_n\}$ exists and $t\leq s$. For the  opposite inequality, observe that 
$\di\sum_{n=1}^{\infty}a_{n}$ is a rearrangement of $\di\sum_{n=1}^{\infty}a_{\pi(n)}$ induced by the bijection $\pi^{-1}:\mathbb{Z}^+\to\mathbb{Z}^+$. Hence, the same argument above shows that $s\leq t$. Combine together, we have $t=s$, thus proving that any rearrangement of the series $\di \sum_{n=1}^{\infty}a_{n}$ has the same sum.

\end{myproof}

Now we can prove that any rearrangemnt of an absolutely convergent series converge to the same sum.
\begin{theorem}{Rearrangement of Absolutely Convergent Series}
If  the series $\di\sum_{n=1}^{\infty}a_n$   converges absolutely, then any rearrangement of the series has the same sum. Namely, for any bijecion $\pi:\mathbb{Z}^+\to\mathbb{Z}^+$, 
\[\sum_{n=1}^{\infty}a_{\pi(n)}=\sum_{n=1}^{\infty}a_n.\]
\end{theorem}\begin{myproof}{Proof}
Define the nonnegative series $\di\sum_{n=1}^{\infty}p_n$ and $\di\sum_{n=1}^{\infty}q_n$ by 
\[p_n=\frac{|a_n|+a_n}{2},\hspace{1cm}q_n=\frac{|a_n|-a_n}{2}.\]
Then
\[a_n=p_n-q_n.\]
Since  the series $\di\sum_{n=1}^{\infty}a_n$   converges absolutely, Theorem \ref{230228_10} says that the series  $\di\sum_{n=1}^{\infty}p_n$ and $\di\sum_{n=1}^{\infty}q_n$  are convergent.

 Lemma \ref{230228_9} says that  for any bijecion $\pi:\mathbb{Z}^+\to\mathbb{Z}^+$, the  series $\di\sum_{n=1}^{\infty}p_{\pi(n)}$ and $\di\sum_{n=1}^{\infty}q_{\pi(n)}$ are convergent, and
\[\sum_{n=1}^{\infty}p_{\pi(n)}=\sum_{n=1}^{\infty}p_n,\hspace{1cm} \sum_{n=1}^{\infty}q_{\pi(n)}=\sum_{n=1}^{\infty}q_n.\]Therefore, the series $\di\sum_{n=1}^{\infty}a_{\pi(n)}$ is convergent and
\[\sum_{n=1}^{\infty}a_{\pi(n)}=\sum_{n=1}^{\infty}p_{\pi(n)}-\sum_{n=1}^{\infty}q_{\pi(n)}=\sum_{n=1}^{\infty}p_{n}-\sum_{n=1}^{\infty}q_{n}=\sum_{n=1}^{\infty}a_n.\]
\end{myproof}

Finally, we come to the celebrated Riemann's theorem for series that converges conditionally.
\begin{theorem}{\\Riemann's Theorem for Conditionally Convergent Series}
Let $\di\sum_{n=1}^{\infty}a_n$ be a series that converges conditionally, and let $b$ and $c$ be two extended real numbers with $b\leq c$. There exists a bijection $\pi:\mathbb{Z}^+\to\mathbb{Z}^+$ such that for the series $\di\sum_{n=1}^{\infty}a_{\pi(n)}$ with partial sums $\di t_n=a_{\pi(1)}+\cdots+a_{\pi(n)}$,
\[\liminf_{n\to\infty}t_n=b,\hspace{1cm}\limsup_{n\to\infty}t_n=c.\]\end{theorem} 
Here an extended real number is either an ordinary real number or $\pm\infty$. This theorem implies that one can have a rearrangement of a conditionally convergent series that diverge to $\pm\infty$ or converge to any real number.
\begin{myproof}{Proof}For $n\in\mathbb{Z}^+$, let
\[p_n=\frac{|a_n|+a_n}{2},\hspace{1cm}q_n=\frac{|a_n|-a_n}{2}.\]
Since  the series $\di\sum_{n=1}^{\infty}a_n$   converges conditionally, Theorem \ref{230228_10} says that the series  $\di\sum_{n=1}^{\infty}p_n$ and $\di\sum_{n=1}^{\infty}q_n$  are divergent. 
Let
\[S_+=\left\{n\in\mathbb{Z}^+\,|\, a_n\geq 0\right\}, \hspace{1cm}S_-=\left\{n\in\mathbb{Z}^+\,|\, a_n< 0\right\}.\] 
Then \[S_+\cup S_-=\mathbb{Z}^+,\hspace{1cm}S_+\cap S_-=\emptyset.\]
There are strictly increasing maps $\pi_1:\mathbb{Z}^+\to\mathbb{Z}^+$ and $\pi_2:\mathbb{Z}^+\to\mathbb{Z}^+$, such that $\pi_1(\mathbb{Z}^+)=S_+$ and $\pi_2(\mathbb{Z}^+)=S_-$.

Define the nonnegative series $\di\sum_{n=1}^{\infty}u_n$ and $\di\sum_{n=1}^{\infty}v_n$ by 
\[u_n=a_{\pi_1(n)},\hspace{1cm}v_n=-a_{\pi_2(n)}.\]\bp
Then  the sequences  $\{u_n\}$ and $\{v_n\}$ are obtained from the sequences $\{p_n\}$ and $\{q_n\}$ by removing some zero terms. 
 Hence, both nonnegative series $\di\sum_{n=1}^{\infty}u_n$ and $\di\sum_{n=1}^{\infty}v_n$ are divergent.

Now we start to define the bijection $\pi:\mathbb{Z}^+\to\mathbb{Z}^+$.  
 Construct two sequences of real numbers $\{b_n\}$ and $\{c_n\}$ such that  $c_1>0$, $b_n\leq c_n$ for all $n\in\mathbb{Z}^+$, and 
\[\lim_{n\to \infty}b_n=b,\hspace{1cm}\lim_{n\to \infty}c_n=c.\]
Take $k_1$ to be the smallest positive integer such that
\[C_1=u_1+u_2+\cdots+u_{k_1}>c_1.\]
Then define
\[\pi(1)=\pi_1(1),\ldots,\pi(k_1)=\pi_1(k_1).\]

Take $l_1$ to be the smallest positive integer such that
\[B_1=C_1-(v_1+v_2+\cdots+v_{l_1})<b_1.\]
Then define
\[\pi(k_1+1)=\pi_2(1),\ldots,\pi(k_1+l_1)=\pi_2(l_1).\]
Take $k_2$ to be the smallest positive integer such that
\[C_2=B_1+u_{k_1+1}+\cdots+u_{k_1+k_2}>c_2.\]
Then define
\[\pi(k_1+l_1+1)=\pi_1(k_1+1),\ldots,\pi(k_1+l_1+k_2)=\pi_1(k_1+k_2).\]
Take $l_2$ to be the smallest positive integer such that
\[B_2=C_2-\left(v_{l_1+1}+v_{l_1+2}+\cdots+v_{l_1+l_2}\right)<b_2.\] 
Then define
\[\pi(k_1+l_1+k_2+1)=\pi_2(l_1+1),\ldots,\pi(k_1+l_1+k_2+l_2)=\pi_2(l_1+l_2).\]\bp
Continue this construction inductively. Since $\di\sum_{n=1}^{\infty}u_n$ and $\di\sum_{n=1}^{\infty}v_n$ are nonnegative sequences that diverges to $\infty$, and $b_m\leq c_m$ for all positive integers $m$, the existence of the positive integers $k_m$ and $l_m$ at each step is guaranteed. It is easy to see that the map $\pi:\mathbb{Z}^+\to\mathbb{Z}^+$ is a bijection. For the series $\di\sum_{n=1}^{\infty}a_{\pi(n)}$, let $t_n=a_{\pi(1)}+\cdots+a_{\pi(n)}$ be its $n^{\text{th}}$ partial sum. Set  $\alpha_0=\beta_0=0$, and for $m\geq 1$, let
\begin{align*}
\alpha_m&= k_1  +k_{2} +\cdots+k_m,\hspace{1cm}\beta_m =  l_1 +l_{2} +\cdots +l_m,\\
\delta_m&=\alpha_{m-1}+\beta_{m-1}+k_m,\hspace{1 cm}\lambda_m=\delta_m+l_m=\alpha_m+\beta_m.\end{align*} 
Then
\begin{gather*}
1\leq \alpha_1<\alpha_2<\cdots<\alpha_m<\cdots,\\
1\leq\beta_1<\beta_2<\cdots<\beta_m<\cdots,\\
1\leq\delta_1<\lambda_1<\delta_2<\lambda_2<\cdots<\delta_m<\lambda_m<\cdots.\end{gather*}
By construction,
\begin{gather*}
t_1\leq t_2\leq\cdots\leq t_{\delta_1-1}\leq c_1<t_{\delta_1}\leq c_1+u_{\alpha_1},\\
t_{\delta_1}\geq t_{\delta_1+1}\geq t_{\delta_1+2}\geq\cdots\geq t_{\lambda_1-1}\geq b_1>
t_{\lambda_1}\geq b_1-v_{\beta_1},\\
t_{\lambda_1}\leq t_{\lambda_1+1}\leq t_{\lambda_1+2}\leq \cdots\leq t_{\delta_2-1}\leq c_2<t_{\delta_2}\leq c_2+u_{\alpha_2},\\
t_{\delta_2}\geq t_{\delta_2+1}\geq t_{\delta_2+2}\geq\cdots\geq t_{\lambda_2-1}\geq b_2>
t_{\lambda_2}\geq  b_2-v_{\beta_2},\\
 \vdots
\end{gather*}
Since the series $\di\sum_{n=1}^{\infty}a_n$ is convergent, $\di\lim_{n\to\infty}a_n=0$. This implies that the sequences $\{u_n\}$ and $\{v_n\}$ converge  to $0$.  Therefore, \[\di\lim_{m\to\infty}u_{\alpha_m}=\lim_{m\to\infty}v_{\beta_m}=0.\]\bp
   Given $\varepsilon>0$,   there exists a positive integer $M_1$ so that for all $m\geq M_1$,
\[0\leq u_{\alpha_m}<\frac{\varepsilon}{2},\hspace{1cm} 0\leq v_{\beta_m}<\frac{\varepsilon}{2}.\] There exists a positive integer $M_2$ so that $M_2\geq M_1$ and  for all $m\geq M_2$,
\[b_m>b-\frac{\varepsilon}{2},\hspace{1cm}c_m<c+\frac{\varepsilon}{2}.\]
  Let $N=\max\{ \alpha_{M_1}, \beta_{M_1}, \lambda_{M_2}\}$. If  $n\geq N$,  then $n\geq \lambda_{M_2}>\delta_{M_2}$. Hence, there exists $m\geq M_2$ such that
\[\delta_m\leq n<\delta_{m+1}.\]Then
\[t_n\leq \max\{c_m+u_{\alpha_m}, c_{m+1}+u_{\alpha_{m+1}}\}.\]

Since $m\geq M_2$, $c_m$ and $c_{m+1}$ are less than $c+\varepsilon/2$. Since $m\geq M_2\geq M_1$, $u_{\alpha_m}$ and $u_{\alpha_{m+1}}$ are less than $\varepsilon/2$. These imply that for all $n\geq N$,
\[t_n<c+\varepsilon.\]
Hence,
\[\limsup_{n\to\infty} t_n\leq c.\]
Similarly, we can show that \[\liminf_{n\to\infty}t_n\geq b.\]
For all $m\in\mathbb{Z}^+$,
\[b_m-v_{\beta_m}\leq t_{\lambda_m}<b_m,\hspace{1cm}c_m<t_{\delta_m}\leq c_m+u_{\alpha_m}.\]
Taking $m\to\infty$ limits show that $\{t_{\lambda_m}\}$ is a subsequence of $\{t_n\}$ that converges to $b$, and $\{t_{\delta_m}\}$ is a subseqeunce of $\{t_n\}$ that converges to $c$. This completes the proof that \[\liminf_{n\to\infty}t_n=b,\hspace{1cm}\limsup_{n\to\infty}t_n=c.\]
\end{myproof}
\vp
\noindent
{\bf \large Exercises  \thesection}
\setcounter{myquestion}{1}

 \begin{question}{\themyquestion}Show that the series \[\di\sum_{n=1}^{\infty}a_n=\sum_{n=1}^{\infty}\frac{(-1)^{n-1}}{\sqrt{n+1}}\] is convergent. 
 If $\pi:\mathbb{Z}^+\to\mathbb{Z}^+$ is a bijective correspondence,   consider the rearrangement of the series $\di\sum_{n=1}^{\infty}a_{n}$ given by $\di\sum_{n=1}^{\infty}a_{\pi(n)}$. Does the series  $\di\sum_{n=1}^{\infty}a_{\pi(n)}$ necessarily converge to the same number as the series $\di\sum_{n=1}^{\infty}a_{n}$? 
  
 \end{question}
 \atc
 \begin{question}{\themyquestion}Show that the  series \[\sum_{n=1}^{\infty}a_n=\sum_{n=1}^{\infty}\frac{(-1)^{n-1}\sqrt{n}}{n^2+1}\] is convergent. 
 If $\pi:\mathbb{Z}^+\to\mathbb{Z}^+$ is a bijective correspondence,   consider the rearrangement of the series $\di\sum_{n=1}^{\infty}a_{n}$ given by $\di\sum_{n=1}^{\infty}a_{\pi(n)}$. Does the series  $\di\sum_{n=1}^{\infty}a_{\pi(n)}$ necessarily converge to the same number as the series $\di\sum_{n=1}^{\infty}a_{n}$? 
  
 \end{question}
 \atc
 \begin{question}{\themyquestion}Show that the  series \[\sum_{n=1}^{\infty}a_n=\sum_{n=1}^{\infty}\frac{(-1)^{n-1}\sqrt{n}}{n+1}\] is convergent. 
 If $\pi:\mathbb{Z}^+\to\mathbb{Z}^+$ is a bijective correspondence,   consider the rearrangement of the series $\di\sum_{n=1}^{\infty}a_{n}$ given by $\di\sum_{n=1}^{\infty}a_{\pi(n)}$. Does the series  $\di\sum_{n=1}^{\infty}a_{\pi(n)}$ necessarily converge to the same number as the series $\di\sum_{n=1}^{\infty}a_{n}$? 
  
 \end{question}

\vp

\section{Infinite Products}\label{sec5.4}

In this section, we consider infinite products and study its convergence. An infinite product is a product of the form
\[\prod_{n=1}^{\infty}u_n=u_1u_2\cdots u_n\cdots,\]
where $\{u_n\}$ is an infinite sequence. The definition of convergence of infinite product is slightly more complicated.

\begin{definition}{Convergence of Infinite Product} 
Given a sequence $\{u_n\}$, consider the infinite product
$\di \prod_{n=1}^{\infty}u_n$.
\begin{enumerate}[(a)]
\item
If infinitely many of the terms $u_n$'s are zero, then we say that the infinite product $\di \prod_{n=1}^{\infty}u_n$ is divergent.
\item If only finitely many of the $u_n$'s are zero, there is a positive integer $\ell$ such that $u_n$ is nonzero for all $n\geq \ell$. Form the partial product
\[P[\ell]_{n}=\prod_{k=\ell}^nu_k,\hspace{1cm}\text{for}\;n\geq \ell.\]

\begin{enumerate}[(i)]
\item If the limit $\di \lim_{n\rightarrow \infty}P[\ell]_{n}$ does not exist or the limit is 0, we say that the infinite product $\di \prod_{n=1}^{\infty}u_n$ is divergent.
\item If the limit  $\di \lim_{n\rightarrow \infty}P[\ell]_{n}$ exists and is equal to a nonzero number $P[\ell]$, we say that  the infinite product $\di \prod_{n=1}^{\infty}u_n$ converges to
\[P=P[\ell]\prod_{k=1}^{\ell-1}u_k.\]

\end{enumerate}
\end{enumerate}
\end{definition}
The convergence of infinite product is not affected by finitely many terms in the product. If $u_n\neq 0$ for all $n\geq 1$, we will denote the partial product $\di P[1]_n=\prod_{k=1}^n u_k$ simply as $P_n$. 

\begin{highlight}{} 
By definition, if the infinite product $\di\prod_{n=1}^{\infty}u_n$ converges to 0, then at least one of the $u_n$ is equal to 0, and  there are only finitely many of the $u_n$'s that are equal to 0.
\end{highlight}
Let us look at a few examples.

\begin{example}[label=ex230301_1]{}
Determine the convergence of the infinite product $\di\prod_{n=1}^{\infty}\left(1+\frac{1}{n}\right)$.
\end{example}
\begin{solution}{Solution}
For $n\geq 1$, $\di u_n=1+\frac{1}{n}\neq 0$. Notice that
\[P_n=\prod_{k=1}^n\left(1+\frac{1}{k}\right)=\frac{2}{1}\times\frac{3}{2}\times\cdots\times\frac{n+1}{n}=n+1.\]
Since $\di\lim_{n\to\infty}P_n=\infty$, the infinite product $\di\prod_{n=1}^{\infty}\left(1+\frac{1}{n}\right)$ is divergent.

\end{solution}

\begin{example}[label=ex230301_2]{}
Determine the convergence of the infinite product $\di\prod_{n=1}^{\infty}\left(1-\frac{1}{n}\right)$.
\end{example}
\begin{solution}{Solution}
 For $n\geq 1$, $\di u_n=1-\frac{1}{n}$. We find that $u_1=0$ and $u_n>0$ for all $n\geq 2$.
\[P[2]_n=\prod_{k=2}^n\left(1-\frac{1}{k}\right)=\frac{1}{2}\times\frac{2}{3}\times\cdots\frac{n-1}{n}=\frac{1}{n}.\]
Since $\di\lim_{n\to\infty}P[2]_n=0$, the infinite product $\di\prod_{n=1}^{\infty}\left(1-\frac{1}{n}\right)$ is divergent.
\end{solution}

\begin{example}[label=ex230301_3]{}
Determine the convergence of the infinite product $\di\prod_{n=1}^{\infty}\left(1-\frac{1}{n^2}\right)$.
\end{example}

\begin{solution}{Solution}
 For $n\geq 1$, $\di u_n=1-\frac{1}{n^2}$. We find that $u_1=0$ and $u_n>0$ for all $n\geq 2$.
\begin{align*}
P[2]_n&=\prod_{k=2}^n\left(1-\frac{1}{k^2}\right)=\prod_{k=2}^n\left(1-\frac{1}{k}\right)\prod_{k=2}^n\left(1+\frac{1}{k}\right)\\
&=\frac{1}{2}\times\frac{2}{3}\times\cdots\frac{n-1}{n}\times\frac{3}{2}\times\cdots\times\frac{n+1}{n}=\frac{n+1}{2n}.
\end{align*}Since $\di P[2]=\lim_{n\to\infty}P[2]_n=\frac{1}{2}$, the infinite product $\di\prod_{n=1}^{\infty}\left(1-\frac{1}{n^2}\right)$ is convergent, and it converges to $\di u_1P[2]=0$.
 
\end{solution}

\begin{highlight}{}When the sequence of  partial products $\{P_{n}\}$ converges to 0, we consider the infinite product as divergent. This is so that the infinite product $\di\prod_{n=1}^{\infty} u_n$ is convergent if and only if the infinite product $\di \prod_{n=1}^{\infty} u_n^{-1}$ is convergent. \end{highlight}

The following is obvious. 
\begin{proposition}{} If the infinite product $\di\prod_{n=1}^{\infty}u_n$ is convergent, then 
$\di\lim_{n\rightarrow\infty}u_n=1$.
\end{proposition}
Using this proposition, when we consider convergence of the infinite product $\di\prod_{n=1}^{\infty}u_n$, we can assume that $u_n>0$ for all $n\in\mathbb{Z}^+$.

There is a Cauchy criterion for convergence of infinite product.

\begin{theorem}{Cauchy Criterion for Infinite Product} 
Let $\{u_n\}$ be a sequence of positive  numbers. The infinite profuct $\di\prod_{n=1}^{\infty}u_n$ is convergent if and only if it satifies the Cauchy criterion, which says that for every $\varepsilon>0$, there exists a positive integer $N$ such that for all $m\geq n\geq N$,
\[\left|\left[\prod_{k=n}^m u_k\right]-1\right|<\varepsilon.\]
\end{theorem}
The proof of this is more complicated than its infinite series counterpart.
\begin{myproof}{Proof}
 Let
$\di P_n= \prod_{k=1}^{n} u_k$ be the $n^{\text{th}}$ partial product. Then $P_n>0$ for all $n\in\mathbb{Z}^+$.
If the infinite product $\di\prod_{n=1}^{\infty}u_n$ is convergent, then the sequence $\{P_n\}$ converges to a positive number $P$. This implies that there is a positive integer $N_1$ such that 
\[P_n>\frac{P}{2}\quad  \text{for all}\; n\geq N_1.\] Given $\varepsilon>0$, apply Cauchy criterion to the convergent sequence $\{P_n\}$, we find that there is a positive integer $N_2$ such that for all $m\geq n\geq N_2$,
\[|P_n-P_{m}|<\frac{P\varepsilon}{2}.\] \bp
Let $N=\max\{N_1, N_2\}+1$.
We find that for all $ m\geq n\geq N$,
\begin{align*}
\left|\left[\prod_{k=n}^m u_k\right]-1\right|&=\left|\frac{P_m}{P_{n-1}}-1\right|\\
&=\frac{1}{P_{n-1}}\times\left|P_m-P_{n-1}\right|\\
&<\frac{2}{P}\times\frac{P\varepsilon}{2}=\varepsilon.
\end{align*}
Therefore, the Cauchy criterion for infinite product is satisfied. 

Conversely, assume the Cauchy criterion for infinite product holds. Taking $\varepsilon=1/2$, we find that there is an integer $N_1$ such that for all $m\geq n\geq N_1$,
\[\left|\left[\prod_{k=n}^m u_k\right]-1\right|< \frac{1}{2}.\]This implies that  
\begin{equation}\label{eq220928_3}\frac{1}{2}< \frac{P_n}{P_{N_1}} <\frac{3}{2}\quad\text{for all}\; n\geq N_1. \end{equation} 
Now given $\varepsilon>0$, there is an integer $N\geq N_1$ such that for all $m\geq n\geq N$,
\[\left|\left[\prod_{k=n+1}^m u_k\right]-1\right|< \frac{ 2\varepsilon}{3 P_{N_1} }.\]
This implies that when $m\geq n\geq N$,
\begin{align*}
|P_m-P_{n}|&=P_n\times \left|\left[\prod_{k=n+1}^m u_k\right]-1\right|<\frac{3P_{N_1}}{ 2 }\times  \frac{2 \varepsilon}{3 P_{N_1}}=\varepsilon.
\end{align*} Hence, $\{P_n\}$ is a Cauchy sequence, and thus it is convergent. Eq. \eqref{eq220928_3} then implies that \[\lim_{n\rightarrow \infty}P_n\geq\frac{1}{2}P_{N_1}>0.\]This proves that $\{P_n\}$ does not converge to 0. Therefore, the infinite product $\di\prod_{n=1}^{\infty}u_n$ is convergent.

\end{myproof}

The following  gives a relation between the convergence of the infinite product with the convergence of infinite series. 
\begin{theorem}[label=thm220929_1]{}
Let $\{u_n\}$ be a sequence of positive numbers.  Then the infinite product $\di\prod_{n=1}^{\infty}u_n$ is convergent if and only if the infinite series $\di\sum_{n=1}^{\infty}\ln u_n$ is convergent. 
\end{theorem}
\begin{myproof}{Proof}
First assume that the infinite product $\di\prod_{k=1}^{\infty}u_k$ is convergent.  
Given $\varepsilon>0$, since
$\di \lim_{x\rightarrow 1}\ln x=0$, there exists a $\delta$ such that
$0<\delta<1$ and if $|x-1|<\delta$, then
$\di |\ln x|<\varepsilon$.
By the Cauchy criterion for infinite products,   there is a
 positive integer $N $ such that for all $m\geq n\geq N$, 
\[\left|\left[\prod_{k=n}^m u_k\right]-1\right|<\delta.\] 
It follows that
  for all $m\geq n\geq N$, 
 \[\left| \sum_{k=n}^m \ln u_k \right|=\left|\ln\left[\prod_{k=m}^n u_k\right]\right|<\varepsilon.\]
This proves that the  infinite series $\di\sum_{n=1}^{\infty}\ln u_n$ satisfies the Cauchy criterion. Hence, it is convergent.

Conversely, assume that the infinite series $\di\sum_{n=1}^{\infty}\ln u_n$ is convergent. 
  Given $\varepsilon>0$, since
$\di \lim_{x\rightarrow 0}e^x=1$, there exists $\delta>0$ such that if $|x|<\delta$, then 
\[|e^x-1|<\varepsilon.\]\bp
 Using Cauchy criterion for infinite series, we find that there is a positive integer $N$ such that for all $m\geq n\geq N$,
\[\left|\sum_{k=n}^m\ln u_k\right|<\delta.\] 
It follows that for all $m\geq n\geq N$,
\[\left|\left[\prod_{k=n}^m u_k\right]-1\right|=\left|\exp\left(\sum_{k=n}^m\ln u_k\right)-1\right|<\varepsilon.\]
 This shows that the infinite product $\di\prod_{n=1}^{\infty}u_n$ satisfies the Cauchy criterion. Hence, it is convergent.
\end{myproof}

\begin{example}{}
For any nonzero real number $a$, the infinite product $\di \prod_{n=1}^{\infty} \exp\left(\frac{a}{n}\right)$ is divergent since the infinite series $\di\sum_{n=1}^{\infty}\frac{a}{n}$ is divergent; while the infinite product $\di \prod_{n=1}^{\infty} \exp\left(\frac{a}{n^2}\right)$ is convergent since the infinite series $\di\sum_{n=1}^{\infty}\frac{a}{n^2}$ is convergent. 
\end{example}

Since 
$\di \lim_{a\rightarrow 0}\frac{\ln(1+a)}{ a}=1$,
 it is natural to compare the convergence of the product $\di\prod_{n=1}^{\infty}(1+a_n)$ to the convergence of the series $\di\sum_{n=1}^{\infty}a_n$.  
\begin{theorem}[label=thm220929_2]{}
Let $\{a_n\}$ be a sequence of   real numbers such that $0<a_n<1$ for all $n\in\mathbb{Z}^+$. Then the following three statements are equivalent.
\begin{enumerate}[(a)]
\item  The series $\di \sum_{n=1}^{\infty}a_n$ is convergent.
\item The infinite product $\di\prod_{n=1}^{\infty}(1+a_n)$ is convergent.
\item The infinite product $\di\prod_{n=1}^{\infty}(1-a_n)$ is convergent.
\end{enumerate}
\end{theorem}

\begin{myproof}{Proof}Since $0<a_n<1$ for all $n\in\mathbb{Z}^+$, we find that $1+a_n>0$ and $1-a_n>0$ for all $n\in\mathbb{Z}^+$. A necessary condition for the convergence of either $\di \sum_{n=1}^{\infty}a_n$, or $\di\prod_{n=1}^{\infty}(1+a_n)$, or $\di\prod_{n=1}^{\infty}(1-a_n)$, is
\[\lim_{n\to\infty}a_n=0.\]
v

By Theorem \ref{thm220929_1}, it is then sufficient to prove that if $\{a_n\}$ is a sequence of real numbers satisfying $0< a_n<1$ for all $n\geq 1$, and $\di\lim_{n\rightarrow\infty}a_n=0$, then
the following three statements are equivalent.  
\begin{enumerate}
\item[(a)]  The series $\di \sum_{n=1}^{\infty}a_n$ is convergent.
\item[(b$^{\prime}$)] The series $\di\sum_{n=1}^{\infty}\ln(1+a_n)$ is convergent.
\item[(c$^{\prime}$)] The series $\di\sum_{n=1}^{\infty}\ln(1-a_n)$ is convergent.
\end{enumerate}\bp
Let $b_n=\ln(1+a_n)$ and $c_n=-\ln(1-a_n)$. Notice that $b_n$ and $c_n$ are also positive numbers.

Now since the sequence $\{a_n\}$ converges to 0, we find that
\[\lim_{n\rightarrow\infty}\frac{b_n}{a_n}=\lim_{n\rightarrow\infty}\frac{\ln(1+a_n)}{a_n}=\lim_{x\rightarrow 0}\frac{\ln(1+x)}{x}=1,\]
\[\lim_{n\rightarrow\infty}\frac{c_n}{a_n}=\lim_{n\rightarrow\infty}\frac{-\ln(1-a_n)}{a_n}=\lim_{x\rightarrow 0}\frac{-\ln(1-x)}{x}=1.\]

By limit comparison test for positive series, we find that
$\di\sum_{n=1}^{\infty}a_n$ is convergent if and only if $\di\sum_{n=1}^{\infty}b_n$ is convergent, and $\di\sum_{n=1}^{\infty}a_n$ is convergent if and only if $\di\sum_{n=1}^{\infty}c_n$ is convergent. These establish the equivalence of (a) and (b$^{\prime}$), and the equivalence of (a) and (c$^{\prime}$). 
\end{myproof}

 \begin{example}{}
Theorem \ref{thm220929_2} can be used to deduce the following.
\begin{enumerate}[1.]
\item[1.]  The infinite product $\di\prod_{n=1}^{\infty}\left(1+\frac{1}{n}\right)$ considered in Example \ref{ex230301_1} is divergent since the infinite series $\di\sum_{n=1}^{\infty}\frac{1}{n}$ is divergent.
\item[2.] The infinite product $\di\prod_{n=1}^{\infty}\left(1-\frac{1}{n}\right)$ considered in Example \ref{ex230301_2} is divergent  since the infinite series $\di\sum_{n=1}^{\infty}\frac{1}{n}$ is divergent.\end{enumerate}\be\begin{enumerate}[1.]
\item[3.] The  infinite product $\di\prod_{n=1}^{\infty}\left(1-\frac{1}{n^2}\right)$ considered in Example \ref{ex230301_3} is convergent  since the infinite series $\di\sum_{n=1}^{\infty}\frac{1}{n^2}$ is convergent.
\end{enumerate}
\end{example2}

\begin{theorem}[label=230301_5]{}
 If the infinite product $\di\prod_{n=1}^{\infty}(1+|a_n|)$ is convergent, then the infinite product  $\di\prod_{n=1}^{\infty}(1+a_n)$ is convergent.
\end{theorem}
\begin{myproof}{Proof} Without loss of generality, we can assume that $|a_n|<1$ for all $n\geq 1$. Given $\varepsilon>0$,
since the infinite product $\di\prod_{n=1}^{\infty}(1+|a_n|)$ is convergent, 
  Cauchy criterion says that there is a  positive integer $N$ such that for all $m\geq n\geq N$, 
\[\left|\left[\prod_{k=n}^m(1+|a_k|)\right]-1\right|<\varepsilon.\]  By an inequality  in the exercises, we find that
\[\left|\left[\prod_{k=n}^m(1+ a_k )\right]-1\right|\leq \prod_{k=n}^m(1+|a_k|)-1<\varepsilon.\]
 This proves that the infinite product $\di\prod_{n=1}^{\infty}(1+a_n)$ satisfies the Cauchy criterion. Hence, it is convergent.
\end{myproof}
\begin{definition}{Absolutely Convergent Infinite Products}
 We say that the infinite product
$\di \prod_{n=1}^{\infty}\left(1+a_n\right)$
converges absolutely if  the infinite product
$\di\prod_{n=1}^{\infty}\left(1+|a_n|\right)$
is  convergent. 
\end{definition} Theorem \ref{230301_5} says that  an absolutely convergent infinite product  is convergent.
 \begin{corollary}
{}
Let $\di \sum_{n=1}^{\infty}a_n$ be a series that converges absolutely. Then the infinite product $\di\prod_{n=1}^{\infty}\left(1+a_n\right)$ converges absolutely.
\end{corollary}

\begin{myproof}{Proof}
Since $\di\sum_{n=1}^{\infty}a_n$ converges absolutely, $\di\lim_{n\to\infty}a_n=0$. Without loss of generality, we can assume that $|a_n|<1$ for all $n\geq 1$.   Since $\di\sum_{n=1}^{\infty}|a_n|$ is convergent, Theorem \ref{thm220929_2} implies that the infinite product $\di\prod_{n=1}^{\infty}\left(1+|a_n|\right)$ is convergent. Theorem \ref{230301_5} then implies that the infinite product $\di\prod_{n=1}^{\infty}\left(1+a_n\right)$ converges absolutely.
\end{myproof}
\begin{example}{}
The infinite product $\di\prod_{n=1}^{\infty}\left(1+\frac{(-1)^{n-1}}{n^2}\right)$ is convergent since the series $\di\sum_{n=1}^{\infty}\frac{(-1)^{n-1}}{n^2}$ converges absolutely.
\end{example}
 
Now it is natural to ask the following question.
 Is it true that the infinite product $\di\prod_{n=1}^{\infty}(1+a_n)$ is convergent if and only if the series $\di\sum_{n=1}^{\infty} a_n$ is convergent?
The following two examples show that neither one implies the other.

\begin{example}{}
Let $\{a_n\}$ be the sequence defined by
\[
a_{2n-1}=\frac{1}{\sqrt{n+1}},\quad a_{2n}=-\frac{1}{\sqrt{n+1}}\hspace{1cm}\text{for}\;n\geq 1.
\] If $s_n=\di\sum_{k=1}^na_k$, then
$\di s_{2n-1}=\frac{1}{\sqrt{n+1}}$ and $ s_{2n}=0$  for all $n\geq 1$.
This implies that the series $\di\sum_{n=1}^{\infty}a_n$ converges to 0.

On the other hand,
if $P_n=\di\prod_{k=1}^n(1+a_k)$, we find that
\[P_{2n-1}=\prod_{k=2}^{n}\left(1-\frac{1}{k}\right)\left(1+\frac{1}{\sqrt{n+1}}\right),\hspace{1cm}P_{2n}=\prod_{k=2}^{n+1}\left(1-\frac{1}{k}\right).\]
Since the infinite product
$\di \prod_{n=2}^{\infty}\left(1-\frac{1}{n}\right)$ is divergent, the infinite product $\di\prod_{n=1}^{\infty}(1+a_n)$ is divergent. 

This gives an example where $\di\sum_{n=1}^{\infty}a_n$ is convergent but $\di\prod_{n=1}^{\infty}(1+a_n)$ is divergent. 
\end{example}

\begin{example}{}
Let $\{a_n\}$ be the sequence defined by
\[
a_{2n-1}=\frac{1}{\sqrt{n}},\quad a_{2n}=-\frac{1}{\sqrt{n}+1}\hspace{1cm}\text{for}\;n\geq 1.
\]  
Then 
\[1+a_{2n-1}=\frac{\sqrt{n}+1}{\sqrt{n}},\quad 1+a_{2n}=\frac{\sqrt{n}}{\sqrt{n}+1}.\]
If $P_n=\di\prod_{k=1}^n(1+a_k)$, we find that
$\di P_{2n-1}=\frac{\sqrt{n}+1}{\sqrt{n}}$ and $P_{2n}=1$  for all $n\geq 1$. Hence, the infinite product  $\di\prod_{n=1}^{\infty}(1+a_n)$ converges to 1.

If $s_n=\di\sum_{k=1}^na_k$, then
\[
s_{2n} =\sum_{k=1}^n\left(\frac{1}{\sqrt{k}}-\frac{1}{\sqrt{k}+1}\right)=\sum_{k=1}^n \frac{1}{\sqrt{k}(\sqrt{k}+1)}.
\]
Compare to
the series $\di\sum_{k=1}^{\infty}\frac{1}{k}$, we find  that the series
$\di \sum_{k=1}^{\infty} \frac{1}{\sqrt{k}(\sqrt{k}+1)}$ is divergent. Therefore, $\di\lim_{n\rightarrow\infty}s_{2n}=\infty$, which implies that $\di\lim_{n\rightarrow \infty}s_n$ does not exist.  Hence, the series $\di\sum_{n=1}^{\infty}a_n$ is divergent.

This gives an example where $\di\prod_{n=1}^{\infty}(1+a_n)$ is convergent but $\di\sum_{n=1}^{\infty}a_n$ is divergent. 
\end{example}
\vp
\noindent
{\bf \large Exercises  \thesection}
\setcounter{myquestion}{1}
\begin{question}{\themyquestion}
Given that $\{a_n\}$ is a sequence of numbers with $a_n>-1$ for all $n\in\mathbb{Z}^+$. Prove that for all $n\in \mathbb{Z}^+$,
\[\left|\prod_{k=1}^n(1+a_k)-1\right|\leq \prod_{k=1}^n (1+|a_k|)-1.\]
\end{question}
\atc
 
\begin{question}{\themyquestion}
Let $s$ be a positive number. Show that the infinite product $\di\prod_{n=2}^{\infty}\left(1-\frac{1}{n^s}\right)$  is convergent if and only if $s>1$. 
\end{question}
\atc
\begin{question}{\themyquestion}
For $n\geq 1$, let
\[u_n=\left(1+\frac{1}{n}\right)\exp\left(-\frac{1}{n}\right).\]Show that the infinite product $\di\prod_{n=1}^{\infty}u_n$ is convergent and find its value.
\end{question}
\vp
\section{Double Sequences and Double Series}\label{sec5.5}
In this section, we give a brief discussion about double sequences.
\begin{definition}{Double Sequences}
A double sequence is a function $f:\mathbb{Z}^+\times\mathbb{Z}^+\to\mathbb{R}$ that is defined on the set $\mathbb{Z}^+\times\mathbb{Z}^+$. It is customary to denote a general term  $f(m,n)$
as $a_{m,n}$, and denote the double sequence by $\{f(m,n)\}_{m,n=1}^{\infty}$ or $\{a_{m,n}\}_{m,n=1}^{\infty}$.
\end{definition}
The following gives some examples of double sequences.
\begin{example}[label=ex230301_7]
{}
\begin{enumerate}[(a)]\item 
$\di\left\{\frac{n(m+1)}{m(n+1)}\right\}_{m,n=1}^{\infty}$  
\item
$\di\left\{\frac{mn}{m^2+n^2}\right\}_{m,n=1}^{\infty}$  
\end{enumerate}
\end{example}
\begin{definition}{Convergence of Double Sequence}
We say that a double sequence  $\{a_{m,n}\}_{m,n=1}^{\infty}$ converges to a number $a$, written as
\[a=\lim_{m,n\to\infty}a_{m,n},\]provided that for every $\varepsilon>0$, there is a positive integer $N$ so that for all positive integers $m$ and $n$ with $m\geq N$, $n\geq N$, 
\[\left|a_{m,n}-a\right|<\varepsilon.\]
\end{definition}
If a double sequence converges to a number $a$, this number $a$ is unique, and we say that the sequence is convergent. Otherwise, we say that the sequence is divergent.
\begin{example}[label=ex230301_8]{}
For the double sequence $\di\left\{\frac{n(m+1)}{m(n+1)}\right\}_{m,n=1}^{\infty}$   considered in Example \ref{ex230301_7}, notice that
\[a_{m,n}=\frac{n(m+1)}{m(n+1)}=\left(1-\frac{1}{n+1}\right)\left(1+\frac{1}{m}\right)=1-\frac{1}{n+1}+\frac{1}{m}-\frac{1}{m(n+1)}.\]
Given $\varepsilon>0$, there is a positive integer $N$ so that $3/N<\varepsilon$. Then if $m\geq N$, $n\geq N$,
\[|a_{m,n}-1|<\frac{1}{n+1}+\frac{1}{m}+\frac{1}{m(n+1)}<\frac{1}{N}+\frac{1}{N}+\frac{1}{N}<\varepsilon.\]This proves that
\[\lim_{m,n\to\infty}\frac{n(m+1)}{m(n+1)}=1.\]
\end{example}

Before we study the convergence of the double sequence $\di\left\{\frac{mn}{m^2+n^2}\right\}_{m,n=1}^{\infty}$, let us prove the following lemma, which says that for a   double sequence $\{a_{m,n}\}_{m,n=1}^{\infty}$ to be convergent, it should approach the same limit regardless of how $m$ and $n$ goes to infinity.

\begin{lemma}{}
Let $\{a_{m,n}\}_{m,n=1}^{\infty}$ be a double sequence that converges to a number $a$, and let $g:\mathbb{Z}^+\to\mathbb{Z}^+$ be a function such that $\di\lim_{n\to\infty} g(n)=\infty$. Define the sequence $\{b_n\}_{n=1}^{\infty}$ by
\[b_n=a_{g(n), n}.\]
Then the sequence $\{b_n\}_{n=1}^{\infty}$ also converges to $a$.
\end{lemma}Notice that $\{g(n)\}$ is a sequence of positive integers that diverges to $\infty$.
\begin{myproof}{Proof}
Given $\varepsilon>0$, there is a positive integer $N_1$ such that for all $(m,n)\in\mathbb{Z}^+\times\mathbb{Z}^+$ with $m\geq N_1$ and $n\geq N_1$, 
\[|a_{m,n}-a|<\varepsilon.\]
Since $\di\lim_{n\to\infty} g(n)=\infty$, there is a positive integer $N\geq N_1$ such that $g(n)\geq N_1$ for all $n\geq N$. If $n\geq N$, $g(n)\geq N_1$ and $n\geq N_1$. Therefore,
\[|b_n-a|=|a_{g(n),n}-a|<\varepsilon.\]
This proves that the sequence $\{b_n\}_{n=1}^{\infty}$ converges to $a$.
\end{myproof}

\begin{example}{}
For the double sequence $\di\left\{\frac{mn}{m^2+n^2}\right\}_{m,n=1}^{\infty}$   considered in Example \ref{ex230301_7}, assume that it converges to  $a$. Take $g_1:\mathbb{Z}^+\to\mathbb{Z}^+$ to be the function $g_1(n)=n$. Then we find that
\[a=\lim_{n\to\infty}\frac{n^2}{n^2+n^2}=\frac{1}{2}.\]
Take $g_2:\mathbb{Z}^+\to\mathbb{Z}^+$ to be the function $g_2(n)=2n$. Then we find that
\[a=\lim_{n\to\infty}\frac{2n^2}{4n^2+n^2}=\frac{2}{5}.\] We get two different values of $a$. This is a contradiction. Therefore,   the double sequence $\di\left\{\frac{mn}{m^2+n^2}\right\}_{m,n=1}^{\infty}$  is divergent. 
\end{example}

It is easy to prove that linearity also holds for limits of double sequences.
\begin{proposition}{Linearity}
Assume that the double sequences  $\{a_{m,n}\}_{m,n=1}^{\infty}$ and  $\{b_{m,n}\}_{m,n=1}^{\infty}$ are convergent. Then for any constants $\alpha$ and $\beta$, the double sequence
 \[\{\alpha a_{m,n}+\beta b_{m,n}\}_{m,n=1}^{\infty}\] is also convergent, and
\[\lim_{m,n\to\infty}\left(\alpha a_{m,n}+\beta b_{m,n}\right)=\alpha  \lim_{m,n\to\infty}a_{m,n}+\beta \lim_{m,n\to\infty}b_{m,n}.\]
\end{proposition}
\begin{myproof}{Proof}Let $\di a=\lim_{m,n\to\infty}a_{m,n}$ and $\di b=\lim_{m,n\to\infty}b_{m,n}$. 
Given $\varepsilon>0$, there are positive integers $N_1$ and $N_2$ such that 
\[|a_{m,n}-a|<\frac{\varepsilon}{2(|\alpha|+1)},\hspace{1cm}\text{for all}\;m\geq N_1, n\geq N_1;\]
\[|b_{m,n}-b|<\frac{\varepsilon}{2(|\beta|+1)},\hspace{1cm}\text{for all}\;m\geq N_2, n\geq N_2.\]Let $N=\max\{N_1, N_2\}$. For all positive integers $m$ and $n$ with $m\geq N$ and $n\geq N$, we have
\begin{align*}
\left|\left(\alpha a_{m,n}+\beta b_{m,n}\right)-\left(\alpha a+\beta b\right)\right| &\leq |\alpha||a_{m,n}-a|+|\beta||b_{m,n}-b|\\
&<\frac{|\alpha|}{2(|\alpha|+1)}\varepsilon+\frac{|\beta|}{2(|\beta|+1)}\varepsilon\\&<\frac{\varepsilon}{2}+\frac{\varepsilon}{2}=\varepsilon.
\end{align*}
This proves the assertion.
\end{myproof}In the proof, we divide $\varepsilon/2$ by $|\alpha|+1$ instead of $|\alpha|$, because $\alpha$ can be 0.

For the double sequence we considered in Example \ref{ex230301_8}, notice that 
\[\lim_{m\to\infty}\left(\lim_{n\to\infty}\frac{n(m+1)}{m(n+1)}\right)=\lim_{m\to\infty}\frac{m+1}{m}=1=\lim_{m,n\to\infty}\frac{n(m+1)}{m(n+1)}.\]The question is whether we can find the limit of a double seqeunce $\{a_{m,n}\}$ by taking the limit $n\to\infty$ first, and then take the limit $m\to\infty$, or in the opposite order. 

For the double sequence $\di\left\{\frac{mn}{m^2+n^2}\right\}_{m,n=1}^{\infty}$, for fixed $m\geq 1$, taking the $n\to\infty$ limit, we have
\[\lim_{n\to\infty}\frac{mn}{m^2+n^2}=0.\] Hence,
\[\lim_{m\to\infty}\left(\lim_{n\to\infty}\frac{mn}{m^2+n^2}\right)=0.\]But we have shown  that the double sequence $\di\left\{\frac{mn}{m^2+n^2}\right\}_{m,n=1}^{\infty}$ is divergent.

Therefore, we find that to study the limit of a double sequence, in general we cannot take one limit first before the other. The following theorem says that  if one knows apriori that the double sequence is convergent, one can take iterated limits under some conditions.
\begin{theorem}[label=230301_10]{}
Assume that the double sequence $\{a_{m,n}\}_{m,n}^{\infty}$ converges to $a$, and for each $m\in\mathbb{Z}^+$, the limit
\[b_m=\lim_{n\to\infty}a_{m,n}\] exists. Then the sequence $\{b_m\}$ also converges to $a$. In other words,
\[\lim_{m,n\to\infty}a_{m,n}=a\quad\implies\quad \lim_{m\to\infty}\left(\lim_{n\to\infty}a_{m,n}\right)=a\] provided that the limit $\di\lim_{n\to\infty}a_{m,n}$ exists for all $m\in\mathbb{Z}^+$.
\end{theorem}
\begin{myproof}{Proof}
Given $\varepsilon>0$, there exists a positive integer $N$ such that
for all $(m,n)\in\mathbb{Z}^+\times\mathbb{Z}^+$ with $m\geq N$ and $n\geq N$, 
\[|a_{m,n}-a|<\frac{\varepsilon}{2}.\]
Hence, for fixed $m\geq N$, taking the $n\to\infty$ limit gives
\[|b_m-a|\leq\frac{\varepsilon}{2}<\varepsilon.\]
This proves that $\di\lim_{n\to\infty}b_m=a$.
\end{myproof}

The assumption that the limit $\di\lim_{n\to\infty}a_{m,n}$ exists for each $m\in\mathbb{Z}^+$ in Theorem \ref{230301_10} is needed, as the convergence of the double sequence  $\{a_{m,n}\}_{m,n}^{\infty}$ does not guarantee that the  limit $\di\lim_{n\to\infty}a_{m,n}$ exists. An example is shown below.
\begin{example}{}
Consider the double sequence $\{a_{m,n}\}_{m,n=1}^{\infty}$ with
\[a_{m,n}=\frac{m+(-1)^{n-1}}{m^2 }.\]
For fixed $m\in\mathbb{Z}^+$, the sequence $\di\left\{\frac{m+(-1)^{n-1}}{m^2 }\right\}_{n=1}^{\infty}$ does not have a limit since it is oscillating between $\di\frac{m+1}{m^2}$ and $\di\frac{m-1}{m^2}$. But the double sequence $\{a_{m,n}\}_{m,n=1}^{\infty}$ converges to zero. This can be proved in the following way. Given $\varepsilon>0$, since $\di\lim_{m\to\infty}\frac{m+1}{m^2}=0$, there exists a positive integer $N$ so that for all $m\geq N$.
\[0<\frac{m+1}{m^2}<\varepsilon.\] This implies that if $m\geq N$, $n\geq N$, then
\[0\leq a_{m,n}\leq\frac{m+1}{m^2}<\varepsilon.\]Hence, the double sequence $\{a_{m,n}\}_{m,n=1}^{\infty}$ converges to zero. 

\end{example}

\begin{definition}{Bounded Double Sequence}
We say that a double sequence $\{a_{m,n}\}_{m,n=1}^{\infty}$  is bounded if the set \[\left\{a_{m,n}\,|\,(m,n)\in\mathbb{Z}^+\times \mathbb{Z}^+\right\}\] is bounded.
\end{definition}
\begin{remark}
{}If a double sequence $\{a_{m,n}\}_{m,n=1}^{\infty}$  is convergent, it is not necessarily bounded. For example, consider the double  sequence  $\{a_{m,n}\}_{m,n=1}^{\infty}$ with
\[a_{m,n}=\begin{cases}n,\quad &\text{if}\; m=1,\\\di 1,\quad & \text{if}\; m\geq 2.\end{cases}\]
Obviously, it is not bounded. However, It is not difficult to prove that the double sequence $\{a_{m,n}\}_{m,n=1}^{\infty}$   converges to 1.
\end{remark}

\begin{definition}{Increasing Double Sequence}We say that  a double sequence $\{a_{m,n}\}_{m,n=1}^{\infty}$  is increasing in both indices provided that for fixed $m\in \mathbb{Z}^+$, $\{a_{m,n}\}_{n=1}^{\infty}$ is an increasing sequence in $n$; and for fixed $n\in\mathbb{Z}^+$, $\{a_{m,n}\}_{m=1}^{\infty}$ is an increasing sequence in $m$. \end{definition}

If a  double sequence $\{a_{m,n}\}_{m,n=1}^{\infty}$  is increasing in both indices, for any positive integers $m_1, m_2, n_1, n_2$, if $m_2\geq m_1$ and $n_2\geq n_1$, then
\[a_{m_2, n_2}\geq a_{m_1,n_1}.\]

\begin{figure}[ht]
\centering
\includegraphics[scale=0.2]{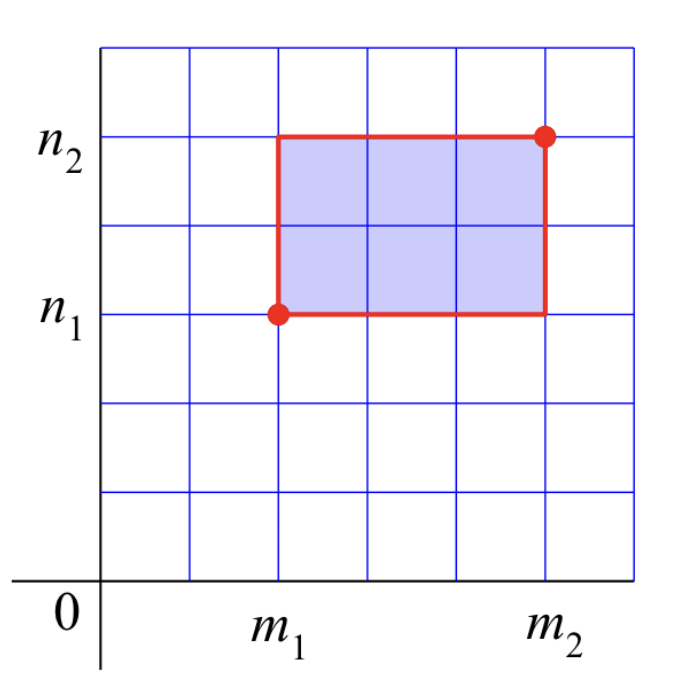}
\caption{An illiustration of the relative positions of $(m_1,n_1)$ and $(m_2, n_2)$ when $m_2>m_1$ and $n_2>n_1$.\fa}\label{figure51}
\end{figure}

 \begin{example}{}
The double sequence $\{a_{m,n}\}_{m,n=1}^{\infty}$ with
\[a_{m,n}=\frac{mn}{(m+1)(n+1)}\] is increasing in both indices.
\end{example}

The following is a counterpart of monotone convergence theorem for double sequences. 
 
\begin{theorem}[label=230301_11]
{Convergence of Increasing Double Sequences}
Let $\{a_{m,n}\}_{m,n=1}^{\infty}$ be a double sequence that is increasing in both indices. Then the double sequence $\{a_{m,n}\}_{m,n=1}^{\infty}$ is convergent if and only if it is bounded above. In case it is convergent, it converges to $\di\sup_{(m,n)\in\mathbb{Z}^+\times\mathbb{Z}^+}\{a_{m,n}\}$.
\end{theorem}
\begin{myproof}{Proof}
If the sequence $\{a_{m,n}\}_{m,n=1}^{\infty}$ converges to $a$, then there is a positive integer $N$ such that for all $(m,n)\in\mathbb{Z}^+\times\mathbb{Z}^+$ with $m\geq N$ and $n\geq N$, 
\[|a_{m,n}-a|<1.\]
This implies that 
\[a_{m,n}<a+1 \quad \text{ for all }\; m\geq N, n\geq N.\]Given $(m,n)\in\mathbb{Z}^+\times \mathbb{Z}^+$, let $k=\max\{m,n,N\}$. Then $k\geq m$, $k \geq n$ and $k\geq N$. Therefore,
\[a_{m,n}\leq  a_{k,k}<a+1.\]\bp
This prove  that the double sequence $\{a_{m,n}\}_{m,n=1}^{\infty}$ is bounded above by $a+1$.
In fact, the same reasoning shows that it is bounded above by $a+\varepsilon$ for any $\varepsilon>0$, but we do not need this.
 
Conversely, if $\{a_{m,n}\}_{m,n=1}^{\infty}$ is bounded above, then 
\[a=\sup_{(m,n)\in\mathbb{Z}^+\times\mathbb{Z}^+}\{a_{m,n}\}\]  exists. Given $\varepsilon>0$, there exists $(m_0, n_0)\in\mathbb{Z}^+\times \mathbb{Z}^+$ such that 
\[a_{m_0,n_0}>a-\varepsilon.\]
 Take $N=\max\{m_0,n_0\}$. Then if $m\geq N\geq m_0$, $n\geq N\geq n_0$,
\[a_{m,n}\geq   a_{m_0,n_0}> a-\varepsilon.\]
By definition $a_{m,n}\leq a$. Therefore, for all $(m,n)\in\mathbb{Z}^+\times\mathbb{Z}^+$ with $m\geq N$ and $n\geq N$, we have
\[|a_{m,n}-a|<\varepsilon.\]This proves that the double sequence $\{a_{m,n}\}_{m,n=1}^{\infty}$ is convergent and it converges to $\di a=\sup_{(m,n)\in\mathbb{Z}^+\times\mathbb{Z}^+}\{a_{m,n}\}$.
\end{myproof}

Now we turn to double series. A double series is a series of the form
\[\sum_{(m,n)\in\mathbb{Z}^+\times\mathbb{Z}^+}a_{m,n},\]where  $\{a_{m,n}\}_{m,n=1}^{\infty}$  is a double sequence. For each $(m,n)\in \mathbb{Z}^+\times\mathbb{Z}^+$, we define the $(m,n)$ partial sum $s_{m,n}$ by
\[s_{m,n}=\sum_{k=1}^m\sum_{l=1}^na_{k,l}.\]

\begin{definition}{Convergence of Double Series}
We say that the double series $\di\sum_{(m,n)\in\mathbb{Z}^+\times\mathbb{Z}^+}a_{m,n}$ is convergent provided that the double sequence of partial sums $\{s_{m,n}\}$ is convergent. In this case, the sum of the double series is
\[\sum_{(m,n)\in\mathbb{Z}^+\times\mathbb{Z}^+}a_{m,n}=s=\lim_{m,n\to\infty}s_{m,n}=\lim_{m,n\to\infty}\sum_{k=1}^m\sum_{l=1}^na_{k,l}.\]
\end{definition}

Notice that for any $(m,n)\in\mathbb{Z}^+\times \mathbb{Z}^+$,
\[s_{m,n}-s_{m,n-1}=\sum_{k=1}^m a_{k,n}.\]
Therefore,
\[s_{m,n}-s_{m,n-1}-s_{m-1,n}+s_{m-1, n-1}=\sum_{k=1}^ma_{k,n}-\sum_{k=1}^{m-1}a_{k,n}=a_{m,n}.\]
From this, we obtain the following immediately.

\begin{proposition}{}
If the double series $\di \sum_{(m,n)\in\mathbb{Z}^+\times\mathbb{Z}^+}a_{m,n}$ is convergent, then the double sequence   $\{a_{m,n}\}_{m,n=1}^{\infty}$  converges to 0.
\end{proposition}
 
\begin{figure}[ht]
\centering
\includegraphics[scale=0.2]{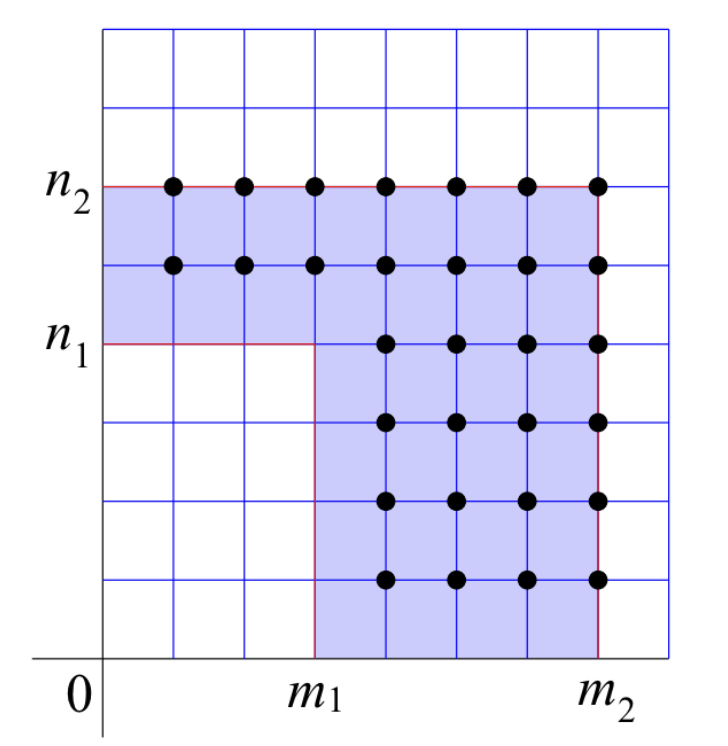}
\caption{An illiustration of those terms $a_{k,l}$  that involved in $s_{m_2, n_2}- s_{m_1, n_1}$ when $m_2>m_1$ and $n_2>n_1$.\fa}\label{figure49}
\end{figure}

If $\di\sum_{(m,n)\in\mathbb{Z}^+\times\mathbb{Z}^+}a_{m,n}$ is a double series with $a_{m,n}\geq 0$ for all $(m,n)\in\mathbb{Z}^+\times\mathbb{Z}^+$, then the double sequence of partial sums $\{s_{m,n}\}$ is a double sequence that is increasing in both indices.
From Theorem \ref{230301_11}, we obtain the following.
\begin{theorem}[label=230302_2]{}
If $\di\sum_{(m,n)\in\mathbb{Z}^+\times\mathbb{Z}^+}a_{m,n}$ is a double series with $a_{m,n}\geq 0$ for all $(m,n)\in\mathbb{Z}^+\times\mathbb{Z}^+$, then it is convergent if and only if the double sequence of partial sums $\{s_{m,n}\}$ is bounded above.
\end{theorem}
\begin{corollary}[label=230302_16]{}
If $\di\sum_{(m,n)\in\mathbb{Z}^+\times\mathbb{Z}^+}a_{m,n}$ is a double series with $a_{m,n}\geq 0$ for all $(m,n)\in\mathbb{Z}^+\times\mathbb{Z}^+$, then it is convergent if and only if the  sequence   $\{s_{n,n}\}_{n=1}^{\infty}$ is convergent. In this case,
\[\sum_{(m,n)\in\mathbb{Z}^+\times\mathbb{Z}^+}a_{m,n}=\lim_{n\to\infty}s_{n,n}=\lim_{n\to\infty}\sum_{k=1}^n\sum_{l=1}^na_{k,l}.\]
\end{corollary}
This says that we can determine the convergence of a nonnegative double series from the sequence $\{s_{n,n}\}_{n=1}^{\infty}$ instead of the double sequence $\{s_{m,n}\}_{m,n=1}^{\infty}$.
\begin{myproof}{Proof}Since $a_{m,n}\geq 0$ for all $(m,n)\in\mathbb{Z}^+\times\mathbb{Z}^+$, the double sequence $\{s_{m,n}\}$ is increasing in both indices, while the sequence $\{s_{n,n}\}$ is increasing.

If the double series $\di\sum_{(m,n)\in\mathbb{Z}^+\times\mathbb{Z}^+}a_{m,n}$ is convergent, Theorem \ref{230302_2} implies that the double sequence of partial sums $\{s_{m,n}\}$ is bounded above. Being a subset, the sequence   $\{s_{n,n}\}_{n=1}^{\infty}$ is also bounded above. By monotone convergence theorem,  the  sequence   $\{s_{n,n}\}_{n=1}^{\infty}$ is convergent. 
\bp
Conversely, assume that  the  sequence   $\{s_{n,n}\}_{n=1}^{\infty}$ is convergent. Then it is bounded above. Let \[t=\sup_{n\in\mathbb{Z}^+}s_{n,n}=\lim_{n\to\infty}s_{n,n}.\]For any positive integers $m$ and $n$, \[s_{m,n}\leq\max\{s_{m,m},s_{n,n}\}\leq t.\]
This implies that the double sequence   $\{s_{m,n}\}$ is bounded above by $t$. Hence, the double series $\di\sum_{(m,n)\in\mathbb{Z}^+\times\mathbb{Z}^+}a_{m,n}$ is convergent.  From the argument above, we also find that
\[\sup_{(m,n)\in\mathbb{Z}^+\times\mathbb{Z}^+}s_{m,n}\leq t=\sup_{n\in\mathbb{Z}^+}s_{n,n}.\]
Since the oppositie inequality is obvious, this is in fact an equality. Hence,
\begin{align*}\sum_{(m,n)\in\mathbb{Z}^+\times\mathbb{Z}^+}a_{m,n}&=\sup_{(m,n)\in\mathbb{Z}^+\times\mathbb{Z}^+}s_{m,n}=\sup_{n\in\mathbb{Z}^+}s_{n,n}\\&=\lim_{n\to\infty}s_{n,n}=\lim_{n\to\infty}\sum_{k=1}^n\sum_{l=1}^na_{k,l}.\end{align*}
\end{myproof}

Let us look at an example.
\begin{example}{}
Show that the double series
\[\sum_{(m,n)\in\mathbb{Z}^+\times\mathbb{Z}^+}\frac{1}{(m^2+n^2)^2}\] is convergent.
\end{example}
\begin{solution}{Solution}
Notice that 
\begin{align*}
s_{n,n}&=\sum_{k=1}^n\sum_{l=1}^n\frac{1}{(k^2+l^2)^2} \leq  \sum_{k=1}^n\sum_{l=1}^k \frac{1}{(k^2+l^2)^2}+\sum_{l=1}^n\sum_{k=1}^l \frac{1}{(k^2+l^2)^2}\\&
\leq 2\sum_{k=1}^n \sum_{l=1}^k\frac{1}{k^4} 
=2\sum_{k=1}^n\frac{k}{k^4}= 2\di\sum_{k=1}^{\infty}\frac{1}{k^3}.
\end{align*}Since the series $\di\sum_{k=1}^{\infty}\frac{1}{k^3}$ is convergent,     the sequence $\{s_{n,n}\}$ is bounded above. Hence,   the double series $\di\sum_{(m,n)\in\mathbb{Z}^+\times\mathbb{Z}^+}\frac{1}{(m^2+n^2)^2}$ is convergent.
\end{solution}

Next, we consider double series that have negative terms. 
Given a double sequence $\{a_{m,n}\}_{m,n=1}^{\infty}$, let $\{p_{m,n}\}_{m,n=1}^{\infty}$ and $\{q_{m,n}\}_{m,n=1}^{\infty}$ be double sequences defined by
\[p_{m,n}=\frac{|a_{m,n}|+a_{m,n}}{2},\hspace{1cm}q_{m,n}=\frac{|a_{m,n}|-a_{m,n}}{2}.\]Then
\[|a_{m,n}|=p_{m,n}+q_{m,n},\hspace{1cm}a_{m,n}=p_{m,n}-q_{m,n}.\]
 $\{p_{m,n}\}_{m,n=1}^{\infty}$ and $\{q_{m,n}\}_{m,n=1}^{\infty}$ are nonnegative double sequences with
\[0\leq p_{m,n}\leq |a_{m,n}|,\hspace{1cm}0\leq q_{m,n}\leq |a_{m,n}|.\]

\begin{definition}{Absolute Convergence of Double Series}We say that the double series $\di \sum_{(m,n)\in\mathbb{Z}^+\times\mathbb{Z}^+}a_{m,n}$
converges absolutely if the double series  $\di \sum_{(m,n)\in\mathbb{Z}^+\times\mathbb{Z}^+}|a_{m,n}|$ is convergent.
\end{definition}
\begin{theorem}[label=230302_12]{}
If the double series $\di \sum_{(m,n)\in\mathbb{Z}^+\times\mathbb{Z}^+}a_{m,n}$ converges absolutely, then it is convergent.
\end{theorem}
\begin{myproof}{Proof}
 For $(m,n)\in\mathbb{Z}^+\times\mathbb{Z}^+$, let
\[p_{m,n}=\frac{|a_{m,n}|+a_{m,n}}{2},\hspace{1cm}q_{m,n}=\frac{|a_{m,n}|-a_{m,n}}{2}.\]Then \[0\leq p_{m,n}\leq |a_{m,n}|,\hspace{1cm}0\leq q_{m,n}\leq |a_{m,n}|.\]Let 
$\{s^+_{m,n}\}$, $\{s^-_{m,n}\}$, $\{t_{m,n}\}$ and $\{s_{m,n}\}$ be respectively the double sequences of partial sums for the double series $\di\sum_{(m,n)\in\mathbb{Z}^+\times\mathbb{Z}^+}p_{m,n}$, $\di\sum_{(m,n)\in\mathbb{Z}^+\times\mathbb{Z}^+}q_{m,n}$, $\di\sum_{(m,n)\in\mathbb{Z}^+\times\mathbb{Z}^+}|a_{m,n}|$ and $\di\sum_{(m,n)\in\mathbb{Z}^+\times\mathbb{Z}^+}a_{m,n}$. Then
\[t_{m,n}=s^+_{m,n}+s^-_{m,n},\hspace{1cm} s_{m,n}=s^+_{m,n}-s^-_{m,n}.\]Moreover,
\begin{equation}\label{eq230302_4}0\leq s^+_{m,n}\leq t_{m,n},\hspace{1cm}0\leq s^-_{m,n}\leq t_{m,n}.\end{equation}
Since $\{p_{m,n}\}$, $\{q_{m,n}\}$ and $\{|a_{m,n}|\}$ are nonnegative double sequences, $\{s^+_{m,n}\}$, $\{s^-_{m,n}\}$ and $\{t_{m,n}\}$ are nonnegative double sequences that are increasing in both indices. By assumption, the double series $\di \sum_{(m,n)\in\mathbb{Z}^+\times\mathbb{Z}^+}|a_{m,n}|$ is convergent. Therefore, the double sequence $\{t_{m,n}\}$ is bounded above. Eq. \eqref{eq230302_4} implies that the double sequences  $\{s^+_{m,n}\}$ and  $\{s^-_{m,n}\}$ are also bounded above. Hence, the   double sequences  $\{s^+_{m,n}\}$ and  $\{s^-_{m,n}\}$ are convergent. By linearity, the double sequence $\{s_{m,n}\}$ is also convergent and
\[\lim_{m,n\to\infty}s_{m,n}=\lim_{m,n\to\infty}s^+_{m,n}-\lim_{m,n\to\infty}s^-_{m,n}.\] This proves that the double series $\di \sum_{(m,n)\in\mathbb{Z}^+\times\mathbb{Z}^+}a_{m,n}$ is convergent, and
\[\sum_{(m,n)\in\mathbb{Z}^+\times\mathbb{Z}^+}a_{m,n}=\sum_{(m,n)\in\mathbb{Z}^+\times\mathbb{Z}^+}p_{m,n}-\sum_{(m,n)\in\mathbb{Z}^+\times\mathbb{Z}^+}q_{m,n}.\]
\end{myproof}There is a simpler proof of this theorem using the same idea as we prove the case for single series. The ideas in the proof that we present above have been used when we prove that any rearrangement of an absolutely convergent single series is convergent and has the same sum. It is a useful technique  for dealing with absolutely convergent series. One should compare this proof to the proof of Theorem \ref{230224_8} for convergence of improper integrals. In fact, infinite series and improper integrals are closely related. An improper integral $\di\int_{-\infty}^{\infty} f(x)dx$ is convergent if and only if the double limit
\[\lim_{a\to-\infty, b\to\infty}\int_a^bf(x)dx\] exists. This can be rephrased as for any two sequences $\{a_m\}$ and $\{b_n\}$ satisfying   $\di\lim_{m\to\infty}a_m=-\infty$ and $\di\lim_{n\to \infty}b_n=\infty$, the double sequence $\{F_{m,n}\}$, with 
\[F_{m,n}=\int_{a_m}^{b_n}f(x)dx\] is convergent and has the same limit.

As we have seen before, we cannot simply compute the limit of a double sequence by taking the limit with respect to one index first before the other. For double series, we cannot find the sum simply by taking the sum with respect to one index first before the other. Let us look at the following example.
\begin{example}[label=ex230301_6]{}
For $(m,n)\in\mathbb{Z}^+\times\mathbb{Z}^+$, let
\[a_{m,n}=\begin{cases} m-n, \quad &\text{if}\;|m-n|=1,\\0,\quad &\text{otherwise},\end{cases}\] and consider the double series $\di\sum_{(m,n)\in\mathbb{Z}^+\times\mathbb{Z}^+}a_{m,n}$. We find that
\begin{align*}
\sum_{n=1}^{\infty}a_{m,n}=\begin{cases}-1,\quad &\text{if} \; m=1,\\0,\quad &\text{if}\;m\geq 2;\end{cases}\hspace{1cm} \sum_{m=1}^{\infty}a_{m,n}=\begin{cases}1,\quad &\text{if} \; n=1,\\0,\quad &\text{if}\;n\geq 2.\end{cases}
\end{align*}\be Therefore,
\[\sum_{m=1}^{\infty}\left(\sum_{n=1}^{\infty}a_{m,n}\right)=-1,\hspace{1cm}\sum_{n=1}^{\infty}\left(\sum_{m=1}^{\infty}a_{m,n}\right)=1.\]We find that changing the orders of summation produces different sums.
\end{example2}

\begin{figure}[ht]
\centering
\includegraphics[scale=0.2]{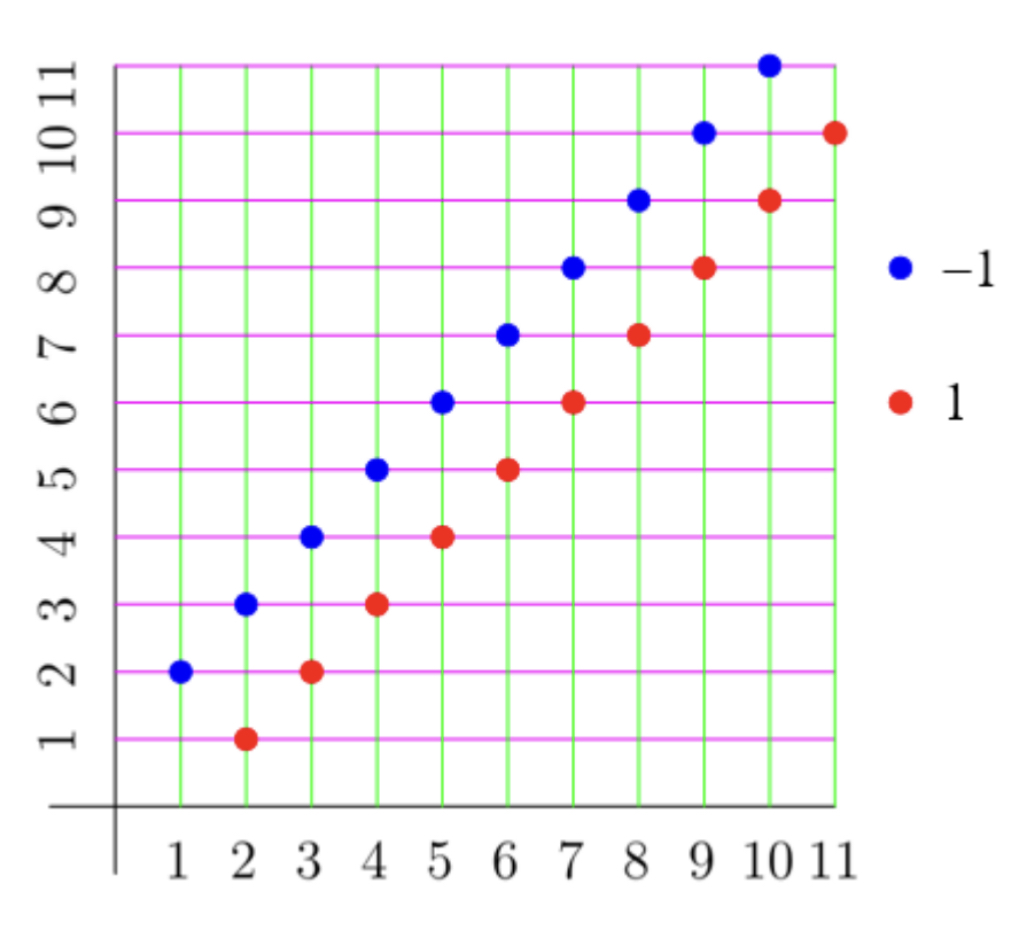}
\caption{An illiustration of the terms in the double series defined in Example \ref{ex230301_6}.\fa}\label{figure50}
\end{figure}

However, we have the following if the double series is convergent.
\begin{theorem}[label=230302_8]
{}
Assume that the double series $\di\sum_{(m,n)\in\mathbb{Z}^+\times\mathbb{Z}^+}a_{m,n}$ converges to $s$, and for every fixed $m\in\mathbb{Z}^+$, the series
$\di\sum_{n=1}^{\infty}a_{m,n}$ is convergent with sum $u_m$. Then the series \[\sum_{m=1}^{\infty}u_m=\sum_{m=1}^{\infty}\left(\sum_{n=1}^{\infty}a_{m,n}\right)\] is convergent and its sum is $s$.
\end{theorem}
\begin{myproof}{Proof}
  Let 
\[s_{m,n}=\sum_{k=1}^m\sum_{l=1}^na_{k,l}.\]  We are given that the double sequence $\{s_{m,n}\}_{m,n=1}^{\infty}$ converges to $s$. Notice that for fixed $m\in\mathbb{Z}^+$,
\[\sum_{k=1}^m u_k=\sum_{k=1}^m \lim_{n\to\infty}\sum_{l=1}^na_{k,l} =\lim_{n\to\infty}s_{m,n}.\] This shows that for fixed $m$, the limit $\di b_m= \lim_{n\to\infty}s_{m,n} $ exists and it equal to $\di\sum_{k=1}^m u_k$.
By Theorem \ref{230301_10}, the sequence $\{b_m\}$ converges to $s$. Therefore, the series $\di\sum_{m=1}^{\infty}u_m$ is convergent and has sum $s$.
\end{myproof}

Let us explore more about
  nonnegative double series first.

\begin{theorem}[label=230302_6]{}
Given that $\di\sum_{(m,n)\in\mathbb{Z}^+\times\mathbb{Z}^+}a_{m,n}$ is a double series with $a_{m,n}\geq 0$ for all $(m,n)\in\mathbb{Z}^+\times\mathbb{Z}^+$, and 
it is convergent with sum $s$.  We have the following.
\begin{enumerate}[(a)]
\item For all $m\in\mathbb{Z}^+$, $\di u_m=\sum_{n=1}^{\infty}a_{m,n}$ is finite.
\item For all $n\in\mathbb{Z}^+$, $\di v_n=\sum_{m=1}^{\infty}a_{m,n}$ is finite.
\item The series $\di\sum_{m=1}^{\infty}u_m$ and the series $\di\sum_{n=1}^{\infty}v_n$ both converge to $s$. Namely,
\[\sum_{m=1}^{\infty} \sum_{n=1}^{\infty} a_{m,n} =\sum_{n=1}^{\infty} \sum_{m=1}^{\infty} a_{m,n}=\sum_{(m,n)\in\mathbb{Z}^+\times\mathbb{Z}^+}a_{m,n}.\]
 \end{enumerate}
\end{theorem}
\begin{myproof}{Proof}
Given positive integers $m$ and $n$, let
\[s_{m,n}=\sum_{k=1}^m\sum_{l=1}^na_{m,n},\quad u_{m,n}=\sum_{l=1}^n a_{m,l},\quad v_{m,n}=\sum_{k=1}^ma_{k,n}.\]Then 
\begin{equation}\label{eq230302_7}
s_{m,n}=\sum_{k=1}^mu_{k,n}=\sum_{l=1}^n v_{m,l}.\end{equation}Since $a_{m,n}\geq 0$ for all $m, n\in\mathbb{Z}^+$, we have
\[u_{m,n}\leq s_{m,n},\hspace{1cm}v_{m,n}\leq s_{m,n}\hspace{1cm}\text{for all}\;(m,n)\in\mathbb{Z}^+\times\mathbb{Z}^+.\]
For fixed $m$, $\{u_{m,n}\}_{n=1}^{\infty}$ and  $\{s_{m,n}\}_{n=1}^{\infty}$ are increasing sequences. For fixed $n$, $\{v_{m,n}\}_{m=1}^{\infty}$  and  $\{s_{m,n}\}_{m=1}^{\infty}$ are  increasing sequences.  

Since the double series $\di\sum_{(m,n)\in\mathbb{Z}^+\times\mathbb{Z}^+}a_{m,n}$ is convergent with sum $s$, $s_{m,n}\leq s$ for all positive integers $m$ and $n$. Therefore,
the sequences $\{u_{m,n}\}_{n=1}^{\infty}$, $\{s_{m,n}\}_{m=1}^{\infty}$, $\{v_{m,n}\}_{m=1}^{\infty}$ and $\{s_{m,n}\}_{m=1}^{\infty}$ are increasing sequences that are bounded above by $s$. 
Therefore, each of these sequences is convergent.  The convergence of the sequences $\{u_{m,n}\}_{n=1}^{\infty}$ and $\{v_{m,n}\}_{m=1}^{\infty}$ are precisely the statements in (a) and (b). By definition,
 \[u_m=\sum_{n=1}^{\infty}a_{m,n}=\lim_{n\to\infty}u_{m,n},\hspace{1cm}v_n=\sum_{m=1}^{\infty}a_{m,n}=\lim_{m\to\infty}v_{m,n}.\]
Now let
\[b_m=\sum_{k=1}^mu_k \hspace{1cm}\text{and}\hspace{1cm}c_n=\sum_{l=1}^nv_l\] be the partial sums of the series $\di\sum_{m=1}^{\infty}u_m$ and $\di\sum_{n=1}^{\infty}v_n$. 
From \eqref{eq230302_7}, we find that
\[\lim_{n\to\infty}s_{m,n}=\sum_{k=1}^mu_k=b_m,\hspace{1cm}\lim_{m\to\infty}s_{m,n}=\sum_{l=1}^nv_l=c_n.\]From these, we find that the sequences $\{b_m\}$ and $\{c_n\}$ are also increasing sequences that are bounded above by $s$.\bp Therefore, 
\[b=\lim_{m\to\infty}b_m\hspace{1cm}\text{and}\hspace{1cm}c=\lim_{n\to\infty}c_n\] exist, and $b\leq s$, $c\leq s$. We are now left to prove that $b=c=s$. It is sufficient to prove that $b=s$. Then $c=s$ follows by interchanging the roles of $m$ and $n$.
Given $\varepsilon>0$, using the fact that $s=\sup\di\left\{s_{n,n}\,|\,n\in\mathbb{Z}^+\right\}$ from Corollary \ref{230302_16}, we find that there is a positive integer $N$ such that
\[s_{N,N}>s-\varepsilon.\]But then
\[s_{N,N}=\sum_{m=1}^N\sum_{n=1}^Na_{m,n}\leq \sum_{m=1}^N\sum_{n=1}^{\infty}a_{m,n}=b_N.\]This shows that
\[b_N>s-\varepsilon.\]
Hence,
\[b=\sup_m b_m> s-\varepsilon.\] Since $\varepsilon>0$ is arbitrary, we conclude that $b\geq s$. Together with $b\leq s$ that is proved earlier, we conclude that $b=s$.

\end{myproof}

\begin{theorem}[label=230302_3]{}
Given that $\di\sum_{(m,n)\in\mathbb{Z}^+\times\mathbb{Z}^+}a_{m,n}$ is a double series with $a_{m,n}\geq 0$ for all $(m,n)\in\mathbb{Z}^+\times\mathbb{Z}^+$.
Assume that for each $m\in\mathbb{Z}^+$, the series $\di\sum_{n=1}^{\infty}a_{m,n}$ converges to $u_m$. If the series
$\di\sum_{m=1}^{\infty}u_m$ is convergent, then the double series $\di\sum_{(m,n)\in\mathbb{Z}^+\times\mathbb{Z}^+}a_{m,n}$ is convergent, and 
\[\sum_{(m,n)\in\mathbb{Z}^+\times\mathbb{Z}^+}a_{m,n}=\sum_{m=1}^{\infty}u_m=\sum_{m=1}^{\infty}\left( \sum_{n=1}^{\infty} a_{m,n}\right).\]

\end{theorem}
\begin{myproof}{Proof}
It is sufficient to prove that the convergence of the series
$\di\sum_{m=1}^{\infty}u_m$ implies the convergence of the double series $\di\sum_{(m,n)\in\mathbb{Z}^+\times\mathbb{Z}^+}a_{m,n}$. The last statement then follows from Theorem \ref{230302_8}. Assume that the series
$\di\sum_{m=1}^{\infty}u_m$ converges to $u$. Using the same notations  as in the proof of Theorem \ref{230302_6}, we find 
that for each positive integer $m$, the sequence $\{u_{m,n}\}_{n=1}^{\infty}$ increases to $u_m$. From \eqref{eq230302_7}, we find  that for any positive integers $m$ and $n$, 
\[s_{m,n}\leq\sum_{k=1}^m u_k\leq u.\]This shows that the double sequence $\{s_{m,n}\}_{m,n=1}^{\infty}$ is bounded above, and hence it is convergent. Therefore,  the double series $\di\sum_{(m,n)\in\mathbb{Z}^+\times\mathbb{Z}^+}a_{m,n}$ is convergent. 

\end{myproof}
\begin{remark}{}
Putting together Theorem \ref{230302_6} and Theorem \ref{230302_3}, we conclude the following. Given a double series $\di\sum_{(m,n)\in\mathbb{Z}^+\times\mathbb{Z}^+}a_{m,n}$    with nonnegative terms $a_{m,n}$, we can determine its convergence and find its sum by first checking whether for each fixed $m$, the series $\di\sum_{n=1}^{\infty}a_{m,n}$ is convergent. If yes,  find the sum, call it as $u_m$, and check whether the series $\di\sum_{m=1}^{\infty}u_m$ is convergent. If yes, then the double series $\di\sum_{(m,n)\in\mathbb{Z}^+\times\mathbb{Z}^+}a_{m,n}$   is convergent and its sum is given by $\di\sum_{m=1}^{\infty}u_m$. Namely, the sum of the double series $\di\sum_{(m,n)\in\mathbb{Z}^+\times\mathbb{Z}^+}a_{m,n}$  can be obtained by iterated summation. \end{remark}
\begin{highlight}{} We can also start with the series $\di\sum_{m=1}^{\infty}a_{m,n}$ for each fixed $n$. This shows that for double series with nonnegative terms,   we can   interchange the orders of summation. In fact, with slightly more effort, one can prove that we can sum in any orders.

If for some integer $m$, the series $\di\sum_{n=1}^{\infty}a_{m,n}$ is divergent, then the double series $\di\sum_{(m,n)\in\mathbb{Z}^+\times\mathbb{Z}^+}a_{m,n}$   is divergent. Even if the series $\di\sum_{n=1}^{\infty}a_{m,n}$ is convergent for all positive integers $m$, the series $\di\sum_{m=1}^{\infty}u_m$ can still be divergent. In this latter case, the double series $\di\sum_{(m,n)\in\mathbb{Z}^+\times\mathbb{Z}^+}a_{m,n}$   is divergent. An example is given by the double series
\[\sum_{(m,n)\in\mathbb{Z}^+\times\mathbb{Z}^+}\frac{1}{m^2+n^2}.\]
For fixed positive integer $m$, comparison with the series $\di \sum_{n=1}^{\infty}\frac{1}{n^2}$ shows that the series $\di\sum_{n=1}^{\infty}\frac{1}{m^2+n^2}$ is convergent.
By integral test, we find that
\[u_m=\sum_{n=1}^{\infty}\frac{1}{m^2+n^2}\geq \int_0^{\infty}\frac{1}{m^2+x^2}dx-\frac{1}{m^2}=\frac{\pi}{2m}-\frac{1}{m^2}>0.\]Since the series $\di\sum_{m=1}^{\infty}\frac{1}{m}$ is divergent but the series  $\di \sum_{n=1}^{\infty}\frac{1}{m^2}$  is convergent, the series $\di \sum_{m=1}^{\infty}\left(\frac{\pi}{2m}-\frac{1}{m^2}\right)$ is divergent. Hence, the series $\di\sum_{m=1}^{\infty}u_m$ is divergent.
\end{highlight}

Finally, we can come back to series with negative terms.
From Theorem \ref{230302_6}, we have the following.
\begin{theorem}[label=230302_10]{}
Given that $\di\sum_{(m,n)\in\mathbb{Z}^+\times\mathbb{Z}^+}a_{m,n}$ is a double series that converges absolutely, and 
it is convergent with sum $s$.  We have the following.
\begin{enumerate}[(a)]
\item For all $m\in\mathbb{Z}^+$, $\di u_m=\sum_{n=1}^{\infty}a_{m,n}$ is finite.
\item For all $n\in\mathbb{Z}^+$, $\di v_n=\sum_{m=1}^{\infty}a_{m,n}$ is finite.
\item The series $\di\sum_{m=1}^{\infty}u_m$ and the series $\di\sum_{n=1}^{\infty}v_n$ both converge to $s$. Namely,
\[\sum_{m=1}^{\infty} \sum_{n=1}^{\infty} a_{m,n} =\sum_{n=1}^{\infty} \sum_{m=1}^{\infty} a_{m,n}=\sum_{(m,n)\in\mathbb{Z}^+\times\mathbb{Z}^+}a_{m,n}.\]
 \end{enumerate}
 
\end{theorem}
\begin{myproof}{Proof}
  Using the same notations as in the proof of Theorem \ref{230302_12}, since the double series $\di\sum_{(m,n)\in\mathbb{Z}^+\times\mathbb{Z}^+}|a_{m,n}|$ is convergent, the double series $\di\sum_{(m,n)\in\mathbb{Z}^+\times\mathbb{Z}^+}p_{m,n}$ and $\di\sum_{(m,n)\in\mathbb{Z}^+\times\mathbb{Z}^+}q_{m,n}$ are convergent. Applying Theorem \ref{230302_6} to the nonnegative series $\di\sum_{(m,n)\in\mathbb{Z}^+\times\mathbb{Z}^+}p_{m,n}$ and $\di\sum_{(m,n)\in\mathbb{Z}^+\times\mathbb{Z}^+}q_{m,n}$, we conclude that for all $m\in\mathbb{Z}^+$ and all $n\in\mathbb{Z}^+$, the series
$\di \sum_{n=1}^{\infty}p_{m,n}$, $\di \sum_{n=1}^{\infty}q_{m,n}$,  $\di\sum_{m=1}^{\infty}p_{m,n}$ and $\di\sum_{m=1}^{\infty}q_{m,n}$  are convergent.  Since
\[a_{m,n}=p_{m,n}-q_{m,n}\hspace{1cm}\text{for all}\;(m,n)\in\mathbb{Z}^+\times\mathbb{Z}^+,\]
we conclude that the series
$\di \sum_{n=1}^{\infty}a_{m,n}$ and the series $\di\sum_{m=1}^{\infty}a_{m,n}$ are convergent. The remaining assertions are concluded using the same arguments.
\end{myproof}
This theorem says that absolutely convergent double series enjoys almost the same privileges as the nonnegative double series.
 The following theorem gives a summary.
\begin{theorem}{}
Given that $\di\sum_{(m,n)\in\mathbb{Z}^+\times\mathbb{Z}^+}a_{m,n}$ is a double series that satisfies the following conditions.
\begin{enumerate}[(i)]
\item
For each fixed $m\in\mathbb{Z}^+$, the series $\di\sum_{n=1}^{\infty}|a_{m,n}|$ is convergent.
\item $\di \sum_{m=1}^{\infty}\sum_{n=1}^{\infty}|a_{m,n}|$ is convergent.
\end{enumerate}
 We have the following.
\begin{enumerate}[(a)]
\item The double series $\di\sum_{(m,n)\in\mathbb{Z}^+\times\mathbb{Z}^+}a_{m,n}$ converges absolutely.
\item For each fixed $n\in\mathbb{Z}^+$, the series $\di\sum_{m=1}^{\infty}a_{m,n}$  converges absolutely.
\item For each fixed $m\in\mathbb{Z}^+$, the series $\di\sum_{n=1}^{\infty}a_{m,n}$  converges absolutely.
\item Both the series  $\di \sum_{m=1}^{\infty}\left|\sum_{n=1}^{\infty} a_{m,n}\right|$ and $\di \sum_{n=1}^{\infty}\left|\sum_{m=1}^{\infty} a_{m,n}\right|$ are convergent.
\item The sum of the double series can be computed by iterated summation. Namely, 
\[\sum_{(m,n)\in\mathbb{Z}^+\times\mathbb{Z}^+}a_{m,n}=\sum_{m=1}^{\infty} \sum_{n=1}^{\infty} a_{m,n} =\sum_{n=1}^{\infty} \sum_{m=1}^{\infty} a_{m,n}.\]
\end{enumerate}
\end{theorem}
\begin{myproof}{Proof}
By Theorem \ref{230302_3}, (i) and (ii) implies that the double series  $\di\sum_{(m,n)\in\mathbb{Z}^+\times\mathbb{Z}^+}|a_{m,n}|$ is convergent, which gives (a). \bp By  Theorem \ref{230302_6}, (a) implies (b) and (c). Theorem \ref{230302_6} also implies that
the two series
\[  \sum_{m=1}^{\infty}\left(\sum_{n=1}^{\infty} |a_{m,n}|\right)\quad\text{and }\quad   \sum_{n=1}^{\infty}\left(\sum_{m=1}^{\infty} |a_{m,n}|\right)\]are convergent.
 Since
\[\left|\sum_{n=1}^{\infty} a_{m,n}\right|\leq \sum_{n=1}^{\infty}\left| a_{m,n}\right|,\hspace{1cm}\left|\sum_{m=1}^{\infty} a_{m,n}\right|\leq \sum_{m=1}^{\infty}\left|   a_{m,n}\right|,\]comparison test shows that  the series  $\di \sum_{m=1}^{\infty}\left|\sum_{n=1}^{\infty} a_{m,n}\right|$ and $\di \sum_{n=1}^{\infty}\left|\sum_{m=1}^{\infty} a_{m,n}\right|$ are convergent. This gives   (d). The statement (e) follows from (a) and Theorem \ref{230302_10}.
\end{myproof}

\vp
\noindent
{\bf \large Exercises  \thesection}
\setcounter{myquestion}{1}
 
\begin{question}{\themyquestion}
If $a$ and $b$ are positive constants,  show that the double series 
\[\sum_{(m,n)\in \mathbb{Z}^+\times\mathbb{Z}^+}\frac{1}{am^2+bn^2}\] is divergent.
\end{question}
\atc

\begin{question}{\themyquestion}
Given that $a$ and $b$ are positive constants,  $u$ and $v$ are   real numbers, and $\alpha$ is a number larger than 1. Show that the double series 
\[\sum_{(m,n)\in \mathbb{Z}^+\times\mathbb{Z}^+}\frac{\sin(mu+nv)}{(am^2+bn^2)^{\alpha}}\] is convergent.
\end{question}

\chapter{Sequences and Series of Functions}\label{ch6}

 In this chapter, we   study sequences and series whose terms depend on a variable. 
 
\section{Convergence of Sequences and Series of Functions}\label{sec6.1}
Let $D$ be a subset of real numbers. For each positive integer $n$, let $f_n:D\to\mathbb{R}$ be a function defined on $D$. Then $\{f_n\}_{n=1}^{\infty}$ is a sequence of functions defined on $D$. Sometimes we will write $\{f_n:D\to\mathbb{R}\}$ or  $\{f_n:D\to\mathbb{R}\}_{n=1}^{\infty}$ to make it explicit that each $f_n$ is a function defined on $D$. 

Given a sequence of functions $\{f_n:D\to\mathbb{R}\}$ that are defined on $D$, for each $x\in D$, $\{f_n(x)\}$ is a sequence of real numbers. We can determine whether such a sequence is convergent.

\begin{definition}{Pointwise Convergence of Sequence of Functions}
Given a sequence of functions $\{f_n:D\to\mathbb{R}\}$ that are defined on $D$, we say that it {\bf converges pointwise} to the function $f:D\to\mathbb{R}$ provided that  for every $x\in D$, the sequence $\{f_n(x)\}$ converges to $f(x)$. Namely,
\[f(x)=\lim_{n\to\infty}f_n(x)\hspace{1cm}\text{for all}\;x\in D.\]In this case, we also say that the function $f:D\to\mathbb{R}$ is the pointwise limit of the sequence of functions $\{f_n:D\to\mathbb{R}\}$.
\end{definition}
 
 Let us look at some examples.
\begin{example}[label=230303_1]{}
 For each positive integer $n$, let $f_n:[0,1]\to\mathbb{R}$ be the function $f_n(x)=x^n$. Study the pointwise convergence of the sequence of functions $\{f_n\}$.
\end{example}
\begin{solution}{Solution}Notice that
 \[\lim_{n\to\infty}x^n=\begin{cases}0,\quad &\text{if}\;0\leq x<1,\\1,\quad &\text{if}\;\quad x=1.\end{cases}\] 
Therefore, the sequence of functions $\{f_n\}$ converges pointwise to the function $f:[0,1]\to\mathbb{R}$, where
 \[f(x)=\begin{cases}0,\quad &\text{if}\;0\leq x<1,\\1,\quad &\text{if}\;\quad x=1.\end{cases}\]
\end{solution}

\begin{figure}[ht]
\centering
\includegraphics[scale=0.2]{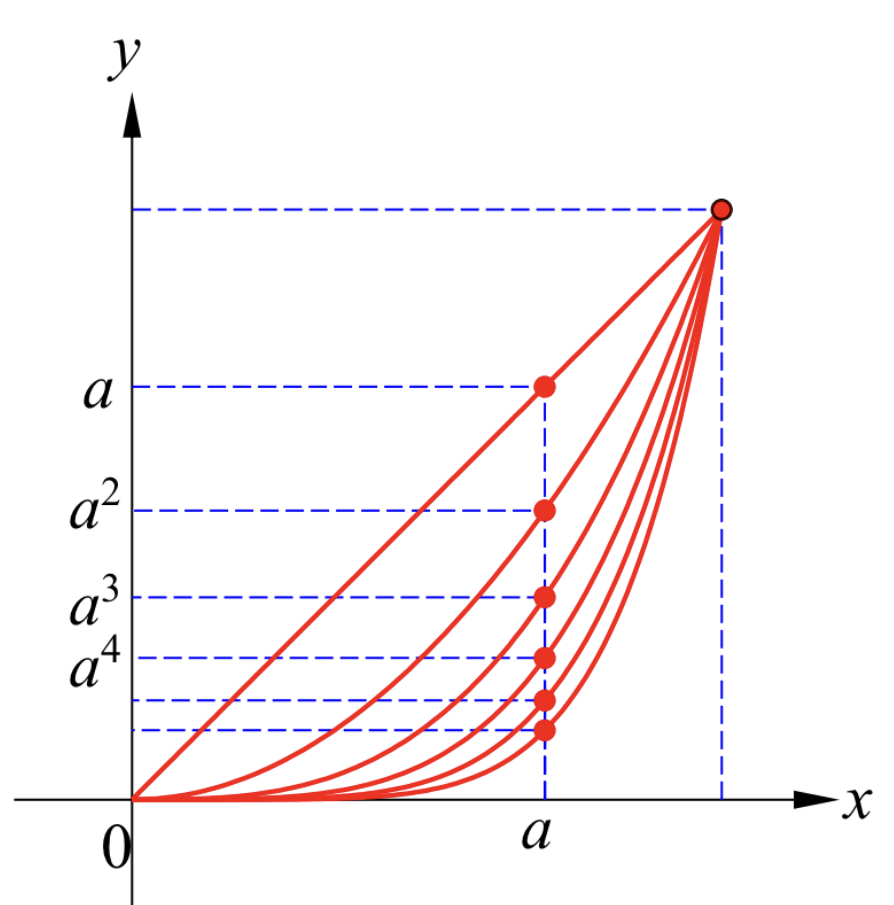}
\caption{The sequence of functions $\{f_n\}$ defined in Example \ref{230303_1}.\fa}\label{figure52}
\end{figure}

\begin{example}[label=230303_2]{}
 For each positive integer $n$, let $f_n:[0,2]\to\mathbb{R}$ be the function $f_n(x)=x^n$. Study the pointwise convergence of the sequence of functions $\{f_n\}$.\end{example}
\begin{solution}{Solution}
  For each $x\in [0,1]$, the sequence $\{f_n(x)\}$ is convergent. For any $x\in (1, 2]$, the sequence $\{f_n(x)\}$ is divergent. Hence, the sequence of functions $\{f_n\}$ does not converge pointwise.
\end{solution}
\begin{highlight}{}
In Example \ref{230303_1}, notice that each   $f_n:[0,1]\to\mathbb{R}$  is a  continuous function, but the limit  $f:[0,1]\to\mathbb{R}$ is not  a continuous function.

 Given that $\{f_n:D\to\mathbb{R}\}$ is a sequence of functions that converges pointwise to the function $f:D\to\mathbb{R}$.  We will consider the following questions.
\begin{enumerate}[1.]
\item If each $f_n$ is continuous, is $f$ continuous?
\item If each $f_n$ is a differentiable function defined on an open interval $I$, is $f$ differentiable on $I$? If yes, does the sequence $\{f_n':I\to\mathbb{R}\}$ converge to $f':I\to\mathbb{R}$?
\item If each $f_n$ is Riemann  integrable on a closed and bounded interval $I$, is $f$ Riemann integrable on $I$? If yes, does 
the sequence of integrals $\di\left\{\int_If_n \right\}$ converge to the integral $\di\int_I f$?
\end{enumerate}We have seen that the answer to the first question is no, as given by Example \ref{230303_1}. The answers to the second  and third questions are also no. We will look at some examples.
\end{highlight}

\begin{example}[label=230303_3]{}
 For each positive integer $n$, let $f_n:\mathbb{R}\to\mathbb{R}$ be the function \[f_n(x)=\frac{1}{1+nx^2}.\]
  Study the pointwise convergence of the sequence of functions $\{f_n\}$.

\end{example}
\begin{solution}{Solution}Since $f_n(0)=1$ for all $n\in\mathbb{Z}^+$,
\[
\lim_{n\to\infty}f_n(0)=1.
\]If $x\neq 0$, 
\[0\leq f_n(x)\leq \frac{1}{nx^2}.\] \bs
By squeeze theorem,
\[\lim_{n\to\infty}f_n(x)=0\hspace{1cm}\text{when}\;x\neq 0.\]
Hence, the sequence  of functions $\{f_n\}$ converges pointwise to the function $f:\mathbb{R}\to\mathbb{R}$, where
 \[f(x)=\begin{cases}0,\quad &\text{if}\;x\neq 0,\\1,\quad &\text{if}\; x=0.\end{cases}\]
\end{solution}

\begin{figure}[ht]
\centering
\includegraphics[scale=0.18]{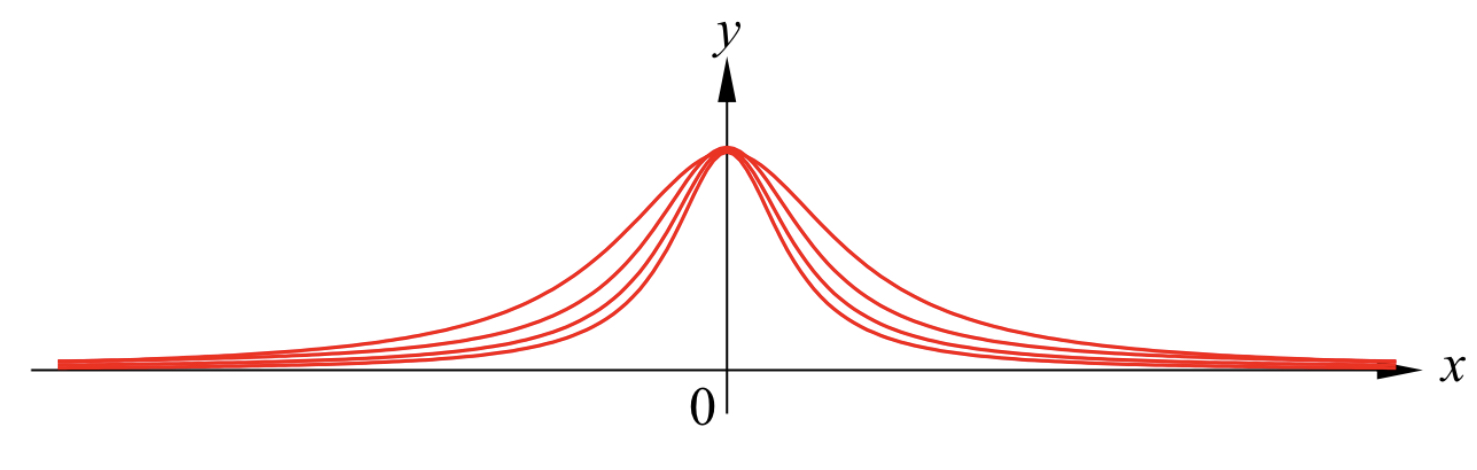}
\caption{The sequence of functions $\{f_n\}$ defined in Example \ref{230303_3}.\fa}\label{figure53}
\end{figure}

\begin{highlight}{}
In Example \ref{230303_3}, each of the functions $f_n$ is differentiable. But the function $f$ is not differentiable at $x=0$ since it is not continuous at $x=0$.

\end{highlight}

\begin{example}[label=230303_4]{}
 For each positive integer $n$, let $f_n:\mathbb{R}\to\mathbb{R}$ be the differentiable function \[f_n(x)=xe^{-nx^2}.\]
 \begin{enumerate}[(a)]
 \item  Study the pointwise convergence of the sequence of functions $\{f_n\}$.
 
 \item Study the pointwise convergence of the sequence of functions $\{f_n'\}$.
 \end{enumerate}
\end{example}
\begin{solution}{Solution}
\begin{enumerate} [(a)]
 
\item Since $f_n(0)=0$ for all $n\in\mathbb{Z}^+$,
\[
\lim_{n\to\infty}f_n(0)=0.
\]\end{enumerate}\bs \begin{enumerate}[]\item If $x\neq 0$, since $\di\lim_{u\to\infty}e^{-u}=0$, we find that
\[\lim_{n\to\infty}f_n(x)=x\lim_{n\to\infty}e^{-nx^2}=x\lim_{u\to \infty}e^{-u}=0.\]
Hence, the sequence of functions $\{f_n\}$ converges pointwise to the function $f:\mathbb{R}\to\mathbb{R}$, where $f(x)=0$ for all $x\in\mathbb{R}$.\end{enumerate}\begin{enumerate}[(b)]
\item For $n\in\mathbb{Z}^+$, 
\[f_n'(x)=(1-2nx^2)e^{-nx^2}.\]
Since $f_n'(0)=1$ for all $n\in\mathbb{Z}^+$,
\[
\lim_{n\to\infty}f_n'(0)=1.
\]If $x\neq 0$, since  $\di\lim_{u\to\infty} e^{-u}=0$ and $\di\lim_{u\to\infty}ue^{-u}=0$, we find that
\[\lim_{n\to\infty}f_n'(x)= \lim_{n\to\infty}(1-2nx^2)e^{-nx^2}= \lim_{u\to \infty}(1-2u)e^{-u}=0.\] 
Hence, the sequence of functions $\{f_n'\}$ converges pointwise to the function $g:\mathbb{R}\to\mathbb{R}$, where  \[g(x)=\begin{cases}0,\quad &\text{if}\;x\neq 0,\\1,\quad &\text{if}\; x=0.\end{cases}\]
\end{enumerate}
\end{solution}

\begin{figure}[ht]
\centering
\includegraphics[scale=0.18]{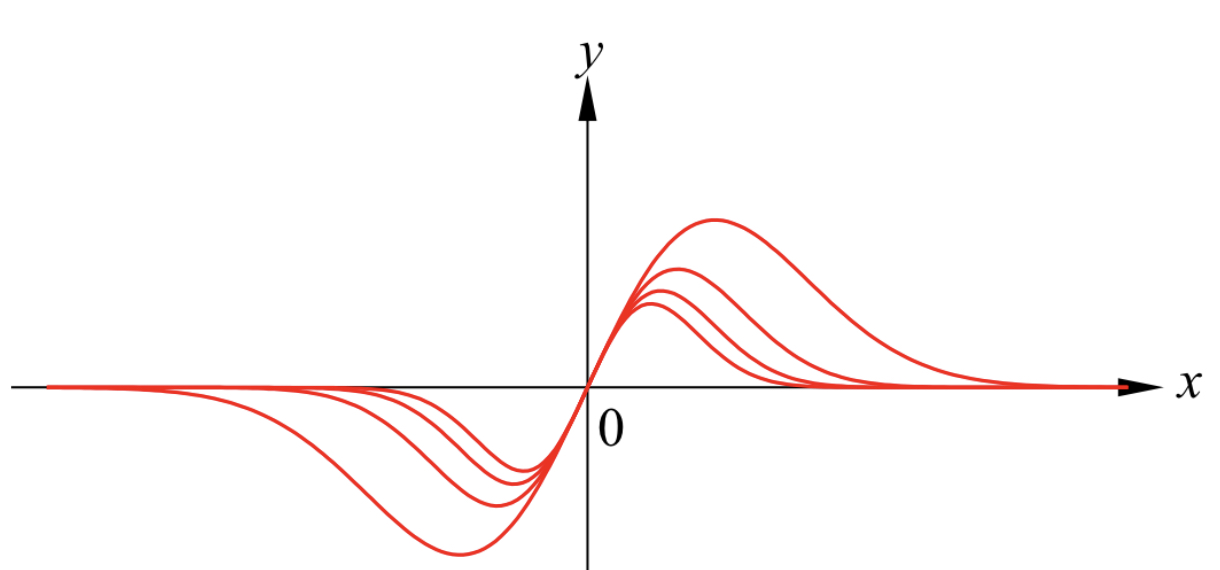}
\caption{The sequence of functions $\{f_n\}$ defined in Example \ref{230303_4}.\fa}\label{figure54}
\end{figure}

\begin{figure}[ht]
\centering
\includegraphics[scale=0.18]{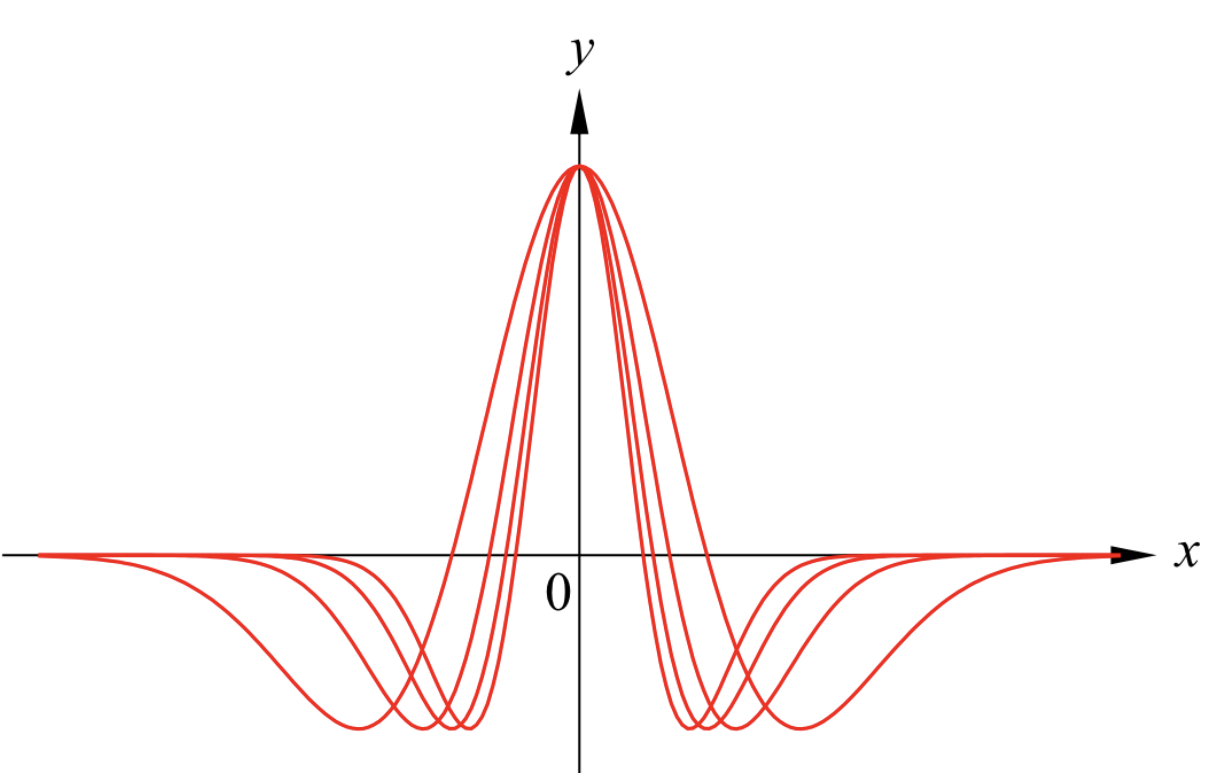}
\caption{The sequence of functions $\{f_n'\}$  in Example \ref{230303_4}.\fa}\label{figure55}
\end{figure}

\begin{highlight}{}
In Example \ref{230303_4}, each of the functions $f_n$ is differentiable and the  function $f$ is also differentiable. The sequence $\{f_n'\}$ also converges pointwise, but it does not converge to the function $f'$.

\end{highlight}

\begin{example}[label=230303_5]{}
 For each positive integer $n$, let 
 \[S_n=\left\{\left.\frac{p}{q}\,\right|\, p, q\in\mathbb{Z},\,0\leq p\leq q\leq n, q\geq 1\right\}.\]
 Define the function $f_n:[0,1]\to\mathbb{R}$ by
 \[f_n(x)=\begin{cases} 1,\quad &\text{if}\; x\in S_n,\\0,\quad &\text{if}\; x\notin S_n.\end{cases}\]Study the pointwise convergence of the sequence of functions $\{f_n\}$.
\end{example}
\begin{solution}{Solution}
If $x$ is a rational number in $[0,1]$, there exists a nonnegative integer  $p$ and a positive integer $q$ such that $0\leq p\leq q$ and $x=p/q$. Therefore, $x\in S_n$ for all $n\geq q$. This implies that 
$f_n(x)=1$ for all $n\geq q$. Hence,
\[\lim_{n\to\infty}f_n(x)=1\hspace{1cm}\text{if $x$ is rational}.\]\bs 
If $x$ is not a rational number, then $x\notin S_n$ for any $n\in\mathbb{Z}^+$. Therefore, $f_n(x)=0$ for all $n\in\mathbb{Z}^+$. Hence,
\[\lim_{n\to\infty}f_n(x)=0\hspace{1cm}\text{if $x$ is irrational}.\]These show that the sequence of functions $\{f_n\}$ converges pointwise to the Dirichlet function $f:\mathbb{R}\to\mathbb{R}$,
\[f(x)=\begin{cases}1,\quad &\text{if $x$ is rational},\\
0,\quad & \text{if $x$ is irrational}.\end{cases}.\]
\end{solution}
\begin{highlight}{}
For each $n\in\mathbb{Z}^+$, the set $S_n$ which  $f_n(x)\neq 0$ is  finite. Thus the function $f_n:[0,1]\to\mathbb{R}$ is  Riemann integrable. But the Dirichlet function $f$ is not Riemann integrable.
\end{highlight}

\begin{example}[label=230303_6]{}
 For each positive integer $n$, let 
 \[f_n(x)=\begin{cases} n^2x(1-nx),\quad &\text{if}\;0\leq x\leq \frac{1}{n},\\0,\quad  &\text{otherwise}.\end{cases}\] Notice that $f_n$ is integrable on $[0,1]$. Let
$\di c_n=\int_0^1f_n(x)dx$. 
 \begin{enumerate}[(a)]
 \item Study the pointwise convergence of the sequence of functions $\{f_n\}$.
\item Determine the limit of the sequence $\{c_n\}$.
\end{enumerate}
\end{example}
\begin{solution}{Solution}
\begin{enumerate}[(a)]
\item  Since $f_n(0)=0$ for all $n\in\mathbb{Z}^+$,
\[
\lim_{n\to\infty}f_n(0)=0.
\]
\end{enumerate}\bs \begin{enumerate}[]\item
If $x>0$, there is  a positive integer $N$ so that $x>1/N$. This implies that $f_n(x)=0$ for all $n\geq N$. Hence, we also have 
\[
\lim_{n\to\infty}f_n(x)=0.
\]Thus, the sequence of functions $\{f_n\}$ converges pointwise to the function $f:\mathbb{R}\to\mathbb{R}$, where $f(x)=0$ for all $x\in\mathbb{R}$.\end{enumerate}\begin{enumerate}[(b)]
\item We compute $c_n$ directly. For $n\in\mathbb{Z}^+$,
\[c_n=n^2\left[\frac{x^2}{2}-\frac{nx^3}{3}\right]_0^{1/n}=\frac{1}{6}.\]
Hence, the sequence $\{c_n\}$ converges to $1/6$.
\end{enumerate}
\end{solution}
\begin{highlight}{}
In Example \ref{230303_6},   each of the functions $f_n$ is integrable and the  function $f$ is also integrable. However, $\di\left\{\int_0^1f_n\right\}$ does not converge to $\di\int_0^1f$. 

\end{highlight}
\begin{figure}[ht]
\centering
\includegraphics[scale=0.2]{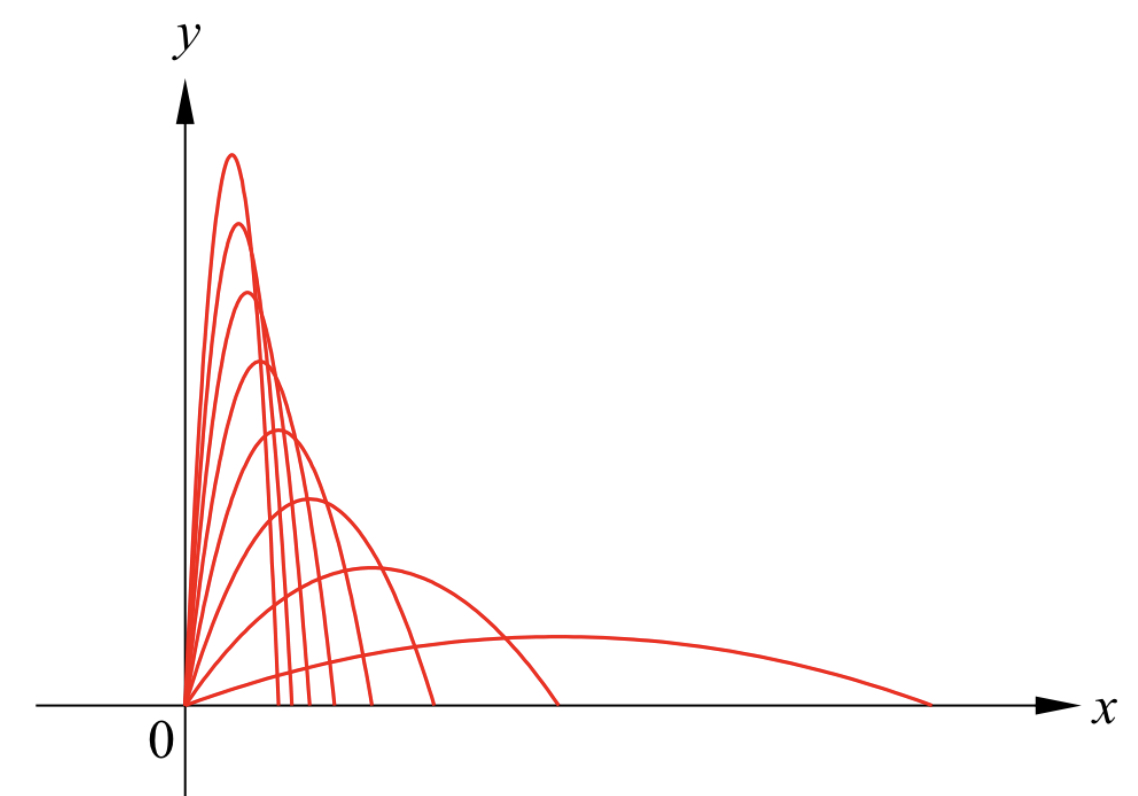}
\caption{The sequence of functions $\{f_n\}$ defined in Example \ref{230303_6}.\fa}\label{figure56}
\end{figure}

Now let us consider series of functions.
\begin{definition}{Pointwise Convergence of Series of Functions}
A series of functions defined on a set $A$ is a series of the form
\[\sum_{n=1}^{\infty}f_n(x),\]
where $\{f_n:A\to\mathbb{R}\}$ is a sequence of functions defined on $A$. For such a series, we form the partial sum
\[s_n(x)=\sum_{k=1}^n f_k(x)\hspace{1cm}\text{for}\;n\geq 1.\]
Then $\{s_n:A\to\mathbb{R}\}$ is a sequence of functons defined on $A$. The {\bf domain of convergence} of the series $\di\sum_{n=1}^{\infty}f_n(x)$ is the set
\[D=\left\{x\in A\,|\, \text{the sequence $\{s_n(x)\}$ is convergent}\right\}.\]
It is the largest subset $D$ of $A$ such that the sequence of functions
$\{s_n:D\to\mathbb{R}\}$ converges pointwise.  
For each $x$ in   $D$, let \[s(x)=\sum_{n=1}^{\infty}f_n(x)=\lim_{n\to\infty}s_n(x)\]
be the sum of the series. Then the  sequence of functions
$\{s_n:D\to\mathbb{R}\}$ converges pointwise to the function $s(x)$.
\end{definition}

 Let us reformulate Theorem \ref{230227_2} using series of functions.
\begin{example}[label=230305_16]{Geometric Series}
For the series $\di\sum_{n=0}^{\infty}x^n$, the terms are the functions $f_n(x)=x^n$, $n\geq 0$. They  are defined on $A=\mathbb{R}$. The partial sums are
\[s_n(x)=1+x+\cdots+x^n=\begin{cases}\di\frac{1-x^{n+1}}{1-x},\quad &\text{if}\;x\neq 1,\\n+1,\quad &\text{if}\;x=1.\end{cases}\]\be
The domain of convergence is the set $D=(-1,1)$. For $x\in D$, 
\[s(x)=\sum_{n=0}^{\infty}x^n=\frac{1}{1-x}.\]
Hence, the series of functons $\di\sum_{n=0}^{\infty}x^n$ converges pointwise on the interval $(-1,1)$ to the function $\di s(x)=\frac{1}{1-x}$.
\end{example2}

\begin{example}[label=230304_4]{}
Determine the domain of convergence of the series $\di\sum_{n=1}^{\infty}e^{-n^2 x}$.
\end{example}
\begin{solution}{Solution}
If $x\leq 0$, $-n^2x\geq 0$ for all $n\in\mathbb{Z}^+$. Therefore, $\di\lim_{n\to\infty}e^{-n^2 x}\neq 0$, and so the series $\di\sum_{n=1}^{\infty}e^{-n^2x}$ is divergent.

If $x>0$,  $\di\lim_{n\to\infty}e^{-n^2 x}=0$. In this case, notice that $n^2\geq n$  for all $n\in\mathbb{Z}^+$ implies that
\[0\leq e^{-n^2x}\leq e^{-nx}\hspace{1cm}\text{for all}\;n\in\mathbb{Z}^+.\]
Since the series $\di\sum_{n=1}^{\infty}e^{-nx}$ is a geometric series with positive constant ratio $r=e^{-x}<1$, it is convergent. By the comparison test, the series  $\di\sum_{n=1}^{\infty}e^{-n^2x}$ is also convergent.

 Therefore, the domain of convergence of the series  $\di\sum_{n=1}^{\infty}e^{-n^2 x}$ is $D=(0, \infty)$.
\end{solution}
\vp

\noindent
{\bf \large Exercises  \thesection}
\setcounter{myquestion}{1}
\begin{question}{\themyquestion}
 For each positive integer $n$, let $f_n:\mathbb{R}\to\mathbb{R}$ be the function defined by
 \[f_n(x)=e^{-nx^2}.\]  
   \begin{enumerate}[(a)]
   \item
  Determine the pointwise convergence of the sequence  of functions $\{f_n\}$.
 \item   Determine the pointwise convergence of the sequence of functions  $\{f_n'\}$.
 \end{enumerate}
\end{question}

\atc
\begin{question}{\themyquestion}
 For each positive integer $n$, let $f_n:(0,\infty)\to\mathbb{R}$ be the function defined by
 \[f_n(x)=\frac{1}{1+x^n}.\]  
  \begin{enumerate}[(a)]
   \item
  Determine the pointwise convergence of the sequence  of functions $\{f_n\}$.
 \item   Determine the pointwise convergence of the sequence of functions  $\{f_n'\}$.
 \end{enumerate}
\end{question}

\atc
\begin{question}{\themyquestion}
For each positive integer $n$, let $f_n:[0,\infty)\to\mathbb{R}$ be the function defined by
 \[f_n(x)=\frac{x}{1+nx},\]and  let
$\di c_n=\int_0^1f_n(x)dx$. 
\begin{enumerate}[(a)]
   \item  Determine the pointwise  convergence of the sequence of functions $\{f_n\}$.
   \item Determine the  convergence of the sequence $\{c_n\}$.
 \end{enumerate}
\end{question}

\atc
\begin{question}{\themyquestion}
For each positive integer $n$, let $f_n:\mathbb{R}\to\mathbb{R}$ be the function defined by
 \[f_n(x)=n\sin\left(\frac{x}{n}\right),\]  
 and let
$\di c_n=\int_0^1f_n(x)dx$.
\begin{enumerate}[(a)]
 \item  Study the  convergence of the sequence  of functions $\{f_n\}$.
\item  Study the  convergence of the sequence  of functions  $\{f_n'\}$.
\item Determine the limit of the sequence $\{c_n\}$.
 \end{enumerate}
\end{question}

 \atc
\begin{question}{\themyquestion}Find the domain of convergence of the series of functions $\di\sum_{n=1}^{\infty}ne^{-nx^2}$.
\end{question}
  \atc
\begin{question}{\themyquestion}Find the domain of convergence of the series of functions $\di\sum_{n=1}^{\infty}ne^{-n^2x}$.
\end{question}
 
\vp

\section{Uniform Convergence of Sequences and Series of Functions}\label{sec6.2}
In Section \ref{sec6.1}, we have seen  examples where a sequence of functions $\{f_n:D\to\mathbb{R}\}$ converges pointwise to a function $f:D\to\mathbb{R}$, but some properties of the sequence $\{f_n\}$, such as continuity, differentiability, or integrability, are lost in the limit function $f$. We also see an example where differentiability is preserved, but the derivative of $f$ is not the limit of the derivatives of the sequence $\{f_n\}$. There is  also an example where each function $f_n$ is integrable over an interval $I$, $f$ is also integrable over $I$, but the limit of the sequence $\di\int_If_n$ is not $\di\int_If$. 

\begin{highlight}{}
Given that $\{f_n:D\to\mathbb{R}\}$ is a sequence of functions  that converges pointwise to the function $f:D\to\mathbb{R}$.
\begin{enumerate}[1.]
\item  Continuity of each $f_n$ does not imply the continuity of $f$.
\[\lim_{x\to x_0}\lim_{n\to\infty}f_n(x) \quad\text{does not necessary equal to}\quad 
\lim_{n\to\infty}\lim_{x\to x_0}f_n(x).\]\end{enumerate} 
\begin{enumerate}[2.]
\item The derivative of $f$ does not necessary equal to the limit of $\{f_n'\}$
\[\frac{d}{dx}\lim_{n\to\infty}f_n(x) \quad\text{does not necessary equal to}\quad \lim_{n\to\infty}\frac{d}{dx} f_n(x).\] \end{enumerate}\begin{enumerate}[3.]
\item  The integral of $f$ over an interval $[a,b]$ does not necessary equal to the limit the  integrals of $f_n$ over $[a,b]$.
\[\int_a^b\lim_{n\to\infty}f_n(x)dx\quad\text{does not necessary equal to}\quad\lim_{n\to\infty}\int_a^b f_n(x) dx .\]
\end{enumerate}
Since derivatives and integrals are also limits, all these pathological behaviors have the same root. Namely, one cannot simply interchange the orders of two limits, as have been shown in Section \ref{sec5.5}. 
\end{highlight}

In this section, we are going to introduce the concept of uniform convergence. We are going to see in next section  how this extra condition can help to remedy some of the pathological behaviors mentioned above.

Let us   review Example \ref{230303_1}. The sequence $f_n:[0,1]\to\mathbb{R}$, $f_n(x)=x^n$ is found to converge pointwise to the function $f:[0,1]\to\mathbb{R}$ given by  \[f(x)=\begin{cases}0,\quad &\text{if}\;0\leq x<1,\\1,\quad &\text{if}\;\quad x=1.\end{cases}\]For the point $x=1$,  $\{f_n(1)\}$ converges to $f(1)=1$. For any $\varepsilon>0$, we can take $N=1$. Then for all $n\geq N$, \[|f_n(1)-f(1)|=0<\varepsilon.\] The same goes for the point  $x=0$. For any other $x$ in the interval $(0,1)$, $\{f_n(x)\}$ converges to $f(x)=0$. Given $\varepsilon>0$, if $\varepsilon<1$, the smallest $N$ such that
\[|f_n(x)-f(x)|=x^n<\varepsilon\hspace{1cm}\text{for all}\;n\geq N\] is the smallest positive integer $N$ such  that 
\[N>\frac{\ln\varepsilon}{\ln x}.\]
One see that this number $N$ would become larger and larger when $x$ approaches 1. The idea of uniform convergence is to say that $N$ can be chosen to be independent of the point $x$ in the domain.

\begin{definition}{Uniform Convergence of Sequences of Functions}
Let $D$ be a subset of real numbers. We say that a sequence of functions $\{f_n:D\to\mathbb{R}\}$ converges uniformly to the function $f:D\to\mathbb{R}$, provided that for every $\varepsilon>0$, there is a positive integer $N$ such that for all $n\geq N$, and all $x\in D$, 
\[|f_n(x)-f(x)|<\varepsilon.\]
\end{definition}

\begin{figure}[ht]
\centering
\includegraphics[scale=0.2]{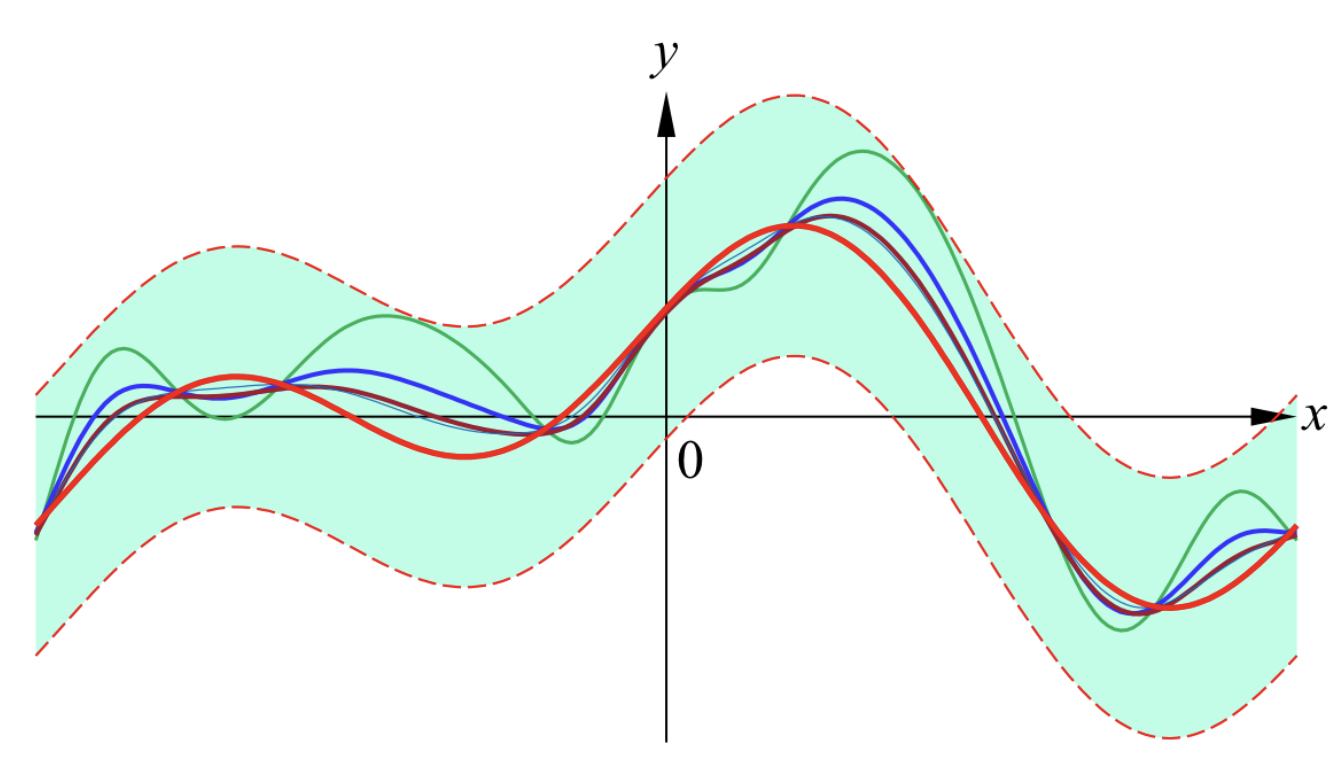}
\caption{Uniform convergence of a sequence of functions.\fa}\label{figure57}
\end{figure}

Obviously,  we have the following.
\begin{proposition}{}
Let $D$ be a subset of real numbers. If  $\{f_n:D\to\mathbb{R}\}$ is a sequence of functions that converges uniformly to the function $f:D\to\mathbb{R}$, then the sequence $\{f_n:D\to\mathbb{R}\}$  converges pointwise to  $f:D\to\mathbb{R}$.
\end{proposition}
\begin{highlight}{Uniform Limit and Pointwise Limit}
If a sequence of functions $\{f_n:D\to\mathbb{R}\}$  converges uniformly, the uniform limit is the same as the pointwise limit.
\end{highlight}

Let us compare the definitions of pointwise  and uniform convergence using logical expressions.
\begin{highlight}{Pointwise Convergence versus Uniform Convergence}
\begin{enumerate}[$\bullet$\;\;]
\item
The sequence of functions $\{f_n:D\to\mathbb{R}\}$ converges pointwise to the function $f:D\to\mathbb{R}$.
\[\forall\;x\in D,\; \forall \varepsilon>0,\;\exists N\in\mathbb{Z}^+,\;\forall \,n\geq N,\;|f_n(x)-f(x)|<\varepsilon.\]
\item
The sequence of functions $\{f_n:D\to\mathbb{R}\}$ converges uniformly to the function $f:D\to\mathbb{R}$.
\[\forall \varepsilon>0,\;\exists N\in\mathbb{Z}^+,\;\forall\;x\in D,\; \forall \,n\geq N,\;|f_n(x)-f(x)|<\varepsilon.\]
\end{enumerate}
\end{highlight}
One sees that it is a matter of the ordering of the quantifiers, but it makes a significant difference when we interchange the orders of a universal quantifier with  a existential quantifier.

 One should also compare the definition of uniform convergence to uniform continuity that we discussed in Section \ref{sec2.5}. In both cases, the uniformity is with respect to the domain $D$.

Before looking at some examples, let us highlight the negation of uniform continuity.
\begin{highlight}{Non-Uniform Convergence}
In logical expressions, the sequence of functions $\{f_n:D\to\mathbb{R}\}$ does not converge uniformly to the function $f:D\to\mathbb{R}$ is expressed by
\begin{equation}\label{eq230303_1}\exists \,\varepsilon>0,\;\forall N\in\mathbb{Z}^+,\;\exists\;x\in D,\; \exists \,n\geq N,\;|f_n(x)-f(x)|\geq \varepsilon.\end{equation}
\end{highlight}

The following gives a prelimary test for uniform convergence.
\begin{proposition}{}
If a sequence of functions $\{f_n:D\to\mathbb{R}\}$ does not converge pointwise, then it does not converge uniformly. 
\end{proposition}
If the sequence $\{f_n:D\to\mathbb{R}\}$ does converge pointwise, to show that it does not converge uniformly, we only need to establish the statement \eqref{eq230303_1} with $f:D\to\mathbb{R}$ the pointwise limit of the  sequence $\{f_n:D\to\mathbb{R}\}$.

\begin{example}[label=230303_7]{}For $n\geq 1$, let $f_n:[0,1]\to\mathbb{R}$ be the function $f_n(x)=x^n$.
Show that the sequence $\{f_n:[0,1]\to\mathbb{R} \}$ does not converge uniformly. 
\end{example}
\begin{solution}{Solution}In Example \ref{230303_1}, we have seen that the sequence $\{f_n \}$ converges  pointwise to the function $f:[0,1]\to\mathbb{R}$, where $f(x)=0$ for $x\in [0,1)$ and $f(1)=1$. If $\{f_n:[0,1]\to\mathbb{R}\}$  converges  uniformly, it must converge to the same function $f::[0,1]\to\mathbb{R}$. Take $\varepsilon=\frac{1}{2}$. There must be a positive integer $N$ such that for all $n\geq N$, for all $x\in [0,1]$,
\[|f_n(x)-f(x)|<\frac{1}{2}.\]\bs
In particular, this says that for all $x\in [0,1)$,
\[x^N<\frac{1}{2}.\]
This is absurd since $\di\lim_{x\to 1^-}x^N=1$. Hence, the  sequence $\{f_n:[0,1]\to\mathbb{R} \}$ does not converge uniformly. 
\end{solution}

\begin{example}[label=230303_8]{}For $n\geq 1$, let $f_n:\mathbb{R}\to\mathbb{R}$ be the function $f_n(x)=xe^{-nx^2}$.
Show that the sequence $\{f_n \}$  converges uniformly. 
\end{example}
\begin{solution}{Solution}In Example \ref{230303_4}, we have seen that the sequence $\{f_n\}$ converges pointwise to the function $f:\mathbb{R}\to\mathbb{R}$ that is identically 0.
Notice that 
\[f_n'(x)=(1-2nx^2)e^{-nx^2}.\]
This shows that $f'_n(x)>0$ for $|x|<1/\sqrt{2n}$, and $f'_n(x)<0$ for $|x|>1/\sqrt{2n}$. Since
\[\lim_{n\to-\infty}f_n(x)=0\quad\text{and}\quad\lim_{n\to\infty}f_n(x)=0,\]we find that $f_n(x)$  decreases from $0$ to $f_n(-1/\sqrt{2n})$ when $x$ goes from $-\infty$ to $ -1/\sqrt{2n}$, $f_n(x)$ increases from $f_n(-1/\sqrt{2n})$  to $f_n(1/\sqrt{2n})$ when $x$ goes from $-1/\sqrt{2n}$ to $1/\sqrt{2n}$, and $f_n(x)$ decreases from $f_n(1/\sqrt{2n})$ to 0 when $x$ goes from $1/\sqrt{2n}$ to $\infty$. Hence, the minimum and maximum values of $f_n$ are  $f_n(-1/\sqrt{2n})$ and $f_n(1/\sqrt{2n})$ respectively.
This shows that
\[|f_n(x)|\leq f_n(1/\sqrt{2n})=\frac{1}{\sqrt{2n}}e^{-\frac{1}{2}}\leq \frac{1}{\sqrt{2n}}\hspace{1cm}\text{for all}\;x\in\mathbb{R}.\]

Given $\varepsilon>0$, there is a positive integer $N$ such that $1/\sqrt{2N}<\varepsilon$. For all $n\in N$, for all $x\in\mathbb{R}$, we find that
\bs
\[|f_n(x)-f(x)|=|f_n(x)|\leq \frac{1}{\sqrt{2n}}\leq\frac{1}{\sqrt{2N}}<\varepsilon.\]
This proves that  the sequence of functions $\{f_n\}$ converges uniformly to the function $f$ that is identically 0.

\end{solution}

 By definition, if a sequence of functions $\{f_n:D\to\mathbb{R}\}$ converges uniformly to the function $f:D\to\mathbb{R}$, then there is a positive integer $N_0$ such that for all $n\geq N_0$, 
\[|f_n(x)-f(x)|<1\hspace{1cm}\text{for all}\;x\in D.\]
This implies that for all $n\geq N_0$, the function $(f_n-f):D\to\mathbb{R}$ is bounded above, and thus $\di M_n=\sup_{x\in D}|f_n(x)-f(x)|$ exists. 

The following theorem provides a systematic way to determine whether a sequence of functions $\{f_n:D\to\mathbb{R}\}$ converges uniformly.
\begin{theorem}[label=230303_9]{}
Let $D$ be a subset of real numbers, and let $\{f_n:D\to\mathbb{R}\}$ be a sequence of functions defined on $D$.
\begin{enumerate}[I.]
\item
If the sequence of functions $\{f_n:D\to\mathbb{R}\}$ does not converge pointwise to a function, then it does not converge uniformly.
\item If the sequence of functions $\{f_n:D\to\mathbb{R}\}$  converges pointwise to a function $f:D\to\mathbb{R}$, for each $n\in\mathbb{Z}^+$, define the function $g_n:D\to\mathbb{R}$ by $g_n(x)=f_n(x)-f(x)$.
\begin{enumerate}[(a)]
\item
If $g_n$ is not bounded for infinitely many $n$, then  the sequence of functions $\{f_n:D\to\mathbb{R}\}$ does not converge uniformly.
\item If only finitely many of the functions $g_n$ are not bounded, there is a positive integer $N_0$ such that $g_n$ is bounded for all $n\geq N_0$. For $n\geq N_0$, let $\di M_n=\sup_{x\in D}|g_n(x)|$. Then the  sequence of functions $\{f_n:D\to\mathbb{R}\}$  converges uniformly to the function $f:D\to\mathbb{R}$ if and only if $\di\lim_{n\to\infty}M_n=0$.
\end{enumerate}
\end{enumerate}
\end{theorem}
\begin{myproof}{Proof}
We have addressed I. and II. (a). Let us now consider II. (b). If the    sequence of functions $\{f_n:D\to\mathbb{R}\}$  converges uniformly to the function $f:D\to\mathbb{R}$, given $\varepsilon>0$, there is a positive integer $N\geq N_0$ such that for all $n\geq N$ and for all $x\in D$,
\[|g_n(x)|=|f_n(x)-f(x)|<\frac{\varepsilon}{2}.\]
This gives 
\[0\leq M_n=\sup_{x\in D}|g_n(x)|\leq\frac{\varepsilon}{2}<\varepsilon\hspace{1cm}\text{for all}\;n\geq N.\]
Therefore, $\di\lim_{n\to\infty}M_n=0$.
 
Conversely, if $\di\lim_{n\to\infty}M_n=0$, given $\varepsilon>0$, there is a positive integer $N\geq N_0$ such that 
\[M_n<\varepsilon\hspace{1cm}\text{for all}\;n\geq N.\]
It follows that for all $n\geq N$, for all $x\in D$,
\[|f_n(x)-f(x)|=|g_n(x)|\leq\sup_{x\in D}|g_n(x)|=M_n<\varepsilon.\]
This proves that the    sequence of functions $\{f_n:D\to\mathbb{R}\}$  converges uniformly to the function $f:D\to\mathbb{R}$.
\end{myproof}
\begin{example}{}
For  the sequence of functions discussed in Example \ref{230303_7}, \[f_n(x)-f(x)=\begin{cases}x^n,\quad &\text{if}\;0\leq x<1,\\0,\quad &\text{if}\;\quad x=1.\end{cases}\]Therefore,
\[M_n=\sup_{0\leq x\leq 1}|f_n(x)-f(x)|= 1.\]
Since $\di\lim_{n\to\infty}M_n=1\neq 0$,  Theorem \ref{230303_9} implies that the sequence of functions $\{f_n:[0,1]\to\mathbb{R}\}$ with $f_n(x)=x^n$ does not converge uniformly.
\end{example}

\begin{example}{}
For  the sequence of functions discussed in Example \ref{230303_8}, $f_n(x)-f(x)=f_n(x)=xe^{-nx^2}$. We have shown that
\[M_n=\sup_{x\in\mathbb{R}}|f_n(x)-f(x)|\leq\frac{ 1}{\sqrt{2n}}.\]
This implies that $\di\lim_{n\to\infty}M_n= 0$. Hence,  Theorem \ref{230303_9}  says that the sequence of functions $\{f_n:\mathbb{R}\to\mathbb{R}\}$ with $f_n(x)=xe^{-nx^2}$   converges uniformly.
\end{example}

To apply Theorem \ref{230303_9}, we need to know \emph{apriori} the pointwise limit of the sequence of functions $\{f_n\}$ to be able to conclude the uniform convergence of the sequence. Sometimes it could be difficult to find the limit function. To circumvent this problem, we introduce the concept of uniformly Cauchy.

\begin{definition}
{Uniformly Cauchy Sequence of Functions}Let $D$ be a subset of real numbers.
A sequence of functions $\{f_n:D\to\mathbb{R}\}$ is {\bf uniformly Cauchy} provided that for every $\varepsilon>0$, there is a positive integer $N$ such that for all $m\geq n\geq N$,
\[|f_m(x)-f_n(x)|<\varepsilon\hspace{1cm}\text{for all}\;x\in D.\]

\end{definition}
We have the following.
\begin{theorem}[label=230303_10]{~\\Cauchy Criterion for Uniform Convergence of Sequences of Functions}A sequence of functions $\{f_n:D\to\mathbb{R}\}$ converges uniformly if and only if it is uniformly Cauchy.
\end{theorem}
\begin{myproof}{Proof}
If the sequence of functions $\{f_n:D\to\mathbb{R}\}$ converges uniformly to $f:D\to\mathbb{R}$, given $\varepsilon>0$, there is a positive integer $N$ such that for all $n\geq N$,
\[|f_n(x)-f(x)|<\frac{\varepsilon}{2}\hspace{1cm}\text{for all}\;x\in D.\]
\bp
Using triangle inequality, this proves that for all $m\geq n\geq N$,
\[|f_m(x)-f_n(x)|<\varepsilon\hspace{1cm}\text{for all}\;x\in D.\]This proves that the sequence   $\{f_n:D\to\mathbb{R}\}$ is  uniformly Cauchy.

Conversely, if the sequence of functions $\{f_n:D\to\mathbb{R}\}$ is  uniformly Cauchy, then for each $x\in D$, the sequence $\{f_n(x)\}$ is a Cauchy sequence. Hence, it converges to a number $f(x)$. This shows that the sequence of functions $\{f_n:D\to\mathbb{R}\}$ converges pointwise to a function $f:D\to\mathbb{R}$. To show that the convergence is uniform, given $\varepsilon>0$, there is a positive integer $N$ such that for all $m\geq n\geq N$, 
\[|f_m(x)-f_n(x)|<\frac{\varepsilon}{2}\hspace{1cm}\text{for all}\;x\in D.\]   For each $x\in D$,  fixed $n\geq N$ and take the limit $m\to\infty$, we find that
\[|f_n(x)-f(x)|\leq\frac{\varepsilon}{2}.\]
This proves that for all $n\geq N$, for all $x\in D$, \[|f_n(x)-f(x)|<\varepsilon.\]Therefore, the sequence of functions $\{f_n:D\to\mathbb{R}\}$ converges uniformly.
\end{myproof}

Using Theorem \ref{230303_10},
Theorem \ref{230303_9} can be finetuned as follows. The proof is straightforward and we leave it to the exercises.
\begin{theorem}[label=230303_11]{}
Given that $\{f_n:D\to\mathbb{R}\}$ is a sequence of functions defined on the subset $D$ of real numbers, for each pair of $(m,n)\in\mathbb{Z}^+\times\mathbb{Z}^+$, define the extended real number $M_{m,n}$ as
\[M_{m,n}=\sup_{x\in D}|f_m(x)-f_n(x)|.\]Then the sequence of functions  $\{f_n:D\to\mathbb{R}\}$ converges uniformly if and only if the double sequence $\{M_{m,n}\}$ converges to 0.
\end{theorem}

Next we turn to series of functions.
\begin{definition}{Uniform Convergence of Series of Functions}
Let $D$ be a subset of real numbers and let $\{f_n:D\to\mathbb{R}\}$ be a sequence of functions defined on $D$. We say that  the series of functions $\di\sum_{n=1}^{\infty}f_n(x)$  converges uniformly to the function $s:D\to\mathbb{R}$ provided that the sequence of partial sums $\di\{s_n:D\rightarrow\mathbb{R}\}$ with $\di s_n(x)=\sum_{k=1}^nf_k(x)$ converges uniformly to the function $s(x)$.
 
\end{definition}

\begin{highlight}{Uniform Convergence of Series of Functions}A necessary condition for a series of functions to converge uniformly is that it converges pointwise.

When the series of functions  $\di\sum_{n=1}^{\infty}f_n(x)$ converges pointwise to a function $s(x)$ on a set $D$, then for any  $n\in \mathbb{Z}^+$ and any $x\in D$, the series
\[\sum_{k=n}^{\infty}f_k(x)\] converges pointwise to the function $s(x)-s_{n-1}(x)$, where  $\di s_n(x)=\sum_{k=1}^nf_k(x)$ is the $n^{\text{th}}$ partial sum, and $s_0(x)=0$ by default.

Therefore, we can reformulate the definition of uniform convergence of series of functions as follows.   The series $\di\sum_{n=1}^{\infty}f_n(x)$  converges uniformly to the function $s(x)$ on the set $D$ provided that for any $\varepsilon>0$, there is a positive integer $N$ such that for all $n\geq N$, 
\[\left|\sum_{k=n}^{\infty} f_k(x)\right|<\varepsilon\hspace{1cm}\text{for all}\;x\in D.\]
\end{highlight}

In most cases, such as Example \ref{230304_4}, we can only justify a series of functions converges pointwise, but we cannot find an explicit close form for the sum $s(x)$. In this case, a Cauchy criterion becomes useful.

From Theorem \ref{230303_10}, we obtain the following immediately. 
\begin{theorem}[label=230304_7]{~\\Cauchy Criterion for Uniform Convergence of Series of Functions}A series of functions $\di\sum_{n=1}^{\infty}f_n(x)$ converges uniformly on a set $D$ if and only if for every $\varepsilon>0$, there is a positive integer $N$ such that for all $m\geq n\geq N$,
\[\left|\sum_{k=n}^mf_k(x)\right|<\varepsilon\hspace{1cm}\text{for all}\;x\in D.\]
\end{theorem}

\begin{example}[label=230304_9]{}
For the series $\di\sum_{n=1}^{\infty}e^{-n^2x}$ considered in Example \ref{230304_4}, we have shown that it converges pointwise on the interval $(0,\infty)$. Let us prove that the convergence is uniform on any set $D$ of the form $D=[a, \infty)$, where $a$ is a positive constant.

We notice that if $m\geq n$ and $x\geq a$,
\[0\leq \sum_{k=n}^me^{-k^2x}\leq\sum_{k=n}^me^{-kx}\leq\sum_{k=n}^{\infty}e^{-kx}\leq \sum_{k=n}^{\infty}e^{-ka}= \frac{e^{-na}}{1-e^{-a}}.
\]
Given $\varepsilon>0$, since \[\lim_{n\to\infty}\frac{e^{-na}}{1-e^{-a}}=0,\]there exists a positive integer $N$ such that for all $n\geq N$,
\[0\leq \frac{e^{-na}}{1-e^{-a}}<\varepsilon.\]\be
It follows that for all $m\geq n\geq N$, and for all $x\in[a,\infty)$,
\[\left|\sum_{k=n}^me^{-k^2x}\right|\leq \frac{e^{-na}}{1-e^{-a}}<\varepsilon.\]
By Theorem \ref{230304_7}, the series $\di\sum_{n=1}^{\infty}e^{-n^2x}$ converges uniformly on $[a, \infty)$.
\end{example2}
The readers are invited  to show that the series $\di\sum_{n=1}^{\infty}e^{-n^2x}$ does not converge uniformly on the set $(0,\infty)$. 
It is  a typical situation that allthough the series converges pointwise on a set $A$,   it fails to converge uniformly on $A$, but it converges uniformly on  subsets of $A$. Most of the time, we do not need uniform convergence on $A$, but uniform convergence on a collection of subsets of $A$ whose union  is $A$ is enough. In the example above, $\mathscr{C}=\left\{[a, \infty)\,|\, a\geq 0\right\}$ is a collection of subsets of $A=(0,\infty)$ whose union is $A$. 

\begin{definition}{Absolute  Convergence of Series of Functions}
 A series of  functions $\di\sum_{n=1}^{\infty}f_n(x)$ is said to converge absolutely on a set $D$ if the series $\di\sum_{n=1}^{\infty}|f_n(x)|$ converges pointwise on $D$. In this case, the series $\di\sum_{n=1}^{\infty}f_n(x)$ also converges pointwise on $D$.  
\end{definition}

Now we   present a useful test to show that a series of functions converges absolutely and uniformly on a set $D$.
\begin{theorem}[label=230617_1]{}Let $\{f_n:D\to \mathbb{R}\}$ be a sequence of functions defined on $D$. If the series
 $\di\sum_{n=1}^{\infty}|f_n(x)|$ converges uniformly on $D$, then the series $\di\sum_{n=1}^{\infty}f_n(x)$ converges absolutely and  uniformly on $D$.
\end{theorem}
\begin{myproof}{Proof}
Since  the series
 $\di\sum_{n=1}^{\infty}|f_n(x)|$ converges uniformly on $D$, it also converges pointwise. Hence, the series $\di\sum_{n=1}^{\infty}f_n(x)$ converges absolutely on $D$.
 
Since   $\di\sum_{n=1}^{\infty}|f_n(x)|$ converges uniformly on $D$, it is uniformly Cauchy. Given $\varepsilon>0$, there is a positive integer $N$ such that for all $m\geq n\geq N$,
\[\sum_{k=n}^m|f_k(x)|<\varepsilon \hspace{1cm}\text{for all}\;x\in D.\]
Triangle inequality implies that
\[\left|\sum_{k=n}^mf_k(x)\right|\leq \sum_{k=n}^m|f_k(x)|<\varepsilon \hspace{1cm}\text{for all}\;x\in D.\]
 
Hence, the series $\di\sum_{n=1}^{\infty}f_n(x)$ is also uniformly Cauchy on $D$. Therefore, it also converges uniformly.
\end{myproof}
\begin{theorem}[label=230305_8]{Weiertrass M-Test}
Let $\{f_n:D\to \mathbb{R}\}$ be a sequence of functions defined on $D$. Assume that the following conditions are satisfied.
\begin{enumerate}[(i)]
\item
For each $n\in\mathbb{Z}^+$, there is a positive constant $M_n$ such that $|f_n(x)|\leq M_n$ for all $x\in D$.
\item The series $\di\sum_{n=1}^{\infty}M_n$ is convergent.
\end{enumerate}Then the series $\di\sum_{n=1}^{\infty}f_n(x)$ converges absolutely and uniformly on $D$. 
\end{theorem}
\begin{myproof}{Proof}By Theorem \ref{230617_1}, we only
 need to show that the series $\di\sum_{n=1}^{\infty}|f_n(x)|$ converges uniformly on $D$.
Given $\varepsilon>0$, since the series $\di\sum_{n=1}^{\infty}M_n$ is convergent, there is a positive integer $N$ such that for all $m\geq n\geq N$,
\[ \sum_{k=n}^m M_k<\varepsilon.\]This implies that
\[\sum_{k=n}^m|f_k(x)|\leq\sum_{k=n}^mM_k<\varepsilon\hspace{1cm}\text{for all}\;x\in D.\]
By Theorem \ref{230304_7}, the series  $\di\sum_{n=1}^{\infty}|f_n(x)|$ converges uniformly on $D$.
\end{myproof}

\begin{example}{}
Let $a$ be a positive number. Show that the series  $\di\sum_{n=1}^{\infty}(-1)^{n-1}e^{-n^2x}$ converges absolutely and  uniformly on $[a,\infty)$.  
\end{example}
\begin{solution}{Solution}
For $n\in\mathbb{Z}^+$, let $f_n(x)=\di (-1)^{n-1}e^{-n^2x}$. For $x\in [a, \infty)$, $x\geq a$. Hence, for $n\in\mathbb{Z}^+$,
\[|f_n(x)|=e^{-n^2x}\leq e^{-nx}\leq e^{-na}.\]
Since $r=e^{-na}<1$, the geometric series $\di\sum_{n=1}^{\infty} e^{-na}$ is convergent. By Weierstrass $M$-test, the series $\di\sum_{n=1}^{\infty}(-1)^{n-1}e^{-n^2x}$ converges absolutely and  uniformly on $[a,\infty)$. 
\end{solution}

\vp
\noindent
{\bf \large Exercises  \thesection}
\setcounter{myquestion}{1}
\begin{question}{\themyquestion}
For $n\geq 1$, let $f_n:[0,1]\to\mathbb{R}$ be the function $f_n(x)=e^{-nx^2}$.
Show that the sequence of functions $\{f_n \}$ does not converge uniformly. 
\end{question}
\atc
\begin{question}{\themyquestion}
For $n\geq 1$, let $f_n:\mathbb{R}\to\mathbb{R}$ be the function $\di f_n(x)=n\sin\left(\frac{x}{n}\right)$. Show that the sequence of functions $\{f_n \}$  does not converge uniformly. 
\end{question}

\atc
\begin{question}{\themyquestion}
For $n\geq 1$, let $f_n:[0,2\pi]\to\mathbb{R}$ be the function $\di f_n(x)=n\sin\left(\frac{x}{n}\right)$. Show that the sequence of functions $\{f_n \}$  converges uniformly. 
\end{question}
\atc
\begin{question}{\themyquestion}
For $n\geq 1$, let $f_n:[0,\infty)\to\mathbb{R}$ be the function $\di f_n(x)=\frac{x}{1+nx}$. Determine whether the sequence of functions $\{f_n \}$  converges uniformly. 
\end{question}
 
 \atc
\begin{question}{\themyquestion}
Let $a$ be a positive constant. Show that the series $\di\sum_{n=1}^{\infty}(-1)^{n-1}ne^{-n^2x}$ converges absolutely and uniformly on the set $[a, \infty)$.
\end{question}
\vp
\section{Properties of  Uniform Limits of Functions}\label{sec6.3}

In this section, we are going to see how uniform convergence can avoid the pathological behaviours we mentioned in the beginning of Section \ref{sec6.2}. First we show that uniform limit of continuous functions is continuous. This is a very important result in mathematical analysis.
\begin{theorem}[label=230304_1]{Uniform Limit of Continuous Functions is Continuous}
Given that $D$ is a subset of real numbers, and $\{f_n:D\to\mathbb{R}\}$ is a sequence of continuous functions that converges uniformly to the function $f:D\to\mathbb{R}$. Then the function  $f:D\to\mathbb{R}$ is continuous. 
\end{theorem}\begin{myproof}{Proof}
The proof is a standard $1/3$ argument. Given $x_0\in D$, we want to show that $f$ is continuous at $x_0$ using the $\varepsilon-\delta$ argument. Given $\varepsilon>0$, there is a positive integer $N$ such that for all $n\geq N$, 
\[|f_n(x)-f(x)|<\frac{\varepsilon}{3}\hspace{1cm}\text{for all}\;x\in D.\]We are only going to use this statement when $n=N$.
Since $f_N$ is continuous at $x_0$, there is a $\delta>0$ such that for all $x\in D$, if $|x-x_0|<\delta$, then
\[|f_N(x)-f_N(x_0)|<\frac{\varepsilon}{3}.\]

From these, we find that if $x$ is in $D$ and $|x-x_0|<\delta$, then
\begin{align*}
|f(x)-f(x_0)|&\leq |f(x)-f_N(x)|+|f_N(x)-f_N(x_0)|+|f_N(x_0)-f(x_0)|\\&<\frac{\varepsilon}{3}+\frac{\varepsilon}{3}+\frac{\varepsilon}{3}=\varepsilon.\end{align*}This proves that $f$ is continuous at $x_0$.
 
\end{myproof}
\begin{example}{}
 For the sequence of functions $\{f_n:[0,1]\to\mathbb{R}\}$ with $f_n(x)=x^n$, its pointwise limit $f:[0,1]\to\mathbb{R}$,  \[f(x)=\begin{cases}0,\quad &\text{if}\;x\neq 0,\\1,\quad &\text{if}\; x=0,\end{cases}\]
is  not continuous. Since  each $f_n$, $n\in\mathbb{Z}^+$ is a continuous function,
  Theorem \ref{230304_1} can be used to infer that the sequence of functions  $\{f_n:[0,1]\to\mathbb{R}\}$ with $f_n(x)=x^n$  does not converge uniformly.

\end{example}
Applying Theorem \ref{230304_1} to series of functions, we have the following.
\begin{corollary}[label=230305_10]{}
Given that $\{f_n:D\rightarrow \mathbb{R}\}$ is a sequence of continuous functions defined on $D$. If the series of functions $\di\sum_{n=1}^{\infty}f_n(x)$ converges uniformly on  $D$, then it defines a continuous function $s:D\to\mathbb{R}$ by
\[s(x)=\sum_{n=1}^{\infty}f_n(x).\]
\end{corollary}\begin{myproof}{Proof}
We apply  Theorem \ref{230304_1} to the sequence of partial sums $\{s_n(x)\}$. Since $\di s_n(x)=\di\sum_{k=1}^nf_k(x)$ is a finite sum of continuous functions, it is continuous. By Theorem \ref{230304_1}, $\di s(x)=\lim_{n\to\infty}s_n(x)$ is continuous.
\end{myproof}

Next, we turn to integration.
\begin{theorem}[label=230304_2]{}
 Assume that for  each $n\in\mathbb{Z}^+$, the funtion $f_n:[a,b]\to\mathbb{R}$ is Riemann integrable. If the sequence of functions $\{f_n:[a,b]\to\mathbb{R}\}$   converges uniformly to the function $f:[a,b]\to\mathbb{R}$, then $f:[a,b]\to\mathbb{R}$ is also Riemann integrable, and the orders of the limit operation and the integration operation can be interchanged. Namely,
\begin{equation}\label{eq230304_3}
\lim_{n\to\infty}\int_a^b f_n(x)dx=\int_a^bf(x)dx=\int_a^b\lim_{n\to\infty}f_n(x)dx.
\end{equation}
\end{theorem}
Notice that we only assume that each $f_n$ is Riemann integrable. We do not need to assume that it is continuous.
\begin{myproof}{Proof}
Given $\varepsilon>0$, there is a positive integer $N$ such that for all $n\geq N$, 
\begin{equation}\label{eq230304_5}|f_n(x)-f(x)|<\frac{\varepsilon}{3(b-a)}\hspace{1cm}\text{for all}\;x\in [a,b].\end{equation} First we take $n=N$. Since $f_N:[a,b]\to\mathbb{R}$ is Riemann integrable, there is a partition  $P=\{x_i\}_{i=0}^k$  of $[a,b]$ such that
\[U(f_N,P)-L(f_N,P)<\frac{\varepsilon}{3}.\]From \eqref{eq230304_5}, we have
\[f_N(x)-\frac{\varepsilon}{3(b-a)}<f(x)<f_N(x)+\frac{\varepsilon}{3(b-a)}\hspace{1cm}\text{for all}\;x\in [a,b].\]
  For any $1\leq i\leq k$, if $x\in [x_{i-1}, x_i]$,
\begin{align*}
\inf_{x_{i-1}\leq x\leq x_i}f_N(x)-\frac{\varepsilon}{3(b-a)}\leq f(x)\leq \sup_{x_{i-1}\leq x\leq x_i}f_N(x)+\frac{\varepsilon}{3(b-a)}.
\end{align*}This implies that
\begin{align*}\inf_{x_{i-1}\leq x\leq x_i}f_N(x)&-\frac{\varepsilon}{3(b-a)} \leq \inf_{x_{i-1}\leq x\leq x_i}f(x)\\ &\leq  \sup_{x_{i-1}\leq x\leq x_i}f(x)\leq \sup_{x_{i-1}\leq x\leq x_i}f_N(x)+\frac{\varepsilon}{3(b-a)}.\end{align*}
 \bp
Therefore,
\[U(f,P)\leq U(f_N,P)+\frac{\varepsilon}{3},\hspace{1cm}L(f,P)\geq L(f_N,P)-\frac{\varepsilon}{3}.\] 
Thus,
\[U(f,P)-L(f,P)\leq U(f_N,P)-L(f_N,P)+\frac{2\varepsilon}{3}<\frac{\varepsilon}{3}+\frac{2\varepsilon}{3}=\varepsilon.\]
From this, we conclude  that $f$ is Riemann integrable.
 This in turn implies that for any $n\in\mathbb{Z}^+$, the function $f_n-f$ is Riemann integrable on $[a,b]$, and so is the function $|f_n-f|$.
With the same $\varepsilon>0$, we find from \eqref{eq230304_5} that for any $n\geq N$,
\begin{align*}
\left|\int_a^bf_n(x)dx-\int_a^b f(x)dx\right|&=\left|\int_a^b(f_n(x)-f(x))dx\right|\\&\leq \int_a^b\left|f_n(x)-f(x) \right|dx\leq \frac{\varepsilon}{3}<\varepsilon.
\end{align*}This proves that \eqref{eq230304_3} holds.
 
\end{myproof}

Applying Theorem \ref{230304_2} to series of functions, we have the following.
\begin{corollary}[label=230305_11]{}
Given that $\{f_n:[a,b]\rightarrow \mathbb{R}\}$ is a sequence of Riemann integrable functions. If the series $\di\sum_{n=1}^{\infty}f_n(x)$ converges uniformly, then the function $\di s(x)=\sum_{n=1}^{\infty}f_n(x)$ is Riemann integrable, the series $\di \sum_{n=1}^{\infty}\int_a^bf_n(x)dx$ is convergent, and we can integrate term by term. Namely,
\[\int_a^bs(x)dx=\int_a^b \sum_{n=1}^{\infty}f_n(x)dx=\sum_{n=1}^{\infty}\int_a^bf_n(x)dx.\] 
\end{corollary}
\begin{myproof}{Proof}
We apply  Theorem \ref{230304_2} to the sequence of partial sums $\{s_n(x)\}$. Since $\di s_n(x)=\di\sum_{k=1}^nf_k(x)$ is a finite sum of Riemann integrable functions, it is Riemann integrable. The rest follows from Theorem \ref{230304_2}.
\end{myproof}
\begin{example}{}
In Example \ref{230303_6}, the sequence of functions $\{f_n:[0,1]\to\mathbb{R}\}$ converges pointwise to the function $f:[0,1]\to\mathbb{R}$ that is identically 0. However, since
$\di\int_0^1f_n(x)dx=1/6$, the sequence $\di\left\{\int_0^1f_n(x)dx\right\}$ does not converge to $\di\int_0^1f(x)dx=0$. Theorem \ref{230304_2} can be used to deduce that $\{f_n:[0,1]\to\mathbb{R}\}$ does not converge to $f:[0,1]\to\mathbb{R}$ uniformly. 

In fact, one can verify that 
\[M_n=\sup_{0\leq x\leq 1}|f_n(x)-f(x)|=\sup_{0\leq x\leq 1/n}n^2x\left(1-nx\right)=\frac{n}{4}.\]Since $\di\lim_{n\to\infty}M_n\neq 0$,    $\{f_n:[0,1]\to\mathbb{R}\}$ does not converge to $f:[0,1]\to\mathbb{R}$ uniformly. 
\end{example}

Now we consider differentiation. In Example \ref{230303_8}, we have shown that the sequence $\{f_n:\mathbb{R}\to\mathbb{R}\}$ defined by $f_n(x)=xe^{-nx^2}$ converges uniformly to the function $f:\mathbb{R}\to\mathbb{R}$ that is identically zero. In Example \ref{230303_4}, we have seen that the derivative sequence $\{f_n'\}$ converges to the function $g:\mathbb{R}\to\mathbb{R}$ given  by
\[g(x)=\begin{cases}0,\quad &\text{if}\;x\neq 0,\\1,\quad &\text{if}\; x=0,\end{cases}\]We  find that
\[\left.\frac{d}{dx} \right|_{x=0}\lim_{n\to\infty} f_n(x)=0 \quad\text{which is not equal to}\quad\lim_{n\to\infty} \left.\frac{d}{dx} \right|_{x=0}f_n(x)=1.\]
Hence, even though the seqeunce of functions $\{f_n\}$ converges uniformly, we cannot interchange limit with differentiation.

 The following theorem gives a sufficient condition for interchanging limit with differentiation.

\begin{theorem}[label=230304_8]{}
Given that   $\{f_n:(a,b)\to\mathbb{R}\}$ is a sequence of   functions which satisfies the following conditions.
\begin{enumerate}[(i)]
\item There is a point $x_0$ in the interval $(a,b)$ such that the sequence   $\{f_n(x_0)\}$ converges to a number   $y_0$.
\item For each $n\in\mathbb{Z}^+$, $f_n:(a,b)\to\mathbb{R}$ is   differentiable.
\item The sequence of derivative functions $\{f_n':(a,b)\to\mathbb{R}\}$ converges uniformly to a function $g:(a,b)\to\mathbb{R}$.

\end{enumerate}
Then we have the following.
\begin{enumerate}[(a)]
\item The sequence of functions $\{f_n:(a,b)\to\mathbb{R}\}$ converges uniformly to a function $f:(a,b)\to\mathbb{R}$.
\item The function $f:(a,b)\to\mathbb{R}$ is differentiable.
\item  We can interchange differentiation and limits. Namely, for any $x\in (a,b)$, \[f'(x)= \frac{d}{dx}\lim_{n\to\infty}f_n(x)=\lim_{n\to\infty}\frac{d}{dx}f_n(x)= g(x).\]
\end{enumerate}
\end{theorem}

\begin{myproof}{Proof}For each $n\in\mathbb{Z}^+$,  since $f_n:I\to\mathbb{R}$ is   differentiable, it is continuous.
Given a point $c$ in the interval $(a,b)$, let $\{h_{n,c}:(a,b)\to\mathbb{R}\}$ be a sequence of functions defined by
\begin{equation}\label{eq230304_13}h_{n,c}(x)=\begin{cases} \di\frac{f_n(x)-f_n(c)}{x-c},\quad &\text{if}\;x\neq c,\\f_n'(c),\quad &\text{if}\;x=c.\end{cases}\end{equation}\bp
Then $h_{n,c}:(a,b)\to\mathbb{R}$ is a continuous function. 
For any positive integers $m$ and $n$, we have
\[h_{m,c}(x)-h_{n,c}(x)=\begin{cases} \di\frac{(f_m(x)-f_n(x))-(f_m(c)-f_n(c))}{x-c},\quad &\text{if}\;x\neq c,\\f_m'(c)-f_n'(c),\quad &\text{if}\;x=c.\end{cases}\]
Applying  mean value theorem to the differentiable function $f_m(x)-f_n(x)$, we find that  for any $x\in (a,b)\setminus\{c\}$, there is a point $\xi_x$ in between $x$ and $c$ such that
\[h_{n,c}(x)= f_{m}'(\xi_x)-f_n'(\xi_{x}).\]

Thus, we find that for any $x\in (a,b)$,
\[|h_{m,c}(x)-h_{n,c}(x)|\leq \sup_{a<x<b}|f_m'(x)-f_n'(x)|.\]
This implies that
\begin{equation}\label{eq230304_10}\sup_{a<x<b}|h_{m,c}(x)-h_{n,c}(x)|\leq \sup_{a<x<b}|f_m'(x)-f_n'(x)|.\end{equation}
Since the sequence of functions $\{f_n'\}$ converges uniformly, Theorem \ref{230303_11} implies that
\[\lim_{m,n\to\infty} \sup_{a<x<b}|f_m'(x)-f_n'(x)|=0.\]
Eq. \eqref{eq230304_10} then implies that 
\[\lim_{m,n\to\infty} \sup_{a<x<b}|h_{m,c}(x)-h_{n,c}(x)|=0.\]
By Theorem  \ref{230303_11} again, we find that the sequence of functions $\{h_{n,c}:(a,b)\to\mathbb{R}\}$ converges uniformly.

Now we specialize to $c=x_0$. Notice that by definition, 
\begin{equation}
\label{eq230304_12} f_m(x)-f_n(x)=f_m(x_0)-f_n(x_0)+(x-x_0)\left(h_{m,x_0}(x)-h_{n,x_0}(x)\right).\end{equation}Given $\varepsilon>0$, since the sequence $\{f_n(x_0)\}$ is convergent, there is a positive integer $N_1$ such that for all $m\geq n\geq N_1$, 
\[|f_m(x_0)-f_n(x_0)|<\frac{\varepsilon}{2}.\]\bp
Since the sequence of functions $\{h_{n,x_0}(x)\}$ converges uniformly,   Theorem \ref{230303_10} implies that there is a positive integer $N\geq N_1$ such that for all $m\geq n\geq N$,
\[|h_{m,x_0}(x)-h_{n,x_0}(x)|<\frac{\varepsilon}{2(b-a)}\hspace{1cm}\text{for all}\;x\in (a,b).\]
Eq. \eqref{eq230304_12}   implies that for all $m\geq n\geq N$, and for all $x\in (a,b)$,
\begin{align*}
|f_m(x)-f_n(x)|&\leq |f_m(x_0)-f_n(x_0)|+|x-x_0||h_{m,x_0}(x)-h_{n,x_0}(x)|
\\
&<\frac{\varepsilon}{2}+(b-a)\times \frac{\varepsilon}{2(b-a)}=\varepsilon.
\end{align*}

By Theorem \ref{230303_10}, this proves that the sequence of functions $\{f_n:(a,b)\to\mathbb{R}\}$ converges uniformly.   Let
\[f(x)=\lim_{n\to \infty}f_n(x)\] be the limit function. Being the limit of a sequence of continuous functions that converges uniformly, Theorem \ref{230304_1} says hat the function $f:(a,b)\to\mathbb{R}$ is continuous. 

 Now we want to prove that $f$ is differentiable and $f'(x)=g(x)$ for each $x\in (a,b)$. 
For any fixed $c\in (a,b)$, since the sequence of continuous functions $\{h_{n,c}(x)\}$ converges uniformly, it also converges pointwise. Taking $n\to\infty$ limits  in \eqref{eq230304_13}, we find that
\[h_c(x)=\lim_{n\to \infty}h_{n,c}(x)=\begin{cases}\di\frac{f(x)-f(c)}{x-c}, \quad &\text{if}\;x\neq c,\\g(c),\quad &\text{if}\;x=c.\end{cases}.\]

Since $\{h_{n,c}(x)\}$ is a sequence of continuous functions that converges uniformly,  Theorem \ref{230304_1} says that the limit function $h_c:(a,b)\to\mathbb{R}$ is continuous. Therefore,
\[g(c)=h_c(c)=\lim_{x\to c}h_c(x)=\lim_{x\to c}\frac{f(x)-f(c)}{x-c}.\]
This shows that $f$ is differentiable at $x=c$ and $f'(c)=g(c)$.
\end{myproof}

\begin{remark}{}
 \begin{enumerate}[1.]
 \item 
 In   Theorem \ref{230304_8}, we do not need to assume that the sequence of functions $\{f_n\}$ converges uniformly. It is a consequence of uniform convergence of the sequence $\{f_n'\}$. The condition that there is a point $x_0$ in $(a,b)$ so that  the sequence   $\{f_n(x_0)\}$ converges is  necessary. For otherwise  if we let $\widetilde{f}_n(x)=f_n(x)+n$ for  $n\in\mathbb{Z}^+$, then $\widetilde{f}_n'(x)=f_n'(x) $. But the sequence $\{\widetilde{f}_n\}$ does not converge if the sequence $\{f_n\}$  is convergent.
 \item
If we assume that  for all $n\in\mathbb{Z}^+$, the function  $f_n:(a,b)\to\mathbb{R}$ is continuously differentiable, there is an easier proof for the conclusions in Theorem \ref{230304_8}. \end{enumerate}
\end{remark}

Applying Theorem \ref{230304_8} to series of functions, we have the following.
\begin{corollary}[label=230304_14]
{}Let $\{f_n:(a,b)\to\mathbb{R}\}$ be a sequence of differentiable functions. Assume that there is a  $x_0\in (a, b)$ such that the series $\di\sum_{n=1}^{\infty}f_n(x_0)$ is convergent, and the series $\di\sum_{n=1}^{\infty}f_n'(x)$ converges uniformly on $(a,b)$, then the series  $\di\sum_{n=1}^{\infty}f_n(x)$ converges uniformly on $(a,b)$ to a differentiable function whose derivative is given by
\[\frac{d}{dx}\sum_{n=1}^{\infty}f_n(x)=\sum_{n=1}^{\infty}f_n'(x)\hspace{1cm}\text{for all}\;x\in (a,b).\]
\end{corollary}
\begin{myproof}{Proof}
Applying Theorem \ref{230304_8} to the sequence of partial sums $\{s_n(x)\}$. Since $\di s_n(x)=\di\sum_{k=1}^nf_k(x)$ is a finite sum of differentiable functions, it is differentiable. The rest follows from Theorem \ref{230304_8}.
\end{myproof}

\begin{example}{}
Consider the series $\di\sum_{n=1}^{\infty} e^{-n^2x}$ discussed in Example \ref{230304_9}. We have shown that it converges uniformly on $[a,\infty)$ when $a$ is a positive number.  For each $n\in\mathbb{Z}^+$, $f_n(x)=\di-\frac{1}{n^2}e^{-n^2 x}$ is a differentiable function with derivative 
\[\frac{d}{dx}\left(-\frac{1}{n^2}e^{-n^2 x}\right)=e^{-n^2 x}.\]
Notice that  
\[|f_n(x)|\leq \frac{1}{n^2}\hspace{1cm}\text{for all}\;x\in [0,\infty).\]
 
Since the series $\di\sum_{n=1}^{\infty}\frac{1}{n^2}$ is convergent, 
 Weierstrass $M$-test shows that the series $\di \sum_{n=1}^{\infty} f_n(x)=\sum_{n=1}^{\infty}-\frac{1}{n^2}e^{-n^2x}$ converges absolutely and uniformly on $[0, \infty)$. Corollary \ref{230304_14} shows that for any $x\in (a, \infty)$, we can do term by term differentiation and obtain
 \begin{equation}\label{eq230304_11}\frac{d}{dx}\sum_{n=1}^{\infty}-\frac{1}{n^2}e^{-n^2x}=\sum_{n=1}^{\infty} e^{-n^2x}.\end{equation}
 Since $a>0$ is arbitrary, eq. \eqref{eq230304_11} holds for any $x>0$. However, this is not true for $x=0$ even if we only consider right derivatives, as the right hand side of the equation is divergent when $x=0$.
\end{example}

\vp
\noindent
{\bf \large Exercises  \thesection}
\setcounter{myquestion}{1}
\begin{question}{\themyquestion}
\begin{enumerate}[(a)]\item Show that 
the series $\di\sum_{n=1}^{\infty} n^2e^{-n^2x}$ defines a continuous function on $(0,\infty)$.
\item 
Show that 
the series $\di\sum_{n=1}^{\infty}  e^{-n^2x}$ defines a differentiable function on $(0,\infty)$, and for each $x>0$,
\[\frac{d}{dx}\sum_{n=1}^{\infty}  e^{-n^2x}=-\sum_{n=1}^{\infty} n^2e^{-n^2x}.\]\end{enumerate}
\end{question}
\atc
\begin{question}{\themyquestion}
 Let $\{f_n:(a,b)\to\mathbb{R}\}$ be a sequence of continuously differentiable  functions. Assume that there is a point $x_0\in [a,b]$ such that the sequence $\{f_n(x_0)\}$ converges to a point $y_0$, and the sequence of functions $\{f_n':(a,b)\to\mathbb{R}\}$ converges uniformly to a function $g:(a,b)\to\mathbb{R}$. Use the fundamental theorems of calculus and Theorem \ref{230304_2} to prove that the sequence of functions $\{f_n:(a,b)\to\mathbb{R}\}$ converges uniformly to a differentiable function $f:(a,b)\to\mathbb{R}$, and $f'(x)=g(x)$ for all $x\in (a,b)$.
\end{question}
 
\vp
\section{Power Series}\label{sec6.4}
In this section, we turn to consider a special class of series of functions  called {\it power series}. The partial sums of a power series are polynomial functions. Hence, power series are limits of polynomial sequences. They play important roles in analysis.
\begin{definition}{Power series}A power series in the variable $x$ is a series of the form
\[\sum_{n=0}^{\infty}c_n(x-x_0)^n,\]
where $x_0$ is a fixed real number, and $c_0, c_1, c_2, \ldots$ are the coefficients. 

\end{definition}
Each term in a power series $\di \sum_{n=0}^{\infty}c_n(x-x_0)^n$ is a simple polynomial $c_n(x-x_0)^n$ which is infinitely differentiable. However, as an infinite series, we need to address the convergence issue. Obviously, the power series $\di \sum_{n=0}^{\infty}c_n(x-x_0)^n$  converges when $x=x_0$. 

Recall that in Chapter \ref{ch5}, we have discussed the ratio test in Theorem \ref{230305_4}. Given $\di\sum_{n=1}^{\infty}a_n$ is a series with $a_n\neq 0$ for all $n\in\mathbb{Z}^+$, let
\[r=\liminf_{n\to\infty}\left|\frac{a_{n+1}}{a_n}\right|\hspace{1cm}\text{and}\hspace{1cm} R=\limsup_{n\to\infty}\left|\frac{a_{n+1}}{a_n}\right|.\]
Then the series $\di\sum_{n=1}^{\infty}a_n$  is divergent if $r>1$, convergent if $R<1$, but inconclusive if $r\leq 1\leq R$.
This test is useful if the limit $\di\lim_{n\to\infty}\left|\frac{a_{n+1}}{a_n}\right|$ exists. For then $r=R$ and we only left with finitely many points which we cannot conclude the convergence of the power series.
Let us look at some examples.
\begin{example}[label=230305_1]{}
Find the domain of convergence of the power series $\di \sum_{n=0}^{\infty}\frac{x^n}{n!}$.
\end{example}
\begin{solution}{Solution}
The power series is convergent when $x=0$. When $x\neq 0$,  using ratio test with $a_n=\di \frac{x^n}{n!}$, we find that
\[ \lim_{n\to\infty}\left|\frac{a_{n+1}}{a_n}\right|=\lim_{n\to\infty}\frac{|x|}{n+1}=0.\] 
Therefore, the series is convergent for all real numbers $x$. The domain of convergence is $\mathbb{R}$, the set of real numbers.
\end{solution}

\begin{example}[label=230305_2]{}
Find the domain of convergence of the power series $\di \sum_{n=0}^{\infty}n!x^n $.
\end{example}
\begin{solution}{Solution}
The power series is convergent when $x=0$. When $|x|\neq 0$, using ratio test with $a_n=\di  n!x^n$, we have 
\[ \lim_{n\to\infty}\left|\frac{a_{n+1}}{a_n}\right|=\lim_{n\to\infty}(n+1)|x|=\infty.\] Hence, the series is divergent if $x\neq 0$.
We conclude that the  series $\di \sum_{n=0}^{\infty}n!x^n $ is only convergent when $x=0$. The domain of convergence is the set $\{0\}$.
\end{solution}

\begin{example}[label=230305_3]{}
Find the domain of convergence of the power series $\di \sum_{n=1}^{\infty}\frac{x^n}{n^2}$.
\end{example}
\begin{solution}{Solution}
The power series is convergent when $x=0$. When $x\neq 0$,
  using  ratio test with $a_n=\di \frac{x^n}{n^2}$, we have
\[ \lim_{n\to\infty}\left|\frac{a_{n+1}}{a_n}\right|=|x|\lim_{n\to\infty}\frac{n^2}{(n+1)^2}=|x|\lim_{n\to\infty}\left(\frac{n}{n+1}\right)^2=|x|.\]

Therefore, the series is convergent if $|x|<1$, and divergent if $ |x|>1$. When $|x|=1$, the test is inconclusive.

But we know that the series $\di\sum_{n=1}^{\infty}\frac{1}{n^2}$ and the series $\di\sum_{n=1}^{\infty}\frac{(-1)^n}{n^2}$ are convergent. Therefore, the series $\di \sum_{n=1}^{\infty}\frac{x^n}{n^2}$ is convergent if and only if $|x|\leq 1$. The domain of convergence is the set $[-1, 1]$.

\end{solution}

In the examples above, we apply the ratio test to determine  the domain of convergence.  This   works fine when all  the coefficients $c_n$ in the power series $
\di\sum_{n=0}^{\infty}c_n(x-x_0)^n$ are nonzero, or only finitely many of them are zero.
We need to find the limit inferior and limit superior of the sequence $\di\left\{\left|\frac{a_{n+1}}{a_n}\right|\right\}$, where $a_n$ is the $n^{\text{th}}$ term $c_n(x-x_0)^n$ in the power series. Since
\[\left|\frac{a_{n+1}}{a_n}\right|=\left|\frac{c_{n+1}(x-x_0)^{n+1}}{c_n(x-x_0)^n}\right|=|x-x_0|\left|\frac{c_{n+1}}{c_n}\right|,\]  essentially we need to find the limit inferior and limit superior of the sequence $\di\left\{\left|\frac{c_{n+1}}{c_n}\right|\right\}$, then multiply by $|x-x_0|$. If the limit of the sequence $\di\left\{\left|\frac{c_{n+1}}{c_n}\right|\right\}$ exists, the limit inferior and limit superior of this sequence are the same, and the domain of convergence can be determined up to the end points of an interval. We apply other convergence test to check the convergence at these end points.

There are two problems with using the ratio test for determining the domain of convergence.
\begin{enumerate}[1.]
\item
If the limit inferior and limit superior of the sequence  $\di\left\{\left|\frac{c_{n+1}}{c_n}\right|\right\}$ are not the same, the ratio test is inconclusive for $x$ in an interval. 
\item When infinitely many of the coefficients $c_n$ in  the power series $\di\sum_{n=0}^{\infty}c_n(x-x_0)^n$ are zero, the ratio test cannot be applied. This problem can be circumvented if there is some patterns on the indices $n$ for which $c_n$ is 0. For example, if $c_{2n}=0$ for all $n\in\mathbb{Z}^+$, the series only contains the odd terms, and it can be written as
\[\sum_{n=1}^{\infty}c_{2n-1}(x-x_0)^{2n-1}.\]
In this case, we can apply the ratio test with $\di a_n=c_{2n-1}(x-x_0)^{2n-1}$. However, the first problem might still be present.
\end{enumerate}

To resolve these problems, we find that the root test (Theorem \ref{230227_23}) is better from the theoretical point of view. Given a  series $\di\sum_{n=1}^{\infty}a_n$, let 
\[\widetilde{\rho}=\limsup_{n\to\infty}\sqrt[n]{|a_n|}.\]The root test says that the series  $\di\sum_{n=1}^{\infty}a_n$ is convergent if $\widetilde{\rho}<1$, divergent if $\widetilde{\rho}>1$, and inconclusive if $\widetilde{\rho}=1$.

Applying the root test to a power series, we have the following.

\begin{theorem}[label=230305_5]{Convergence of Power Series}
Given a power series $\di\sum_{n=0}^{\infty}c_n(x-x_0)^n$, let
\[\rho=\limsup_{n\to\infty}\sqrt[n]{|c_n|}.\]
\begin{enumerate}[1.]
\item If $\rho=0$, then the power series converges  for all real numbers $x$.
\item If $\rho=\infty$, then the power series only converges at the point $x=x_0$.
\item If $\rho$ is a finite positive number, let $R=1/\rho$. Then $R$ is a positive number. The power series is convergent for all $x$ satisfying $|x-x_0|<R$, and divergent for all $x$ satisfying $|x-x_0|>R$. 

\end{enumerate}
\end{theorem}

\begin{myproof}{\linkt Proof of Theorem \ref{230305_5}}
For the power series $\di\sum_{n=0}^{\infty}c_n(x-x_0)^n$, the $n^{\text{th}}$ term is
 $a_n=\di c_n(x-x_0)^n$. 
 \[\widetilde{\rho}=\limsup_{n\to\infty}\sqrt[n]{|a_n|}=|x-x_0|\limsup_{n\to\infty}\sqrt[n]{|c_n|}=\rho|x-x_0|.\]
Now we apply root test as stipulated in Theorem \ref{230227_23}.
\begin{enumerate}[1.]
\item If $\rho=0$, then $\widetilde{\rho}=0$, and so the power series converges  for all real numbers $x$.\end{enumerate}\begin{enumerate}[2.]
\item If $\rho=\infty$, then $\widetilde{\rho}=\infty$ if $x\neq x_0$. Hence, the power series is divergent if $x\neq x_0$. Therefore, the power series only converges at the point $x=x_0$.\end{enumerate} \begin{enumerate}[3.]
\item If $\rho$ is a finite positive number and $R=1/\rho$,  then when $|x-x_0|<R$, $\widetilde{\rho}=|x-x_0|\rho<R\rho=1$; when $|x-x_0|>R$, $\widetilde{\rho}=|x-x_0|\rho>R\rho=1$. Therefore, the power series is convergent when $|x-x_0|<R$, divergent when $|x-x_0|>R$.
\end{enumerate}
\end{myproof}

\begin{corollary}{}
Given a power series $\di\sum_{n=0}^{\infty}c_n(x-x_0)^n$ such that $c_n\neq 0$ for all $n$, assume that the limit
\[\rho=\lim_{n\to\infty}\left|\frac{c_{n+1}}{c_n}\right|\]exists in the extended sense.
\begin{enumerate}[1.]
\item If $\rho=0$, then the power series converges  for all real numbers $x$.
\item If $\rho=\infty$, then the power series only converges at the point $x=x_0$.
\item If $\rho$ is a finite positive number, let $R=1/\rho$. Then $R$ is a positive number. The power series is convergent for all $x$ satisfying $|x-x_0|<R$, and divergent for all $x$ satisfying $|x-x_0|>R$. 

\end{enumerate}
\end{corollary}
\begin{myproof}{Proof}
By Theorem \ref{230227_22}, we find that  $\di\lim_{n\to\infty}\left|\frac{c_{n+1}}{c_n}\right|$ exists implies that
\[\limsup_{n\to\infty}\sqrt[n]{|c_n|}=\lim_{n\to\infty}\sqrt[n]{|c_n|}=\lim_{n\to\infty}\left|\frac{c_{n+1}}{c_n}\right|.\] The rest follows from Theorem \ref{230305_5}.
\end{myproof}

\begin{highlight}{Domain of Convergence}
Theorem \ref{230305_5} shows that the domain of convergence of a power series centered at $x_0$ can only be one of the following cases:
\[\begin{array}{ccp{1cm}cc} 1.&\mathbb{R}  && 2.&
   \{x_0\} \\
   3.&(x_0-R,x_0+R)  &&4.&
   [x_0-R, x_0+R) \\
5.&   (x_0-R, x_0+R]  &&6.&
   [x_0-R, x_0+R] 
\end{array}\]Here $R$ is a positive number.
\end{highlight}
\begin{definition}{Radius of Convergence}
Given a power series $\di\sum_{n=0}^{\infty}c_n(x-x_0)^n$, let
\[\rho=\limsup_{n\to\infty}\sqrt[n]{|c_n|}\]as an extended real number. Then $\rho\geq 0$. Let $R=1/\rho$  in the extended sense. Namely, $R=\infty$ if $\rho=0$, and $R=0$ if $\rho=\infty$. 
This  number $R$ is called the {\bf radius of convergence} of the power series $\di\sum_{n=0}^{\infty}c_n(x-x_0)^n$. The power series is convergent when $|x-x_0|<R$, and divergent when $|x-x_0|>R$. 

\end{definition}

Let us look at the following example.
\begin{example}{}
Let $\{c_n\}$ be the sequence defined by
\[c_n=\begin{cases}n,\quad &\text{if $n$ is even},\\1,\quad  &\text{if $n$ is odd}.\end{cases}\]
Find the domain of convergence of the power series
$\di\sum_{n=1}^{\infty}c_nx^n$.
\end{example}
\begin{solution}{Solution}
Notice that 
\[\frac{c_{n+1}}{c_n}=\begin{cases} n+1,\quad  &\text{if $n$ is odd},\\\di \frac{1}{n},\quad  &\text{if $n$ is even}.\end{cases}\]Applying ratio test with
$a_n=c_nx^n$, we find that if $x\neq 0$,
\[\liminf_{n\to\infty}\left|\frac{a_{n+1}}{a_n}\right| =|x|\liminf_{n\to\infty}\frac{c_{n+1}}{c_n}=0\] \[
\limsup_{n\to\infty}\left|\frac{a_{n+1}}{a_n}\right| =|
x|\limsup_{n\to\infty}\frac{c_{n+1}}{c_n}=\infty.\]
This shows that the ratio test is inconclusive for any $x$ except $x=0$.

Let us turn to root test.
By \eqref{eq230305_7}, we have
\[\lim_{n\to\infty}\sqrt[n]{n}=1.\]  This implies that
\[\lim_{n\to\infty}\sqrt[n]{|c_n|}=1.\]
Therefore, the power series $\di\sum_{n=1}^{\infty}c_nx^n$ is convergent when $|x|<1$, divergent when $|x|>1$. When $x=1$ or $-1$, we have the series
$\di \sum_{n=0}^{\infty}(\pm 1)^nc_n$. Since $\di\lim_{n\to\infty}c_n\neq 0$, we conclude that the power series is divergent when $x=1$ or $x=-1$.
Hence, the domain of convergence of the power series is $(-1,1)$.
\end{solution}
This example shows that  applying ratio test naively will leads to inconclusive scenario, but the root test has rescued the problem. In practice, we always want to avoid applying the root test  because it is difficult to find the limit superior of the sequence $\{\sqrt[n]{|c_n|}\}$ when the coefficients $c_n$. In the example above, we can avoid using root test by writing the power series as a sum of two power series, and apply the ratio test to the two power series separately.  In any case, the root test has given a theoretical  decisive conclusion about the possible types of domain of convergence for a power series.

 For a power series whose radius of convergence $R$ is 0, it only converges at a single point $x=x_0$. So there is no point to consider such power series.
If the radius of convergence $R$ of a power series   $\di \sum_{n=0}^{\infty}c_n(x-x_0)^n$  is positive, the power series defines a function on the open interval 
$(x_0-R, x_0+R)$. We want to study the continuity, differentiability and integrability of such a power series. Therefore, we need to determine whether the power series converges uniformly.

Unfortunately, in general, a power series $\di\sum_{n=0}^{\infty}c_n(x-x_0)^n$  does not converge uniformly on the interval $(x_0-R, x_0+R)$.   For example, consider the  series $\di s(x)=\sum_{n=0}^{\infty}x^n$. In Example \ref{230305_16}, we have seen that   $\di\sum_{n=0}^{\infty} x^n$ is convergent when $|x|<1$, and divergent when $|x|>1$. Hence, its radius of convergence is $R=1$. When $|x|<1$, the power series $\di\sum_{n=0}^{\infty} x^n$ defines the function
\[s(x)=\sum_{n=0}^{\infty} x^n=\frac{1}{1-x}.\] The $n^{\text{th}}$ partial sum of the series is
\[s_n(x)=1+x+\cdots+x^n=\frac{1-x^{n+1}}{1-x}.\]
Therefore, when $x\in (-1,1)$,
\[s(x)-s_n(x)=\frac{1}{1-x}-\frac{1-x^{n+1}}{1-x}=\frac{x^{n+1}}{1-x}.\]Since
\[\lim_{x\to 1^-}\frac{x^{n+1}}{1-x}=\infty,\]we find  that
\[\sup_{|x|<1}|s(x)-s_n(x)|=\infty.\]
Hence, the series $\di\sum_{n=0}^{\infty}x^n$ does not converge uniformly on $(-1,1)$. However, if $a$ is a number such that $0<a<1$, then for $|x|\leq a$,
\[\left|\frac{x^{n+1}}{1-x}\right|\leq \frac{a^{n+1}}{1-a}.\]
Therefore,
\[ \sup_{|x|\leq a}|s(x)-s_n(x)|\leq\frac{a^{n+1}}{1-a}.\]
This implies that
\[\lim_{n\to\infty}\sup_{|x|\leq a}|s(x)-s_n(x)|=0.\]
Hence, the series $\di\sum_{n=0}^{\infty}x^n$  converges uniformly on $[-a,a]$.

A general power series  also have similar behavior.
\begin{theorem}[label=230305_9]{Absolute and Uniform  Convergence of a Power Series}
Given that   $\di \sum_{n=0}^{\infty}c_n(x-x_0)^n$ is a power series whose radius of convergence $R$ is positive. If $R_1$ is any number satisfying $0<R_1<R$, then the power series  $\di \sum_{n=0}^{\infty}c_n(x-x_0)^n$ converges absolutely and uniformly on the set $D_1=\left\{x\,|\,|x-x_0|\leq R_1\right\}$.
\end{theorem}
\begin{myproof}{Proof}
Let $\di R_2=\di \frac{R+R_1}{2}$. Then $R_1<R_2<R$. Hence, the series $\di \sum_{n=0}^{\infty}c_n(x-x_0)^n$  is convergent when $|x-x_0|=R_2$. Let $x_2=x_0+R_2$. Then $\di \sum_{n=0}^{\infty}c_n(x_2-x_0)^n=\sum_{n=0}^{\infty}c_nR_2^n$ is convergent. This implies that
$\di \lim_{n\to\infty}c_nR_2^n=0$.
\bp
In particular, the sequence $\{c_nR_2^n\}$ is bounded. Let $M$ be a positive number such that 
\[|c_nR_2^n|\leq M\hspace{1cm}\text{for all}\;n\geq 0.\]
We apply the Weiertrass $M$-test (Theorem \ref{230305_8}) with $f_n(x)=c_n(x-x_0)^n$. We find that
\[|c_n(x-x_0)^n|\leq |c_n|R_1^n \leq M\left(\frac{R_1}{R_2}\right)^n=Mr^n\hspace{1cm}\text{when}\;|x-x_0|\leq R_1.\]
Here $r=R_1/R_2$. Since $0<r<1$, the geometric series $\di\sum_{n=0}^{\infty} Mr^n$ is convergent. By Weierstrass $M$-test, the power series 
$\di \sum_{n=0}^{\infty}c_n(x-x_0)^n$ converges absolutely and uniformly on the set $D_1=\left\{x\,|\,|x-x_0|\leq R_1\right\}$.\end{myproof}
 
 \begin{remark}{Radius of Convergence Revisited}
In the proof of Theorem \ref{230305_9}, essentially we show that if the power series $\di\sum_{n=0}^{\infty}c_n(x-x_0)^n$ is convergent when $x=x_2$, then it is convergent for all $x$ in the interval $(x_0-R_2, x_0+R_2)$, where $R_2=|x_2-x_0|$. The contrapositive of this statement says that if the power series $\di\sum_{n=0}^{\infty}c_n(x-x_0)^n$  is divergent when $x=x_3$, then it is divergent for all $x$ satisfying $|x-x_0|>R_3$, where $R_3=|x_3-x_0|$. Hence, if 
$S$ is the set
\[S=\left\{|x_1-x_0|\,\left|\, \sum_{n=0}^{\infty}c_n(x-x_0)^n\;\text{is convergent when} \; x=x_1\right.\right\},\]
then $S$ contains only nonnegative numbers. Obviously, 0 is in $S$. If $R_1$ is in $S$,   any positive number $r$ that is less than $R_1$ is also in $S$. This implies that if $R=\sup S$, then $[0, R)\subset S$ and $(R, \infty)$ is disjoint from $S$. This provides an alternative way to define the radius of convergence of the power series without using the root test. Namely, the radius of convergence $R$ is defined as the supremum of the set $S$.
 \end{remark}
 
 From Theorem \ref{230305_9}, we obtain the following.
 \begin{theorem}{Continuity of a Power Series}
Given that   $\di \sum_{n=0}^{\infty}c_n(x-x_0)^n$ is a power series whose radius of convergence $R$ is positive.  It defines a function 
\[f(x)=\sum_{n=0}^{\infty}c_n(x-x_0)^n\] that is continuous on the set $D=\left\{ x\,|\,|x-x_0|< R \right\}$.
\end{theorem}
\begin{myproof}{Proof}
Given any $x_1\in D=\left\{ x\,|\,|x-x_0|< R \right\}$, $R_1=|x_1-x_0|<R$.  Theorem \ref{230305_9} says that  the power series converges uniformly on the set  $D_1=\left\{ x\,|\,|x-x_0|\leq R_1 \right\}$, which contains the point $x_1$. 

 For $n\geq 0$,  the function $f_n(x)=c_n(x-x_0)^n$ is continuous. By Corollary \ref{230305_10}, the power series  $\di \sum_{n=0}^{\infty}c_n(x-x_0)^n$ is continuous at $x=x_1$.  
\end{myproof}

The next is about term by term integration of a power series.
 \begin{theorem}[label=230305_21]{Term by Term Integration of a Power Series}
Given that   $\di \sum_{n=0}^{\infty}c_n(x-x_0)^n$ is a power series whose radius of convergence $R$ is positive. If $[a,b]$ is a closed interval that is contained in the interval $(x_0-R,x_0+R)$, then the function
\[f(x)=\sum_{n=0}^{\infty}c_n(x-x_0)^n\] is Riemann integrable on $[a,b]$, and
we can integrate  term by term. Namely,
\begin{equation}\label{eq230305_12}\int_a^bf(x)dx=\int_a^b\sum_{n=0}^{\infty}c_n(x-x_0)^ndx=\sum_{n=0}^{\infty}c_n\int_a^b (x-x_0)^ndx.\end{equation}
\end{theorem}
\begin{myproof}{Proof}
Let $\di R_1=\max\{|a-x_0|, |b-x_0|\}$. Then $0<R_1<R$ and $[a,b]$ is contained in $[x_0-R_1, x_0+R_1]$. By 
  Theorem \ref{230305_9}, the power series $\di \sum_{n=0}^{\infty}c_n(x-x_0)^n$ converges uniformly on  $[x_0-R_1, x_0+R_1]$, and hence on $[a,b]$. 
    For any $n\geq 0$,  the function $f_n(x)=c_n(x-x_0)^n$ is Riemann integrable.
   By Corollary \ref{230305_11}, the function $\di f(x)= \sum_{n=0}^{\infty}c_n(x-x_0)^n$  is Riemann integrable on $[a,b]$ and \eqref{eq230305_12} holds.
\end{myproof}

Before we discuss term by term differentiation, we need to prove the uniform convergence of the derivative series. We will first prove the following lemma.
\begin{lemma}[label=230618_2]{}
Given that $\{a_n\}$ is a sequence of nonnegative numbers, 
\[\limsup_{n\to \infty}\left(\sqrt[n]{n}\, a_n\right)=\limsup_{n\to\infty}a_n.\]

\end{lemma}
\begin{myproof}{Proof}For all $n\in\mathbb{Z}^+$, $\sqrt[n]{n}\geq 1$. 
By \eqref{eq230305_7}, we have
$\di\lim_{n\to\infty}\sqrt[n]{n}=1$. Hence, given $\varepsilon>0$, there is a positive integer $N$ such that for all $n\geq N$,
\[ 1\leq\sqrt[n]{n}<1+\varepsilon.\]Therefore, for all $n\geq N$,
\[ a_n\leq \sqrt[n]{n}\;a_n\leq \left(1+\varepsilon\right)a_n.\]This implies that
\[ \limsup_{n\to\infty}a_n\leq\limsup_{n\to\infty}\left(\sqrt[n]{n}\;a_n\right)\leq \left(1+\varepsilon\right)\limsup_{n\to\infty}a_n.\]Since $\varepsilon$ can be any positive number, we conclude that \[\limsup_{n\to \infty}\left(\sqrt[n]{n}\, a_n\right)=\limsup_{n\to\infty}a_n.\]
\end{myproof}
\begin{theorem}[label=230305_15]{}
Let $\di\sum_{n=0}^{\infty}c_n(x-x_0)^n$ be a power series with a positive radius of convergence $R$. Then the radius of convergence of the  derived series $\di \sum_{n=1}^{\infty}nc_n(x-x_0)^{n-1} $ is also $R$.
\end{theorem}
\begin{myproof}{Proof}
Let $R'$ be the radius of convergence of the  derived  series $\di \sum_{n=1}^{\infty}nc_n(x-x_0)^{n-1}$. It is not difficult to see that it is the same as the radius of convergence of the series $\di \sum_{n=1}^{\infty}nc_n(x-x_0)^{n}$.  By Lemma \ref{230618_2},
\[\frac{1}{R'}=\limsup_{n\to\infty}\sqrt[n]{|nc_{n}|}= \limsup_{n\to\infty}\sqrt[n]{|c_n|}=\frac{1}{R}.\]
This proves that $R'=R$. 
\end{myproof}

Notice that
if $k\in\mathbb{Z}^+$,
\[\frac{d^k}{dx^k}(x-x_0)^n=n(n-1)\cdots (n-k+1)(x-x_0)^{n-k}.\]
By induction, we can deduce the following.
\begin{corollary}{}
Given that   $\di \sum_{n=0}^{\infty}c_n(x-x_0)^n$ is a power series whose radius of convergence $R$ is positive. For any $k\in\mathbb{Z}^+$, the series \[ \sum_{n=k}^{\infty}n(n-1)\cdots (n-k+1)c_n(x-x_0)^{n-k}\] has radius of convergence $R$.
\end{corollary}
\begin{myproof}{Proof}
The $k=1$ case, which says that the series $\di\sum_{n=1}^{\infty}nc_n(x-x_0)^{n-1}$ has radius of convergence $R$,  is given by Theorem \ref{230305_15}. Applying Theorem  \ref{230305_15} to the series $\di\sum_{n=1}^{\infty}nc_n(x-x_0)^{n-1}$, we find that the series  $\di\sum_{n=2}^{\infty}n(n-1)c_n(x-x_0)^{n-2}$ also has radius of convergence $R$. This is the statement we need to prove for the $k=2$ case. For general $k\in\mathbb{Z}^+$, we  proceed by induction.
\end{myproof}

The next theorem says  that we can differentiate a power series term by term.
\begin{theorem}{Term by Term Differentiation of a Power Series}
Given that   $\di \sum_{n=0}^{\infty}c_n(x-x_0)^n$ is a power series whose radius of convergence $R$ is positive. 
Then the   function
\[f(x)=\sum_{n=0}^{\infty}c_n(x-x_0)^n\] is   differentiable on $(x_0-R, x_0+R)$. When $x\in (x_0-R, x_0+R)$,
  we can differentiate the power series  $\di \sum_{n=0}^{\infty}c_n(x-x_0)^n$ term by term to obtain
\begin{equation}\label{eq230305_13}f'(x)=\frac{d}{dx}\sum_{n=0}^{\infty}c_n(x-x_0)^n =\sum_{n=1}^{\infty}nc_n (x-x_0)^{n-1}.\end{equation}
 
\end{theorem}

\begin{myproof}{Proof}
By Theorem \ref{230305_15}, the radius of convergence of the derived series $\di \sum_{n=1}^{\infty}nc_n (x-x_0)^{n-1}$ is also $R$. By Theorem \ref{230305_9}, the series  $\di\sum_{n=1}^{\infty}nc_n (x-x_0)^{n-1}$ converges absolutely and uniformly on $[x_0-R_1, x_0+R_1]$ if $R_1<R$. Given $x_1\in (x_0-R, x_0+R)$, $|x_1-x_0|<R$. Choose $R_1$ such that $|x_1-x_0|<R_1<R$. Then $x_1\in (x_0-R_1, x_0+R_1)$. Corollary \ref{230304_14} implies that the function\bp
\[ f(x)=\sum_{n=0}^{\infty}c_n(x-x_0)^n\]   is differentiable on  $(x_0-R_1, x_0+R_1)$, and we can perform term by term differentiation to obtain \eqref{eq230305_13}.  Since $x_1$ is any point in $(x_0-R, x_0+R)$, this proves the statement of the theorem.
\end{myproof}

By induction, we have the following.
\begin{corollary}[label=230305_19]{}
Let $\di\sum_{n=0}^{\infty}c_n(x-x_0)^n$ be a power series with a positive radius of convergence $R$. Then the function
\[f(x)=\sum_{n=0}^{\infty}c_n(x-x_0)^n\] is infinitely differentiable on $(x_0-R, x_0+R)$.
 For any $k\geq 1$ and $x\in (x_0-R, x_0+R)$,
\begin{equation}\label{eq230305_14}f^{(k)}(x) = \sum_{n=k}^{\infty}n(n-1)\cdots (n-k+1)c_n(x-x_0)^{n-k}.\end{equation}
In particular,
\[f^{(k)}(x_0)=k!c_k.\]
\end{corollary}
Let us summarize what we have learned about power series.
\begin{highlight}{Functions Defined by Power Series}

A power series  $\di\sum_{n=0}^{\infty}c_n(x-x_0)^n$ is convergent when $x=x_0$. If the series is convergent for some $x_1\neq x_0$, then it has a positive radius of convergence $R$. The series is convergent for all $x$ satisfying $|x-x_0|<R$, and divergent for all $x$ satisfying $|x-x_0|>R$.\end{highlight}\begin{highlight}{}

The power series defines a function
\[f(x)=\sum_{n=0}^{\infty}c_n(x-x_0)^n\] on the interval $(x_0-R, x_0+R)$. This function $f(x)$ is infinitely differentiable, and we can perform term by term differentiation and term by term integration.

Functions that are representable by power series are called {\bf analytic functions}. Their  domains can be naturally  extended to complex numbers. This is the main topic that is discussed in a course in complex analysis.
\end{highlight}

\begin{definition}{Power Series Expansion of a Function}If a power series $\di\sum_{n=0}^{\infty}c_n(x-x_0)^n$ has positive radius of convergence $R$, it defines an analytic function
$f:(x_0-R,x_0+R)\to\mathbb{R}$ by
\[f(x)= \sum_{n=0}^{\infty}c_n(x-x_0)^n.\]
 We say that $\di\sum_{n=0}^{\infty}c_n(x-x_0)^n$ is a {\bf power series expansion} or {\bf power series representation} of the function $f(x)$ on the interval $(x_0-R, x_0+R)$.
\end{definition}

\begin{example}{} When $|x|<1$,  the function $\di f(x)=\frac{1}{1-x}$ has a power series expansion given by
\begin{equation}\label{eq230305_17}\frac{1}{1-x}=\sum_{n=0}^{\infty}x^n=1+x+x^2+\cdots+x^n+\cdots.\end{equation}
\end{example}

Applying term by term differentiation to \eqref{eq230305_17}, we obtain the following.

\begin{theorem}[label=230307_1]{}
Let $k$ be a nonnegative integer. Then for $|x|<1$, 
\begin{equation}\label{eq230305_20}\frac{1}{(1-x)^{k+1}} =\sum_{n=k}^{\infty}\binom{n}{k}x^{n-k}.\end{equation}Here
$\di \binom{n}{k}=\di \frac{n(n-1)\cdots(n-k+1)}{k!}$ are the binomial coefficients.
\end{theorem}
\begin{myproof}{Proof}
The $k=0$ case is just the formula \eqref{eq230305_17}. By Corollary \ref{230305_19}, we can differentiate term by term $k$ times and \eqref{eq230305_14} gives
\[\frac{k!}{(1-x)^{k+1}}=\sum_{n=k}^{\infty}n(n-1)\cdots (n-k+1)x^{n-k}\hspace{1cm}\text{when}\;|x|<1.\]
Dividing by $k!$ on both sides gives \eqref{eq230305_20}.
\end{myproof}

The formula \eqref{eq230305_20} is very useful. It   has applications in probability theory.
\begin{example}{}
In probability theory, a geometric random variable is a random variable $X$ that depends on a parameter $p$ where $0<p<1$. If one performs a series of identical and independent Bernoulli trials, each has a probability $p$ to be a success, then $X$ is the number of these Bernoulli trials need to be performed until the first success occurs.  For any $n\in\mathbb{Z}^+$, the probability that $X$ is equal to $n$ is 
\[P(X=n)=(1-p)^{n-1}p.\]\be The expected number of Bernoulli trials need to be performed until the first success is
\[E(X)=\sum_{n=1}^{\infty}nP(X=n)=p\sum_{n=1}^{\infty}n(1-p)^{n-1}.\]
Using \eqref{eq230305_20} with $k=1$ and $x=1-p$, we find that
\[E(X)=\frac{p}{(1-(1-p))^2}=\frac{1}{p}.\] 
The variance of $X$ is given by $\text{Var}\,(X)=E(X^2)-E(X)^2$. To find this, we compute $E(X^2)$ first.
\begin{align*}
E(X^2)&=\sum_{n=1}^{\infty}n^2P(X=n)=p\sum_{n=1}^{\infty}n^2(1-p)^{n-1}.
\end{align*}Using \eqref{eq230305_20} with $k=2$ and $x=1-p$, we find that
\[\sum_{n=1}^{\infty}n(n-1)(1-p)^{n-2}=\frac{2}{(1-(1-p))^3}=\frac{2}{p^3}.\]
Thus,
\begin{align*}
\sum_{n=1}^{\infty}n^2(1-p)^{n-1}&=\sum_{n=1}^{\infty}n(n-1)(1-p)^{n-1}+\sum_{n=1}^{\infty}n(1-p)^{n-1}\\&=\frac{2(1-p)}{p^3}+\frac{1}{p^2}=\frac{2-p}{p^3}.\end{align*}
Therefore, the variance of $X$ is
\[\text{Var}\,(X)=E(X^2)-E(X)^2=\frac{2-p}{p^2}-\frac{1}{p^2}=\frac{1-p}{p^2}.\]

\end{example2}

Recall that the logarithm function $f(x)=\ln x$ is defined so that
$\di f'(x)=\frac{1}{x}$.
This gives
\[\frac{d}{dx}\ln (1+x)=\frac{1}{1+x}.\]
Using term by term integration, we can obtain power series representation for the logarithm function.
\begin{theorem}{Power Series Expansion of Logarithm Function }
For $|x|<1$, the power series $\di\sum_{n=1}^{\infty}(-1)^{n-1}\frac{x^n}{n}$ is convergent and
\[\ln(1+x)=\sum_{n=1}^{\infty}(-1)^{n-1}\frac{x^n}{n}=x-\frac{x^2}{2}+\frac{x^3}{3}+\cdots+(-1)^{n-1}\frac{x^n}{n}+\cdots.\]
\end{theorem}
\begin{myproof}{Proof}
Given any $x_1$ with $|x_1|<1$, let $R_1=|x_1|$. Since the geometric series $\di \sum_{n=0}^{\infty}x^n$ has radius of convergence 1, Theorem \ref{230305_21} says that we can integrate \eqref{eq230305_17} term by term over the interval with 0 and $x_1$ as end points.
\[\int_0^{x_1}\frac{1}{1-x}dx=\int_0^{x_1}\sum_{n=0}^{\infty}x^ndx=\sum_{n=0}^{\infty}\int_0^{x_1}x^ndx.\]
This gives
\[-\ln(1-x_1)=\sum_{n=0}^{\infty}\frac{x_1^{n+1}}{n+1}=\sum_{n=1}^{\infty}\frac{x_1^n}{n}.\]
Replacing $x_1$ by $-x$, we find that if $|x|<1$, then
\begin{equation}\label{eq230305_22}\ln(1+x)=-\sum_{n=1}^{\infty}(-1)^n\frac{x^n}{n}=\sum_{n=1}^{\infty}(-1)^{n-1}\frac{x^n}{n}.\end{equation}
\end{myproof}
The theories that we have deveoped so far do not allow us to take the limit $x\to 1^-$ term by term on the right hand side of \eqref{eq230305_22}.   However, we can go around this problem in another way.
\begin{example}{}
Show that 
\begin{equation}\label{eq230305_23}1-\frac{1}{2}+\frac{1}{3}-\frac{1}{4}+\cdots=\sum_{n=1}^{\infty}\frac{(-1)^{n-1}}{n}=\ln 2.\end{equation}
\end{example}
\begin{solution}{Solution}
Notice that if $x\neq -1$, then for any $n\in\mathbb{Z}^+$,
\[1-x+x^2-x^3+\cdots+(-1)^{n-1}x^{n-1}=\frac{1-(-x)^{n}}{1+x}=\frac{1}{1+x}+(-1)^{n-1}\frac{x^{n}}{1+x}.\] 
Each of the functions is continuous on $[0,1]$. Therefore, 
\begin{align*}
&\int_0^1\left(1-x+x^2-x^3+\cdots+(-1)^{n-1}x^{n-1}\right)dx\\&=\int_0^1\frac{1}{1+x}dx+\int_0^1(-1)^{n-1}\frac{x^{n}}{1+x}dx.\end{align*}
This gives
\[ \sum_{k=1}^n\frac{(-1)^{k-1}}{k}=1-\frac{1}{2}+\frac{1}{3}-\frac{1}{4}+\cdots +(-1)^{n-1}\frac{1}{n}=\ln 2+R_n,\]
where
\[R_n=(-1)^{n-1}\int_0^1\frac{x^{n}}{1+x}dx.\]
  When  $0\leq x\leq 1$, $1+x\geq 1$, and so
\[\left| \frac{x^{n}}{1+x}\right|\leq x^n\hspace{1cm}\text{when}\;0\leq x\leq 1.\]
Therefore,
\[|R_n|\leq\int_0^1\left| \frac{x^{n}}{1+x}\right|dx\leq \int_0^1x^ndx=\frac{1}{n+1}.\]
This implies that
$\di \lim_{n\to\infty}R_n=0$. Thereofore,
\[1-\frac{1}{2}+\frac{1}{3}-\frac{1}{4}+\cdots=\lim_{n\to\infty}\sum_{k=1}^n\frac{(-1)^{k-1}}{k}=\ln2.\] 
\end{solution}
 
Next, we give the power series that represents the exponential function.
\begin{theorem}[label=230307_6]{Power Series Expansion of Exponential Function}
For any real numbers $x$,
\[e^x=\sum_{n=0}^{\infty}\frac{x^n}{n!}=1+\frac{x}{1!}+\frac{x^2}{2!}+\frac{x^3}{3!}+\frac{x^4}{4!}+\cdots.\]
\end{theorem}
\begin{myproof}{Proof}
We have shown in Example \ref{230305_1} that the power series $\di\sum_{n=0}^{\infty}\frac{x^n}{n!}$ is convergent for all real numbers $x$.  Corollary \ref{230305_19} says that it defines an infinitely differentiable function $f:\mathbb{R}\to\mathbb{R}$ by
\[f(x)=\sum_{n=0}^{\infty}\frac{x^n}{n!}=1+\frac{x}{1!}+\frac{x^2}{2!}+\frac{x^3}{3!}+\frac{x^4}{4!}+\cdots.\]From this, we have $f(0)=1$. Term by term differentiation gives
\[f'(x)=\sum_{n=1}^{\infty}\frac{nx^{n-1}}{n!}=\sum_{n=1}^{\infty}\frac{x^{n-1}}{(n-1)!}=\sum_{n=0}^{\infty}\frac{x^n}{n!}=f(x).\]
Let $g:\mathbb{R}\to\mathbb{R}$ be the function defined by
\[g(x)=e^{-x}f(x).\]
Then we find that
\[g'(x)=e^{-x}f'(x)-e^{-x}f(x)=0.\]
This shows that there is a constant $C$ such that
\[g(x)=C\hspace{1cm}\text{for all}\;x\in\mathbb{R}.\]
Set $x=0$, we find that $C=g(0)=e^{0}f(0)=1$. Hence, $e^{-x}f(x)=1$ for all real numbers $x$, which implies that $f(x)=e^x$ for all real numbers $x$.
\end{myproof}

Now we want to return to address an existence problem in Chapter \ref{ch3}. In Theorem \ref{thm230218_3}, we claim that  there is a twice differentiable function $f(x)$ that satisfies the equation
\[f''(x)+f(x)=0\] and the initial conditions
\[f(0)=0,\quad f'(0)=1.\]
We define this function as $\sin x$. 
We can now prove the existence. This is actually the power series method for solving differential equations.
Assume that $f(x)$ can be written as a power series 
\[f(x)=\sum_{n=0}^{\infty}c_nx^n.\]
Then $f(0)=0$ and $f'(0)=1$ implies that $c_0=0$ and $c_1=1$. Differentiate two times, we have
\[f''(x)=\sum_{n=2}^{\infty}n(n-1)c_nx^{n-2}=\sum_{n=0}^{\infty}(n+2)(n+1)c_{n+2}x^n.\]
Substitute into the equation $f''(x)+f(x)=0$, we find that
\[\sum_{n=0}^{\infty}\left[(n+2)(n+1)c_{n+2}+c_n\right]x^n=0.\]
Hence, we find that if $\{c_n\}$ is defined recursively by $c_0=0$, $c_1=1$, and for all $n\geq 0$,
\[c_{n+2}=-\frac{c_n}{(n+1)(n+2)},\] we get a candidate solution for our problem.  The recursive formula for $\{c_n\}$ can be easily solved to give
\[c_{2n-1}=\frac{(-1)^{n-1}}{(2n-1)!},\hspace{1cm}c_{2n}=0\hspace{1cm}\text{for all}\;n\in\mathbb{Z}^+.\]Now we are left to justify this is indeed the solution to our problem.
\begin{theorem}{}
The power series \[\di\sum_{n=1}^{\infty} \frac{(-1)^{n-1}}{(2n-1)!}x^{2n-1}=x-\frac{x^3}{3!}+\frac{x^5}{5!}-\frac{x^7}{7!}+\cdots\] defines an infinitely differentiable function $f:\mathbb{R}\to\mathbb{R}$ that satisfies  
\[f''(x)+f(x)=0, \hspace{1cm}f(0)=0, \;f'(0)=1.\] 
\end{theorem}
\begin{myproof}{Proof}
First, we need to show that the power series is convergent everywhere. We can use the ratio test with $\di a_n=\frac{(-1)^{n-1}}{(2n-1)!}x^{2n-1}$. Then if $x\neq 0$, 
\[\lim_{n\to\infty}\left|\frac{a_{n+1}}{a_n}\right|=x^2\lim_{n\to\infty}\frac{1}{2n(2n+1)}=0.\]
This shows that the power series is convergent for all $x\in\mathbb{R}$. By Corollary \ref{230305_19}, it defines an infinitely differentiable function $f:\mathbb{R}\to\mathbb{R}$,
\[f(x)=\sum_{n=1}^{\infty} \frac{(-1)^{n-1}}{(2n-1)!}x^{2n-1}=x-\frac{x^3}{3!}+\frac{x^5}{5!}-\frac{x^7}{7!}+\cdots.\]From here, it is straightforward to find that $f(0)=0$.
 We can differentiate term by term to obtain
\[f'(x)=1-\frac{x^2}{2!}+\frac{x^4}{4}-\frac{x^6}{6!}+\cdots.\]This gives $f'(0)=1$. Differentiate term by term again, we have
\begin{align*}f''(x)&=\sum_{n=2}^{\infty} (2n-1)(2n-2)\frac{(-1)^{n-1}}{(2n-1)!}x^{2n-3}\\&=\sum_{n=2}^{\infty}  \frac{(-1)^{n-1}}{(2n-3)!}x^{2n-3}=- \sum_{n=1}^{\infty} \frac{(-1)^{n-1}}{(2n-1)!}x^{2n-1}=-f(x).\end{align*}This proves that $f''(x)+f(x)=0$, and thus the proof is completed.
\end{myproof}

As a byproduct, we obtain the power series expansion for the functions $\sin x$ and $\cos x$.
\begin{theorem}{Power Series Expansion of Sine and Cosine Functions}
For any real numbers $x$,
\begin{align*}\sin x&=\sum_{n=1}^{\infty} \frac{(-1)^{n-1}}{(2n-1)!}x^{2n-1}=x-\frac{x^3}{3!}+\frac{x^5}{5!}-\frac{x^7}{7!}+\cdots,\\
\cos x&=\sum_{n=0}^{\infty} \frac{(-1)^{n}}{(2n)!}x^{2n}=1-\frac{x^2}{2!}+\frac{x^4}{4}-\frac{x^6}{6!}+\cdots.\end{align*}
\end{theorem}The power series for $\cos x$ is obtained by term by term differentiating the power series for $\sin x$. 

Finally, we want to consider the multiplication of two power series. Given that $p(x)$ and $q(x)$ are polynomials of degree $k$ and $l$ respectively,   with
\begin{align*}
p(x)=a_0+a_1x+\cdots+a_kx^k\quad\text{and}\quad q(x)=b_0+b_1x+\cdots+b_lx^l.
\end{align*}The product $p(x)q(x)$ is a polynomial of degree $k+l$, with
\begin{align*}p(x)q(x)&=c_0+c_1x+\cdots+c_{k+1}x^{k+1}\\
&=(a_0b_0)+(a_0b_1+a_1b_0)x+\cdots  +a_kb_lx^{k+l}.\end{align*}For $0\leq n\leq \max\{k, l\}$, we find that
\[c_n=a_0b_n+a_1b_{n-1}+\cdots+a_{n-1}b_1+a_nb_0=\sum_{m=0}^{n}a_{m}b_{n-m}.\]This motivates the following.
\begin{definition}{Cauchy Product of Two Series}
Given the two infinite series $\di\sum_{n=0}^{\infty}a_n$ and $\di\sum_{n=0}^{\infty}b_n$, their {\bf Cauchy product} is the infinite series $\di\sum_{n=0}^{\infty}c_n$, where
\[c_n=\sum_{m=0}^{n}a_{m}b_{n-m}.\]
\end{definition}The following theorem is a special case of the Merten's theorem on Cauchy products.  
\begin{theorem}[label=230306_2]{Term by Term Multiplication of Power Series}
Let $\di\sum_{n=0}^{\infty}a_n(x-x_0)^n$ and $\di\sum_{n=0}^{\infty}b_n(x-x_0)^n$ be two power series with positive radii of convergence $R_a$ and $R_b$ respectively. Define the sequence $\{c_n\}_{n=0}^{\infty}$ by
\[c_n=\sum_{m=0}^n a_mb_{n-m}.\]Then the power series  $\di\sum_{n=0}^{\infty}c_n(x-x_0)^n$ has radius of convergence $R\geq R_c$, where $R_c= \min\{R_a, R_b\}$. If
\[f(x)=\sum_{n=0}^{\infty}c_n(x-x_0)^n,\quad g(x)=\sum_{n=0}^{\infty}a_n(x-x_0)^n,\quad h(x)=\sum_{n=0}^{\infty}b_n(x-x_0)^n \]are the functions defined by each of the power series on $(x_0-R_c, x_0+R_c)$, then
 we have
\[f(x)=g(x)h(x).\]

\end{theorem}
\begin{myproof}{Proof}Without loss of generality, assume that   $x_0=0$.  

It is suficient to prove that for any $x_1$ satisfying $0<R_1=|x_1|<R_c$, the series $\di\sum_{n=0}^{\infty}c_nx_1^n$ is convergent, and it converges to $AB$, where
\[A=\sum_{n=0}^{\infty}a_nx_1^n,\hspace{1cm}B=\sum_{n=0}^{\infty}b_nx_1^n.\]
This would imply that the series $\di\sum_{n=0}^{\infty}c_nx^n$ is convergent on $(-R_c, R_c)$, which proves that its radius of convergence $R$ is  at least $R_c$.\bp
Take an $R_2$ such that $R_1<R_2<R_c$. Then $R_2<R_a$ and $R_2<R_b$. Using the same reasoning as in the proof of Theorem \ref{230305_9}, we find that there is a positive constant $M$ such that
\[|a_nx_1^n|\leq Mr^n\quad\text{and}\quad |b_nx_1^n|\leq Mr^n\hspace{1cm}\text{for all}\;n\geq 0,\]
where $r=R_1/R_2$ is a number satisfying $0<r<1$.

For a positive integer $n$, let 
\[C_n =\sum_{k=0}^{n}c_kx_1^k,\quad A_n =\sum_{k=0}^{n}a_kx_1^k,\quad B_n =\sum_{k=0}^{n}b_kx_1^k \]be the partial sums of each series.
For any $n\in\mathbb{Z}^+$, 
\[|B-B_n|=\left|\sum_{k=n+1}^{\infty}b_kx_1^k\right|\leq \sum_{k=n+1}^{\infty}\left|b_kx_1^k\right|\leq\sum_{k=n+1}^{\infty}Mr^{k}=\frac{Mr^{n+1}}{1-r}.\]
By definitions of the sequence $\{c_n\}$, we find that for $n\in\mathbb{Z}^+$,
\begin{align*}
C_n&=\sum_{l=0}^n\sum_{k=0}^la_kb_{l-k}x_1^l=\sum_{k=0}^n a_kx_1^k\sum_{l=k}^n b_{l-k}x_1^{l-k}\\&=\sum_{k=0}^n a_kx_1^k\sum_{l=0}^{n-k} b_{l}x_1^{l}=\sum_{k=0}^na_kx_1^k B_{n-k} \\
&=BA_n -\sum_{k=0}^{n}a_kx_1^k(B-B_{n-k}).
\end{align*} Therefore, for any $n\in\mathbb{Z}^+$,
\begin{align*}
|C_n-BA_n|&\leq\sum_{k=0}^n |a_kx_1^k||B-B_{n-k}|\leq\sum_{k=0}^n Mr^k\times \frac{Mr^{n-k+1}}{1-r}\\
&=\frac{M^2}{1-r}(n+1)r^{n+1}.
\end{align*}By Theorem \ref{230306_1}, $\di \lim_{n\to\infty}(n+1)r^{n+1}=0$. By squeeze theorem, we find that
\[\lim_{n\to\infty}C_n=B\lim_{n\to\infty}A_n=AB,\]which completes the proof of the theorem.
 
\end{myproof}

Let us look at an example.
\begin{example}{}
Consider the function $f:\mathbb{R}\to\mathbb{R}$ defined by $f(x)=e^x\sin x$. Find  the power series expansion of $f(x)$ up to the $x^5$ term, and find $f^{(5)}(0)$. 
\end{example}
\begin{solution}{Solution} 
We know that \begin{align*}
e^x&=\sum_{n=0}^{\infty}\frac{x^n}{n!}=1+x+\frac{x^2}{2}+\frac{x^3}{6}+\frac{x^4}{24}+\frac{x^5}{120}+\cdots\hspace{1cm}\text{for all}\;x\in\mathbb{R},\\
 \sin x&=\sum_{n=1}^{\infty}(-1)^{n-1}\frac{x^{2n-1}}{(2n-1)!}=x-\frac{x^3}{6}+\frac{x^5}{120}+\cdots\hspace{1cm}\text{for all}\;x\in\mathbb{R}.\end{align*}
By Theorem \ref{230306_2},
\begin{align*}
e^x\sin x&=\left(1+x+\frac{x^2}{2}+\frac{x^3}{6}+\frac{x^4}{24}+\frac{x^5}{120}+\cdots\right)\left(x-\frac{x^3}{6} +\frac{x^5}{120}+\cdots\right)\\
&=x+x^2+\frac{x^3}{2}+\frac{x^4}{6}+\frac{x^5}{24}-\frac{x^3}{6}-\frac{x^4}{6}-\frac{x^5}{12}+\frac{x^5}{120}+\cdots\\
&=x+x^2+\frac{x^3}{3}-\frac{x^5}{30}+\cdots.
\end{align*}This gives the power series expansion of $f(x)$ up to the  $x^5$ term.  From this, we find that
\[f^{(5)}(0)=5!\times\left(-\frac{1}{30}\right)=-4.\]
\end{solution}
\vp
\noindent
{\bf \large Exercises  \thesection}
\setcounter{myquestion}{1}
\begin{question}{\themyquestion}
Let $p$ be a positive number. Determine the  domain of convergence of the power series $\di \sum_{n=0}^{\infty}\frac{x^n}{n^p}$.
\end{question}
\atc
 \begin{question}{\themyquestion}Show that when $|x|<1$,
\[\tan^{-1}x=\sum_{n=1}^{\infty}(-1)^{n-1}\frac{x^{2n-1}}{2n-1}=x-\frac{x^3}{3}+\frac{x^5}{5}+\cdots.\]
\end{question}
\atc
\begin{question}{ \themyquestion\; [The Newton-Gregory Formula]}
Show that 
\[1-\frac{1}{3}+\frac{1}{5}-\frac{1}{7}+\cdots=\sum_{n=1}^{\infty}\frac{ (-1)^{n-1}}{2n-1}=\frac{\pi}{4}.\]
\end{question}
\atc
\begin{question}{\themyquestion}
Find a closed form formula for the sum of the series
$\di \sum_{n=1}^{\infty}n^3x^n$ when $|x|<1$.
\end{question}

\atc
\begin{question}{\themyquestion}
Consider the function $f:\mathbb{R}\to\mathbb{R}$ defined by $f(x)=e^x\cos x$. Find  the power series expansion of $f(x)$ up to the  $x^5$ term, and find $f^{(5)}(0)$. 
\end{question}
\vp

\section{Taylor Series and Taylor Polynomials}\label{sec6.5}
In Section \ref{sec6.4}, we have seen that  the exponential function, logarithm function, sine and cosine functions have power series representations that are valid on its domain or a subset of its domain. Power series are limits of sequences of polynomials. They are infinitely differentiable, and they can be differentiated term by term and integrated term by term. Thus, they are very useful. Hence, we can ask the following two questions.

\begin{highlight}{}
\begin{enumerate}[1.]
\item
If $I$ is an open interval that contains the point $x_0$, and the function $f:I\to\mathbb{R}$ is infinitely differentiable on $I$, does there exist a positive constant $R$ and a power series $\di\sum_{n=0}^{\infty}c_n(x-x_0)^n$ such that
\[f(x)=\sum_{n=0}^{\infty}c_n(x-x_0)^n\hspace{1cm}\text{when}\;|x-x_0|<R.\]
\item If the power series expansion exists, what is the error when we approximate $f(x)$ by the partial sum $\di s_n(x)=\sum_{k=0}^{n}c_k(x-x_0)^k$?

\end{enumerate}
\end{highlight}
For the first question, Corollary \ref{230305_19} says that if such a representation exists, then we must have
\[f^{(n)}(x_0)=n!c_n\hspace{1cm}\text{for all}\; n\in\mathbb{Z}^+.\]

This leads us to the following definition.
\begin{definition}{Taylor Series and Maclaurin Series}
If $I$ is an open interval that contains the point $x_0$, and the function $f:I\to\mathbb{R}$ is infinitely differentiable on $I$,  the {\bf Taylor series} of $f(x)$ at $x_0$ is the series
\[\sum_{n=0}^{\infty}\frac{f^{(n)}(x_0)}{n!}(x-x_0)^n=f(x_0)+\frac{f'(x_0)}{1!}(x-x_0) +\cdots+\frac{f^{(n)}(x_0)}{n!}(x-x_0)^n+\cdots.\]When $x_0=0$, the Taylor series at $0$ is also called a {\bf Maclaurin series}.
\end{definition}Here the Taylor series is defined as a power series as long as the function is infinitely differentiable in an open interval $I$ that contains the point $x_0$. We do not assume any convergence. Even though the Taylor series is convergent, we cannot assume that it converges to the function $f(x)$ itself. In Section \ref{sec6.6.3}, we are going to see a classical example of an infinitely differentiable function whose Taylor series converges but to a different function.

Nevertheless, for functions that are defined by a power series centered at $x_0$, Corollary \ref{230305_19} gives the following.
\begin{theorem}[label=230306_3]{}
Assume that the power series $\di\sum_{n=0}^{\infty}c_n(x-x_0)^n$ has positive radius of convergence $R$. If $f(x)$ is the function defined by the power series $\di\sum_{n=0}^{\infty}c_n(x-x_0)^n$ on the interval $(x_0-R, x_0+R)$, then the Taylor series of $f(x)$ at $x_0$ is $\di\sum_{n=0}^{\infty}c_n(x-x_0)^n$. Namely,
\[f(x)=\sum_{n=0}^{\infty}\frac{f^{(n)}(x_0)}{n!}(x-x_0)^n\hspace{1cm}\text{when}\;x\in (x_0-R, x_0+R).\]This shows that the Taylor series of $f(x)$ converges to the function $f(x)$.
It also says that the power series expansion of a function at a point $x_0$, if exists,  is unique, which is the Taylor series of the function at $x_0$.
\end{theorem}
We have the following list of Maclaurin series from Section \ref{sec6.4}.

\begin{highlight}{Useful Maclaurin Series}
\begin{enumerate}[1.]
\item $\di \frac{1}{1-x}=\sum_{n=0}^{\infty}x^n=1+x+x^2+x^3+\cdots$ when $|x|<1$.
\item $\di e^x=\sum_{n=0}^{\infty}\frac{x^n}{n!}=1+x+\frac{x^2}{2!}+\frac{x^3}{3!}+\cdots$ for all $x\in\mathbb{R}$.
\item $\di \sin x=\sum_{n=1}^{\infty}(-1)^{n-1}\frac{x^{2n-1}}{(2n-1)!}=x-\frac{x^3}{3!}+\frac{x^5}{5!}-\frac{x^7}{7!}+\cdots$ for all $x\in\mathbb{R}$.
\item $\di \cos x= \sum_{n=0}^{\infty}(-1)^n\frac{x^{2n }}{(2n)!}=1-\frac{x^2}{2!}+\frac{x^4}{4!}-\frac{x^6}{6!} +\cdots$ for all $x\in\mathbb{R}$.
\item $\di \ln(1+x)=\sum_{n=1}^{\infty}(-1)^{n-1}\frac{x^n}{n}=x-\frac{x^2}{2}+\frac{x^3}{3}-\frac{x^4}{4}+\cdots$ when $|x|<1$.
\end{enumerate}
\end{highlight}

\begin{remark}{Maclaurin Series for Odd Functions and Even Functions}
Let $a$ be a positive number and let $f:(-a,a)\to\mathbb{R}$ be an infinitely differentiable function. Since the derivative of an odd function is even, and the derivative of an even function is odd, the following holds.
\begin{enumerate}[1.]
\item If $f(x)$ is an odd function, the Taylor series of $f(x)$ at $x=0$ has the form
\[\sum_{n=1}^{\infty}\frac{f^{(2n-1)}(0)}{(2n-1)!}x^{2n-1},\]which only contains the odd power terms. 
\item If $f(x)$ is a even function,   the Taylor series of $f(x)$ at $x=0$ has the form
\[\sum_{n=0}^{\infty}\frac{f^{(2n)}(0)}{(2n)!}x^{2n},\]which only contains the even power terms.
\end{enumerate}
\end{remark}

By the uniqueness of power series expansion asserted in Theorem \ref{230306_3}, and the results proved in Section \ref{sec6.4}, we can use term by term addition, multiplication, differentiation and integration   to obtain the power series for new functions from old ones. This is a useful tactic to find Taylor series of a large class of functions from the few elementary ones listed above.

\begin{example}{}
Find the power series expansion of the function $f(x)=\di\frac{x+1}{4+x^2}$ at $x=0$, and find the largest open interval where this series is convergent.

\end{example}
\begin{solution}{Solution}
Applying the formula
\[\frac{1}{1-x}=\sum_{n=0}^{\infty}x^n,\hspace{1cm}|x|<1 \]
gives
\begin{align*}
\frac{1}{4+x^2}= \frac{1}{4\left(1+\di\frac{x^2}{4}\right)}=\frac{1}{4}\sum_{n=0}^{\infty}(-1)^n\frac{x^{2n}}{2^{2n}},\hspace{1cm}|x|<2.
\end{align*}Multiply by $1+x$, we find that when $|x|<2$, 
\begin{align*}
\frac{x+1}{4+x^2}=\frac{1}{4}\sum_{n=0}^{\infty}(-1)^n\frac{x^{2n+1}}{2^{2n}}+\frac{1}{4}\sum_{n=0}^{\infty}(-1)^n\frac{x^{2n}}{2^{2n}}.
\end{align*} The largest open interval where this series is convergent is $(-2,2)$.
\end{solution}

\begin{example}{}
Let $f:\mathbb{R}\to\mathbb{R}$ be the function defined by
\[f(x)=\begin{cases} \di\frac{\sin x}{x},\quad &\text{if}\;x\neq 0,\\1,\quad & \text{if} \;x=0.\end{cases}\]Show that  $f$ is infinitely differentiable, and find $f^{(n)}(0) $ for all $n\geq 0$.
\end{example}
\begin{solution}{Solution}
For any real number $x$, the series
\[\sum_{n=1}^{\infty}(-1)^{n-1}\frac{x^{2n-1}}{(2n-1)!}=x-\frac{x^3}{3!}+\frac{x^5}{5!}-\frac{x^7}{7!}+\cdots\] 

converges to $\sin x$. When $x\neq 0$, dividing by $x$, we find that the series 
\[\sum_{n=1}^{\infty}(-1)^{n-1}\frac{x^{2n-2}}{(2n-1)!}=\sum_{n=0}^{\infty}(-1)^n\frac{x^{2n}}{(2n+1)!}=1-\frac{x^2}{3!}+\frac{x^4}{5!}-\frac{x^6}{7!}+\cdots\]

 converges to $\di \frac{\sin x}{x}$.  Therefore, the power series 
$\di \sum_{n=0}^{\infty}(-1)^n\frac{x^{2n}}{(2n+1)!}$ converges for all $x$, and when $x\neq 0$, it is equal to $f(x)$. When $x=0$, it has value 1, which is equal to $f(0)$. This proves that 
\[f(x)=\sum_{n=0}^{\infty}(-1)^n\frac{x^{2n}}{(2n+1)!}\hspace{1cm}\text{for all}\;x\in\mathbb{R}.\]
Since the function $f(x)$ has a power series expansion that converges everywhere, it is an infinitely differentiable function. From the power series expansion, we find that
\[f^{(2n+1)}(0)=0,\hspace{1cm}f^{(2n)}(0)=\frac{(-1)^n}{2n+1}\hspace{1cm}\text{for all}\;n\geq 0.\]
\end{solution}

An important power series expansion that cannot be derived from the list of Taylor series for elementary functions is the binomial series. Recall that if $n$ is a positive integer, 
the binomial expansion of $(1+x)^n$ is given by
\[(1+x)^n=\sum_{k=0}^n\binom{n}{k}x^k.\]If $n$ is a negative integer, let $m=-n-1$. Then $m$ is a nonnegative integer. By Theorem \ref{230307_1}, we find that when $|x|<1$.
\[
(1+x)^n =\frac{1}{(1+x)^{m+1}}=\sum_{k=m}^{\infty}\binom{k}{m}(-x)^{k-m} =\sum_{k=0}^{\infty}(-1)^{k}\binom{k+m}{m}x^k.\]
Notice that for $k\geq 0$,
\begin{align*}
(-1)^{k}\binom{k+m}{m}
&=(-1)^{k}\frac{(k+m)!}{k!m!}\\
&=(-1)^{k}\frac{(m+k)(m+k-1)\cdots (m+1)}{k!}\\
&=\frac{(-m-1)(-m-2)\cdots (-m-k)}{k!}\\
&=\frac{n(n-1)\cdots (n-k+1)}{k!}.
\end{align*}Thus, when $n$ is a negative integer, we find that $(1+x)^n$ has a power series expansion on the interval $(-1,1)$,  which can be written as
\[(1+x)^n=\sum_{k=0}^{\infty}\frac{n(n-1)\cdots (n-k+1)}{k!}x^k.\]
This motivates us to extend the definition of the binomial coefficients.
\begin{definition}{Generalized Binomial Coefficients}
For any real number $\alpha$ and any nonnegative integer $k$, we define the generalized binomial coefficient $\di\binom{\alpha}{k}$ by
\[\binom{\alpha}{0}=1,\] and for $k\geq 1$,
\[\binom{\alpha}{k}=\frac{\alpha(\alpha-1)\cdots(\alpha-k+1)}{k!}.\]
\end{definition}

\begin{example}[label=230307_9]{}
Let $\alpha$ be a real number. Show that the Maclaurin series of the function $f:(-1,1)\to\mathbb{R}$, $f(x)=(1+x)^{\alpha}$ is
\[\sum_{k=0}^{\infty}\binom{\alpha}{k}x^k.\]
When $\alpha$ is not a nonnegative integer, show that the radius of convergence of this power series is 1.
\end{example}\begin{solution}{Solution}
The function $f$ is infinitely differentiable. By straightforward computation, we have
\[f^{(k)}(x)=\alpha(\alpha-1)\cdots (\alpha-k+1)(1+x)^{\alpha-k}.\]
This gives, 
\[f^{(k)}(0)=\alpha(\alpha-1)\cdots (\alpha-k+1).\]
Therefore, the Maclaurin series of $f$ is
\[\sum_{k=0}^{\infty}\frac{f^{(k)}(0)}{k!}x^k=\sum_{k=0}^{\infty}\frac{\alpha(\alpha-1)\cdots (\alpha-k+1)}{k!}x^k=\sum_{k=0}^{\infty}\binom{\alpha}{k}x^k.\]
For the radius of convergence, we note that
\[c_k=\binom{\alpha}{k}=\frac{\alpha(\alpha-1)\cdots(\alpha-k+1)}{k!}\] is nonzero for all $k\geq 0$ when $\alpha$ is not a nonnegative integer. Thus, we can apply ratio test. Since
\[\lim_{k\to\infty}\left|\frac{c_{k+1}}{c_k}\right|=\lim_{k\to\infty}\left|\frac{\alpha -k}{k+1}\right|=1,\]
we find that the radius of convergence of the power series is 1.
\end{solution}
In the example above, we have shown that the Maclaurin  series of $f(x)=(1+x)^{\alpha}$ is $\di \sum_{k=0}^{\infty}\binom{\alpha}{k}x^k$, which is a power series that converges on $(-1,1)$. But we have not shown that the Maclaurin series converges to $f(x)$ on $(-1,1)$, except when $\alpha$ is an integer. To prove the convergence of the Maclaurin series to the function, we will study the convergence of the sequence of partial sums.

The  partial sums of Taylor series are called Taylor polynomials. They are important in their own right.

\begin{definition}{Taylor Polynomials}
Let $I$ be an open interval that contains the point $x_0$, and let $n$ be a positive integer.  If  the function $f:I\to\mathbb{R}$ is $n$ times differentiable on $I$,  the $n^{\text{th}}$  {\bf Taylor polynomial} of $f(x)$ at $x_0$ is the polynomial
\begin{align*}T_n(x)&=\sum_{k=0}^{n}\frac{f^{(k)}(x_0)}{k!}(x-x_0)^k\\&=f(x_0)+\frac{f'(x_0)}{1!}(x-x_0) +\cdots+\frac{f^{(n)}(x_0)}{n!}(x-x_0)^n.\end{align*}
\end{definition}
In particular,
\begin{align*}
T_1(x)&=f(x_0)+f'(x_0)(x-x_0),\\
T_2(x)&=f(x_0)+f'(x_0)(x-x_0)+\frac{f''(x_0)}{2}(x-x_0)^2,\\
T_3(x)&=f(x_0)+f'(x_0)(x-x_0)+\frac{f''(x_0)}{2}(x-x_0)^2+\frac{f'''(x_0)}{6}(x-x_0)^3, 
\end{align*}and so on.  

 Notice that to define Taylor polynomials of degree $n$ for a function $f$, we do not need to assume that $f$ is infinitely differentiable. We just need to assume that $f$ is $n$ times diferentiable. 

\begin{example}
{} For the function $f(x)=x\cos x$, its Taylor series is 
\[f(x)=x\sum_{n=0}^{\infty}(-1)^n\frac{x^{2n}}{(2n)!}=x\left(1-\frac{x^2}{2}+\frac{x^4}{24}+\cdots\right)=x-\frac{x^3}{2}+\frac{x^5}{24}+\cdots.\]
If $T_n(x)$ is the $n^{\text{th}}$ Taylor polynomial for $f(x)$ at $x=0$, then
\begin{align*}T_1(x) &=T_2(x)=x,\\
T_3(x) &=T_4(x)=x-\frac{x^3}{2},\\
T_5(x)&=T_6(x)=x-\frac{x^3}{2}+\frac{x^5}{24},
\end{align*}and so on.
\end{example}
In this example, we notice that $T_{2n-1}(x)=T_{2n}(x)$ for all $n\in\mathbb{Z}^+$.  This is because   $f(x)$ is an odd function.

We have noticed that the first Taylor polynomial \[T_1(x)=f(x_0)+f'(x_0)(x-x_0)\] of a function $f(x)$ at $x=x_0$ is related to the tangent line to the graph of the function. 
In fact, by definition of derivatives, we have
\[\lim_{x\to x_0}\frac{f(x)-T_1(x)}{x-x_0}=\lim_{x\to x_0}\frac{f(x)-f(x_0)-(x-x_0)f'(x_0)}{x-x_0}=0.\]We say that $T_1(x)$ is a first order approximation of $f(x)$ at $x=x_0$. 
In general, we can define the following concept.
\begin{definition}{Order of Approximation}
Let $I$ be an open interval that contains the point $x_0$, and let $n$ be a positive integer. We say that two functions $f:I\to\mathbb{R}$ and $g: I\to\mathbb{R}$ are $n^{\text{th}}$-order approximations of each other at the point $x_0$ if
\[\lim_{x\to x_0}\frac{f(x)-g(x)}{(x-x_0)^n}=0.\]
\end{definition}
We will show that the $n^{\text{th}}$ Taylor polynomial of a function $f(x)$ at $x_0$ is an $n^{\text{th}}$-order approximation of the function at $x_0$. First, we prove 
the following lemma which says that for any real number $x_0$, any polynomial of degree $n$ can be written in the form$\di\sum_{k=0}^n c_k(x-x_0)^k$. 

\begin{lemma}[label=230307_4]{}
Given a real number $x_0$, and a polynomial $p(x)$ of degree $n$, we have
\[p(x)=\sum_{k=0}^n\frac{p^{(k)}(x_0)}{k!}(x-x_0)^k.\]In other words, the $n^{\text{th}}$ Taylor polynomial of $p(x)$ is $p(x)$ itself, and the Taylor series of $p(x)$ is also $p(x)$.
\end{lemma}
\begin{myproof}{Proof}
Let
\[p(x)=a_0+a_1x+\cdots+a_nx^n,\]
and let $h=x-x_0$. Then   $x=h+x_0$. Substitute $x$ by $x_0+h$, we have
\[p(x)=a_0+a_1(h+x_0)+\cdots+a_n(h+x_0)^n.\]
For $0\leq k\leq n$, $(h+x_0)^k$ is a polynomial of degree $k$ in $h$. Thus,
\[a_0+a_1(h+x_0)+\cdots+a_n(h+x_0)^n\] 
 is a polynomial of degree $n$ in $h$. This implies that there are constants $c_0$, $c_1$, $\ldots$, $c_n$ such that
\[p(x)= \sum_{k=0}^n c_kh^k=\sum_{k=0}^nc_k(x-x_0)^k.\]
Differentiate both sides $k$ times and set $x=x_0$ gives
\[p^{(k)}(x_0)=k!c_k.\]
This proves that
\[p(x)=\sum_{k=0}^n\frac{p^{(k)}(x_0)}{k!}(x-x_0)^k.\]

\end{myproof}

As  a corollary, we have the following, which can be deduced from Theorem \ref{thm230216_11}.
\begin{corollary}[label=230307_18]{}
If $p(x)$ is a polynomial of degree at most $n$, and there is a point $x_0$ such that
\[p(x_0)=p'(x_0)=\cdots =p^{(n)}(x_0)=0,\]
then $p(x)$ is identically zero.
\end{corollary} 

We would also like to emphasize again the following.
\begin{corollary}[label=230307_5]{}
Let $I$ be an interval that contains the point $x_0$, and let $n$ be a positive integer. Given that $f:I\to\mathbb{R}$   is a function that is $n$ times differentiable, let \[T_n(x)=\di\sum_{k=0}^n\frac{f^{(k)}(x_0)}{k!}(x-x_0)^k\]be its $n^{\text{th}}$ Taylor polynomial  at $x_0$. For $0\leq k\leq n$, we have
\[T_n^{(k)}(x_0)=f^{(k)}(x_0).\]
\end{corollary}
\begin{myproof}{Proof}
By Lemma \ref{230307_4}, we have
\[ T_n(x)=\di\sum_{k=0}^n\frac{T_n^{(k)}(x_0)}{k!}(x-x_0)^k.\]
The result follows by comparing coefficients.
\end{myproof}

Now we prove the approximation theorem.
\begin{theorem}[label=230307_2]{}Let $I$ be an open interval that contains the point $x_0$, and let $n$ be a positive integer. Assume that the function $f:I\to\mathbb{R}$   is $n$ times differentiable.
\begin{enumerate}[(a)]
\item The $n^{\text{th}}$ Taylor polynomial \[T_n(x)=\sum_{k=0}^n\frac{f^{(k)}(x_0)}{k!}(x-x_0)^k\] of $f(x)$ at $x_0$ is an $n^{\text{th}}$-order approximation of $f(x)$ at $x_0$.
\item If $p(x)$ is a polynomial of degree at most $n$, and $p(x)$  is an $n^{\text{th}}$-order approximation of $f(x)$ at $x_0$, then $p(x)=T_n(x)$.

\end{enumerate}
\end{theorem}
\begin{myproof}{Proof}Let us  consider (a) first.
If $n=1$, we need to show that
\[\lim_{x\to x_0}\frac{f(x)-T_1(x)}{x-x_0}=\lim_{x\to x_0}\frac{f(x)-f(x_0)-f'(x_0)(x-x_0)}{x-x_0}=0.\] 
But this is just the definition of $f'(x_0)$.
  Assume that we have proved the statement for the $n-1$ case. Now we look at
\[\lim_{x\to x_0}\frac{f(x)-T_n(x)}{(x-x_0)^n}.\]This is a limit of the indeterminate form $0/0$. 

Notice that \[T_n'(x)=\sum_{k=1}^n\frac{f^{(k)}(x_0)}{(k-1)!}(x-x_0)^{k-1}=\sum_{k=0}^{n-1}\frac{(f')^{(k)}(x_0)}{k!}(x-x_0)^k\] 
is the $(n-1)^{\text{th}}$ Taylor polynomial for $f'$, and $f'$ is $(n-1)$ times differentiable. By inductive hypothesis,
\[\lim_{x\to x_0}\frac{f'(x)-T_n'(x)}{(x-x_0)^{n-1}}=0.\]
By   l' H$\hat{\text{o}}$pital's rule,
\[\lim_{x\to x_0}\frac{f(x)-T_n(x)}{(x-x_0)^n}=\lim_{ x\to x_0}\frac{f'( x)-T_n'(x)}{n(x-x_0)^{n-1}}=0.\]This finishes the induction for (a).
 
Now we consider (b). 
Let $\di p(x)=\sum_{k=0}^nc_k(x-x_0)^k$ be a polynomial of degree at most $n$ which  is an $n^{\text{th}}$-order approximation of $f(x)$ at $x_0$. Then
\[\lim_{x\to x_0}\frac{f(x)-p(x)}{(x-x_0)^n}=0.\]
We have proved in part (a) that
\[\lim_{x\to x_0}\frac{f(x)-T_n(x)}{(x-x_0)^n}=0.\] \bp
These give
\[\lim_{x\to x_0}\frac{T_n(x)-p(x)}{(x-x_0)^n}=\lim_{x\to x_0}\frac{f(x)-p(x)}{(x-x_0)^n}-\lim_{x\to x_0}\frac{f(x)-T_n(x)}{(x-x_0)^n}=0.\]
It follows that for all $0\leq k\leq n$, 
\begin{equation}\label{eq230306_9}\lim_{x\to x_0}\frac{T_n(x)-p(x)}{(x-x_0)^k}=0.\end{equation}
Notice that
\[T_n(x)-p(x)=\sum_{k=0}^n \left(\frac{f^{(k)}(x_0)}{k!}-c_k\right)(x-x_0)^k.\]
Take $k=0$ in \eqref{eq230306_9}, we find that
$c_0=f(x_0)$. Then take $k=1$ shows that $c_1=f'(x_0)$. Inductively, we show that $c_k=\di\frac{f^{(k)}(x_0)}{k!}$ for all $0\leq k\leq n$. This completes the proof of the theorem.
\end{myproof}

\begin{highlight}{}
Theorem \ref{230307_2} says that the $n^{\text{th}}$ Taylor polynomial  of a function $f(x)$ at a point $x_0$ is the unique polynomial of degree at most $n$ which is an $n^{\text{th}}$-order approximation of $f(x)$ at $x_0$. Hence, we also called $T_n(x)$ the Taylor polynomial of $f(x)$ of order $n$ at $x_0$. One should avoid calling it the $n^{\text{th}}$ degree Taylor polynomial as we have seen that $T_n(x)$ does not necessary have degree $n$. 
\end{highlight}

Given that $I$ is an interval that contains the point $x_0$, and  $f:I\to\mathbb{R}$   is an $n$ times differentiable function, the Taylor polynomial $T_n(x)$ is well defined. By Theorem \ref{230307_2},
\[\lim_{x\to x_0}\frac{f(x)-T_n(x)}{(x-x_0)^n}=0.\]
This implies that for any $\varepsilon>0$, there is a $\delta>0$ such that $(x_0-\delta, x_0+\delta)\subset I$, and for all $x\in (x_0-\delta, x_0+\delta)$,
\begin{equation}\label{eq230307_3}|f(x)-T_n(x)|\leq \varepsilon |x-x_0|^n.\end{equation}
The function 
\[R_n(x)=f(x)-T_n(x)\] is called the remainder when we approximate the function $f(x)$ by its $n^{\text{th}}$ Taylor polynomial $T_n(x)$ at $x_0$. Eq \eqref{eq230307_3} says that when $x$ approaches $x_0$, the order of $R_n(x)$ is smaller than the order of $|x-x_0|^n$. If we assume that $f$ has one more derivative, we can   say more. 

We will first prove the Lagrange  remainder theorem which  assumes that $f:I\to\mathbb{R}$ is $(n+1)$ times differentiable.  

\begin{theorem}{The Lagrange Remainder Theorem}
Let $I$ be an open interval that contains the point $x_0$, and let $n$ be a positive integer. Given that $f:I\to\mathbb{R}$   is a function that is $(n+1)$ times differentiable, let \[T_n(x)=\di\sum_{k=0}^n\frac{f^{(k)}(x_0)}{k!}(x-x_0)^k\]be its Taylor polynomial of order $n$ at $x_0$. For any $x\in I\setminus \{x_0\}$, there is a number $c\in (0,1)$ such that
\[f(x)-T_n(x)=\frac{f^{(n+1)}(\xi)}{(n+1)!}(x-x_0)^{n+1},\quad\text{where}\;\xi=x_0+c(x-x_0).\]
\end{theorem}Recall that $\xi=x_0+c(x-x_0)$ with $c\in (0,1)$ means that $\xi$ is a point strictly between $x_0$ and $x$. 
\begin{myproof}{Proof}
The proof   is   a straightforward application of Theorem \ref{thm230216_11}, which is a consequence of Cauchy mean value theorem.
Let $g:I\to\mathbb{R}$ be the function defined by
\[g(x)=f(x)-T_n(x)=f(x)-\sum_{k=0}^n\frac{f^{(k)}(x_0)}{k!}(x-x_0)^k.\]Then $g$ is $(n+1)$ times differentiable. Corollary \ref{230307_5} implies that
\[g(x_0)=g'(x_0)=\cdots=g^{(n)}(x_0)=0.\]\bp
Applying Theorem \ref{thm230216_11} to the function $g$, we find that for any $x\in I\setminus \{x_0\}$, there is a number $c\in (0,1)$ such that
\[f(x)-T_n(x)=g(x)=\frac{g^{(n+1)}(\xi)}{(n+1)!}(x-x_0)^{n+1},\quad\text{where}\;\xi=x_0+c(x-x_0).\]  Since $T_n(x)$ is a polynomial of degree $n$, $T_n^{(n+1)}(x)=0$ for all $x\in I$. Therefore, $g^{(n+1)}(x)=f^{(n+1)}(x)$ for all $x\in I$. This concludes the proof.
\end{myproof}
The Lagrange remainder theorem also holds in the $n=0$ case. This is just 
the
Lagrange mean value theorem. Thus Lagrange remainder theorem is an extension of the Lagrange mean value theorem. It gives   useful estimates on the error term in approximating a function by its Taylor polynomial, especially if $f^{(n+1)}(x)$ is always positive or always negative in a neighbourhood of $x_0$.

\begin{example}{}In this example, we demonstrate how we can use the Lagrange remainder theorem to show that the Taylor series of the function $
 f:\mathbb{R}\to\mathbb{R}$, $f(x)=e^x$, converges to $f(x)$ for all real numbers $x$.  The $n^{\text{th}}$ Taylor polynomial of $f(x)=e^x$ at $x=0$ is
\[T_n(x)=1+x+\frac{x^2}{2}+\cdots+\frac{x^n}{n!}.\]Since $f^{(n)}(x)=e^x$ for any $n\in\mathbb{Z}^+$, Lagrange remainder theorem says that for any $n\geq 0$, for any real number $x\neq 0$, there is a number $\xi$ strictly between $0$ and $x$ such that
\begin{equation}\label{eq230307_8}e^x-T_n(x)=\frac{f^{(n+1)}(\xi)}{(n+1)!}x^{n+1}=\frac{e^{\xi}}{(n+1)!}x^{n+1}.\end{equation}
For fixed $x$, $\xi$ depends on $n$ but we can use  $\xi<|x|$ to get the estimate $e^{\xi}\leq e^{|x|}$ that is independent of $n$.
This implies that
\[\left|e^x-T_n(x)\right|\leq  e^{|x|}\frac{|x|^{n+1}}{(n+1)!}.\]\be
We have proved in Example \ref{230305_1} that the power series $\di\sum_{n=0}^{\infty}\frac{|x|^n}{n!}$ is convergent. Therefore, $\di \lim_{n\to\infty} \frac{|x|^{n+1}}{(n+1)!}=0$. This allows us to conclude that
\begin{equation}\label{eq230307_7}e^x=\lim_{n\to\infty}T_n(x)=1+x+\frac{x^2}{2}+\cdots+\frac{x^n}{n!}+\cdots.\end{equation}
This is an alternative way to prove that the Taylor series of $e^x$ converges to $e^x$, instead of the aproach used in the proof of Theorem \ref{230307_6}.
From the series expansion \eqref{eq230307_7}, it is easy to deduce that for all $x>0$, and all $n\geq 1$, 
\[e^x>1+x +\frac{x^2}{2}+\cdots+\frac{x^n}{n!}.\]In particular, we have
\begin{align*}
e^x&>1+x,\\
e^x&>1+x+\frac{x^2}{2},\\
e^x&>1+x+\frac{x^2}{2}+\frac{x^3}{6},\\
e^x&>1+x+\frac{x^2}{2}+\frac{x^3}{6}+\frac{x^4}{24}.
\end{align*}  For $x<0$, the series  \eqref{eq230307_7} is alternating. The sequence $\{b_n\}$ with $\di b_n=\frac{|x|^n}{n!}$ is not decreasing. However, since $x^{2n-1}<0$ and $x^{2n}>0$ for all $n\geq 1$, and $e^{\xi}>0$ for all $\xi$, we can use \eqref{eq230307_8} to conclude that for   $x<0$,
\begin{align*}
e^x&>1+x,\\
e^x&<1+x+\frac{x^2}{2},\\
e^x&>1+x+\frac{x^2}{2}+\frac{x^3}{6},\\
e^x&<1+x+\frac{x^2}{2}+\frac{x^3}{6}+\frac{x^4}{24}.
\end{align*}
\end{example2}

\begin{figure}[ht]
\centering
\includegraphics[scale=0.2]{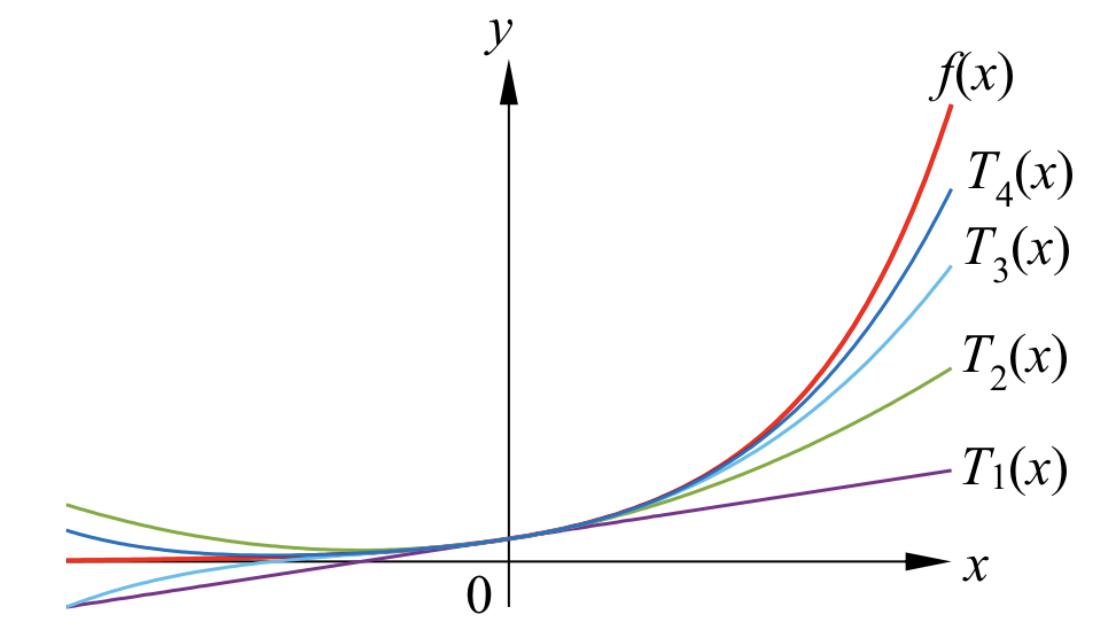}
\caption{The function $f(x)=e^x$ and its Taylor polynomials at $x=0$.\fa}\label{figure59}
\end{figure}

In Theorem \ref{230307_8}, we apply  mean value theorem to prove that $|\sin x|\leq |x|$ for all real numbers $x$. In the following example, we extend this result partially.
\begin{example}{}
Show that
for $x\in (0, \pi)$,  $\sin x>x-\di\frac{x^3}{6}$.
\end{example}
\begin{solution}{Solution}
Let $f(x)=\sin x$. Then $f(x)$ is infinitely differentiable, with the third  Taylor polynomial at $x=0$ given by
\[ T_3(x)= x-\frac{x^3}{6}.\]
Apply the Lagrange remainder theorem, we find that for any $x\in (0,\pi)$, there is a $\xi \in (0,x)\subset (0,\pi)$ so that
\[\sin x-x+\di\frac{x^3}{6}=\frac{f^{(4)}(\xi)}{24}x^4=\frac{\sin \xi}{24}x^4.\] Since $\sin\xi>0$ for $\xi\in (0,\pi)$,
 this proves that \[\sin x>x-\di\frac{x^3}{6}\hspace{1cm}\text{for}\;x\in (0,\pi).\]
\end{solution}

\begin{figure}[ht]
\centering
\includegraphics[scale=0.2]{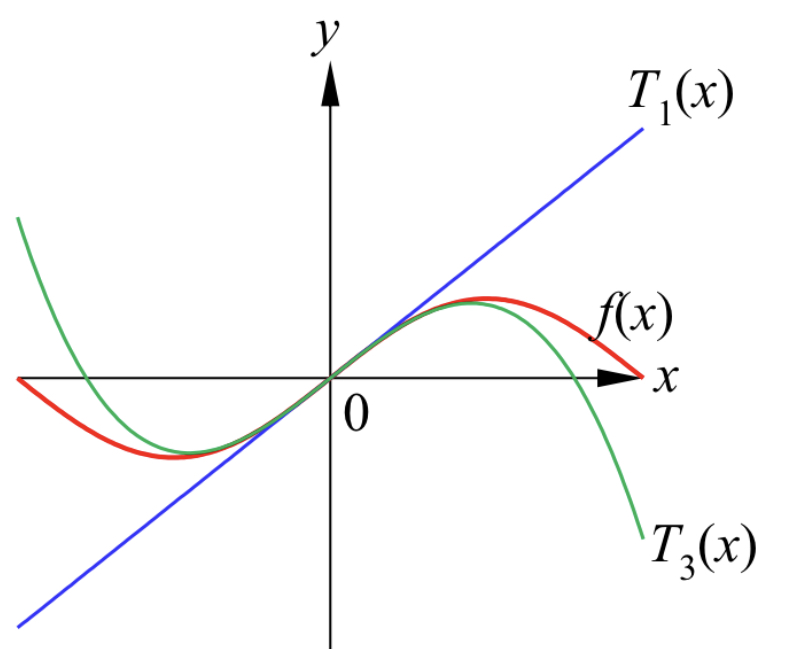}
\caption{The function $f(x)=\sin x$ and its Taylor polynomials at $x=0$.\fa}\label{figure60}
\end{figure}

  We have repeatedly used the fact that if $f:I\to\mathbb{R}$ is a differentiable function defined on an open interval $I$, and $f'(x)=0$ for all $x\in I$, then $f(x)$ is a constant function. The next theorem extends this result.
\begin{theorem}[label=230307_17]{}
Let $I$ be an open interval, and let $n$ be a positive integer. Assume that the function $f:I\to\mathbb{R}$ is $(n+1)$ times differentiable, and
$f^{(n+1)}(x)=0$ for all $x\in I$. Then $f(x)$ is a polynomial of degree at most $n$.
\end{theorem} 
\begin{myproof}{Proof}If $f^{(n+1)}(x)=0$ for all $x\in I$, take any point $x_0$  in $I$, and let $T_n(x)$ be the $n^{\text{th}}$ Taylor polynomial of $f$ at $x_0$. By definition, $f(x_0)=T_n(x_0)$. Given $x\in I\setminus\{x_0\}$, the Lagrange remainder theorem implies that there is a point $\xi\in I$ such that
\[f(x)-T_n(x)=\frac{f^{(n+1)}(\xi)}{(n+1)!}(x-x_0)^{n+1}.\]
Since $f^{(n+1)}(x)$ is identically 0, we find that for all $x\in I$, 
$f(x)=T_n(x)$. This proves that $f$ is a polynomial of degree at most $n$.

\end{myproof}

As a corollary, we have the following.
\begin{corollary}{}
Let $I$ be an open interval, and let $n$ be a positive integer. Assume that $f:I\to\mathbb{R}$ and $g:I\to\mathbb{R}$ are $(n+1)$ times differentiable functions such that
\[f^{(n+1)}(x)=g^{(n+1)}(x)\hspace{1cm}\text{for all}\;x\in I,\] then there is a polynomial $p(x)$ of degree at most $n$ such that
\[f(x)=g(x)+p(x).\]
\end{corollary}

Next we turn to the Cauchy remainder theorem. In Example \ref{230307_10}, we have shown that if $g:I\to\mathbb{R}$ is a continuous function, $x_0$ is a point in $I$,  $n$ is a positive integer, then  the function $G:I\to\mathbb{R}$ defined by
\[G(x)=\frac{1}{n!}\int_{x_0}^x(x-t)^ng(t)dt \]
  is $(n+1)$ times continuously differentiable,  
\[G(x_0)=G'(x_0)=\ldots=G^{(n)}(x_0)=0,\]
and
\[G^{(n+1)}(x)=g(x)\hspace{1cm}\text{for all}\;x\in I.\]
\begin{theorem}{The Cauchy Remainder Formula}
Let $I$ be an open interval that contains the point $x_0$, and let $n$ be a positive integer. Given that $f:I\to\mathbb{R}$   is a function that is $(n+1)$ times continuously differentiable, let \[T_n(x)=\di\sum_{k=0}^n\frac{f^{(k)}(x_0)}{k!}(x-x_0)^k\]be its Taylor polynomial of order $n$ at $x_0$. For any $x\in I$,  
\[f(x)-T_n(x)=\frac{1}{n!}\int_{x_0}^x (x-t)^n f^{(n+1)}(t)dt.\]

\end{theorem}
\begin{myproof}{Proof}
Let $H: I\to\mathbb{R}$ be the function defined by
\[H(x)=f(x)-T_n(x)-\frac{1}{n!}\int_{x_0}^x (x-t)^n f^{(n+1)}(t)dt.\]
Then by the result proved in Example \ref{230307_10}, we find that $H$ is a function that is $(n+1)$ times continuously differentiable,
\[H(x_0)=H'(x_0)=\cdots=H^{(n)}(x_0)=0,\]
and
\[H^{(n+1)}(x)=f^{(n+1)}(x)-f^{(n+1)}(x)=0\hspace{1cm}\text{for all}\;x\in I.\]
By Theorem \ref{230307_17} and Corollary \ref{230307_18}, $H(x)=0$ for all $x\in I$. This completes the proof of the assertion.
\end{myproof}

In Cauchy remainder formula, the error term is expressed as a precise integral, although in practice it might not be possible to evaluate such an integral. Let us now apply the Cauchy remainder formula to prove that the Maclaurin series of the function $f(x)=(1+x)^{\alpha}$ converges to $f(x)$ when $x\in (-1,1)$.
\begin{theorem}
{}Let $\alpha$ be a real number. For $|x|<1$,
\[(1+x)^{\alpha}= \sum_{k=0}^{\infty}\binom{\alpha}{k}x^k.\]
\end{theorem}
\begin{myproof}{Proof}Let $f(x)=(1+x)^{\alpha}$, $-1<x<1$. We have seen that $\di \sum_{k=0}^{\infty}\binom{\alpha}{k}x^k$ is the Maclaurin series of $f(x)$.
The $n^{\text{th}}$ Taylor polynomial of $f(x)$ at $x=0$ is 
\[T_n(x)=\sum_{k=0}^{n}\binom{\alpha}{k}x^k.\]\bp
We need to show that $\di \lim_{n\to\infty}T_n(x)=f(x)$ for all $|x|<1$.
It is easy to verify that
\[f^{(n+1)}(x)= (n+1)! \binom{\alpha}{n+1}(1+x)^{\alpha-n-1}.\]
By  Cauchy remainder formula,  for $x\in (-1,1)$,
\begin{align*}f(x)-T_n(x)&=\frac{1}{n!}\int_0^x (x-t)^nf^{(n+1)}(t)dt\\&=(n+1)\binom{\alpha}{n+1}\int_0^x(x-t)^n(1+t)^{\alpha-n-1}dt.\end{align*}
We  need to  estimate this last integral for $x\neq 0$. Making a change of variables $t=x\tau$, we find that
\[\int_0^x(x-t)^n(1+t)^{\alpha-n-1}dt=x^{n+1}\int_0^1 (1-\tau)^n(1+x\tau)^{\alpha-n-1}d\tau.\]
Notice that since $x\in (-1,1)$, when $\tau\in [0,1]$, 
\[(1-\tau)^n(1+x\tau)^{\alpha-n-1}\geq 0.\]  For fixed $x\in (-1,1)$, the function $g:[0,1]\to\mathbb{R}$, $g(\tau)=(1+x\tau)^{\alpha-1}$ is continuous.   Therefore, there is a constant $M$ such that 
\[0\leq (1+x\tau)^{\alpha-1}\leq M\hspace{1cm}\text{for all}\;\tau\in [0,1].\]
This implies that
\[0\leq \int_0^1 (1-\tau)^n(1+x\tau)^{\alpha-n-1}d\tau\leq M\int_0^1\left(\frac{1-\tau}{1+x\tau}\right)^{n}d\tau.\]
For any $x\in (-1,1)$ and $\tau\in [0,1]$,
\[1+x\tau\geq 1-\tau\geq 0.\]
 This  implies that for any $x\in (-1,1)$, 
\[0\leq \left(\frac{1-\tau}{1+x\tau}\right)^{n} \leq 1\hspace{1cm}\text{for all}\;\tau\in [0,1].\]\bp
Therefore, we find that
\[0\leq \int_0^1\left(\frac{1-\tau}{1+x\tau}\right)^{n}d\tau\leq 1\hspace{1cm}\text{for all}\;|x|<1.\]Hence,
\begin{equation}\label{eq230307_13}|f(x)-T_n(x)|\leq M (n+1)\left|\binom{\alpha}{n+1}x\right|^{n+1}.\end{equation}
In Example \ref{230307_9}, we have proved that the series $\di\sum_{n=0}^{\infty}\binom{\alpha}{n}x^n$ is convergent when $|x|<1$. It follows from Theorem \ref{230305_15} that the derived series $\di\sum_{n=0}^{\infty}(n+1)\binom{\alpha}{n+1}x^n$ is also convergent when $|x|<1$. This implies that 
\[\lim_{n\to\infty}(n+1)\left|\binom{\alpha}{n+1}x\right|^n=0.\]
Using squeeze theorem, we deduce from \eqref{eq230307_13} that
\[\lim_{n\to\infty}T_n(x)=f(x).\] This completes the proof.
\end{myproof}

\vp
\noindent
{\bf \large Exercises  \thesection}
\setcounter{myquestion}{1}
\begin{question}{\themyquestion}
Let $f:\mathbb{R}\to\mathbb{R}$ be the function defined by
\[f(x)=\begin{cases} \di\frac{2-2\cos x}{x^2},\quad &\text{if}\;x\neq 0,\\1,\quad & \text{if} \;x=0.\end{cases}\]Show that  $f$ is infinitely differentiable, and find $f^{(n)}(0) $ for all $n\geq 0$.
\end{question}

\atc
\begin{question}{\themyquestion}
Show that for all $x>0$,
\[x-\frac{x^2}{2}<\ln (1+x)<x.\]
\end{question}
 
\atc
\begin{question}{\themyquestion}
Show that for all $x>0$,
\[1+\frac{x}{2}-\frac{x^2}{8}<\sqrt{1+x}<1+\frac{x}{2}-\frac{x^2}{8}+\frac{x^3}{16}.\]
\end{question}

\atc
\begin{question}{\themyquestion}
Show that for all $x\in (-\pi, \pi)$,
\[\cos x\geq 1-\frac{x^2}{2}.\]
\end{question}

\atc
\begin{question}{\themyquestion}
Let $\alpha$ be a real number. Assume that $\alpha$ is not an integer. In Example \ref{230307_9}, we have shown that the power series $\di\sum_{k=0}^{\infty}\binom{\alpha}{k}x^k$,
which is the Maclaurin series of the function $f(x)=(1+x)^{\alpha}$, has radius of convergence 1. Define the function 
$g:(-1,1)\to\mathbb{R}$ by
\[g(x)=\sum_{k=0}^{\infty}\binom{\alpha}{k}x^k.\]
In this question, you are asked to show that $g(x)=(1+x)^{\alpha}$ for $x\in (-1,1)$, without using the Cauchy remainder formula.
\begin{enumerate}[(a)]\item
Show that $(1+x)g'(x)=\alpha g(x)$ for all $x\in (-1,1)$.
\item Let $h:(-1,1)\to\mathbb{R}$ be the function defined by $h(x)=g(x)(1+x)^{-\alpha}$. Prove that $h$ is a constant function.
\item Conclude that $g(x)=(1+x)^{\alpha}$ for all $x\in (-1,1)$.
\end{enumerate}
\end{question}

\vp

\section{Examples and Applications}\label{sec6.6}
In this section, we discuss some examples and applications.

The number $e$ and the number $\pi$ are two important numbers in mathematics. 
In Section \ref{sec6.6.1} and Section \ref{sec6.6.2}, we prove respectively that these two numbers are irrational. 

In Section \ref{sec6.6.3}, we prove that there is an infinitely differentiable function whose Taylor series at a point does not converge to the function itself. We also briefly discuss the applications of such functions, despite its non-analyticity.

In Section \ref{sec6.6.4}, we construct a continuous function that is differentiable nowhere. It uses Theorem \ref{230304_1} which says that uniform limit of continuous functions is continuous.

In Section \ref{sec6.6.5}, we prove the Weierstrass approximation theorem, which says that any continuous function defined on a closed and bounded interval can be uniformly approximated by a polynomial. We give a proof that uses Bernstein's approach. It uses the fact that a continuous function defined on a closed and bounded interval is bounded and uniformly continuous.  Later when we study Fourier series, we are going to prove this important theorem again using the theory of Fourier series.

\subsection[The Irrationality of $e$]{The Irrationality of $\pmb{e}$}\label{sec6.6.1}
In Example \ref{23020507}, we have defined the number $e$ as the limit of the increasing sequence $\{a_n\}$, where $\di a_n=\left(1+\frac{1}{n}\right)^n$. We have proved that $a_n\leq 3$ for all $n\in\mathbb{Z}^+$. This implies that $e\leq 3$. In Theorem \ref{230307_6}, we proved that
\[e=\sum_{n=0}^{\infty}\frac{1}{n!}=1+\frac{1}{1!}+\frac{1}{2!}+\cdots+\frac{1}{n!}+\cdots.\]

\begin{theorem}{Irrationality of $\pmb{e}$}
The number $e$ is irrational.
\end{theorem}
\begin{myproof}{Proof}Assume to the contrary that $e$ is rational. Then since $e$ is positive, there are positive integers $a$ and $b$ such that
\[e=\frac{a}{b}.\]
For any positive integer $n$, we apply  Lagrange remainder theorem to the $n^{\text{th}}$ Taylor polynomial of $e^x$ at the point $x_0=0$. With $x=1$, we find that there is a number $c_n$ in the interval $(0,1)$ such that
\begin{equation}\label{eq230307_20}e=1+\frac{1}{1!}+\frac{1}{2!}+\cdots+\frac{1}{n!}+\frac{e^{c_n}}{(n+1)!}.\end{equation} 
 
For $n\geq b$, we find that
$n!$ is divisible by $b$, and so $n!e$ is an integer. From \eqref{eq230307_20}, we have
\[0<n!e-\left(n!+n!+\frac{n!}{2!}+\cdots+\frac{n!}{n!}\right)=\frac{e^{c_n}}{n+1}<\frac{e}{n+1}\leq  \frac{3}{n+1}.\]
Notice that for each $1\leq k\leq n$, $n!/k!$ is an integer. Hence, for $n\geq b$, 
\[n!e-\left(n!+n!+\frac{n!}{2!}+\cdots+\frac{n!}{n!}\right)\] is a positive integer that is less than $3/(n+1)$. For $n>3$, $3/n+1$ is less than 1. This gives a contradiction. Hence, $e$ must be irrational.
\end{myproof}

\bigskip
\subsection[The Irrationality of $\pi$]{The Irrationality of $\pmb{\pi}$}\label{sec6.6.2}

As in the case of the number $e$, we will show that $\pi$ is an irrational number using proof by contradiction. 
We begin by two lemmas.
\begin{lemma}[label=230309_1]{}
Given that $f:\mathbb{R}\to\mathbb{R}$ and $g:\mathbb{R}\to\mathbb{R}$ are two infinitely differentiable functions. For any $n\in\mathbb{Z}^+$, and any numbers $\alpha$ and $\beta$,
\begin{equation}\label{eq230308_1}\begin{split}
&\int_{\alpha}^{\beta} f^{(2n+1)}(x)g(x)dx+\int_{\alpha}^{\beta} f(x)g^{(2n+1)}(x)dx\\ &=\sum_{k=0}^{2n}(-1)^k f^{(k)}(\beta)g^{(2n-k)}(\beta)-\sum_{k=0}^{2n}(-1)^k f^{(k)}(\alpha)g^{(2n-k)}(\alpha).
\end{split}\end{equation}

\end{lemma}
\begin{myproof}{Proof}Given  $n\in\mathbb{Z}^+$,
define the function $F:\mathbb{R}\to\mathbb{R}$ by
\[F (x)=\sum_{k=0}^{2n}(-1)^k f^{(k)}(x)g^{(2n-k)}(x).\]
Then
\[
F'(x) =\sum_{k=0}^{2n}(-1)^kf^{(k+1)}(x)g^{(2n-k)}(x)+\sum_{k=0}^{2n}(-1)^kf^{(k)}(x)g^{(2n-k+1)}(x)\]Because of the alternating signs, the $k=0$ to $k=2n-1$ terms in the first sum cancel with the $k=1$ to $k=2n$ terms in the second sum.  This gives
\[
F'(x)=f^{(2n+1)}(x)g(x)+f(x)g^{(2n+1)}(x).
\]  By fundamental theorem of calculus, we find that  
\[
F(\beta)-F(\alpha)=\int_{\alpha}^{\beta} f^{(2n+1)}(x)g(x)dx+\int_{\alpha}^{\beta} f(x)g^{(2n+1)}(x)dx.
\]This proves \eqref{eq230308_1}.
\end{myproof}

\begin{lemma}[label=230309_2]{}
Let $a$, $b$ and $n$ be positive integers. Define the polynomial $p:\mathbb{R}\to\mathbb{R}$ by
\[p(x)=\frac{x^n(a-bx)^n}{n!}.\]
For any  integer $k$ satisfying $0\leq k\leq 2n$, $p^{(k)}(0)$ and $p^{(k)}(a/b)$ are integers.
\end{lemma}
\begin{myproof}{Proof} 
Using binomial expansion, we have
\[p(x)=\frac{x^n}{n!}\sum_{m=0}^n \binom{n}{m}a^{n-m}(-1)^mb^mx^m.\]
By Lemma \ref{230307_4},
\[p(x)=\sum_{k=0}^{2n}\frac{p^{(k)}(0)}{k!}x^k.\] 
Comparing the coeficients, we find that   
\begin{align*}
p^{(k)}(0)=\begin{cases}0,\quad &\text{if}\; \;0\leq k\leq n-1,\\ \di
 \frac{k!}{n!}\binom{n}{k-n}a^{2n-k}(-1)^{k-n}b^{k-n},\quad &\text{if}\;\;n\leq k\leq 2n.\end{cases}\end{align*}Since $k!$ is divisible by $n!$ when $k\geq n$, we find  that $p^{(k)}(0)$ is an integer for all $0\leq k\leq 2n$.
Expanding $p(x)$ in powers of $(x-a/b)$, we find that
\begin{align*}
p(x)=(-1)^n\frac{b^n}{n!}\left(x-\frac{a}{b}\right)^n\sum_{m=0}^n\binom{n}{m}\left(\frac{a}{b}\right)^{n-m}\left(x-\frac{a}{b}\right)^m.
\end{align*}
By Lemma \ref{230307_4},
\[p(x)=\sum_{k=0}^{2n}\frac{p^{(k)}(a/b)}{k!}\left(x-\frac{a}{b}\right)^k.\]  \bp
Comparing the coeficients, 
we find that
\begin{align*}p^{(k)}\left(\frac{a}{b}\right)=\begin{cases}0,\quad &\text{if}\; \;0\leq k\leq n-1,\\ \di (-1)^n\frac{k!}{n!}\binom{n}{k-n}a^{2n-k}b^{k-n},\quad &\text{if}\;\;n\leq k\leq 2n.\end{cases}\end{align*}
Hence, $p^{(k)}(a/b)$ is also an integer  for all $0\leq k\leq 2n$.
\end{myproof}

Now we can prove the theorem.
\begin{theorem}{Irrationality of $\pmb{\pi}$}
The number $\pi$ is irrational.

\end{theorem}
\begin{myproof}{Proof}

Assume that $\pi$ is a rational number. Then there are positive integers $a$ and $b$ such that \[\pi =\frac{a}{b}.\] For   $n\in\mathbb{Z}^+$, define the polynomial $p_{n}(x)$ by
\[p_n(x)=\frac{x^n(a-bx)^n}{n!}=\frac{b^nx^n(\pi-x)^n}{n!},\]and let
\[I_n=\int_0^{\pi}p_n(x)\sin xdx.\]
Take $f(x)=p_n(x)$, $g(x)=\cos x$ and $\alpha=0$, $\beta=\pi$ in Lemma \ref{230309_1}. 
Since $p_n(x)$ is a polynomial of degree $2n$, we find that $p_n^{(2n+1)}(x)=0$. On the other hand, for all $k\geq 0$,
\begin{gather*}
g^{(4k)}(x)=\cos x, \quad g^{(4k+1)}(x)=-\sin x,\\g^{(4k+2)}(x)=-\cos x,\quad g^{(4k+3)}(x)=\sin x.\end{gather*}\bp
In particular,  $g^{(2n+1)}(x)=(-1)^{n-1}\sin x$. 
 From \eqref{eq230308_1}, we have
\begin{equation}\label{eq230309_3}\begin{split}
I_n &= (-1)^{n-1}\left\{\sum_{k=0}^{2n}(-1)^k p_n^{(k)}(\pi)g^{(2n-k)}(\pi)\right.\\&\hspace{3cm}\left.-\sum_{k=0}^{2n}(-1)^k p_n^{(k)}(0)g^{(2n-k)}(0)\right\}.\end{split}
\end{equation}By Lemma \ref{230309_2}, $p_n^{(k)}(0)$ and $p_n^{(k)}(\pi)$ are  integers for all $0\leq k\leq 2n$. Since
\[\sin 0=0,\quad \sin \pi =0, \quad \cos 0=1,\quad \cos \pi =-1\] are   integers, we find that $g^{(k)}(0)$ and $g^{(k)}(\pi)$ are  integers for all $k\geq 0$. The right hand side of \eqref{eq230309_3} shows  that $I_n$ is an integer for all $n\in\mathbb{Z}^+$. 
On the other hand, for all $0\leq x\leq \pi$,
\[0\leq x(\pi -x)\leq \left(\frac{\pi}{2}\right)^2.\]
Therefore, for $0\leq x\leq \pi$, 
\[0\leq p_n(x)\leq \frac{1}{n!}\left(\frac{\pi^2 b}{4}\right)^n.\]
Since we also have $0\leq\sin x\leq 1$ for all $0\leq x\leq \pi$, we conclude that
\[0\leq I_n=\int_0^{\pi}p_n(x)\sin x dx\leq  \frac{\pi}{n!}\left(\frac{\pi^2 b}{4}\right)^n.\] 
Because the series $\di\sum_{n=0}^{\infty} \frac{1}{n!}\left(\frac{\pi^2 b}{4}\right)^n$ is convergent, we find that 
\[\lim_{n\to\infty}  \frac{1}{n!}\left(\frac{\pi^2 b}{4}\right)^n=0.\]Therefore, there is a positive integer $N$ such that for all $n\geq N$,
\[  \frac{1}{n!}\left(\frac{\pi^2 b}{4}\right)^n<\frac{1}{\pi},\]
which gives
\[0\leq I_n<1\hspace{1cm}\text{for all}\;n\geq N.\]\bp
We only need the $n=N$ case now. Since $I_N$ is an integer, we must have $I_N=0$. However, since 
\[p_{N}(x)\sin  x\] is a  continuous function and it is positive on $(0,\pi)$, by Example \ref{ex230222_5}, \[I_N=\di\int_0^{\pi}p_N(x)\sin xdx\] cannot be zero. This gives a contradiction. Hence, $\pi$ must be an irrational number. 

\end{myproof}

\bigskip
\subsection{ Infinitely Differentiable Functions that are Non-Analytic}\label{sec6.6.3}
We consider the function $f:\mathbb{R}\to\mathbb{R}$ defined by
\begin{equation*} f(x)=\begin{cases}\di\exp\left(-\frac{1}{x^2}\right),\quad &\text{if}\;x\neq 0,\\0,\quad &\text{if}\;x=0.\end{cases}\end{equation*}
We will show that this function  is infinitely differentiable and $f^{(n)}(0)=0$ for all $n\geq 0$. 
 
\begin{figure}[ht]
\centering
\includegraphics[scale=0.2]{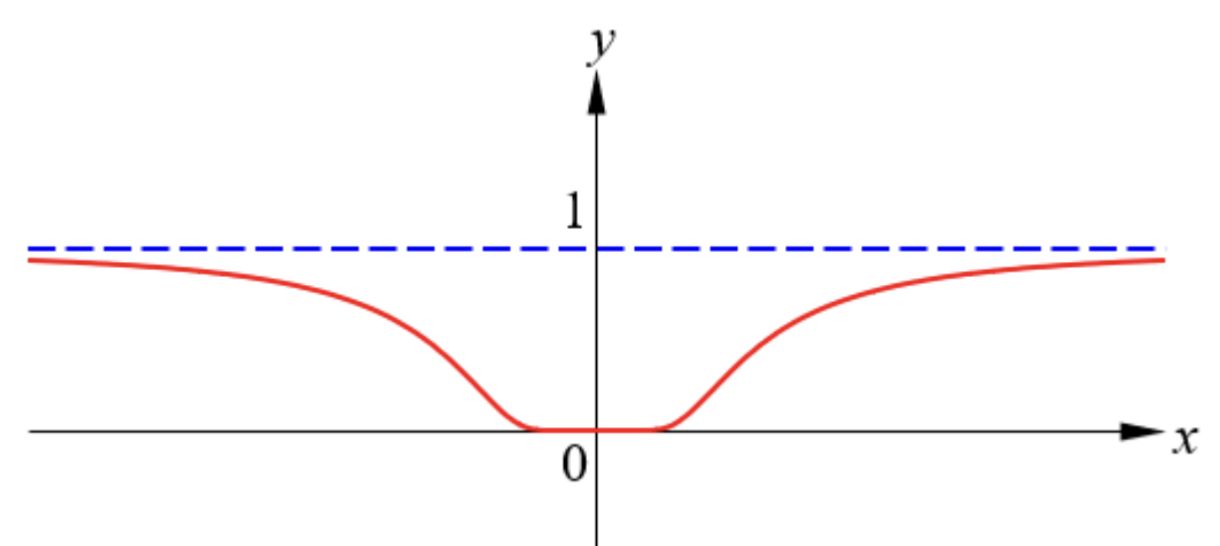}
\caption{The function $\di f(x)=\exp\left(-\frac{1}{x^2}\right)$.\fa}\label{figure58}
\end{figure}

Let us first prove the following lemma.
\begin{lemma}[label=230309_6]{}
If $p(x)$ is a polynomial, then
\begin{subequations}\label{eq230309_4}\begin{align}\lim_{x\to 0^+}p\left(\frac{1}{x}\right)\exp\left(-\frac{1}{x}\right)=0 \label{eq230309_4_1}\\
\lim_{x\to 0}p\left(\frac{1}{x}\right)\exp\left(-\frac{1}{x^2}\right)=0.\label{eq230309_4_2}
\end{align}\end{subequations}
\end{lemma}\begin{myproof}{Proof}
In Example \ref{230307_11}, we have shown that for any  real number $s$,
$\di \lim_{y\to \infty}y^s e^{-y}=0$. From this, we find that
 if $k$ is an integer,
\[\lim_{x\to 0^+}\frac{1}{x^k}\exp\left(-\frac{1}{x}\right)=\lim_{y\to \infty}y^ke^{-y}=0,\] and 
\[\lim_{x\to 0}\frac{1}{x^{2k}}\exp\left(-\frac{1}{x^2}\right)=\lim_{y\to \infty}y^ke^{-y}=0.\]
The latter one implies that for any integer $k$,
\[\lim_{x\to 0}\frac{1}{x^{2k-1}}\exp\left(-\frac{1}{x^2}\right)=\left[\lim_{x\to 0}x\right]\left[\lim_{x\to 0}\frac{1}{x^{2k}}\exp\left(-\frac{1}{x^2}\right)\right]=0.\]
These prove \eqref{eq230309_4}.
\end{myproof}

Next, we prove the following.
\begin{theorem}[label=230309_7]{}
Let  $f:\mathbb{R}\to\mathbb{R}$ be the function defined by
\begin{equation}\label{eq230307_12}f(x)=\begin{cases}\di\exp\left(-\frac{1}{x^2}\right),\quad &\text{if}\;x\neq 0,\\0,\quad &\text{if}\;x=0.\end{cases}\end{equation}
Then $f$ is an infinitely differentiable function with $f^{(n)}(0)=0$ for all $n\geq 0$.
\end{theorem}

\begin{myproof}{Proof} We claim that for each positive integer $n$, there is a polynomial $p_n(x)$ of degree $3n$ such that
 \begin{equation}
\label{eq230309_5}f^{(n)}(x)=\begin{cases}\di p_n\left(\frac{1}{x}\right)\exp\left(-\frac{1}{x^2}\right),\quad &\text{if}\;x\neq 0,\\0,\quad &\text{if}\;x=0.\end{cases}\end{equation}
This will show that $f$ is infinitely differentiable. In fact, Lemma \ref{230309_6} implies that 
\[\lim_{x\to 0}f_n(x)=\lim_{n\to 0}p_n\left(\frac{1}{x}\right)\exp\left(-\frac{1}{x^2}\right)=0=f^{(n)}(0),\]which says that $f^{(n)}(x)$ is continuous at $x=0$.
We will prove  \eqref{eq230309_5}  by induction on $n$.
When $n=1$, we find from the definition \eqref{eq230307_12} that
\[f'(x)=\frac{2}{x^3}\exp\left(-\frac{1}{x^2}\right)\hspace{1cm}\text{when}\;x\neq 0.\] 

When $x=0$, we apply Lemma \ref{230309_6} to get
\[f'(0)=\lim_{x\to 0}\frac{\di \exp\left(-\frac{1}{x^2}\right)-0}{x}=\lim_{x\to 0}\frac{1}{x}\exp\left(-\frac{1}{x^2}\right)=0.\]Therefore, the $n=1$ statement is true with $p_1(x)$ is polynomial of degree 3 given by
\[p_1(x)=2x^3.\]
Assume that the statement is true for the $n-1$ case. This means that there is a polynomial $p_{n-1}(x)$ of degree   $3n -3$ such that
\[f^{(n-1)}(x)=\begin{cases}\di  p_{n-1}\left(\frac{1}{x }\right)\exp\left(-\frac{1}{x^2}\right),\quad &\text{if}\;x\neq 0,\\0,\quad &\text{if}\;x=0.\end{cases}\]When $x\neq 0$,
\[f^{(n)}(x)=\left( -\frac{1}{x^2}p_{n-1}'\left(\frac{1}{x}\right)+\frac{2}{x^3}p_{n-1}\left(\frac{1}{x}\right)\right)\exp\left(-\frac{1}{x^2}\right).\]\bp This shows that when $x\neq 0$,
\[  f^{(n)}(x)=p_{n}\left(\frac{1}{x }\right)\exp\left(-\frac{1}{x^2}\right),\]
where
\[p_n(x)=-x^2 p_{n-1}'(x)+2x^3p_{n-1}(x).\]
By inductive hypothesis, $p_{n-1}'(x)$ is a polynomial of degree $3n-4$. Thus, $x^2 p_{n-1}'(x)$ is a polynomial of degree $3n-2$. Since $2x^3p_{n-1}(x)$ is a polynomial of degree $3n$, $p_n(x)$ is a polynomial of degree $3n$. 
For the derivative at 0, Lemma \ref{230309_6} implies that
\[f^{(n)}(0)=\lim_{x\to 0}\frac{f^{(n-1)}(x)-f^{(n-1)}(0)}{x}=\lim_{x\to 0}\frac{1}{x}p_{n-1}\left(\frac{1}{x}\right)\exp\left(-\frac{1}{x^2}\right)=0.\]
This proves the statement for the $n$ case, and thus completes the induction.\end{myproof}

Now we prove our main theorem in this section.
\begin{theorem}{}Let $I$ be an open interval that contains the point $x_0$. There is an infinitely differentiable function $f:I\to\mathbb{R}$ whose Taylor series at the point $x=x_0$  is convergent pointwise on $I$, but it does not converge to $f(x)$ pointwise on $I$.
\end{theorem}
\begin{myproof}{Proof}Define the function $f:I\to\mathbb{R}$ by
\begin{equation*} \begin{cases}\di\exp\left(-\frac{1}{(x-x_0)^2}\right),\quad &\text{if}\;x\neq x_0,\\0,\quad &\text{if}\;x=x_0.\end{cases}\end{equation*}Theorem \ref{230309_7} implies  that
 the function $f(x)$  is infinitely differentiable on $I$, and 
$f^{(n)}(x_0)=0$ for all $n\geq 0$. Hence, the Taylor series of $f$ at $x=x_0$,\bp
\[\sum_{n=0}^{\infty}\frac{f^{(n)}(x_0)}{n!}(x-x_0)^n,\] is the series that is identically 0. Therefore, it converges everywhere, but it does not converge to $f(x)$ except at the point $x=x_0$.  
\end{myproof}

 Using almost the same proof as for Theorem \ref{230309_7}, we obtain the following.
\begin{theorem}[label=230309_11]{}
Given a real number $x_0$, the function $g:\mathbb{R}\to\mathbb{R}$ defined by
\begin{equation}\label{eq230309_8}g(x)=\begin{cases} \di\exp\left(-\frac{1}{x-x_0}\right),\quad &\text{if}\;x> x_0,\\0,\quad &\text{if}\;x\leq x_0.\end{cases}\end{equation} is infinitely differentiable.
\end{theorem}

The function $g(x)$ defined by \eqref{eq230309_8} is also not analytic. Nevertheless, it has some important applications. It is usually used to "smooth" up a function or truncate a function smoothly.

\begin{theorem}{}
Given two real numbers $a$ and $b$ with $a<b$, define the function $h:\mathbb{R}\to\mathbb{R}$   by
\begin{equation}\label{eq230309_12}\begin{split}
h(x)&=\frac{g(x-a)}{g(x-a)+g(b-x)}\\&=\begin{cases}
0,\quad &\text{if}\; \quad x\leq a,\\ \di \frac{\di\exp\left(-\frac{1}{x-a}\right)}{\di\exp\left(-\frac{1}{x-a}\right)+\exp\left(-\frac{1}{b-x}\right)},\quad &\text{if}\;a<x<b,\\1,\quad &\text{if}\; \quad x\geq b.\end{cases}\end{split}\end{equation}Then $h$ is a functon that is infinitely differentiable.
\end{theorem}
 
\begin{myproof}{Proof}
We just need to show that $g(x-a)+g(b-x)$ is nonzero for all $x\in\mathbb{R}$. The rest follows from Theorem \ref{230309_11} and the definition of $h(x)$. Since the function $g$ is nonnegative, in order for $g(x-a)+g(b-x)=0$, we must have $g(x-a)=g(b-x)=0$. But we know that $g(x-a)=0$ only  when $x\leq a$, and $g(b-x)=0$ only when $x\geq b$. Since the set $\{x\,|\,x\leq a\}$ and the set $\{x\,|\,x\geq b\}$ are disjoint, we conclude that $g(x-a)+g(b-x)$ is never 0.
\end{myproof}
 
\begin{figure}[ht]
\centering
\includegraphics[scale=0.2]{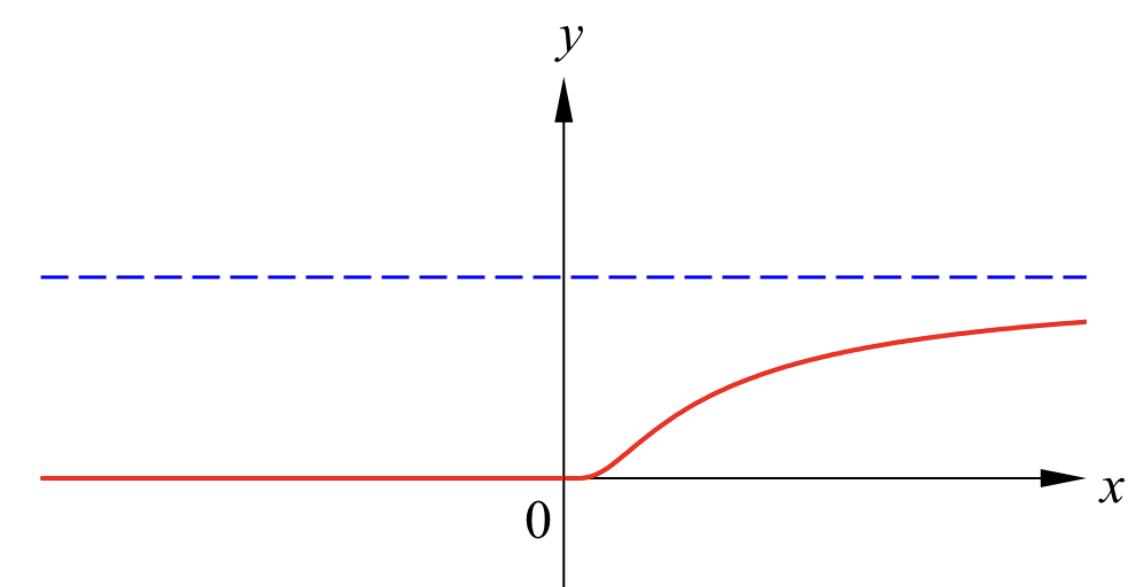}
\caption{The function $\di g(x)$ defined by \eqref{eq230309_8} when $x_0=0$.\fa}\label{figure61}
\end{figure}

\begin{figure}[ht]
\centering
\includegraphics[scale=0.2]{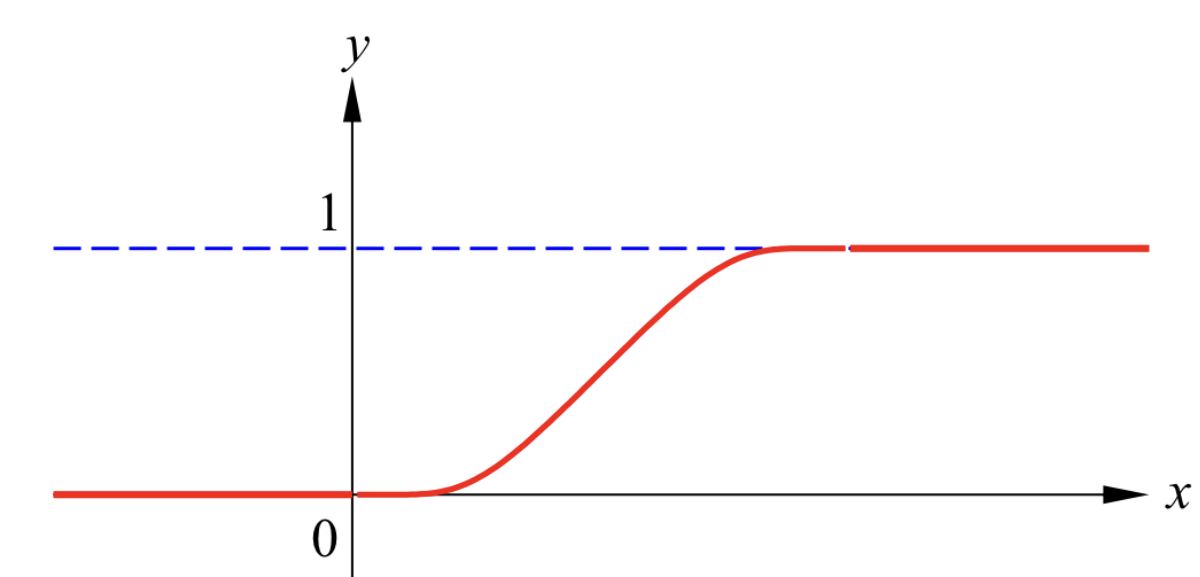}
\caption{The function $\di h(x)$ defined by \eqref{eq230309_12}.\fa}\label{figure62}
\end{figure}

\begin{remark}{}
The function $h(x)$ defined  by \eqref{eq230309_12}  is an example of an infinitely diferentiable function that is increasing but assume constant values outside a bounded interval.
\end{remark}

\bigskip
\subsection{A Continuous Function that is Nowhere Differentiable}\label{sec6.6.4}
In this section, we want to construct a continuous function $f:\mathbb{R}\to\mathbb{R}$ which is not differentiable at any point. The main ingredient in the proof is to note that the function $g:\mathbb{R}\to\mathbb{R}$, $g(x)=|x|$ is continuous, and it is not differentiable at $x=0$. 

\begin{definition}[label=230309_10]{The function $\pmb{h_m}$}
For any positive number $m$, let $h_m:\mathbb{R}\to\mathbb{R}$ be the function defined by
\[h_m(x)=|x|,\hspace{1cm}\text{for all}\;-m\leq x\leq m,\]
and
\[h_m(x+2m)=h(x)\hspace{1cm}\text{for all}\;x\in\mathbb{R}.\]

\end{definition}

\begin{figure}[ht]
\centering
\includegraphics[scale=0.2]{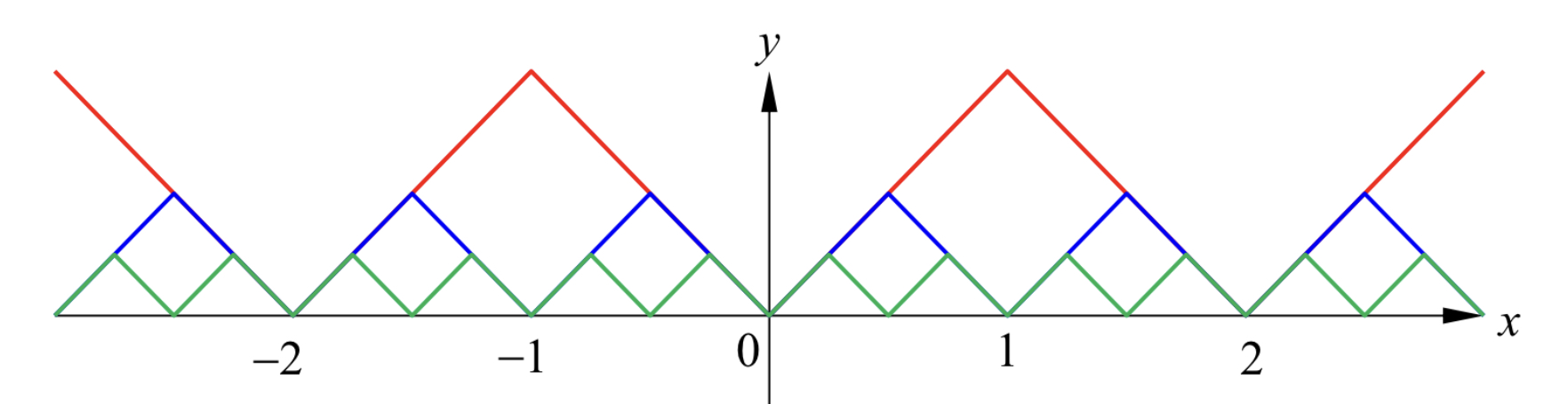}
\caption{The functions $\di h_m(x)$ when $m=1$, $\frac{1}{2}$ and $\frac{1}{4}$.\fa}\label{figure63}
\end{figure}

 Let us first explore the properties of the function $h_m$.
\begin{lemma}[label=230309_9]{}
Given  a positive number $m$, define $x_{m,k}=mk$ for all $k\in\mathbb{Z}^+$. The function $h_m$ defined in Definition \ref{230309_10} has the following properties.
\begin{enumerate}[(a)]
\item $h_m$  is a continuous even function that is periodic of period $2m$.
\item For $k\in \mathbb{Z}$,  the graph of $h_m:[x_{m,2k}, x_{m, 2k+1}]\to\mathbb{R}$ is a straightline segment of slope 1; while the graph of $h_m:[x_{m,2k-1}, x_{m, 2k}]\to\mathbb{R}$ is a straightline segment of slope $-1$. Hence, the graph  of $h_m$ is a union of straightline segments alternatingly having slopes 1 and $-1$.
\item $0\leq h_m(x)\leq m$ for all $x\in\mathbb{R}$.
\end{enumerate}
\end{lemma}
\begin{myproof}
{Proof}
Part (a) follows from $h_m(-m)=h_m(m)$. Part (b) can be proved by  induction on $k\geq 0$, using the periodicity of $h_m$ and the fact that $h_m$ is an even function.  Part (c) follows from the definition of $h_m$ and periodicity.
\end{myproof}

\begin{lemma}[label=230309_13]{}
Given  a positive number $\ell$ and a point  $x\in \mathbb{R}$, let \[U=[x-\ell /2, x]\hspace{1cm} \text{and} \hspace{1cm} V=[x,x+\ell/2].\] For a positive integer $m$, let $h_m:\mathbb{R}\to\mathbb{R}$ be the function defined in Definition \ref{230309_10}. Then  one of the following holds.
 
\begin{enumerate}[(a)]
\item For each nonnegative integer $k$, the graph of $h_{2^k\ell}:U\to\mathbb{R}$ is a line segment of slope 1 or $-1$.

\item  For each nonnegative integer $k$, the graph of $h_{2^k\ell}:V\to\mathbb{R}$ is a line segment of slope 1 or $-1$.

\end{enumerate}
\end{lemma}
\begin{myproof}{Proof} 
The points $n\ell$, $n\in\mathbb{Z}$, partition the real line into subintervals of the form $[n\ell, (n+1)\ell]$, each of length $\ell$. Since $U$ and $V$ are adjacent intervals of length $\ell/2$, one of them must lie entirely inside one of the intervals of the form $[n\ell, (n+1)\ell]$.

 By part (b) in Lemma \ref{230309_9}, the graph of $h_{\ell}:[n\ell, (n+1)\ell]\to\mathbb{R}$ is a line segment of slope 1 or $-1$. This proves the assertion when $k=0$. To prove the assertion for $k\geq 1$, we notice that to obtain the graph of the function $h_{\ell}$ from the graph of the function $h_{2\ell}$, we divide each line segment in the graph of $h_{2\ell}$ into two equal parts, one of the parts change slope from 1 to $-1$ or from $-1$ to $1$. Hence, if $W$ is an interval and the graph of $h_{\ell}:W\to\mathbb{R}$ is a line segment, the graph of $h_{2^k\ell}:W\to\mathbb{R}$ must also be a line segment for any $k\in\mathbb{Z}^+$. This completes the proof the the lemma.
\end{myproof}

Now we can prove the main theorem in this section.

\begin{theorem}
{}For a positive number $m$, let $h_m:\mathbb{R}\to\mathbb{R}$ be the function
\[h_m(x)=|x|\;\;\text{for }\;|x|\leq m, \quad \text{and}\quad h(x+2m)=h(x)\;\;\text{for all}\;x\in\mathbb{R}.\]For $n\geq 0$, let $g_n:\mathbb{R}\to\mathbb{R}$ be the function defined by
$g_n(x)=h_{m_n}(x)$, with $m_n=\di\frac{1}{4^n}$. 
Then the series $\di \sum_{n=0}^{\infty}g_n(x)$ converges uniformly to a function  $f:\mathbb{R}\to\mathbb{R}$,  \[ f(x)=\sum_{n=0}^{\infty}g_n(x).\] $f(x)$ is  a continuous function that is not differentiable at any point.
\end{theorem}
\begin{myproof}{Proof}
 By Lemma \ref{230309_9}, for each $n\in\mathbb{Z}^+$, the function $g_n:\mathbb{R}\to\mathbb{R}$ is continuous and 
\[|g_n(x)|\leq\frac{1}{4^n}\hspace{1cm}\text{for all}\;x\in\mathbb{R}.\]   Since the series $\di\sum_{n=0}^{\infty}\frac{1}{4^n}$ is convergent,
  Weierstrass $M$-test implies that the series $\di \sum_{n=0}^{\infty}g_n(x)$ converges uniformly on $\mathbb{R}$.  Since each $g_n(x)$ is a continuous function,   Corollary \ref{230305_10} implies that  the function  
$\di f(x)=\sum_{n=0}^{\infty}g_n(x)$ is continuous.

 Now we are left to prove that $f(x)$ is not diferentiable at any $x\in \mathbb{R}$. Given $x_0\in \mathbb{R}$, assume that 
\[f'(x_0)=\lim_{x\to x_0}\frac{f(x)-f(x_0)}{x-x_0}\] exists. Then for any sequence $\{x_k\}_{k=0}^{\infty}$ in $\mathbb{R}\setminus\{x_0\}$,  if $\di \lim_{k\to\infty}x_k=x_0$, then the limit  
\[\lim_{k\to \infty}\frac{f(x_k)-f(x_0)}{x-x_0}\]exists and is equal to $f'(x_0)$.\bp

We construt a sequence $\{x_{k}\}_{k=0}^{\infty}$ as follows.  For each $k\in \mathbb{Z}^+$, Lemma \ref{230309_13} implies that either the graph of $
h_{m_k}:[x_0-m_k/2, x_0]\to\mathbb{R}$ or the graph of $h:[x_0, x_0+m_k/2]\to\mathbb{R}$ is a line segment with slope 1 or $-1$.   In the former case, we let $x_{k}=x_0-m_k/2$. In the latter case, we let $x_{k}=x_0+m_k/2$. In any case, we find that
\[|x_{k}-x_0|=\frac{m_k}{2}=\frac{1}{2^{2k+1}}\hspace{1cm}\text{for all}\;k\geq 0.\]
This shows that $\{x_{k}\}$ is a sequence in $\mathbb{R}\setminus\{x_0\}$ that converges to $x_0$. 

For fixed $k\in \mathbb{Z}^+$, 
$m_k/2$ is a multiple of $2m_n$ for all $n>k+1$. By periodicity of $h_{m_n}$, 
\[g_n(x_{k})-g_n(x_0)=h_{m_n}\left(x_0\pm \frac{m_k}{2}\right)-h_{m_n}(x_0)=0\hspace{1cm} \text{for all}\;n>k.\]
This implies that
\[ \frac{f(x_k)-f(x_0)}{x_k-x_0}=\sum_{n=0}^{\infty}\frac{g_n(x_k)-g_n(x_0)}{x_k-x_0}=\sum_{n=0}^{k}\frac{g_n(x_k)-g_n(x_0)}{x_k-x_0}.\]
By the  definition of $x_k$ and Lemma \ref{230309_13}, 
\[\frac{g_n(x_k)-g_n(x_0)}{x_k-x_0}\] is equal to 1 or $-1$ for each $0\leq k\leq n$. The sum of an odd number of 1 or $-1$ must be odd. The sum of an even number of 1 or $-1$ must be even.
Therefore,
\[c_k=\frac{f(x_k)-f(x_0)}{x_k-x_0}\] is odd when $k$ is even, and is even when $k$ is odd. 
This implies that the sequence $\{c_k\}_{k=0}^{\infty}$ is an integer sequence that is alternatingly odd and even. Hence, it does not have a limit. This is a contradiction, which allows us to  conclude that $f$ cannot be differentiable at $x_0$.

\end{myproof}

\begin{figure}[ht]
\centering
\includegraphics[scale=0.18]{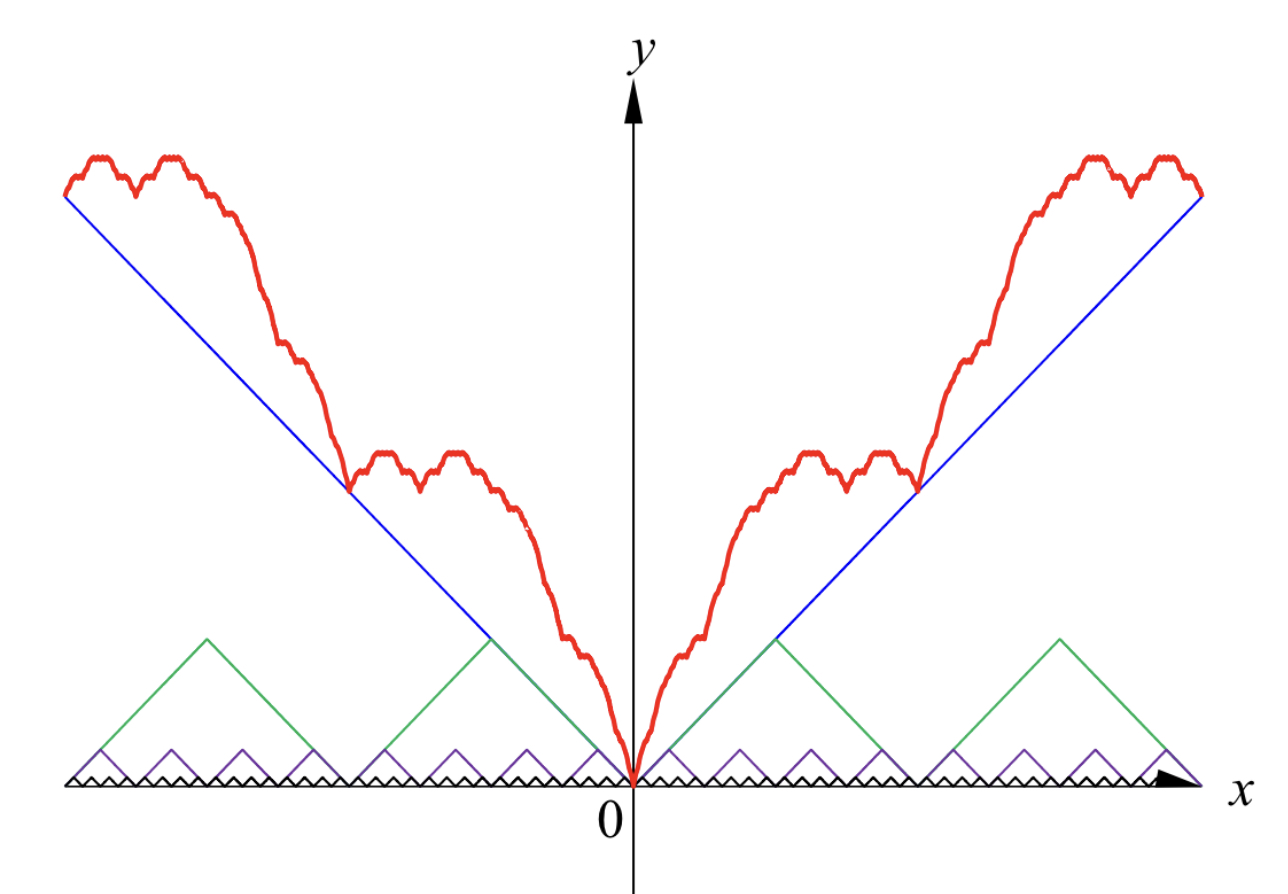}
\caption{The functions $g_n(x)$ for $n=0, 1, 2, 3$, and the function $f(x)$.\fa}\label{figure64}
\end{figure}
 
 \bigskip
\subsection{The Weierstrass Approximation Theorem}\label{sec6.6.5}
In this section, we prove the Weierstrass approximation theorem using Bernstein's ingenious approach. We start with a  lemma.

\begin{lemma}[label=230309_14]{}
The following identities hold.
\begin{enumerate}[(a)]
\item For $n\geq 0$, $\di\sum_{k=0}^n\binom{n}{k}x^k(1-x)^{n-k}=1$.
\item For $n\geq 1$, $\di\sum_{k=1}^n\frac{k}{n}\binom{n}{k}x^k(1-x)^{n-k}=x$.

%\item For $n\geq 2$, $\di\sum_{k=2}^n\frac{k(k-1)}{n(n-1)}\binom{n}{k}x^k(1-x)^{n-k}=x^2$.
\item For $n\geq 2$, $\di\sum_{k=1}^n\frac{k^2}{n^2}\binom{n}{k}x^k(1-x)^{n-k}=x^2+\frac{x(1-x)}{n}$.
\item For $n\geq 2$, $\di\sum_{k=0}^n\left(x-\frac{k}{n}\right)^2\binom{n}{k}x^k(1-x)^{n-k}=\frac{x(1-x)}{n}$.
 \end{enumerate}
\end{lemma}
\begin{myproof}{Proof}
The first  identity  (a)
is just a consequence of the binomial expansion theorem.
\bp
 For the   identity in (b), notice that when $n\geq k\geq 1$,
\[\frac{k}{n}\binom{n}{k}=\frac{(n-1)!}{(k-1)!(n-k)!}=\binom{n-1}{k-1}.\]
Therefore, 
\begin{align*}\sum_{k=1}^n\frac{k}{n}\binom{n}{k}x^k(1-x)^{n-k}&=x\sum_{k=1}^n\binom{n-1}{k-1}x^{k-1}(1-x)^{n-k}\\&=x\sum_{k=0}^{n-1}\binom{n-1}{k}x^k(1-x)^{n-1-k}=x.\end{align*}
 For part (c), we find that when $n\geq k\geq 2$,
\[\frac{k(k-1)}{n(n-1)}\binom{n}{k}=\frac{(n-2)!}{(k-2)!(n-k)!}=\binom{n-2}{k-2}.\]  It follows that
\begin{align*}\sum_{k=2}^n\frac{k(k-1)}{n(n-1)}\binom{n}{k}x^k(1-x)^{n-k} =x^2\sum_{k=0}^{n-2}\binom{n-2}{k}x^k(1-x)^{n-2-k}=x^2.\end{align*}
Writing  $k^2=k(k-1)+k$, we have
\begin{align*}
&\sum_{k=1}^n\frac{k^2}{n^2}\binom{n}{k}x^k(1-x)^{n-k}\\&=\sum_{k=1}^n\frac{k(k-1)}{n^2}\binom{n}{k}x^k(1-x)^{n-k}+\sum_{k=1}^n\frac{k }{n^2}\binom{n}{k}x^k(1-x)^{n-k}\\
&=\frac{n-1}{n}x^2+\frac{1}{n}x=x^2+\frac{x(1-x)}{n}.
\end{align*} 

For the identity in part (d),
  a straightforward computation  gives
\begin{align*}
&\sum_{k=0}^n\left(x-\frac{k}{n}\right)^2\binom{n}{k}x^k(1-x)^{n-k}\\
&= \sum_{k=0}^n\left(x^2- \frac{2k}{n}x+\frac{k^2}{n^2}\right) \binom{n}{k}x^k(1-x)^{n-k}\\
&=x^2-2x^2+x^2+\frac{x(1-x)}{n} =\frac{x(1-x)}{n}.
\end{align*}
\end{myproof}

\begin{definition}{Bernstein Basis Polynomials}
For any positive integer $n$, there are $n+1$ Bernstein basis polynomials given by
\[p_{n,k}(x)=\binom{n}{k}x^k(1-x)^{n-k},\hspace{1cm}0\leq k\leq n.\]
\end{definition}

 \begin{figure}[ht]
\centering
\includegraphics[scale=0.2]{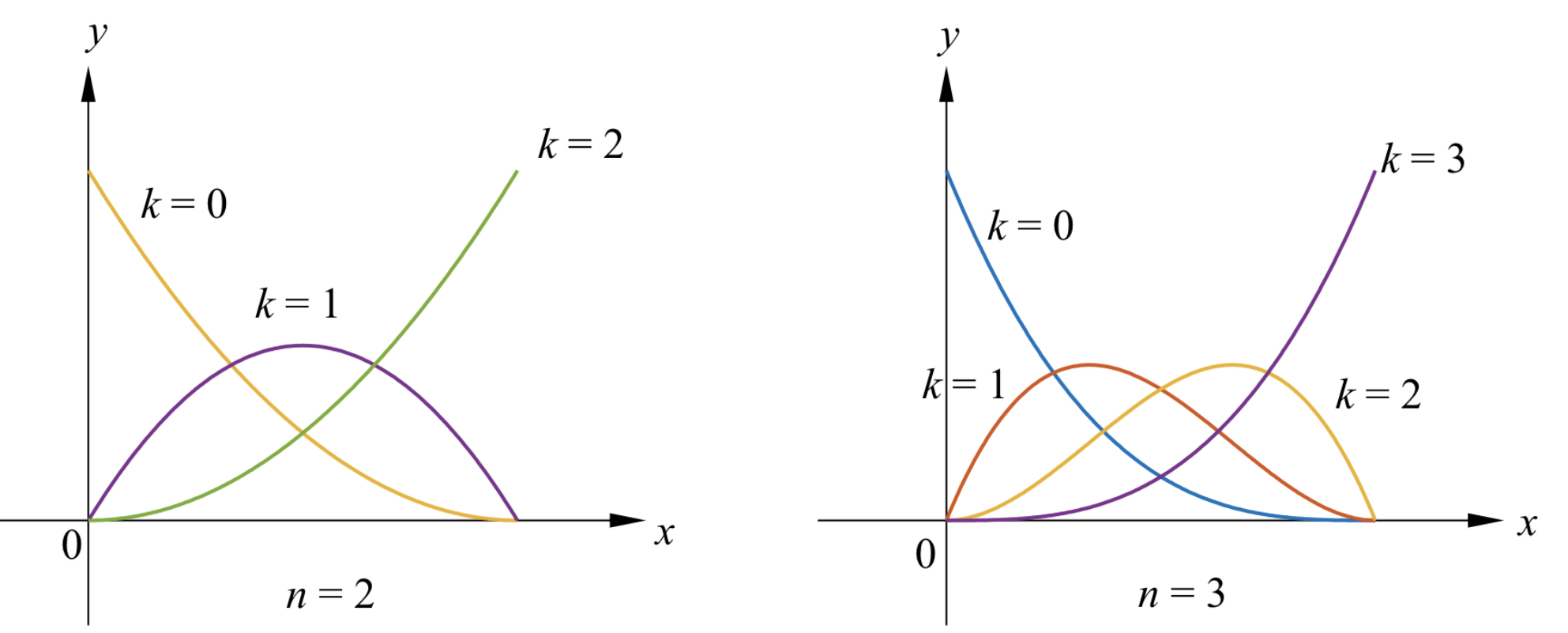}
\caption{The polynomials $\di p_{n,k}(x)=\binom{n}{k}x^k(1-x)^{n-k}$ when $n=2$ and $n=3$, for all $0\leq k\leq n$.\fa}\label{figure65}
 
\includegraphics[scale=0.2]{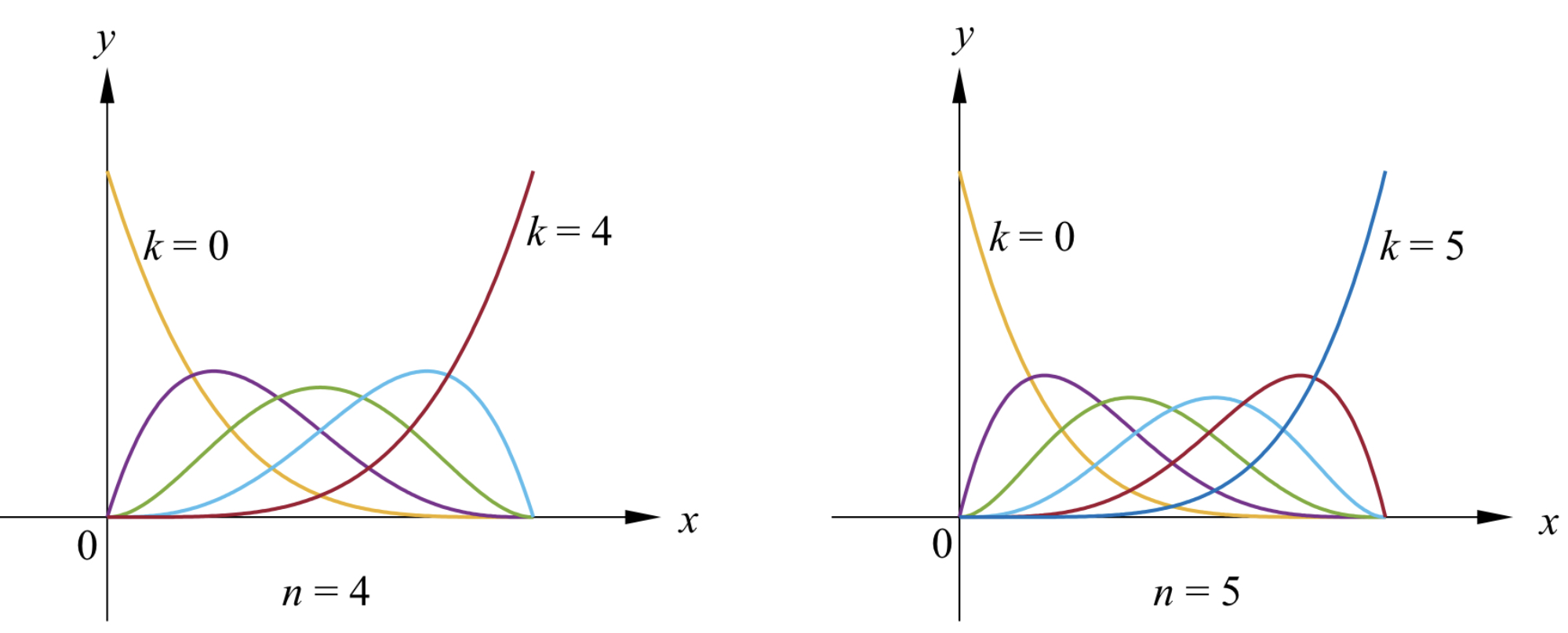}
\caption{The polynomials $\di  p_{n,k}(x)=\binom{n}{k}x^k(1-x)^{n-k}$ when $n=4$ and $n=5$, for all $0\leq k\leq n$.\fa}\label{figure66}
\end{figure}

\begin{figure}[ht]
\centering
\includegraphics[scale=0.2]{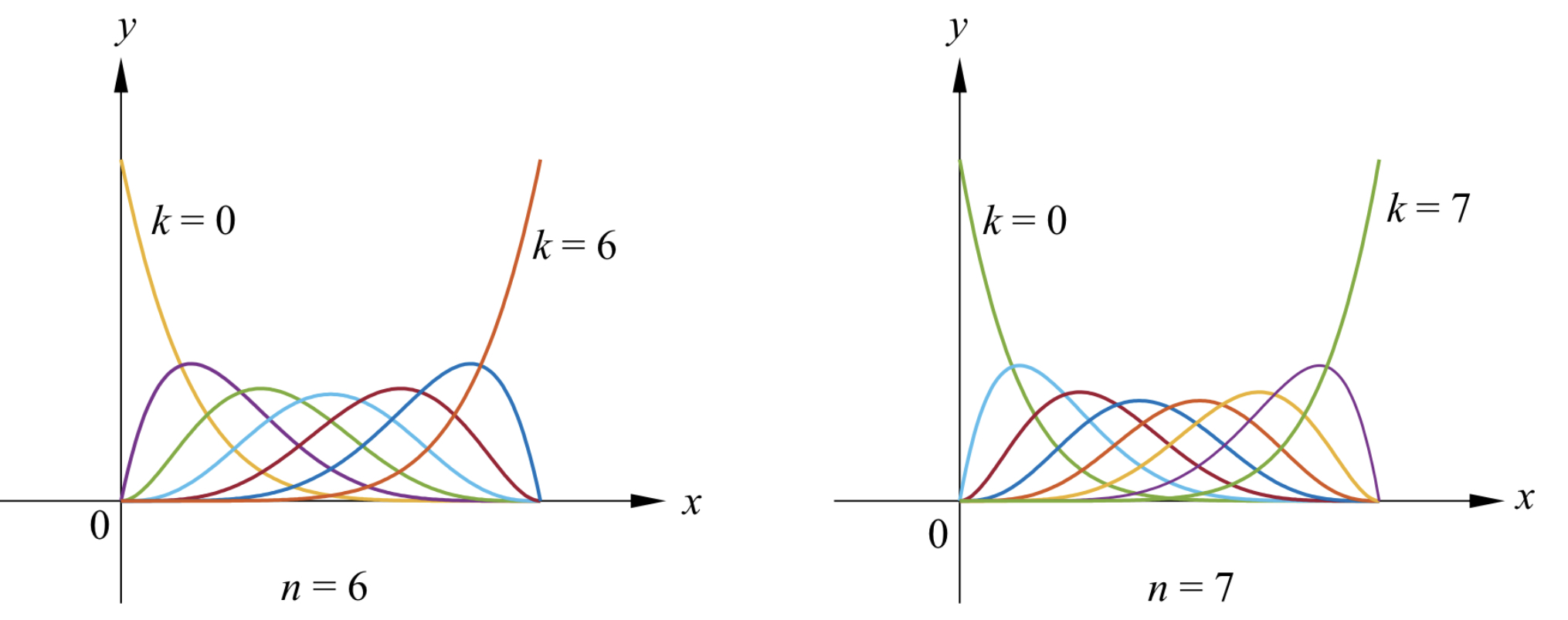}
\caption{The polynomials $\di  p_{n,k}(x)=\binom{n}{k}x^k(1-x)^{n-k}$ when $n=6$ and $n=7$, for all $0\leq k\leq n$.\fa}\label{figure67}
\end{figure}

Now we come to our main theorem.
\begin{theorem}{Weierstrass Approximation Theorem}
Let $f:[a,b]\to\mathbb{R}$ be a continuous function defined on $[a, b]$. Given $\varepsilon>0$, there is a polynomial $p(x)$ such that
\[|f(x)-p(x)|<\varepsilon\hspace{1cm}\text{for all}\;x\in [a,b].\]
\end{theorem}
\begin{myproof}{Proof}
We first consider the case where $[a,b]=[0,1]$. Since $f:[0,1]\to\mathbb{R}$ is continuous on a closed and bounded interval, it is uniformly continuous and bounded. The boundeness of $f$ implies that there is a positive number $M$ such that
\[|f(x)|\leq M\hspace{1cm}\text{for all}\;x\in [0,1].\] Given $\varepsilon>0$, since $f$ is uniformly continuous, there is a $\delta>0$ such that for all $x_1$ and $x_2$ in $[0,1]$, if $|x_1-x_2|<\delta$, then 
\[|f(x_1)-f(x_2)|<\frac{\varepsilon}{2}.\]  

For any positive integer  $n$, we construct a polynomial $p_n(x)$ to be a polynomial of degree at most $n$ given by the following linear combination of Bernstein basis polynomials.\bp

\[p_n(x)=\sum_{k=0}^n f\left(\frac{k}{n}\right)p_{n,k}(x)=\sum_{k=0}^n f\left(\frac{k}{n}\right)\binom{n}{k}x^k(1-x)^{n-k}.\]  Let us estimate the supremum of $|f(x)-p_n(x)|$ on $[0,1]$. For fixed $x\in [0,1]$, part (a) in Lemma \ref{230309_14} implies that
\[
f(x)-p_n(x)=\sum_{k=0}\left(f(x)-f\left(\frac{k}{n}\right)\right)\binom{n}{k}x^k(1-x)^{n-k}.\] 
 
Since $x^k(1-x)^{n-k}\geq 0$ for all $x\in [0, 1]$ and all $n\geq k\geq 0$, triangle inequality gives
\[
\left|f(x)-p_n(x)\right|\leq \sum_{k=0}\left|f(x)-f\left(\frac{k}{n}\right)\right|\binom{n}{k}x^k(1-x)^{n-k}.\]  
For  $0\leq k\leq n$, if  $\di \left|x-\frac{k}{n}\right|<\delta$, then 
\[\left|f(x)-f\left(\frac{k}{n}\right)\right|<\frac{\varepsilon}{2}.\]
If $\di \left|x-\frac{k}{n}\right|\geq \delta$, then
\[\left|f(x)-f\left(\frac{k}{n}\right)\right|\leq \left|f(x)\right|+\left|f\left(\frac{k}{n}\right)\right|\leq 2M\leq \frac{2M}{\delta^2}\left(x-\frac{k}{n}\right)^2.\]In any case, we find that
\[\left|f(x)-f\left(\frac{k}{n}\right)\right|<\frac{\varepsilon}{2}+\frac{2M}{\delta^2}\left(x-\frac{k}{n}\right)^2\hspace{1cm}\text{for all}\;0\leq k\leq n.\]
Therefore,
\[
\left|f(x)-p_n(x)\right|  < \sum_{k=0}^n\left[\frac{\varepsilon}{2}+\frac{2M}{\delta^2}\left(x-\frac{k}{n}\right)^2\right]\binom{n}{k}x^k(1-x)^{n-k}.
\]By Lemma \ref{230309_14}, and the fact that
\[0\leq x(1-x)\leq\frac{1}{4}\hspace{1cm}\text{for all}\;0\leq x\leq 1,\]
\bp
we find that
\[\left|f(x)-p_n(x)\right|  <\frac{\varepsilon}{2}+\frac{2M}{\delta^2}\frac{x(1-x)}{n}\leq \frac{\varepsilon}{2}+\frac{M}{2\delta^2n}.\]
If 
$\di n\geq \frac{M}{\varepsilon\delta^2}$, 
then
$\di \frac{M}{2\delta^2n}\leq\frac{\varepsilon}{2}$.
For any such $n$, we find that
\[\left|f(x)-p_n(x)\right|  <\varepsilon\hspace{1cm}\text{for all}\;0\leq x\leq 1.\]
This completes the proof when $[a,b]=[0,1]$. 

  For general $[a,b]$, let $u:[0,1]\to\mathbb{R}$ be the polynomial function $u(t)=a+t(b-a)$. This is a continuous function mapping $[0,1]$ bijectively onto $[a,b]$. The inverse is the continuous function $u^{-1}(x)=\di\frac{x-a}{b-a}$. The function $g=f\circ u:[0,1]\to\mathbb{R}$, being a composition of continuous functions, is continuous.  By what we have proved above, given $\varepsilon>0$, there is a polynomial $q(t)$ so that 
\[|f(u(t))-q(t)|<\varepsilon\hspace{1cm}\text{for all}\;t\in [0,1].\]
Let 
\[p(x)=q(u^{-1}(x))=q\left(\frac{x-a}{b-a}\right).\] Then $p(x)$ is also a polynomial, and $p(u(t))=q(t)$. Therefore,
\[|f(u(t))-p(u(t))|<\varepsilon\hspace{1cm}\text{for all}\;t\in [0,1],\]which implies that
\[|f(x)-p(x)|<\varepsilon \hspace{1cm}\text{for all}\;x\in [a,b].\]This completes the proof of the Weierstrass approximation theorem for the general  case.

\end{myproof}
\begin{figure}[ht]
\centering
\includegraphics[scale=0.2]{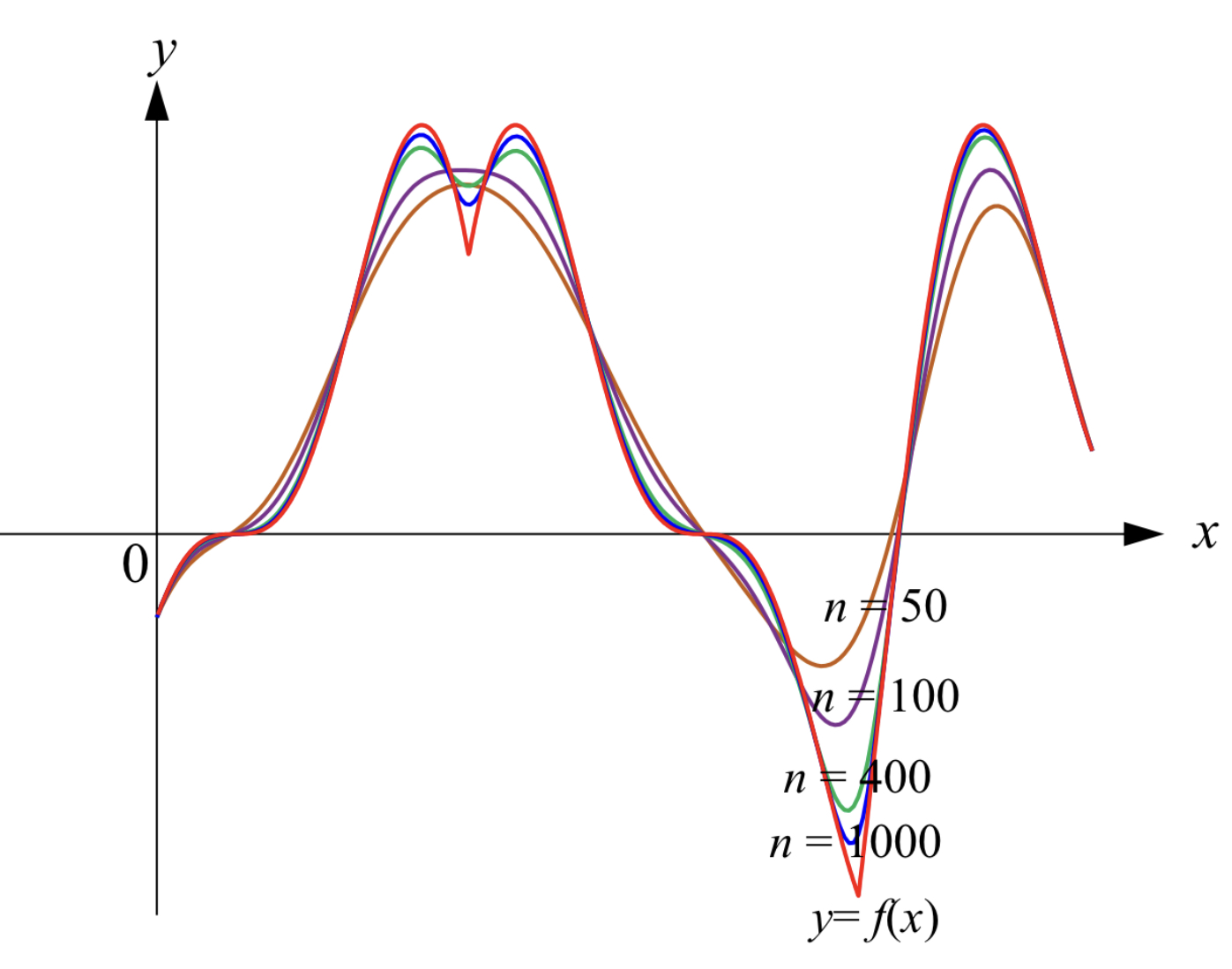}
\caption{Approximations of the continuous function $f(x)$ by the polynomials $p_n(x)$, where $f(x)=\di 3\sin( 4\pi|x-1/3|)+2\sin (6\pi|x-3/4|)$.\fa}\label{figure68}
\end{figure}
~
 
One cannot extend the Weierstrass approximation theorem to the case where $f:I\to\mathbb{R}$ is a continuous function defined on an unbounded interval $I$.  This is because a non-constant polynomial would approach $\infty$ or $-\infty$ when $x$ approaches $\infty$ or $-\infty$. However, there are bounded continuous functions defined on unbounded intervals. For example, the function
\[f(x)=\frac{x}{x^2+1}\] is a bounded continuous function defined on $\mathbb{R}$.

\begin{remark}
{}

In probability theory, a binomial random variable $X$ with parameters $n$ and $p$ counts the number of successes in $n$ independent and identical Bernoulli trials, each has a probability $p\in (0,1)$ of being a success.  $X$ can take integer values between $0$ and $n$. The probability that $X=k$ is
\[P(X=k)=\binom{n}{k}p^k(1-p)^{n-k},\quad 0\leq k\leq n.\]
The identity in (a) of Lemma \ref{230309_14} amounts to
\[\sum_{k=0}^n \binom{n}{k}p^k(1-p)^k=1,\]
\end{remark}\begin{highlight}{}
which reflects that the total probability is 1. The identity in part (b) gives 
\[E(X)=\sum_{k=0}^n k\binom{n}{k}p^k(1-p)^k=np,\]

which is the expected value of a binomial  random variable $X$ with parameters $n$ and $p$. The identity in part (c) gives
\[E(X^2)=\sum_{k=0}^n k^2\binom{n}{k}p^k(1-p)^k=n^2p^2+np(1-p).\] Together with the identity in part (b),  the variance of $X$ is given by
\begin{align*}\text{Var}\,(X)&=E(X^2)-E(X)^2 =n^2p^2+np(1-p)-n^2p^2=np(1-p).\end{align*}In fact, the variance of a random variable $X$ is defined as
\[\text{Var}\,(X)=E([X-E(X)]^2).\]The identity in part (d) of Lemma \ref{230309_14} is just another way of computing the variance.  Using part (d), we have
\begin{align*}\text{Var}\,(X)&= \sum_{k=0}^n (k-np)^2\binom{n}{k}p^k(1-p)^k\\&=n^2 \sum_{k=0}^n \left(\frac{k}{n}-p\right)^2\binom{n}{k}p^k(1-p)^k=np(1-p).\end{align*}

\end{highlight}

\backmatter

\chapter*{References}
\bibliographystyle{amsalpha}
\bibliography{ref}

\providecommand{\bysame}{\leavevmode\hbox to3em{\hrulefill}\thinspace}
\providecommand{\MR}{\relax\ifhmode\unskip\space\fi MR }
\providecommand{\MRhref}[2]{%
  \href{http://www.ams.org/mathscinet-getitem?mr=#1}{#2}
}
\providecommand{\href}[2]{#2}
\begin{thebibliography}{SCW20}

\bibitem[Abb15]{Abbott}
Stephen Abbott, \emph{Understanding analysis}, second ed., Undergraduate Texts
  in Mathematics, Springer, New York, 2015. \MR{3331079}

\bibitem[Apo74]{Apostol}
Tom~M. Apostol, \emph{Mathematical analysis}, second ed., Addison-Wesley
  Publishing Co., Reading, Mass.-London-Don Mills, Ont., 1974. \MR{0344384}

\bibitem[BS92]{Bartle}
Robert~G. Bartle and Donald~R. Sherbert, \emph{Introduction to real analysis},
  second ed., John Wiley \& Sons, Inc., New York, 1992. \MR{1135107}

\bibitem[Fit09]{Fitzpatrick}
Patrick~M. Fitzpatrick, \emph{Advanced calculus}, second ed., American
  Mathematical Society, 2009.

\bibitem[Ros18]{Rosen}
Kenneth Rosen, \emph{Discrete mathematics and its applications}, eighth ed., Mc
  Graw Hill, 2018.

\bibitem[Rud76]{Rudin}
Walter Rudin, \emph{Principles of mathematical analysis}, third ed.,
  International Series in Pure and Applied Mathematics, McGraw-Hill Book Co.,
  New York-Auckland-D\"{u}sseldorf, 1976. \MR{0385023}

\bibitem[SCW20]{Stewart}
James Stewart, Daniel~K. Clegg, and Saleem Watson, \emph{Calculus}, ninth ed.,
  Cengage Learning, 2020.

\bibitem[Tao14]{Tao_2}
Terence Tao, \emph{Analysis. {II}}, third ed., Texts and Readings in
  Mathematics, vol.~38, Hindustan Book Agency, New Delhi, 2014. \MR{3310023}

\bibitem[Tao16]{Tao_1}
\bysame, \emph{Analysis. {I}}, third ed., Texts and Readings in Mathematics,
  vol.~37, Hindustan Book Agency, New Delhi; Springer, Singapore, 2016,
  Edectronic edition of [ MR3309891]. \MR{3728289}

\bibitem[Zor15]{Zorich_1}
Vladimir~A. Zorich, \emph{Mathematical analysis. {I}}, second ed.,
  Universitext, Springer-Verlag, Berlin, 2015, With Appendices A--F and new
  problems translated by Octavio Paniagua T. \MR{3495809}

\bibitem[Zor16]{Zorich_2}
\bysame, \emph{Mathematical analysis. {II}}, second ed., Universitext,
  Springer, Heidelberg, 2016. \MR{3445604}

\end{thebibliography}
 
\begin{coverpage}
~
 \end{coverpage}
\end{document}